# AXIOMATIC THEORY OF DISTRIBUTIONS

— by —

Newton C.A. da Costa

and

J.A. Baêta Segundo

*To Neusa and Téia, with profound gratitude.*

# CONTENTS

# PROLOGUE

Very frequently we find, for instance, in physics books, an incursion on the theory of distributions, especially Dirac's delta distribution, quite unsatisfactory from a mathematical point of view. The main reason for this, in our view, resides in the fact that a constructive approach to distributions, at the heart of mathematical analysis, requires the consideration of many questions of topological nature that, in a certain sense, are far away from those of genuine interest of physicists, engineers, and other mathematics users. However, opting up for a non-constructive, but axiomatic, approach, it is possible to capture the referred concept by means of simple axioms with which, through a process of logical inference, its properties can be established, hence obtaining a rigorous treatment from the mathematical point of view and, at the same time, appealing from a physicist perspective, for instance. Something analogous to consider the system of real numbers as defined by the axioms of complete ordered field, forgetting about the constructions (models) made by Dedekind, Cantor and others.

Regarding distributions, the idea can be expressed, *mutatis mutandis*, by what said Spivak*:

> "It is utterly irrelevant that a real number happens to be a collection of rational numbers, and such a fact should never enter the proof of any important theorem about real numbers. Reasonable proofs should use only the fact that the real numbers are a complete ordered field, because this property of the real numbers characterizes them up to isomorphism, and any significant mathematical property of the real numbers will be true for all isomorphic fields."

José Sebastião e Silva (1914 - 1972), influential Portuguese mathematician with a relevant scientific work integrated into the evolution of functional analysis in the post-war, published in 1954, in the Revista da Faculdade de Ciências de Lisboa, the important work entitled "Sur Une Construction Axiomatique de la Théorie des Distributions"†. In this article, J.S. e Silva, among other things, shows that the statement made in the previous

* Spivak, M. *Calculus.* W. A. Benjamin, 1967.

† Silva, J. S. *Sur Une Construction Axiomatique de la Théorie des Distributions.* Rev. da Fac. de Ciências de Lisboa, 2ª série-A, vol 4, pp. 79–186, (1954–1955).

paragraph, about the possibility of a simple axiomatic capturing the concept of distribution, is true. In the referred article, J.S. e Silva presents an axiomatic foundation for the L. Schwartz' distributions theory, backed by concepts as simple as the ones of continuity and differentiability studied in elementary Calculus courses. Such as the axioms of complete ordered field, J.S. e Silva's ones admit, essentially, a single model, namely, the distributions; more precisely, Silva's axiomatic is categoric and has as a model the classic structure of the distributions as constructed by L. Schwartz.

When elaborating this monograph, we not only used the work of J.S. e Silva above mentioned, but were also motivated and guided by it to introduce certain concepts and to promote some generalizations, described synthetically in (a) to (e) ahead, that allows, at least in principle, new applications. Our approach is composed, basically, of the following points:

**(a)** First, we define an axiomatic structure — $\mathscr{G}(I)$ — that, in honor to Sebastião e Silva, we denominate $S$-space, and we formulate a notion of extension of a $S$-space $\mathscr{G}(I)$ consisting of, roughly speaking, another $S$-space, $\widehat{\mathscr{G}}(I)$, that in a certain way includes $\mathscr{G}(I)$ properly;

**(b)** next, we prove two "theorems of extension" that ensure, for $S$-spaces $\mathscr{G}(I)$ of a certain type, let us say, of type A, the existence of two extensions, $\widetilde{\mathscr{G}}(I)$ and $\overline{\mathscr{G}}(I)$, both essentially unique;

**(c)** starting with the extensions $\widetilde{\mathscr{G}}(I)$ and $\overline{\mathscr{G}}(I)$ (of $\mathscr{G}(I)$) we define the notions of $\widetilde{\mathscr{G}}(I)$ and $\overline{\mathscr{G}}(I)$-distributions and its associated concepts of domain, addition, derivative and restriction;

**(d)** we prove that the $\widetilde{\mathscr{G}}(I)$ and $\overline{\mathscr{G}}(I)$-distributions (and associated concepts) can be defined through simple axiomatics and we formulate, for each case, a pair of equivalent axiomatics;

**(e)** finally, we prove that the two pairs of axiomatics referred to in (d), specialized (particularized) to a specific $S$-space $\mathscr{G}(I)$, i.e., to a given model of the axiomatic structure of $S$-space, a model involving the continuous functions and their partial derivatives, categorically define the (*à la* Schwartz) finite order distributions, one of the pairs, and the distributions (of all orders, finite and infinite), the other pair.

In few words, we elaborated a kind of "axiomatization schema" for the $\widetilde{\mathscr{G}}(I)$ and $\overline{\mathscr{G}}(I)$-distributions: for each model $\mathscr{G}(I)$ of the $S$-space structure of type A, the schema provides categoric axiomatics that have as models the $\widetilde{\mathscr{G}}(I)$ and $\overline{\mathscr{G}}(I)$-distributions correspondent to $\mathscr{G}(I)$. In particular, axiomatics for the Schwartz' distributions are obtained through this schema.

Hence, the structure of $S$-space and its corresponding theory, formulated in elementary algebraic basis, not only allows obtaining the results of J. Sebastião e Silva as a particular case but, mainly, due to its abstract nature, it explicits a structural environment conceptually very simple where other generalizations, beyond the one promoted to differential calculus by the theory of distributions, can, in principle, stand.

Regarding the prerequisites to study this work, nothing beyond some basic notions associated with groups and semigroups theory — homomorphisms, endomorphisms, isomorphism between groups and between semigroups — to linear algebra — vector spaces and subspaces, linear transformations — to analysis in $\mathbb{R}^n$ — open, closed and compact sets, continuous functions, partial derivative and integration of continuous functions — and general notions about axiomatic theories — model, categoricity, equivalence of systems of axioms.

We would like to thank Guilherme Tavares da Silva for the translation of the original text in Portuguese to English, and its LaTeX digitization, as well as all others that helped us with their critics and suggestions, especially to professor and friend Luis Fernando Mello for the priceless support and insightful analysis of the text.

Newton C.A. da Costa
J.A. Baêta Segundo

Itajubá, October 2021.

# 1

# DISTRIBUTIONS AND $S$-GROUPS

## Introduction

### 1.1 Motivation

The mathematical modeling of fragments of physical reality leads, frequently, to differential equations. For instance, under certain conditions, the wave equation

$$\frac{\partial^2 u}{\partial t^2} - k^2 \frac{\partial^2 u}{\partial x^2} = 0,$$

and Laplace's equation,

$$\nabla^2 u = 0,$$

where

$$\nabla^2 u = \frac{\partial^2 u}{\partial x^2} + \frac{\partial^2 u}{\partial y^2} \quad \text{in} \quad \mathbb{R}^2$$

and

$$\nabla^2 u = \frac{\partial^2 u}{\partial x^2} + \frac{\partial^2 u}{\partial y^2} + \frac{\partial^2 u}{\partial z^2} \quad \text{in} \quad \mathbb{R}^3,$$

well describe the physical situations of a vibrant string (tense string, along the $x$-axis, vibrating on the $xy$-plane) and the electrostatic potential of a region free of electrical charges, respectively.

Not rare, the own portion of the reality modeled by a differential equation suggests that the latter admits, as solutions, non-differentiable functions — a surprising situation. For instance[1]: it is known that, being $f(\xi)$, $\xi \in \mathbb{R}$, any real function twice differentiable in $\mathbb{R}$, $u(x,t) = f(x - kt)$ is a solution to the wave equation describing, mathematically, the physical situation of a pulse, whose form is given by the function $f(\xi)$, propagating through the string, without deformation, at a speed $k$. We believe that a physicist would like to consider the propagation of a pulse whose form is described by a continuous function $f(\xi)$, not differentiable in a certain finite number of points. However, regardless of physical

[1] Strichartz, R. S. *A Guide to Distribution Theory and Fourier Transforms.* CRC Press, 1993.

arguments, with such a form $f$, the function $u(x,t) = f(x-kt)$ is not a solution to the wave equation, since it does not satisfy it on points where $f$ is not differentiable. As significant physical possibilities should not be discarded for not meeting requirements related to technical aspects of the model, two paths are immediately glimpsed:

**via A** – the shorter, actually too short as will be seen, consists of modifying the concept of solution to a differential equation in such a way that functions satisfying it except in a certain "reduced" number of points, at which they are not differentiable, are allowed;

**via B** – the second one, longer, modifies the own concept of derivative, extending it to a bigger class of objects than that of differentiable functions, bigger in the sense of proper inclusion, in such a way that the new model, that is, the differential equation with the derivatives now taken in this modified, extended sense, also admits as solutions the physically significative ones which are not allowed in the classical model.

The simple modification, as described in via A, of the concept of solution to a differential equation would accommodate relevant physical situations discarded by the usual concept of solution as, for instance, that of a pulse with form $f$ continuous and not differentiable propagating through a string, but, on the other hand, would also open space for so many others to which physical reasons exist to discard then. In fact, the functions

$$u(x,y) = \log\left(x^2 + y^2\right)$$

and

$$u(x,y,z) = \left(x^2 + y^2 + z^2\right)^{-\frac{1}{2}}$$

which satisfy, respectively, the Laplace equation in two and three dimensions at every point except the origin, where they are not differentiable, would be, therefore, in the sense of via A, considered as solutions for the referred equations, despite the fact of existing physical reasons to reject them: the electrostatic potential in regions free of electrical charges are continuous, while the functions above do have poles.

As for the second path, via B, before start exploring it, it is worth remembering some results from analysis as well as introduce some notation.

## 1.2 A Brief Review

In what follows, we highlight some results from mathematical analysis which we will use in our considerations, as well as introduce some notation.

**(a)** The set of natural numbers, that is, non-negative integers, will be denoted by $\mathbb{N}$, whereas $\mathbb{N}^n$ symbolizes the set of all $n$-tuples of natural numbers;

**(b)** $\mathbb{R}$ represents a complete ordered field, namely, the set of real numbers;

**(c)** $\mathbb{C}$ denotes the field of complex numbers;

**(d)** $\mathbb{R}^n$ represents the vector space of the $n$-tuples of real numbers, equipped with the euclidean metric:

$$d(x, y) = \sqrt{\sum_{i=1}^{n} (x_i - y_i)^2}$$

with $x = (x_1, \dots, x_n)$ and $y = (y_1, \dots, y_n)$ in $\mathbb{R}^n$;

**(e)** $C(\Omega)$, where $\Omega \subseteq \mathbb{R}^n$ is an open set, symbolizes the class of functions with domain $\Omega$, assuming complex values, which are continuous on $\Omega$;

**(f)** $C^k(\Omega)$, with $k \in \mathbb{N}$, denotes the class of functions $f \in C(\Omega)$ such that every partial derivatives of order less or equal to $k$ exist and are continuous functions in $\Omega$. For $k = 0$ we write $C^0(\Omega) = C(\Omega)$. Hence, for instance, $f \in C^2(\Omega)$ if the partial derivatives,

$$\begin{aligned} \frac{\partial f}{\partial x_i} : \Omega &\longrightarrow \mathbb{C} \\ x &\longmapsto \frac{\partial f}{\partial x_i}(x) \end{aligned}$$

for $i \in \{1, \dots, n\}$, are differentiable functions on $\Omega$ and their partial derivatives (the order 2 partial derivatives of $f$),

$$\begin{aligned} \frac{\partial}{\partial x_j}\left(\frac{\partial f}{\partial x_i}\right) : \Omega &\longrightarrow \mathbb{C} \\ x &\longmapsto \left(\frac{\partial}{\partial x_j}\left(\frac{\partial f}{\partial x_i}\right)\right)(x) \end{aligned}$$

for $i, j \in \{1, \dots, n\}$, are continuous on $\Omega$.

The use of multi-indexes, as defined in (h), allows, as we will see, a concise notation for partial derivatives of high orders;

**(g)** $C^\infty(\Omega)$ symbolizes the set of functions $f \in C(\Omega)$ with continuous partial derivatives of every order in $\Omega$. Hence, $f \in C^\infty(\Omega)$ if and only if $f \in C^k(\Omega)$ for every $k \in \mathbb{N}$;

**(h)** $\alpha$ is a **multi-index** or, more precisely, a $n$-multi-index if and only if $\alpha \in \mathbb{N}^n$. For multi-indexes $\alpha = (\alpha_1, \dots, \alpha_n)$ and $\beta = (\beta_1, \dots, \beta_n)$ we define:

$$\alpha + \beta = (\alpha_1 + \beta_1, \dots, \alpha_n + \beta_n)$$

and

$$|\alpha| = \alpha_1 + \cdots + \alpha_n;$$

**(i)** Let $\alpha = (\alpha_1, \ldots, \alpha_n)$ be a multi-index and $f \in C(\Omega)$ be a function for which exist (in $\Omega$) the partial derivative (of order $|\alpha|$)

$$\frac{\partial^{\alpha_1}}{\partial x_1^{\alpha_1}}\left(\frac{\partial^{\alpha_2}}{\partial x_2^{\alpha_2}}\left(\cdots\left(\frac{\partial^{\alpha_n} f}{\partial x_n^{\alpha_n}}\right)\cdots\right)\right).$$

We denote this partial derivative of $f$ as $\partial_\Omega^\alpha(f)$. Clearly, if $f \in C^{|\alpha|}(\Omega)$, then we have $\partial_\Omega^\alpha(f) \in C(\Omega)$ and, hence, we have well defined the function

$$\begin{aligned} \partial_\Omega^\alpha : C^{|\alpha|}(\Omega) &\longrightarrow C(\Omega) \\ f &\longmapsto \partial_\Omega^\alpha(f) = \tfrac{\partial^{\alpha_1}}{\partial x_1^{\alpha_1}}\Big(\cdots\Big(\tfrac{\partial^{\alpha_n} f}{\partial x_n^{\alpha_n}}\Big)\cdots\Big). \end{aligned}$$

For example, $\partial_{\mathbb{R}^3}^{(2,0,3)}$ represents the function

$$\begin{aligned} \partial_\Omega^{(2,0,3)} : C^5(\Omega) &\longrightarrow C(\Omega) \\ f &\longmapsto \partial_\Omega^{(2,0,3)}(f) = \tfrac{\partial^2}{\partial x_1^2}\Big(\tfrac{\partial^3 f}{\partial x_3^3}\Big). \end{aligned}$$

If $\alpha = 0 = (0, \ldots, 0)$ we write $\partial_\Omega^0 = I_{C(\Omega)}$, where $I_{C(\Omega)}$ is the identity function on $C(\Omega)$:

$$\begin{aligned} I_{C(\Omega)} : C(\Omega) &\longrightarrow C(\Omega) \\ f &\longmapsto I_{C(\Omega)}(f) = f; \end{aligned}$$

**(j)** From the mathematical analysis, we have that:

$$\partial_\Omega^\alpha\Big(\partial_\Omega^\beta(f)\Big) = \partial_\Omega^\beta\Big(\partial_\Omega^\alpha(f)\Big) = \partial_\Omega^{\alpha+\beta}(f)$$

for every $f \in C^\infty(\Omega)$, whatever the multi-indexes $\alpha$ and $\beta$; in other words, $\partial_\Omega^\alpha$ and $\partial_\Omega^\beta$ commute in $C^\infty(\Omega)$;

Beyond that, if $k \in \mathbb{N}$ and the multi-indexes $\alpha$ and $\beta$ are such that $|\alpha + \beta| \leqslant k$, then $\partial_\Omega^\alpha$ and $\partial_\Omega^\beta$ also commute in $C^k(\Omega)$;

**(k)** $C^k(\Omega)$, for every $k \in \mathbb{N}$, with the usual operation of addition of functions is an abelian group and, as $C^k(\Omega) \subseteq C(\Omega)$, then $C^k(\Omega)$ is a subgroup of the abelian group $C(\Omega)$. When, in this monograph, we refer to $C^k(\Omega)$ as a group, we will be considering it equipped with the usual operation of addition of functions;

**(l)** We highlight that $\partial_\Omega^\alpha\colon C^{|\alpha|}(\Omega) \longrightarrow C(\Omega)$ is a surjective homomorphism (between groups), that is, if $f \in C(\Omega)$, then there exists $g \in C^{|\alpha|}(\Omega)$ such that

$$\partial_\Omega^\alpha(g) = f;$$

**(m)** Finally, the symbol "$A \coloneqq B$" will be employed to denote that: $A$ is equal, by definition, to $B$.

## 1.3 The via B and the Distributions

The endeavor to free differential calculus from difficulties as those presented in 1.1, which occur due to the existence of continuous functions which are not differentiable, following the route of via B, culminates, as we will see, in the theory of distributions, which extends the notion of differentiability to a class of objects, the distributions, which properly includes the class of continuous functions and, therefore, of the differentiable ones in the usual sense. Naturally, of any such extension it is expected, in order for it to be useful, that:

- every continuous function to be a distribution;
- distributions to be differentiable and, its derivatives, to be also distributions, therefore, distributions are infinitely differentiable;
- the notion of derivative of a distribution, when applied to classically differentiable functions, to coincide with the usual derivative; and
- the usual rules of differential calculus to remain valid.

We can say that the via B poses an extension problem: given the “embryonic universe”

$$\left(C(\Omega), \partial(\Omega) = \left\{\partial^{\alpha}_{\Omega} \colon \ \alpha \in \mathbb{N}^n\right\}\right)$$

constituted of continuous functions and partial derivatives, we ask for expand it so that to obtain another one,

$$\left(\widetilde{C}(\Omega), \widetilde{\partial}(\Omega) = \left\{\widetilde{\partial^{\alpha}_{\Omega}} \colon \ \alpha \in \mathbb{N}^n\right\}\right),$$

at which the elements of $\widetilde{C}(\Omega)$ would be the distributions and those of $\widetilde{\partial}(\Omega)$ the generalized derivatives, satisfying the conditions above rewritten below in a more convenient form.

**(a)** $C(\Omega) \subseteq \widetilde{C}(\Omega)$;

**(b)** For each $\alpha \in \mathbb{N}^n$, $\widetilde{\partial^{\alpha}_{\Omega}}$ is a function with domain and codomain equal to $\widetilde{C}(\Omega)$, $\widetilde{\partial^{\alpha}_{\Omega}} \colon \widetilde{C}(\Omega) \longrightarrow \widetilde{C}(\Omega)$;

**(c)** For each $\alpha \in \mathbb{N}^n$, the restriction of the function $\widetilde{\partial^{\alpha}_{\Omega}}$ to the subset $C^{|\alpha|}(\Omega)$ of its domain $\widetilde{C}(\Omega)$, is the partial derivative $\partial^{\alpha}_{\Omega}$, that is,

$$\widetilde{\partial^{\alpha}_{\Omega}}(f) = \partial^{\alpha}_{\Omega}(f)$$

for every $f \in C^{|\alpha|}(\Omega)$;

**(d)** For any $\alpha, \beta \in \mathbb{N}^n$ and every $\widetilde{f}, \widetilde{g} \in \widetilde{C}(\Omega)$,

$$\widetilde{\partial_\Omega^\alpha}\left(\widetilde{f}+\widetilde{g}\right) = \widetilde{\partial_\Omega^\alpha}\left(\widetilde{f}\right) + \widetilde{\partial_\Omega^\alpha}\left(\widetilde{g}\right)$$

and

$$\widetilde{\partial_\Omega^\alpha}\left(\widetilde{\partial_\Omega^\beta}\left(\widetilde{f}\right)\right) = \widetilde{\partial_\Omega^{\alpha+\beta}}\left(\widetilde{f}\right).$$

It is necessary, in order for conditions (a) to (d) above to make sense, to take $C(\Omega)$ not only as a set but as the abelian group $C(\Omega)$ described in 1.2(k); in the same way, it is necessary to consider the basic set, $\widetilde{C}(\Omega)$, from the universe of distributions, equipped with a binary operation that coincides with the addition of the group $C(\Omega)$ when restricted to it. We will admit that $\widetilde{C}(\Omega)$ with this operation is also an abelian group which, therefore, has $C(\Omega)$ as a subgroup (condition (a)).

And what to say about the components $\partial(\Omega)$ and $\widetilde{\partial}(\Omega)$ of the universes $(C(\Omega), \partial(\Omega))$ and $(\widetilde{C}(\Omega), \widetilde{\partial}(\Omega))$? What type of structure the conditions (a) to (d) impose to $\partial(\Omega)$ and $\widetilde{\partial}(\Omega)$? We remark that, for both $\partial(\Omega)$ and $\widetilde{\partial}(\Omega)$, once met the conditions (a) to (d), multiplication operations are then defined as

$$\left(\partial_\Omega^\alpha \partial_\Omega^\beta\right)(f) := \partial_\Omega^\alpha\left(\partial_\Omega^\beta(f)\right) \quad \text{for every} \quad f \in C^{|\alpha+\beta|}(\Omega)$$

and

$$\left(\widetilde{\partial_\Omega^\alpha}\widetilde{\partial_\Omega^\beta}\right)(\widetilde{f}) := \widetilde{\partial_\Omega^\alpha}\left(\widetilde{\partial_\Omega^\beta}(\widetilde{f})\right) \quad \text{for every} \quad \widetilde{f} \in \widetilde{C}(\Omega)$$

which are binary in $\partial(\Omega)$ and $\widetilde{\partial}(\Omega)$, respectively (since $\partial_\Omega^\alpha \partial_\Omega^\beta = \partial_\Omega^{\alpha+\beta}$ and $\widetilde{\partial_\Omega^\alpha}\widetilde{\partial_\Omega^\beta} = \widetilde{\partial_\Omega^{\alpha+\beta}}$), and associative. Therefore, $\partial(\Omega)$ as well as $\widetilde{\partial}(\Omega)$ are, with its corresponding multiplications, semigroups: the first, $\partial(\Omega)$, a semigroup of surjective homomorphisms of subgroups of $C(\Omega)$ $\left(C^{|\alpha|}(\Omega)\right)$ over $C(\Omega)$, which contains the identity on $C(\Omega)$ ($\partial_\Omega^0 = I_{C(\Omega)}$), while the second, $\widetilde{\partial}(\Omega)$, a semigroup of endomorphisms on $\widetilde{C}(\Omega)$.

In short, the task of improving, through via B, the mathematical apparatus, aiming to obtain mathematical models of physical reality more suitable to the latter, which does not exclude significant physical possibilities, leads us naturally to the structure $(C(\Omega), \partial(\Omega))$ constituted by the abelian group of continuous functions, $C(\Omega)$, and the semigroup $\partial(\Omega)$ of partial derivatives, and to the following associated extension problem:

> To extend the structure $(C(\Omega), \partial(\Omega))$ obtaining another one, $(\widetilde{C}(\Omega), \widetilde{\partial}(\Omega))$, the habitat of the distributions, satisfying the conditions (a) to (d) above formulated.

A precise formulation of the problem above is given in the next section, through concepts whose definitions allow us to formalize the intuitive and informal approach adopted until now.

# $S$-Groups and the Extension Problem

## 1.4 Definition

**(a)** We say that the ordered pair

$$\mathbb{G} = (G, H)$$

is a $S$-**group** if and only if $G$ is an abelian group and $H$, with the multiplication defined ahead, is a semigroup of homomorphisms of subgroups of $G$ into $G$.

**Multiplication in** $H$ — Given the homomorphisms $\Phi, \Psi \in H$, we denote as $G_\Phi \subseteq G$ and $G_\Psi \subseteq G$ the subgroups of $G$ which are, respectively, the domains of $\Phi$ and $\Psi$, and as $G_{\Phi\Psi}$ the subset of $G$ defined as

$$G_{\Phi\Psi} := \Big\{ g \in G \mid g \in G_\Psi \quad \text{and} \quad \Psi(g) \in G_\Phi \Big\}.$$

The multiplication is then defined as an operation that to each ordered pair $(\Phi, \Psi)$ of elements of $H$ associates a function denoted as $\Phi\Psi$, denominated the product of $\Phi$ by $\Psi$, and defined as

$$\begin{aligned} \Phi\Psi : G_{\Phi\Psi} &\longrightarrow G \\ g &\longmapsto (\Phi\Psi)(g) := \Phi\Big(\Psi(g)\Big). \end{aligned}$$

It is easy to verify that $G_{\Phi\Psi}$ is, as $G_\Phi$ and $G_\Psi$, a subgroup of $G$ and that $\Phi\Psi$ is, as $\Phi, \Psi \in H$, a homomorphism from a subgroup of $G$ ($G_{\Phi\Psi}$) into $G$; also, the multiplication as defined is associative and hence $H$ is a semigroup, if and only if, $\Phi\Psi \in H$ whatever $\Phi, \Psi \in H$.

**(b)** We say that the $S$-group $\mathbb{G} = (G, H)$ is **surjective** if and only if each homomorphism $\Phi \in H$ is surjective; $\mathbb{G}$ is **abelian** if and only if the semigroup $H$ is abelian, i.e., $\Phi\Psi = \Psi\Phi$ whatever $\Phi, \Psi \in H$; $\mathbb{G}$ is a $S$-group **with identity** if and only if the homomorphism $I_G$, the identity on $G$, belongs to $H$.

**(c)** Let $\mathbb{G} = (G, H)$ and $\mathbb{E} = (E, F)$ be $S$-groups. We say that $\mathbb{G}$ is **isomorphic** to $\mathbb{E}$ if and only if there exist isomorphisms $\alpha\colon G \longrightarrow E$ and $\beta\colon H \longrightarrow F$ between groups and semigroups, respectively, such that

$$\alpha\Big(\Phi(g)\Big) = \beta(\Phi)\Big(\alpha(g)\Big)$$

for every $\Phi \in H$ and $g \in G_\Phi$ ($G_\Phi \subseteq G$ denotes, as in (a) above, the domain of $\Phi$).

**Remark.** It is easy to verify that if $\mathbb{G}$ is isomorphic to $\mathbb{E}$, then $\mathbb{E}$ is isomorphic to $\mathbb{G}$ and hence, without ambiguity, we can state that $\mathbb{G}$ and $\mathbb{E}$ are isomorphic $S$-groups; it is also easy to verify that if $\mathbb{G}$, $\mathbb{E}$ and $\mathbb{L}$ are $S$-groups such that $\mathbb{G}$ is isomorphic to $\mathbb{E}$ and $\mathbb{E}$ is isomorphic to $\mathbb{L}$, then $\mathbb{G}$ is isomorphic to $\mathbb{L}$. Yet, it can be seen that isomorphic $S$-groups cannot be distinguished through their properties and, hence, we can identify them.

## 1.5 Example: the $S$-group $\mathbb{C}(\Omega)$

Our considerations in 1.3, relative to the question of obtaining a more flexible notion of differentiability than the usual, classical, at which continuous functions would be infinitely differentiable, and which reduces to the classical notion when applied to the the domain of the latter, led us naturally to the following $S$-group, henceforth referred to as **$S$-group of continuous functions**:

$$\mathbb{C}(\Omega) := \Big(C(\Omega), \partial(\Omega)\Big)$$

where $C(\Omega)$ is the abelian group of the complex values functions, whose domain is an open set $\Omega \subseteq \mathbb{R}^n$, which are continuous on $\Omega$, and $\partial(\Omega)$ is the semigroup of partial derivatives, that is (see 1.2)

$$\partial(\Omega) := \Big\{\partial^\alpha_\Omega : \ \alpha \in \mathbb{N}^n\Big\}$$

with the multiplication operation defined by

$$\partial^\alpha_\Omega \partial^\beta_\Omega := \partial^{\alpha+\beta}_\Omega$$

whatever the multi-indexes $\alpha, \beta \in \mathbb{N}^n$. Taking into account the review provided in 1.2, in particular the itens 1.2(i), (j) and (l), we have that $\partial(\Omega)$ is an abelian semigroup ($\partial^\alpha_\Omega \partial^\beta_\Omega = \partial^\beta_\Omega \partial^\alpha_\Omega$), surjective (if $f \in C(\Omega)$, then there exists $g \in C^{|\alpha|}_\Omega$ such that $\partial^\alpha_\Omega(g) = f$) and with identity ($\partial^0_\Omega = I_{C(\Omega)}$). Then, remembering the Definition 1.4(b), $\mathbb{C}(\Omega) = (C(\Omega), \partial(\Omega))$ is an abelian, surjective, and with identity $S$-group.

## 1.6 Definition

**(a)** Let $G, \widehat{G}, \widehat{G'}$ and $\delta$ be such that: $G$ and $\widehat{G}$ are groups, $\widehat{G'} \subseteq \widehat{G}$ is a subgroup of $\widehat{G}$ and $\delta$ is an isomorphism of the group $G$ onto the group $\widehat{G'}$ ($\delta\colon G \longrightarrow \widehat{G'}$). Under these conditions, we say that $\widehat{G}$ **admits** or, more precisely, **$\delta$-admits** $G$ as a subgroup;

**(b)** Let $G$ and $\widehat{G}$ be groups such that $\widehat{G}$ $\delta$-admits $G$ as a subgroup. Let also $H = \{\Phi, \Psi, \ldots\}$ be a semigroup of homomorphisms of subgroups of $G$ into $G$. For $\Phi \in H$, $\Phi\colon G_\Phi \longrightarrow G$, we say that $\widehat{\Phi}\colon \widehat{G} \longrightarrow \widehat{G}$ is a **$\delta$-prolongation of $\Phi$ to $\widehat{G}$** if and only if $\widehat{\Phi}$ is a homomorphism and, therefore, an endomorphism on $\widehat{G}$, such that

$$\widehat{\Phi}\Big(\delta(g)\Big) = \delta\Big(\Phi(g)\Big) \quad \text{for every} \quad g \in G_\Phi.$$

Since $\widehat{G}$ $\delta$-admits $G$ as a subgroup, then

$$\delta(G) := \Big\{\delta(g) : \ g \in G\Big\}$$

and

$$\delta(G_\Phi) := \Big\{\delta(g) : \ g \in G_\Phi\Big\}$$

are subgroups of $\widehat{G}$ isomorphic, respectively, to $G$ and $G_\Phi$, which allow us to make the following identifications[2]:

$$G \equiv \delta(G), \quad G_\Phi \equiv \delta(G_\Phi) \quad \text{and} \quad g \equiv \delta(g) \quad \text{for every} \quad g \in G.$$

With these identifications in mind, we write, as an abuse of language, for example,

$$\widehat{\Phi}(g) = \Phi(g) \quad \text{for every} \quad g \in G_\Phi,$$

instead of the correct form,

$$\widehat{\Phi}\Big(\delta(g)\Big) = \delta\Big(\Phi(g)\Big) \quad \text{for every} \quad g \in G_\Phi,$$

and then, based on this abuse, we say that $\Phi$ is the restriction of $\widehat{\Phi}$ to $G_\Phi$ and that the endomorphism on $\widehat{G}$, $\widehat{\Phi}$, is a prolongation of $\Phi$ to $\widehat{G}$ if and only if $\Phi$ is its restriction to $G_\Phi$, without any reference to the isomorphism $\delta$;

**(c)** With the same hypotheses of (b), we say that $\widehat{H} = \{\widehat{\Phi}, \widehat{\Psi}, \cdots\}$ is a **$\delta$-prolongation of $H = \{\Phi, \Psi, \cdots\}$ to $\widehat{G}$** if and only if its elements are, exact and precisely, $\delta$-prolongations to $\widehat{G}$ of each $\Phi \in H$ such that:

$$\widehat{\Phi\Psi} = \widehat{\Phi}\widehat{\Psi} \quad \text{for every} \quad \Phi, \Psi \in H,$$

with the multiplication $\widehat{\Phi}\widehat{\Psi}$ in $\widehat{H}$ taken as the usual composition of functions. It immediately results that every $\delta$-prolongation of the semigroup $H$ is also a semigroup (of endomorphisms on $\widehat{G}$), and if $H$ is abelian then its $\delta$-prolongations also are.

In other terms, a $\delta$-prolongation of the semigroup $H$ to $\widehat{G}$ is a semigroup $\widehat{H}$ of endomorphisms on $\widehat{G}$ such that:

**(c-1)** is isomorphic to $H$, i.e., there exist a bijective function from $H$ onto $\widehat{H}$, denoted

$$\widehat{\ } : H \longrightarrow \widehat{H}$$
$$\Phi \longmapsto \widehat{\Phi}$$

such that

$$\widehat{\Phi\Psi} = \widehat{\Phi}\widehat{\Psi} \quad \text{for every} \quad \Phi, \Psi \in H;$$

**(c-2)** $\widehat{\Phi}\Big(\delta(g)\Big) = \delta\Big(\Phi(g)\Big)$ for every $\Phi \in H$ and $g \in G_\Phi$;

**(d)** Let $\mathbb{G} = (G, H)$ and $\widehat{\mathbb{G}} = (\widehat{G}, \widehat{H})$ be $S$-groups. We say that $\widehat{\mathbb{G}}$ is a **$\delta$-extension of $\mathbb{G}$** if and only if $\widehat{G}$ $\delta$-admits $G$ as a subgroup and $\widehat{H}$ is a $\delta$-prolongation of $H$ to $\widehat{G}$. It then results that $\delta$-extensions of an abelian $S$-group are, themselves, abelian $S$-groups.

[2] The symbol $A \equiv B$ informs that, in a certain sense, the object $A$ is identified to the object $B$.

## 1.7 Remark

The abuse of language mentioned in Definition 1.6(b), at which instead of "$\widehat{\Phi}$ is a $\delta$-prolongation of $\Phi$ to $\widehat{G}$", we omit the reference to the isomorphism $\delta$ and say, simply, "$\widehat{\Phi}$ is a prolongation of $\Phi$ to $\widehat{G}$", will be frequently used in this monograph, not only with respect to the concept of $\delta$-prolongation, but also relative to all those which, as this, encompass the isomorphism $\delta$. Hence, being $\mathbb{G} = (G, H)$ and $\widehat{\mathbb{G}} = (\widehat{G}, \widehat{H})$ $S$-groups, we say, for example, that "$\widehat{G}$ admits $G$", "$\widehat{H}$ is a prolongation of $H$ to $\widehat{G}$", "$\widehat{\mathbb{G}}$ is an extension of $\mathbb{G}$", instead of, respectively, "$\widehat{G}$ $\delta$-admits $G$", "$\widehat{H}$ is a $\delta$-prolongation of $H$ to $\widehat{G}$", "$\widehat{\mathbb{G}}$ is a $\delta$-extension of $\mathbb{G}$". However, in some future developments it will be convenient to get back to the precise language.

## 1.8 The Extension Problem

The problem which via B led us, described in 1.3, of constructing the universe of the distributions,

$$\widetilde{\mathbb{C}}(\Omega) = \left(\widetilde{C}(\Omega), \widetilde{\partial}(\Omega) = \left\{\widetilde{\partial_\Omega^\alpha} : \ \alpha \in \mathbb{N}^n\right\}\right),$$

starting with

$$\mathbb{C}(\Omega) = \left(C(\Omega), \partial(\Omega) = \left\{\partial_\Omega^\alpha : \ \alpha \in \mathbb{N}^n\right\}\right),$$

can now, with the definitions 1.4 and 1.6 in mind, be stated in the following simple and precise form:

> **Extension Problem** — Given the abelian, surjective, and with identity $S$-group
>
> $$\mathbb{C}(\Omega) = \left(C(\Omega), \partial(\Omega) = \left\{\partial_\Omega^\alpha : \ \alpha \in \mathbb{N}^n\right\}\right),$$
>
> to construct an extension
>
> $$\widetilde{\mathbb{C}}(\Omega) = \left(\widetilde{C}(\Omega), \widetilde{\partial}(\Omega) = \left\{\widetilde{\partial_\Omega^\alpha} : \ \alpha \in \mathbb{N}^n\right\}\right)$$
>
> of $\mathbb{C}(\Omega)$.

Clearly, any solution to this problem, i.e., every extension of the $S$-group of continuous functions (described in Example 1.5) which can be constructed, will attend the conditions (a) to (d) formulated in 1.3; by the way, the Definitions 1.4 and 1.6 were elaborated for, among others, this goal.

# Strict Extension

## 1.9 Strict Growth

What other conditions, beyond those described in 1.3 (items (a) to (d)) as essential to any habitat suitable for distributions and that, therefore, relies bundled in the definition of $S$-groups extension, should we demand of possible solutions to the extension problem? Some direction in this sense can be obtained considering, again, the question of modeling physical systems with differential equations. For the sake of simplicity, we will limit ourselves to ordinary differential equations, whose (classical) solutions are, therefore, differentiable functions of a single real variable. Hence, our discussion will concern the group of continuous functions $C(\Omega) = C(J)$, where $\Omega = J \subseteq \mathbb{R}$ is an open interval of the real line, and the derivatives $\partial_\Omega^\alpha = \mathrm{d}_J^n$ being

$$\begin{aligned} \mathrm{d}_J^n : C^n(J) &\longrightarrow C(J) \\ f &\longmapsto \mathrm{d}_J^n(f) = \tfrac{\mathrm{d}^n f}{\mathrm{d}x^n} \end{aligned}$$

for each $n \in \mathbb{N}$.

As previously remarked, the differential equations that model physical systems, when taken in its extended versions on the universe of the distributions,

$$\widetilde{\mathbb{C}}(J) = \left(\widetilde{C}(J); \widetilde{\mathrm{d}}(J) = \left\{\widetilde{\mathrm{d}_J^n} : \ n \in \mathbb{N}\right\}\right),$$

should admit, as solutions, the classical ones, i.e., those attending the usual solution sense, not extended, as well as the ones representing physically relevant situations but, due to technicalities of the classical model (not extended), are not accepted as solutions. Hence, no classical solution should be discarded by extended models: “expansion without exclusion” is then a condition which should be demanded of possible solutions to the extension problem.

As we know, the concept of derivative of order $n$ is well-defined (classically) not only for functions $f \in C^n(J)$, whose derivatives of order $n$ are continuous functions on $J$, but also for certain functions $f \in C(J)$ which do not belong to $C^n(J)$. These functions have, therefore, a $n$-th derivative which is not continuous on $J$ and, as mathematical analysis teach us, their derivatives of order $n$ only have discontinuities of second kind[3]. We denote by $C'^n(J)$ the set of functions $f \in C(J)$ with classical order $n$ derivative not necessarily continuous on $J$, and by $C'(J)$ the functions of complex values with domain $J$ which are continuous or present second kind discontinuities. Clearly, $C^n(J)$ and $C(J)$ are, respectively, proper subsets of $C'^n(J)$ and $C'(J)$ and, if we denote by $\mathrm{d}'^n_J$ the function from $C'^n(J)$ into $C'(J)$ which associates to each $f \in C'^n(J)$ its $n$-th classical derivative, $\frac{\mathrm{d}^n f}{\mathrm{d}x^n}$, in $C'(J)$, i.e.,

$$\begin{aligned} \mathrm{d}'^n_J : C'^n(J) &\longrightarrow C'(J) \\ f &\longmapsto \mathrm{d}'^n_J(f) := \tfrac{\mathrm{d}^n f}{\mathrm{d}x^n}, \end{aligned}$$

we have that:

$$\mathrm{d}'^n_J(f) = \mathrm{d}^n_J(f) \quad \text{for every} \quad f \in C^n(J),$$

namely, $\mathrm{d}'^n_J$ is a prolongation of $\mathrm{d}^n_J$ to $C'^n(J)$.

[3] If $f \in C(J)$ is differentiable in $J$, then its derivative, $f'$, certainly does not present simple discontinuities in $J$; however, discontinuities of second kind may occur (Rudin, W. *Principles of Mathematical Analysis*. 3rd ed. McGraw-Hill Book Company, 1976).

It would then be reasonable to restrict the search for solutions to the extension problem to those,

$$\widetilde{\mathbb{C}}(J) = \left(\widetilde{C}(J), \widetilde{\mathrm{d}}(J) = \left\{\widetilde{\mathrm{d}^n_J} : \; n \in \mathbb{N}\right\}\right),$$

whose "extended derivatives", $\widetilde{\mathrm{d}^n_J} \in \widetilde{\mathrm{d}}(J)$, were prolongations not only of $d^n_J \in \mathrm{d}(J)$, but also of the derivatives $\mathrm{d}'^n_J$, i.e., such that:

$$\widetilde{\mathrm{d}^n_J}(f) = \mathrm{d}'^n_J(f) \quad \text{for every} \quad f \in C'^n_J.$$

The extension to such a universe $\widetilde{\mathbb{C}}(J)$ of the differential equations that describe physical situations, would result in new models which not only would keep accommodating the classical solutions given by functions on $C^n(J)$ (the domain of $\mathrm{d}^n_J$), but also those others (if existing), yet classical, represented by functions on $C'^n(J)$ (the domain of $\mathrm{d}'^n_J$) which does not belong to $C^n(J)$ and, possibly, other non-classical ones — expansion without exclusion.

Observe that the prolongation of $\mathrm{d}^n_J \colon C^n(J) \to C(J)$ to $C'^n(J)$, i.e., $\mathrm{d}'^n_J \colon C'^n(J) \to C'(J)$, requires the enlargement of the basic set $C(J)$, the image of $d^n_J$, to $C'(J)$: we are obliged to consider new objects, in the sense of not being in $C(J)$, namely, the functions with second kind discontinuities, since, for $n = 1, 2, \dots,$

$$\text{if} \quad f \in C'^n(J) \quad \text{is such that} \quad f \notin C^n(J), \quad \text{then} \quad \mathrm{d}'^n_J(f) \notin C(J).$$

Hence, in the search for expansion without exclusion, we must select, among the possible solutions to the extension problem, those satisfying the condition above. One possibility would be to require of the solutions $\widetilde{\mathbb{C}}(J)$ that, for $n = 1, 2, \dots,$

$$\text{if} \quad f \in C(J) \quad \text{is such that} \quad f \notin C^n(J), \quad \text{then} \quad \widetilde{\mathrm{d}^n_J}(f) \notin C(J),$$

or, more generally, for the case of multivariable functions defined on an open set $\Omega \subseteq \mathbb{R}^n$, require of the solutions $\widetilde{\mathbb{C}}(\Omega)$ to the extension problem that

$$\text{if} \quad f \in C(\Omega) \quad \text{and} \quad f \notin C^{|\alpha|}(\Omega), \quad \text{then} \quad \widetilde{\partial^\alpha_\Omega}(f) \notin C(\Omega),$$

for every multi-index $\alpha \in \mathbb{N}^n$, $\alpha \neq 0$.

This would not only imply the proper inclusion of $C(\Omega)$ in $\widetilde{C}(\Omega)$, but also that generalized derivatives of continuous non-differentiable functions in the usual sense to be new entities, i.e., elements of $\widetilde{C}(\Omega)$ not present in the set $C(\Omega)$ of the $S$-group $\mathbb{C}(\Omega) = (C(\Omega), \partial(\Omega))$. We can say that the referred condition required to the solutions $\widetilde{\mathbb{C}}(\Omega)$ of the extension problem, demands a proper, strict, growing of the universe $\mathbb{C}(\Omega)$. These considerations motivate the definition presented below.

## 1.10 Definition

Let $\mathbb{G} = (G, H)$ be a $S$-group and $\widehat{\mathbb{G}} = (\widehat{G}, \widehat{H})$ be an extension of $\mathbb{G}$.

**(a)** Given $\Phi \in H$, we say that its prolongation $\widehat{\Phi} \in \widehat{H}$ to $\widehat{G}$ is a **strict prolongation** if and only if

$$\widehat{\Phi}(g) \notin G \quad \text{if} \quad g \in G \quad \text{and} \quad g \notin G_\Phi \quad \text{(the domain of } \Phi)$$

or, equivalently,

$$g' = \widehat{\Phi}(g) \quad \text{with} \quad g, g' \in G \quad \text{if and only if} \quad g' = \Phi(g);$$

**(b)** We say that $\widehat{\mathbb{G}}$ is a **strict extension** of the $S$-group $\mathbb{G}$ if and only if for each $\Phi \in H$, its prolongation $\widehat{\Phi} \in \widehat{H}$ to $\widehat{G}$ is a strict prolongation.

Observe the use, in the definition above, of the abuse of language referred to in Remark 1.7. Formally, we should define, for example, based on the hypothesis of $\widehat{\mathbb{G}}$ being a $\delta$-extension of $\mathbb{G}$, the concept of $\delta$-strict prolongation of $\Phi$ to $\widehat{G}$ through the following requirement on $\widehat{\Phi}$:

$$\widehat{\Phi}\Big(\delta(g)\Big) \notin \delta(G) \quad \text{if} \quad g \in G \quad \text{and} \quad g \notin G_\Phi.$$

## 1.11 The Strict Extension Problem

In terms of Definition 1.10, the strict growing condition that we judged, by the reasons exposed in 1.9, reasonable to impose on possible solutions to the extension problem, equals the requirement on solutions $\widetilde{\mathbb{C}}(\Omega)$ of the referred problem to be strict extensions of the $S$-group $\mathbb{C}(\Omega)$. We will then divert our attention, so far directed to the extension problem, to concentrate it on the strict extension problem formulated ahead.

**Strict Extension Problem** — Given the abelian, surjective, and with identity $S$-group

$$\mathbb{C}(\Omega) = \left(C(\Omega), \partial(\Omega) = \left\{\partial^\alpha_\Omega : \ \alpha \in \mathbb{N}^n\right\}\right),$$

to construct a strict extension

$$\widetilde{\mathbb{C}}(\Omega) = \left(\widetilde{C}(\Omega), \widetilde{\partial}(\Omega) = \left\{\widetilde{\partial_\Omega^\alpha}\colon\ \alpha \in \mathbb{N}^n\right\}\right)$$

of $\mathbb{C}(\Omega)$.

Future considerations will require some properties associated with strict extensions of $S$-groups, from which we highlight, for convenience, the one presented in the proposition below.

## 1.12 Proposition

*Let $G$ be a group, $\Phi$ be a surjective homomorphism from a subgroup $G_\Phi$ of $G$ onto $G$, $\widehat{G}$ be a group which admits $G$ as subgroup and $\widehat{\Phi}$ be a prolongation of $\Phi$ to $\widehat{G}$. Denoting by $N(\Phi)$ and $N(\widehat{\Phi})$ the kernels*[4] *of $\Phi$ and $\widehat{\Phi}$, respectively, then $\widehat{\Phi}$ is a strict prolongation of $\Phi$ to $\widehat{G}$ if and only if*

$$G \cap N(\widehat{\Phi}) = N(\Phi).$$

*Proof.* First, suppose that $\widehat{\Phi}$ is a strict prolongation of $\Phi$ to $\widehat{G}$. In this case, from Definition 1.10(a) we have that, for $g, g' \in G$,

$$g' = \widehat{\Phi}(g) \quad \text{if and only if} \quad g' = \Phi(g).$$

Hence, for $g' = 0 \in G$, we get that, for $g \in G$

$$\widehat{\Phi}(g) = 0 \quad \text{if and only if} \quad \Phi(g) = 0,$$

that is,

$$G \cap N(\widehat{\Phi}) = N(\Phi).$$

Lets now prove the converse, i.e., we will demonstrate that the prolongation $\widehat{\Phi}$ of $\Phi$ to $\widehat{G}$ is a strict prolongation if the last equation above is true. In order to do so, let $g, g' \in G$. Clearly, if $g' = \Phi(g)$ then $g' = \widehat{\Phi}(g)$ since $\widehat{\Phi}$ is a prolongation of $\Phi$. Suppose now that $g' = \widehat{\Phi}(g)$. Since $\Phi\colon G_\Phi \longrightarrow G$ is surjective, then there exist $g'' \in G_\Phi$ such that $g' = \Phi(g'') = \widehat{\Phi}(g'')$. Hence, $\widehat{\Phi}(g'') = \widehat{\Phi}(g)$ and, therefore, $g - g'' \in N(\widehat{\Phi})$ and $g - g'' \in G$, i.e.,

$$g - g'' \in G \cap N(\widehat{\Phi}) = N(\Phi) \subseteq G_\Phi.$$

[4] Let $G_1$ and $G_2$ be groups, $\varphi\colon G_1 \longrightarrow G_2$ a homomorphism and $0 \in G_2$ the neutral element of $G_2$. The kernel of $\varphi$, denoted by $N(\varphi)$, is defined by:

$$N(\varphi) := \{g \in G_1\colon\ \varphi(g) = 0\}.$$

Hence, we have $g - g'' \in G_\Phi$ and $g'' \in G_\Phi$ which allow us to conclude that $(g - g'') + g'' = g \in G_\Phi$ since $G_\Phi$ is a group. Thus, $g' = \widehat{\Phi}(g) = \Phi(g)$. In short, we proved that if $G \cap N(\widehat{\Phi}) = N(\Phi)$ then, for $g, g' \in G$, $g' = \widehat{\Phi}(g)$ if and only if $g' = \Phi(g)$, which means that the prolongation $\widehat{\Phi}$ is strict.[5] ■

# Closed Extension

## 1.13 Moderate Growth

As we have seen, the considerations made in 1.9 suggest to consider, not the extension problem stated in 1.8, of constructing an extension of the $S$-group $\mathbb{C}(\Omega) = (C(\Omega), \partial(\Omega))$, but the Strict Extension Problem, formulated in 1.11, of constructing a strict extension of $\mathbb{C}(\Omega)$, $\widetilde{\mathbb{C}}(\Omega) = (\widetilde{C}(\Omega), \widetilde{\partial}(\Omega))$, from which we would define the distributions and their derivatives as being, respectively, the elements of the sets $\widetilde{C}(\Omega)$ and $\widetilde{\partial}(\Omega)$.

On the other hand, it is reasonable to take $\widetilde{\mathbb{C}}(\Omega)$ as the most economic strict extension of $\mathbb{C}(\Omega)$ that meets our goal: to free differential calculus from difficulties originated from the existence of continuous and non-differentiable functions. We would like to extend $\mathbb{C}(\Omega)$ not more than the necessary and sufficient for this purpose.

What conditions, beyond that of being strict, could be reasonably required of an extension $\widetilde{\mathbb{C}}(\Omega)$ of $\mathbb{C}(\Omega)$ aiming to attend this "tight measure growth" question?

In $\widetilde{\mathbb{C}}(\Omega) = (\widetilde{C}(\Omega), \widetilde{\partial}(\Omega) = \{\widetilde{\partial^\alpha_\Omega} \colon \ \alpha \in \mathbb{N}^n\})$ will reside the derivatives of the continuous functions, i.e., there exist, as distributions, all derivatives of any continuous function or, more precisely, if $f \in C(\Omega)$ and $\alpha \in \mathbb{N}^n$, $\widetilde{\partial^\alpha_\Omega}(f) \in \widetilde{C}(\Omega)$. On the other hand, in view of the basic goal of every continuous function to be differentiable, it is reasonable to require the converse, i.e., that every distribution to be the derivative of a continuous function or, more precisely: if $\widetilde{f} \in \widetilde{C(\Omega)}$, then exist a multi-index $\alpha$ and a function $f \in C(\Omega)$ such that $\widetilde{f} = \widetilde{\partial^\alpha_\Omega}(f)$. This requirement looks reasonable insofar as, being possible to obtain the class of distributions, $\widetilde{C}(\Omega)$, from $C(\Omega)$ adding up to the latter only the generalized derivatives of the continuous and non-differentiable functions in the usual sense, $C(\Omega)$ would be extended the necessary and the sufficient, i.e., tightly measured. In fact, since the generalized derivatives, $\widetilde{\partial^\alpha_\Omega}$, are endomorphisms on $\widetilde{C}(\Omega)$, then every distribution is infinitely differentiable and, as $C(\Omega) \subseteq \widetilde{C}(\Omega)$, so are the continuous functions. Hence, on the hypothesis of $\widetilde{\mathbb{C}}(\Omega)$ fulfil the condition in consideration, i.e., if each distribution $\widetilde{f} \in \widetilde{C}(\Omega)$ is the generalized derivative of some continuous functions, no proper subclass of $\widetilde{C}(\Omega)$ would be adequate; in this sense, such extension is so economic as it could possibly be.

[5] The end of each proof will be indicated by the symbol ■, a notation due to P. Halmos.

In what follows, as a definition, we introduce in a precise way the concept of "tight measure growth".

## 1.14 Definition

Let $\widehat{\mathbb{G}} = (\widehat{G}, \widehat{H})$ be an extension of the $S$-group $\mathbb{G} = (G, H)$. We say that $\widehat{\mathbb{G}}$ is a **closed extension** of $\mathbb{G}$ if and only if for each $\widehat{g} \in \widehat{G}$, there exist $g \in G$ and $\Phi \in H$ such that $\widehat{g} = \widehat{\Phi}(g)$, where $\widehat{\Phi} \in \widehat{H}$ is the prolongation of $\Phi$ to $\widehat{G}$.

## 1.15 The Strict and Closed Extension Problem

Its clear, with Definition 1.14 in mind, that our economic extension condition discussed in 1.13, is equivalent to require of $\mathbb{C}(\Omega)$ an extension $\widetilde{\mathbb{C}}(\Omega)$ that, beyond strict as previously examined, to be also closed. We change therefore, again, our attention, so far focused on the strict extension problem, to concentrate it now on the strict and closed extension problem proposed below.

**Strict and Closed Extension Problem** — Given the abelian, surjective, and with identity $S$-group

$$\mathbb{C}(\Omega) = \left( C(\Omega), \partial(\Omega) = \left\{ \partial_\Omega^\alpha : \ \alpha \in \mathbb{N}^n \right\} \right),$$

to construct a strict and closed extension

$$\widetilde{\mathbb{C}}(\Omega) = \left( \widetilde{C}(\Omega), \widetilde{\partial}(\Omega) = \left\{ \widetilde{\partial_\Omega^\alpha} : \ \alpha \in \mathbb{N}^n \right\} \right)$$

of $\mathbb{C}(\Omega)$.

Would the problem above be well-posed, in the sense of the conditions required of $\widetilde{\mathbb{C}}(\Omega)$, enclosed on the definitions of $S$-group extension, strict extension and closed extension, be consistent? Being consistent would exist, for $\mathbb{C}(\Omega)$, strict and closed extensions not isomorphic? These are questions we are going to ponder. However, instead of treating them relative to the particular $S$-group of continuous functions, $\mathbb{C}(\Omega)$, why not in the context of an abstract version of the problem where, instead of $\mathbb{C}(\Omega)$, we take an arbitrary $S$-group, $\mathbb{G}$, satisfying conditions obtained by abstractions of the main properties of $\mathbb{C}(\Omega)$. After all, as *J. Sebastião e Silva*[6] would say,

[6] Silva, J. S. *Sur Une Construction Axiomatique de la Théorie des Distributions.* Rev. da Fac. de Ciências de Lisboa, 2ª série-A, vol 4, pp. 79–186, (1954–1955).

*"... qu'on ne sait jamais si d'autres applications ne seront pas possibles."*

Hence, through the path of abstraction, we are led to one of the central questions of this monograph, the abstract version of the strict and closed extension problem, presented below:

> **Strict and Closed Extension Problem (abstract version)** — Given an abelian, surjective, and with identity $S$-group, $\mathbb{G} = (G, H)$, to construct a strict and closed extension $\widetilde{\mathbb{G}} = (\widetilde{G}, \widetilde{H})$ of $\mathbb{G}$.

It seems opportune to present some properties associated with strict and closed extensions, as done through the proposition below.

## 1.16 Proposition

*Let $\mathbb{G} = (G, H)$ be a $S$-group and $\widehat{\mathbb{G}} = (\widehat{G}, \widehat{H})$ be an extension of $\mathbb{G}$:*

**(a)** *if $\mathbb{E} = (E, F)$ is a closed extension of $\widehat{\mathbb{G}}$, then $\mathbb{E} = \widehat{\mathbb{G}}$;*

**(b)** *if $\widehat{\mathbb{G}}$ is a strict and closed extension of $\mathbb{G}$ and $\mathbb{G}$ is a surjective $S$-group, then, for every $\Phi \in H$,*

$$N(\widehat{\Phi}) = N(\Phi)$$

*being $\widehat{\Phi} \in \widehat{H}$ the prolongation of $\Phi$ to $\widehat{G}$.*

*Proof.* (a) Since $\mathbb{E} = (E, F)$ is, by hypothesis, an extension of $\widehat{\mathbb{G}} = (\widehat{G}, \widehat{H})$, then $\widehat{G}$ is a subgroup of $E$, therefore $\widehat{G} \subseteq E$, and there exists an isomorphism from the semigroup $\widehat{H}$ onto the semigroup $F$, let us say

$$\begin{aligned} p : \widehat{H} &\longrightarrow F \\ \widehat{\Phi} &\longmapsto p(\widehat{\Phi}) = \Phi' \end{aligned}$$

where $\Phi'$ is a prolongation of the endomorphism $\widehat{\Phi} \colon \widehat{G} \longrightarrow \widehat{G}$ to $E$. Hence, $\Phi' \colon E \longrightarrow E$ is an endomorphism on $E$ such that

$$\Phi'(\widehat{g}) = \widehat{\Phi}(\widehat{g}) \quad \text{for every} \quad \widehat{g} \in \widehat{G}.$$

Now, let $e \in E$ be arbitrarily fixed. By hypothesis, $\mathbb{E}$ is a closed extension of $\widehat{\mathbb{G}}$, then there exist $\widehat{\Phi} \in \widehat{H}$ and $\widehat{g} \in \widehat{G}$ such that

$$e = \left(p(\widehat{\Phi})\right)(\widehat{g}) = \Phi'(\widehat{g}).$$

But, as $\widehat{g} \in \widehat{G}$, then $\Phi'(\widehat{g}) = \widehat{\Phi}(\widehat{g})$ and, hence,

$$e = \widehat{\Phi}(\widehat{g}) \in \widehat{G}.$$

Since $e \in E$ is arbitrary, we conclude that $E \subseteq \widehat{G}$ and, hence, as $\widehat{G} \subseteq E$, we get $E = \widehat{G}$ and, consequently, $p(\widehat{\Phi}) = \widehat{\Phi}$ for every $\widehat{\Phi} \in \widehat{H}$, which leads us, being $p\colon \widehat{H} \longrightarrow F$ an isomorphism, to $\widehat{H} = F$. With $E = \widehat{G}$ and $F = \widehat{H}$ we conclude that $\mathbb{E} = \widehat{\mathbb{G}}$.

(b) Since $\widehat{\mathbb{G}} = (\widehat{G}, \widehat{H})$ is, by hypothesis, a strict and closed extension of the $S$-group $\mathbb{G} = (G, H)$ which, in turn, also by hypothesis, is surjective, we have, by Proposition 1.12, that

$$G \cap N(\widehat{\Phi}) = N(\Phi)$$

for every $\Phi \in H$. Hence, to prove (b) it is enough to verify that

$$N(\widehat{\Phi}) \subseteq G,$$

that is:

$$\text{if} \quad \widehat{g} \in \widehat{G} \quad \text{and} \quad \widehat{g} \notin G, \quad \text{then} \quad \widehat{g} \notin N(\widehat{\Phi}).$$

Suppose then that $\widehat{g} \in \widehat{G}$ is such that

$$\widehat{g} \notin G.$$

By *reductio ad absurdum*, let us admit that

$$\widehat{g} \in N(\widehat{\Phi}).$$

Now, $\widehat{\mathbb{G}}$ is a closed extension of $\mathbb{G}$ and hence there exist $\Psi \in H$ and $g \in G$ such that

$$\widehat{g} = \widehat{\Psi}(g).$$

Since $\widehat{g} = \widehat{\Psi}(g) \in N(\widehat{\Phi})$, then

$$\widehat{\Phi}\left(\widehat{\Psi}(g)\right) = \left(\widehat{\Phi\Psi}\right)(g) = 0,$$

that is,

$$g \in G \cap N\left(\widehat{\Phi\Psi}\right).$$

But $\widehat{\Phi\Psi}$ is a strict prolongation of $\Phi\Psi$, since $\widehat{\mathbb{G}}$ is a strict extension of $\mathbb{G}$, and beyond that, $\Phi\Psi \in H$ is a surjective homomorphism, since $\mathbb{G}$ is a surjective $S$-group. Thus, by Proposition 1.12,

$$G \cap N\left(\widehat{\Phi\Psi}\right) = N(\Phi\Psi)$$

and hence we get that

$$g \in N(\Phi\Psi) \subseteq G_{\Phi\Psi} = \left\{g \in G\colon\ g \in G_\Psi \quad \text{and} \quad \Psi(g) \in G_\Phi\right\}.$$

From where $g \in G_\Psi$ and, therefore,

$$\widehat{g} = \widehat{\Psi}(g) = \Psi(g) \in G$$

what contradicts the hypothesis of $\widehat{g} \notin G$. ■

# Consistency of the Strict and Closed Extension Problem

## 1.17 Remarks

To justify the abstract formulation given, in 1.15, to the strict and closed extension problem, it is necessary that exists at least one abelian, surjective and with identity $S$-group not isomorphic to that of the continuous functions, $\mathbb{C}(\Omega)$. We present ahead, in 1.18, one such $S$-group[7]. Furthermore, for the referred example we prove the existence of a strict and closed extension. Hence, not only the abstract approach proposed in 1.15 is justified but, above all, the consistency of the conditions gathered on the definitions of $S$-group extension, strict extension, and closed extension are ensured.

## 1.18 The $S$-group $\mathbb{Z}$

It is known that the set of integer numbers, $\{0, \pm 1, \pm 2, \ldots\}$, with the usual addition operation is an abelian group and that, this group, with the usual integer multiplication is a ring[8]. Here, the letter $Z$ will be used to denote both the group and the ring referred.

For each integer $n > 0$ we define:

$$nZ := \Big\{ nm : \ m \in Z \Big\}$$

and

$$\begin{aligned} f_n : nZ &\longrightarrow Z \\ m &\longmapsto f_n(m) := \tfrac{m}{n}. \end{aligned}$$

It results that, for each integer $n > 0$, $nZ$ is a subgroup of $Z$ and $f_n$ is a surjective homomorphism from the subgroup $nZ$ over $Z$. Beyond that, $f_1 = I_Z$, i.e., $f_1(m) = m$ for every $m \in 1Z = Z$.

Observe now that the set $\{f_n : \ n = 1, 2, \ldots\}$ equipped with the multiplication defined by

$$f_{n_1} f_{n_2} := f_{n_1 n_2} \quad \text{with} \quad n_1, n_2 \in \Big\{1, 2, \ldots\Big\},$$

is an abelian semigroup. Hence, $\{f_n : \ n = 1, 2, \ldots\}$ with the multiplication above is a semigroup of surjective homomorphisms of subgroups of $Z$ (the $nZ$'s) over $Z$, which is abelian and contains the identity $I_Z = f_1$. Therefore, denoting this semigroup as $H_Z$, we get that

$$\mathbb{Z} := (Z, H_Z)$$

---

[7] The content of this item was gently suggested by Leandro Gustavo Gomes.

[8] To the reader not familiarized with the algebraical content of this section we suggest: Birkhoff, G. and MacLane, S. *A Survey of Modern Algebra*. 4th ed. New York: Macmillan, 1977.

is an abelian, surjective and with identity $S$-group, not isomorphic to $\mathbb{C}(\Omega) = (C(\Omega), \partial(\Omega))$ since $Z$ and $C(\Omega)$ have distinct cardinalities.

Let now be denoted by $Q$ the ring with unity given by the set of rational numbers with the usual operations of addition and multiplication. Let also $\widetilde{Q}$ be a ring that has $Q$ as a subring. Hence, for instance, $\widetilde{Q}$ can denote the ring of real numbers, $\mathbb{R}$, the ring of complex numbers, $C = \mathbb{R}^2$, or $Q$ itself. In this context, we define $\widetilde{f}_n$, for each integer $n > 0$, as:

$$\begin{aligned} \widetilde{f}_n : \widetilde{Q} &\longrightarrow \widetilde{Q} \\ \widetilde{q} &\longmapsto \widetilde{f}_n(\widetilde{q}) := \tfrac{\widetilde{q}}{n}. \end{aligned}$$

Since every ring is also a group with respect to its addition, then, $\widetilde{Q}$ is a group (with respect to its addition) which has as subgroups, among others, $Q$ with the rational addition, the group $Z$ and its subgroups $nZ$. Furthermore, it is easily verifiable that, for each integer $n > 0$, $\widetilde{f}_n$ is an endomorphism on the group $\widetilde{Q}$ such that

$$\widetilde{f}_n(m) = f_n(m) \quad \text{for every} \quad m \in nZ \subseteq \widetilde{Q}.$$

Defining now, for the members of the set $\{\widetilde{f}_n : n = 1, 2, \ldots\}$, the multiplication

$$\widetilde{f}_{n_1}\widetilde{f}_{n_2} := \widetilde{f}_{n_1 n_2} \quad \text{with} \quad n_1, n_2 \in \Big\{1, 2, \ldots\Big\},$$

we can see that the mentioned set equipped with this multiplication is an abelian semigroup of endomorphisms on $\widetilde{Q}$, henceforth denoted $\widetilde{H}_Z$. Beyond that, the function $\widetilde{p}$ defined by

$$\begin{aligned} \widetilde{p} : H_Z &\longrightarrow \widetilde{H}_Z \\ f_n &\longmapsto \widetilde{p}(f_n) := \widetilde{f}_n, \end{aligned}$$

is a bijection such that

$$\widetilde{p}(f_{n_1} f_{n_2}) = \widetilde{p}(f_{n_1})\widetilde{p}(f_{n_2}),$$

that is,

$$\widetilde{f_{n_1} f_{n_2}} = \widetilde{f}_{n_1}\widetilde{f}_{n_2}.$$

In other words, the semigroups $H_Z$ and $\widetilde{H}_Z$ are isomorphic, with $\widetilde{p}$ an isomorphism from $H_Z$ onto $\widetilde{H}_Z$: $\widetilde{H}_Z$ is composed exact and precisely by the prolongations $\widetilde{p}(f_n) = \widetilde{f}_n$ of the elements $f_n \in H_Z$ to $\widetilde{Q}$. Now, with the Definition 1.6 ((c) and (d)) in mind, the conclusions above show us that

$$\widetilde{H}_Z = \Big\{\widetilde{f}_n : \ n = 1, 2, \ldots\Big\}$$

is a prolongation of the semigroup $H_Z = \{f_n : \ n = 1, 2, \ldots\}$ to $\widetilde{Q}$ and that

$$\widetilde{\mathbb{Q}} := \Big(\widetilde{Q}, \widetilde{H}_Z\Big)$$

is an extension of the $S$-group $\mathbb{Z} = (Z, H_Z)$, with $\widetilde{Q}$ being any ring with unity that has $Q$ as one of its subrings.

Finally, the question about the consistency of the strict and closed extension problem is answered positively with the proof of the proposition below.

## 1.19 Proposition

*The extension* $\widetilde{\mathbb{Q}} = (\widetilde{Q}, \widetilde{H}_Z)$ *of the abelian, surjective and with identity S-group* $\mathbb{Z} = (Z, H_Z)$, *is a strict extension. Furthermore,* $\widetilde{\mathbb{Q}}$ *is a closed extension of* $\mathbb{Z}$ *if and only if* $\widetilde{Q} = Q$.

*Proof.* Let us first show that the extension is strict. In order to do so, take $n \in \{1, 2, \ldots\}$ and $m \in Z$ such that $m \notin nZ$ ($= Z_{f_n}$, the domain of $f_n$, according to the notation introduced in Definition 1.4(a)) In this case, $m = nk + l$ with $k \in Z$ and $l \in Z$ such that $1 \leqslant |l| < n$. Hence, $\widetilde{f}_n(m) = \frac{m}{n} = k + \frac{l}{n}$ with $k \in Z$ and $0 < \left|\frac{l}{n}\right| < 1$. Thus, $\widetilde{f}_n(m) \notin Z$ and, therefore, the extension $\widetilde{\mathbb{Q}}$ of $\mathbb{Z}$ is strict.

Observe now that, if $q \in Q \subseteq \widetilde{Q}$, then there exist $m_0, n_0 \in Z$, with $n_0 > 0$, such that $q = \frac{m_0}{n_0}$. Hence, $q = \frac{m_0}{n_0} = f_{n_0}(m_0) = \widetilde{f}_{n_0}(m_0)$, what allows us to conclude that $q \in \widetilde{f}_{n_0}(Z)$ and, therefore, $q \in \bigcup_{n=1}^{\infty} \widetilde{f}_n(Z)$. Thus, $Q \subseteq \bigcup_{n=1}^{\infty} \widetilde{f}_n(Z)$.

Conversely, if $q \in \bigcup_{n=1}^{\infty} \widetilde{f}_n(Z)$, then $q \in \widetilde{f}_{n_0}(Z)$ for some integer $n_0 > 0$ and, therefore, $q = \widetilde{f}_{n_0}(m_0)$ for some $m_0 \in Z$, i.e., $q = \frac{m_0}{n_0} \in Q$. Hence, we see that $\bigcup_{n=1}^{\infty} \widetilde{f}_n(Z) \subseteq Q$. Consequently,

$$\bigcup_{n=1}^{\infty} \widetilde{f}_n(Z) = Q.$$

The equality above completes the proof, since it informs us that: given $\widetilde{q} \in \widetilde{Q}$, then there exist $\widetilde{f}_n \in \widetilde{H}_Z$ and $m \in Z$ such that $\widetilde{f}_n(m) = \widetilde{q}$ if and only if $\widetilde{Q} = Q$. ∎

# 2

# THEOREM OF EXTENSION OF $S$-GROUPS AND THE $\widetilde{\mathbb{G}}$-DISTRIBUTIONS

## The Equivalence Relation $\sim$

### 2.1 Preliminaries

Of the questions raised in 1.15, right after the formulation of the strict and closed extension problem, namely, if the referred problem would be well-posed, that is, if consistent regarding the required conditions, and the one relative to the existence, for an abelian, surjective, and with identity $S$-group (such as, for instance, the $S$-group $\mathbb{C}(\Omega)$ of continuous functions), of non-isomorphic strict and closed extensions, only the first one was completely answered through Proposition 1.19. This chapter intends to answer the second question. Hence, in a more explicit and precise manner, one of our goals in this chapter consists of answering the following questions: Being $\mathbb{G}$ an abelian, surjective, and with identity $S$-group, is there a strict and closed extension of $\mathbb{G}$? Existing one, is it unique unless isomorphism?

As we will see, the hypothesis on the existence of such an extension of $\mathbb{G}$, that is, that the strict and closed extension problem (in its abstract version) has a solution, will inform us, in a very natural manner, how, step-by-step, we can construct it. Each step will culminate in a necessary proposition that, if true, allows the next step. Hence, proposition by proposition, at the end of the path we will have provided the proof of a theorem, the Theorem of Extension of $S$-Groups, that ensures the existence of a unique solution (unless isomorphism) of the referred problem, in addition to displaying it.

## 2.2 The Relation $\sim$

The problem stated in 1.15, about constructing a strict and closed extension, $\widetilde{\mathbb{G}} = (\widetilde{G}, \widetilde{H})$, of a given abelian, surjective, and with identity $S$-group, $\mathbb{G} = (G, H)$, reflects the purpose of giving citizenship to the symbol $\Phi(g)$, where $\Phi \in H$ and $g \in G$, when $g$ does not belong to the domain $G_\Phi \subseteq G$ of $\Phi$; lets remember that the $S$-group $\mathbb{G} = (G, H)$ is the abstract version of $\mathbb{C}(\Omega) = (C(\Omega), \partial(\Omega) = \{\partial_\Omega^\alpha\colon \alpha \in \mathbb{N}^n\})$, and that our basic goal is to eliminate difficulties caused by the existence of non-differentiable continuous functions, that is, of elements $f \in C(\Omega)$ (elements $g \in G$ in the abstract version) for which the expression $\partial_\Omega^\alpha(f)$ ($\Phi(g)$ in the abstract version) does not make sense.

We want an "environment" at which $\Phi(g)$ makes sense, exists, even if $g \notin G_\Phi$. However, this environment, $\widetilde{\mathbb{G}}$, must be built from $\mathbb{G}$, that is, the objects of $\widetilde{\mathbb{G}}$ must be built only and solely by material made available in $\mathbb{G}$. Well, since $\mathbb{G}$ is a structure endowed with two sets, $G$ and $H$, with two types of "raw material", sounds reasonable to take the cartesian product $H \times G$ as a first mix of the available ingredients to produce something new: the ordered pair $(\Phi, g) \in H \times G$.

This mix, $H \times G$, shelter the symbol $(\Phi, g)$ with $\Phi \in H$ and $g \in G$, even if $g \notin G_\Phi$; hence, when we fell bothered with the presence of $\Phi \in H$ (respectively, of $\partial_\Omega^\alpha \in \partial(\Omega)$) juxtaposed to $g \in G$ (respectively, to $f \in C(\Omega)$), that is, with $\Phi(g)$ (respectively, with $\partial_\Omega^\alpha(f)$), by the fact that $g \notin G_\Phi$ (respectively, that $f \notin C^{|\alpha|}(\Omega)$), we can undo the hassle by distancing $\Phi$ from $g$, separating them with a comma and, to give unity to this triplet, enclosing them with parenthesis, obtaining in this way, from the "homeless" $\Phi(g)$, a respectable citizen of $H \times G$, namely, the ordered pair $(\Phi, g)$. Could we then construct the objects of $\widetilde{G}$ with elements of $H \times G$, for example, taking for $\widetilde{G}$ the "mixture" $H \times G$?

Let us suppose that there exists a solution to our problem, i.e., lets say that, $\widehat{\mathbb{G}} = (\widehat{G}, \widehat{H})$ is a strict and closed extension of $\mathbb{G} = (G, H)$. In this case, the following function on $H \times G$ into $\widehat{G}$ will be well-defined:

$$p: H \times G \longrightarrow \widehat{G}$$
$$(\Phi, g) \longmapsto p(\Phi, g) := \widehat{\Phi}(g)$$

where $\widehat{\Phi} \in \widehat{H}$ is the prolongation of $\Phi \in H$ to $\widehat{G}$. Note that the Definition 1.6(c), of semigroup prolongation, requires that in $\widehat{H}$ there exists a unique prolongation, $\widehat{\Phi}$, of $\Phi$ to $\widehat{G}$, which supports the above definition for the function $p$. Observe now that $p$ is surjective, since if $\widehat{g} \in \widehat{G}$, then, as $\widehat{\mathbb{G}}$ is a closed extension of $\mathbb{G}$, there exits $(\Phi, g) \in H \times G$ such that $\widehat{g} = \widehat{\Phi}(g) = p(\Phi, g)$.

Now, for each $\widehat{g} = \widehat{\Phi}(g) \in \widehat{G}$, let $[\Phi, g] \subseteq H \times G$ be defined by:

$$[\Phi, g] \coloneqq p^{-1}\Big(\widehat{g} = \widehat{\Phi}(g)\Big) = \\ = \Big\{(\varphi, g') \in H \times G \colon p(\varphi, g') = \widehat{g}\Big\} = \\ = \Big\{(\varphi, g') \in H \times G \colon \widehat{\varphi}(g') = \widehat{\Phi}(g)\Big\}.$$

Clearly, for $\widehat{g} = \widehat{\Phi}(g)$ and $\widehat{h} = \widehat{\Psi}(h)$ in $\widehat{G}$, we have that:

**(a)** $[\Phi, g] = [\Psi, h]$ if and only if $\widehat{g} = \widehat{h}$;

**(b)** $[\Phi, g] \cap [\Psi, h] = \varnothing$ if and only if $\widehat{g} \neq \widehat{h}$;

**(c)** $[\Phi, g] \neq \varnothing$;

**(d)** $\bigcup_{(\Phi,g)\in H\times G} [\Phi, g] = H \times G$.

From (b), (c) and (d) result that the set $(H \times G)_p$ defined by

$$(H \times G)_p \coloneqq \Big\{[\Phi, g] \colon (\Phi, g) \in H \times G\Big\},$$

is a partition of $H \times G$, while (a) legitimates the definition of the following bijective function, $\alpha$, on $\widehat{G}$ onto $(H \times G)_p$:

$$\begin{aligned} \alpha : \quad \widehat{G} &\longrightarrow (H \times G)_p \\ \widehat{g} = \widehat{\Phi}(g) &\longmapsto \alpha(\widehat{g}) \coloneqq [\Phi, g]. \end{aligned}$$

Therefore, the hypothesis on the existence of $\widehat{\mathbb{G}}$ implies, through the bijective function $\alpha$, that we can identify each object $\widehat{g} \in \widehat{G}$ with a collection of ordered pairs in $H \times G$, $\alpha(\widehat{g})$, which, somehow, supports the previously considered idea of constructing the "distributions", that is, the members of the extension searched for $\mathbb{G}$, from the elements of $H \times G$.

As it is known, every partition $\mathcal{P}$ of a set $A \neq \varnothing$ allows to define an equivalence relation in $A$, denominated the relation induced by $\mathcal{P}$, whose equivalence classes are, exactly, the elements of the partition. In particular, for the partition of $H \times G$ determined by the function $p$, namely, $(H \times G)_p$, the induced equivalence relation (in $H \times G$), henceforth denoted $\sim$, is given by: whatever $(\Phi, g)$ and $(\Psi, h)$ in $H \times G$,

$$(\Phi, g) \sim (\Psi, h) \quad \text{if and only if} \quad [\Phi, g] = [\Psi, h].$$

Since $[\Phi, g] = [\Psi, h]$ is equivalent to $\alpha^{-1}([\Phi, g]) = \alpha^{-1}([\Psi, h])$, that is, $\widehat{\Phi}(g) = \widehat{\Psi}(h)$, we have that:

$$(\Phi, g) \sim (\Psi, h) \quad \text{if and only if} \quad \widehat{\Phi}(g) = \widehat{\Psi}(h).$$

Taking now into account that the $S$-group $\mathbb{G} = (G, H)$ is surjective, then the sets

$$\Psi^{-1}(g) := \left\{ g' \in G_\Psi : \Psi(g') = g \right\}$$

and

$$\Phi^{-1}(h) := \left\{ g' \in G_\Phi : \Phi(g') = h \right\}$$

are not empty. Hence, for every $g^* \in \Psi^{-1}(g)$ and $h^* \in \Phi^{-1}(h)$ we have:

$$\widehat{\Phi}(g) = \widehat{\Phi}\Big(\Psi(g^*)\Big) = \widehat{\Phi}\Big(\widehat{\Psi}(g^*)\Big) = \Big(\widehat{\Phi\Psi}\Big)(g^*)$$

and

$$\widehat{\Psi}(h) = \widehat{\Psi}\Big(\Phi(h^*)\Big) = \widehat{\Psi}\Big(\widehat{\Phi}(h^*)\Big) = \widehat{\Phi}\Big(\widehat{\Psi}(h^*)\Big) = \Big(\widehat{\Phi\Psi}\Big)(h^*).$$

Consequently, $\widehat{\Phi}(g) = \widehat{\Psi}(h)$ if and only if

$$g^* - h^* \in N\left(\widehat{\Phi\Psi}\right) \quad \text{for every} \quad g^* \in \Psi^{-1}(g) \quad \text{and every} \quad h^* \in \Phi^{-1}(h)$$

or, equivalently[9],

$$\widehat{\Phi}(g) = \widehat{\Psi}(h) \quad \text{if and only if} \quad \Psi^{-1}(g) - \Phi^{-1}(h) \subseteq N\left(\widehat{\Phi\Psi}\right).$$

Thus, the equivalence relation $\sim$ in $H \times G$, induced by the partition $(H \times G)_p$, can be expressed in the following form:

$$(\Phi, g) \sim (\Psi, h) \quad \text{if and only if} \quad \Psi^{-1}(g) - \Phi^{-1}(h) \subseteq N\left(\widehat{\Phi\Psi}\right).$$

In short, existing a solution $\widehat{\mathbb{G}}$ for the strict and closed extension problem, then, necessarily, $\sim$, as defined above, is an equivalence relation in $H \times G$ such that the partition of $H \times G$ determined by it, that is, the set of its equivalence classes, henceforth denoted $[H \times G]_\sim$, is equal to the partition $(H \times G)_p$ determined by the function $p$, that is,

$$[H \times G]_\sim = (H \times G)_p$$

and, therefore, through the bijective function

$$\begin{aligned} \alpha : \qquad \widehat{G} &\longrightarrow (H \times G)_p = [H \times G]_\sim \\ \widehat{g} = \widehat{\Phi}(g) &\longmapsto \alpha(\widehat{g}) = [\Phi, g], \end{aligned}$$

we identify $\widehat{G}$ with $[H \times G]_\sim$.

Observe now that the memory of the hypothesis of the existence of $\widehat{\mathbb{G}}$ in the definition of $\sim$, is fully concentrated in the term $N(\widehat{\Phi\Psi})$. But, from Proposition 1.16(b) we

9 Being $A$ and $B$ subsets of the same group $C$, we define:

$$A \pm B := \Big\{ c \in C : c = a \pm b \quad \text{with} \quad a \in A \quad \text{and} \quad b \in B \Big\}.$$

have that $N(\widehat{\Phi\Psi}) = N(\Phi\Psi)$ and thus we can rewrite the definition of $\sim$ in the following form which does not mention $\widehat{\mathbb{G}}$, that is, independent of the existence of $\widehat{\mathbb{G}}$:

$$(\Phi, g) \sim (\Psi, h) \quad \text{if and only if} \quad \Psi^{-1}(g) - \Phi^{-1}(h) \subseteq N\left(\Phi\Psi\right).$$

Therefore, the hypothesis of the existence of $\widehat{\mathbb{G}}$ suggests us the definition of a binary relation in $H \times G$, $\sim$, through, exclusively, the elements of $\mathbb{G} = (G, H)$, without any reference to $\widehat{\mathbb{G}}$, which, being the hypothesis true, shall be an equivalence relation. Hence, if $\sim$ is not an equivalence relation, our problem does not have any solution. On the other hand, being it an equivalence relation, even not being this sufficient to ensure the existence of a solution $\widehat{\mathbb{G}}$ for the strict and closed extension problem, if there exists any, whatever it is, the elements of $\widehat{G}$ will be identified to those of $[H \times G]_{\sim}$ by the bijective function $\alpha$.

The following proposition imposes itself in this point and its proof, beyond strengthening our belief in the existence of a solution $\widehat{\mathbb{G}} = (\widehat{G}, \widehat{H})$ for our problem, it will also suggest to take $\widehat{G} = [H \times G]_{\sim}$.

## 2.3 Proposition

*Let $\mathbb{G} = (G, H)$ be an abelian, surjective, and with identity $S$-group. The binary relation, $\sim$, in $H \times G$, defined ahead, is an equivalence relation.*

**Definition of $\sim$:** for $(\Phi, g)$ and $(\Psi, h)$ in $H \times G$,

$$(\Phi, g) \sim (\Psi, h) \quad \text{if and only if} \quad \Psi^{-1}(g) - \Phi^{-1}(h) \subseteq N(\Phi\Psi).$$

**Remark.** Before demonstrating the proposition above, we list a set of results collected in Lemma 2.4 below, proved in the appendix at the end of this chapter, that will be used in the proof of several propositions, among which Proposition 2.3 above. Whenever one or more results provided by the referred lemma are employed in a proof, we will cite the used items or the lemma without specifying the items, or yet, in some cases, not even the lemma will be mentioned, at which cases the identification is left to the reader.

## 2.4 Lemma

*Let $\mathbb{G} = (G, H)$ be an abelian, surjective, and with identity $S$-group, $g, h \in G$, $\Phi, \Psi, \varphi \in H$ and $F, K, L, M \subseteq G$. Then:*

**(a)** *if $F \subseteq K$ and $L \subseteq M$, then $F \pm L \subseteq K \pm M$; in particular, since $M \subseteq M$, then $F \pm M \subseteq K \pm M$;*

**(b)** *if $F \subseteq K$, then $\Phi^{-1}(F) \subseteq \Phi^{-1}(K)$;*

**(c)** $\Phi^{-1}\Big(\Psi^{-1}(g)\Big) = \Psi^{-1}\Big(\Phi^{-1}(g)\Big)$ *and* $\Big(\Phi\Psi\Big)^{-1}(g) = \Psi^{-1}\Big(\Phi^{-1}(g)\Big)$;

**(d)** $\varphi^{-1}\Big(\Psi^{-1}(g)\Big) \pm \varphi^{-1}\Big(\Phi^{-1}(h)\Big) \subseteq \varphi^{-1}\Big(\Psi^{-1}(g) \pm \Phi^{-1}(h)\Big)$;

**(e)** $\varphi^{-1}\Big(\Psi^{-1}(g) \pm \Phi^{-1}(h)\Big) \subseteq \varphi^{-1}\Big(\Psi^{-1}(g)\Big) \pm \varphi^{-1}\Big(\Phi^{-1}(h)\Big) + N(\varphi)$;

**(f)** $N\Big(\Phi\Psi\Big) = \Phi^{-1}\Big(N(\Psi)\Big) + \Psi^{-1}\Big(N(\Phi)\Big) + N(\Phi) + N(\Psi)$;

**(g)** $N(\Psi\varphi\Phi) + N(\Psi) + N(\Phi) \subseteq N(\Psi\varphi\Phi)$;

**(h)** *if* $h = \Phi(g)$ *and* $g \in N(\Phi\Psi)$, *then,* $h \in N(\Psi)$.

*Proof (of Proposition 2.3).* Let us show first that $\sim$ is reflexive, that is, $(\Phi, g) \sim (\Phi, g)$ for every $(\Phi, g) \in H \times G$. Hence, we must prove that

$$\Phi^{-1}(g) - \Phi^{-1}(g) \subseteq N(\Phi\Phi) = N(\Phi^2).$$

Let then $h \in \Phi^{-1}(g) - \Phi^{-1}(g)$ be arbitrarily chosen, that is, $h = g' - g''$ where $g', g'' \in G_\Phi$ are such that $\Phi(g') = \Phi(g'') = g$. Therefore, $h \in G_\Phi$ and $\Phi(h) = \Phi(g') - \Phi(g'') = 0$. Consequently, $\Phi\left(\Phi(h)\right) = \Phi(0) = 0$ and hence $h \in N(\Phi\Phi) = N(\Phi^2)$. Since $h \in \Phi^{-1}(g) - \Phi^{-1}(g)$ is arbitrary, we conclude then that

$$\Phi^{-1}(g) - \Phi^{-1}(g) \subseteq N(\Phi\Phi) = N(\Phi^2).$$

The symmetry of $\sim$, that is, for $(\Phi, g)$ and $(\Psi, h)$ in $H \times G$, if $(\Phi, g) \sim (\Psi, h)$, then, $(\Psi, h) \sim (\Phi, g)$, follows from the following obvious fact:

$$\Psi^{-1}(g) - \Phi^{-1}(h) \subseteq N(\Phi\Psi) \quad \text{implies} \quad \Phi^{-1}(h) - \Psi^{-1}(g) \subseteq N(\Psi\Phi).$$

Remains then to prove that $\sim$ is transitive, that is, in the hypothesis of $(\Phi, g)$, $(\Psi, h)$ and $(\varphi, h')$ in $H \times G$ be such that

$$(\Phi, g) \sim (\Psi, h) \quad \text{and} \quad (\Psi, h) \sim (\varphi, h'),$$

then,

$$(\Phi, g) \sim (\varphi, h').$$

Suppose then that $(\Phi, g) \sim (\Psi, h)$ and $(\Psi, h) \sim (\varphi, h')$, that is:

$$\Psi^{-1}(g) - \phi^{-1}(h) \subseteq N(\Phi\Psi) \quad \text{and} \quad \varphi^{-1}(h) - \Psi^{-1}(h') \subseteq N(\Psi\varphi).$$

From this hypothesis and items (b), (a) and (d) of Lemma 2.4, result that:

$$\varphi^{-1}\Big(\Psi^{-1}(g)\Big) - \varphi^{-1}\Big(\Phi^{-1}(h)\Big) + \Phi^{-1}\Big(\varphi^{-1}(h)\Big) - \Phi^{-1}\Big(\Psi^{-1}(h')\Big) \subseteq$$
$$\subseteq \varphi^{-1}\Big(N(\Phi\Psi)\Big) + \Phi^{-1}\Big(N(\Psi\varphi)\Big).$$

From this last relation and items (c), (a) and (e) of Lemma 2.4 we obtain:

$$\Psi^{-1}\Big(\varphi^{-1}(g) - \Phi^{-1}(h')\Big) + \Phi^{-1}\Big(\varphi^{-1}(h) - \varphi^{-1}(h)\Big) \subseteq$$
$$\subseteq \varphi^{-1}\Big(N(\Phi\Psi)\Big) + \Phi^{-1}\Big(N(\Psi\varphi)\Big) + N(\Psi) + N(\Phi).$$

Now, it is easy to see that

$$\Psi^{-1}\Big(\varphi^{-1}(g) - \Phi^{-1}(h')\Big) \subseteq \Psi^{-1}\Big(\varphi^{-1}(g) - \Phi^{-1}(h')\Big) + \Phi^{-1}\Big(\varphi^{-1}(h) - \varphi^{-1}(h)\Big),$$

and thus we have

$$\Psi^{-1}\Big(\varphi^{-1}(g) - \Phi^{-1}(h')\Big) \subseteq \varphi^{-1}\Big(N(\Phi\Psi)\Big) + \Phi^{-1}\Big(N(\Psi\varphi)\Big) + N(\Psi) + N(\Phi).$$

On the other hand, clearly, we have

$$\varphi^{-1}\Big(N(\Phi\Psi)\Big) + \Phi^{-1}\Big(N(\Psi\varphi)\Big) \subseteq N(\varphi\Phi\Psi)$$

and hence we can write that

$$\Psi^{-1}\Big(\varphi^{-1}(g) - \Phi^{-1}(h')\Big) \subseteq N(\varphi\Phi\Psi) + N(\Psi) + N(\Phi)$$

and since, by item (g) of Lemma 2.4,

$$N(\varphi\Phi\Psi) + N(\Psi) + N(\Phi) \subseteq N(\varphi\Psi\Phi),$$

results that

$$\Psi^{-1}\Big(\varphi^{-1}(g) - \Phi^{-1}(h')\Big) \subseteq N(\varphi\Psi\Phi). \tag{2.3-1}$$

Let now $g' \in \varphi^{-1}(g) - \Phi^{-1}(h')$ be arbitrarily chosen. Since $\Psi\colon G_\Psi \longrightarrow G$ is a homomorphism over $G$, that is, $\Psi$ is surjective, then there exists $g'' \in G_\Psi$ such that $\Psi(g'') = g'$. Hence,

$$g'' \in \Psi^{-1}(g') \subseteq \Psi^{-1}\Big(\varphi^{-1}(g) - \Phi^{-1}(h')\Big)$$

and therefore, taking (2.3-1) into account,

$$g'' \in N(\varphi\Phi\Psi).$$

Hence, $g' = \Psi(g'')$ and $g'' \in N(\varphi\Phi\Psi)$ which, in turn, taking into account the item (h) of Lemma 2.4, informs us that $g' \in N(\varphi\Phi)$. Therefore,

$$\varphi^{-1}(g) - \Phi^{-1}(h') \subseteq N(\varphi\Phi)$$

that is,

$$(\Phi, g) \sim (\varphi, h'). \qquad \blacksquare$$

# The Abelian Group $[H \times G]_\sim$

## 2.5 Addition in $[H \times G]_\sim$

What information, beyond those obtained previously, namely, that $\sim$ is an equivalence relation in $H \times G$ and the function

$$\begin{aligned} \alpha : \qquad \widehat{G} &\longrightarrow [H \times G]_\sim \\ \widehat{g} = \widehat{\Phi}(g) &\longmapsto \alpha(\widehat{g}) = [\Phi, g] \end{aligned}$$

is bijective, could the hypothesis on the existence of a solution $\widehat{\mathbb{G}} = (\widehat{G}, \widehat{H})$ for the strict and closed extension problem provides us?

Well, if we intend, as suggest the considerations in 2.2, to make $[H \times G]_\sim$ a solution to the problem, that is, for it to play the role of $\widehat{G}$, then we need to define an operation in $[H \times G]_\sim$ that transforms it into an abelian group. Existing $\widehat{\mathbb{G}} = (\widehat{G}, \widehat{H})$, then, clearly, the bijectivity of the function $\alpha \colon \widehat{G} \longrightarrow [H \times G]_\sim$ assures not only the well-definition of the addition operation in $[H \times G]_\sim$ given by, for $\widetilde{g}$ and $\widetilde{h}$ in $[H \times G]_\sim$,

$$\widetilde{g} + \widetilde{h} \coloneqq \alpha(\widehat{g} + \widehat{h})$$

where

$$\widehat{g} = \alpha^{-1}(\widetilde{g}) \quad \text{and} \quad \widehat{h} = \alpha^{-1}(\widetilde{h}),$$

but also that $[H \times G]_\sim$ with this operation is an abelian group isomorphic to the group $\widehat{G}$ with $\alpha$ as an isomorphism.

Now, since $\widehat{g} = \widehat{\Phi}(g)$ and $\widehat{h} = \widehat{\Psi}(h)$ with $\Phi, \Psi \in H$ and $g, h \in G$ (remember that $\widehat{\mathbb{G}}$ is a closed extension of $\mathbb{G}$), then

$$\widehat{g} + \widehat{h} = \widehat{\Phi}(g) + \widehat{\Psi}(h),$$

and if $g^* \in \Psi^{-1}(g) \quad \text{and} \quad h^* \in \Phi^{-1}(h)$,

$$\widehat{g} + \widehat{h} = \widehat{\Phi}\Big(\Psi(g^*)\Big) + \widehat{\Psi}\Big(\Phi(h^*)\Big) = \widehat{\Phi}\Big(\widehat{\Psi}(g^*)\Big) + \widehat{\Psi}\Big(\widehat{\Phi}(h^*)\Big)$$

that is,

$$\widehat{g} + \widehat{h} = \widehat{\Phi\Psi}(g^* + h^*)$$

and from this results that:

$$\widetilde{g} + \widetilde{h} = \alpha(\widehat{g} + \widehat{h}) = [\Phi\Psi, g^* + h^*]$$

or yet, since $\widetilde{g} = \alpha(\widehat{g}) = \alpha(\widehat{\Phi}(g)) = [\Phi, g]$ and $\widetilde{h} = \alpha(\widehat{h}) = \alpha(\widehat{\Psi}(h)) = [\Psi, h]$,

$$[\Phi, g] + [\Psi, h] = [\Phi\Psi, g^* + h^*].$$

These considerations suggest to define addition in $[H \times G]_\sim$ as following.

**Addition in $[H \times G]_\sim$:** For $[\Phi, g]$ and $[\Psi, h]$ in $[H \times G]_\sim$ we define:

$$[\Phi, g] + [\Psi, h] \coloneqq [\Phi\Psi, g^* + h^*]$$

where $g^* \in \Psi^{-1}(g)$ and $h^* \in \Phi^{-1}(h)$.

Observing now that the definition above does not involve elements of $\widehat{\mathbb{G}} = (\widehat{G}, \widehat{H})$, but only those in $\mathbb{G} = (G, H)$, and that, as we saw, existing $\widehat{\mathbb{G}}$, $[H \times G]_\sim$ with the operation introduced by the referred definition, is an abelian group isomorphic to $\widehat{G}$, then, just like we were led to Proposition 2.3, about the relation $\sim$, here we are also led to question if the alluded addition in $[H \times G]_\sim$ is, itself, well-defined and, case affirmative, if $[H \times G]_\sim$ with this operation is in fact an abelian group. Not being this the case, then, the strict and closed extension problem will be solved since, as goes the popular Portuguese saying, "o que não tem solução resolvido está" (something like "what does not have a solution is solved"). It is then time and place for the following proposition.

## 2.6 Proposition

*Let $\mathbb{G} = (G, H)$ be an abelian, surjective, and with identity $S$-group. Let also $\sim$ be the equivalence relation in $H \times G$ defined in Proposition 2.3 and $[H \times G]_\sim$ the partition of $H \times G$ determined by $\sim$. Then, defining, for $[\Phi, g]$ and $[\Psi, h]$ in $[H \times G]_\sim$,*

$$[\Phi, g] + [\Psi, h] \coloneqq [\Phi\Psi, g^* + h^*]$$

*where*

$$g^* \in \Psi^{-1}(g) \quad \text{and} \quad h^* \in \Phi^{-1}(h),$$

*the (addition) operation is well-defined, that is, it is independent of the representatives taken in the equivalence classes $[\Phi, g]$ and $[\Psi, h]$, and, yet, $[H \times G]_\sim$ with this operation is an abelian group.*

*Proof.*

**(a)** The operation is well-defined.

We must prove that if $(\Phi_1, g_1) \in [\Phi, g]$ and $(\Psi_1, h_1) \in [\Psi, h]$, then $[\Phi_1\Psi_1, g_1^* + h_1^*] = [\Phi\Psi, g^* + h^*]$ where $g^* \in \Psi^{-1}(g)$, $h^* \in \Phi^{-1}(h)$, $g_1^* \in \Psi_1^{-1}(g_1)$ and $h_1^* \in \Phi_1^{-1}(h_1)$. Suppose then that

$$(\Phi_1, g_1) \in [\Phi, g] \quad \text{and} \quad (\Psi_1, h_1) \in [\Psi, h],$$

i.e., that

$$\Phi^{-1}(g_1) - \Phi_1^{-1}(g) \subseteq N(\Phi\Phi_1) \quad \text{and} \quad \Psi^{-1}(h_1) - \Psi_1^{-1}(h) \subseteq N(\Psi\Psi_1), \tag{2.6-1}$$

and, based on this hypothesis, let us prove that

$$[\Phi_1\Psi_1, g_1^* + h_1^*] = [\Phi\Psi, g^* + h^*]$$

or, equivalently, that

$$(\Phi_1\Psi_1)^{-1}(g^* + h^*) - (\Phi\Psi)^{-1}(g_1^* + h_1^*) \subseteq N(\Phi_1\Psi_1\Phi\Psi) \tag{2.6-2}$$

where

$$g^* \in \Psi^{-1}(g), \quad h^* \in \Phi^{-1}(h), \quad g_1^* \in \Psi_1^{-1}(g_1) \quad \text{and} \quad h_1^* \in \Phi_1^{-1}(h_1).$$

From (2.6-1), taking into account the items (b) and (a) of Lemma 2.4, results that

$$(\Psi\Psi_1)^{-1}\Big(\Phi^{-1}(g_1) - \Phi_1^{-1}(g)\Big) + (\Phi\Phi_1)^{-1}\Big(\Psi^{-1}(h_1) - \Psi_1^{-1}(h)\Big) \subseteq \\ \subseteq (\Psi\Psi_1)^{-1}\Big(N(\Phi\Phi_1)\Big) + (\Phi\Phi_1)^{-1}\Big(N(\Psi\Psi_1)\Big)$$

and hence, since by item (f) of the same lemma

$$(\Phi_1\Phi)^{-1}\Big(N(\Psi_1\Psi)\Big) + (\Psi_1\Psi)^{-1}\Big(N(\Phi_1\Phi)\Big) \subseteq N(\Phi_1\Phi\Psi_1\Psi),$$

comes that

$$(\Psi\Psi_1)^{-1}\Big(\Phi^{-1}(g_1) - \Phi_1^{-1}(g)\Big) + (\Phi\Phi_1)^{-1}\Big(\Psi^{-1}(h_1) - \Psi_1^{-1}(h)\Big) \subseteq N(\Phi_1\Phi\Psi\Psi_1),$$

which along with the item (d) of Lemma 2.4 allow us to conclude that

$$(\Psi\Psi_1)^{-1}\Big(\Phi^{-1}(g_1)\Big) + (\Phi\Phi_1)^{-1}\Big(\Psi^{-1}(h_1)\Big) - \left[(\Psi\Psi_1)^{-1}\Big(\Phi_1^{-1}(g)\Big) + (\Phi\Phi_1)^{-1}\Big(\Psi_1^{-1}(h)\Big)\right] \subseteq \\ \subseteq N(\Phi_1\Phi\Psi_1\Psi),$$

or yet, with item (c) of Lemma 2.4 in mind,

$$(\Phi\Psi)^{-1}\Big(\Psi_1^{-1}(g_1)\Big) + (\Phi\Psi)^{-1}\Big(\Phi_1^{-1}(h_1)\Big) - \left[(\Phi_1\Psi_1)^{-1}\Big(\Psi^{-1}(g)\Big) + (\Phi_1\Psi_1)^{-1}\Big(\Phi^{-1}(h)\Big)\right] \subseteq \\ \subseteq N(\Phi_1\Phi\Psi_1\Psi).$$

Taking

$$A := (\Phi\Psi)^{-1}\Big(\Psi_1^{-1}(g_1)\Big) + (\Phi\Psi)^{-1}\Big(\Phi_1^{-1}(h_1)\Big)$$

and

$$B := (\Phi_1\Psi_1)^{-1}\Big(\Psi^{-1}(g)\Big) + (\Phi_1\Psi_1)^{-1}\Big(\Phi^{-1}(h)\Big),$$

the last expression assumes the shorter form

$$A - B \subseteq N(\Phi_1\Phi\Psi_1\Psi),$$

from where we conclude (item (a) of Lemma 2.4) that

$$A + N(\Phi\Psi) - [B + N(\Phi_1\Psi_1)] \subseteq N(\Phi_1\Phi\Psi_1\Psi) + N(\Phi\Psi) - N(\Phi_1\Psi_1).$$

Now, taking into account the definitions of $A$ and $B$ given above, as well as item (e) of Lemma 2.4,

$$(\Phi\Psi)^{-1}\Big(\Psi_1^{-1}(g_1) + \Phi_1^{-1}(h_1)\Big) \subseteq A + N(\Phi\Psi)$$

and

$$(\Phi_1\Psi_1)^{-1}\Big(\Psi^{-1}(g) + \Phi^{-1}(h)\Big) \subseteq B + N(\Phi_1\Psi_1),$$

which, along with the last inclusion relation, allows us to conclude that:

$$(\Phi\Psi)^{-1}\Big(\Psi_1^{-1}(g_1) + \Phi_1^{-1}(h_1)\Big) - (\Phi_1\Psi_1)^{-1}\Big(\Psi^{-1}(g) + \Phi^{-1}(h)\Big) \subseteq$$
$$\subseteq N(\Phi_1\Phi\Psi_1\Psi) + N(\Phi\Psi) - N(\Phi_1\Psi_1).$$

Finally, since

$$N(\Phi_1\Phi\Psi_1\Psi) + N(\Phi\Psi) - N(\Phi_1\Psi_1) \subseteq N(\Phi_1\Phi\Psi_1\Psi),$$

we obtain the relation

$$(\Phi\Psi)^{-1}\Big(\Psi_1^{-1}(g_1) + \Phi_1^{-1}(h_1)\Big) - (\Phi_1\Psi_1)^{-1}\Big(\Psi^{-1}(g) + \Phi^{-1}(h)\Big) \subseteq N(\Phi_1\Phi\Psi_1\Psi),$$

from where we conclude that:

$$(\Phi_1\Psi_1)^{-1}(g^* + h^*) - (\Phi\Psi)^{-1}(g_1^* + h_1^*) \subseteq N(\Phi_1\Phi\Psi_1\Psi)$$

where

$$g^* \in \Psi^{-1}(g), \quad h^* \in \Phi^{-1}(h), \quad g_1^* \in \Psi_1^{-1}(g_1) \quad \text{and} \quad h_1^* \in \Phi_1^{-1}(h_1),$$

which is the relation (2.6-2).

**(b)** The operation is associative.

Let $\widetilde{f} := [\varphi, f]$, $\widetilde{g} := [\Phi, g]$ and $\widetilde{h} := [\Psi, h]$ be arbitrarily chosen elements in $[H \times G]_\sim$. We want to prove that

$$(\widetilde{f} + \widetilde{g}) + \widetilde{h} = \widetilde{f} + (\widetilde{g} + \widetilde{h}).$$

Let us calculate first $(\widetilde{f} + \widetilde{g}) + \widetilde{h}$. We have, from the definition of addition in $[H \times G]_\sim$, that

$$\widetilde{f} + \widetilde{g} = [\varphi, f] + [\Phi, g] = [\varphi\Phi, f^* + g^*]$$

where

$$f^* \in \Phi^{-1}(f) \quad \text{and} \quad g^* \in \varphi^{-1}(g),$$

that is,

$$\widetilde{f} + \widetilde{g} = [\varphi\Phi, l^*]$$

with

$$l^* := f^* + g^* \in \Phi^{-1}(f) + \varphi^{-1}(g)$$

arbitrarily chosen.

Resorting again to the definition of addition, we get that

$$\left(\tilde{f} + \tilde{g}\right) + \tilde{h} = [\varphi\Phi, l^*] + [\Psi, h] = [\varphi\Phi\Psi, m^* + h^*]$$

being

$$m^* \in \Psi^{-1}(l^*) \quad \text{and} \quad h^* \in (\varphi\Phi)^{-1}(h),$$

that is,

$$\left(\tilde{f} + \tilde{g}\right) + \tilde{h} = [\varphi\Phi\Psi, k^*]$$

with

$$k^* := m^* + h^* \in \Psi^{-1}(l^*) + (\varphi\Phi)^{-1}(h),$$

or yet, since $l^*$ is an arbitrarily chosen element of the set $\Phi^{-1}(f) + \varphi^{-1}(g)$, with

$$k^* \in \Psi^{-1}\left(\Phi^{-1}(f) + \varphi^{-1}(g)\right) + (\varphi\Phi)^{-1}(h)$$

arbitrarily chosen.

Analogously, we obtain that

$$\tilde{f} + \left(\tilde{g} + \tilde{h}\right) = [\varphi\Phi\Psi, u^*]$$

with

$$u^* \in (\Phi\Psi)^{-1}(f) + \varphi^{-1}\left(\Psi^{-1}(g) + \Phi^{-1}(h)\right)$$

arbitrarily chosen.

Taking now into account the item (d) of Lemma 2.4, we obtain the following inclusion relations:

$$\Psi^{-1}\left(\Phi^{-1}(f) + \varphi^{-1}(g)\right) + (\varphi\Phi)^{-1}(h) \supseteq \Psi^{-1}\left(\Phi^{-1}(f)\right) + \Psi^{-1}\left(\varphi^{-1}(g)\right) + (\varphi\Phi)^{-1}(h)$$

and

$$(\Phi\Psi)^{-1}(f) + \varphi^{-1}\left(\Psi^{-1}(g) + \Phi^{-1}(h)\right) \supseteq (\Phi\Psi)^{-1}(f) + \varphi^{-1}\left(\Psi^{-1}(g)\right) + \varphi^{-1}\left(\Phi^{-1}(h)\right).$$

On the other hand, the item (c) of the same lemma shows that

$$\Psi^{-1}\left(\Phi^{-1}(f)\right) + \Psi^{-1}\left(\varphi^{-1}(g)\right) + (\varphi\Phi)^{-1}(h) = (\Phi\Psi)^{-1}(f) + \varphi^{-1}\left(\Psi^{-1}(g)\right) + \varphi^{-1}\left(\Phi^{-1}(h)\right),$$

and hence we obtain

$$\Psi^{-1}\left(\Phi^{-1}(f) + \varphi^{-1}(g)\right) + (\varphi\Phi)^{-1}(h) \supseteq A$$

and

$$(\Phi\Psi)^{-1}(f) + \varphi^{-1}\Big(\Psi^{-1}(g) + \Phi^{-1}(h)\Big) \supseteq A$$

where

$$A := \Psi^{-1}\Big(\Phi^{-1}(f)\Big) + \Psi^{-1}\Big(\varphi^{-1}(g)\Big) + (\varphi\Phi)^{-1}(h) \neq \varnothing.$$

Since

$$k^* \in \Psi^{-1}\Big(\Phi^{-1}(f) + \varphi^{-1}(g)\Big) + (\varphi\Phi)^{-1}(h) \supseteq A$$

and

$$u^* \in (\Phi\Psi)^{-1}(f) + \varphi^{-1}\Big(\Psi^{-1}(g) + \Phi^{-1}(h)\Big) \supseteq A$$

are arbitrarily chosen, we can take them equal to any element of $A$, that is, $k^* = u^* = a \in A$ and hence conclude that

$$[\varphi\Phi\Psi, k^*] = [\varphi\Phi\Psi, u^*] = [\varphi\Phi\Psi, a],$$

that is,

$$\left(\widetilde{f} + \widetilde{g}\right) + \widetilde{h} = \widetilde{f} + \left(\widetilde{g} + \widetilde{h}\right).$$

**(c)** The operation is commutative.

The commutativity of the operation is a direct consequence of its definition and of $\mathbb{G} = (G, H)$ being an abelian $S$-group.

**(d)** The existence of an additive neutral element.

Since $\mathbb{G} = (G, H)$ is a $S$-group with identity, we have that $[I_G, 0]$ is an element of $[H \times G]_\sim$ and, in fact, an additive neutral, that is, if $[\Phi, g] \in [H \times G]_\sim$,

$$[\Phi, g] + [I_G, 0] = [\Phi, g],$$

since, from the definition of addition,

$$[\Phi, g] + [I_G, 0] = [\Phi, g + h^*]$$

where $h^* \in N(\Phi)$. However, it is easy to see that,

$$(\Phi, g) \sim (\Phi, g + h^*),$$

that is,

$$[\Phi, g] = [\Phi, g + h^*],$$

and, hence, $[\Phi, g] + [I_G, 0] = [\Phi, g]$.

**(e)** The existence of inverse elements.

We must prove that if $\widetilde{g} = [\Phi, g] \in [H \times G]_\sim$, then there exists $-\widetilde{g} \in [H \times G]_\sim$ such that

$$\widetilde{g} + (-\widetilde{g}) = [I_G, 0].$$

Let then $\widetilde{g} = [\Phi, g] \in [H \times G]_\sim$ be arbitrarily fixed. Now define $-\widetilde{g}$ by:

$$-\widetilde{g} \coloneqq [\Phi, -g].$$

Hence, we have that

$$\widetilde{g} + (-\widetilde{g}) = [\Phi\Phi, g^* + h^*]$$

where

$$g^* \in \Phi^{-1}(g) \quad \text{and} \quad h^* \in \Phi^{-1}(-g).$$

It is easy to verify that

$$(\Phi\Phi, g^* + h^*) \sim (I_G, 0),$$

and hence we have

$$[\Phi\Phi, g^* + h^*] = [I_G, 0].$$

Therefore,

$$\widetilde{g} + (-\widetilde{g}) = [\Phi\Phi, g^* + h^*] = [I_G, 0].$$ ■

## 2.7 $G$ as a Subgroup of $[H \times G]_\sim$

Concerning the strict and closed extension problem, we find ourselves in the following stage: regardless of whether the problem has or not a solution, the set $[H \times G]_\sim$, described in Proposition 2.6, along with the addition operation defined in the same proposition, is an abelian group that, if there exists a solution $\widehat{\mathbb{G}} = (\widehat{G}, \widehat{H})$ of the referred problem, it will be isomorphic to the group $\widehat{G}$, with the function

$$\begin{aligned} \alpha : \qquad \widehat{G} &\longrightarrow [H \times G]_\sim \\ \widehat{g} = \widehat{\Phi}(g) &\longmapsto \alpha(\widehat{g}) = [\Phi, g] \end{aligned}$$

as an isomorphism.

Well, in this case, of the existence of $\widehat{\mathbb{G}}$, $G$ is a subgroup of $\widehat{G}$ and, then,

$$\alpha(G) = \Big\{\alpha(g) \colon g \in G\Big\} \subseteq [H \times G]_\sim$$

is a subgroup of $[H \times G]_\sim$ isomorphic to the group $G$, with the restriction of $\alpha$ to $G$, $\alpha|_G$, being an isomorphism from $G$ onto $\alpha(G)$. On the other hand, for every $g \in G \subseteq \widehat{G}$ we have $g = I_G(g) = \widehat{I}_G(g)$ and, therefore, $\alpha(g) = [I_G, g]$, that is

$$\alpha(G) = \Big\{\alpha(g) \colon g \in G\Big\} = \Big\{[I_G, g] \colon g \in G\Big\}$$

and hence we have a description of the set $\alpha(G) \subseteq [H \times G]_\sim$ independent of the function $\alpha$, that is, independent of the existence of $\widehat{\mathbb{G}}$, given exclusively in terms of $[H \times G]_\sim$. Consequently, the restriction of $\alpha$ to $G$, $\alpha|_G$, also becomes independent of $\widehat{\mathbb{G}}$, since, both its domain, $G$, and its image, $\alpha(G)$, and the association rule, $g \in G \longmapsto [I_G, g] \in \alpha(G)$, makes no allusion to $\widehat{\mathbb{G}}$. Hence, again, we have “ou pipoca ou piruá” (something like “life or death”): either the set $\alpha(G)$, henceforth denoted by the more suggestive symbol $[I_G \times G]_\sim$, that is,

$$[I_G \times G]_\sim := \alpha(G) = \Big\{[I_G, g] \colon g \in G\Big\} \subseteq [H \times G]_\sim,$$

is a subgroup of $[H \times G]_\sim$ isomorphic to the group $G$ with

$$\begin{aligned} \alpha|_G : G &\longrightarrow [I_G \times G]_\sim \\ g &\longmapsto \alpha|_G(g) = [I_G, g] \end{aligned}$$

as an isomorphism from $G$ onto $[I_G \times G]_\sim$ and, in this case, the hope of existence of a solution $\widehat{\mathbb{G}}$ for our problem increase, or, otherwise, “piruá” (“death”). Thus, the next step of our walk is established: to decide the truthiness of Proposition 2.8 stated below.

## 2.8 Proposition

*Let $[H \times G]_\sim$ be the abelian group described in Proposition 2.6 and $[I_G \times G]_\sim$ be the subset of $[H \times G]_\sim$ defined by:*

$$[I_G \times G]_\sim := \Big\{[I_G, g] \colon g \in G\Big\}$$

*where $I_G \in H$ is the identity function on $G$. Then, $[I_G \times G]_\sim$ is a subgroup of $[H \times G]_\sim$ isomorphic to the group $G$ with the function $\gamma$ defined below as an isomorphism:*

$$\begin{aligned} \gamma : G &\longrightarrow [I_G \times G]_\sim \\ g &\longmapsto \gamma(g) := [I_G, g]. \end{aligned}$$

*In other words, $[H \times G]_\sim$ is an abelian group which admits, more precisely, $\gamma$-admits $G$ as a subgroup (see Definition 1.6(a)).*

*Proof.* Note first that, for any $g$ and $h$ in $G$,

$$\gamma(g) + \gamma(h) = [I_G, g] + [I_G, h] = [I_G, g + h] = \gamma(g + h),$$

that is, considering $\gamma$ as a function from the group $G$ into the group

$$[H \times G]_\sim \supseteq [I_G \times G]_\sim,$$

$\gamma$ is a homomorphism and, therefore, its image $\gamma(G) = [I_G \times G]_\sim \subseteq [H \times G]_\sim$ is a subgroup of $[H \times G]_\sim$. Hence, $\gamma$ is a surjective homomorphism from the group $G$ onto the group $[I_G \times G]_\sim$. Beyond that, for $g, h \in G$ such that $\gamma(g) = \gamma(h)$, that is, $[I_G, g] = [I_G, h]$, we have $(I_G, g) \sim (I_G, h)$, that is, $g - h \in N(I_G)$ which is equivalent to $g = h$. Thus, $\gamma$ is injective and, therefore, an isomorphism from the group $G$ onto the subgroup $[I_G \times G]_\sim$ of $[H \times G]_\sim$. ■

# $\widetilde{\mathbb{G}}$, an Extension of $\mathbb{G}$

## 2.9 The Prolongation $\widetilde{H}$ of $H$

Existing or not a strict and closed extension of a given abelian, surjective, and with identity $S$-group $\mathbb{G} = (G, H)$, we saw until now that: the relation $\sim$ defined in Proposition 2.3 is an equivalence relation in $H \times G$, the set of the equivalence classes determined by $\sim$, $[H \times G]_\sim$, with the addition operation given in Proposition 2.6, is an abelian group and that this group admits $G$ as a subgroup or, more specifically, $G$ is isomorphic to the subgroup $[I_G \times G]_\sim = \{[I_G, g] \colon g \in G\}$ of the group $[H \times G]_\sim$, with

$$\begin{aligned} \gamma : G &\longrightarrow [I_G \times G]_\sim \\ g &\longmapsto \gamma(g) = [I_G, g] \end{aligned}$$

as an isomorphism.

On the other hand, we also saw that, if there exists a strict and closed extension, $\widehat{\mathbb{G}} = (\widehat{G}, \widehat{H})$, of $\mathbb{G} = (G, H)$, then the group $\widehat{G}$ is isomorphic to the group $[H \times G]_\sim$, with

$$\begin{aligned} \alpha : \qquad \widehat{G} &\longrightarrow [H \times G]_\sim \\ \widehat{g} = \widehat{\Phi}(g) &\longmapsto \alpha(\widehat{g}) = [\Phi, g] \end{aligned}$$

as an isomorphism such that its restriction to $G \subseteq \widehat{G}$ is $\alpha\big|_G = \gamma$ and, therefore, keeps fixed the elements of $G$ since, if $g \in G$,

$$\alpha(g) = \alpha\big|_G(g) = \gamma(g) = [I_G, g] \equiv g$$

($g \in G$ is identified to $[I_G, g] \in [I_G \times G]_\sim$, $g \equiv [I_G, g]$, via the isomorphism $\gamma \colon G \longrightarrow [I_G \times G]_\sim$).

Let us now, again admitting the existence of a solution $\widehat{\mathbb{G}} = (\widehat{G}, \widehat{H})$ for our problem, verify if this hypothesis can lead us, as it led, inevitably, to the abelian group $[H \times G]_\sim$ isomorphic to $\widehat{G}$, also to necessary prolongations $\widetilde{\Phi}$ to $[H \times G]_\sim$ of the homomorphisms $\Phi \in H$ (see Definition 1.6(b)); in other words, we will investigate if the hypothesis of existence of a solution $\widehat{\mathbb{G}} = (\widehat{G}, \widehat{H})$ suggests any semigroup $\widetilde{H}$ of endomorphisms on $[H \times G]_\sim$, which can play the role of the semigroup $\widehat{H}$ of the hypothetical solution $\widehat{\mathbb{G}}$, just like it suggested the group $[H \times G]_\sim$ as a protagonist for $\widehat{G}$.

It is worth noting that the identification $g \equiv [I_G, g]$, of $g \in G$ with $[I_G, g] \in [I_G \times G]_\sim$ or, more generally, of a subset $G'$ of $G$ with

$$[I_G \times G']_\sim := \Big\{[I_G, g] \colon g \in G'\Big\} \subseteq [I_G \times G]_\sim,$$

realized by the isomorphism

$$\begin{aligned} \gamma : G &\longrightarrow [I_G \times G]_\sim \\ g &\longmapsto \gamma(g) = [I_G, g], \end{aligned}$$

allows to express (identify) each homomorphism

$$\begin{aligned}\Phi : G_\Phi &\longrightarrow G\\ g &\longmapsto \Phi(g)\end{aligned}$$

of $H$ in the following form:

$$\begin{aligned}\Phi : [I_G \times G_\Phi]_\sim \subseteq [I_G \times G]_\sim &\longrightarrow [I_G \times G]_\sim\\ [I_G, g] &\longmapsto \Phi\Big([I_G, g]\Big) := \Big[I_G, \Phi(g)\Big].\end{aligned}$$

Thus, $\widetilde{\Phi}\colon [H \times G]_\sim \longrightarrow [H \times G]_\sim$ will be a prolongation of $\Phi \in H$ to $[H \times G]_\sim$ (see Definition 1.6(b)) if and only if it is an endomorphism on $[H \times G]_\sim$ and

$$\widetilde{\Phi}\Big([I_G, g]\Big) = \Phi\Big([I_G, g]\Big) = \Big[I_G, \Phi(g)\Big]$$

for every $g \in G_\Phi$.

In the hypothesis of existence of $\widehat{\mathbb{G}} = (\widehat{G}, \widehat{H})$ we have, for each $\Phi \in H$, a prolongation, $\widehat{\Phi}$, to the group $\widehat{G}$, with this group being isomorphic to $[H \times G]_\sim$ with

$$\begin{aligned}\alpha : \quad \widehat{G} &\longrightarrow [H \times G]_\sim\\ \widehat{g} = \widehat{\Phi}(g) &\longmapsto \alpha(\widehat{g}) = [\Phi, g]\end{aligned}$$

as an isomorphism. Thus, through the isomorphism $\alpha$, a natural definition for a prolongation $\widetilde{\Phi}$ of $\Phi\colon [I_G \times G_\Phi]_\sim \equiv G_\Phi \longrightarrow [I_G \times G]_\sim \equiv G$ to $[H \times G]_\sim$, which imposes itself in this context, represented in the diagram of the figure below, is the following:

$$\widetilde{\Phi} := \alpha\widehat{\Phi}\alpha^{-1},$$

that is,

$$\begin{aligned}\widetilde{\Phi} : [H \times G]_\sim &\longrightarrow [H \times G]_\sim\\ [\Psi, g] &\longmapsto \widetilde{\Phi}\Big([\Psi, g]\Big) := [\Phi\Psi, g].\end{aligned}$$

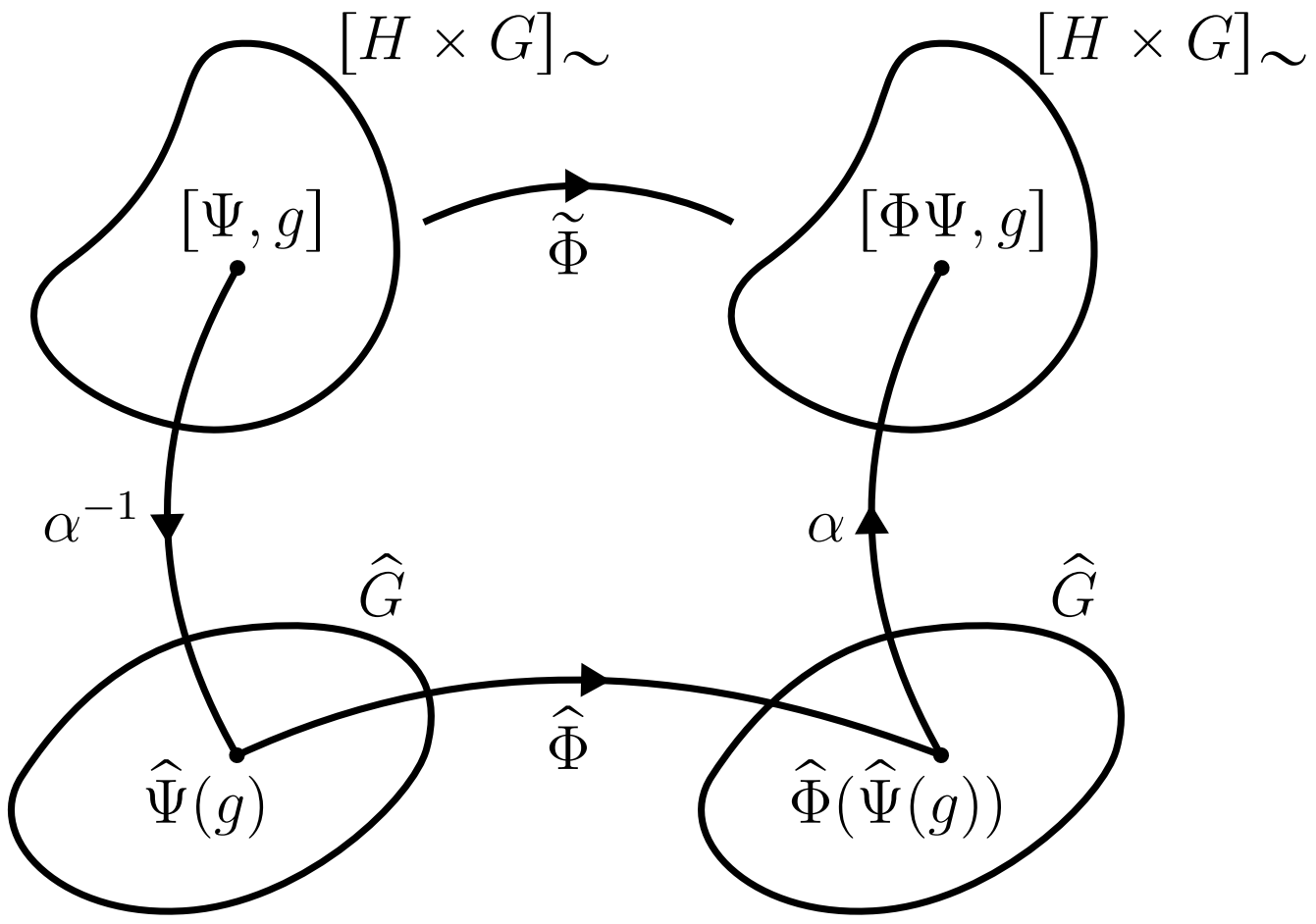

Note that $\widetilde{\Phi}$ is in fact a prolongation of $\Phi$, since:

**(a)** $\widetilde{\Phi}$ is an endomorphism on $[H\times G]_\sim$ since, being $[\Psi,g]$ and $[\varphi,h]$ any elements of $[H\times G]_\sim$,

$$\begin{aligned}\widetilde{\Phi}\Big([\Psi,g]+[\varphi,h]\Big) &= \alpha\left(\widehat{\Phi}\left(\alpha^{-1}\Big([\Psi,g]+[\varphi,h]\Big)\right)\right)=\\ &= \alpha\left(\widehat{\Phi}\left(\alpha^{-1}\Big([\Psi,g]\Big)+\alpha^{-1}\Big([\varphi,h]\Big)\right)\right) =\\ &= \alpha\left(\widehat{\Phi}\Big(\widehat{\Psi}(g)+\widehat{\varphi}(h)\Big)\right)=\\ &= \alpha\Big(\widehat{\Phi\Psi}(g)+\widehat{\Phi\varphi}(h)\Big)=\\ &= \alpha\Big(\widehat{\Phi\Psi}(g)\Big)+\alpha\Big(\widehat{\Phi\varphi}(h)\Big) =\\ &= [\Phi\Psi,g]+[\Phi\varphi,h] =\\ &= \widetilde{\Phi}\Big([\Psi,g]\Big)+\widetilde{\Phi}\Big([\varphi,h]\Big);\end{aligned}$$

**(b)** $\Phi$ is the restriction of $\widetilde{\Phi}$ to $G_\Phi \equiv [I_G\times G_\Phi]_\sim$, since for $g\equiv[I_G,g]\in[I_G\times G_\Phi]_\sim\equiv G_\Phi$ we have:

$$\begin{aligned}\widetilde{\Phi}\Big([I_G,g]\Big) &= \alpha\left(\widehat{\Phi}\left(\alpha^{-1}\Big([I_G,g]\Big)\right)\right) = \alpha\left(\widehat{\Phi}\Big(\widehat{I}_G(g)\Big)\right) = \alpha\left(\widehat{I}_G\Big(\widehat{\Phi}(g)\Big)\right) =\\ &= \alpha\Big(\widehat{I}_G\Big(\Phi(g)\Big)\Big) = \Big[I_G,\Phi(g)\Big] \equiv\\ &\equiv \Phi(g);\end{aligned}$$

and finally

**(c)** $\widetilde{\Phi\Psi}=\widetilde{\Phi}\widetilde{\Psi}$. In fact, we have

$$\widetilde{\Phi\Psi}=\alpha(\widehat{\Phi\Psi})\alpha^{-1}=\alpha(\widehat{\Phi}\widehat{\Psi})\alpha^{-1}=(\alpha\widehat{\Phi}\alpha^{-1})(\alpha\widehat{\Psi}\alpha^{-1})=\widetilde{\Phi}\widetilde{\Psi}.$$

Taking into account that $H$ is a semigroup, (a), (b) and (c) above ensure that the set

$$\widetilde{H}:=\left\{\widetilde{\Phi}=\alpha\widehat{\Phi}\alpha^{-1} : \Phi\in H\right\}$$

is a prolongation of $H$ to $[H\times G]_\sim$ (see Definition 1.6(c)) and, thus, the existence of $\widehat{\mathbb{G}}$ implies

$$\widetilde{\mathbb{G}}:=\Big([H\times G]_\sim,\widetilde{H}\Big)$$

to be an extension of the $S$-group $\mathbb{G}=(G,H)$.

On the other hand, since $\widetilde{\Phi}$ has a definition independent of $\widehat{\mathbb{G}}$, namely,

$$\widetilde{\Phi}\Big([\Psi,g]\Big):=[\Phi\Psi,g] \quad \text{with} \quad [\Psi,g]\in[H\times G]_\sim,$$

fully expressed in $\mathbb{G} = (G, H)$, becomes important to known if, without the hypothesis of existence of $\widehat{\mathbb{G}}$, $\widetilde{\Phi}$ is intrinsically well-defined and, case affirmative, if it is really an endomorphism on $[H \times G]_{\sim}$ such that $\widetilde{H} = \{\widetilde{\Phi}\colon \Phi \in H\}$ be a prolongation of $H$ to $[H \times G]_{\sim}$. Here is the content of our next proposition.

## 2.10 Proposition

*Let $\mathbb{G} = (G, H)$ and $[H \times G]_{\sim}$ be as in Proposition 2.6. Let also, for each $\Phi \in H$, $\widetilde{\Phi}$ be defined by:*

$$\begin{aligned} \widetilde{\Phi} : [H \times G]_{\sim} &\longrightarrow [H \times G]_{\sim} \\ [\Psi, g] &\longmapsto \widetilde{\Phi}\Big([\Psi, g]\Big) := [\Phi\Psi, g]. \end{aligned}$$

*Then, $\widetilde{\Phi}$ is well-defined (it is independent of the representative taken into the equivalence class $[\Psi, g]$), the set*

$$\widetilde{H} := \Big\{\widetilde{\Phi}\colon \Phi \in H\Big\},$$

*with the usual operation of composition of functions, is a prolongation of $H$ to $[H \times G]_{\sim}$ (see Definition 1.6(c)) and, therefore,*

$$\widetilde{\mathbb{G}} := \Big([H \times G]_{\sim}, \widetilde{H}\Big)$$

*is an extension of the $S$-group $\mathbb{G} = (G, H)$.*

*Proof.* **(a)** $\widetilde{\Phi}$ is well defined.

In fact, if $(\varphi, h) \in [\Psi, g] \in [H \times G]_{\sim}$, that is, if $(\varphi, h) \sim (\Psi, g)$, then,

$$\Psi^{-1}(h) - \varphi^{-1}(g) \subseteq N(\Psi\varphi)$$

and, therefore, taking into account Lemma 2.4,

$$(\Phi\Psi)^{-1}(h) - (\Phi\varphi)^{-1}(g) \subseteq \Phi^{-1}\Big(N(\Psi\varphi)\Big) \subseteq N(\Phi\Psi\Phi\varphi)$$

and hence we see that

$$(\Phi\varphi, h) \sim (\Phi\Psi, g).$$

Therefore,

$$[\Phi\varphi, h] = [\Phi\Psi, g].$$

**(b)** $\widetilde{\Phi}$ is an endomorphism on $[H \times G]_{\sim}$.

In the proof of this claim we will use the following result, whose prove is immediate:

$$\text{if} \quad \chi \in H, (\Delta, u) \in H \times G \quad \text{and} \quad u^* \in \chi^{-1}(u), \quad \text{then} \quad (\Delta, u) \sim (\Delta\chi, u^*). \tag{2.10-1}$$

Let now $\widetilde{g} = [\Psi, g]$ and $\widetilde{h} = [\varphi, h]$ be elements of $[H \times G]_\sim$ arbitrarily fixed. We then have that

$$\widetilde{g} + \widetilde{h} = [\Psi\varphi, g^* + h^*]$$

where

$$g^* \in \varphi^{-1}(g) \quad \text{and} \quad h^* \in \Psi^{-1}(h).$$

Therefore,

$$\widetilde{\Phi}(\widetilde{g} + \widetilde{h}) = [\Phi\Psi\varphi, g^* + h^*].$$

On the other hand, taking into account (2.10-1),

$$\widetilde{g} = [\Psi, g] = [\Psi\varphi, g^*]$$

and

$$\widetilde{h} = [\varphi, h] = [\Psi\varphi, h^*].$$

Thus,

$$\widetilde{\Phi}(\widetilde{g}) + \widetilde{\Phi}(\widetilde{h}) = \Big[(\Phi\Psi\varphi)^2, g^{**} + h^{**}\Big],$$

where

$$g^{**} \in (\Phi\Psi\varphi)^{-1}(g^*) \quad \text{and} \quad h^{**} \in (\Phi\Psi\varphi)^{-1}(h^*)$$

and, therefore, resorting to Lemma 2.4, we get that

$$g^{**} + h^{**} \in (\Phi\Psi\varphi)^{-1}(g^* + h^*).$$

From this, we can apply (2.10-1) again, now with $\Phi\Psi\varphi$ both as $\chi$ and $\Delta$, to obtain:

$$\Big(\Phi\Psi\varphi, g^* + h^*\Big) \sim \Big((\Phi\Psi\varphi)^2, g^{**} + h^{**}\Big)$$

or, equivalently,

$$\widetilde{\Phi}\Big(\widetilde{g} + \widetilde{h}\Big) = \widetilde{\Phi}\Big(\widetilde{g}\Big) + \widetilde{\Phi}\Big(\widetilde{h}\Big).$$

**(c)** $\Phi$ is the restriction of $\widetilde{\Phi}$ to $G_\Phi$.

Considering the isomorphism

$$\begin{aligned} \gamma : G &\longrightarrow [I_G \times G]_\sim \\ g &\longmapsto \gamma(g) = [I_G, g] \end{aligned}$$

from the group $G$ onto the subgroup $[I_G \times G]_\sim$ of $[H \times G]_\sim$ (see Proposition 2.8), which realizes the identifications

$$g \equiv [I_G, g] \quad \text{and} \quad G_\Phi \equiv \gamma(G_\Phi) = \Big\{[I_G, g] : g \in G_\Phi\Big\},$$

to state that "$\Phi$ is the restriction of $\widetilde{\Phi}$ to $G_\Phi$" means, precisely, to say that

$$\widetilde{\Phi}\Big(\gamma(g)\Big) = \gamma\Big(\Phi(g)\Big)$$

that is,

$$\widetilde{\Phi}\Big([I_G, g]\Big) = \Big[I_G, \Phi(g)\Big] \quad \text{for every} \quad g \in G_\Phi,$$

or yet, taking into account the definition of $\widetilde{\Phi}$, that

$$[\Phi, g] = \Big[I_G, \Phi(g)\Big] \quad \text{for every} \quad g \in G_\Phi.$$

We must then prove that, if $g \in G_\Phi$,

$$(\Phi, g) \sim \Big(I_G, \Phi(g)\Big).$$

Now, this last relation can be obtained from (2.10-1) taking $\chi = \Phi$, $\Delta = I_G$, $u = \Phi(g)$ and $u^* = g$ (which is possible since $g \in \Phi^{-1}(\Phi(g))$).

**(d)** $\widetilde{\Phi\Psi} = \widetilde{\Phi}\widetilde{\Psi}$ for every $\Phi, \Psi \in H$.

In fact, if $\widetilde{g} = [\varphi, g] \in [H \times G]_\sim$ is arbitrarily fixed, we have

$$\begin{aligned}\widetilde{\Phi\Psi}\Big(\widetilde{g}\Big) &= \Big[(\Phi\Psi)\varphi, g\Big] = \Big[\Phi(\Psi\varphi), g\Big] = \\ &= \widetilde{\Phi}\Bigg([\Psi\varphi, g]\Bigg) = \widetilde{\Phi}\Bigg(\widetilde{\Psi}\Big([\varphi, g]\Big)\Bigg) = \\ &= \widetilde{\Phi}\Big(\widetilde{\Psi}(\widetilde{g})\Big) = \Big(\widetilde{\Phi}\widetilde{\Psi}\Big)(\widetilde{g})\end{aligned}$$

and hence, since $\widetilde{g}$ is arbitrary, we obtain

$$\widetilde{\Phi\Psi} = \widetilde{\Phi}\widetilde{\Psi}.$$

It results from (b), (c) and (d) above that

$$\widetilde{H} = \Big\{\widetilde{\Phi} \colon \Phi \in H\Big\}$$

is a prolongation of $H$ to $[H \times G]_\sim$ (see Definition 1.6(c)) and, therefore, that

$$\widetilde{\mathbb{G}} = \Big([H \times G]_\sim, \widetilde{H}\Big)$$

is an extension of the $S$-group $\mathbb{G} = (G, H)$ (see Definition 1.6(d)). ∎

# $\widetilde{\mathbb{G}}$, a Strict Extension of $\mathbb{G}$

## 2.11 The Strict Prolongation, $\widetilde{H}$, of $H$

The considerations made in 2.9, which led us to formulate Proposition 2.10, can be summarized as follows: in the hypothesis of existence of a strict and closed extension, $\widehat{\mathbb{G}} = (\widehat{G}, \widehat{H})$, of the abelian, surjective, and with identity $S$-group $\mathbb{G} = (G, H)$, and being $\widetilde{H} = \left\{\widetilde{\Phi} \colon \Phi \in H\right\}$ the prolongation of $H$ to $[H \times G]_{\sim}$ described in the referred proposition, then, the function

$$\begin{aligned} \beta : \widehat{H} &\longrightarrow \widetilde{H} \\ \widehat{\Phi} &\longmapsto \beta(\widehat{\Phi}) := \alpha\widehat{\Phi}\alpha^{-1} = \widetilde{\Phi} \end{aligned}$$

is an isomorphism from the semigroup $\widehat{H}$ onto the semigroup $\widetilde{H}$, where

$$\begin{aligned} \alpha : \qquad \widehat{G} &\longrightarrow [H \times G]_{\sim} \\ \widehat{g} = \widehat{\Phi}(g) &\longmapsto \alpha(\widehat{g}) = [\Phi, g] \end{aligned}$$

is, as we saw, an isomorphism (between groups) that keeps fixed the elements of $G \subseteq \widehat{G}$ (that is, $\alpha(g) = [I_G, g] \equiv g$ for every $g \in G$).

Now we ask: Is the property "to be a strict prolongation" preserved by the isomorphism $\beta$? More precisely, we inquire: For $\widehat{\Phi} \in \widehat{H}$, $\beta(\widehat{\Phi}) = \widetilde{\Phi}$ is, such as $\widehat{\Phi}$, a strict prolongation of $\Phi$? The answer is affirmative, since given $g \in G$ such that $g \notin G_\Phi \subseteq G$, as $\alpha \colon \widehat{G} \longrightarrow [H \times G]_{\sim}$ keeps fixed the elements of $G$, then, $\alpha^{-1}([I_G, g] \equiv g) \notin G_\Phi$ and, therefore, we have $\widehat{\Phi}(\alpha^{-1}([I_G, g])) \notin G$ since $\widehat{\Phi}$ is a strict prolongation of $\Phi$. Hence,

$$\alpha\left(\widehat{\Phi}\left(\alpha^{-1}\Big([I_G, g]\Big)\right)\right) \notin [I_G \times G]_{\sim} \equiv G,$$

that is

$$\widetilde{\Phi}\Big([I_G, g]\Big) \notin [I_G \times G]_{\sim}$$

which means that $\widetilde{\Phi}$ is a strict prolongation of $\Phi$.

Without the hypothesis of the existence of $\widehat{\mathbb{G}}$, in particular of the existence of the isomorphism $\beta \colon \widehat{H} \longrightarrow \widetilde{H}$ through which we proved above that $\widetilde{\Phi}$ is a strict prolongation, but resorting only to the definition of $\widetilde{\Phi}$, not the one given by $\widetilde{\Phi} = \beta(\widehat{\Phi})$, but its version independent of $\beta$ (and thus of $\widehat{\mathbb{G}}$), namely,

$$\widetilde{\Phi}\,([\Psi, u]) = [\Phi\Psi, u] \quad \text{for} \quad [\Psi, u] \in [H \times G]_{\sim},$$

will we be able to construct a proof that "to be a strict prolongation" (of $\Phi$) is a property of $\widetilde{\Phi}$? The proposition below provides the answer.

## 2.12 Proposition

*The extension $\widetilde{\mathbb{G}} = ([H \times G]_{\sim}, \widetilde{H})$, described in Proposition 2.10, of the abelian, surjective, and with identity $S$-group $\mathbb{G} = (G, H)$, is a strict extension of $\mathbb{G}$.*

*Proof.* Proposition 1.12 establishes that $\widetilde{\mathbb{G}}$ will be a strict extension of $\mathbb{G}$ if and only if, for each $\Phi \in H$,

$$G \cap N(\widetilde{\Phi}) = N(\Phi)$$

or, more precisely, if and only if

$$\gamma(G) \cap N(\widetilde{\Phi}) = \gamma\Big(N(\Phi)\Big),$$

being $\gamma \colon G \longrightarrow [I_G \times G]_\sim$ the isomorphism described in Proposition 2.8 through which each element $g \in G$ is identified to $\gamma(g) = [I_G, g]$. Let then $\widetilde{g} \in \gamma(G) \cap N(\widetilde{\Phi})$, that is,

$$\widetilde{g} = [I_G, g] \in \gamma(G) \cap N(\widetilde{\Phi})$$

be arbitrarily chosen. In this case, $g \in G$ and, remembering that $[I_G, 0]$ is the neutral element of the group $[H \times G]_\sim$,

$$\widetilde{\Phi}(\widetilde{g}) = [\Phi, g] = [I_G, 0],$$

and hence we have

$$(\Phi, g) \sim (I_G, 0)$$

or, equivalently,

$$g - N(\Phi) \subseteq N(\Phi),$$

from where we obtain that $g \in N(\Phi)$ and, therefore,

$$\widetilde{g} = [I_G, g] \in \gamma\Big(N(\Phi)\Big).$$

Since $\widetilde{g} \in \gamma(G) \cap N(\widetilde{\Phi})$ is arbitrary, we have

$$\gamma(G) \cap N\left(\widetilde{\Phi}\right) \subseteq \gamma\Big(N(\Phi)\Big).$$

Conversely, let now $\widetilde{g} \in \gamma\Big(N(\Phi)\Big)$ be arbitrarily fixed. Then we have

$$\widetilde{g} = [I_G, g] \in \gamma(G) \quad \text{and} \quad g \in N(\Phi).$$

With $g \in N(\Phi)$, then,

$$g - N(\Phi) \subseteq N(\Phi),$$

that is,

$$I_G^{-1}(g) - \Phi^{-1}(0) \subseteq N(\Phi),$$

which means

$$(\Phi, g) \sim (I_G, 0)$$

and, therefore,

$$[\Phi, g] = [I_G, 0],$$

or yet,

$$\widetilde{\Phi}\Big([I_G, g]\Big) = \widetilde{\Phi}(\widetilde{g}) = [I_G, 0],$$

that is, $\widetilde{g} \in N(\widetilde{\Phi})$. Since $\widetilde{g} \in \gamma(N(\Phi))$ is arbitrary we conclude that

$$\gamma\Big(N(\Phi)\Big) \subseteq \gamma(G) \cap N(\widetilde{\Phi}).$$

From this, since $\gamma(G) \cap N(\widetilde{\Phi}) \subseteq \gamma(N(\Phi))$, we conclude that

$$\gamma(G) \cap N(\widetilde{\Phi}) = \gamma\Big(N(\Phi)\Big). \qquad \blacksquare$$

# $\widetilde{\mathbb{G}}$, a Strict and Closed Extension of $\mathbb{G}$

## 2.13 Is $\widetilde{\mathbb{G}}$ a Closed Extension of $\mathbb{G}$?

From Propositions 2.10 and 2.12 we have that $\widetilde{\mathbb{G}} = ([H \times G]_{\sim}, \widetilde{H})$ is a strict extension of $\mathbb{G} = (G, H)$. Is $\widetilde{\mathbb{G}}$ a closed extension of $\mathbb{G}$?

Let us note that, if there exists a strict and closed extension, $\widehat{\mathbb{G}} = (\widehat{G}, \widehat{H})$, of $\mathbb{G}$, then, necessarily, $\widetilde{\mathbb{G}}$ is also a closed extension of $\mathbb{G}$. In fact, taking into account that, in this case (of the existence of $\widetilde{\mathbb{G}}$),

$$\begin{aligned} \alpha : \quad \widehat{G} &\longrightarrow [H \times G]_{\sim} \\ \widehat{g} = \widehat{\Phi}(g) &\longmapsto \alpha(\widehat{g}) = [\Phi, g] \end{aligned}$$

and

$$\begin{aligned} \beta : \widehat{H} &\longrightarrow \widetilde{H} \\ \widehat{\Phi} &\longmapsto \beta(\widehat{\Phi}) = \widetilde{\Phi} = \alpha\widehat{\Phi}\alpha^{-1} \end{aligned}$$

are isomorphisms between groups and semigroups, respectively, then, for $\widetilde{g} \in [H \times G]_{\sim}$ arbitrarily fixed, we have

$$\widetilde{g} = \alpha(\widehat{g})$$

with $\widehat{g} \in \widehat{G}$, and since $\widehat{\mathbb{G}}$ is a closed extension of $\mathbb{G}$,

$$\widehat{g} = \widehat{\Phi}(g)$$

with $\Phi \in H$ and $g \in G$, which leads us to

$$\widetilde{g} = \alpha\left(\widehat{\Phi}(g)\right).$$

Now we calculate $\widetilde{\Phi}(\alpha(g))$:

$$\widetilde{\Phi}\Big(\alpha(g)\Big) = \Big(\alpha\widehat{\Phi}\alpha^{-1}\Big)\Big(\alpha(g)\Big) = \alpha\Big(\widehat{\Phi}(g)\Big) = \widetilde{g}.$$

Since $\alpha(g) = \alpha(\widehat{I}_G(g)) = [I_G, g] = \gamma(g)$, with $\gamma\colon G \longrightarrow [I_G \times G]_\sim$ the isomorphism (described in Proposition 2.8) which identifies $g \in G$ with $[I_G, g] \in [I_G \times G]_\sim \equiv G$, we conclude that

$$\widetilde{g} = \widetilde{\Phi}\Big(\gamma(g)\Big),$$

which shows us that $\widetilde{\mathbb{G}}$ is a closed extension of $\mathbb{G}$.

We are, therefore, regarding the strict and closed extension problem, in the following stage:

**(a)** either we can prove, without using the hypothesis of existence of a solution $\widehat{\mathbb{G}}$, that the strict extension $\widetilde{\mathbb{G}} = ([H \times G]_\sim, \widetilde{H})$ of $\mathbb{G}$ is closed and, in this case, we will have $\widetilde{\mathbb{G}}$ as a solution, the only one unless isomorphism (see Definition 1.4(c)) since, thus, our previous considerations show that, if $\widehat{\mathbb{G}} = (\widehat{G}, \widehat{H})$ is another solution, then, the isomorphisms $\alpha\colon \widehat{G} \longrightarrow [H \times G]_\sim$ and $\beta\colon \widehat{H} \longrightarrow \widetilde{H}$ are such that

$$\alpha\Big(\widehat{\Phi}(\widehat{g})\Big) = \beta\Big(\widehat{\Phi}\Big)\Big(\alpha(\widehat{g})\Big),$$

that is, $\widehat{\mathbb{G}}$ is isomorphic to $\mathbb{G}$ (see Definition 1.4(c)), or, otherwise;

**(b)** the problem does not have any solution.

Fortunately, we can provide a proof for the proposition below.

### 2.14 Proposition

*The extension $\widetilde{\mathbb{G}} = ([H \times G]_\sim, \widetilde{H})$, described in Proposition 2.10, of the abelian, surjective, and with identity $S$-group $\mathbb{G} = (G, H)$, is a closed extension of $\mathbb{G}$.*

*Proof.* In fact, if $\widetilde{g} = [\Phi, g] \in [H \times G]_\sim$ is arbitrarily chosen in $[H \times G]_\sim$, we have, by the definition of $\widetilde{\Phi}$ given in Proposition 2.10 and by the identification of $G$ with $[I_G \times G]_\sim \subseteq [H \times G]_\sim$ realized by the isomorphism $\gamma$ defined in Proposition 2.8, that

$$\widetilde{g} = [\Phi, g] = \widetilde{\Phi}\Big([I_G, g]\Big) = \widetilde{\Phi}\Big(\gamma(g)\Big).$$ ∎

# Theorem of Extension of $S$-Groups

### 2.15 Remark

The established in Propositions 2.3, 2.6, 2.8, 2.10, 2.12 and 2.14 shows us that the strict and closed extension problem, in its abstract version (see 1.15), has, unless isomorphism, a unique solution. Furthermore, the referred propositions provide a specific solution to the problem. The Theorem of Extension of $S$-Groups, stated below, gathers the results provided by the above referred propositions.

## 2.16 Theorem of Extension of $S$-Groups

*Let $\mathbb{G} = (G, H)$ be an abelian, surjective, and with identity $S$-group. Let also $\widetilde{\mathbb{G}} = ([H \times G]_{\sim}, \widetilde{H})$, where $[H \times G]_{\sim}$ is the abelian group described in Proposition 2.6 and $\widetilde{H}$ is the semigroup defined in Proposition 2.10. Then, $\widetilde{\mathbb{G}}$ is a strict and closed extension of $\mathbb{G}$. If $\widehat{\mathbb{G}} = (\widehat{G}, \widehat{H})$ is another strict and closed extension of $\mathbb{G}$, then the functions*

$$\begin{aligned} \alpha : \quad \widehat{G} &\longrightarrow [H \times G]_{\sim} \\ \widehat{g} = \widehat{\Phi}(g) &\longmapsto \alpha(\widehat{g}) = [\Phi, g] \end{aligned}$$

*and*

$$\begin{aligned} \beta : \widehat{H} &\longrightarrow \widetilde{H} \\ \widehat{\Phi} &\longmapsto \beta(\widehat{\Phi}) = \widetilde{\Phi} = \alpha\widehat{\Phi}\alpha^{-1} \end{aligned}$$

*are isomorphisms (between groups and semigroups, respectively), such that*

$$\alpha\Big(\widehat{\Phi}\,(\widehat{g})\Big) = \beta(\widehat{\Phi})\Big(\alpha(\widehat{g})\Big)$$

*for every $\widehat{\Phi} \in \widehat{H}$ and every $\widehat{g} \in \widehat{G}$, that is, $\widehat{\mathbb{G}}$ is isomorphic to $\widetilde{\mathbb{G}}$.*

*Furthermore, the isomorphism $\alpha$ keeps fixed the elements of the group $G$ in the following sense: $\alpha(g) = \gamma(g) \equiv g$, where*

$$\gamma \colon G \longrightarrow [I_G \times G]_{\sim}$$

*is the isomorphism defined in Proposition 2.8.* ■

# $\widetilde{\mathbb{G}}$-Distributions and its Axiomatics

## 2.17 Preliminaries

It is convenient, at this point, to remember our basic motivation described in the items 1.1 and 1.3, which led us to formulate, first, the extension problem (item 1.8), next, the strict extension problem (item 1.11) and, finally, the strict and closed extension problem for the particular $S$-group

$$\mathbb{C}(\Omega) = \Big(C(\Omega), \partial(\Omega)\Big)$$

of the continuous functions and, then, for an abstract $S$-group,

$$\mathbb{G} = (G, H)\,,$$

that, as $\mathbb{C}(\Omega)$, is abelian, surjective, and with identity. In the context of that motivation, a strict and closed extension of $\mathbb{C}(\Omega)$,

$$\widetilde{\mathbb{C}}(\Omega) = \left(\widetilde{C}(\Omega), \widetilde{\partial}(\Omega)\right),$$

constitutes a "habitat of the distributions"; a mathematical structure that encloses our desired extensions (generalizations) of the notions of continuous function and differentiability, attending the requirements described in 1.3, namely,

- every continuous function is a distribution;
- distributions are differentiable, and its derivatives are also distributions, therefore, distributions are infinitely differentiable;
- the notion of derivative of a distribution, when applied to classically differentiable functions coincide with the usual derivative and;
- the usual rules of differential calculus remain valid.

The Theorem of Extension of $S$-Groups, particularized to the case where $\mathbb{G}$ is taken as $\mathbb{C}(\Omega)$, the $S$-group of the continuous functions, provides such a mathematical structure — the wanted "habitat" — besides informing us about its uniqueness, unless isomorphism. Hence, nothing more natural now than formalizing the informally introduced notions of distributions and its derivatives, through a definition where the distributions are taken as being the elements of the group $\widetilde{C}(\Omega)$ and, its derivatives, the members of the semigroup $\widetilde{\partial}(\Omega)$ of the strict and closed extension (described in the referred theorem),

$$\widetilde{\mathbb{C}}(\Omega) = \left(\widetilde{C}(\Omega), \widetilde{\partial}(\Omega)\right),$$

of the $S$-group $\mathbb{C}(\Omega) = (C(\Omega), \partial(\Omega))$.

Since the established in the Theorem of Extension of $S$-Groups applies to any $S$-group $\mathbb{G}$ that is abelian, surjective, and with identity, and not only to the $S$-group $\mathbb{C}(\Omega)$ of the continuous functions, we can define an abstract notion of distribution and its derivatives, lets say, the $\widetilde{\mathbb{G}}$-distributions and respective derivatives, as being the elements of the group $\widetilde{G}$ and the semigroup $\widetilde{H}$, respectively, of the unique (unless isomorphism) strict and closed extension $\widetilde{\mathbb{G}} = (\widetilde{G}, \widetilde{H})$, of $\mathbb{G}$; hence, the distributions as generalizations of the concept of continuous function, and its derivatives as an extension of the classical notion of differentiability, attending the conditions above remembered, are formalized, mathematically defined, as the $\widetilde{\mathbb{C}}(\Omega)$-distributions and its derivatives.

The considerations above motivate the definition formulated below.

## 2.18 Definition

**(a)** We say that the ordered pair $(\mathbb{G}, \widehat{\mathbb{G}})$ is a **domain of distribution** if and only if $\mathbb{G} = (G, H)$ is an abelian, surjective, and with identity $S$-group, and $\widehat{\mathbb{G}} = (\widehat{G}, \widehat{H})$ is an extension of $\mathbb{G}$ isomorphic to the strict and closed extension, $\widetilde{\mathbb{G}} = (\widetilde{G}, \widetilde{H})$, of $\mathbb{G}$, defined in Theorem 2.16;

**(b)** Let $(\mathbb{G} = (G, H), \widehat{\mathbb{G}} = (\widehat{G}, \widehat{H}))$ be a domain of distribution. We say that $\widehat{g}$ is a $\widetilde{\mathbb{G}}$**-distribution** (or $\widehat{\mathbb{G}}$-distribution) if and only if $\widehat{g} \in \widehat{G}$; $\widehat{\Phi}$ is a **derivative** of the $\widetilde{\mathbb{G}}$-distributions if and only if $\widehat{\Phi} \in \widehat{H}$.

## 2.19 Remark

The Theorem of Extension of $S$-Groups, as we saw, establishes the existence and uniqueness (unless isomorphism) of a strict and closed extension, $\widetilde{\mathbb{G}} = (\widetilde{G}, \widetilde{H})$, for a given abelian, surjective, and with identity $S$-group, $\mathbb{G} = (G, H)$. In other terms, the referred theorem says, or suggests, that, roughly speaking, if we postulate the existence of two new objects (new relative to those that integrate the $S$-group $\mathbb{G} = (G, H)$) which compose two new classes, lets say, $\widehat{G}$ and $\widehat{H}$, satisfying the conditions implicit in the definition of strict and closed extension of a $S$-group, then, these classes, $\widehat{G}$ and $\widehat{H}$, in fact exist and, in a certain sense, are unique.

Hence, it is reasonable to admit that the $\widetilde{\mathbb{G}}$-distributions and its derivatives, as introduced in Definition 2.18, can be defined, equivalently, through a categoric axiomatic, that is, more explicitly, that the referred concepts, $\widetilde{\mathbb{G}}$-distributions and its derivatives, can be taken as primitive terms (concepts) defined, implicitly, by the properties that characterize the condition of being a strict and closed extension of $\mathbb{G} = (G, H)$, taken as axioms; the Theorem of Extension of $S$-Groups would then ensure that such axiomatic is consistent, that is, has a model and, yet, that any two models are isomorphic, that is, it is a categoric axiomatic.

As we know, in the formulation of an axiomatic theory we must, first of all, state in a complete form the initial assumptions, that is:

- the concepts (terms) that will be taken as primitives;
- the properties that, with reference to the primitive concepts, we want to admit, that is, the axioms and;
- the precedent theories, i.e., those already axiomatized and that will be used in the axiomatical structuring of the topic at hand.

Next, in 2.20, following the above prescriptions and guided by the previous considerations, we formulate an axiomatic that, as will be proved in 2.22 after obtained, in 2.21, some elementary consequences of the axioms, admit as a model the objects that compose a strict and closed extension of the $S$-group $\mathbb{G}$, thus being a categoric axiomatic.

## 2.20 $\widetilde{\mathbb{G}}$-distributions Axiomatic

The primitive terms and the precedent theories of the $\widetilde{\mathbb{G}}$-distributions axiomatic are stated below.

- Primitive terms: $\widetilde{\mathbb{G}}$-distribution, derivative (of $\widetilde{\mathbb{G}}$-distribution) and addition (of $\widetilde{\mathbb{G}}$-distribution);
- Precedent theories: Classical logic, Set theory and $S$-groups theory.

It is important to observe that, with the above choice of primitive terms, we started to attribute two meanings to each of the terms $\widetilde{\mathbb{G}}$-distribution and $\widetilde{\mathbb{G}}$-distribution derivative, namely: those instituted by Definition 2.18(b) and those that, as primitive terms, the axioms of the theory listed ahead, implicitly, attribute to them. In principle, these two meanings may not coincide; however, our choice of axioms, as suggest our considerations in Remark 2.19, will be such that they coincide, as we will see.

The use of primitive symbols, that is, symbols associated with primitive concepts, allows a more economic formulation of the axioms. With this intent we introduce the following primitive symbols:

$\widehat{G}$, for the class of the $\widetilde{\mathbb{G}}$-distributions;

$\widehat{H}$, represents the set of the $\widetilde{\mathbb{G}}$-distribution derivatives and;

$+$, denotes the addition of $\widetilde{\mathbb{G}}$-distributions.

The axioms, stated in relation to an arbitrarily fixed $S$-group, $\mathbb{G} = (G, H)$, abelian, surjective, and with identity, are the following:

**Axiom 1** $G \subseteq \widehat{G}$.

**Axiom 2** The addition, $+$, is a binary operation in $\widehat{G}$ whose restriction to $G \times G$ is the addition of the abelian group $G$.

**Axiom 3** The derivatives (of $\widetilde{\mathbb{G}}$-distribution), that is, the members $\widehat{\Phi} \in \widehat{H}$, are functions with domain and codomain equal to $\widehat{G}$, such that:

**(a)** $\widehat{\Phi}(\widehat{g}+\widehat{h}) = \widehat{\Phi}(\widehat{g}) + \widehat{\Phi}(\widehat{h})$ whatever $\widehat{g}, \widehat{h} \in \widehat{G}$;

**(b)** each $\widehat{\Phi} \in \widehat{H}$ is an extension to $\widehat{G}$ of an unique homomorphism $\Phi \in H$, that is, $\widehat{\Phi}(g) = \Phi(g)$ for every $g \in G_\Phi$ (domain of $\Phi$), for an unique $\Phi \in H$;

**(c)** the class $\widehat{H}$ with the usual composition of functions, is a semigroup isomorphic to the semigroup $H$ and the function $\widehat{\ }$ defined by

$$\begin{aligned} \widehat{\ } : H &\longrightarrow \widehat{H} \\ \Phi &\longmapsto \widehat{\Phi} := \text{the extension of } \Phi \text{ to } \widehat{G}, \end{aligned}$$

is an isomorphism (between semigroups).

**Axiom 4** For each $\Phi \in H$, if $g \in G$ and $g \notin G_\Phi$ (the domain of $\Phi$), then, $\widehat{\Phi}(g) \notin G$.

**Axiom 5** For each $\widehat{g} \in \widehat{G}$, there exist $\Phi \in H$ and $g \in G$ such that $\widehat{g} = \widehat{\Phi}(g)$.

**Remark.** The symbol $\widehat{\Phi}$ denotes in the axioms above, as it will denote in all our discussion about these axioms, the (unique) extension of $\Phi$ to $\widehat{G}$ that belongs to $\widehat{H}$, in other terms, $\widehat{\Phi}$ is the image of $\Phi \in H$ by the isomorphism $\widehat{\ }$. It is convenient to remark also that the function $\widehat{\ }$ introduced in Axiom 3(c) has its "definition" authorized by Axiom 3(b). Nevertheless, this function, $\widehat{\ }$, can not be taken as a defined concept of our axiomatic, since it contains information about the primitive concept of "derivative" that does not appear in the other axioms: in fact, with the function $\widehat{\ }: H \longrightarrow \widehat{H}$ as postulated, to **every** $\Phi \in H$ corresponds one derivative $\widehat{\Phi} \in \widehat{H}$, an information that, in the absence of Axiom 3(c), that is, without the function $\widehat{\ }$, can not be found explicit nor implicitly in the other axioms.

## 2.21 Elementary Consequences of the Axioms

Among the logical consequences of the axioms formulated in 2.20, we highlight the following:

**(a)** $\widehat{g}+\widehat{h} = \widehat{h}+\widehat{g}$ for every $\widehat{g}, \widehat{h} \in \widehat{G}$.

In fact, if $\widehat{g}, \widehat{h} \in \widehat{G}$, then, by Axiom 5, there exist $\Phi, \Psi \in H$ and $g, h \in G$ such that

$$\widehat{g} = \widehat{\Phi}(g) \quad \text{and} \quad \widehat{h} = \widehat{\Psi}(h).$$

Since $\Phi$ and $\Psi$ are surjective homomorphisms (for the $S$-group $\mathbb{G} = (G, H)$ is surjective), there exist $g^* \in G_\Psi$ (domain of $\Psi$) and $h^* \in G_\Phi$ (domain of $\Phi$) such that

$$g = \Psi(g^*) \quad \text{and} \quad h = \Phi(h^*),$$

and hence

$$\widehat{g} = \widehat{\Phi}\Big(\Psi(g^*)\Big) \quad \text{and} \quad \widehat{h} = \widehat{\Psi}\Big(\Phi(h^*)\Big).$$

However, by Axiom 3(b), we have

$$\Psi(g^*) = \widehat{\Psi}(g^*) \quad \text{and} \quad \Phi(h^*) = \widehat{\Phi}(h^*),$$

leading us to

$$\widehat{g} = \widehat{\Phi}\left(\widehat{\Psi}(g^*)\right) \quad \text{and} \quad \widehat{h} = \widehat{\Psi}\left(\widehat{\Phi}(h^*)\right).$$

Now, we know from Axiom 3(c) that $\widehat{\ }: H \longrightarrow \widehat{H}$ is an isomorphism and, thus,

$$\widehat{\Phi\Psi} = \widehat{\Phi}\widehat{\Psi},$$

and since $\Phi\Psi = \Psi\Phi$, for $H$ is an abelian semigroup (remember that the $S$-group $\mathbb{G} = (G, H)$ is abelian), we get

$$\widehat{\Phi}\widehat{\Psi} = \widehat{\Phi\Psi} = \widehat{\Psi\Phi} = \widehat{\Psi}\widehat{\Phi}.$$

Therefore,

$$\widehat{g} + \widehat{h} = \widehat{\Phi}\left(\widehat{\Psi}(g^*)\right) + \widehat{\Psi}\left(\widehat{\Phi}(h^*)\right) = \left(\widehat{\Phi}\widehat{\Psi}\right)(g^*) + \left(\widehat{\Phi}\widehat{\Psi}\right)(h^*)$$

that, taking into account Axiom 3(a), allows us to conclude that

$$\widehat{g} + \widehat{h} = \widehat{\Phi}\widehat{\Psi}(g^* + h^*).$$

But $g^*, h^* \in G$ and $G$ is an abelian group and, since, by Axiom 2, the addition restricted to $G$ is the own group addition of $G$, from the last expression we obtain

$$\widehat{g} + \widehat{h} = \widehat{\Phi}\widehat{\Psi}(h^* + g^*).$$

Resorting once more to Axiom 3(a) and to the commutativity of the derivatives $\widehat{\Phi}$ and $\widehat{\Psi}$, results from the expression above that

$$\widehat{g} + \widehat{h} = \widehat{\Phi}\widehat{\Psi}(h^*) + \widehat{\Phi}\widehat{\Psi}(h^*) = \widehat{\Psi}\left(\widehat{\Phi}(h^*)\right) + \widehat{\Phi}\left(\widehat{\Psi}(g^*)\right).$$

However, as we saw above,

$$\widehat{\Phi}(h^*) = \Phi(h^*) = h$$

and

$$\widehat{\Psi}(g^*) = \Psi(g^*) = g,$$

therefore,

$$\widehat{g} + \widehat{h} = \widehat{\Psi}(h) + \widehat{\Phi}(g) = \widehat{h} + \widehat{g}.$$

**(b)** $(\widehat{f} + \widehat{g}) + \widehat{h} = \widehat{f} + (\widehat{g} + \widehat{h})$ for every $\widehat{f}, \widehat{g}, \widehat{h} \in \widehat{G}$.

With a procedure analogous to that of the proof of the above item (a), we can also prove (b).

**(c)** There exists $\widehat{0} \in \widehat{G}$ such that $\widehat{g} + \widehat{0} = \widehat{g}$ for every $\widehat{g} \in \widehat{G}$.

Let $\widehat{g} \in \widehat{G}$ arbitrarily chosen. According to Axiom 5,

$$\widehat{g} = \widehat{\Phi}(g)$$

for some $\Phi \in H$ and $g \in G$. But $\Phi \colon G_\Phi \longrightarrow G$ is a homomorphisms and hence $\Phi(0) = 0$, where $0 \in G_\Phi \subseteq G$ is the neutral element of the group $G$. Now, from Axiom 3(b),

$$\widehat{\Phi}(0) = \Phi(0) = 0.$$

Hence, and taking into account Axiom 3(a),

$$\widehat{g} + 0 = \widehat{\Phi}(g) + \widehat{\Phi}(0) = \widehat{\Phi}(g + 0).$$

But, by Axiom 2, the addition of $\widetilde{\mathbb{G}}$-distributions, when restricted to elements of the group $G$, coincides with the addition of the group $G$ for which $g + 0 = g$. Hence, we have that

$$\widehat{g} + 0 = \widehat{\Phi}(g).$$

But, as above, $\widehat{g} = \widehat{\Phi}(g)$, and hence, defining $\widehat{0} \coloneqq 0$, we get:

$$\widehat{g} + \widehat{0} = \widehat{g}.$$

**(d)** For each $\widehat{g} \in \widehat{G}$, there exist $-\widehat{g} \in \widehat{G}$ such that $\widehat{g} + (-\widehat{g}) = \widehat{0}$.

In fact, let $\widehat{g} \in \widehat{G}$ be arbitrarily chosen and let also $\Phi \in H$ and $g \in G$ be such that

$$\widehat{g} = \widehat{\Phi}(g),$$

with the existence of $\Phi$ and $g$ ensured by Axiom 5. Now, we define:

$$-\widehat{g} \coloneqq \widehat{\Phi}(-g),$$

where $-g \in G$ is the opposite element (in the group $G$) of $g$. Hence, we have:

$$\widehat{g} + (-\widehat{g}) = \widehat{\Phi}(g) + \widehat{\Phi}(-g),$$

or yet, resorting to Axiom 3(a),

$$\widehat{g} + (-\widehat{g}) = \widehat{\Phi}\Big(g + (-g)\Big).$$

However, by Axiom 2, the addition in $\widehat{G}$ coincides with the addition of the group $G$ when applied to elements of the latter. Hence, we have that $g + (-g) = 0$ and, from this, we obtain:

$$\widehat{g} + (-\widehat{g}) = \widehat{\Phi}(0).$$

But, according to Axiom 3(b), $\widehat{\Phi}(0) = \Phi(0) = 0$, and hence

$$\widehat{g} + (-\widehat{g}) = 0 = \widehat{0}.$$

**(e)** $\widehat{G}$ with the addition of $\widetilde{\mathbb{G}}$-distribution is an abelian group that has the group $G$ as a subgroup.

In fact, the items (a) to (d) above show us that $\widehat{G}$ equipped with the addition (of $\widetilde{\mathbb{G}}$-distribution) is an abelian group, while from Axiom 1 and Axiom 2 we have that the group $G \subseteq \widehat{G}$ is a subgroup of the group $\widehat{G}$.

**(f)** The derivatives $\widehat{\Phi} \in \widehat{H}$ of the $\widetilde{\mathbb{G}}$-distributions are endomorphisms on the group $\widehat{G}$; $\widehat{H}$, with the usual composition of functions, is a semigroup (of endomorphisms on $\widehat{G}$) and a prolongation of the semigroup $H$ to $\widehat{G}$ (see Definition 1.6(c)).

In fact, from (e) above together with Axiom 3(a) one concludes that the derivatives $\widehat{\Phi} \in \widehat{H}$ are endomorphisms on the group $\widehat{G}$, while by Axiom 3(b) and Axiom 3(c) one obtain that $\widehat{H}$ is a prolongation of the semigroup $H$ to $\widehat{G}$ (according to Definition 1.6(c)).

**(g)** $\widehat{\mathbb{G}} := (\widehat{G}, \widehat{H})$ is an $S$-group that is an extension of $S$-group $\mathbb{G} = (G, H)$.

This comes trivially from (e) and (f) above, in accordance with the definitions of $S$-group (Definition 1.4(a)) and $S$-group extension (Definition 1.6(d)).

**(h)** $\widehat{\mathbb{G}} = (\widehat{G}, \widehat{H})$ is a strict and closed extension of the $S$-group $\mathbb{G} = (G, H)$.

Immediate consequence of (g) above, Axiom 4 and Axiom 5, taking into account the pertinent definitions (1.10(b) and 1.14).

## 2.22 Consistency and Categoricity

A simple inspection of the axioms stated in 2.20 is enough to convince us that: if $\widehat{\mathbb{G}} = (\widehat{G}, \widehat{H})$ is a strict and closed extension of the $S$-group $\mathbb{G} = (G, H)$ associated to the axiomatic in appreciation (that is, the $S$-group $\mathbb{G}$ the axioms make reference to), then, the members of the classes $\widehat{G}$ and $\widehat{H}$ constitute a model of this axiomatic.

Conversely, the logical consequences of these same axioms, obtained in 2.21, show us that: if the axiomatic of the $\widetilde{\mathbb{G}}$-distributions does have a model, then, the sets $\widehat{G}$ and $\widehat{H}$ of the $\widetilde{\mathbb{G}}$-distributions and its derivatives, respectively, correspondents to a model, are such that the ordered pair $(\widehat{G}, \widehat{H})$ is a strict and closed extension of the abelian, surjective, and with identity $S$-group $\mathbb{G} = (G, H)$ associated to the axiomatic.

Since, by the Theorem of Extension of $S$-Groups (Theorem 2.16), the $S$-group $\widetilde{\mathbb{G}} = (\widetilde{G}, \widetilde{H})$ there defined is a strict and closed extension of $\mathbb{G}$, and, by the same theorem, any other extension of $\mathbb{G}$ is isomorphic to $\widetilde{\mathbb{G}}$, we conclude that the $\widetilde{\mathbb{G}}$-distributions axiomatic does have a model, that is, it is consistent and, further, it is categoric, that is, it has (essentially) an unique model.

From the conclusion above, it results that the two meanings attributed to the terms "$\widetilde{\mathbb{G}}$-distribution" and "$\widetilde{\mathbb{G}}$-distribution derivative" mentioned in 2.20 (the ones from Definition 2.18(b) and the ones implicitly determined by the axioms stated in 2.20), are equivalent.

## 2.23 A Special Set of Properties

Let $\mathbb{G} = (G, H)$ be an abelian, surjective, and with identity $S$-group, and $\widehat{\mathbb{G}} = (\widehat{G}, \widehat{H})$ be a strict and closed extension of $\mathbb{G}$. In this case, $\widehat{\mathbb{G}}$ is a model of the $\widetilde{\mathbb{G}}$-distributions axiomatic and, therefore, taking into account Axiom 1, Axiom 3 and Axiom 5 of the referred axiomatic, we have that:

**(a)** $G \subseteq \widehat{G}$;

**(b)** the function

$$\begin{aligned} \widehat{\ } : H &\longrightarrow \widehat{H} \\ \Phi &\longmapsto \widehat{\Phi}, \end{aligned}$$

where $\widehat{\Phi}(g) = \Phi(g)$ for every $g \in G_\Phi$, is an isomorphism;

**(c)** for each $\widehat{g} \in \widehat{G}$, there exists $\Phi \in H$ and $g \in G$ such that

$$\widehat{g} = \widehat{\Phi}(g).$$

Furthermore, by Proposition 1.16(b), we also have that:

**(d)** if $g, h \in G$ and $\Phi \in H$, then,

$$\widehat{\Phi}(g) = \widehat{\Phi}(h) \quad \text{if and only if} \quad g - h \in N(\Phi).$$

The highlight given to the properties (a) to (d) above, that is, the special interest that lead us to make them explicit, is justified by what establishes Proposition 2.24 (stated below) regarding a set of conditions mirroring the properties (a) to (d); as we will see, the cited proposition allows the formulation of a axiomatic equivalent to the $\widetilde{\mathbb{G}}$-distributions one, with a smaller number of primitive terms and axioms. In this new axiomatic, stated in 2.25, with only two primitive terms, "$\widetilde{\mathbb{G}}$-distribution" and "derivative (of $\widetilde{\mathbb{G}}$-distribution)", the term "addition" assumes the status of defined term.

## 2.24 Proposition

*Let $\mathbb{G} = (G, H)$ be an abelian, surjective, and with identity $S$-group. Suppose that there exists a set $\widehat{G}$ and, for each $\Phi\colon G_\Phi \longrightarrow G$ in $H$ there exists a corresponding function $\widehat{\Phi}\colon \widehat{G} \longrightarrow \widehat{G}$ (with the set $\widehat{G}$ as domain and codomain), such that:*

**(a)** *$G \subseteq \widehat{G}$;*

**(b)** *$\widehat{\Phi}(g) = \Phi(g)$ for every $g \in G_\Phi$ and the set $\widehat{H} := \{\widehat{\Phi}\colon \Phi \in H\}$, with the usual composition of functions, is a semigroup isomorphic to the semigroup $H$ with the function*

$$\begin{aligned} \widehat{\ }\colon H &\longrightarrow \widehat{H} \\ \Phi &\longmapsto \widehat{\Phi} \end{aligned}$$

*as an isomorphism;*

**(c)** *for each $\widehat{g} \in \widehat{G}$, there exist $\Phi \in H$ and $g \in G$ such that*

$$\widehat{g} = \widehat{\Phi}(g);$$

**(d)** *if $g, h \in G$ and $\Phi \in H$, then,*

$$\widehat{\Phi}(g) = \widehat{\Phi}(h) \quad \textit{if and only if} \quad g - h \in N(\Phi).$$

*Satisfied the above conditions, there exists an unique manner to extend the addition of the group $G$ to the set $\widehat{G}$, so that $\widehat{G}$ becomes an abelian group and $\widehat{\Phi}\colon \widehat{G} \longrightarrow \widehat{G}$ an endomorphism on $\widehat{G}$.*

*Proof.* The proof provided ahead is divided into six lemmas, at which the strict and closed extension, $\widetilde{\mathbb{G}} = (\widetilde{G}, \widetilde{H})$, of the $S$-group $\mathbb{G} = (G, H)$, defined by the Theorem of Extension of $S$-Groups (Theorem 2.16), plays a major role. Hence, it is appropriate to remember that

$$\widetilde{G} = [H \times G]_\sim$$

is the abelian group described in Proposition 2.6, whose members, $[\Phi, g]$ with $\Phi \in H$ and $g \in G$, are the equivalence classes of the set $H \times G$ determined by the equivalence relation $\sim$ defined in Proposition 2.3.

**Lemma 1.** *For any $\Phi, \Psi \in H$ and $g, h \in G$ we have:*

$$\widehat{\Phi}(g) = \widehat{\Psi}(h) \quad \textit{if and only if} \quad \Phi^{-1}(h) - \Psi^{-1}(g) \subseteq N(\Phi\Psi).$$

*Proof.* Let $\Phi, \Psi \in H$ and $g, h \in G$ be arbitrarily fixed. Since the $S$-group $\mathbb{G} = (G, H)$ is surjective, $\Phi\colon G_\Phi \longrightarrow G$ and $\Psi\colon G_\Psi \longrightarrow G$ are surjective homomorphisms and hence $\Phi^{-1}(h) \subseteq G_\Phi$ and $\Psi^{-1}(g) \subseteq G_\Psi$ are non-empty sets. Let then $g^*$ and $h^*$ be arbitrarily chosen elements in $\Psi^{-1}(g)$ and $\Phi^{-1}(h)$, respectively, that is,

$$g^* \in \Psi^{-1}(g) \subseteq G_\Psi \quad \text{and} \quad h^* \in \Phi^{-1}(h) \subseteq G_\Phi.$$

Thus, with the hypothesis (b) in mind,

$$g = \Psi(g^*) = \widehat{\Psi}(g^*) \quad \text{and} \quad h = \Phi(h^*) = \widehat{\Phi}(h^*)$$

and, therefore,

$$\widehat{\Phi}(g) = \widehat{\Phi}\left(\widehat{\Psi}(g^*)\right) \quad \text{and} \quad \widehat{\Psi}(h) = \widehat{\Psi}\left(\widehat{\Phi}(h^*)\right).$$

By hypothesis (b), $\widehat{\ }: H \longrightarrow \widehat{H}$ is an isomorphism and hence

$$\widehat{\Phi\Psi} = \widehat{\Phi}\widehat{\Psi} \quad \text{and} \quad \widehat{\Psi\Phi} = \widehat{\Psi}\widehat{\Phi}.$$

Now $\Phi\Psi = \Psi\Phi$, since $H$ is an abelian semigroup (for $\mathbb{G} = (G, H)$ is abelian), and hence

$$\widehat{\Phi}\widehat{\Psi} = \widehat{\Psi}\widehat{\Phi}.$$

Therefore,

$$\widehat{\Phi}(g) = \widehat{\Phi\Psi}(g^*) \quad \text{and} \quad \widehat{\Psi}(h) = \widehat{\Phi\Psi}(h^*)$$

and, with that, and by what establishes the hypothesis (d), we obtain that:

$$\widehat{\Phi\Psi}(g^*) = \widehat{\Phi\Psi}(h^*) \quad \text{if and only if} \quad g^* - h^* \in N(\Phi\Psi)$$

that is,

$$\widehat{\Phi}(g) = \widehat{\Psi}(h) \quad \text{if and only if} \quad g^* - h^* \in N(\Phi\Psi).$$

Since $g^*$ and $h^*$ are arbitrary elements of $\Psi^{-1}(g)$ and $\Phi^{-1}(h)$, respectively, the above equivalence assumes the form:

$$\widehat{\Phi}(g) = \widehat{\Psi}(h) \quad \text{if and only if} \quad \Psi^{-1}(g) - \Phi^{-1}(h) \subseteq N(\Phi\Psi),$$

or yet, taking into account that $N(\Phi\Psi) \subseteq G$ is a subgroup of $G$,

$$\widehat{\Phi}(g) = \widehat{\Psi}(h) \quad \text{if and only if} \quad \Phi^{-1}(h) - \Psi^{-1}(g) \in N(\Phi\Psi).$$

**Lemma 2.** *The function $\alpha$ ahead is well-defined and is a bijection.*

$$\begin{aligned} \alpha : \quad & \widehat{G} \longrightarrow \widetilde{G} = [H \times G]_{\sim} \\ & \widehat{g} = \widehat{\Phi}(g) \longmapsto \alpha(\widehat{g}) := [\Phi, g]. \end{aligned}$$

*(Remark — By the hypothesis (c), each $\widehat{g} \in \widehat{G}$ can be written in the form $\widehat{g} = \widehat{\Phi}(g)$ for some $\Phi \in H$ and $g \in G$).*

*Proof.* From Proposition 2.3 we have, for $(\Phi, g)$ and $(\Psi, h)$ in $H \times G$, that

$$(\Phi, g) \sim (\Psi, h) \quad \text{if and only if} \quad \Psi^{-1}(g) - \Phi^{-1}(h) \subseteq N(\Phi\Psi).$$

Hence, taking into account Lemma 1,

$$\widehat{\Phi}(g) = \widehat{\Psi}(h) \quad \text{if and only if} \quad (\Phi, g) \sim (\Psi, h).$$

But, $(\Phi, g) \sim (\Psi, h)$ if and only if the equivalence classes determined by $(\Phi, g)$ and $(\Psi, h)$, namely, $[\Phi, g]$ and $[\Psi, h]$, respectively, are equal. Hence,

$$\widehat{\Phi}(g) = \widehat{\Psi}(h) \quad \text{if and only if} \quad [\Phi, g] = [\Psi, h].$$

This equivalence, among with the hypothesis (c) which allows to write each $\widehat{g} \in \widehat{G}$ in the form $\widehat{g} = \widehat{\Phi}(g)$, with $\Phi \in H$ and $g \in G$, show us not only that $\alpha$ is well-defined but also that it is a bijection.

**Lemma 3.** *Let $\widehat{g}, \widehat{h} \in \widehat{G}$ be arbitrarily chosen. Let also $\Phi, \Psi \in H$ and $g, h \in G$ be such that (according to hypothesis (c))*

$$\widehat{g} = \widehat{\Phi}(g) \quad \textit{and} \quad \widehat{h} = \widehat{\Psi}(h).$$

*Defining $\widehat{g} + \widehat{h}$ (addition in $\widehat{G}$) by*

$$\widehat{g} + \widehat{h} := \widehat{\Phi\Psi}(g^* + h^*),$$

*where $g^*$ and $h^*$ are arbitrarily chosen elements in $\Psi^{-1}(g)$ and $\Phi^{-1}(h)$, respectively, that is,*

$$g^* \in \Psi^{-1}(g) \quad \textit{and} \quad h^* \in \Phi^{-1}(h),$$

*then:*

- *$\widehat{g} + \widehat{h}$ is well-defined for each $\widehat{g}$ and $\widehat{h}$ of $\widehat{G}$;*
- *$\widehat{G}$ with the addition above defined is an abelian group isomorphic to the group $\widetilde{G} = [H \times G]_\sim$;*
- *the function $\alpha\colon \widehat{G} \longrightarrow \widetilde{G}$ defined in Lemma 2 is an isomorphism from the group $\widehat{G}$ onto the group $\widetilde{G} = [H \times G]_\sim$.*

*Proof.* From Proposition 2.6 we have, for $[\Phi, g]$ and $[\Psi, h]$ in $[H \times G]_\sim$, that

$$[\Phi, g] + [\Psi, h] := [\Phi\Psi, g^* + h^*]$$

where $g^* \in \Psi^{-1}(g)$ and $h^* \in \Phi^{-1}(h)$, is a well-defined binary operation in $[H \times G]_\sim$, and this set with this operation is an abelian group. Hence, and with the bijection $\alpha\colon \widehat{G} \longrightarrow [H \times G]_\sim$ defined in Lemma 2, it is also well-defined the following operation (addition) in $\widehat{G}$: for $\widehat{g}, \widehat{h} \in \widehat{G}$ arbitrarily fixed we define

$$\widehat{g} + \widehat{h} := \alpha^{-1}\left[\alpha(\widehat{g}) + \alpha(\widehat{h})\right].$$

It is easy to show that $\widehat{G}$ with this operation is an abelian group and that $\alpha$ is an isomorphism from this group onto the group $[H \times G]_\sim$. On the other hand, by hypothesis (c), $\widehat{g}$ and $\widehat{h}$ can be written in the form:

$$\widehat{g} = \widehat{\Phi}(g) \quad \text{and} \quad \widehat{h} = \widehat{\Psi}(h)$$

where $\Phi, \Psi \in H$ and $g, h \in G$. Hence, with the association rule of the function $\alpha$ in mind,

$$\alpha(\widehat{g}) = [\Phi, g] \quad \text{and} \quad \alpha(\widehat{h}) = [\Psi, h]$$

and with that, through the addition in $[H \times G]_\sim$,

$$\alpha(\widehat{g}) + \alpha(\widehat{h}) = [\Phi, g] + [\Psi, h] = [\Phi\Psi, g^* + h^*]$$

with

$$g^* \in \Psi^{-1}(g) \quad \text{and} \quad h^* \in \Phi^{-1}(h).$$

Therefore,

$$\widehat{g} + \widehat{h} = \alpha^{-1}\Big(\alpha(\widehat{g}) + \alpha(\widehat{h})\Big) = \alpha^{-1}\Big([\Phi\Psi, g^* + h^*]\Big)$$

and hence, taking into account the definition of $\alpha$ (in Lemma 2) and that $\alpha$ is bijective, we obtain:

$$\widehat{g} + \widehat{h} = \widehat{\Phi\Psi}(g^* + h^*)$$

which concludes the proof of the lemma.

**Lemma 4.** *The abelian group $\widehat{G}$, described in Lemma 3, has $G$ as a subgroup, and the function $\alpha\colon \widehat{G} \longrightarrow [H \times G]_\sim$, defined in Lemma 2, restricted to $G$, $\alpha|_G$, is an isomorphism from the subgroup $G$ of $\widehat{G}$ onto the subgroup*

$$[I_G \times G]_\sim := \Big\{[I_G, g]\colon g \in G\Big\} \subseteq [H \times G]_\sim$$

*of the group $[H \times G]_\sim$.*

*Proof.* Lets considerate the abelian group $\widehat{G}$ described in Lemma 3 and, for now, denote by "$\oplus$" its addition operation. Since $G \subseteq \widehat{G}$ (by the hypothesis (a)) and $G$ is also an abelian group (for $\mathbb{G} = (G, H)$ is an $S$-group), then, the group $G$ is a subgroup of $\widehat{G}$ if and only if the addition of $\widehat{G}$, $\oplus$, when restricted to the elements of $G \subseteq \widehat{G}$, coincides with the addition, $+$, of $G$. We must then prove that

$$g \oplus h = g + h$$

whatever $g, h \in G \subseteq \widehat{G}$. Let then $g, h \in G$ be arbitrarily fixed. As we know, the identity homomorphism on $G$, $I_G\colon G \longrightarrow G$ where $I_G(g) = g$ for every $g \in G$, is an element of the semigroup $H$, for $\mathbb{G} = (G, H)$ is a $S$-group with identity. Hence, by hypothesis (b) we have that $\widehat{I}_G \in \widehat{H}$ is such that[10]

$$\widehat{I}_G(g') = I_G(g') = g' \quad \text{for every} \quad g' \in G.$$

[10] Actually, $\widehat{I}_G(\widehat{g}) = \widehat{g}$ for every $\widehat{g} \in \widehat{G} \supseteq G$, that is, $\widehat{I}_G = I_{\widehat{G}}$. In fact, for $\Phi \in H$ arbitrarily chosen we have that $I_G\Phi = \Phi$ and hence, by hypothesis (b), $\widehat{I_G\Phi} = \widehat{I}_G\widehat{\Phi} = \widehat{\Phi}$. If now $\widehat{g} \in \widehat{G}$ is arbitrarily fixed, then, by hypothesis (c), $\widehat{g} = \widehat{\Psi}(g)$ with $\Psi \in H$ and $g \in G$ and, therefore, $\widehat{I}_G(\widehat{g}) = \widehat{I}_G(\widehat{\Psi}(g)) = \Big(\widehat{I}_G\widehat{\Psi}\Big)(g) = \widehat{\Psi}(g) = \widehat{g}$.

Therefore,

$$g = \widehat{I}_G(g) \quad \text{and} \quad h = \widehat{I}_G(h)$$

and, thus, by the definition of the function $\alpha$ in Lemma 2, comes that:

$$\alpha(g) = [I_G, g] \quad \text{and} \quad \alpha(h) = [I_G, h].$$

Now, from the definition of addition of the group $[H \times G]_\sim$ (given in Proposition 2.6),

$$\alpha(g) + \alpha(h) = [I_G, g] + [I_G, h] = [I_G, g^* + h^*]$$

where

$$g^* \in I_G^{-1}(g) \quad \text{and} \quad h^* \in I_G^{-1}(h),$$

that is,

$$g^* = g \quad \text{and} \quad h^* = h.$$

Hence,

$$\alpha(g) + \alpha(h) = [I_G, g + h].$$

On the other hand, the addition in $\widehat{G}$, $\oplus$, defined in Lemma 3, can be written, as we saw in the proof of the referred lemma, through the function $\alpha$ in the following form:

$$\widehat{g} \oplus \widehat{h} = \alpha^{-1}\left(\alpha(\widehat{g}) + \alpha(\widehat{h})\right)$$

for every $\widehat{g}, \widehat{h} \in \widehat{G}$. Hence,

$$g \oplus h = \alpha^{-1}\Big(\alpha(g) + \alpha(h)\Big) = \alpha^{-1}\Big([I_G, g + h]\Big).$$

Since $\alpha$ is a bijection and taking into account its definition, results that

$$\alpha^{-1}\Big([I_G, g + h]\Big) = \widehat{I}_G(g + h) = g + h.$$

Thus,

$$g \oplus h = g + h$$

which allows us to conclude that the group $G$ is a subgroup of the group $\widehat{G}$.

From Proposition 2.8 we have that: $[I_G \times G]_\sim$ is a subgroup of the group $[H \times G]_\sim$, isomorphic to the group $G$, and the function

$$\begin{aligned} \gamma : G &\longrightarrow [I_G \times G]_\sim \\ g &\longmapsto \gamma(g) := [I_G, g] \end{aligned}$$

is an isomorphism. Now, if $g \in G \subseteq \widehat{G}$, then, $g = \widehat{I}_G(g)$ and hence

$$\alpha(g) = [I_G, g] = \gamma(g),$$

that is, $\alpha|_G = \gamma$, which shows the restriction of $\alpha$ to $G$ to be an isomorphism from the subgroup $G \subseteq \widehat{G}$ of $\widehat{G}$ onto the subgroup $[I_G \times G]_\sim$ of the group $[H \times G]_\sim$.

**Lemma 5.** *Let $\widehat{G}$ be the abelian group described in Lemma 3. For each $\Phi \in H$, the corresponding function $\widehat{\Phi}\colon \widehat{G} \longrightarrow \widehat{G}$ in $\widehat{H}$ is an endomorphism on $\widehat{G}$.*

*Proof.* Let $\widehat{g}, \widehat{h} \in \widehat{G}$ and $\widehat{\Phi} \in \widehat{H}$ arbitrarily fixed. From hypothesis (c), there exists $g, h \in G$ and $\varphi, \Psi \in H$ such that

$$\widehat{g} = \widehat{\varphi}(g) \quad \text{and} \quad \widehat{h} = \widehat{\Psi}(h).$$

Since $\Phi\colon G_\Phi \longrightarrow G$ is surjective, we have that

$$g = \Phi(g_1) \quad \text{and} \quad h = \Phi(h_1)$$

with $g_1, h_1 \in G_\Phi$. Hence we have

$$\widehat{g} = \widehat{\varphi}(g) = \widehat{\varphi}\Big(\Phi(g_1)\Big) = \widehat{\varphi}\Big(\widehat{\Phi}(g_1)\Big) = \widehat{\varphi\Phi}(g_1)$$

and

$$\widehat{h} = \widehat{\Psi}(h) = \widehat{\Psi}\Big(\Phi(h_1)\Big) = \widehat{\Psi}\Big(\widehat{\Phi}(h_1)\Big) = \widehat{\Psi\Phi}(g_1),$$

where hypothesis (b) was used. Now, with the addition in $\widehat{G}$ given in Lemma 3 in mind,

$$\widehat{g} + \widehat{h} = \widehat{\varphi\Phi}(g_1) + \widehat{\Psi\Phi}(h_1) = \Big(\widehat{\varphi\Phi\Psi\Phi}\Big)\,(g_1^* + h_1^*)$$

where

$$g_1^* \in (\Psi\Phi)^{-1}(g_1) = \Big\{\xi^* \in G_{\Psi\Phi}\colon\ (\Psi\Phi)\,(\xi^*) = g_1\Big\}$$

and

$$h_1^* \in (\varphi\Phi)^{-1}(h_1) = \Big\{\eta^* \in G_{\varphi\Phi}\colon\ (\varphi\Phi)\,(\eta^*) = h_1\Big\}.$$

Hence, we have that

$$\begin{aligned}
\widehat{\Phi}\Big(\widehat{g} + \widehat{h}\Big) &= \widehat{\Phi}\left(\Big(\widehat{\varphi\Phi\Psi\Phi}\Big)(g_1^* + h_1^*)\right) = \\
&= \Big(\widehat{\varphi\Phi\Psi\Phi}\Big)\,\Big(\widehat{\Phi}(g_1^* + h_1^*)\Big) = \\
&= \Big(\widehat{\varphi\Phi\Psi\Phi}\Big)\,\Big(\Phi(g_1^* + h_1^*)\Big) = \\
&= \Big(\widehat{\varphi\Phi\Psi\Phi}\Big)\,\Big(\Phi(g_1^*) + \Phi(h_1^*)\Big),
\end{aligned}$$

where, again, we used the hypothesis (b) and the fact that $\Phi \in H$ is a homomorphism. Since $g_1^* \in (\Psi\Phi)^{-1}(g_1)$ and $h_1^* \in (\varphi\Phi)^{-1}(h_1)$, then

$$(\Psi\Phi)\,(g_1^*) = g_1 \quad \text{and} \quad (\varphi\Phi)\,(h_1^*) = h_1.$$

Thus,

$$\Phi\Big((\Psi\Phi)(g_1^*)\Big) = \Phi(g_1) = g$$

and

$$\Phi\Big((\varphi\Phi)(h_1^*)\Big) = \Phi(h_1) = h,$$

from where we conclude that

$$\Phi(g_1^*) \in (\Phi\Psi)^{-1}(g) \quad \text{and} \quad \Phi(h_1^*) \in (\Phi\varphi)^{-1}(h).$$

Now, these pertinence relations and the definition of addition of the group $\widehat{G}$ allow us to conclude that

$$\left(\widehat{\Phi\varphi}\right)(g) + \left(\widehat{\Psi\Phi}\right)(h) = \left(\widehat{\Phi\varphi\Psi\Phi}\right)\left(\Phi(g_1^*) + \Phi(h_1^*)\right).$$

Using this result in the expression obtained above for $\widehat{\Phi}(\widehat{g} + \widehat{h})$, we get:

$$\begin{aligned}\widehat{\Phi}(\widehat{g} + \widehat{h}) &= \left(\widehat{\Phi\varphi}\right)(g) + \left(\widehat{\Psi\Phi}\right)(h) = \\ &= \widehat{\Phi}\left(\widehat{\varphi}(g)\right) + \widehat{\Phi}\left(\widehat{\Psi}(h)\right) = \\ &= \widehat{\Phi}(\widehat{g}) + \widehat{\Phi}(\widehat{h}).\end{aligned}$$

**Lemma 6.** *The addition in $\widehat{G}$ defined in Lemma 3, is the only possible extension of the addition of the group $G$ to $\widehat{G}$, that turns $\widehat{G}$ into an abelian group and $\widehat{\Phi}$ (for each $\Phi \in H$) an endomorphism on $\widehat{G}$.*

*Proof.* Lemmas 3, 4 and 5 show us that the operation in $\widehat{G}$ defined in Lemma 3 is an extension to $\widehat{G}$ of the addition of the group $G$, that turns $\widehat{G}$ into an abelian group and $\widehat{\Phi}$ an endomorphism on $\widehat{G}$. Suppose now that exists another extension of the addition, $+$, of the group $G$ to $\widehat{G}$, lets say, $\boxplus$, such that $\widehat{G}$ equipped with $\boxplus$ is an abelian group and $\widehat{\Phi}$ an endomorphism on the group $(\widehat{G}, \boxplus)$. Let then $\widehat{g}, \widehat{h} \in \widehat{G}$ be arbitrarily chosen and, consonant to hypothesis (c), $\Phi, \Psi \in H$ and $g, h \in G$ such that

$$\widehat{g} = \widehat{\Phi}(g) \quad \text{and} \quad \widehat{h} = \widehat{\Psi}(h).$$

As we know, $\Phi\colon G_\Phi \longrightarrow G$ and $\Psi\colon G_\Psi \longrightarrow G$ are surjective homomorphisms and hence, since $g$ and $h$ belong to $G$, we have

$$g = \Psi(g^*) \quad \text{and} \quad h = \Phi(h^*)$$

with

$$g^* \in \Psi^{-1}(g) \quad \text{and} \quad h^* \in \Phi^{-1}(h)$$

arbitrarily chosen. Thus

$$\widehat{g} = \widehat{\Phi}(g) = \widehat{\Phi}\left(\Psi(g^*)\right) = \widehat{\Phi}\left(\widehat{\Psi}(g^*)\right) = \left(\widehat{\Phi\Psi}\right)(g^*)$$

and

$$\widehat{h} = \widehat{\Psi}(h) = \widehat{\Psi}\left(\Phi(h^*)\right) = \widehat{\Psi}\left(\widehat{\Phi}(h^*)\right) = \left(\widehat{\Phi\Psi}\right)(h^*)$$

where we resorted to hypothesis (b). Hence,

$$\widehat{g} \boxplus \widehat{h} = \left(\widehat{\Phi\Psi}\right)(g^*) \boxplus \left(\widehat{\Phi\Psi}\right)(g^*),$$

and since, by hypothesis, $\widehat{\Phi\Psi}$ is an endomorphism on the group $(\widehat{G}, \boxplus)$, we have

$$\widehat{g} \boxplus \widehat{h} = \left(\widehat{\Phi\Psi}\right)(g^* \boxplus h^*).$$

But $g^* \boxplus h^* = g^* + h^*$ since, by hypothesis, $\boxplus$ is an extension of $+$. Thus,

$$\widehat{g} \boxplus \widehat{h} = \left(\widehat{\Phi\Psi}\right)(g^* + h^*)$$

where, as we saw above,

$$g^* \in \Psi^{-1}(g) \quad \text{and} \quad h^* \in \Phi^{-1}(h).$$

Recurring now to the definition of addition, $+$, given in Lemma 3, we see that

$$\left(\widehat{\Phi\Psi}\right)(g^* + h^*) = \widehat{g} + \widehat{h}$$

and hence we conclude that

$$\widehat{g} \boxplus \widehat{h} = \widehat{g} + \widehat{h},$$

so finishing the proof of this lemma and, also, of the proposition. ■

## 2.25 $\widetilde{\mathbb{G}}$-distributions Axiomatic — A Simplified Version

The conditions (a) to (d) of Proposition 2.24, stated relative to an abelian, surjective, and with identity $S$-group $\mathbb{G} = (G, H)$, are, as we saw in 2.23, attended by the sets $\widehat{G}$ and $\widehat{H}$ that compose any strict and closed extension, $\widehat{\mathbb{G}} = (\widehat{G}, \widehat{H})$, of the $S$-group $\mathbb{G}$. On the other hand, the referred proposition tells us that if $\widehat{G}$ and $\widehat{H}$ are any sets, not necessarily the components of a strict and closed extension of $\mathbb{G}$, that satisfy the conditions (a) to (d), then, apparentely, we can transform them, through the definition of an addition for the elements of the set $\widehat{G}$, in the components of such an extension of $\mathbb{G} = (G, H)$. In fact, once attended the referred conditions (a) to (d) by the sets $\widehat{G}$ and $\widehat{H}$, then, the addition defined in Lemma 3 of the proof of proposition at hand, transforms $\widehat{G}$ into an abelian group that has $G$ as a subgroup and $\widehat{H}$ into a semigroup of endomorphisms on $\widehat{G}$ which is a prolongation to $\widehat{G}$ of the semigroup $H$, turning then $(\widehat{G}, \widehat{H})$ into a $S$-group that is an extension of the $S$-group $\mathbb{G} = (G, H)$, and, even more, as the condition (c) is attended, into a closed extension of $\mathbb{G}$. Further, from condition (d) we see that $(\widehat{G}, \widehat{H})$ also is a strict extension of $\mathbb{G}$: in fact, from (d) results that $N(\widehat{\Phi}) = N(\Phi)$ and, then, $G \cap N(\widehat{\Phi}) = N(\Phi)$, that, in turn, according to Proposition 1.12, is a necessary and sufficient condition for the prolongation $\widehat{\Phi}$ of $\Phi$ to $\widehat{G}$ to be strict.

Now, the addition defined in Lemma 3 above referred, responsible for the metamorphoses of the sets $\widehat{G}$ and $\widehat{H}$ into the strict and closed extension $(\widehat{G}, \widehat{H})$ of the $S$-group $\mathbb{G} = (G, H)$, is, as establishes the proposition at hand, the only one capable of such transformation. Since the pair $(\widehat{G}, \widehat{H})$ thus metamorphized is a model of the $\widetilde{\mathbb{G}}$-distributions axiomatic, unique unless isomorphism since the referred axiomatic, as we saw in 2.22, is categoric, seems then more than reasonable that: taking the conditions (a) to (d) of Proposition 2.24 as axioms defining implicitly the primitive terms "$\widetilde{\mathbb{G}}$-distribution" (the members of $\widehat{G}$) and "$\widetilde{\mathbb{G}}$-distribution derivative" (the elements of $\widehat{H}$), we obtain an equivalent axiomatic to the one defined in 2.20.

The above observations lead us then to formulate a new axiomatic, now with only two primitive terms, "$\widetilde{\mathbb{G}}$-distribution" and "derivative", with the same precedent theories of the axiomatic formulated in 2.20, namely, the classical logic, the set theory, and the $S$-groups theory. The axioms, stated relative to an arbitrary, yet fixed, abelian, surjective, and with identity $S$-group, $\mathbb{G} = (G, H)$, are presented ahead using the primitive symbols $\widehat{G}$ and $\widehat{H}$ which denotes, respectively, the class of the $\widetilde{\mathbb{G}}$-distributions and the set of its derivatives.

**Axiom 1** $G \subseteq \widehat{G}$.

**Axiom 2** The derivatives, that is, the members $\widehat{\Phi} \in \widehat{H}$, are functions with domain and codomain equal to $\widehat{G}$, such that:

- **(a)** each $\widehat{\Phi} \in \widehat{H}$ is an extension to $\widehat{G}$ of an unique homomorphism $\Phi \in H$, that is, $\widehat{\Phi}(g) = \Phi(g)$ for every $g \in G_\Phi$, for an unique $\Phi \in H$;
- **(b)** the class $\widehat{H}$, with the usual composition of functions, is a semigroup isomorphic to the semigroup $H$ and the function $\widehat{\ }$ defined by

$$\begin{aligned} \widehat{\ } : H &\longrightarrow \widehat{H} \\ \Phi &\longmapsto \widehat{\Phi} := \text{the extension of } \Phi \text{ to } \widehat{G}, \end{aligned}$$

is an isomorphism (between semigroups).

**Axiom 3** For each $\widehat{g} \in \widehat{G}$, there exist $\Phi \in H$ and $g \in G$ such that $\widehat{g} = \widehat{\Phi}(g)$.

**Axiom 4** If $g, h \in G$ and $\Phi \in H$, then $\widehat{\Phi}(g) = \widehat{\Phi}(h)$ if and only if $g - h \in N(\Phi)$.

**Remark.** Such as for the axiomatic in 2.20, here too, the symbol $\widehat{\Phi}$ denotes the unique extension of $\Phi$ to $\widehat{G}$ that belongs to $\widehat{H}$, that is, $\widehat{\Phi}$ is the image of $\Phi \in H$ by the isomorphism $\widehat{\ }: H \longrightarrow \widehat{H}$.

## 2.26 The Equivalence of the Axiomatics

The remarks in the previous item, which precede the presentation of the axioms, strongly suggest that this new axiomatic, defined in the referred item, is equivalent to that formulated in 2.20, or, more explicitly, that: every axiom of any one of the axiomatics at hand is an axiom or theorem of the other and, also, every primitive term of one is a primitive term or defined term of the other. In fact, this is what happens, as ahead, in an explicit form, we prove.

First we observe that in the axiomatic of the item 2.25, henceforth denominated "axiomatic 2" to distinguish it from that formulated in 2.20 to whose we will refer as "axiomatic 1", the axioms are, exactly, the conditions (a) to (d) of Proposition 2.24 and, as occurs with these conditions, the axioms are stated relative to an abelian, surjective, and with identity $S$-group $\mathbb{G} = (G, H)$. Therefore, the six lemmas that compose the proof of the cited proposition, that have its proofs exclusively backed on these four conditions and the properties of the $S$-group $\mathbb{G}$, are theorems of the axiomatic 2. In particular, Lemma 3 allows to introduce in the axiomatic 2, as a defined term, the "addition", defining it as the operation described in the alluded lemma. Furthermore, with the "addition" incorporated via definition to axiomatic 2, we conclude, by simple inspection of the lemmas, now theorems of axiomatic 2, that the axioms of the axiomatic 1 that does not figure as axioms in the axiomatic 2, namely, Axiom 2, Axiom 3(a) and Axiom 4, only the Axiom 4 does not find itself among the "lemmas ≡ theorems" referred. However, also this Axiom 4 (from the axiomatic 1) is a theorem of the axiomatic 2 as we show below.

We want then to prove, starting with the axioms and theorems of axiomatic 2, that for each $\Phi \in H$,

$$\text{if} \quad g \in G \quad \text{and} \quad g \notin G_\Phi, \quad \text{then} \quad \widehat{\Phi}(g) \notin G.$$

It is easy to see that the implication above is equivalent to

$$\text{if} \quad g \in G \quad \text{and} \quad \widehat{\Phi}(g) \in G, \quad \text{then} \quad g \in G_\Phi.$$

We will prove then this last implication and, for such, let $\Phi \in H$ and $g \in G$ arbitrarily fixed and such that $\widehat{\Phi}(g) \in G$. Since $\Phi\colon G_\Phi \longrightarrow G$ is surjective, there exists $g' \in G_\Phi$ such that $\Phi(g') = \widehat{\Phi}(g)$, from where we conclude, taking into account the Axiom 2 (from axiomatic 2), that

$$\widehat{\Phi}(g') = \widehat{\Phi}(g),$$

and hence, resorting now to Axiom 4 (from axiomatic 2), we obtain that $g' - g \in N(\Phi) \subseteq G_\Phi$. Consequently, $g' - g \in G_\Phi$ and, since $G_\Phi$ is a group and $g' \in G_\Phi$, then $g' - (g' - g) = g \in G_\Phi$.

From this we can now affirm that: all primitive terms of axiomatic 1 are primitive ($\widetilde{\mathbb{G}}$-distribution and derivative) or defined (addition) terms of the axiomatic 2 and all axioms of axiomatic 1 are axioms or theorems of axiomatic 2. Regarding the converse we have, clearly, that the primitive terms of the axiomatic 2 are also primitive terms of the axiomatic 1; about the axioms, only Axiom 4 (of axiomatic 2) does not figure in the list of axioms of axiomatic 1. However, as proved ahead, the referred axiom is a theorem of axiomatic 1.

We must then prove that: for $g, h \in G$ and $\Phi \in H$,

$$\widehat{\Phi}(g) = \widehat{\Phi}(h) \quad \text{if and only if} \quad g - h \in N(\Phi).$$

Hence, let $g, h \in G$ and $\Phi \in H$ be arbitrarily fixed. If $g - h \in N(\Phi)$, then, $g - h \in G_\Phi$ and, thus, by Axiom 3(b) of axiomatic 1,

$$\widehat{\Phi}(g - h) = \Phi(g - h) = 0.$$

Now, among the logical consequences (theorems) of the axiomatic 1, obtained in 2.21, we have that $\widehat{G}$ is an abelian group that has $G \subseteq \widehat{G}$ as a subgroup (2.21(e)) and $\widehat{H}$ is a semigroup of endomorphisms on $\widehat{G}$ (2.21(f)). Hence, from $\widehat{\Phi}(g - h) = 0$ results that $\widehat{\Phi}(g) = \widehat{\Phi}(h)$. Conversely, suppose now that $\widehat{\Phi}(g) = \widehat{\Phi}(h)$, that is, taking into account the above cited consequences of axiomatic 1, that $g - h \in N(\widehat{\Phi})$. But, it is also a theorem of axiomatic 1 that $(\widehat{G}, \widehat{H})$ is a strict and closed extension of the $S$-group (abelian, surjective, and with identity) $\mathbb{G} = (G, H)$ (see 2.21(h)). Hence being, and taking into account Proposition 1.16(b) (which is also a theorem of axiomatic 1 since elapsed of fulfilled hypothesis of this axiomatic), we have that $N(\widehat{\Phi}) = N(\Phi)$ and hence, $g - h \in N(\Phi)$, concluding the proof of the converse.

The axiomatics 1 and 2, as we anticipated in the item 2.25, are, in fact, equivalent.

# Schwartz' Distributions

## 2.27 Preliminaries

The idea of generalize the notions of continuous function and differentiability, motivated by questions associated with solutions for partial differential equations, leading to the notion of distributions and their derivatives, has a long and interesting history which well illustrate the important role of the strong connection between the development of mathematical concepts and theories and applied questions, especially those in the domain of natural sciences. Regarding this history, we suggest the Jesper Lützen's book[11] at which, among several mathematicians referred to as important for the development of the concept of distribution, the alluded author highlights Sergej L'vovic Sobolev and Laurent Schwartz as the ones with higher relevance writing that:

[11] Lützen, J. *The Prehistory of the Theory of Distributions.* Springer-Verlag, New York, 1982.

"Thus Sobolev invented distributions, but it was Schwartz who created the theory of distributions."

In this section we present, schematically, a few aspects of the L. Schwartz' distributions theory[12] that will prove themselves relevant for us to establish the relation between the Schwartz' distributions and the $\widetilde{\mathbb{C}}(\Omega)$-distributions (the $\widetilde{\mathbb{G}}$-distributions specialized (particularized) to the $S$-group $\mathbb{C}(\Omega) = (C(\Omega), \partial(\Omega))$ of the continuous functions). We begin, through the definition below, with a recapitulation of some concepts, as well as the introduction of some notation.

## 2.28 Definition

It is convenient, before introducing additional notation, to remember those defined in the item 1.2: $C^\infty(\Omega)$ denotes the set of functions $f \in C(\Omega)$ with partial derivatives of all orders continuous on $\Omega$, being $\Omega \subseteq \mathbb{R}^n$ a non-empty open set; $C(\Omega)$, in turn, is the set of functions with domain $\Omega$, assuming complex values, that are continuous on the open set $\Omega \subseteq \mathbb{R}^n$. Hence, in particular, $C^\infty(\mathbb{R}^n)$ is the set of functions of complex values, defined in $\mathbb{R}^n$, with partial derivatives of all orders continuous on $\mathbb{R}^n$, that is,

$$C^\infty(\mathbb{R}^n) = \Big\{\phi \in C(\mathbb{R}^n)\colon \partial^\alpha_{\mathbb{R}^n}(\phi) \in C(\mathbb{R}^n) \quad \text{for every} \quad \alpha \in \mathbb{N}^n\Big\}.$$

Still from 1.2, for $n \in \mathbb{N}$, $C^n(\Omega)$ (respectively, $C^n(\mathbb{R}^n)$) is the class of functions of $C(\Omega)$ (respectively, $C(\mathbb{R}^n)$) with partial derivatives of order less than or equal to $n$ continuous on $\Omega$ (respectively, $\mathbb{R}^n$).

In the definitions given ahead, we freely use the notions of closure of a set, compact set, vector space, (vector) subspace, linear transformation, among others; to the reader unfamiliar with these concepts we suggest, for analytical nature ones the text of W. Rudin cited in footnote 3 (page 11), and for those in the domain of algebra, the book of G. Birkhoff and S. MacLane indicated in footnote 8 (page 19).

**(a)** Let $\Omega \subseteq \mathbb{R}^n$ be a non-empty open set and $f \in C(\Omega)$. We define the **support of $f$** as the closure of the following set:

$$\Big\{x \in \Omega\colon f(x) \neq 0\Big\}.$$

**(b)** Let $K \subseteq \mathbb{R}^n$, $K \neq \varnothing$, be a compact set. We denote by $C^\infty_K(\mathbb{R}^n)$ the set of the functions $\phi \in C^\infty(\mathbb{R}^n)$ whose support is included in $K$, that is,

$$C^\infty_K(\mathbb{R}^n) \coloneqq \Big\{\phi \in C^\infty(\mathbb{R}^n)\colon \text{support of } \phi \subseteq K\Big\}.$$

[12] Schwartz, L. *Théorie des Distribuitions.* Tome I, Hermann, Paris, 1950; Tome II, Hermann, Paris, 1951.

We know that $C(\Omega)$ with usual operations of function addition and function multiplication by a complex number is a vector space and that $C^\infty(\Omega) \subseteq C(\Omega)$ and $C^n(\Omega) \subseteq C(\Omega)$, $n \in \mathbb{N}$, are subspaces of $C(\Omega)$. Results immediately that $C_K^\infty(\mathbb{R}^n) \subseteq C^\infty(\mathbb{R}^n)$ is a subspace of the vector space $C^\infty(\mathbb{R}^n)$. Furthermore, if the compact $K$ is included in an open set $\Omega$ ($K \subseteq \Omega$), then, $C_K^\infty(\mathbb{R}^n)$ is isomorphic to a subspace of $C^\infty(\Omega)$, namely,

$$\Big\{\phi \in C^\infty(\Omega) \colon \text{support of } \phi \subseteq K \subseteq \Omega\Big\}.$$

In this case, it is common to refer to $C_K^\infty(\mathbb{R}^n)$ through the notation $C_K^\infty(\Omega)$ that, however, we will not use.[13]

**(c)** Let $\Omega \subseteq \mathbb{R}^n$ a non-empty open set. We define now the set $C_\Omega^\infty(\mathbb{R}^n)$ as the union of the sets $C_K^\infty(\mathbb{R}^n)$ for which the compact $K$ is included in $\Omega$, that is,

$$C_\Omega^\infty(\mathbb{R}^n) \coloneqq \bigcup_{\substack{K \subseteq \Omega \\ K \text{ compact}}} C_K^\infty(\mathbb{R}^n).$$

Hence, we have that, $\phi \in C_\Omega^\infty(\mathbb{R}^n)$ if and only if $\phi \in C^\infty(\mathbb{R}^n)$ and the support of $\phi$ is a compact subset of $\Omega$.

$C_\Omega^\infty(\mathbb{R}^n)$ also is a subspace of the vector space $C^\infty(\mathbb{R}^n)$, which Schwartz denominated space of **test functions**.

The space of test functions is also isomorphic, and thus being able to be identified, to a subspace of $C^\infty(\Omega)$, namely, that of the functions $\phi \in C^\infty(\Omega)$ whose support is a compact included in $\Omega$.

**(d)** A **linear functional** on $C_\Omega^\infty(\mathbb{R}^n)$ is, by definition, a (any) linear transformation from the space of test functions into the (complex) vector space of the complex numbers $C$, that is, a function $\Lambda \colon C_\Omega^\infty(\mathbb{R}^n) \longrightarrow C$ such that

$$\Lambda(a\phi + b\varphi) = a\Lambda(\phi) + b\Lambda(\varphi)$$

whatever $\phi, \varphi \in C_\Omega^\infty(\mathbb{R}^n)$ and $a, b \in C$.

## 2.29 Remark

The distributions, as defined by L. Schwartz, are the linear functionals on $C_\Omega^\infty(\mathbb{R}^n)$ that are continuous relative to a specific topology defined for the vector space of the test functions, which makes it a topological vector space. However, taking into account that our purpose here consists only of establishing the relation between the $\widetilde{\mathbb{C}}(\Omega)$-distributions and the Schwartz ones, it is not necessary to involve ourselves with this subjacent topological vector space once that the Definition 2.30 formulated ahead, equivalent (as proved in the literature[14]) to that of topological nature referred above, is, as we will see, appropriate to our goals and does not involve us with subtle topological aspects.

[13] The reason we highlight this isomorphism and the referred notation justified by it, is to conciliate our notation with the one used by many other texts.

## 2.30 Definition

**(a)** A **distribution in** $\Omega$ ($\Omega \subseteq \mathbb{R}^n$ an open set), or a **distribution of domain** $\Omega$, is a linear functional on the space of test functions, $\Lambda\colon C^{\infty}_{\Omega}(\mathbb{R}^n) \longrightarrow C$, that satisfies the following condition: for each compact $K \subseteq \Omega$ there exist a constant $C_K < \infty$ and an integer $N_K \geqslant 0$, generally dependents of $K$, such that

$$|\Lambda(\phi)| \leqslant C_K \max\left\{\left|\Big(\partial^{\alpha}_{\mathbb{R}^n}(\phi)\Big)(x)\right| : x \in \Omega \quad \text{and} \quad |\alpha| \leqslant N_K\right\} \tag{2.30-1}$$

for every $\phi \in C^{\infty}_{K}(\mathbb{R}^n)$.

We denote by $D'(\Omega)$ the class of the distributions in $\Omega$ (of domain $\Omega$).

The condition (2.30-1) is equivalent to affirm that the linear functional $\Lambda$ is continuous relative to the previously mentioned "specific topology" which makes the vector space of the test functions, $C^{\infty}_{\Omega}(\mathbb{R}^n)$, a topological vector space. Hence, the distributions of domain $\Omega$, in the context of this topology, are the continuous linear functionals on $C^{\infty}_{\Omega}(\mathbb{R}^n)$.

**(b)** Let $\Lambda$ be a distribution in $\Omega$, that is, $\Lambda \in D'(\Omega)$. If $\Lambda$ is such that the condition (2.30-1) is satisfied for a choice of $N_K = N$ independent of the compact $K$, that is, so that a single $N$ can be chosen for every $K$, then, the smaller of these integers $N$ is denominated the **order** of $\Lambda$, and $\Lambda$ is said to be a distribution of **finite order**; if no such $N$ can be so chosen (independent of $K$) we say that $\Lambda \in D'(\Omega)$ is of **infinite order**. Hence, $\Lambda \in D'(\Omega)$ is of finite order if and only if there exists an integer $N \geqslant 0$ such that, for each compact $K \subseteq \Omega$ corresponds a constant $C_K < \infty$ (which generally depends on $K$) in such a way that

$$|\Lambda(\phi)| \leqslant C_K \max\left\{\left|\Big(\partial^{\alpha}_{\mathbb{R}^n}(\phi)\Big)(x)\right| : x \in \Omega \quad \text{and} \quad |\alpha| \leqslant N\right\}$$

for every $\phi \in C^{\infty}_{K}(\mathbb{R}^n)$.

We will employ the symbol $D'_f(\Omega)$ to denote the class of distributions in $\Omega$ of finite order.

[14] A nice and concise presentation of the Schwartz' distributions theory can be found in the following book by W. Rudin, where the underlying topology is explicitly presented: Rudin, W. *Functional Analysis*. 2nd ed. McGraw-Hill Book Company, 1991.

## 2.31 The Group $D'(\Omega)$ and the Subgroups $D'_f(\Omega)$ and $C'(\Omega)$

It can be verified that the class $D'(\Omega)$ of the distributions of domain $\Omega$, with the usual operations of addition of linear functionals and multiplication of a linear functional by a complex number, is a vector space. Beyond that, the class of distributions in $\Omega$ of finite order, $D'_f(\Omega) \subseteq D'(\Omega)$, is a subspace of the vector space $D'(\Omega)$, that is, the sum of distributions in $\Omega$ of finite order, as well as the product of such a distribution by a complex number, also are distributions of finite order and same domain $\Omega$. In particular, $D'(\Omega)$ with the referred addition is an abelian group that has $D'_f(\Omega) \subseteq D'(\Omega)$ as a subgroup.

Let now $f$ be a function arbitrarily chosen in the group $C(\Omega)$ of continuous functions. Since the function

$$h(x) \coloneqq f(x)\phi(x) \quad \text{for} \quad x \in \Omega$$

is integrable in $\Omega$ whatever $\phi \in C^\infty_\Omega(\mathbb{R}^n)$, we can associate to $f$ the linear functional on $C^\infty_\Omega(\mathbb{R}^n)$ defined by:

$$\begin{aligned} \Lambda_f : C^\infty_\Omega(\mathbb{R}^n) &\longrightarrow C \\ \phi &\longmapsto \Lambda_f(\Phi) \coloneqq \textstyle\int_\Omega f(x)\phi(x)\,\mathrm{d}x\,. \end{aligned}$$

It results that $\Lambda_f$ is a distribution in $\Omega$, more specifically, $\Lambda_f$ is a distribution in $\Omega$ of finite order: $\Lambda_f \in D'_f(\Omega)$. In fact, for if $K \subseteq \Omega$ is a compact, then, for every $\phi \in C^\infty_K(\mathbb{R}^n)$, we have:

$$|\Lambda_f(\phi)| \leqslant \int_K |f(x)||\phi(x)|\,\mathrm{d}x \leqslant \left(\int_K |f(x)|\,\mathrm{d}x\right) \max\left\{|\phi(x)| : x \in \Omega\right\}$$

that is,

$$|\Lambda_f(\phi)| \leqslant C_K \max\left\{\left|\left(\partial^\alpha_{\mathbb{R}^n}(\phi)\right)(x)\right| : x \in \Omega \quad \text{and} \quad |\alpha| \leqslant 0\right\}$$

where we take

$$C_K \coloneqq \int_K |f(x)|\,\mathrm{d}x\,.$$

From this last inequality, with Definition 2.30(b) in mind, we conclude that $\Lambda_f$ is a distribution in $\Omega$ of order zero.

We say that $\Lambda_f$, as above defined, is the distribution in $\Omega$ **induced** by the function $f \in C(\Omega)$. Let us consider now the following function:

$$\begin{aligned} \eta : C(\Omega) &\longrightarrow D'_f(\Omega) \\ f &\longmapsto \eta(f) \coloneqq \Lambda_f \end{aligned}$$

where $\Lambda_f$ is the distribution induced by $f$. It results immediately that:

$$\eta(f+g) = \eta(f) + \eta(g) \quad \text{for any} \quad f, g \in C(\Omega),$$

that is, $\eta$ is a homomorphism on the abelian group $C(\Omega)$ of the continuous functions on $\Omega$, into the group $D'_f(\Omega)$ of the distributions in $\Omega$ of finite order. Beyond that, if for $f, g \in C(\Omega)$, $\eta(f) = \eta(g)$, that is, $\eta(f-g) = 0$, then

$$\Lambda_{f-g}(\phi) = \int_\Omega (f-g)(x)\phi(x)\,\mathrm{d}x = 0 \quad \text{for every} \quad \phi \in C_\Omega^\infty(\mathbb{R}^n)$$

which only occurs if $f = g$. Hence, $\eta\colon C(\Omega) \longrightarrow D'_f(\Omega)$ is an injective homomorphism and, therefore, an isomorphism on the group $C(\Omega)$ onto the subgroup

$$C'(\Omega) := \eta\Big(C(\Omega)\Big) = \Big\{\Lambda_f \colon f \in C(\Omega)\Big\}$$

of the group $D'_f(\Omega)$. From this, we conclude that the group $D'_f(\Omega)$, as well as the group $D'(\Omega)$, admit the group $C(\Omega)$ as a subgroup (see Definition 1.6(a)). Hence, through the isomorphism $\eta\colon C(\Omega) \longrightarrow C'(\Omega)$, we identify the functions $f \in C(\Omega)$ with the distributions induced by them,

$$f \equiv \Lambda_f,$$

as well as we identify $C(\Omega)$ with $C'(\Omega)$,

$$C(\Omega) \equiv C'(\Omega).$$

## 2.32 Derivatives of Distributions

As we saw in the previous item, through the isomorphism $\eta\colon C(\Omega) \longrightarrow C'(\Omega)$, where $\eta(f) = \Lambda_f$ is the distribution induced by $f \in C(\Omega)$, we identify the functions $f \in C(\Omega)$ with the distributions $\Lambda_f$ of $C'(\Omega)$ (the group of the distributions induced by the continuous functions on $\Omega$): $f \equiv \Lambda_f$ if $f \in C(\Omega)$. Let us suppose now that $f \in C^{|\alpha|}(\Omega)$ with $\alpha$ a multi-index (see 1.2(f) and (h)). In this case we have $\partial_\Omega^\alpha(f) \in C(\Omega)$ and thus $\partial_\Omega^\alpha(f) \equiv \Lambda_{\partial_\Omega^\alpha(f)}$. Now, $\Lambda_{\partial_\Omega^\alpha(f)}\colon C_\Omega^\infty(\mathbb{R}^n) \longrightarrow C$ is, as we know, the distribution defined by:

$$\Lambda_{\partial_\Omega^\alpha(f)}(\phi) = \int_\Omega \Big(\partial_\Omega^\alpha(f)\Big)(x)\phi(x)\,\mathrm{d}x$$

for every $\phi \in C_\Omega^\infty(\mathbb{R}^n)$, or yet, after $|\alpha|$ integrations by parts (at which the compactness of the support of $\phi$ is used),

$$\Lambda_{\partial_\Omega^\alpha(f)}(\phi) = (-1)^{|\alpha|}\int_\Omega f(x)\Big(\partial_{\mathbb{R}^n}^\alpha(\phi)\Big)(x)\,\mathrm{d}x = (-1)^{|\alpha|}\Lambda_f\Big(\partial_{\mathbb{R}^n}^\alpha(\phi)\Big)$$

for every $\phi \in C_\Omega^\infty(\mathbb{R}^n)$.

In short, we have, for $f \in C^{|\alpha|}(\Omega)$, that:

$$f \equiv \Lambda_f \quad \text{where} \quad \Lambda_f(\phi) = \int_\Omega f(x)\phi(x)\,\mathrm{d}x$$

and

$$\partial^{\alpha}_{\Omega}(f) \equiv \Lambda_{\partial^{\alpha}_{\Omega}(f)} \quad \text{where} \quad \Lambda_{\partial^{\alpha}_{\Omega}(f)}(\phi) = (-1)^{|\alpha|}\Lambda_f\Big(\partial^{\alpha}_{\mathbb{R}^n}(\phi)\Big),$$

for every $\phi \in C^{\infty}_{\Omega}(\mathbb{R}^n)$.

Thus, it is reasonable to admit that a definition for the derivative of a distribution in $\Omega$, with the pretension to be taken as a generalization of the derivative $\partial^{\alpha}_{\Omega}$, if denoted by, lets say, $D^{\alpha}_{\Omega}$, should attend the following condition:

$$D^{\alpha}_{\Omega}(\Lambda_f) = \Lambda_{\partial^{\alpha}_{\Omega}(f)},$$

i.e., that $D^{\alpha}_{\Omega}(\Lambda_f)$ should be a distribution in $\Omega$ given by

$$\Big(D^{\alpha}_{\Omega}(\Lambda_f)\Big)(\phi) = (-1)^{|\alpha|}\Lambda_f\Big(\partial^{\alpha}_{\mathbb{R}^n}(\phi)\Big) \quad \text{for every} \quad \phi \in C^{\infty}_{\Omega}(\mathbb{R}^n).$$

We observe now that the second member of the above equality is well-defined (for every $\phi \in C^{\infty}_{\Omega}(\mathbb{R}^n)$), not only for the distribution induced by the function $f \in C(\Omega)$, $\Lambda_f$, but also for any distribution $\Lambda$ of domain $\Omega$. Hence, it is reasonable that, given a multi-index $\alpha \in \mathbb{N}^n$ and a distribution $\Lambda \in D'(\Omega)$, we take for the derivative of $\Lambda$, $D^{\alpha}_{\Omega}(\Lambda)$, the following linear functional on $C^{\infty}_{\Omega}(\mathbb{R}^n)$:

$$\begin{aligned} D^{\alpha}_{\Omega}(\Lambda) : C^{\infty}_{\Omega}(\mathbb{R}^n) &\longrightarrow C \\ \phi &\longmapsto \Big(D^{\alpha}_{\Omega}(\Lambda)\Big)(\phi) := (-1)^{|\alpha|}\Lambda\Big(\partial^{\alpha}_{\mathbb{R}^n}(\phi)\Big). \end{aligned}$$

However, even if $D^{\alpha}_{\Omega}(\Lambda)$ as given above were well-defined, as it is, and be, as it is, a linear functional on $C^{\infty}_{\Omega}(\mathbb{R}^n)$, for us to accept it as a derivative of $\Lambda$, it is essential for this functional to be, also, a distribution in $\Omega$. But, this is in fact the case as we will prove ahead.

Since $\Lambda$ is a distribution in $\Omega$, then, for each compact $K \subseteq \Omega$ there exist a constant $C_K < \infty$ and an integer $N_K \geqslant 0$ such that

$$|\Lambda(\phi)| \leqslant C_K \max\left\{\left|\Big(\partial^{\alpha'}_{\mathbb{R}^n}(\phi)\Big)(x)\right| : x \in \Omega \quad \text{and} \quad |\alpha'| \leqslant N_K\right\}$$

for every $\phi \in C^{\infty}_{K}(\mathbb{R}^n)$, and hence (taking into account that $\partial^{\alpha}_{\mathbb{R}^n}(\phi) \in C^{\infty}_{K}(\mathbb{R}^n)$ for every $\phi \in C^{\infty}_{K}(\mathbb{R}^n)$),

$$\left|\Big(D^{\alpha}_{\Omega}(\Lambda)\Big)(\phi)\right| = \left|\Lambda\Big(\partial^{\alpha}_{\mathbb{R}^n}(\phi)\Big)\right| \leqslant C_K \max\left\{\left|\Big(\partial^{\alpha'+\alpha}_{\mathbb{R}^n}(\phi)\Big)(x)\right| : x \in \Omega \quad \text{and} \quad |\alpha'| \leqslant N_K\right\},$$

or yet, making $\beta := \alpha' + \alpha$,

$$\left|\Big(D^{\alpha}_{\Omega}(\Lambda)\Big)(\phi)\right| \leqslant C_K \max\left\{\left|\Big(\partial^{\beta}_{\mathbb{R}^n}(\phi)\Big)(x)\right| : x \in \Omega \quad \text{and} \quad |\beta| \leqslant N_K + |\alpha|\right\}$$

from where we conclude that $D^{\alpha}_{\Omega}(\Lambda) \in D'(\Omega)$.

The definition given below formalizes the ideas discussed above.

## 2.33 Definition

Let $\alpha \in \mathbb{N}^n$ be a multi-index and $\Omega \subseteq \mathbb{R}^n$ be a non-empty open set. $D_\Omega^\alpha$ is then the function with domain and codomain equal to $D'(\Omega)$ defined by:

$$\begin{aligned} D_\Omega^\alpha : D'(\Omega) &\longrightarrow D'(\Omega) \\ \Lambda &\longmapsto D_\Omega^\alpha(\Lambda) \end{aligned}$$

with $D_\Omega^\alpha(\Lambda)$ being the following distribution in $\Omega$,

$$\begin{aligned} D_\Omega^\alpha(\Lambda) : C_\Omega^\infty(\mathbb{R}^n) &\longrightarrow C \\ \phi &\longmapsto \Big(D_\Omega^\alpha(\Lambda)\Big)(\phi) := (-1)^{|\alpha|}\Lambda\Big(\partial_{\mathbb{R}^n}^\alpha(\phi)\Big). \end{aligned}$$

Fixed a non-empty open set $\Omega \subseteq \mathbb{R}^n$, we say that the functions $D_\Omega^\alpha$ above defined, one for each multi-index $\alpha \in \mathbb{N}^n$, are the **derivatives** of the distributions of domain $\Omega$ and we denote by $D(\Omega)$ the class of these derivatives, that is,

$$D(\Omega) := \Big\{D_\Omega^\alpha \colon \alpha \in \mathbb{N}^n\Big\}.$$

We highlight that $D_\Omega^0 = I_{D'(\Omega)} \in D(\Omega)$, with $I_{D'(\Omega)}$ being the identity function on $D'(\Omega)$: $I_{D'(\Omega)}(\Lambda) = \Lambda$ for every $\Lambda \in D'(\Omega)$. In fact, from the definition of $D_\Omega^\alpha(\Lambda)$ and taking into account that $\partial_{\mathbb{R}^n}^0 = I_{C(\mathbb{R}^n)}$ (see 1.2(i)), we have, for $\alpha = 0 = (0, \ldots, 0)$ and $\Lambda \in D'(\Omega)$ arbitrarily fixed, that:

$$\Big(D_\Omega^0(\Lambda)\Big)(\phi) = (-1)^0\Lambda\Big(\partial_{\mathbb{R}^n}^0(\phi)\Big) = \Lambda(\phi) \quad \text{for every} \quad \phi \in C_\Omega^\infty(\mathbb{R}^n)$$

and hence, $D_\Omega^0(\Lambda) = \Lambda$.

## 2.34 The Semigroup $D(\Omega)$

Let us consider now the class $D(\Omega)$ of the derivatives of the distributions in $\Omega$ with the usual operation of composition of functions. Clearly, the function d defined ahead is a "natural" bijection between the set $\partial(\Omega) = \{\partial_\Omega^\alpha \colon \alpha \in \mathbb{N}^n\}$ of the usual derivatives and the set of the derivatives of distributions in $\Omega$, $D(\Omega) = \{D_\Omega^\alpha \colon \alpha \in \mathbb{N}^n\}$:

$$\begin{aligned} \mathrm{d} : \partial(\Omega) &\longrightarrow D(\Omega) \\ \partial_\Omega^\alpha &\longmapsto \mathrm{d}\left(\partial_\Omega^\alpha\right) := D_\Omega^\alpha. \end{aligned}$$

As we know (see the items 1.2 and 1.5), $\partial(\Omega)$ with the multiplication given by

$$\partial_\Omega^\alpha \partial_\Omega^\beta = \partial_\Omega^{\alpha+\beta},$$

is an abelian semigroup. We will see ahead that $D(\Omega)$, with the operation of composition of functions, is also an abelian semigroup and that the bijection $\mathrm{d} \colon \partial(\Omega) \longrightarrow D(\Omega)$ above is an isomorphism from the semigroup $\partial(\Omega)$ onto the semigroup $D(\Omega)$. In order to do so, let $\alpha, \beta \in \mathbb{N}^n$ be multi-indexes and $\Lambda \in D'(\Omega)$ be arbitrarily fixed. Remembering that $\partial^{\alpha}_{\mathbb{R}^n}$ and $\partial^{\beta}_{\mathbb{R}^n}$ commute in $C^{\infty}(\mathbb{R}^n)$ (see 1.2(j)), let us calculate $D^{\alpha}_{\Omega}(D^{\beta}_{\Omega}(\Lambda))$ applying the Definition 2.33 of derivative of a distribution, to the distribution $D^{\beta}_{\Omega}(\Lambda)$:

$$\begin{aligned}
\left(D^{\alpha}_{\Omega}\left(D^{\beta}_{\Omega}(\Lambda)\right)\right)(\phi) &= (-1)^{|\alpha|}\left(D^{\beta}_{\Omega}(\Lambda)\right)\left(\partial^{\alpha}_{\mathbb{R}^n}(\phi)\right) = \\
&= (-1)^{|\alpha|}(-1)^{|\beta|}\Lambda\left(\partial^{\beta}_{\mathbb{R}^n}\left(\partial^{\alpha}_{\mathbb{R}^n}(\phi)\right)\right) = \\
&= (-1)^{|\alpha+\beta|}\Lambda\left(\partial^{\alpha+\beta}_{\mathbb{R}^n}(\phi)\right) = \\
&= \left(D^{\alpha+\beta}_{\Omega}(\Lambda)\right)(\phi) \quad \text{for every} \quad \phi \in C^{\infty}_{\Omega}(\mathbb{R}^n).
\end{aligned}$$

Therefore,

$$D^{\alpha}_{\Omega}D^{\beta}_{\Omega} = D^{\alpha+\beta}_{\Omega}$$

and hence, since $D^{\alpha+\beta}_{\Omega} \in D(\Omega)$, we conclude that the composition of functions of $D(\Omega)$ ($D^{\alpha}_{\Omega}D^{\beta}_{\Omega}$) is also a function of $D(\Omega)$ ($D^{\alpha}_{\Omega}D^{\beta}_{\Omega} = D^{\alpha+\beta}_{\Omega} \in D(\Omega)$) that, along with the fact of composition being an associative operation, show us that $D(\Omega)$ with the composition is a semigroup. Now, from $D^{\alpha}_{\Omega}D^{\beta}_{\Omega} = D^{\alpha+\beta}_{\Omega}$ and the "tautology" $D^{\alpha+\beta}_{\Omega} = D^{\beta+\alpha}_{\Omega}$, results that

$$D^{\alpha}_{\Omega}D^{\beta}_{\Omega} = D^{\beta}_{\Omega}D^{\alpha}_{\Omega},$$

that is, $D(\Omega)$ is an abelian semigroup. Moreover, from the properties above and the definition of the bijection $\mathrm{d}\colon \partial(\Omega) \longrightarrow D(\Omega)$, we have that

$$\mathrm{d}\left(\partial^{\alpha}_{\Omega}\partial^{\beta}_{\Omega}\right) = \mathrm{d}\left(\partial^{\alpha}_{\Omega}\right)\mathrm{d}\left(\partial^{\beta}_{\Omega}\right)$$

whatever $\partial^{\alpha}_{\Omega}, \partial^{\beta}_{\Omega} \in \partial(\Omega)$, that is, the function

$$\begin{aligned}
\mathrm{d} : \partial(\Omega) &\longrightarrow D(\Omega) \\
\partial^{\alpha}_{\Omega} &\longmapsto \mathrm{d}\left(\partial^{\alpha}_{\Omega}\right) := D^{\alpha}_{\Omega}
\end{aligned}$$

is an isomorphism from the semigroup $\partial(\Omega)$ onto the semigroup $D(\Omega)$.

Another logical consequence of the Definition 2.33 consists in the derivatives of distributions in $\Omega$ being endomorphisms on the abelian group $D'(\Omega)$. In fact, let $D^{\alpha}_{\Omega}\colon D'(\Omega) \longrightarrow D'(\Omega)$ be arbitrarily chosen in the semigroup $D(\Omega)$: we have, for any $\Lambda, \Lambda' \in D'(\Omega)$, that:

$$\begin{aligned}
\left(D^{\alpha}_{\Omega}(\Lambda + \Lambda')\right)(\phi) &= (-1)^{|\alpha|}(\Lambda + \Lambda')\left(\partial^{\alpha}_{\mathbb{R}^n}(\phi)\right) = \\
&= (-1)^{|\alpha|}\Lambda\left(\partial^{\alpha}_{\mathbb{R}^n}(\phi)\right) + (-1)^{|\alpha|}\Lambda'\left(\partial^{\alpha}_{\mathbb{R}^n}(\phi)\right) = \\
&= \left(D^{\alpha}_{\Omega}(\Lambda)\right)(\phi) + \left(D^{\alpha}_{\Omega}(\Lambda')\right)(\phi) = \\
&= \left(D^{\alpha}_{\Omega}(\Lambda) + D^{\alpha}_{\Omega}(\Lambda')\right)(\phi)
\end{aligned}$$

for every $\phi \in C^{\infty}_{\Omega}(\mathbb{R}^n)$, and hence we conclude that

$$D^{\alpha}_{\Omega}(\Lambda + \Lambda') = D^{\alpha}_{\Omega}(\Lambda) + D^{\alpha}_{\Omega}(\Lambda').$$

In short, the set $D(\Omega)$ of derivatives of the distributions in $\Omega$, with the operation of composition of functions, is an abelian semigroup of endomorphisms on the group $D'(\Omega)$ of the distributions in $\Omega$, isomorphic to the semigroup $\partial(\Omega)$ of the classic partial derivatives, and the function

$$\begin{aligned} \mathrm{d} : \partial(\Omega) &\longrightarrow D(\Omega) \\ \partial^{\alpha}_{\Omega} &\longmapsto \mathrm{d}\left(\partial^{\alpha}_{\Omega}\right) \coloneqq D^{\alpha}_{\Omega} \end{aligned}$$

is an isomorphism between the referred semigroups.

Another important aspect regarding our purpose of revealing the relation between the $\widetilde{\mathbb{C}}(\Omega)$-distributions and the L. Schwartz' distributions in $\Omega$, that easily follows from the definition of derivative of distributions in $\Omega$ (which, by the way, was the motivation, presented in 2.32, that led us to the definition given in 2.33) is that

$$D^{\alpha}_{\Omega}(\Lambda_f) = \Lambda_{\partial^{\alpha}_{\Omega}(f)}$$

for every $f \in C^{|\alpha|}(\Omega)$ (the domain of the homomorphism $\partial^{\alpha}_{\Omega} \colon C^{|\alpha|}(\Omega) \longrightarrow C(\Omega)$). In fact, for $\phi \in C^{\infty}_{\Omega}(\mathbb{R}^n)$ arbitrarily chosen,

$$\Big(D^{\alpha}_{\Omega}(\Lambda_f)\Big)(\phi) = (-1)^{|\alpha|}\Lambda_f\Big(\partial^{\alpha}_{\mathbb{R}^n}(\phi)\Big) = (-1)^{|\alpha|}\int_{\Omega} f(x)\Big(\partial^{\alpha}_{\mathbb{R}^n}(\phi)\Big)(x)\,\mathrm{d}x$$

from where, after $|\alpha|$ integrations by parts, we obtain that

$$\Big(D^{\alpha}_{\Omega}(\Lambda_f)\Big)(\phi) = (-1)^{2|\alpha|}\int_{\Omega}\Big(\partial^{\alpha}_{\Omega}(f)\Big)(x)\phi(x)\,\mathrm{d}x = \Lambda_{\partial^{\alpha}_{\Omega}(f)}(\phi),$$

and since $\phi \in C^{\infty}_{\Omega}(\mathbb{R}^n)$ is arbitrary we conclude that

$$D^{\alpha}_{\Omega}(\Lambda_f) = \Lambda_{\partial^{\alpha}_{\Omega}(f)}.$$

Resorting now to the isomorphism $\eta \colon C(\Omega) \longrightarrow C'(\Omega)$ from the group of continuous functions, $C(\Omega)$, onto the group of the distributions in $\Omega$ induced by the functions $f \in C(\Omega)$, $C'(\Omega)$, defined by (see 2.31)

$$\eta(f) = \Lambda_f \quad \text{for every} \quad f \in C(\Omega),$$

the last identity can be written in the following form:

$$D^{\alpha}_{\Omega}\Big(\eta(f)\Big) = \eta\Big(\partial^{\alpha}_{\Omega}(f)\Big)$$

for every $f \in C^{|\alpha|}(\Omega)$, or yet, in an abuse of language (described in the Definition 1.6(b) and, also, in the Remark 1.7)

$$D^{\alpha}_{\Omega}(f) = \partial^{\alpha}_{\Omega}(f) \quad \text{for every} \quad f \in C^{|\alpha|}(\Omega),$$

which informs us that the endomorphism $D^{\alpha}_{\Omega}\colon D'(\Omega) \longrightarrow D'(\Omega)$ is a prolongation, more precisely, an $\eta$-prolongation, of the homomorphism $\partial^{\alpha}_{\Omega}\colon C^{|\alpha|}(\Omega) \longrightarrow C(\Omega)$ to $D'(\Omega)$ (according to Definition 1.6(b)). Thus, and given the fact that the function

$$\begin{aligned} \mathrm{d} : \partial(\Omega) &\longrightarrow D(\Omega) \\ \partial^{\alpha}_{\Omega} &\longmapsto \mathrm{d}\left(\partial^{\alpha}_{\Omega}\right) = D^{\alpha}_{\Omega} \end{aligned}$$

is an isomorphism from the semigroup $\partial(\Omega)$ onto the semigroup $D(\Omega)$, we conclude, endorsed by the Definition 1.6(c), that the semigroup $D(\Omega)$ of the derivatives of the distributions in $\Omega$ is a prolongation, more precisely, an $\eta$-prolongation, to $D'(\Omega)$, of the semigroup $\partial(\Omega)$ of classic derivatives.

## 2.35 The Semigroup $D_f(\Omega)$

Let us now occupy ourselves with the subset of $D'(\Omega)$ constituted by the distributions in $\Omega$ of finite order, that is, the set $D'_f(\Omega) \subseteq D'(\Omega)$ which, as we saw in 2.31, is a subgroup of the abelian group $D'(\Omega)$ of the distributions in $\Omega$. An important property of $D'_f(\Omega)$ is the one of being closed with respect to the operation of differentiation, that is, if $\Lambda \in D'_f(\Omega)$, then, for any $D^{\alpha}_{\Omega} \in D(\Omega)$, $D^{\alpha}_{\Omega}(\Lambda) \in D'_f(\Omega)$. In fact, let $\Lambda \in D'_f(\Omega)$ and $D^{\alpha}_{\Omega} \in D(\Omega)$ be arbitrarily fixed. In this case, since $\Lambda$ is of finite order, there exists an integer $N \geqslant 0$ such that, for each compact $K \subseteq \Omega$ corresponds a constant $C_K < \infty$ and

$$|\Lambda(\phi)| \leqslant C_K \max\left\{\left|\left(\partial^{\alpha'}_{\mathbb{R}^n}(\phi)\right)(x)\right| : x \in \Omega \quad \text{and} \quad |\alpha'| \leqslant N\right\}$$

for every $\phi \in C^{\infty}_K(\mathbb{R}^n)$. But $\partial^{\alpha}_{\mathbb{R}^n}(\phi) \in C^{\infty}_K(\mathbb{R}^n)$ for every $\phi \in C^{\infty}_K(\mathbb{R}^n)$ and hence results from the inequality above that

$$\left|\Lambda\left(\partial^{\alpha}_{\mathbb{R}^n}(\phi)\right)\right| \leqslant C_K \max\left\{\left|\left(\partial^{\alpha'+\alpha}_{\mathbb{R}^n}(\phi)\right)(x)\right| : x \in \Omega \quad \text{and} \quad |\alpha'| \leqslant N\right\},$$

that is,

$$\left|\Lambda\left(\partial^{\alpha}_{\mathbb{R}^n}(\phi)\right)\right| \leqslant C_K \max\left\{\left|\left(\partial^{\beta}_{\mathbb{R}^n}(\phi)\right)(x)\right| : x \in \Omega \quad \text{and} \quad |\beta| \leqslant N + |\alpha|\right\}.$$

On the other hand, from the definition of $D^{\alpha}_{\Omega}(\Lambda)$ we have that

$$\left(D^{\alpha}_{\Omega}(\Lambda)\right)(\phi) = (-1)^{|\alpha|}\Lambda\left(\partial^{\alpha}_{\mathbb{R}^n}(\phi)\right)$$

and hence

$$\left|\left(D^{\alpha}_{\Omega}(\Lambda)\right)(\phi)\right| \leqslant C_K \max\left\{\left|\left(\partial^{\beta}_{\mathbb{R}^n}(\phi)\right)(x)\right| : x \in \Omega \quad \text{and} \quad |\beta| \leqslant N + |\alpha|\right\}$$

for every $\phi \in C_K^\infty(\mathbb{R}^n)$, which allows us to conclude that $D_\Omega^\alpha(\Lambda)$ is a distribution in $\Omega$ of finite order, that is, $D_\Omega^\alpha(\Lambda) \in D'_f(\Omega)$.

Now, knowing that the derivatives of distributions of finite order also are distributions of finite order we get that, for each derivative $D_\Omega^\alpha\colon D'(\Omega) \longrightarrow D'(\Omega)$ in the semigroup $D(\Omega)$, the restriction of $D_\Omega^\alpha$ to the subset $D'_f(\Omega)$ of its domain $D'(\Omega)$, $D_\Omega^\alpha|_{D'_f(\Omega)}$, is a function assuming values in $D'_f(\Omega)$. It is convenient to introduce a notation for these restrictions: we will use the symbol $D_{\Omega-f}^\alpha$ to denote the restriction $D_\Omega^\alpha|_{D'_f(\Omega)}$, that is, $D_{\Omega-f}^\alpha$ is the following function:

$$\begin{aligned} D_{\Omega-f}^\alpha : D'_f(\Omega) &\longrightarrow D'_f(\Omega) \\ \Lambda &\longmapsto D_{\Omega-f}^\alpha(\Lambda) := D_\Omega^\alpha(\Lambda). \end{aligned}$$

The set of the restrictions of the derivatives $D_\Omega^\alpha \in D(\Omega)$ to $D'_f(\Omega)$, that is, the set of the functions $D_{\Omega-f}^\alpha$ (one for each multi-index $\alpha \in \mathbb{N}^n$) will be denoted by $D_f(\Omega)$. Hence,

$$D_f(\Omega) := \left\{ D_{\Omega-f}^\alpha := D_\Omega^\alpha|_{D'_f(\Omega)} : D_\Omega^\alpha \in D(\Omega) \right\}$$

or, more briefly,

$$D_f(\Omega) = \left\{ D_{\Omega-f}^\alpha : \alpha \in \mathbb{N}^n \right\},$$

and we will refer to the members $D_{\Omega-f}^\alpha \in D_f(\Omega)$ as finite order distributions derivatives.

It results immediately from the properties of the semigroup $D(\Omega)$ and of its members, the derivatives $D_\Omega^\alpha$, obtained in item 2.34, along with the above definitions of the functions $D_{\Omega-f}^\alpha$ and its set $D_f(\Omega)$, that:

**(a)** $D_{\Omega-f}^\alpha \in D_f(\Omega)$ is an endomorphism on the group $D'_f(\Omega)$ of the distributions (in $\Omega$) of finite order;

**(b)** $D_{\Omega-f}^\alpha \in D_f(\Omega)$ is an $\eta$-prolongation to $D'_f(\Omega)$ of the homomorphism $\partial_\Omega^\alpha \in \partial(\Omega)$, that is, $D_{\Omega-f}^\alpha$ is an endomorphism on $D'_f(\Omega)$ such that

$$D_{\Omega-f}^\alpha\Big(\eta(f)\Big) = \eta\Big(\partial_\Omega^\alpha(f)\Big),$$

that is,

$$D_{\Omega-f}^\alpha\left(\Lambda_f\right) = \Lambda_{\partial_\Omega^\alpha(f)}$$

for every $f \in C^{|\alpha|}(\Omega)$;

**(c)** $D_f(\Omega)$, with the composition of functions, is an abelian semigroup, i.e.,

$$D_{\Omega-f}^\alpha D_{\Omega-f}^\beta = D_{\Omega-f}^{\alpha+\beta} = D_{\Omega-f}^{\beta+\alpha} = D_{\Omega-f}^\beta D_{\Omega-f}^\alpha$$

for every $D_{\Omega-f}^\alpha, D_{\Omega-f}^\beta \in D_f(\Omega)$;

**(d)** the semigroup $D_f(\Omega)$ is isomorphic to the semigroup $\partial(\Omega)$ and the function

$$\begin{aligned} \mathrm{d}_f : \partial(\Omega) &\longrightarrow D_f(\Omega) \\ \partial_\Omega^\alpha &\longmapsto \mathrm{d}_f(\partial_\Omega^\alpha) \coloneqq D_{\Omega-f}^\alpha \end{aligned}$$

is an isomorphism between the referred semigroups;

**(e)** the isomorphism $\mathrm{d}_f$ of item (d) is such that

$$\Big(\mathrm{d}_f\,(\partial_\Omega^\alpha)\Big)\Big(\eta(f)\Big) = \eta\Big(\partial_\Omega^\alpha(f)\Big),$$

that is,

$$D_{\Omega-f}^\alpha\,(\Lambda_f) = \Lambda_{\partial_\Omega^\alpha(f)}$$

for every $f \in C^{|\alpha|}(\Omega)$, and thus, according to Definition 1.6(c), the semigroup $D_f(\Omega)$ is a prolongation (more precisely, an $\eta$-prolongation) to $D'_f(\Omega)$ of the semigroup $\partial(\Omega)$;

**(f)** such as $D(\Omega)$, which contains an identity, that is, $I_{D'(\Omega)} \in D(\Omega)$ (see Definition 2.33), the semigroup $D_f(\Omega)$ also contains an identity, that is, $I_{D'_f(\Omega)} \in D_f(\Omega)$.

## 2.36 The Extensions $\mathbb{D}'(\Omega)$ and $\mathbb{D}'_f(\Omega)$ of the $S$-group $\mathbb{C}(\Omega)$

Regarding the $S$-group of the continuous functions

$$\mathbb{C}(\Omega) = \left(C(\Omega), \partial(\Omega) = \Big\{\partial_\Omega^\alpha \colon \alpha \in \mathbb{N}^n\Big\}\right),$$

the results so far obtained about distributions in $\Omega$ allow us to affirm that:

**(a)** the abelian group $D'(\Omega)$ of the distributions in $\Omega$, as well as its subgroup $D'_f(\Omega)$ of the distributions in $\Omega$ of finite order, both $\eta$-admits the abelian group of continuous functions, $C(\Omega)$, as a subgroup, being $\eta \colon C(\Omega) \longrightarrow C'(\Omega)$ the isomorphism given by

$$\eta(f) = \Lambda_f \quad \text{for every} \quad f \in C(\Omega);$$

**(b)** the semigroup $D_f(\Omega)$ is an $\eta$-prolongation to $D'_f(\Omega)$ of the semigroup $\partial(\Omega)$;

**(c)** the semigroup $D(\Omega)$ is an $\eta$-prolongation to $D'(\Omega)$ of the semigroup $\partial(\Omega)$;

**(d)** the semigroup $D(\Omega)$ is isomorphic to the semigroup $D_f(\Omega)$, the function

$$\begin{aligned} \mathrm{e} : D_f(\Omega) &\longrightarrow D(\Omega) \\ D_{\Omega-f}^\alpha &\longmapsto \mathrm{e}\big(D_{\Omega-f}^\alpha\big) \coloneqq D_\Omega^\alpha \end{aligned}$$

is an isomorphism between the respective semigroups and, hence, taking into account the Definition 1.6(c), $D(\Omega)$ is a prolongation to $D'(\Omega)$ of the semigroup $D_f(\Omega)$.

Defining now

$$\mathbb{D}'(\Omega) := \Big(D'(\Omega), D(\Omega)\Big)$$

and

$$\mathbb{D}'_f(\Omega) := \Big(D'_f(\Omega), D_f(\Omega)\Big)$$

and with the definitions of $S$-group (Definition 1.4(a)) and $S$-group extension (Definition 1.6(d)) in mind, we conclude that:

**(e)** from (a) and (b), and (a) and (c), respectively, that $\mathbb{D}'_f(\Omega)$ and $\mathbb{D}'(\Omega)$ are extensions (more precisely, $\eta$-extensions) of the $S$-group $\mathbb{C}(\Omega)$;

**(f)** from (d) and that $D'_f(\Omega) \subseteq D'(\Omega)$ is a subgroup of $D'(\Omega)$, result that $\mathbb{D}'(\Omega)$ is an extension of the $S$-group $\mathbb{D}'_f(\Omega)$.

Among these two extensions, $\mathbb{D}'(\Omega)$ and $\mathbb{D}'_f(\Omega)$, of the $S$-group $\mathbb{C}(\Omega)$, does any of them is strict, closed or even strict and closed? The next item answers this question.

## 2.37 $\mathbb{D}'_f(\Omega)$, a Strict and Closed Extension of $\mathbb{C}(\Omega)$

Taking advantage of the abuse of language described in Definition 1.6(b) (see also the Remark 1.7), we get that $\mathbb{D}'(\Omega)$ is a strict extension of the $S$-group $\mathbb{C}(\Omega)$ if and only if, for every $\partial^\alpha_\Omega \in \partial(\Omega)$,

$$\text{if} \quad f \in C(\Omega) \quad \text{is such that} \quad f \notin C^{|\alpha|}(\Omega), \quad \text{then} \quad D^\alpha_\Omega(f) \notin C(\Omega). \tag{2.37-1}$$

In the same way, $\mathbb{D}'_f(\Omega)$ is a strict extension of $\mathbb{C}(\Omega)$ if and only if, for each $\partial^\alpha_\Omega \in \partial(\Omega)$,

$$\text{if} \quad f \in C(\Omega) \quad \text{is such that} \quad f \notin C^{|\alpha|}(\Omega), \quad \text{then} \quad D^\alpha_{\Omega-f}(f) \notin C(\Omega). \tag{2.37-2}$$

But, since every function $f \in C(\Omega)$ is a finite order distribution (more precisely, it can be identified with such a distribution, $f \equiv \Lambda_f \in D'_f(\Omega)$), we have that $D^\alpha_\Omega(f) = D^\alpha_{\Omega-f}(f)$ if $f \in C(\Omega)$ and hence the conditions (2.37-1) and (2.37-2) to be satisfied by $\mathbb{D}'(\Omega)$ and $\mathbb{D}'_f(\Omega)$, respectively, for them to be strict extensions of the $S$-group $\mathbb{C}(\Omega)$, are equal. Therefore, either $\mathbb{D}'(\Omega)$ and $\mathbb{D}'_f(\Omega)$ are, both, strict extensions of $\mathbb{C}(\Omega)$ or none of them is.

Now, without the abuse of language used, let us say, in (2.37-1), the referred condition assumes the following precise form: for each $\partial^\alpha_\Omega \in \partial(\Omega)$, if $\Lambda \in \eta(C(\Omega)) = C'(\Omega)$ (that is, if $\Lambda = \Lambda_f$ for some $f \in C(\Omega)$) and $\Lambda = \Lambda_f \notin \eta\left(C^{|\alpha|}(\Omega)\right)$, then, $D^\alpha_\Omega(\Lambda_f) \notin \eta(C(\Omega)) = C'(\Omega)$.

Below we demonstrate this last implication. For such, let $\Lambda \in C'(\Omega)$ be arbitrarily chosen in such a way that

$$\Lambda \notin \eta\Big(C^{|\alpha|}(\Omega)\Big).$$

In this case, there exists $f \in C(\Omega)$ such that $f \notin C^{|\alpha|}(\Omega)$ and $\Lambda = \Lambda_f$. Suppose now (*reductio ad absurdum*) that

$$D^{\alpha}_{\Omega}(\Lambda_f) \in C'(\Omega).$$

Hence, since $\eta\colon C(\Omega) \longrightarrow C'(\Omega)$ is an isomorphism, there exists $g \in C(\Omega)$ such that

$$\eta(g) = \Lambda_g = D^{\alpha}_{\Omega}(\Lambda_f).$$

Therefore, for every $\phi \in C^{\infty}_{\Omega}(\mathbb{R}^n)$ we have

$$\Big(D^{\alpha}_{\Omega}(\Lambda_f)\Big)(\phi) = \Lambda_g(\phi) = \int_{\Omega} g(x)\phi(x)\,\mathrm{d}x\,.$$

On the other hand, from the definition of derivative of distribution,

$$\Big(D^{\alpha}_{\Omega}(\Lambda_f)\Big)(\phi) = (-1)^{|\alpha|}\Lambda_f\Big(\partial^{\alpha}_{\mathbb{R}^n}(\phi)\Big) = (-1)^{|\alpha|}\int_{\Omega} f(x)\Big(\partial^{\alpha}_{\mathbb{R}^n}(\phi)\Big)(x)\,\mathrm{d}x$$

for every $\phi \in C^{\infty}_{\Omega}(\mathbb{R}^n)$. Thus, from the last two equations we obtain:

$$\int_{\Omega}\Big[(-1)^{|\alpha|}f(x)\Big(\partial^{\alpha}_{\mathbb{R}^n}(\phi)\Big)(x) - g(x)\phi(x)\Big]\,\mathrm{d}x = 0$$

for every $\phi \in C^{\infty}_{\Omega}(\mathbb{R}^n)$.

Now, $\partial^{\alpha}_{\Omega}\colon C^{|\alpha|}(\Omega) \longrightarrow C(\Omega)$ is a surjective homomorphism and hence, since $g \in C(\Omega)$, there exists $h \in C^{|\alpha|}(\Omega)$ such that $\partial^{\alpha}_{\Omega}(h) = g$ and, from this, the last equation then assumes the form

$$\int_{\Omega}\Big[(-1)^{|\alpha|}f(x)\Big(\partial^{\alpha}_{\mathbb{R}^n}(\phi)\Big)(x) - \Big(\partial^{\alpha}_{\Omega}(h)\Big)(x)\phi(x)\Big]\,\mathrm{d}x = 0$$

for every $\phi \in C^{\infty}_{\Omega}(\mathbb{R}^n)$. However, $|\alpha|$ integrations by parts show that

$$\int_{\Omega}\Big(\partial^{\alpha}_{\Omega}(h)\Big)(x)\phi(x)\,\mathrm{d}x = (-1)^{|\alpha|}\int_{\Omega} h(x)\Big(\partial^{\alpha}_{\mathbb{R}^n}(\phi)\Big)(x)\,\mathrm{d}x\,,$$

a result that, when taken into the previous equality, leads to

$$\int_{\Omega}\Big(f(x) - h(x)\Big)\Big(\partial^{\alpha}_{\mathbb{R}^n}(\phi)\Big)(x)\,\mathrm{d}x = 0$$

for every $\phi \in C^{\infty}_{\Omega}(\mathbb{R}^n)$, from where we conclude that

$$f = h,$$

which is contrary to our initial hypothesis that $f \notin C^{|\alpha|}(\Omega)$, since $h \in C^{|\alpha|}(\Omega)$.

We know now that the $S$-groups $\mathbb{D}'(\Omega)$ and $\mathbb{D}'_f(\Omega)$, henceforth denominated, respectively, $S$-group of the distributions (in $\Omega$) and $S$-group of the distributions (in $\Omega$) of finite order, are strict extensions of the continuous functions $S$-group, $\mathbb{C}(\Omega)$. It can also be proved that the $S$-group of the distributions of finite order is a closed extension of $\mathbb{C}(\Omega)$, that is: every distribution in $\Omega$ of finite order is a derivative of some continuous function or, more precisely, without abuses of language, that, if $\Lambda \in D'_f(\Omega)$, then, there exists $f \in C(\Omega)$ and a multi-index $\alpha \in \mathbb{N}^n$ such that

$$\Lambda = D^{\alpha}_{\Omega - f}(\Lambda_f)$$

or, equivalently, there exists $f \in C(\Omega)$ and $\partial^{\alpha}_{\Omega} \in \partial(\Omega)$ such that

$$\Lambda = \Big(\mathrm{d}_f(\partial^{\alpha}_{\Omega})\Big)\Big(\eta(f)\Big)$$

($\mathrm{d}_f \colon \partial(\Omega) \longrightarrow D_f(\Omega)$ is the isomorphism defined in 2.35(d)).

The infinite order distributions can not be expressed as derivatives of continuous functions, since the continuous functions are (identified to) finite order distributions and its derivatives also are finite order distributions. Therefore, if there exist infinite order distributions, that is, if $D'_f(\Omega)$ is a proper subset of $D'(\Omega)$, then, the $S$-group $\mathbb{D}'(\Omega)$ of the distributions in $\Omega$ is not a closed extension of the $S$-group $\mathbb{C}(\Omega)$; the infinite order distributions do exist, an example is given by the functional

$$\begin{aligned} \Lambda : C^{\infty}_{(0,\infty)}(\mathbb{R}) &\longrightarrow C \\ \phi &\longmapsto \Lambda(\phi) := \textstyle\sum_{m=1}^{\infty} \left(\frac{\mathrm{d}^m \phi}{\mathrm{d}x^m}\right)\left(\frac{1}{m}\right) \end{aligned}$$

which the reader can, as an exercise, verify.

In short, we have that the $S$-group of the distributions in $\Omega$ of finite order,

$$\mathbb{D}'_f(\Omega) = \Big(D'_f(\Omega), D_f(\Omega)\Big),$$

is a strict and closed extension of the $S$-group of continuous functions,

$$\mathbb{C}(\Omega) = \Big(C(\Omega), \partial(\Omega)\Big),$$

while the $S$-group of distributions in $\Omega$,

$$\mathbb{D}'(\Omega) = \Big(D'(\Omega), D(\Omega)\Big),$$

is a strict (but not closed) extension of $\mathbb{C}(\Omega)$ and, also, an extension of $\mathbb{D}'_f(\Omega)$.

## 2.38 The $\widetilde{\mathbb{C}}(\Omega)$-Distributions and the Distributions in $\Omega$ of Finite Order

According to the Definition 2.18, a domain of distribution is an ordered pair of $S$-groups, $(\mathbb{G} = (G, H), \widehat{\mathbb{G}} = (\widehat{G}, \widehat{H}))$, where $\mathbb{G}$ is an abelian, surjective, and with identity $S$-group and $\widehat{\mathbb{G}}$ is a strict and closed extension of $\mathbb{G}$. Hence, given an abelian, surjective, and with identity $S$-group, $\mathbb{G}$, the Theorem of Extension of $S$-Groups (Theorem 2.16) ensures the existence of a unique (unless isomorphism) $S$-group $\widehat{\mathbb{G}}$ that is a strict and closed extension of $\mathbb{G}$; to each $S$-group $\mathbb{G}$ as above, corresponds, then, an essentially unique domain of distribution $(\mathbb{G}, \widehat{\mathbb{G}})$. The $\widetilde{\mathbb{G}}$-distributions are, then, according to the referred Definition 2.18, the members of the group $\widehat{G}$ of the $S$-group $\widehat{\mathbb{G}} = (\widehat{G}, \widehat{H})$, while the derivatives of the $\widetilde{\mathbb{G}}$-distributions are the elements of the semigroup $\widehat{H}$ of the $S$-group $\widehat{\mathbb{G}}$. Thus, taking the group $\mathbb{G}$ as the continuous functions one (which is abelian, surjective, and with identity), that is, for

$$\mathbb{G} = \mathbb{C}(\Omega) = \Big(C(\Omega), \partial(\Omega)\Big),$$

and taking into account that the $S$-group

$$\mathbb{D}'_f(\Omega) = \Big(D'_f(\Omega), D_f(\Omega)\Big),$$

is a strict and closed extension of $\mathbb{C}(\Omega)$, we get that the ordered pair

$$\left(\mathbb{C}(\Omega) = \Big(C(\Omega), \partial(\Omega)\Big), \mathbb{D}'_f(\Omega) = \Big(D'_f(\Omega), D_f(\Omega)\Big)\right)$$

is a domain of distribution — the domain of the $\widetilde{\mathbb{C}}(\Omega)$-distributions — at which the $\widetilde{\mathbb{C}}(\Omega)$-distributions are the L. Schwartz' finite order distributions in $\Omega$, the elements of the group $D'_f(\Omega)$, and the derivatives of the $\widetilde{\mathbb{C}}(\Omega)$-distributions are the elements of the semigroup $D_f(\Omega)$, that is, the derivatives of the finite order distributions as defined by the eminent French mathematician.

On the other hand, in items 2.20 to 2.26, we formulated two axiomatics, the $\widetilde{\mathbb{G}}$-distributions one (in item 2.20) and its simplified version (in item 2.25), and we proved that these axiomatics are equivalent to each other, categoric and both admit as a model, exactly, the $\widetilde{\mathbb{G}}$-distributions and its derivatives. Therefore, these axiomatics translated to the case where $\mathbb{G} = \mathbb{C}(\Omega)$, that is, to the domain of distributions $(\mathbb{C}(\Omega), \mathbb{D}'_f(\Omega))$, constitute, each one of them, a categoric axiomatic definition of the finite order distributions in $\Omega$ and its derivatives. Next, in item 2.39, we present this translation to the domain $(\mathbb{C}(\Omega), \mathbb{D}'_f(\Omega))$ of the referred axiomatics, with a natural substitution of the primitive terms "$\widetilde{\mathbb{C}}(\Omega)$-distribution" and "derivative of $\widetilde{\mathbb{C}}(\Omega)$-distribution" by "distribution in $\Omega$ of finite order" and "derivative of distribution in $\Omega$ of finite order", respectively.

## 2.39 Finite Order Distributions Axiomatic

The axioms of the two axiomatics formulated ahead are stated relative to the $S$-group of continuous functions, $\mathbb{C}(\Omega) = (C(\Omega), \partial(\Omega) = \{\partial_\Omega^\alpha \colon \alpha \in \mathbb{N}^n\})$, and define implicitly the following primitive terms:

- distribution in $\Omega$ of finite order, addition and derivative (of distributions of finite order), for the first axiomatic;
- distribution in $\Omega$ of finite order and derivative (of distributions of finite order), for the second axiomatic.

The precedent theories for each of the referred axiomatics, as was established in 2.20 and 2.25, are Classical Logic, Set Theory, and $S$-Groups Theory.

According to the considerations in 2.38, these are categoric and equivalent axiomatics, with model, essentially unique, constituted of the distributions (in $\Omega$) of finite order and its derivatives as formulated by L. Schwartz.

To allow a more concise formulation of the axioms, we will employ the following primitive symbols which, surely, is the ones a watchful reader would choose:

$D'_f(\Omega)$, for the class of the distributions in $\Omega$ of finite order;

$D_f(\Omega)$, for the class of the derivatives of distributions in $\Omega$ of finite order;

$+$, for the addition in $D'_f(\Omega)$.

We pass now to the presentation of the axioms.

- **Finite order distributions axiomatic:**

  **Axiom 1** $C(\Omega) \subseteq D'_f(\Omega)$.

  **Axiom 2** The addition, $+$, is a binary operation in $D'_f(\Omega)$ whose restriction to $C(\Omega) \times C(\Omega)$ is the addition of the group $C(\Omega)$; for $\Lambda, \Lambda' \in D'_f(\Omega)$, the corresponding element in $D'_f(\Omega)$ associated to the pair $(\Lambda, \Lambda')$ by $+$, is denoted by $\Lambda + \Lambda'$ and denominated the sum of $\Lambda$ with $\Lambda'$.

  **Axiom 3** The derivatives (of distributions of finite order), that is, the members of the class $D_f(\Omega)$, are functions with domain and codomain equal to $D'_f(\Omega)$, such that:

  **(a)** the derivative of the sum is the sum of the derivatives;

  **(b)** each element of $D_f(\Omega)$ is an extension to $D'_f(\Omega)$ of a single $\partial^\alpha_\Omega \in \partial(\Omega)$, that is, if $\mathrm{d} \in D_f(\Omega)$, then, there exists an unique $\partial^\alpha_\Omega \in \partial(\Omega)$ such that

  $$\mathrm{d}(f) = \partial^\alpha_\Omega(f)$$

  for every $f \in C^{|\alpha|}(\Omega)$ (the domain of $\partial^\alpha_\Omega$);

**(c)** the class of the derivatives, $D_f(\Omega)$, with the usual operation of composition of functions, is a semigroup isomorphic to the semigroup $\partial(\Omega)$ and the function $\mathrm{d}_f$ defined by

$$\begin{aligned} \mathrm{d}_f : \partial(\Omega) &\longrightarrow D_f(\Omega) \\ \partial^\alpha_\Omega &\longmapsto \mathrm{d}_f(\partial^\alpha_\Omega) \coloneqq D^\alpha_{\Omega-f}, \end{aligned}$$

where $D^\alpha_{\Omega-f}$ is the extension of $\partial^\alpha_\Omega$ to $D'_f(\Omega)$, is an isomorphism (between semigroups).

**Axiom 4** For each $\partial^\alpha_\Omega \in \partial(\Omega)$, if $f \in C(\Omega)$ and $f \notin C^{|\alpha|}(\Omega)$ (the domain of $\partial^\alpha_\Omega$), then,

$$\Big(\mathrm{d}_f(\partial^\alpha_\Omega)\Big)(f) = D^\alpha_{\Omega-f}(f) \notin C(\Omega).$$

**Axiom 5** For each $\Lambda \in D'_f(\Omega)$, there exist $\partial^\alpha_\Omega \in \partial(\Omega)$ and $f \in C(\Omega)$ such that

$$\Lambda = \Big(\mathrm{d}_f(\partial^\alpha_\Omega)\Big)(f) = D^\alpha_{\Omega-f}(f).$$

- **Finite order distributions axiomatic — simplified version:**

**Axiom 1** $C(\Omega) \subseteq D'_f(\Omega)$.

**Axiom 2** The derivatives, that is, the members of the class $D_f(\Omega)$, are functions with domain and codomain equal to $D'_f(\Omega)$, such that:

**(a)** each element of $D_f(\Omega)$ is an extension to $D'_f(\Omega)$ of a single $\partial^\alpha_\Omega \in \partial(\Omega)$, that is, if $\mathrm{d} \in D_f(\Omega)$, then, there exists an unique $\partial^\alpha_\Omega \in \partial(\Omega)$ such that

$$\mathrm{d}(f) = \partial^\alpha_\Omega(f)$$

for every $f \in C^{|\alpha|}(\Omega)$ (the domain of $\partial^\alpha_\Omega$);

**(b)** the class of the derivatives, $D_f(\Omega)$, with the usual operation of composition of functions, is a semigroup isomorphic to the semigroup $\partial(\Omega)$ and the function $\mathrm{d}_f$ defined by

$$\begin{aligned} \mathrm{d}_f : \partial(\Omega) &\longrightarrow D_f(\Omega) \\ \partial^\alpha_\Omega &\longmapsto \mathrm{d}_f(\partial^\alpha_\Omega) \coloneqq D^\alpha_{\Omega-f}, \end{aligned}$$

where $D^\alpha_{\Omega-f}$ is the extension of $\partial^\alpha_\Omega$ to $D'_f(\Omega)$, is an isomorphism (between semigroups).

**Axiom 3** For each $\Lambda \in D'_f(\Omega)$, there exist $\partial^\alpha_\Omega \in \partial(\Omega)$ and $f \in C(\Omega)$ such that

$$\Lambda = \Big(\mathrm{d}_f(\partial^\alpha_\Omega)\Big)(f) = D^\alpha_{\Omega-f}(f).$$

**Axiom 4** If $f, g \in C(\Omega)$ and $\partial^\alpha_\Omega \in \partial(\Omega)$, then, $(\mathrm{d}_f(\partial^\alpha_\Omega))(f) = (\mathrm{d}_f(\partial^\alpha_\Omega))(g)$ if and only if $f - g \in N(\partial^\alpha_\Omega)$, that is, $D^\alpha_{\Omega-f}(f) = D^\alpha_{\Omega-f}(g)$ if and only if $f - g \in N(\partial^\alpha_\Omega)$.

In this point of our considerations about the notions of distributions and its derivatives, it is opportune to highlight an aspect of the definitions provided for these concepts through the categoric (and equivalent) axiomatics above, which we consider of the highest importance, namely, the simplicity of the axioms; resorting to very simple mathematical concepts, such as group, semigroup, the classical notion of derivative, among others, it is defined, through a set of only four axioms (of the "simplified version") formulated from these elementary notions, categorically, the concepts of finite order distribution and derivative of a finite order distribution taken as primitives, definitions which coincide (in the sense of being isomorphic to) with the corresponding constructions of these concepts, in the heart of mathematical analysis, provided by L. Schwartz which, in turn, requires a mathematical basis that, compared with that required in the elaboration of the referred axiomatics, can be classified as highly sophisticated. To the interested reader, for a better appreciation of the simplicity above highlighted, we recommend to "browse" the chapter destined to the theory of distributions in the text by W. Rudin referred to in page 69 (footnote 14).

Sir Winston Churchill, in a speech to the British (see "Blood, Sweat and Tears" of Winston S. Churchill) by occasion of the Second World War, referred to the debt of his people to the members of their air force saying that:

> *"Never was so much owed by so many to so few."*

Maybe we should, taking into account that the axiomatic method and mathematics as we understand it today are creations of thinkers from ancient Greece, to say, following the example of the great British statesman, that: never was so much owned by so many (actually the whole humanity) to so few philosophers from ancient Greece.

## 2.40 Remark

The concepts of distribution and derivative of a distribution seen as the integrating elements of the strict and closed extension of the $S$-group $\mathbb{C}(\Omega) = (C(\Omega), \partial(\Omega))$ of continuous functions, endorsed by the Theorem of Extension of $S$-Groups (Theorem 2.16), do not include, as seen in this chapter, the infinite order distributions and its derivatives. On the other hand, the extensions promoted by the L. Schwartz' distributions theory to the concepts of continuous functions and differentiability, not only satisfy the conditions (a) to (d) listed in 1.3, attended by the strict and closed extension of the $S$-group $\mathbb{C}(\Omega)$, that any such generalization should attend, but it is also a wider generalization of that offered by the Theorem of Extension of $S$-Groups, insofar as it houses the infinite order distributions which, in turn, does not figure in the strict and closed extension of $\mathbb{C}(\Omega)$.

At this point, naturally, the following question imposes itself: Something analogous to what we accomplished for the finite order distributions can be obtained for the class $D'(\Omega)$ of all the distributions? More specifically, we inquire: Through a simple categoric axiomatic, formulated with elementary concepts, taking the ones used in the formulation of the finite order distributions axiomatics as examples, would not be possible to capture that wider extension promoted to the concepts of continuous function and differentiability by L. Schwartz' distributions theory?

Clearly, if we pretend to obtain one such axiomatization starting with those described in 2.39 for the finite order distributions, we must, necessarily, modify the Axiom 5 (respectively, the Axiom 3 of the simplified version) which requires every distribution in $\Omega$ to be the derivative of some function in $C(\Omega)$. Maybe we could take a "local" version of the referred axiom, requiring every distribution to be at least "locally" the derivative of a continuous function, that, in turn, would require other primitive notions beyond those of distribution, addition and derivative (respectively, of distribution and derivative in the simplified version), for example, the ones of domain and restriction of a distribution (in some sensible sense), through which the notion of local equality could be defined.

The next chapters of this monograph are dedicated to answer the question above; in Chapter 3 we define a structure — $S$-space — that plays an analogous role to that of a $S$-group concerning our axiomatic purposes and, next, in the fourth chapter, we prove two theorems, the (1st and 2nd) Theorems of Extension of $S$-Spaces, that, such as the Theorem of Extension of $S$-Groups which sustains the finite order distributions axiomatics, provide the base on which lean the axiomatics formulated in the fifth and last chapter, which meet the desired conditions, thus answering, positively, the question posed above.

# Appendix: Proof of Lemma 2.4

**Lemma 2.4.** *Let $\mathbb{G} = (G, H)$ be an abelian, surjective, and with identity $S$-group, $g, h \in G$, $\Phi, \Psi, \varphi \in H$ and $F, K, L, M \subseteq G$. Then:*

**(a)** *if $F \subseteq K$ and $L \subseteq M$, then $F \pm L \subseteq K \pm M$; in particular, since $M \subseteq M$, then $F \pm M \subseteq K \pm M$;*

**(b)** *if $F \subseteq K$, then $\Phi^{-1}(F) \subseteq \Phi^{-1}(K)$;*

**(c)** $\Phi^{-1}\Big(\Psi^{-1}(g)\Big) = \Psi^{-1}\Big(\Phi^{-1}(g)\Big)$ *and* $\Big(\Phi\Psi\Big)^{-1}(g) = \Psi^{-1}\Big(\Phi^{-1}(g)\Big)$;

**(d)** $\varphi^{-1}\Big(\Psi^{-1}(g)\Big) \pm \varphi^{-1}\Big(\Phi^{-1}(h)\Big) \subseteq \varphi^{-1}\Big(\Psi^{-1}(g) \pm \Phi^{-1}(h)\Big)$;

**(e)** $\varphi^{-1}\Big(\Psi^{-1}(g) \pm \Phi^{-1}(h)\Big) \subseteq \varphi^{-1}\Big(\Psi^{-1}(g)\Big) \pm \varphi^{-1}\Big(\Phi^{-1}(h)\Big) + N(\varphi)$;

**(f)** $N\Big(\Phi\Psi\Big) = \Phi^{-1}\Big(N(\Psi)\Big) + \Psi^{-1}\Big(N(\Phi)\Big) + N(\Phi) + N(\Psi)$;

**(g)** $N(\Psi\varphi\Phi) + N(\Psi) + N(\Phi) \subseteq N(\Psi\varphi\Phi)$;

**(h)** *if $h = \Phi(g)$ and $g \in N(\Phi\Psi)$, then $h \in N(\Psi)$.*

*Proof.*

**(a)** Let $x \in F \pm L$. In this case, $x = x_1 \pm x_2$, where $x_1 \in F$ and $x_2 \in L$. Thus, if $F \subseteq K$ and $L \subseteq M$, then $x_1 \in K$ and $x_2 \in M$ and, therefore, the element $x = x_1 \pm x_2$ is a member of $K \pm M$;

**(b)** Let $x \in \Phi^{-1}(F)$. Hence, $\Phi(x) \in F \subseteq K$ and, therefore, $\Phi(x) \in K$, that is, $x \in \Phi^{-1}(K)$;

**(c)** Let $x \in \Phi^{-1}(\Psi^{-1}(g))$. In this case, $\Phi(x) \in \Psi^{-1}(g)$ and hence $\Psi(\Phi(x)) = g$. Since $\Phi\Psi = \Psi\Phi$, comes that $\Phi(\Psi(x)) = g$ and hence $\Psi(x) \in \Phi^{-1}(g)$ and $x \in \Psi^{-1}(\Phi^{-1}(g))$. So we conclude that $\Phi^{-1}(\Psi^{-1}(g)) \subseteq \Psi^{-1}(\Phi^{-1}(g))$. Analogously one can prove that $\Psi^{-1}(\Phi^{-1}(g)) \subseteq \Phi^{-1}(\Psi^{-1}(g))$. It is even easier to prove that $(\Phi\Psi)^{-1}(g) = \Psi^{-1}(\Phi^{-1}(g))$;

**(d)** Let $x \in \varphi^{-1}(\Psi^{-1}(g)) \pm \varphi^{-1}(\Phi^{-1}(h))$. In this case: $x = x_g \pm x_h$ where $x_g \in \varphi^{-1}(\Psi^{-1}(g))$ and $x_h \in \varphi^{-1}(\Phi^{-1}(h))$. Thus, $\varphi(x_g) \in \Psi^{-1}(g)$ and $\varphi(x_h) \in \Phi^{-1}(h)$ and hence, since $\varphi(x) = \varphi(x_g \pm x_h) = \varphi(x_g) \pm \varphi(x_h)$, then comes that $x \in \varphi^{-1}(\Psi^{-1}(g) \pm \Phi^{-1}(h))$;

**(e)** Let $x \in \varphi^{-1}(\Psi^{-1}(g) \pm \Phi^{-1}(h))$. Hence, we have that

$$\varphi(x) \in \Psi^{-1}(g) \pm \Phi^{-1}(h)$$

and, therefore,

$$\varphi(x) = g^* \pm h^*$$

where $g^* \in \Psi^{-1}(g)$ and $h^* \in \Phi^{-1}(h)$, that is, $g^*$ and $h^*$ are such that

$$\Psi(g^*) = g \quad \text{and} \quad \Phi(h^*) = h.$$

Since $\varphi$ is surjective, then there exist $x_g, x_h \in G_\varphi$ such that

$$\varphi(x_g) = g^* \quad \text{and} \quad \varphi(x_h) = h^*.$$

Now, we define:

$$x' := x_g \pm x_h.$$

Observe that $x_g \in \varphi^{-1}(g^*)$ and $x_h \in \varphi^{-1}(h^*)$, and yet, that $g^* \in \Psi^{-1}(g)$ and $h^* \in \Phi^{-1}(h)$. Thus,

$$x_g \in \varphi^{-1}\Big(\Psi^{-1}(g)\Big) \quad \text{and} \quad x_h \in \varphi^{-1}\Big(\Phi^{-1}(h)\Big).$$

Beyond that, we have:

$$\varphi(x') = \varphi(x_g \pm x_h) = \varphi(x_g) \pm \varphi(x_h) = g^* \pm h^* = \varphi(x),$$

that is,

$$\varphi(x - x') = 0$$

and hence $x - x' = n$, where $n \in N(\varphi)$. Therefore,

$$x = x' + n = x_g \pm x_h + n,$$

where $x_g \in \varphi^{-1}\Big(\Psi^{-1}(g)\Big)$, $x_h \in \varphi^{-1}\Big(\Phi^{-1}(h)\Big)$ and $n \in N(\varphi)$, that is,

$$x \in \varphi^{-1}\Big(\Psi^{-1}(g)\Big) \pm \varphi^{-1}\Big(\Phi^{-1}(h)\Big) + N(\varphi);$$

**(f)** Let $x \in N(\Phi\Psi)$. In this case, we have that:

$$\Phi\Big(\Psi(x)\Big) = \Psi\Big(\Phi(x)\Big) = 0$$

and hence, $\Psi(x) \in N(\Phi)$ and $\Phi(x) \in N(\Psi)$, which allows us to conclude that

$$x \in \Psi^{-1}\Big(N(\Phi)\Big) \quad \text{and} \quad x \in \Phi^{-1}\Big(N(\Psi)\Big).$$

Analogously we obtain that

$$x + x \in \Phi^{-1}\Big(N(\Psi)\Big)$$

and

$$-x \in \Psi^{-1}\Big(N(\Phi)\Big).$$

Since $x$ can be written as

$$x = (x + x) + (-x) + 0 + 0$$

where $x + x \in \Phi^{-1}(N(\Psi))$, $-x \in \Psi^{-1}(N(\Phi))$, $0 \in N(\Phi)$ and $0 \in N(\Psi)$, then,

$$x \in \Phi^{-1}\Big(N(\Psi)\Big) + \Psi^{-1}\Big(N(\Phi)\Big) + N(\Phi) + N(\Psi).$$

Since $x \in N(\Phi\Psi)$ was arbitrarily fixed, then we have

$$N(\Phi\Psi) \subseteq \Phi^{-1}\Big(N(\Psi)\Big) + \Psi^{-1}\Big(N(\Phi)\Big) + N(\Phi) + N(\Psi).$$

Let now $x \in \Phi^{-1}(N(\Psi)) + \Psi^{-1}(N(\Phi)) + N(\Phi) + N(\Psi)$ be arbitrarily chosen. Hence,

$$x = x_1 + x_2 + n_1 + n_2,$$

where

$$x_1 \in \Phi^{-1}\Big(N(\Psi)\Big), \quad x_2 \in \Psi^{-1}\Big(N(\Phi)\Big), \quad n_1 \in N(\Phi) \quad \text{and} \quad n_2 \in N(\Psi).$$

Thus:

$\Phi(x_1) \in N(\Psi)$ and, therefore, $\Psi\Big(\Phi(x_1)\Big) = 0$, that is, $x_1 \in N(\Phi\Psi)$;

$\Psi(x_2) \in N(\Phi)$ and, therefore, $\Phi\Big(\Psi(x_2)\Big) = 0$, that is, $x_2 \in N(\Phi\Psi)$;

$\Phi(n_1) = 0$ and, therefore, $\Psi\Big(\Phi(n_1)\Big) = 0$, that is, $n_1 \in N(\Phi\Psi)$;

$\Psi(n_2) = 0$ and, therefore, $\Phi\Big(\Psi(n_2)\Big) = 0$, that is, $n_2 \in N(\Phi\Psi)$.

Since $N(\Phi\Psi)$ is a group, we conclude that

$$x = x_1 + x_2 + n_1 + n_2 \in N(\Phi\Psi)$$

and, since $x$ was arbitrarily chosen in $\Phi^{-1}(N(\Psi)) + \Psi^{-1}(N(\Phi)) + N(\Phi) + N(\Psi)$, we conclude that

$$\Phi^{-1}\Big(N(\Psi)\Big) + \Psi^{-1}\Big(N(\Phi)\Big) + N(\Phi) + N(\Psi) \subseteq N(\Phi\Psi);$$

**(g)** Let $x \in N(\Psi\varphi\Phi) + N(\Psi) + N(\Phi)$. Hence,

$$x = n_1 + n_2 + n_3$$

where $n_1 \in N(\Psi\varphi\Phi)$, $n_2 \in N(\Psi)$ and $n_3 \in N(\Phi)$. Thus,

$$\Psi(n_2) = \phi(n_3) = 0$$

and hence

$$(\Phi\varphi)\Big(\Psi(n_2)\Big) = 0 = (\Psi\varphi)\Big(\Phi(n_3)\Big),$$

which shows us that

$$n_2, n_3 \in N\left(\Phi\varphi\Psi\right).$$

Since $n_1 \in N(\Phi\varphi\Psi)$ and $N(\Phi\varphi\Psi)$ is a group, we conclude that

$$x = n_1 + n_2 + n_3 \in N(\Phi\varphi\Psi),$$

that is, since $x \in N(\Psi\varphi\Phi) + N(\Psi) + N(\Phi)$ was arbitrarily chosen,

$$N(\Psi\varphi\Phi) + N(\Psi) + N(\Phi) \subseteq N(\Phi\varphi\Psi) = N(\Psi\varphi\Phi);$$

**(h)** Suppose that $h = \Phi(g)$ with $g \in N(\Phi\Psi)$. In this case,

$$(\Phi\Psi)(g) = \Psi\Big(\Phi(g)\Big) = \Psi(h) = 0$$

and, therefore, $h \in N(\Psi)$. ■

# 3

# $S$-SPACES

## Introduction

### 3.1 Motivation

Let us consider the possibility suggested at the end of Chapter 2 of, starting with one of the finite order distributions axiomatics, let us say, the first version presented in 2.39, modifying it conveniently, to obtain another one that is categoric (like the latter) and admits as a model the pair $(D'(\Omega), D(\Omega))$ of the distributions in $\Omega$ and its derivatives.

Among the axioms in 2.39, as can be verified by inspection, except for the fifth, suppressing of them the saying "finite order" one obtain true claims about $(D'(\Omega), D(\Omega))$, and hence, with this modification they could be part of an axiomatic definition for the elements of $(D'(\Omega), D(\Omega))$. Regarding the referred Axiom 5, the suggested modification leads to a false claim about the class of all distributions since, as we know, there exist elements in $D'(\Omega)$, namely, the infinite order distributions, that can not be expressed as derivatives of continuous functions.

On the other hand, the L. Schwartz' distributions theory has a well-defined notion of local equality of distributions and, associated with this concept, a theorem establishing that: locally, every distribution is a derivative of some continuous function.[15] This suggests a local version of Axiom 5 to integrate a possible axiomatic for the class of (all) the distributions, something like:

- Every member $\Lambda \in D'(\Omega)$ is locally a finite order distribution.

But proceeding in this way, a new term, "local equality", becomes part of the axiomatic and, therefore, as a primitive term, requires new axioms to define it implicitly.

[15] In this initial considerations, of motivational nature, it is not relevant which, exactly, are the mentioned definition and theorem; in Chapter 5, when revisiting the L. Schwartz' distributions theory, we will be more specific about this point.

As a rule, regarding the concept of function, the notion of locality is associated with the restriction of a function to a subset of its domain and hence, taking into account that the distributions in $\Omega$ are linear functionals on $C^{\infty}_{\Omega}(\mathbb{R}^n)$, that is, functions, the natural choice seems to be the one in which "local equality" is considered not as a primitive concept, but as a concept defined from "domain" and "restriction" taken as primitives.

If now we remember that every continuous function with domain $\Omega$, that is, every $f \in C(\Omega)$, is a distribution in $\Omega$ in the sense of being identified, through isomorphism $\eta\colon C(\Omega) \longrightarrow C'(\Omega) \subseteq D'(\Omega)$ (see item 2.31), to the distribution in $\Omega$

$$\begin{aligned} \Lambda_f : C^{\infty}_{\Omega}(\mathbb{R}^n) &\longrightarrow C \\ \phi &\longmapsto \Lambda_f(\phi) = \textstyle\int_{\Omega} f(x)\phi(x)\,\mathrm{d}x\,, \end{aligned}$$

it appears then natural to define the domain of $\Lambda_f$, as of any distribution $\Lambda \in D'(\Omega)$, as the open set $\Omega$. In this case, such as to each pair $(\Omega' \subseteq \Omega, \Omega)$ of open sets of $\mathbb{R}^n$ corresponds the function

$$\begin{aligned} \rho_{(\Omega' \subseteq \Omega, \Omega)} : C(\Omega) &\longrightarrow C(\Omega') \\ f &\longmapsto \rho_{(\Omega' \subseteq \Omega, \Omega)}(f), \end{aligned}$$

where $\rho_{(\Omega' \subseteq \Omega, \Omega)}(f)$ is the usual restriction of the function $f \in C(\Omega)$ to the subset $\Omega'$ of its domain $\Omega$, we could considerate the possibility of defining, for each pair $(\Omega' \subseteq \Omega, \Omega)$ of open sets of $\mathbb{R}^n$, a function $\sigma_{(\Omega' \subseteq \Omega, \Omega)}\colon D'(\Omega) \longrightarrow D'(\Omega')$ that would play, regarding distributions, an analogous role to that of the function $\rho_{(\Omega' \subseteq \Omega, \Omega)}\colon C(\Omega) \longrightarrow C(\Omega')$ does to continuous functions. This being possible, and considering that the class of finite order distributions in $\Omega$, $D'_f(\Omega)$, is a subclass of the class $D'(\Omega)$ of all distributions in $\Omega$, $D'_f(\Omega) \subseteq D'(\Omega)$, then, the "restriction functions" $\sigma_{(\Omega', \Omega)}\colon D'(\Omega) \longrightarrow D'(\Omega')$ (one for each $(\Omega' \subseteq \Omega, \Omega)$) would endorse the following definition:

- We say that $\Lambda \in D'(\Omega)$ is locally a finite order distribution, if and only if, for each $x \in \Omega$ there exists an open set $\Omega' \subseteq \Omega$ such that $x \in \Omega'$ and $\sigma_{(\Omega', \Omega)}(\Lambda) \in D'_f(\Omega')$.

Hence, in short, the above considerations indicate to be reasonable, with regards to the possibility of an axiomatization of the class of (all) the distributions, to take as primitive terms, beyond those of distribution, addition and derivative, the ones of domain and restriction, defined implicitly by a set of axioms that incorporate those of the finite order distributions axiomatic in 2.39 (from which the saying "finite order" must be suppressed), except for the Axiom 5, beyond others providing a sensible characterization for the terms "domain" and "restriction" of a distribution, sensible in the sense of mirroring themselves in properties of these same concepts regarding the notion of function; in this environment we would then have the possibility to introduce, as a defined term, the notion of local equality and, with it, postulate that: locally, every distribution is a finite order distribution.

## 3.2 The Spaces $\mathscr{C}(\mathbb{R}^n)$, $\mathscr{D}'_f(\mathbb{R}^n)$ and $\mathscr{D}'(\mathbb{R}^n)$

It is important to observe that the considerations presented in the previous item, 3.1, as far as they refer, albeit in a very vague way, to concepts such as domain, restriction and local equality of distributions, transport us to an environment in which are present not only the three $S$-groups $\mathbb{C}(\Omega) = (C(\Omega), \partial(\Omega))$, $\mathbb{D}'_f(\Omega) = (D'_f(\Omega), D_f(\Omega))$ and $\mathbb{D}'(\Omega) = (D'(\Omega), D(\Omega))$ associated to an arbitrary, yet fixed, open set $\Omega \subseteq \mathbb{R}^n$, but, instead, a multitude of such triplet, one for each open set $\Omega \subseteq \mathbb{R}^n$. In fact, the notions of restriction and domain of distribution, as intuitively described in 3.1, through which the concept of local equality was, roughly, sketched, lift us to such an environment as far as the restrictions $\rho_{(\Omega',\Omega)} \colon C(\Omega) \longrightarrow C(\Omega')$ and $\sigma_{(\Omega',\Omega)} \colon D'(\Omega) \longrightarrow D'(\Omega')$ for each pair $(\Omega' \subseteq \Omega, \Omega)$ of open sets in $\mathbb{R}^n$, establish correspondences between functions and between distributions, respectively, with distinct domains.

More specifically, with the sets $\Gamma(\mathbb{R}^n)$ and $\Delta(\mathbb{R}^n)$ defined by

$$\Gamma(\mathbb{R}^n) \coloneqq \Big\{\Omega \subseteq \mathbb{R}^n \colon \Omega \text{ is an open set of } \mathbb{R}^n\Big\}$$

and

$$\Delta(\mathbb{R}^n) \coloneqq \Big\{(\Omega', \Omega) \in \Gamma(\mathbb{R}^n) \times \Gamma(\mathbb{R}^n) \colon \Omega' \subseteq \Omega\Big\},$$

the ideas suggested in Motivation 3.1 involve us with the following three families of $S$-groups indexed by $\Gamma(\mathbb{R}^n)$,

$$\mathbb{C}\Big(\Gamma(\mathbb{R}^n)\Big) \coloneqq \Big\{\mathbb{C}(\Omega) = \Big(C(\Omega), \partial(\Omega)\Big)\Big\}_{\Omega \in \Gamma(\mathbb{R}^n)},$$

$$\mathbb{D}'_f\Big(\Gamma(\mathbb{R}^n)\Big) \coloneqq \Big\{\mathbb{D}'_f(\Omega) = \Big(D'_f(\Omega), D_f(\Omega)\Big)\Big\}_{\Omega \in \Gamma(\mathbb{R}^n)}$$

and

$$\mathbb{D}'\Big(\Gamma(\mathbb{R}^n)\Big) \coloneqq \Big\{\mathbb{D}'(\Omega) = \Big(D'(\Omega), D(\Omega)\Big)\Big\}_{\Omega \in \Gamma(\mathbb{R}^n)},$$

as well as with the families constituted of the restrictions $\rho_{(\Omega',\Omega)} \colon C(\Omega) \to C(\Omega')$ and $\sigma_{(\Omega',\Omega)} \colon$
$D'(\Omega) \longrightarrow D'(\Omega')$, one for each pair $(\Omega', \Omega) \in \Delta(\mathbb{R}^n)$, that is,

$$\rho\Big(\Delta(\mathbb{R}^n)\Big) \coloneqq \Big\{\rho_{(\Omega',\Omega)}\Big\}_{(\Omega',\Omega) \in \Delta(\mathbb{R}^n)}$$

and

$$\sigma\Big(\Delta(\mathbb{R}^n)\Big) \coloneqq \Big\{\sigma_{(\Omega',\Omega)}\Big\}_{(\Omega',\Omega) \in \Delta(\mathbb{R}^n)}.$$

It is also involved, lets say, in a masked way, a family of restrictions associated to finite order distributions once that, by the fact that $D'_f(\Omega) \subseteq D'(\Omega)$ for each $\Omega \subseteq \mathbb{R}^n$, it is reasonable to admit that the restriction of the restriction $\sigma_{(\Omega',\Omega)} \colon D'(\Omega) \longrightarrow D'(\Omega')$ to the subset $D'_f(\Omega)$ of its domain $D'(\Omega)$, that is, the function

$$\begin{aligned} \tau_{(\Omega',\Omega)} \colon D'_f(\Omega) &\longrightarrow D'_f(\Omega') \\ \Lambda &\longmapsto \tau_{(\Omega',\Omega)}(\Lambda) \coloneqq \sigma_{(\Omega',\Omega)}(\Lambda), \end{aligned}$$

is also a restriction and, hence, the family

$$\tau\Big(\Delta(\mathbb{R}^n)\Big) := \Big\{\tau_{(\Omega',\Omega)}\Big\}_{(\Omega',\Omega)\in\Delta(\mathbb{R}^n)}$$

makes itself present.

In summary, the idea of formulating for the class of all the distributions a categorical axiomatic in which, as an axiom or theorem, every distribution is locally a finite order distribution, bonds us, that is, inserts us in a context where the "spaces"

$$\mathscr{C}(\mathbb{R}^n) := \Bigg(\mathbb{C}\Big(\Gamma(\mathbb{R}^n)\Big), \rho\Big(\Delta(\mathbb{R}^n)\Big)\Bigg),$$

$$\mathscr{D}'_f(\mathbb{R}^n) := \Bigg(\mathbb{D}'_f\Big(\Gamma(\mathbb{R}^n)\Big), \tau\Big(\Delta(\mathbb{R}^n)\Big)\Bigg),$$

and

$$\mathscr{D}'(\mathbb{R}^n) := \Bigg(\mathbb{D}'\Big(\Gamma(\mathbb{R}^n)\Big), \sigma\Big(\Delta(\mathbb{R}^n)\Big)\Bigg)$$

are necessary and naturally involved.

Observe now that the spaces $\mathscr{C}(\mathbb{R}^n)$, $\mathscr{D}'_f(\mathbb{R}^n)$ and $\mathscr{D}'(\mathbb{R}^n)$ keep between them the following relationship:

**(a)** the $S$-groups family of $\mathscr{D}'_f(\mathbb{R}^n)$ is the strict and closed extension of the corresponding family of $\mathscr{C}(\mathbb{R}^n)$, in the sense that, for each $\Omega \in \Gamma(\mathbb{R}^n)$, the $S$-group $\mathbb{D}'_f(\Omega) = (D'_f(\Omega), D_f(\Omega))$ is the (unique unless isomorphism) strict and closed extension of the $S$-group $\mathbb{C}(\Omega) = (C(\Omega), \partial(\Omega))$;

**(b)** the $S$-groups family of $\mathscr{D}'(\mathbb{R}^n)$ is an extension of the $S$-groups family of $\mathscr{D}'_f(\mathbb{R}^n)$, in the sense that, for each $\Omega \in \Gamma(\mathbb{R}^n)$, the $S$-group $\mathbb{D}'(\Omega) = (D'(\Omega), D(\Omega))$ is an extension of the $S$-group $\mathbb{D}'_f(\Omega) = (D'_f(\Omega), D_f(\Omega))$;

**(c)** the $S$-groups family of $\mathscr{D}'(\mathbb{R}^n)$ is a strict, but not closed, extension of the corresponding family of $\mathscr{C}(\mathbb{R}^n)$, in an analogous sense to the above described;

**(d)** the family of restrictions of $\mathscr{D}'(\mathbb{R}^n)$, $\sigma(\Delta(\mathbb{R}^n))$, is an extension of the family of restrictions $\tau(\Delta(\mathbb{R}^n))$ of the space $\mathscr{D}'_f(\mathbb{R}^n)$, that is, for each pair $(\Omega', \Omega) \in \Delta(\mathbb{R}^n)$ the function $\sigma_{(\Omega',\Omega)}\colon D'(\Omega) \longrightarrow D'(\Omega')$ is an extension to $D'(\Omega)$ of the function $\tau_{(\Omega',\Omega)}\colon D'_f(\Omega) \longrightarrow D'_f(\Omega')$ ($\sigma_{(\Omega',\Omega)}(\Lambda) = \tau_{(\Omega',\Omega)}(\Lambda)$ if $\Lambda \in D'_f(\Omega)$).

Based now in these observations, allied to an obvious (informal) notion of extension for the spaces $\mathscr{C}(\mathbb{R}^n)$, $\mathscr{D}'_f(\mathbb{R}^n)$ and $\mathscr{D}'(\mathbb{R}^n)$, we have:

**(e)** from (b) and (d), that $\mathscr{D}'(\mathbb{R}^n)$ is an extension of $\mathscr{D}'_f(\mathbb{R}^n)$;

**(f)** from (c) and in the hypothesis of the functions $\sigma_{(\Omega',\Omega)} \in \sigma(\Delta(\mathbb{R}^n))$ being extensions to $D'(\Omega)$ of the corresponding restrictions $\rho_{(\Omega',\Omega)} \in \rho(\Delta(\mathbb{R}^n))$, that $\mathscr{D}'(\mathbb{R}^n)$ is a strict, but not closed, extension of $\mathscr{C}(\mathbb{R}^n)$.

## 3.3 Work Agenda

In the scenario described in the previous item, 3.2, the following questions seem inevitable:

**(a)** Is there a sensible definition that can be given to the functions $\sigma_{(\Omega',\Omega)}$, sensible in the sense of being aligned with the classical notion of restriction of functions, in such a way that the space $\mathscr{D}'(\mathbb{R}^n) = (\mathbb{D}'(\Gamma(\mathbb{R}^n)), \sigma(\Delta(\mathbb{R}^n)))$ becomes a locally closed extension of the space $\mathscr{D}'_f(\mathbb{R}^n) = (\mathbb{D}'_f(\Gamma(\mathbb{R}^n)), \tau(\Delta(\mathbb{R}^n)))$, that is, such that: for each $\Omega \in \Gamma(\mathbb{R}^n)$ and each $x \in \Omega$, if $\Lambda \in D'(\Omega)$, there exists $\Omega' \in \Gamma(\mathbb{R}^n)$ such that $\Omega' \subseteq \Omega$, $x \in \Omega'$ and $\sigma_{(\Omega',\Omega)}(\Lambda) \in D'_f(\Omega')$?

**(b)** In the hypothesis of existing a family $\sigma(\Delta(\mathbb{R}^n))$ such that its members, $\sigma_{(\Omega',\Omega)}$, satisfy the condition described in (a), i.e., that the space $\mathscr{D}'(\mathbb{R}^n) = (\mathbb{D}'(\Gamma(\mathbb{R}^n)), \sigma(\Delta(\mathbb{R}^n)))$ becomes a locally closed extension of $\mathscr{D}'_f(\mathbb{R}^n) = (\mathbb{D}'_f(\Gamma(\mathbb{R}^n)), \tau(\Delta(\mathbb{R}^n)))$, would $\mathscr{D}'(\mathbb{R}^n)$ be the only such space?

It is opportune to remember now that the finite order distributions in $\Omega$ are, exactly, the $\widetilde{\mathbb{C}}(\Omega)$-distributions, that is, the $\widetilde{\mathbb{G}}$-distributions specialized to the case where $\mathbb{G} = \mathbb{C}(\Omega)$ (see 2.38); the axiomatic of the finite order distributions in $\Omega$, is obtained by specializing the $\widetilde{\mathbb{G}}$-distributions axioms to the case where the $S$-group $\mathbb{G}$ is that of the continuous functions on $\Omega$, $\mathbb{C}(\Omega)$ (see 2.39).

On the other hand, the $\widetilde{\mathbb{G}}$-distributions axioms were extracted from the Theorem of Extension of $S$-Groups (Theorem 2.16), which establishes the existence of a single strict and closed extension of the $S$-group $\mathbb{G}$; the defining conditions of the concept of strict and closed extension, to which the theorem refers to, was literally taken as axioms (see Remark 2.19 as well as item 2.20).

That said, i.e., that the axiomatic of the finite order distributions in $\Omega$ was obtained from the Theorem of Extension of $S$-Groups, and taking into account that the questions (a) and (b) above formulated inquire about the existence and uniqueness of a locally closed extension of the space $\mathscr{D}'_f(\mathbb{R}^n)$ that, in turn, can be taken (considering 3.2(a) and on the hypothesis of the functions $\tau_{(\Omega',\Omega)}$ being extensions to $D'_f(\Omega)$ of the restrictions $\rho_{(\Omega',\Omega)}$) as a strict and closed extension of the space $\mathscr{C}(\mathbb{R}^n)$, it appears reasonable the following program to obtain a categoric axiomatic for the class of all the distributions (finite and infinite order):

- To formulate in a precise way the following extension problem:

  "To obtain a locally closed extension of the strict and closed extension, $\mathscr{D}'_f(\mathbb{R}^n)$, of the space $\mathscr{C}(\mathbb{R}^n)$."

- To state and prove an extension theorem relative to the problem above formulated, something like:

  "There exists essentially an unique locally closed extension of the extension $\mathscr{D}'_f(\mathbb{R}^n)$ of the space $\mathscr{C}(\mathbb{R}^n)$."

- To extract from the extension theorem obtained, as was done with the Theorem of Extension of $S$-Groups to obtain axioms for the finite order distributions, an axiomatic where the axioms express the defining properties of the type of extension whose existence and uniqueness are ensured in the obtained theorem.

- Finally, to prove that the axiomatic formulated in the previous stage does have a model, and that any of its models constitutes an extension whose existence and uniqueness are ensured by the extension theorem obtained, establishing then the consistency and categoricity of the axiomatic.

## 3.4 Work Agenda — Abstract Version

As previously observed, the axioms of the finite order distributions were obtained from the "abstract" axiomatic of the $\widetilde{\mathbb{G}}$-distributions, abstract in the sense of its axioms referring not to a particular $S$-group but any abelian, surjective, and with identity $S$-group; the simple translation of the $\widetilde{\mathbb{G}}$-distributions axioms to the case where $\mathbb{G}$ is the continuous functions $S$-group, $\mathbb{C}(\Omega)$, leads to the finite order distributions axiomatic.

It has also been highlighted that the $\widetilde{\mathbb{G}}$-distributions axioms are extracted from the Theorem of Extension of $S$-Groups which, in turn, refers not to a particular $S$-group but anyone among the abelian, surjective, and with identity ones. Case we did not have proved the referred theorem in its abstract version, that is, with relation to any $S$-group $\mathbb{G}$ (abelian, surjective, and with identity), if instead we knew what establishes this theorem only for the $S$-group $\mathbb{C}(\Omega)$, from it we would extract the axiomatic of the finite order distributions but not the $\widetilde{\mathbb{G}}$-distributions one.

Now we inquire if something analogous would not be possible with respect to our actual intention. More precisely, we question if the work agenda established in item 3.3, in the figure of the program there proposed to obtain an axiomatic for the class of all the distributions, could not be formulated and, above all, executed, not only for the particular space $\mathscr{C}(\mathbb{R}^n) = (\mathbb{C}(\Gamma(\mathbb{R}^n)), \rho(\Delta(\mathbb{R}^n)))$ and its strict and closed "extension" $\mathscr{D}'_f(\mathbb{R}^n) = (\mathbb{D}'_f(\Gamma(\mathbb{R}^n)), \tau(\Delta(\mathbb{R}^n)))$, but with respect to a pair of abstract spaces,

$$\mathscr{G}(I) = \Big(\mathbb{G}\Big(\Gamma(I)\Big), \Theta\Big(\Delta(I)\Big)\Big)$$

and

$$\widehat{\mathscr{G}}(I) = \Big(\widehat{\mathbb{G}}\Big(\Gamma(I)\Big), \widehat{\Theta}\Big(\Delta(I)\Big)\Big),$$

where $\mathbb{G}(\Gamma(I))$ is a family of abstract $S$-groups that (as the ones from the family $\mathbb{C}(\Gamma(\mathbb{R}^n)))$ are abelian, surjective, and with identity, and $\widehat{\mathbb{G}}(\Gamma(I))$ is the family of its corresponding strict and closed extensions (like the family $\mathbb{D}'_f(\Gamma(\mathbb{R}^n))$ in its relation with $\mathbb{C}(\Gamma(\mathbb{R}^n)))$.

A big part of the rest of this monograph is intended to execute this abstract version of the program formulated at the end of item 3.3, which, roughly speaking, is obtained (from the one in 3.3) exchanging the spaces $\mathscr{C}(\mathbb{R}^n)$ and $\mathscr{D}'_f(\mathbb{R}^n)$ by its respective abstract versions $\mathscr{G}(I)$ and $\widehat{\mathscr{G}}(I)$.

In this chapter, in particular, we will occupy ourselves with the formulation of a structure — the $S$-space structure — that admits as models the spaces $\mathscr{C}(\mathbb{R}^n)$, $\mathscr{D}'_f(\mathbb{R}^n)$ and $\mathscr{D}'(\mathbb{R}^n)$ informally characterized in our previous considerations.

# Indexer Semigroup

## 3.5 The Semigroup $\Gamma(\mathbb{R}^n)$

Our first task in order to define the concept of $S$-space, that, above mentioned, that should play, in the abstract context of the work agenda proposed in 3.4, the role of the spaces $\mathscr{C}(\mathbb{R}^n)$, $\mathscr{D}'_f(\mathbb{R}^n)$ and $\mathscr{D}'(\mathbb{R}^n)$, consists in the characterization of the set $\Gamma(\mathbb{R}^n)$ used to index the families of $S$-groups of these spaces. As we saw,

$$\Gamma(\mathbb{R}^n) := \Big\{\Omega \in \mathcal{P}(\mathbb{R}^n) \colon \Omega \text{ is an open set of } \mathbb{R}^n\Big\},$$

where $\mathcal{P}(I)$ is the set of the subsets of a given set $I$. It results from the definition of the set $\Gamma(\mathbb{R}^n)$ that:

- $\Gamma(\mathbb{R}^n) \subseteq \mathcal{P}(\mathbb{R}^n)$;
- $\mathbb{R}^n \in \Gamma(\mathbb{R}^n)$;
- if $\Omega, \Omega' \in \Gamma(\mathbb{R}^n)$, then, $\Omega \cap \Omega' \in \Gamma(\mathbb{R}^n)$, that is, $\Gamma(\mathbb{R}^n)$ is closed relative to the operation of intersection of sets;

- if $\Omega \in \Gamma(\mathbb{R}^n)$, then, the set

$$\Gamma(\Omega) := \Gamma(\mathbb{R}^n) \cap \mathcal{P}(\Omega),$$

  such as $\Gamma(\mathbb{R}^n)$, satisfies: $\Gamma(\Omega) \subseteq \mathcal{P}(\Omega)$, $\Omega \in \Gamma(\Omega)$ and $\Gamma(\Omega)$ is closed relative to the operation of intersection of sets.

In other words, $\Gamma(\mathbb{R}^n)$ is, with the operation of intersection of sets restricted to its elements, an abelian semigroup with identity ($\mathbb{R}^n \in \Gamma(\mathbb{R}^n)$) such that, for each $\Omega \in \Gamma(\mathbb{R}^n)$,

$$\Gamma(\Omega) := \Gamma(\mathbb{R}^n) \cap \mathcal{P}(\Omega) \subseteq \mathcal{P}(\Omega)$$

is also an abelian semigroup with identity ($\Omega \in \Gamma(\Omega)$).

The families of abstract $S$-groups that will replace, in the abstract approach outlined in 3.4, the corresponding ones in the spaces $\mathscr{C}(\mathbb{R}^n)$, $\mathscr{D}'_f(\mathbb{R}^n)$ and $\mathscr{D}'(\mathbb{R}^n)$, will be indexed by sets whose (axiomatic) definition is inspired in the properties of the set $\Gamma(\mathbb{R}^n)$ above highlighted.

## 3.6 Definition

Let $I$ be an arbitrarily fixed non-empty set and $\mathcal{P}(I)$ be the powerset of $I$ (whose members are the subsets of $I$). Let also $\Gamma(I) \subseteq \mathcal{P}(I)$. We say that $\Gamma(I)$ is an ***I*-indexer**, or simply an **indexer** (when the mention to the set $I$ is not necessary), if and only if it satisfies the following conditions:

**(a)** $I \in \Gamma(I)$;

**(b)** whatever $\gamma, \gamma' \in \Gamma(I)$, one has that $\gamma \cap \gamma' \in \Gamma(I)$.

Take now an $I$-indexer, $\Gamma(I)$, and $\gamma \in \Gamma(I)$ such that $\gamma \neq \varnothing$. It results directly from the definition above that:

- the set $\Gamma(\gamma)$ defined by

$$\Gamma(\gamma) := \Gamma(I) \cap \mathcal{P}(\gamma)$$

  is a $\gamma$-indexer, i.e., it is such that $\gamma \in \Gamma(\gamma)$ and is closed relative to the intersection of sets (if $\gamma', \gamma'' \in \Gamma(\gamma)$, then $\gamma' \cap \gamma'' \in \Gamma(\gamma)$);

- if $\gamma' \in \Gamma(I)$ is such that $\gamma' \subseteq \gamma$,

$$\Gamma(\gamma') \subseteq \Gamma(\gamma);$$

- $\Gamma(\gamma)$ with the operation of intersection of sets restricted to its elements, is an abelian semigroup with identity ($\gamma \in \Gamma(\gamma)$).

# Family and Subfamily of $S$-Groups

## 3.7 Definition

The term **family of $S$-groups** will be used here, unless otherwise explicitly informed, to refer to a family whose set of indexes is an $I$-indexer, $\Gamma(I)$, and whose members are $S$-groups; if for each index $\gamma \in \Gamma(I)$ the corresponding $S$-group of such a family is denoted by $\mathbb{G}(\gamma)$, the family will then be denoted by $\mathbb{G}(\Gamma(I))$, that is,

$$\mathbb{G}\Big(\Gamma(I)\Big) = \Big\{\mathbb{G}(\gamma)\Big\}_{\gamma\in\Gamma(I)}.$$

Now, as we know, for each $\gamma \in \Gamma(I)$, the set $\Gamma(\gamma) = \Gamma(I) \cap \mathcal{P}(\gamma)$ is a subset of $\Gamma(I)$, $\Gamma(\gamma) \subseteq \Gamma(I)$, which is also an indexer ($\gamma$-indexer). Hence, being $\mathbb{G}(\Gamma(I))$, as above, a family of $S$-groups, then, $\mathbb{G}(\Gamma(\gamma)) \subseteq \mathbb{G}(\Gamma(I))$ and $\mathbb{G}(\Gamma(\gamma))$ also is a family of $S$-groups; we will say that $\mathbb{G}(\Gamma(\gamma))$, in its relationship with $\mathbb{G}(\Gamma(I))$, is a **subfamily** of $\mathbb{G}(\Gamma(I))$ and, unless otherwise explicitly said, the term "subfamily of $S$-groups" will be used with this meaning.

Taking into account the properties of the set

$$\Gamma(\mathbb{R}^n) = \Big\{\Omega \in \mathcal{P}(\mathbb{R}^n) : \Omega \text{ is an open set of } \mathbb{R}^n\Big\}$$

presented in 3.5, and the Definition 3.6 of $I$-indexer, one concludes that $\Gamma(\mathbb{R}^n)$ is an indexer. Hence, as important examples of families of $S$-groups we have

$$\mathbb{C}\Big(\Gamma(\mathbb{R}^n)\Big) := \Big\{\mathbb{C}(\Omega) = \Big(C(\Omega), \partial(\Omega)\Big)\Big\}_{\Omega\in\Gamma(\mathbb{R}^n)},$$

$$\mathbb{D}'_f\Big(\Gamma(\mathbb{R}^n)\Big) := \Big\{\mathbb{D}'_f(\Omega) = \Big(D'_f(\Omega), D_f(\Omega)\Big)\Big\}_{\Omega\in\Gamma(\mathbb{R}^n)}$$

and

$$\mathbb{D}'\Big(\Gamma(\mathbb{R}^n)\Big) := \Big\{\mathbb{D}'(\Omega) = \Big(D'(\Omega), D(\Omega)\Big)\Big\}_{\Omega\in\Gamma(\mathbb{R}^n)},$$

to which we will refer, respectively, as the family of $S$-groups of the continuous functions, of the finite order distributions and of the distributions.

## 3.8 Definition

Let $\mathbb{G}(\Gamma(I)) = \{\mathbb{G}(\gamma)\}_{\gamma\in\Gamma(I)}$ and $\mathbb{E}(\Gamma(I)) = \{\mathbb{E}(\gamma)\}_{\gamma\in\Gamma(I)}$ be families of $S$-groups indexed by the same indexer $\Gamma(I)$. We will say that $\mathbb{E}(\Gamma(I))$ is an **extension of the family** $\mathbb{G}(\Gamma(I))$ if and only if, for each $\gamma \in \Gamma(I)$, the $S$-group $\mathbb{E}(\gamma)$ is an extension (in the sense of Definition 1.6(d)) of the $S$-group $\mathbb{G}(\gamma)$ and, if so, the family $\mathbb{E}(\Gamma(I))$ will be said an **strict** (respectively **closed**) **extension** of the family $\mathbb{G}(\Gamma(I))$ if and only if, for each $\gamma \in \Gamma(I)$, the $S$-group $\mathbb{E}(\gamma)$ is an strict (respectively closed) extension of the $S$-group $\mathbb{G}(\gamma)$.

We will say yet that a family of $S$-groups, let us say $\mathbb{G}(\Gamma(I)) = \{\mathbb{G}(\gamma)\}_{\gamma\in\Gamma(I)}$, is **abelian**, **surjective** or **with identity**, if and only if, for each $\gamma \in \Gamma(I)$, the $S$-group $\mathbb{G}(\gamma)$ is abelian, surjective or with identity, respectively.

Regarding the families of $S$-groups of the continuous functions, $\mathbb{C}(\Gamma(\mathbb{R}^n))$, of the finite order distributions, $\mathbb{D}'_f(\Gamma(\mathbb{R}^n))$ and of the distributions, $\mathbb{D}'(\Gamma(\mathbb{R}^n))$, referred in Definition 3.7, we have:

- $\mathbb{C}\Big(\Gamma(\mathbb{R}^n)\Big)$ is an abelian and surjective family with identity;
- $\mathbb{D}'_f\Big(\Gamma(\mathbb{R}^n)\Big)$ is a strict and closed extension of the family $\mathbb{C}\Big(\Gamma(\mathbb{R}^n)\Big)$;
- $\mathbb{D}'\Big(\Gamma(\mathbb{R}^n)\Big)$ is an extension of $\mathbb{D}'_f\Big(\Gamma(\mathbb{R}^n)\Big)$;
- $\mathbb{D}'\Big(\Gamma(\mathbb{R}^n)\Big)$ is a strict extension of $\mathbb{C}\Big(\Gamma(\mathbb{R}^n)\Big)$ that is not closed.

# Bonded Family

## 3.9 Remarks

It is opportune to remember that our goal, motivated and described informally in the Introduction of this chapter (items 3.1 to 3.4) and synthesized in 3.4 as a work program, consists, essentially, of the following (recap item 3.4): to prove, with relation to a pair of abstract (axiomatic) structures $\mathscr{G}(I) = (\mathbb{G}(\Gamma(I)), \Theta(\Delta(I)))$ and $\widehat{\mathscr{G}}(I) = (\widehat{\mathbb{G}}(\Gamma(I)), \widehat{\Theta}(\Delta(I)))$ of which are models, respectively, the spaces $\mathscr{C}(\mathbb{R}^n) = (\mathbb{C}(\Gamma(\mathbb{R}^n)), \rho(\Delta(\mathbb{R}^n)))$ and $\mathscr{D}'_f(\mathbb{R}^n) = (\mathbb{D}'_f(\Gamma(\mathbb{R}^n)), \tau(\Delta(\mathbb{R}^n)))$, a theorem that ensures the existence and uniqueness of a locally closed extension of $\widehat{\mathscr{G}}(I)$.

Hence, it is more than reasonable that the defining axioms of this abstract structure suitable to our purposes, mirror themselves in the relevant properties of those objects, with the spaces $\mathscr{C}(\mathbb{R}^n)$ and $\mathscr{D}'_f(\mathbb{R}^n)$ as examples, that motivate its definition and, for this reason, will be models of it. In particular, the concept of restriction for the family of abstract $S$-groups, $\mathbb{G}(\Gamma(I))$, of the structure $\mathscr{G}(I) = (\mathbb{G}(\Gamma(I)), \Theta(\Delta(I)))$, expressed through the family $\Theta(\Delta(I))$, should be defined (implicitly) by a set of properties shared by the members of the family $\rho(\Delta(\mathbb{R}^n))$ which, in turn, express, through the usual notion of restriction of a function to a subset of its domain, the concept of restriction for the family of $S$-groups of the continuos functions, $\mathbb{C}(\Gamma(\mathbb{R}^n))$, of the space $\mathscr{C}(\mathbb{R}^n) = (\mathbb{C}(\Gamma(\mathbb{R}^n)), \rho(\Delta(\mathbb{R}^n)))$. In the same way and by the same reasons, the characterization of the family of abstract $S$-groups, $\mathbb{G}(\Gamma(I))$, of the structure $\mathscr{G}(I)$, should be inspired in the distinctive aspects of the family of $S$-groups of the space $\mathscr{C}(\mathbb{R}^n)$, namely, the family $\mathbb{C}(\Gamma(\mathbb{R}^n)) \coloneqq \{\mathbb{C}(\Omega) = (C(\Omega), \partial(\Omega))\}_{\Omega \in \Gamma(\mathbb{R}^n)}$ of the $S$-groups of the continuous functions. In the next item, 3.10, we will highlight a property of the family $\mathbb{C}(\Gamma(\mathbb{R}^n))$, more specifically of the semigroups $\partial(\Omega)$ of the $S$-groups $\mathbb{C}(\Omega) = (C(\Omega), \partial(\Omega))$, which should appear in the definition of the structure $\mathscr{G}(I)$. Regarding the notion of restriction, an analogous investigation will be made in 3.12.

## 3.10 $\partial(\Omega)$ and $\partial(\Omega')$: Isomorphic Semigroups

Let us consider the family of the $S$-groups of the continuous functions,

$$\mathbb{C}\Big(\Gamma(\mathbb{R}^n)\Big) = \Big\{\mathbb{C}(\Omega) = \Big(C(\Omega), \partial(\Omega)\Big)\Big\}_{\Omega \in \Gamma(\mathbb{R}^n)},$$

and let us highlight from it the semigroups $\partial(\Omega)$ and $\partial(\Omega')$ correspondents to the indexes $\Omega, \Omega' \in \Gamma(\mathbb{R}^n)$, arbitrarily fixed. As we know (recap itens 1.2 and 1.5),

$$\partial(\Omega) = \Big\{\partial^\alpha_\Omega \colon \alpha \in \mathbb{N}^n\Big\}$$

and

$$\partial(\Omega') = \Big\{\partial^\alpha_{\Omega'} \colon \alpha \in \mathbb{N}^n\Big\},$$

with $\partial^\alpha_\Omega \colon C^{|\alpha|}(\Omega) \longrightarrow C(\Omega)$ and $\partial^\alpha_{\Omega'} \colon C^{|\alpha|}(\Omega') \longrightarrow C(\Omega')$ being the classical partial derivatives. Clearly, $\partial(\Omega)$ and $\partial(\Omega')$ are sets of the same cardinality and, beyond that, the "natural" bijection between them, namely,

$$j_{(\Omega',\Omega)} \colon \partial(\Omega) \longrightarrow \partial(\Omega')$$

where

$$j_{(\Omega',\Omega)}\Big(\partial^\alpha_\Omega\Big) \coloneqq \partial^\alpha_{\Omega'} \quad \text{for every} \quad \alpha \in \mathbb{N}^n,$$

is an isomorphism from the semigroup $\partial(\Omega)$ onto the semigroup $\partial(\Omega')$. In fact, for any multi-indexes $\alpha, \beta \in \mathbb{N}^n$ we have:

$$j_{(\Omega',\Omega)}\Big(\partial^\alpha_\Omega \partial^\beta_\Omega\Big) = j_{(\Omega',\Omega)}\Big(\partial^{\alpha+\beta}_\Omega\Big) = \partial^{\alpha+\beta}_{\Omega'} = \partial^\alpha_{\Omega'}\partial^\beta_{\Omega'} = j_{(\Omega',\Omega)}\Big(\partial^\alpha_\Omega\Big) j_{(\Omega',\Omega)}\Big(\partial^\beta_\Omega\Big).$$

Furthermore, we have that $\big(j_{(\Omega',\Omega)}\big)^{-1} = j_{(\Omega,\Omega')}$ and the composition of $j_{(\Omega'',\Omega')}$ with $j_{(\Omega',\Omega)}$, in this order, results in $j_{(\Omega'',\Omega)}$, that is, $j_{(\Omega'',\Omega')} j_{(\Omega',\Omega)} = j_{(\Omega'',\Omega)}$.

So we see that the semigroups of the family of the $S$-groups of the continuous functions, are related to each other through a bond expressed by the isomorphisms $j_{(\Omega',\Omega)}$, one for each pair $(\Omega', \Omega) \in \Gamma(\mathbb{R}^n) \times \Gamma(\mathbb{R}^n)$, that is, by the family

$$j\Big(\Gamma^2(\mathbb{R}^n)\Big) \coloneqq \Big\{j_{(\Omega',\Omega)}\Big\}_{(\Omega',\Omega) \in \Gamma^2(\mathbb{R}^n)}$$

indexed by the set $\Gamma^2(\mathbb{R}^n) := \Gamma(\mathbb{R}^n) \times \Gamma(\mathbb{R}^n)$.

This bond between the semigroups of the family of $S$-groups $\mathbb{C}(\Gamma(\mathbb{R}^n))$ is, no doubt, a distinctive aspect of this family and, hence being, it is reasonable for it to participate of the characterization of the abstract structure, $\mathscr{G}(I) = (\mathbb{G}(\Gamma(I)), \Theta(\Delta(I)))$, that we look forward to define, particularly regarding the defining properties of its family of $S$-groups $\mathbb{G}(\Gamma(I))$. With this purpose in view, we formulate the next definition.

## 3.11 Definition

Let

$$\mathbb{G}\Big(\Gamma(I)\Big) = \Big\{\mathbb{G}(\gamma) = \Big(G(\gamma), H(\gamma)\Big)\Big\}_{\gamma \in \Gamma(I)}$$

be a family of $S$-groups for which there is a family

$$i\Big(\Gamma^2(I)\Big) := \Big\{i_{(\gamma',\gamma)}\Big\}_{(\gamma',\gamma)\in\Gamma^2(I)}$$

indexed by the set $\Gamma^2(I) := \Gamma(I) \times \Gamma(I)$, such that:

**(a)** for each $(\gamma', \gamma) \in \Gamma^2(I)$, $i_{(\gamma',\gamma)}$ is an isomorphism from the semigroup $H(\gamma)$ onto the semigroup $H(\gamma')$, that is, a function from $H(\gamma)$ onto $H(\gamma')$,

$$i_{(\gamma',\gamma)} \colon H(\gamma) \longrightarrow H(\gamma'),$$

bijective and such that

$$i_{(\gamma',\gamma)}\left(\Phi\Psi\right) = i_{(\gamma',\gamma)}\left(\Phi\right) i_{(\gamma',\gamma)}\left(\Psi\right),$$

whatever $\Phi, \Psi \in H(\gamma)$;

**(b)** $\Big(i_{(\gamma',\gamma)}\Big)^{-1} = i_{(\gamma,\gamma')}$ for each $(\gamma', \gamma) \in \Gamma^2(I)$;

**(c)** $i_{(\gamma'',\gamma')} i_{(\gamma',\gamma)} = i_{(\gamma'',\gamma)}$ whatever $\gamma, \gamma', \gamma'' \in \Gamma(I)$.

In this case, we will say that $i(\Gamma^2(I))$ is a **bonding family** or a **bonding** of the family of $S$-groups $\mathbb{G}(\Gamma(I))$, and also that the ordered pair

$$\left(\mathbb{G}\Big(\Gamma(I)\Big), i\Big(\Gamma^2(I)\Big)\right)$$

is a **bonded family**.

**Remark.** It results that $i_{(\gamma,\gamma)} = I_{H(\gamma)}$, the identity function on $H(\gamma)$. In fact, from 3.11(c), $i_{(\gamma,\gamma)} i_{(\gamma,\gamma)} = i_{(\gamma,\gamma)}$ and if $\Phi \in H(\gamma)$ is arbitrarily chosen, then, by 3.11(a) there exists $\Psi \in H(\gamma)$ such that $\Phi = i_{(\gamma,\gamma)}(\Psi)$, and hence we have that

$$i_{(\gamma,\gamma)}(\Phi) = i_{(\gamma,\gamma)}\Big(i_{(\gamma,\gamma)}(\Psi)\Big) = i_{(\gamma,\gamma)}(\Psi) = \Phi.$$

# The Concept of Restriction

## 3.12 Preliminaries

As suggested by Remarks 3.9, in order to obtain indications about the properties that could compose, in the form of axioms, a sensible notion of restriction, the one represented by the family $\Theta(\Delta(I))$ of the abstract structure $\mathscr{G}(I) = (\mathbb{G}(\Gamma(I)), \Theta(\Delta(I)))$ we want to define accordingly to our purposes (defined in the work program formulated in 3.4), we must inspire ourselves in the properties of the family of restrictions $\rho(\Delta(\mathbb{R}^n))$ of the space $\mathscr{C}(\mathbb{R}^n) = (\mathbb{C}(\Gamma(\mathbb{R}^n)), \rho(\Delta(\mathbb{R}^n)))$. In turn, the family $\rho(\Delta(\mathbb{R}^n))$ have as members the functions $\rho_{(\Omega',\Omega)}$, one for each $(\Omega', \Omega) \in \Delta(\mathbb{R}^n)$, defined by:

$$\begin{aligned} \rho_{(\Omega',\Omega)} : C(\Omega) &\longrightarrow C(\Omega') \\ f &\longmapsto \rho_{(\Omega',\Omega)}(f) := f|_{\Omega'} \end{aligned}$$

where $f|_{\Omega'}$ is the usual restriction of the function $f\colon \Omega \longrightarrow C$ to the subset $\Omega' \subseteq \Omega$ of its domain $\Omega$, that is,

$$f|_{\Omega'}(x) = f(x) \quad \text{for every} \quad x \in \Omega' \subseteq \Omega.$$

It is convenient to remember here that

$$\Gamma(\mathbb{R}^n) = \Big\{\Omega \subseteq \mathbb{R}^n : \Omega \text{ is an open set of } \mathbb{R}^n\Big\}$$

and

$$\Delta(\mathbb{R}^n) = \Big\{(\Omega', \Omega) \in \Gamma^2(\mathbb{R}^n) : \Omega' \subseteq \Omega\Big\},$$

and hence, if $(\Omega', \Omega) \in \Delta(\mathbb{R}^n)$, then, $\Omega'$ and $\Omega$ are open sets of $\mathbb{R}^n$ such that $\Omega' \subseteq \Omega$, which validates the definition of $\rho_{(\Omega',\Omega)}$ given above.

Now, as can easily be proved, the family $\rho(\Delta(\mathbb{R}^n))$ has the following properties:

**(a)** for each $(\Omega', \Omega) \in \Delta(\mathbb{R}^n)$,

$$\rho_{(\Omega',\Omega)}(f + g) = \rho_{(\Omega',\Omega)}(f) + \rho_{(\Omega',\Omega)}(g)$$

whatever $f, g \in C(\Omega)$, that is, $\rho_{(\Omega',\Omega)} : C(\Omega) \longrightarrow C(\Omega')$ is a homomorphism from the group $C(\Omega)$ into the group $C(\Omega')$;

**(b)** $\rho_{(\Omega'',\Omega')}\rho_{(\Omega',\Omega)} = \rho_{(\Omega'',\Omega)}$, that is,

$$\rho_{(\Omega'',\Omega')}\Big(\rho_{(\Omega',\Omega)}(f)\Big) = \rho_{(\Omega'',\Omega)}(f)$$

for every $f \in C(\Omega)$ and every $\Omega, \Omega', \Omega'' \in \Gamma(\mathbb{R}^n)$ such that $\Omega'' \subseteq \Omega' \subseteq \Omega$.

As we can see from the properties above, the family $\rho(\Delta(\mathbb{R}^n))$ establishes a relation, a bond, between the groups of the $S$-groups of the family

$$\mathbb{C}\Big(\Gamma(\mathbb{R}^n)\Big) = \Big\{\mathbb{C}(\Omega) = \Big(C(\Omega), \partial(\Omega)\Big)\Big\}_{\Omega \in \Gamma(\mathbb{R}^n)}.$$

On the other hand, the family $j(\Gamma^2(\mathbb{R}^n))$ defined in 3.10 relates, bonds, the corresponding semigroups. In turn, these two families, $\rho(\Delta(\mathbb{R}^n))$ and $j(\Gamma^2(\mathbb{R}^n))$, are also related (bonded) since, as can be easily proved:

**(c)** if $(\Omega', \Omega) \in \Delta(\mathbb{R}^n)$, $\partial^\alpha_\Omega \in \partial(\Omega)$ and $f \in C^{|\alpha|}(\Omega)$ (the domain of $\partial^\alpha_\Omega$), then

$$\rho_{(\Omega',\Omega)}\Big(\partial^\alpha_\Omega(f)\Big) = j_{(\Omega',\Omega)}(\partial^\alpha_\Omega)\Big(\rho_{(\Omega',\Omega)}(f)\Big).$$

This last property establishes a kind of commutativity between the families $\rho(\Delta(\mathbb{R}^n))$ and $j(\Gamma^2(\mathbb{R}^n))$, hence expliciting an important correlation between them. In turn, this entanglement suggests that it is more appropriate to refer to $\rho(\Delta(\mathbb{R}^n))$, in its relation with the family of $S$-groups $\mathbb{C}(\Gamma(\mathbb{R}^n))$, not as "restriction for $\mathbb{C}(\Gamma(\mathbb{R}^n))$" but "restriction for the bonded family $(\mathbb{C}(\Gamma(\mathbb{R}^n)), j(\Gamma^2(\mathbb{R}^n)))$".

Finally, we highlight another property associated to the family $\rho(\Delta(\mathbb{R}^n))$ that, as the one in (c), involves not only the groups but also the semigroups of the family of $S$-groups $\mathbb{C}(\Gamma(\mathbb{R}^n))$. Here too, as we did in (c), we omit its proof.

- Let $f \in C(\Omega)$, $\partial^\alpha_\Omega \in \partial(\Omega)$, $\Gamma(\Omega) = \Gamma(\mathbb{R}^n) \cap \mathcal{P}(\Omega)$ (an $\Omega$-indexer) and $\overline{\Gamma}(\Omega) \subseteq \Gamma(\Omega)$ a cover of $\Omega$, that is, $\overline{\Gamma}(\Omega)$ is such that $\bigcup_{\Omega' \in \overline{\Gamma}(\Omega)} \Omega' = \Omega$. If $\rho_{(\Omega',\Omega)}(f) \in C^{|\alpha|}(\Omega')$ for every $\Omega' \in \overline{\Gamma}(\Omega)$, then, $f \in C^{|\alpha|}(\Omega)$ (the domain of $\partial^\alpha_\Omega$).

In particular, it follows that:

**(d)** For $f \in C(\Omega)$, $\partial^\alpha_\Omega \in \partial(\Omega)$ and $\overline{\Gamma}(\Omega) \subseteq \Gamma(\Omega) = \Gamma(\mathbb{R}^n) \cap \mathcal{P}(\Omega)$ a cover of $\Omega$ one has: if $\partial^\alpha_{\Omega'}(\rho_{(\Omega',\Omega)}(f)) = 0$ for every $\Omega' \in \overline{\Gamma}(\Omega)$, then, $f \in C^{|\alpha|}(\Omega)$.

Taking into account that $\partial^\alpha_{\Omega'} = j_{(\Omega',\Omega)}(\partial^\alpha_\Omega)$, an equivalent version of this property, which explicits the participation of the bonding family $j(\Gamma^2(\mathbb{R}^n))$, establishes (with the same hypothesis) that:

If $j_{(\Omega',\Omega)}(\partial^\alpha_\Omega)(\rho_{(\Omega',\Omega)}(f)) = 0$ for every $\Omega' \in \overline{\Gamma}(\Omega)$, then, $f$ belongs to the domain of $\partial^\alpha_\Omega$.

The considerations of this item, specially the properties (a) to (d), motivate and suggest the definition presented below.

## 3.13 Definition

Let $\Gamma(I)$ be an $I$-indexer and $\Delta(I)$ the set defined by:

$$\Delta(I) \coloneqq \Big\{(\gamma', \gamma) \in \Gamma(I) \times \Gamma(I) \colon \gamma' \subseteq \gamma\Big\}.$$

Let also $(\mathbb{G}(\Gamma(I)), i(\Gamma^2(I)))$, where

$$\mathbb{G}\Big(\Gamma(I)\Big) = \Big\{\mathbb{G}(\gamma) = \Big(G(\gamma), H(\gamma)\Big)\Big\}_{\gamma \in \Gamma(I)}$$

is a family of $S$-groups and

$$i\Big(\Gamma^2(I)\Big) = \Big\{i_{(\gamma',\gamma)}\Big\}_{(\gamma',\gamma)\in\Gamma^2(I)}$$

a bonding of the family $\mathbb{G}\Big(\Gamma(I)\Big)$. We will say that a family indexed by the set $\Delta(I)$, let us say

$$\Theta\Big(\Delta(I)\Big) = \Big\{\Theta_{(\gamma',\gamma)}\Big\}_{(\gamma',\gamma)\in\Delta(I)},$$

is a **restriction** for the bonded family $(\mathbb{G}(\Gamma(I)), i(\Gamma^2(I)))$ if and only if:

**(a)** for each $(\gamma', \gamma) \in \Delta(I)$, $\Theta_{(\gamma',\gamma)}$ is a homomorphism from the group $G(\gamma)$ into the group $G(\gamma')$ ($\Theta_{(\gamma',\gamma)} \colon G(\gamma) \longrightarrow G(\gamma')$ such that $\Theta_{(\gamma',\gamma)}(g+g') = \Theta_{(\gamma',\gamma)}(g) + \Theta_{(\gamma',\gamma)}(g')$ for every $g, g' \in G(\gamma)$);

**(b)** $\Theta_{(\gamma'',\gamma')}\Big(\Theta_{(\gamma',\gamma)}(g)\Big) = \Theta_{(\gamma'',\gamma)}(g)$ for every $g \in G(\gamma)$ and every $\gamma, \gamma', \gamma'' \in \Gamma(I)$ such that $\gamma'' \subseteq \gamma' \subseteq \gamma$;

**(c)** for every $\Phi \in H(\gamma)$ and every $g \in G(\gamma)_\Phi$[16],

$$\Theta_{(\gamma',\gamma)}\Big(\Phi(g)\Big) = i_{(\gamma',\gamma)}(\Phi)\Big(\Theta_{(\gamma',\gamma)}(g)\Big)$$

for every $(\gamma', \gamma) \in \Delta(I)$;

**(d)** for $\gamma \in \Gamma(I)$, $\Gamma(\gamma) = \Gamma(I) \cap \mathcal{P}(\gamma)$ (a $\gamma$-indexer), $\overline{\Gamma}(\gamma) \subseteq \Gamma(\gamma)$ a cover of $\gamma$, $g \in G(\gamma)$ and $\Phi \in H(\gamma)$, one has that: if

$$i_{(\gamma',\gamma)}(\Phi)\Big(\Theta_{(\gamma',\gamma)}(g)\Big) = 0$$

for every $\gamma' \in \overline{\Gamma}(\gamma)$, then, $g \in G(\gamma)_\Phi$.

[16] $G(\gamma)_\Phi$, according to the notation introduced in Definition 1.4, denotes the domain of the homomorphism $\Phi$.

# $S$-Space and $S$-Subspace

## 3.14 Remark

In possession of the definitions provided for the concepts of indexer, family of $S$-groups, bonding, bonded family, and restriction, we can now formulate a definition that makes precise the notion of "abstract spaces $\mathscr{G}(I)$ and $\widehat{\mathscr{G}}(I)$", informally introduced when we formulated the abstract version of our work agenda in 3.4, to which, since then, we have frequently been referring to.

## 3.15 Definition

We will say that the triplet

$$\mathscr{G}(I) = \Big(\mathbb{G}\big(\Gamma(I)\big), i\big(\Gamma^2(I)\big), \Theta\big(\Delta(I)\big)\Big)$$

is an $S$-**space** if and only if:

**(a)** $\Gamma(I)$ is an $I$-indexer;

**(b)** $\Delta(I) = \Big\{(\gamma', \gamma) \in \Gamma(I) \times \Gamma(I) \colon \gamma' \subseteq \gamma\Big\}$;

**(c)** $\mathbb{G}\big(\Gamma(I)\big)$ is a family of $S$-groups indexed by $\Gamma(I)$;

**(d)** $i\big(\Gamma^2(I)\big)$ is a bonding of the family $\mathbb{G}\big(\Gamma(I)\big)$ and, therefore, $\big(\mathbb{G}(\Gamma(I)), i(\Gamma^2(I))\big)$ is a bonded family;

**(e)** $\Theta\big(\Delta(I)\big)$ is a restriction for the bonded family $\big(\mathbb{G}(\Gamma(I)), i(\Gamma^2(I))\big)$.

We will say that a $S$-space is **abelian**, **surjective**, or **with identity** if and only if its family of $S$-groups is abelian, surjective, or with identity, respectively.

## 3.16 The $S$-Space $\mathscr{C}(\mathbb{R}^n)$

As we have already seen, the set

$$\Gamma(\mathbb{R}^n) = \Big\{\Omega \subseteq \mathbb{R}^n \colon \Omega \text{ is an open set of } \mathbb{R}^n\Big\}$$

is an indexer. Beyond that, the properties described in 3.10 of the family

$$j\Big(\Gamma^2(\mathbb{R}^n)\Big) = \Big\{j_{(\Omega',\Omega)}\Big\}_{(\Omega',\Omega)\in\Gamma^2(\mathbb{R}^n)},$$

where, for each $(\Omega', \Omega) \in \Gamma^2(\mathbb{R}^n)$,

$$j_{(\Omega',\Omega)} \colon \partial(\Omega) \longrightarrow \partial(\Omega')$$

is the function defined by

$$j_{(\Omega',\Omega)}(\partial_{\Omega}^{\alpha}) = \partial_{\Omega'}^{\alpha} \quad \text{for every} \quad \partial_{\Omega}^{\alpha} \in \partial(\Omega),$$

allow us to affirm that this family is a bonding of the family of $S$-groups of the continuous functions,

$$\mathbb{C}\Big(\Gamma(\mathbb{R}^n)\Big) = \Big\{\mathbb{C}(\Omega) = \Big(C(\Omega), \partial(\Omega)\Big)\Big\}_{\Omega \in \Gamma(\mathbb{R}^n)}.$$

We can also affirm that the family

$$\rho\Big(\Delta(\mathbb{R}^n)\Big) = \Big\{\rho_{(\Omega',\Omega)}\Big\}_{(\Omega',\Omega) \in \Delta(\mathbb{R}^n)},$$

where, for each $(\Omega', \Omega) \in \Delta(\mathbb{R}^n)$,

$$\rho_{(\Omega',\Omega)} \colon C(\Omega) \longrightarrow C(\Omega')$$

is the function defined by

$$\rho_{(\Omega',\Omega)}(f) = f|_{\Omega'} \quad \text{for every} \quad f \in C(\Omega),$$

being $f|_{\Omega'}$ the usual restriction of the function $f \colon \Omega \longrightarrow C$ to the subset $\Omega' \subseteq \Omega$ of its domain $\Omega$, is, taking into account the properties given in (a) to (d) of item 3.12, a restriction for the bonded family

$$\left(\mathbb{C}\Big(\Gamma(\mathbb{R}^n)\Big), j\Big(\Gamma^2(\mathbb{R}^n)\Big)\right).$$

Hence, the triplet $\mathscr{C}(\mathbb{R}^n)$ defined by

$$\mathscr{C}(\mathbb{R}^n) \coloneqq \left(\mathbb{C}\Big(\Gamma(\mathbb{R}^n)\Big), j\Big(\Gamma^2(\mathbb{R}^n)\Big), \rho\Big(\Delta(\mathbb{R}^n)\Big)\right)$$

is, according to Definition 3.15, an $S$-space, to which we will refer as **$S$-space of the continuous functions**, which is abelian, surjective, and with identity.

## 3.17 $S$-Subspace

Let us remember that, being $\Gamma(I)$ an indexer, for each $\gamma \in \Gamma(I)$ one has that:

$$\Gamma(\gamma) = \Gamma(I) \cap \mathcal{P}(\gamma)$$

is also an indexer. Furthermore, taking into account the definition of subfamilies of a family $\mathbb{G}(\Gamma(I))$ of $S$-groups, the subfamilies of $\mathbb{G}(\Gamma(I))$ are, exactly, the families $\mathbb{G}(\Gamma(\gamma)) \subseteq \mathbb{G}(\Gamma(I))$ with $\gamma \in \Gamma(I)$.

Let now

$$\mathscr{G}(I) = \Big(\mathbb{G}\big(\Gamma(I)\big), i\big(\Gamma^2(I)\big), \Theta\big(\Delta(I)\big)\Big)$$

be a $S$-space and $\gamma \in \Gamma(I)$ be arbitrarily fixed. It clearly results that, with

$$\Gamma^2(\gamma) = \Gamma(\gamma) \times \Gamma(\gamma) \quad \text{and} \quad \Delta(\gamma) = \Big\{(\gamma'', \gamma') \in \Gamma^2(\gamma) \colon \gamma'' \subseteq \gamma'\Big\},$$

$$\mathscr{G}(\gamma) := \Big(\mathbb{G}\big(\Gamma(\gamma)\big), i\big(\Gamma^2(\gamma)\big), \Theta\big(\Delta(\gamma)\big)\Big)$$

also is a $S$-space that, in its relationship with $\mathscr{G}(I)$, will be said to be an **$S$-subspace** of $\mathscr{G}(I)$. Unless otherwise explicitly said, the term $S$-subspace will be thus understood. Hence, the $S$-subspaces of a given $S$-space $\mathscr{G}(I)$ as above (with indexer $\Gamma(I)$) are, exactly, the $S$-spaces $\mathscr{G}(\gamma)$ with $\gamma \in \Gamma(I)$. It results that, if $\mathscr{G}(I)$ is abelian, surjective or with identity, then, the $S$-subspaces of $\mathscr{G}(I)$ also are abelian, surjective or with identity, respectively.

# Coherent $S$-Space

## 3.18 The Property of Coherency in $\mathscr{C}(\mathbb{R}^n)$

Let us resume to the $S$-space of continuous functions,

$$\mathscr{C}(\mathbb{R}^n) = \Big(\mathbb{C}\big(\Gamma(\mathbb{R}^n)\big), j\big(\Gamma^2(\mathbb{R}^n)\big), \rho\big(\Delta(\mathbb{R}^n)\big)\Big),$$

described in 3.16. Consider, for a $\Omega \in \Gamma(\mathbb{R}^n)$, $\Omega \neq \varnothing$, arbitrarily fixed, the corresponding $S$-subspace of $\mathscr{C}(\mathbb{R}^n)$, namely:

$$\mathscr{C}(\Omega) = \Big(\mathbb{C}\big(\Gamma(\Omega)\big), j\big(\Gamma^2(\Omega)\big), \rho\big(\Delta(\Omega)\big)\Big).$$

Let now $\overline{\Gamma}(\Omega) \subseteq \Gamma(\Omega)$ be a cover of $\Omega$ and $f(\overline{\Gamma}(\Omega)) = \{f_{\Omega'}\}_{\Omega' \in \overline{\Gamma}(\Omega)}$ be a family (indexed by $\overline{\Gamma}(\Omega)$) whose members, $f_{\Omega'}$, are elements of the corresponding groups, $C(\Omega')$, of the $S$-groups $\mathbb{C}(\Omega') = (C(\Omega'), \partial(\Omega'))$ of the family $\mathbb{C}(\Gamma(\Omega))$, that is, $f_{\Omega'} \in C(\Omega')$ for each $\Omega' \in \overline{\Gamma}(\Omega)$. Suppose now that the family $f(\overline{\Gamma}(\Omega))$ is such that

$$\rho_{(\Omega' \cap \Omega'', \Omega')}(f_{\Omega'}) = \rho_{(\Omega' \cap \Omega'', \Omega'')}(f_{\Omega''})$$

for every $\Omega', \Omega'' \in \overline{\Gamma}(\Omega)$ such that $\Omega' \cap \Omega'' \neq \varnothing$. In this case, there exists an unique function $f \in C(\Omega)$ such that

$$\rho_{(\Omega', \Omega)}(f) = f_{\Omega'} \quad \text{for every} \quad \Omega' \in \overline{\Gamma}(\Omega).$$

In fact, let $f$ be defined by:

$$\begin{aligned} f : \Omega &\longrightarrow C \\ x &\longmapsto f(x) := f_{\Omega'}(x) \quad \text{if} \quad x \in \Omega' \in \overline{\Gamma}(\Omega). \end{aligned}$$

As $\overline{\Gamma}(\Omega)$ is a cover of $\Omega$, then, for each $x \in \Omega$ there exists $\Omega' \in \overline{\Gamma}(\Omega)$ such that $x \in \Omega'$. Case exists $\Omega'' \in \overline{\Gamma}(\Omega)$, $\Omega'' \neq \Omega'$, such that $x \in \Omega''$, then, by the definition of $f$, $f(x) = f_{\Omega'}(x)$ and $f(x) = f_{\Omega''}(x)$. However, as

$$\rho_{(\Omega' \cap \Omega'', \Omega')}(f_{\Omega'}) = \rho_{(\Omega' \cap \Omega'', \Omega'')}(f_{\Omega''}),$$

we have $f_{\Omega'}(x) = f_{\Omega''}(x)$, which shows us that $f$ is well-defined. Now, it is clear from the definition of $f$ that $f \in C(\Omega)$ and $\rho_{(\Omega', \Omega)}(f) = f_{\Omega'}$ for every $\Omega' \in \overline{\Gamma}(\Omega)$; furthermore, $f$ is the only such function, for if $g \in C(\Omega)$ is such that $\rho_{(\Omega', \Omega)}(g) = f_{\Omega'}$ for every $\Omega' \in \overline{\Gamma}(\Omega)$, then, if $x \in \Omega$, $x \in \Omega'$ for some $\Omega' \in \overline{\Gamma}(\Omega)$ and, hence, $g(x) = f_{\Omega'}(x) = f(x)$, that is, $g = f$.

If now we baptize the families

$$f\Big(\overline{\Gamma}(\Omega)\Big) = \Big\{f_{\Omega'}\Big\}_{\Omega' \in \overline{\Gamma}(\Omega)},$$

where $\overline{\Gamma}(\Omega) \subseteq \Gamma(\Omega)$ is a cover of $\Omega$, $f_{\Omega'} \in C(\Omega')$ for each $\Omega' \in \overline{\Gamma}(\Omega)$ and

$$\rho_{(\Omega' \cap \Omega'', \Omega')}(f_{\Omega'}) = \rho_{(\Omega' \cap \Omega'', \Omega'')}(f_{\Omega''})$$

for every $\Omega', \Omega'' \in \overline{\Gamma}(\Omega)$ such that $\Omega' \cap \Omega'' \neq \varnothing$, of **coherent families in** $\mathscr{C}(\Omega)$, then, what we proved above can be expressed in the following form:

- Let $\mathscr{C}(\Omega)$ be a $S$-subspace of $\mathscr{C}(\mathbb{R}^n)$ and $\{f_{\Omega'}\}_{\Omega' \in \overline{\Gamma}(\Omega)}$ be a coherent family in $\mathscr{C}(\Omega)$, both arbitrarily fixed. There exists a single $f \in C(\Omega)$ such that

  $$\rho_{(\Omega', \Omega)}(f) = f_{\Omega'}$$

  for each $\Omega' \in \overline{\Gamma}(\Omega)$.

When referring to this property of $\mathscr{C}(\mathbb{R}^n)$, we will say that $\mathscr{C}(\mathbb{R}^n)$ is a **coherent** $S$**-space**.

It is convenient to extend the concept of coherency, defining it relative to any $S$-space.

## 3.19 Definition

Let

$$\mathscr{G}(I) = \left(\mathbb{G}\Big(\Gamma(I)\Big) = \Big\{\mathbb{G}(\gamma) = \Big(G(\gamma), H(\gamma)\Big)\Big\}_{\gamma \in \Gamma(I)}, i\Big(\Gamma^2(I)\Big), \Theta\Big(\Delta(I)\Big)\right)$$

be a $S$-space and $\mathscr{G}(\gamma)$ be a $S$-subspace of $\mathscr{G}(I)$, both arbitrarily fixed.

**(a)** We will say that a family $\{g_{\gamma'}\}_{\gamma' \in \overline{\Gamma}(\gamma)}$ is **coherent in** $\mathscr{G}(\gamma)$ if and only if,

**(a-1)** $\overline{\Gamma}(\gamma) \subseteq \Gamma(\gamma)$ and $\overline{\Gamma}(\gamma)$ is a cover of $\gamma$;

**(a-2)** $g_{\gamma'} \in G(\gamma')$ for each $\gamma' \in \overline{\Gamma}(\gamma)$;

**(a-3)** $\Theta_{(\gamma' \cap \gamma'', \gamma')}(g_{\gamma'}) = \Theta_{(\gamma' \cap \gamma'', \gamma'')}(g_{\gamma''})$ for any $\gamma', \gamma'' \in \overline{\Gamma}(\gamma)$ such that $\gamma' \cap \gamma'' \neq \varnothing$.

**(b)** We will say that the $S$-space $\mathscr{G}(I)$ is **coherent** if and if, for each coherent family $\{g_{\gamma'}\}_{\gamma' \in \overline{\Gamma}(\Omega)}$ in $\mathscr{G}(\gamma)$, there exists a single $g \in G(\gamma)$ such that

$$\Theta_{(\gamma', \gamma)}(g) = g_{\gamma'} \quad \text{for every} \quad \gamma' \in \overline{\Gamma}(\gamma).$$

## 3.20 Notation

Let $\Gamma(I)$ be an $I$-indexer. Covers of the set $I$ included in $\Gamma(I)$ will be systematically denoted by greek letters, capitalized or not, overlined, "$^{-}$", and followed by the symbol "$I$" between parenthesis. Hence, $\overline{\Gamma}(I)$, $\overline{\xi}(I)$, $\overline{\chi}(I)$, for example, will denote subsets of $\Gamma(I)$ such that $\bigcup \overline{\Gamma}(I) = \bigcup \overline{\xi}(I) = \bigcup \overline{\chi}(I) = I$.

Since (recap Definition 3.6), if $\gamma \in \Gamma(I)$ and $\gamma \neq \varnothing$, $\Gamma(\gamma) = \Gamma(I) \cap \mathcal{P}(\gamma)$ is also an indexer (more precisely a $\gamma$-indexer), symbols such as, for instance, $\overline{\eta}(\gamma)$ or $\overline{\nu}(\gamma)$, will denote, in line with the notation above introduced, covers of $\gamma$ included in $\Gamma(\gamma)$.

Generic members of covers as, for instance, $\overline{\xi}(\gamma) \subseteq \Gamma(\gamma)$ or $\overline{\nu}(I) \subseteq \Gamma(I)$, will be, generally but not necessarily, denoted by the same greek letters (without the overline) used in the symbols of the covers: $\xi \in \overline{\xi}(\gamma)$ and $\nu \in \overline{\nu}(I)$. Beyond the notation commonly employed to represent a family indexed by the members of a set, let us say $\overline{\xi}(\gamma)$, as, for example,

$$\left\{g_\xi\right\}_{\xi \in \overline{\xi}(\gamma)},$$

we will also use, frequently, the simpler notation, $g(\overline{\xi}(\gamma))$. Hence,

$$g\Big(\overline{\xi}(\gamma)\Big) = \left\{g_\xi\right\}_{\xi \in \overline{\xi}(\gamma)}.$$

# Bonded Family Extension

## 3.21 Bonding Extension

Let $\mathbb{G}(\Gamma(I)) = \{\mathbb{G}(\gamma) = (G(\gamma), H(\gamma))\}_{\gamma \in \Gamma(I)}$ be a family of $S$-groups and let also $\widehat{\mathbb{G}}(\Gamma(I)) = \{\widehat{\mathbb{G}}(\gamma) = (\widehat{G}(\gamma), \widehat{H}(\gamma))\}_{\gamma \in \Gamma(I)}$ be an extension of $\mathbb{G}(\Gamma(I))$ (see Definition 3.8). Hence, for each $\gamma \in \Gamma(I)$, the $S$-group $\widehat{\mathbb{G}}(\gamma) = (\widehat{G}(\gamma), \widehat{H}(\gamma))$ is an extension of the $S$-group $\mathbb{G}(\gamma) = (G(\gamma), H(\gamma))$ and, therefore, (see Definition 1.6) the semigroup $\widehat{H}(\gamma)$ is constituted, exactly, by prolongations to $\widehat{G}(\gamma)$ of each member $\Phi \in H(\gamma)$. More precisely (see Definition 1.6(c)), $\widehat{H}(\gamma)$ and $H(\gamma)$ are, for each $\gamma \in \Gamma(I)$, isomorphic semigroups and the function $\widehat{\gamma}\colon H(\gamma) \longrightarrow \widehat{H}(\gamma)$, as defined ahead, is an isomorphism:

$$\begin{aligned} \widehat{\gamma} : H(\gamma) &\longrightarrow \widehat{H}(\gamma) \\ \Phi &\longmapsto \widehat{\gamma}(\Phi) = \widehat{\Phi} \end{aligned}$$

being $\widehat{\Phi}$ the prolongation $(\in \widehat{H}(\gamma))$ of the homomorphism $\Phi(\in H(\gamma))$ to $\widehat{G}(\gamma)$.

Worth noting that the notation above employed for the image of $\Phi \in H(\gamma)$ by the isomorphism $\widehat{\gamma}$, namely, $\widehat{\gamma}(\Phi) = \widehat{\Phi}$, is consistent with the notation $\widehat{H}(\gamma)$ for the semigroup of the $S$-group $\widehat{\mathbb{G}}(\gamma) = (\widehat{G}(\gamma), \widehat{H}(\gamma))$, once this semigroup is, exactly, the image of $H(\gamma)$ by the isomorphism $\widehat{\gamma}$, that is, $\widehat{\gamma}(H(\gamma)) = \widehat{H}(\gamma)$.

Let us suppose now that $i(\Gamma^2(I)) = \{i_{(\gamma',\gamma)}\}_{(\gamma',\gamma)\in\Gamma^2(I)}$ is a bonding of the family of $S$-groups $\mathbb{G}(\Gamma(I))$ (see Definition 3.11). Hence, for each $(\gamma',\gamma) \in \Gamma^2(I)$, $i_{(\gamma',\gamma)}\colon H(\gamma) \longrightarrow H(\gamma')$ is an isomorphism from the semigroup $H(\gamma)$ onto the semigroup $H(\gamma')$ and, therefore, since

$$\begin{aligned} \widehat{\gamma} : H(\gamma) &\longrightarrow \widehat{H}(\gamma) \\ \Phi &\longmapsto \widehat{\gamma}(\Phi) = \widehat{\Phi}\Big(\text{prolongation of } \Phi \text{ to } \widehat{G}(\gamma)\Big) \end{aligned}$$

and

$$\begin{aligned} \widehat{\gamma'} : H(\gamma') &\longrightarrow \widehat{H}(\gamma') \\ \Psi &\longmapsto \widehat{\gamma'}(\Psi) = \widehat{\Psi}\Big(\text{prolongation of } \Psi \text{ to } \widehat{G}(\gamma')\Big) \end{aligned}$$

are isomorphisms from $H(\gamma)$ onto $\widehat{H}(\gamma)$ and from $H(\gamma')$ onto $\widehat{H}(\gamma')$, respectively, it results that the function $\widehat{i}_{(\gamma',\gamma)}\colon \widehat{H}(\gamma) \longrightarrow \widehat{H}(\gamma')$ defined by

$$\widehat{i}_{(\gamma',\gamma)} := \widehat{\gamma'} i_{(\gamma',\gamma)} (\widehat{\gamma})^{-1},$$

that is,

$$\widehat{i}_{(\gamma',\gamma)}(\widehat{\Phi}) := \widehat{\gamma'}\left( i_{(\gamma',\gamma)}\Big( (\widehat{\gamma})^{-1}(\widehat{\Phi})\Big)\right) \quad \text{for every} \quad \widehat{\Phi} \in \widehat{H}(\gamma),$$

is an isomorphism from the semigroup $\widehat{H}(\gamma)$ onto the semigroup $\widehat{H}(\gamma')$. Thus, the family $\widehat{i}(\Gamma^2(I))$ defined by

$$\widehat{i}\Big(\Gamma^2(I)\Big) := \Big\{\widehat{i}_{(\gamma',\gamma)} = \widehat{\gamma'} i_{(\gamma',\gamma)} (\widehat{\gamma})^{-1}\Big\}_{(\gamma',\gamma)\in\Gamma^2(I)},$$

is a family of isomorphisms between the semigroups of the extension $\widehat{\mathbb{G}}(\Gamma(I))$ of the family $\mathbb{G}(\Gamma(I))$. Furthermore, it results immediately that the family $\widehat{i}(\Gamma^2(I))$ is such that:

- $\Big(\widehat{i}_{(\gamma',\gamma)}\Big)^{-1} = \widehat{i}_{(\gamma,\gamma')}$ for every $(\gamma',\gamma) \in \Gamma^2(I)$ and
- $\widehat{i}_{(\gamma'',\gamma')}\widehat{i}_{(\gamma',\gamma)} = \widehat{i}_{(\gamma'',\gamma)}$ for every $\gamma, \gamma', \gamma'' \in \Gamma(I)$.

From this, and taking Definition 3.11 into account, we conclude that $\widehat{i}(\Gamma^2(I))$ is a bonding of the family of $S$-groups $\widehat{\mathbb{G}}(\Gamma(I))$ and, therefore, that $(\widehat{G}(\Gamma(I)), \widehat{i}(\Gamma^2(I)))$ is a bonded family. This result endorses and motivates the definition of the concept of extension of a bonded family, formulated below.

## 3.22 Definition

Let

$$\left(\mathbb{G}\Big(\Gamma(I)\Big) = \Big\{\mathbb{G}(\gamma) = \Big(G(\gamma), H(\gamma)\Big)\Big\}_{\gamma\in\Gamma(I)}, i\Big(\Gamma^2(I)\Big) = \Big\{i_{(\gamma',\gamma)}\Big\}_{(\gamma',\gamma)\in\Gamma^2(I)}\right)$$

be a bonded family of $S$-groups and

$$\widehat{\mathbb{G}}\Big(\Gamma(I)\Big) = \Big\{\widehat{\mathbb{G}}(\gamma) = \Big(\widehat{G}(\gamma), \widehat{H}(\gamma)\Big)\Big\}_{\gamma\in\Gamma(I)}$$

an extension of the family $\mathbb{G}(\Gamma(I))$. Let also the family

$$\widehat{i}\Big(\Gamma^2(I)\Big) := \Big\{\widehat{i}_{(\gamma',\gamma)}\Big\}_{(\gamma',\gamma)\in\Gamma^2(I)}$$

where, for each $(\gamma', \gamma) \in \Gamma^2(I)$, $\widehat{i}_{(\gamma',\gamma)}$ is the function defined by

$$\widehat{i}_{(\gamma',\gamma)} : \widehat{H}(\gamma) \longrightarrow \widehat{H}(\gamma')$$

$$\widehat{\Phi} \longmapsto \widehat{i}_{(\gamma',\gamma)}\Big(\widehat{\Phi}\Big) := \widehat{\gamma}'\left(i_{(\gamma',\gamma)}\left((\widehat{\gamma})^{-1}\Big(\widehat{\Phi}\Big)\right)\right),$$

being, for each $\delta \in \Gamma(I)$, $\widehat{\delta}$ the function given by

$$\widehat{\delta} : H(\delta) \longrightarrow \widehat{H}(\delta)$$

$$\Psi \longmapsto \widehat{\delta}(\Psi) := \widehat{\Psi}$$

where $\widehat{\Psi}$ is the prolongation $(\in \widehat{H}(\delta))$ of $\Psi(\in H(\delta))$ to $\widehat{G}(\delta)$. As we saw in 3.21, the family $\widehat{i}(\Gamma^2(I))$ is a bonding of the family $\widehat{\mathbb{G}}(\Gamma(I))$. We will say that $\widehat{i}(\Gamma^2(I))$ is the **extension to $\widehat{\mathbb{G}}(\Gamma(I))$ of the bonding** $i(\Gamma^2(I))$ of the bonded family

$$\left(\mathbb{G}\Big(\Gamma(I)\Big), i\Big(\Gamma^2(I)\Big)\right).$$

We will also say that a bonded family of $S$-groups, let us say

$$\left(\mathbb{E}\Big(\Gamma(I)\Big), k\Big(\Gamma^2(I)\Big)\right),$$

is an **extension of the bonded family**

$$\left(\mathbb{G}\Big(\Gamma(I)\Big), i\Big(\Gamma^2(I)\Big)\right),$$

if and only if, $\mathbb{E}(\Gamma(I))$ is an extension of the family $\mathbb{G}(\Gamma(I))$ and $k(\Gamma^2(I))$ is the extension to $\mathbb{E}(\Gamma(I))$ of the bonding $i(\Gamma^2(I))$.

# $S$-Space Extension

## 3.23 Definition

**(a)** Let $\mathscr{G}(I) = (\mathbb{G}(\Gamma(I)), i(\Gamma^2(I), \Theta(\Delta(I))))$ be a $S$-space, $(\widehat{\mathbb{G}}(\Gamma(I)), \widehat{i}(\Gamma^2(I)))$ be an extension of the bonded family $(\mathbb{G}(\Gamma(I)), i(\Gamma^2(I)))$ of $\mathscr{G}(I)$ and $\widehat{\Theta}(\Delta(I))$ be a restriction for $(\widehat{\mathbb{G}}(\Gamma(I)), \widehat{i}(\Gamma^2(I)))$. We will say that $\widehat{\Theta}(\Delta(I))$ is a **prolongation** of $\Theta(\Delta(I))$ to the bonded family $(\widehat{\mathbb{G}}(\Gamma(I)), \widehat{i}(\Gamma^2(I)))$ if and only if, for every $(\gamma', \gamma) \in \Delta(I)$,

$$\widehat{\Theta}_{(\gamma',\gamma)}(g) = \Theta_{(\gamma',\gamma)}(g)$$

whatever $g \in G(\gamma) \subseteq \widehat{G}(\gamma)$.

**(b)** We will say that a $S$-space

$$\mathscr{E}(I) = \Big(\mathbb{E}\big(\Gamma(I)\big), k\big(\Gamma^2(I)\big), \chi\big(\Delta(I)\big)\Big)$$

is an **extension** of the $S$-space

$$\mathscr{G}(I) = \Big(\mathbb{G}\big(\Gamma(I)\big), i\big(\Gamma^2(I)\big), \Theta\big(\Delta(I)\big)\Big),$$

if and only if,

**(b-1)** the bonded family $(\mathbb{E}(\Gamma(I)), k(\Gamma^2(I)))$ is an extension of the bonded family $(\mathbb{G}(\Gamma(I)), i(\Gamma^2(I)))$ and

**(b-2)** $\chi(\Delta(I))$ is a prolongation of the restriction $\Theta(\Delta(I))$ to the bonded family $(\mathbb{E}(\Gamma(I)), k(\Gamma^2(I)))$.

**(c)** Let

$$\widehat{\mathscr{G}}(I) = \Big(\widehat{\mathbb{G}}\big(\Gamma(I)\big) = \Big\{\widehat{\mathbb{G}}(\gamma) = \big(\widehat{G}(\gamma), \widehat{H}(\gamma)\big)\Big\}_{\gamma\in\Gamma(I)}, \widehat{i}\big(\Gamma^2(I)\big), \widehat{\Theta}\big(\Delta(I)\big)\Big)$$

be an extension of the $S$-space

$$\mathscr{G}(I) = \Big(\mathbb{G}(\Gamma(I)) = \Big\{\mathbb{G}(\gamma) = \big(G(\gamma), H(\gamma)\big)\Big\}_{\gamma\in\Gamma(I)}, i\big(\Gamma^2(I)\big), \Theta\big(\Delta(I)\big)\Big).$$

We will say that:

**(c-1)** $\widehat{\mathscr{G}}(I)$ is a **strict** (respectively, **closed**) **extension** of $\mathscr{G}(I)$ if and only if $\widehat{\mathbb{G}}(\Gamma(I))$ is a strict (respectively, closed) extension of $\mathbb{G}(\Gamma(I))$;

**(c-2)** $\widehat{\mathscr{G}}(I)$ is a **locally closed extension** of $\mathscr{G}(I)$ if and only if, for each $\gamma \in \Gamma(I)$, if $\widehat{g} \in \widehat{G}(\gamma)$ and if $x \in \gamma$, there exists $\gamma' \in \Gamma(\gamma) = \Gamma(I) \cap \mathcal{P}(\gamma)$ such that $x \in \gamma'$ and

$$\widehat{\Theta}_{(\gamma',\gamma)}(\widehat{g}) \in G(\gamma');$$

**(c-3)** $\widehat{\mathscr{G}}(I)$ is a **coherent extension** of $\mathscr{G}(I)$ if and only if $\widehat{\mathscr{G}}(I)$ is a coherent $S$-space.

## 3.24 Warning

When affirming that a family of $S$-groups,

$$\widehat{\mathbb{G}}\Big(\Gamma(I)\Big) = \Big\{\widehat{\mathbb{G}}(\gamma) = \Big(\widehat{G}(\gamma), \widehat{H}(\gamma)\Big)\Big\}_{\gamma\in\Gamma(I)},$$

is an extension of the family

$$\mathbb{G}\Big(\Gamma(I)\Big) = \Big\{\mathbb{G}(\gamma) = (G(\gamma), H(\gamma))\Big\}_{\gamma\in\Gamma(I)},$$

we are saying, according to Definition 3.8, that, for each $\gamma \in \Gamma(I)$, the $S$-group $\widehat{\mathbb{G}}(\gamma) = (\widehat{G}(\gamma), \widehat{H}(\gamma))$ is an extension of the $S$-group $\mathbb{G}(\gamma) = (G(\gamma), H(\gamma))$ which, in turn, with Definition 1.6(d) in mind, means, among other things, that the group $\widehat{G}(\gamma)$ $\delta_\gamma$-admits the group $G(\gamma)$ as a subgroup, that is, that exists an injective homomorphism, $\delta_\gamma\colon G(\gamma) \longrightarrow \widehat{G}(\gamma)$, from $G(\gamma)$ into $\widehat{G}(\gamma)$ or, equivalently, that the group $G(\gamma)$ is isomorphic to the subgroup $\delta_\gamma(G(\gamma)) = \{\delta_\gamma(g)\colon g \in G(\gamma)\}$ of $\widehat{G}(\gamma)$. Hence, through $\delta_\gamma$ we identify the members $g \in G(\gamma)$ with the elements $\delta_\gamma(g) \in \widehat{G}(\gamma)$, as well as the group $G(\gamma)$ with the subgroup $\delta_\gamma(G(\gamma))$ of $\widehat{G}(\gamma)$. From these identifications originate the abuse of language referred in Definition 1.6(b) and Remark 1.7 (recap these itens), that we have, since then, used. Hence, for example, in Definition 3.23(a) the condition (necessary and sufficient) for the family $\widehat{\Theta}(\Delta(I))$ to be a prolongation of $\Theta(\Delta(I))$, when expressed without the referred abuse of language, assumes the form

$$\widehat{\Theta}_{(\gamma',\gamma)}\Big(\delta_\gamma(g)\Big) = \delta_{\gamma'}\Big(\Theta_{(\gamma',\gamma)}(g)\Big) \quad \text{for every} \quad g \in G(\gamma).$$

By the same token, in item (c-2) of the same Definition 3.23, instead of

$$\widehat{\Theta}_{(\gamma',\gamma)}(\widehat{g}) \in G(\gamma'),$$

we have, in the version without the abuse of language,

$$\widehat{\Theta}_{(\gamma',\gamma)}(\widehat{g}) \in \delta_{\gamma'}\Big(G(\gamma')\Big).$$

Clearly, the referred abuse of language is undone if the groups $G(\gamma)$ of the $S$-groups $\mathbb{G}(\gamma) = (G(\gamma), H(\gamma))$ were, in fact, subgroups of the groups $\widehat{G}(\gamma)$ of the extensions $\widehat{\mathbb{G}}(\gamma) = (\widehat{G}(\gamma), \widehat{H}(\gamma))$; that is, $G(\gamma)$ is set-like included in $\widehat{G}(\gamma)$, $G(\gamma) \subseteq \widehat{G}(\gamma)$, and its addition is the addition of the group $\widehat{G}(\gamma)$ restricted to its elements. In this case, the injective homomorphisms $\delta_\gamma \colon G(\gamma) \longrightarrow \widehat{G}(\gamma)$ will be the identities in $G(\gamma)$, that is, for each $\gamma \in \Gamma(I)$,

$$\begin{aligned} \delta_\gamma = I_{G(\gamma)} : G(\gamma) &\longrightarrow \widehat{G}(\gamma) \\ g &\longmapsto \delta_\gamma(g) := g \end{aligned}$$

As we proceeded so far, we will not prevent ourselves from this abuse of language in what follows. Hence, we will admit that everything happens as if the injective homomorphisms $\delta_\gamma \colon G(\gamma) \longrightarrow \widehat{G}(\gamma)$ were the identities in $G(\gamma)$, $\delta_\gamma = I_{G(\gamma)}$; when convenient or if necessary we will proceed to the proper insertions of $\delta_\gamma$ to eliminate the abuse.

## 3.25 Revisited Goals

Reviewing now the considerations of items 3.1 to 3.4 which, informally, introduced the "spaces" $\mathscr{C}(\mathbb{R}^n)$, $\mathscr{D}'_f(\mathbb{R}^n)$ and $\mathscr{D}'(\mathbb{R}^n)$ and motivated the elaboration of a work agenda, described in item 3.3, and its abstract version, presented in 3.4, we can affirm that we are now, given the definitions presented so far, able to formulate a precise, rigorous, version for the questions informally introduced there. In fact, essentially, the questions loosely proposed in 3.4 assume, in terms of the referred definitions, the following precise, rigorous, formulation.

Given a $S$-space, $\mathscr{G}(I)$, abelian, surjective, and with identity, one asks:

**(a)** Is there a strict and closed extension, $\widetilde{\mathscr{G}}(I)$, of $\mathscr{G}(I)$?

**(b)** In the hypothesis of existence of $\widetilde{\mathscr{G}}(I)$ as in (a), is there a locally closed extension, $\overline{\mathscr{G}}(I)$, of $\widetilde{\mathscr{G}}(I)$?

**(c)** If exist $\widetilde{\mathscr{G}}(I)$ and $\overline{\mathscr{G}}(I)$ as in (a) and (b), respectively, are they unique?

In reality, the question (c) above is still informal in nature, that is, it is not backed by the definitions so far formulated, since we have not yet defined the notion of isomorphism for $S$-spaces, necessary to a characterization of uniqueness. This gap will be filled shortly.

The above questions will be considered in detail in the next chapter. For now, it is opportune to approach here some few aspects of the questions (a) and (b), for what we dedicate the next item.

## 3.26 About the Questions (a) and (b) in 3.25

Let us consider a $S$-space, say

$$\mathscr{G}(I)=\Bigg(\mathbb{G}\Big(\Gamma(I)\Big)=\Big\{\mathbb{G}(\gamma)=\Big(G(\gamma),H(\gamma)\Big)\Big\}_{\gamma\in\Gamma(I)},\; i\Big(\Gamma^2(I)\Big)=\Big\{i_{(\gamma',\gamma)}\Big\}_{(\gamma',\gamma)\in\Gamma^2(I)},\; \Theta\Big(\Delta(I)\Big)\Bigg),$$

which is abelian, surjective, and with identity, as, for instance, the $S$-space of the continuous functions

$$\mathscr{C}(\mathbb{R}^n)=\Bigg(\mathbb{C}\Big(\Gamma(\mathbb{R}^n)\Big),j\Big(\Gamma^2(\mathbb{R}^n)\Big),\rho\Big(\Delta(\mathbb{R}^n)\Big)\Bigg)$$

described in item 3.16. By the Theorem of Extension of $S$-groups (Theorem 2.16), we have that the family of $S$-groups

$$\widetilde{\mathbb{G}}\Big(\Gamma(I)\Big)=\Big\{\widetilde{\mathbb{G}}(\gamma)=\Big(\widetilde{G}(\gamma),\widetilde{H}(\gamma)\Big)\Big\}_{\gamma\in\Gamma(I)}$$

where, for each $\gamma\in\Gamma(I)$, the $S$-group

$$\widetilde{\mathbb{G}}(\gamma)=\Big(\widetilde{G}(\gamma),\widetilde{H}(\gamma)\Big)$$

is the strict and closed extension of the $S$-group

$$\mathbb{G}(\gamma)=\Big(G(\gamma),H(\gamma)\Big),$$

as defined in the referred theorem, which is unique unless isomorphism, is, taking into account Definition 3.8, a strict and closed extension of the family of $S$-groups $\mathbb{G}(\Gamma(I))$ of $\mathscr{G}(I)$. Now, since $\widetilde{\mathbb{G}}(\Gamma(I))$ is an extension of the family of $S$-groups $\mathbb{G}(\Gamma(I))$ of $\mathscr{G}(I)$, we have, as establishes the Definition 3.22, that the family

$$\widetilde{i}\Big(\Gamma^2(I)\Big)=\Big\{\widetilde{i}_{(\gamma',\gamma)}\Big\}_{(\gamma',\gamma)\in\Gamma^2(I)}$$

where, for each $(\gamma',\gamma)\in\Gamma^2(I)$, $\widetilde{i}_{(\gamma',\gamma)}$ is the function

$$\begin{aligned}\widetilde{i}_{(\gamma',\gamma)}:\widetilde{H}(\gamma)&\longrightarrow\widetilde{H}(\gamma')\\ \widetilde{\Phi}&\longmapsto\widetilde{i}_{(\gamma',\gamma)}(\widetilde{\Phi}):=\widetilde{\gamma}'\Bigg(i_{(\gamma',\gamma)}\Big((\widetilde{\gamma})^{-1}(\widetilde{\Phi})\Big)\Bigg)\end{aligned}$$

being $\widetilde{\delta}$, for each $\delta\in\Gamma(I)$, the function

$$\begin{aligned}\widetilde{\delta}:H(\delta)&\longrightarrow\widetilde{H}(\delta)\\ \Psi&\longmapsto\widetilde{\delta}(\Psi)=\widetilde{\Psi}\Big(\text{the prolongation of }\Psi\text{ to }\widetilde{G}(\delta)\Big),\end{aligned}$$

is the extension of the bonding $i(\Gamma^2(I))$ (of $\mathscr{G}(I)$) to $\widetilde{\mathbb{G}}(\Gamma(I))$ and, therefore, by the same Definition 3.22,

$$\left(\widetilde{G}\Big(\Gamma(I)\Big),\widetilde{i}\Big(\Gamma^2(I)\Big)\right)$$

is an extension of the bonded family

$$\left(G\Big(\Gamma(I)\Big),i\Big(\Gamma^2(I)\Big)\right)$$

of the $S$-space $\mathscr{G}(I)$.

From this, and taking into account the definition of $S$-space extension (Definition 3.23), we can conclude that the triplet

$$\widetilde{\mathscr{G}}(I) := \left(\widetilde{\mathbb{G}}\Big(\Gamma(I)\Big),\widetilde{i}\Big(\Gamma^2(I)\Big),\widetilde{\Theta}\Big(\Delta(I)\Big)\right)$$

is a strict and closed extension of the $S$-space (abelian, surjective, and with identity)

$$\mathscr{G}(I) = \left(\mathbb{G}\Big(\Gamma(I)\Big),i\Big(\Gamma^2(I)\Big),\Theta\Big(\Delta(I)\Big)\right),$$

if and only if, $\widetilde{\Theta}\Big(\Delta(I)\Big)$ is a restriction for the bonded family

$$\left(\widetilde{\mathbb{G}}\Big(\Gamma(I)\Big),\widetilde{i}\Big(\Gamma^2(I)\Big)\right)$$

and, also, a prolongation of $\Theta(\Delta(I))$ to this bonded family.

We will see in the next chapter that there is only one such family $\widetilde{\Theta}(\Delta(I))$ and that, according to a sensible definition (given ahead in the present chapter) of the concept of isomorphism for $S$-spaces, the $S$-space

$$\widetilde{\mathscr{G}}(I) = \left(\widetilde{\mathbb{G}}\Big(\Gamma(I)\Big),\widetilde{i}\Big(\Gamma^2(I)\Big),\widetilde{\Theta}\Big(\Delta(I)\Big)\right)$$

is, unless isomorphism, the only strict and closed extension of the abelian and surjective $S$-space with identity

$$\mathscr{G}(I) = \left(\mathbb{G}\Big(\Gamma(I)\Big),i\Big(\Gamma^2(I)\Big),\Theta\Big(\Delta(I)\Big)\right),$$

that answers affirmatively the question (a) and, partially, (c). This answer, precisely formulated, constitutes the statement of the first of two essential theorems to our "axiomatic goals" — the (1st and 2nd) $S$-spaces extension theorems — both proved in Chapter 4.

With enough material, for now, about the question (a), let us turn back our attention to the question (b). Let then $\widetilde{\mathscr{G}}(I)$ be the above described strict and closed extension of the $S$-space (abelian, surjective, and with identity) $\mathscr{G}(I)$, and let us suppose that the question (b) has an affirmative answer, that is, let us admit that there exists a locally closed extension,

$$\overline{\mathscr{G}}(I) = \Big(\overline{\mathbb{G}}\big(\Gamma(I)\big), \overline{i}\big(\Gamma^2(I)\big), \overline{\Theta}\big(\Delta(I)\big)\Big),$$

of the $S$-space $\widetilde{\mathscr{G}}(I)$.

As we know, the family $\widetilde{i}(\Gamma^2(I))$ defined on page 116 is the extension to the family $\widetilde{\mathbb{G}}(\Gamma(I))$ of the bonding $i(\Gamma^2(I))$ of the $S$-space $\mathscr{G}(I)$. Similarly, since, by hypothesis, $\overline{\mathscr{G}}(I)$ is an extension of $\widetilde{\mathscr{G}}(I)$, $\overline{i}(\Gamma^2(I))$ is the extension to $\overline{\mathbb{G}}(\Gamma(I))$ of the bonding $\widetilde{i}(\Gamma^2(I))$ of the $S$-space $\widetilde{\mathcal{G}}(I)$ and is, therefore, according to Definition 3.22, the family defined ahead:

$$\overline{i}\big(\Gamma^2(I)\big) := \Big\{\overline{i}_{(\gamma',\gamma)}\Big\}_{(\gamma',\gamma)\in\Gamma^2(I)},$$

where $\overline{i}_{(\gamma',\gamma)}$, for each $(\gamma',\gamma) \in \Gamma^2(I)$, is the function

$$\begin{aligned} \overline{i}_{(\gamma',\gamma)} : \overline{H}(\gamma) &\longrightarrow \overline{H}(\gamma') \\ \overline{\Phi} &\longmapsto \overline{i}_{(\gamma',\gamma)}(\overline{\Phi}) := \overline{\gamma}'\Big(\widetilde{i}_{(\gamma',\gamma)}\Big(\big(\overline{\gamma}\big)^{-1}(\overline{\Phi})\Big)\Big) \end{aligned}$$

being $\overline{\delta}$, for each $\delta \in \Gamma(I)$, the function

$$\begin{aligned} \overline{\delta} : \widetilde{H}(\delta) &\longrightarrow \overline{H}(\delta) \\ \widetilde{\Psi} &\longmapsto \overline{\delta}(\widetilde{\Psi}) := \overline{\Psi}\Big(\text{the prolongation of } \widetilde{\Psi} \text{ to } \overline{G}(\delta)\Big). \end{aligned}$$

Remembering now the definitions of extension of a family of $S$-groups (Definition 3.8) and extension of $S$-groups (Definition 1.6(d)), then, since $\overline{\mathbb{G}}(\Gamma(I))$ is an extension of the family of $S$-groups $\widetilde{\mathbb{G}}(\Gamma(I))$ (for, by hypothesis, $\overline{\mathscr{G}}(I)$ is an extension of $\widetilde{\mathscr{G}}(I)$) which, in turn, is an extension of the family $\mathbb{G}(\Gamma(I))$ (since $\widetilde{\mathscr{G}}(I)$ is an extension of $\mathscr{G}(I)$), we conclude that $\overline{\mathbb{G}}(\Gamma(I))$ is also an extension of the family of $S$-groups $\mathbb{G}(\Gamma(I))$. Therefore, beyond the bondings $\widetilde{i}(\Gamma^2(I))$ (the extension of $i(\Gamma^2(I))$ to $\widetilde{\mathbb{G}}(\Gamma(I))$) and $\overline{i}(\Gamma^2(I))$ (the extension of $\widetilde{i}(\Gamma^2(I))$ to $\overline{\mathbb{G}}(\Gamma(I))$) we also have, since $\overline{\mathbb{G}}(\Gamma(I))$ is an extension of the family $\mathbb{G}(\Gamma(I))$, a third bonding that we will denote by $\widehat{i}(\Gamma^2(I))$, namely, the extension of $i(\Gamma^2(I))$ to the family $\overline{\mathbb{G}}(\Gamma(I))$ that, according to Defnition 3.22, is given by:

$$\widehat{i}\big(\Gamma^2(I)\big) := \Big\{\widehat{i}_{(\gamma',\gamma)}\Big\}_{(\gamma',\gamma)\in\Gamma^2(I)}$$

where, for each $(\gamma',\gamma) \in \Gamma^2(I)$, $\widehat{i}_{(\gamma',\gamma)}$ is the function

$$\begin{aligned} \widehat{i}_{(\gamma',\gamma)} : \overline{H}(\gamma) &\longrightarrow \overline{H}(\gamma') \\ \overline{\Phi} &\longmapsto \widehat{i}_{(\gamma',\gamma)}(\overline{\Phi}) := \widehat{\gamma}'\Big(i_{(\gamma',\gamma)}\Big(\big(\widehat{\gamma}\big)^{-1}(\overline{\Phi})\Big)\Big) \end{aligned}$$

being $\widehat{\delta}$, for each $\delta \in \Gamma(I)$, the function

$$\widehat{\delta} : H(\delta) \longrightarrow \overline{H}(\delta)$$
$$\Phi \longmapsto \widehat{\delta}(\Phi) := \text{the prolongation of } \Phi \text{ to } \overline{G}(\delta).$$

Observe that if the families $\overline{i}(\Gamma^2(I))$ and $\widehat{i}(\Gamma^2(I))$ are not equal, then, the $S$-space $\overline{\mathscr{G}}(I) = (\overline{\mathbb{G}}(\Gamma(I)), \overline{i}(\Gamma^2(I)), \overline{\Theta}(\Delta(I)))$ will not be an extension of the $S$-space $\mathscr{G}(I) = (\mathbb{G}(\Gamma(I)), i(\Gamma^2(I)), \Theta(\Delta(I)))$ since, by Definition 3.23, for $\overline{\mathscr{G}}(I)$ to be an extension of $\mathscr{G}(I)$ is necessary that the bonded family $(\overline{\mathbb{G}}(\Gamma(I)), \overline{i}(\Gamma^2(I)))$ to be an extension of the bonded family $(\mathbb{G}(\Gamma(I)), i(\Gamma^2(I)))$ which, in turn (according to Definition 3.23) only occurs if $\overline{i}(\Gamma^2(I))$ is the extension to the family $\overline{\mathbb{G}}(\Gamma(I))$ of $i(\Gamma^2(I))$, that is,

$$\overline{i}\Big(\Gamma^2(I)\Big) = \widehat{i}\Big(\Gamma^2(I)\Big).$$

However, it is easy to see that this occurs. In fact, reporting ourselves to the definitions of $S$-group extension (Definition 1.6(d)) and homomorphism prolongation (Definition 1.6(b)), and taking into account that, for each $\gamma \in \Gamma(I)$, the $S$-group $\overline{\mathbb{G}}(\gamma) = (\overline{G}(\gamma), \overline{H}(\gamma))$ is an extension of the $S$-group $\widetilde{\mathbb{G}}(\gamma) = (\widetilde{G}(\gamma), \widetilde{H}(\gamma))$ which, in turn, is an extension of the $S$-group $\mathbb{G}(\gamma) = (G(\gamma), H(\gamma))$, it results that, for each $\gamma \in \Gamma(I)$,

$$\widehat{\gamma}(\Phi) = \overline{\gamma}\Big(\widetilde{\gamma}(\Phi)\Big) \quad \text{for every} \quad \Phi \in H(\gamma),$$

that is, the prolongation of $\Phi \in H(\gamma)$ to $\overline{G}(\gamma)$ $(\widehat{\gamma}(\Phi) \in \overline{H}(\gamma))$ is the prolongation to $\overline{G}(\gamma)$ of the prolongation $\widetilde{\Phi} = \widetilde{\gamma}(\Phi)$ of $\Phi$ to $\widetilde{G}(\gamma)$. Hence,

$$\widehat{\gamma} = \overline{\gamma}\widetilde{\gamma} \quad \text{and} \quad (\widehat{\gamma})^{-1} = (\widetilde{\gamma})^{-1}(\overline{\gamma})^{-1},$$

and, thus, taking into account the definition of the function $\widehat{i}_{(\gamma',\gamma)}$ given above, we can conclude that:

$$\widehat{i}_{(\gamma',\gamma)}(\overline{\Phi}) = \overline{\gamma'}\left(\widetilde{\gamma'}\left(i_{(\gamma',\gamma)}\left((\widetilde{\gamma})^{-1}\left((\overline{\gamma})^{-1}\left(\overline{\Phi}\right)\right)\right)\right)\right) \quad \text{for every} \quad \overline{\Phi} \in \overline{H}(\gamma).$$

On the other hand, from the definition of $\widetilde{i}_{(\gamma',\gamma)}$ presented on page 116 we have

$$\widetilde{i}_{(\gamma',\gamma)} = \widetilde{\gamma}' i_{(\gamma',\gamma)}(\widetilde{\gamma})^{-1},$$

that allows us to write the above expression of $\widehat{i}_{(\gamma',\gamma)}(\overline{\Phi})$ in the form

$$\widehat{i}_{(\gamma',\gamma)}(\overline{\Phi}) = \overline{\gamma'}\left(\widetilde{i}_{(\gamma',\gamma)}\left((\overline{\gamma})^{-1}\left(\overline{\Phi}\right)\right)\right)$$

for every $\overline{\Phi} \in \overline{H}(\gamma)$. Comparing now this last expression of $\widehat{i}_{(\gamma',\gamma)}$ with the one for $\overline{i}_{(\gamma',\gamma)}$ given on page 118, we obtain that, for any $(\gamma', \gamma) \in \Gamma^2(I)$,

$$\widehat{i}_{(\gamma',\gamma)}\left(\overline{\Phi}\right) = \overline{i}_{(\gamma',\gamma)}\left(\overline{\Phi}\right) \quad \text{for every} \quad \overline{\Phi} \in \overline{H}(\gamma),$$

that is,

$$\widehat{\imath}\Big(\Gamma^2(I)\Big) = \overline{\imath}\Big(\Gamma^2(I)\Big).$$

Finally, let us note that the arguments that led us to the result above does not depend on $\mathscr{G}(I)$ being abelian, surjective, and with identity, as well as they are independent of $\overline{\mathscr{G}}(I)$ being a locally closed extension of $\widetilde{\mathscr{G}}(I)$ and from the latter being a strict and closed extension of $\mathscr{G}(I)$; they are only backed on the hypothesis of $\overline{\mathscr{G}}(I)$ being an extension of $\widetilde{\mathscr{G}}(I)$ and the latter being a extension of $\mathscr{G}(I)$. Hence, taking this into account, we obtain the following result:

- Let

$$\overline{\mathscr{G}}(I) = \bigg(\overline{\mathbb{G}}\Big(\Gamma(I)\Big), \overline{\imath}\Big(\Gamma^2(I)\Big), \overline{\Theta}\Big(\Delta(I)\Big)\bigg),$$
$$\widetilde{\mathscr{G}}(I) = \bigg(\widetilde{\mathbb{G}}\Big(\Gamma(I)\Big), \widetilde{\imath}\Big(\Gamma^2(I)\Big), \widetilde{\Theta}\Big(\Delta(I)\Big)\bigg) \quad \text{and}$$
$$\mathscr{G}(I) = \bigg(\mathbb{G}\Big(\Gamma(I)\Big), i\Big(\Gamma^2(I)\Big), \Theta\Big(\Delta(I)\Big)\bigg)$$

  be $S$-spaces such that $\overline{\mathscr{G}}(I)$ is an extension of $\widetilde{\mathscr{G}}(I)$ which, in turn, is an extension of $\mathscr{G}(I)$. Then, $(\overline{\mathbb{G}}(\Gamma(I)), \overline{\imath}(\Gamma^2(I)))$ is an extension of the bonded family $(\mathbb{G}(\Gamma(I)), i(\Gamma^2(I)))$, that is, the family of $S$-groups $\overline{\mathbb{G}}(\Gamma(I))$ is an extension of the family of $S$-groups $\mathbb{G}(\Gamma(I))$ and $\overline{\imath}(\Gamma^2(I)) = \widehat{\imath}(\Gamma^2(I))$, being $\widehat{\imath}(\Gamma^2(I))$ as previously defined (at page 118), that is, the extension of the bonding $i(\Gamma^2(I))$ to the family $\overline{\mathbb{G}}(\Gamma(I))$.

With this result, and taking into account that, under the same hypothesis for $\overline{\mathscr{G}}(I)$, $\widetilde{\mathscr{G}}(I)$ and $\mathscr{G}(I)$, we trivially get that

$$\overline{\Theta}_{(\gamma',\gamma)}(g) = \Theta_{(\gamma',\gamma)}(g)$$

for any $(\gamma', \gamma) \in \Delta(I)$ and every $g \in G(\gamma)$, that is, the restriction $\overline{\Theta}(\Delta(I))$ for the bonded family $(\overline{\mathbb{G}}(\Gamma(I)), \overline{\imath}(\Gamma^2(I)))$ is a prolongation of the restriction $\Theta(\Delta(I))$ to the bonded family in question (Definition 3.23(a)), we conclude that:

- If $\overline{\mathscr{G}}(I)$, $\widetilde{\mathscr{G}}(I)$ and $\mathscr{G}(I)$ are $S$-spaces such that $\overline{\mathscr{G}}(I)$ is an extension of $\widetilde{\mathscr{G}}(I)$ which, in turn, is an extension of $\mathscr{G}(I)$, then, $\overline{\mathscr{G}}(I)$ also is an extension of $\mathscr{G}(I)$.

# $S$-Space Isomorphism

## 3.27 Motivation

Our considerations in this item will refer to the $S$-spaces

$$\mathscr{G}(I) = \Bigg( \mathbb{G}\Big(\Gamma(I)\Big) = \Big\{ \mathbb{G}(\gamma) = \Big(G(\gamma), H(\gamma)\Big) \Big\}_{\gamma \in \Gamma(I)}, \quad i\Big(\Gamma^2(I)\Big) = \Big\{ i_{(\gamma',\gamma)} \Big\}_{(\gamma',\gamma) \in \Gamma^2(I)}, \quad \Theta\Big(\Delta(I)\Big) = \Big\{ \Theta_{(\gamma',\gamma)} \Big\}_{(\gamma',\gamma) \in \Delta(I)} \Bigg)$$

and

$$\mathscr{E}(I) = \Bigg( \mathbb{E}\Big(\Gamma(I)\Big) = \Big\{ \mathbb{E}(\gamma) = (E(\gamma), F(\gamma)) \Big\}_{\gamma \in \Gamma(I)}, \quad k\Big(\Gamma^2(I)\Big) = \Big\{ k_{(\gamma',\gamma)} \Big\}_{(\gamma',\gamma) \in \Gamma^2(I)}, \quad \sigma\Big(\Delta(I)\Big) = \Big\{ \sigma_{(\gamma',\gamma)} \Big\}_{(\gamma',\gamma) \in \Delta(I)} \Bigg)$$

arbitrarily fixed, unless the requirement for the same indexer $\Gamma(I)$.

For $S$-groups we already have a well-defined concept of isomorphism, formulated in Definition 1.4(c). Since a $S$-space, for example $\mathscr{G}(I)$ above, has a family of $S$-groups, $\mathbb{G}(\Gamma(I))$, as one of its constituents, the most fundamental of them in the sense that the other two, the families $i(\Gamma^2(I))$ and $\Theta(\Delta(I))$, establish bonds between the elements of the $S$-groups of the family $\mathbb{G}(\Gamma(I))$, then, it is more than reasonable that the notion of $S$-group isomorphism to be incorporated into the homonymous concept we seek for $S$-space, in the sense that the families of $S$-groups $\mathbb{G}(\Gamma(I))$ and $\mathbb{E}(\Gamma(I))$ are such that, for each $\gamma \in \Gamma(I)$, the respective $S$-groups, $\mathbb{G}(\gamma) = (G(\gamma), H(\gamma))$ and $\mathbb{E}(\gamma) = (E(\gamma), F(\gamma))$, are isomorphic by Definition 1.4(c). However, clearly is not enough for the families $\mathbb{G}(\Gamma(I))$ and $\mathbb{E}(\Gamma(I))$ to be isomorphic in the sense above specified so that one can claim, in a sensible manner, that the $S$-spaces $\mathscr{G}(I)$ and $\mathscr{E}(I)$ are isomorphic.

In fact, as the families $i(\Gamma^2(I))$ and $k(\Gamma^2(I))$, as well as $\Theta(\Delta(I))$ and $\sigma(\Delta(I))$, establish bonds between the elements of the families of $S$-groups $\mathbb{G}(\Gamma(I))$ and $\mathbb{E}(\Gamma(I))$ ($i(\Gamma^2(I))$ and $k(\Gamma^2(I))$ between the semigroups and $\Theta(\Delta(I))$ and $\sigma(\Delta(I))$ between the groups, of the respective $S$-groups), it is necessary for the isomorphism between $\mathbb{G}(\Gamma(I))$ and $\mathbb{E}(\Gamma(I))$ to "extend" itself to the families $i(\Gamma^2(I))$ and $k(\Gamma^2(I))$, as well as to $\Theta(\Delta(I))$ and $\sigma(\Delta(I))$, somehow isomorphic, that is, compatible with the general notion of isomorphism.

To better clarify this last statement, let us suppose that the families

$$\mathbb{G}\Big(\Gamma(I)\Big) = \Big\{\mathbb{G}(\gamma) = \Big(G(\gamma), H(\gamma)\Big)\Big\}_{\gamma\in\Gamma(I)}$$

and

$$\mathbb{E}\Big(\Gamma(I)\Big) = \Big\{\mathbb{E}(\gamma) = \Big(E(\gamma), F(\gamma)\Big)\Big\}_{\gamma\in\Gamma(I)}$$

are isomorphic, that is, for each $\gamma \in \Gamma(I)$ the $S$-groups $\mathbb{G}(\gamma) = (G(\gamma), H(\gamma))$ and $\mathbb{E}(\gamma) = (E(\gamma), F(\gamma))$ are isomorphic. Hence, taking into account the definition of isomorphism of $S$-groups (Definition 1.4(c)), for each $\gamma \in \Gamma(I)$, there exist an isomorphism $\alpha_\gamma \colon G(\gamma) \longrightarrow E(\gamma)$ from the group $G(\gamma)$ onto the group $E(\gamma)$ and an isomorphism $\beta\gamma \colon H(\gamma) \longrightarrow F(\gamma)$ from the semigroup $H(\gamma)$ onto the semigroup $F(\gamma)$ such that:

$$\alpha_\gamma\Big(\Phi(g)\Big) = \beta_\gamma(\Phi)\Big(\alpha_\gamma(g)\Big)$$

for every $\Phi \in H(\gamma)$ and every $g \in G(\gamma)_\Phi$ (the domain of $\Phi$).

Now, as we know, for each $(\gamma', \gamma) \in \Delta(I)$, $\Theta_{(\gamma',\gamma)} \colon G(\gamma) \longrightarrow G(\gamma')$ is a homomorphism from the group $G(\gamma)$ into the group $G(\gamma')$ while $\sigma_{(\gamma',\gamma)} \colon E(\gamma) \longrightarrow E(\gamma')$ is a homomorphism from the group $E(\gamma)$ into the group $E(\gamma')$. Hence, since the groups $G(\gamma)$ and $E(\gamma)$, as well as $G(\gamma')$ and $E(\gamma')$, are related through the isomorphisms $\alpha_\gamma \colon G(\gamma) \longrightarrow E(\gamma)$ and $\alpha_{\gamma'} \colon G(\gamma') \longrightarrow E(\gamma')$, any notion of isomorphism for $S$-space, in order to be sensible, should require $\alpha_\gamma$ and $\alpha_{\gamma'}$ to be such that (as illustrates the diagram on Figure 3.1): the image of $g \in G(\gamma)$ by the homomorphism $\Theta_{(\gamma',\gamma)} \colon G(\gamma) \longrightarrow G(\gamma')$ ($\Theta_{(\gamma',\gamma)}(g) \in G(\gamma')$) and the image of $\alpha_\gamma(g) \in E(\gamma)$ by the homomorphism $\sigma_{(\gamma',\gamma)} \colon E(\gamma) \longrightarrow E(\gamma')$ ($\sigma_{(\gamma',\gamma)}(\alpha_\gamma(g)) \in E(\gamma')$) correspond with each other through the isomorphism $\alpha_{\gamma'} \colon G(\gamma') \longrightarrow E(\gamma')$, that is,

$$\alpha_{\gamma'}\Big(\Theta_{(\gamma',\gamma)}(g)\Big) = \sigma_{(\gamma',\gamma)}\Big(\alpha_\gamma(g)\Big).$$

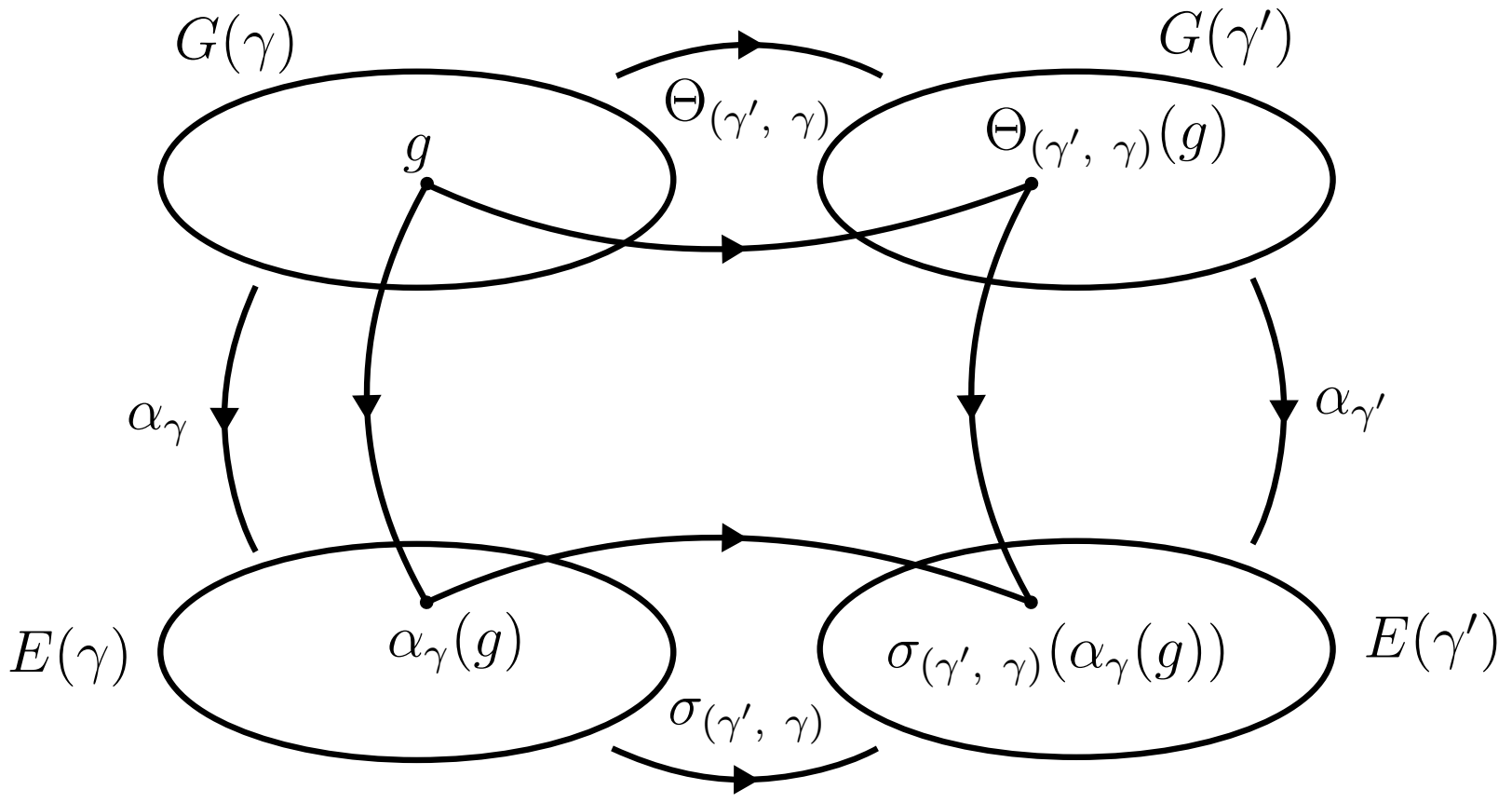

Figure 3.1 – Representation of $\alpha_{\gamma'}\Big(\Theta_{(\gamma',\gamma)}(g)\Big) = \sigma_{(\gamma',\gamma)}\Big(\alpha_\gamma(g)\Big)$.

And about the families $i(\Gamma^2(I))$ and $k(\Gamma^2(I))$ which establish connections between the semigroups of the $S$-groups of the families $\mathbb{G}(\Gamma(I))$ and $\mathbb{E}(\Gamma(I))$, respectively, how the isomorphisms $\beta_\gamma \colon H(\gamma) \longrightarrow F(\gamma)$ should "extend" to them, that is, relate with them?

Analogously, for each $(\gamma', \gamma) \in \Gamma^2(I)$, $i_{(\gamma',\gamma)}\colon H(\gamma) \longrightarrow H(\gamma')$ is an isomorphism between the semigroups $H(\gamma)$ and $H(\gamma')$ while $k_{(\gamma',\gamma)}\colon F(\gamma) \longrightarrow F(\gamma')$ is an isomorphism between the semigroups $F(\gamma)$ and $F(\gamma')$. Therefore, since the semigroups $H(\gamma)$ and $F(\gamma)$, as well as $H(\gamma')$ and $F(\gamma')$, are related through the isomorphisms $\beta_\gamma\colon H(\gamma) \longrightarrow F(\gamma)$ and $\beta_{\gamma'}\colon H(\gamma') \longrightarrow F(\gamma')$, the above mentioned relationship that is sought for isomorphisms $\beta_\gamma$ with respect to the families $i(\Gamma^2(I))$ and $k(\Gamma^2(I))$, seems to be no other if not the following (the diagram in Figure 3.2 illustrates):

$$\beta_{\gamma'}\Big(i_{(\gamma',\gamma)}(\Phi)\Big) = k_{(\gamma',\gamma)}\Big(\beta_\gamma(\Phi)\Big)$$

for every $\Phi \in H(\gamma)$ and every $(\gamma', \gamma) \in \Gamma^2(I)$.

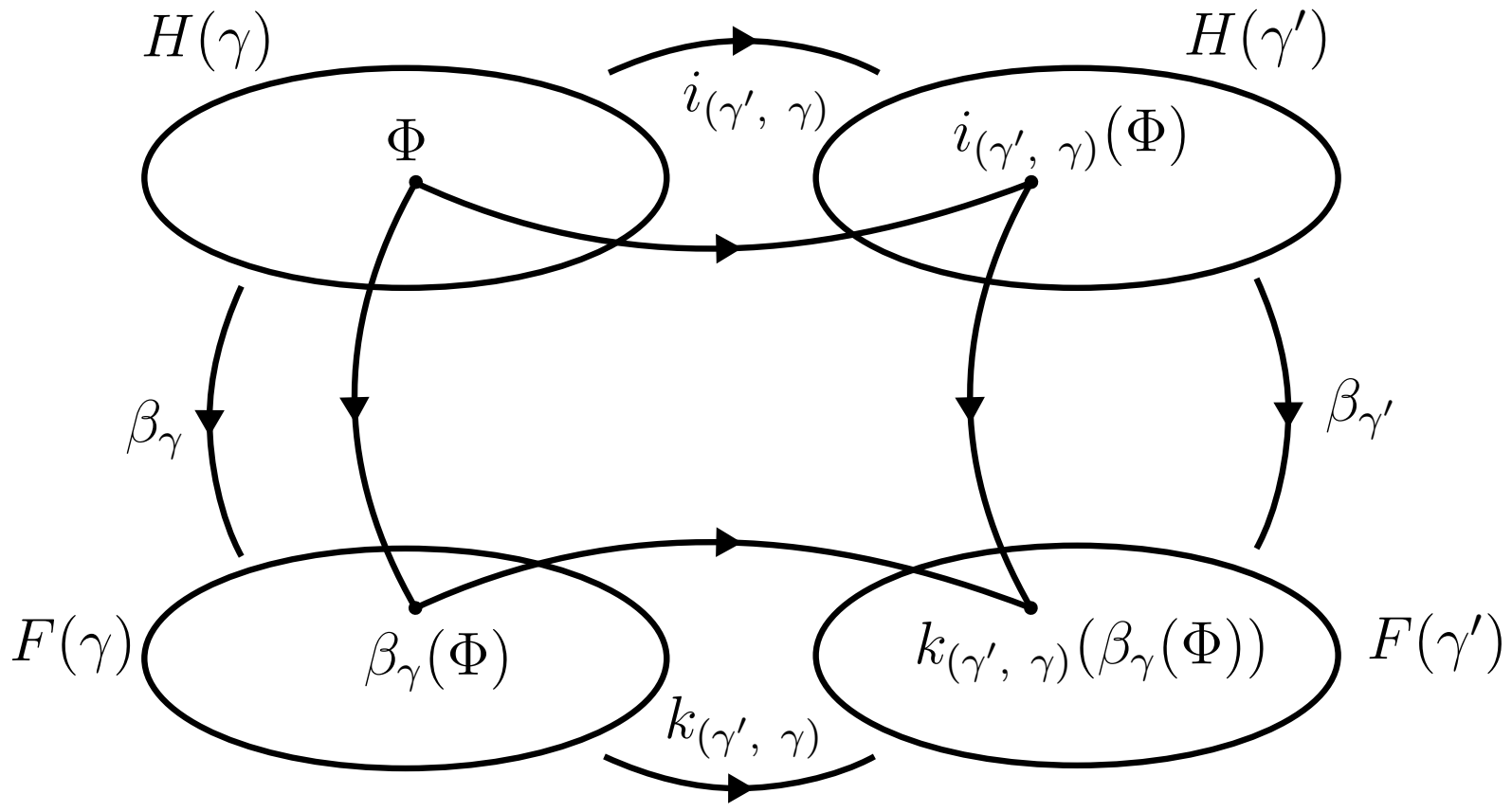

Figure 3.2 – Representation of $\beta_{\gamma'}\Big(i_{(\gamma',\gamma)}(\Phi)\Big) = k_{(\gamma',\gamma)}\Big(\beta_\gamma(\Phi)\Big)$.

The above considerations suggest the definition of an isomorphism for $S$-space given below.

## 3.28 Definition

Let

$$\mathscr{G}(I) = \left(\mathbb{G}\Big(\Gamma(I)\Big) = \Big\{\mathbb{G}(\gamma) = \Big(G(\gamma), H(\gamma)\Big)\Big\}_{\gamma\in\Gamma(I)}, \quad i\Big(\Gamma^2(I)\Big) = \Big\{i_{(\gamma',\gamma)}\Big\}_{(\gamma',\gamma)\in\Gamma^2(I)}, \quad \Theta\Big(\Delta(I)\Big) = \Big\{\Theta_{(\gamma',\gamma)}\Big\}_{(\gamma',\gamma)\in\Delta(I)}\right)$$

and

$$\mathscr{E}(I) = \left( \mathbb{E}\Big(\Gamma(I)\Big) = \Big\{ \mathbb{E}(\gamma) = (E(\gamma), F(\gamma)) \Big\}_{\gamma \in \Gamma(I)}, \quad k\Big(\Gamma^2(I)\Big) = \Big\{ k_{(\gamma',\gamma)} \Big\}_{(\gamma',\gamma) \in \Gamma^2(I)}, \quad \sigma\Big(\Delta(I)\Big) = \Big\{ \sigma_{(\gamma',\gamma)} \Big\}_{(\gamma',\gamma) \in \Delta(I)} \right)$$

be $S$-spaces with the same indexer $\Gamma(I)$. We will say that $\mathscr{G}(I)$ is **isomorphic** to $\mathscr{E}(I)$ if and only if, for each $\gamma \in \Gamma(I)$, there exist an isomorphism $\alpha_\gamma \colon G(\gamma) \longrightarrow E(\gamma)$ from the group $G(\gamma)$ onto the group $E(\gamma)$ and an isomorphism $\beta_\gamma \colon H(\gamma) \longrightarrow F(\gamma)$ from the semigroup $H(\gamma)$ onto the semigroup $F(\gamma)$, such that:

**(a)** $\alpha_\gamma\Big(\Phi(g)\Big) = \beta_\gamma(\Phi)\Big(\alpha_\gamma(g)\Big)$ for every $\Phi \in H(\gamma)$ and every $g \in G(\gamma)_\Phi$;

**(b)** $\alpha_{\gamma'}\Big(\Theta_{(\gamma',\gamma)}(g)\Big) = \sigma_{(\gamma',\gamma)}\Big(\alpha_\gamma(g)\Big)$ for every $g \in G(\gamma)$ and every $(\gamma', \gamma) \in \Delta(I)$;

**(c)** $\beta_{\gamma'}\Big(i_{(\gamma',\gamma)}(\Phi)\Big) = k_{(\gamma',\gamma)}\Big(\beta_\gamma(\Phi)\Big)$ for every $\Phi \in H(\gamma)$ and every $(\gamma', \gamma) \in \Gamma^2(I)$.

**Remark.** It turns out that if $\mathscr{G}(I)$, $\mathscr{E}(I)$ and $\mathscr{F}(I)$ are $S$-spaces such that $\mathscr{G}(I)$ is isomorphic to $\mathscr{E}(I)$ and $\mathscr{E}(I)$ is isomorphic to $\mathscr{F}(I)$, then, $\mathscr{G}(I)$ is isomorphic to $\mathscr{F}(I)$; also, $\mathscr{E}(I)$ is isomorphic to $\mathscr{G}(I)$.

# Inherited Properties

## 3.29 Motivation

Let

$$\mathscr{G}(I) = \left( \mathbb{G}\Big(\Gamma(I)\Big), i\Big(\Gamma^2(I)\Big), \Theta\Big(\Delta(I)\Big) \right)$$

be a $S$-space with indexer $\Gamma(I)$. As we know, if $\gamma \in \Gamma(I)$, then, $\Gamma(\gamma) = \Gamma(I) \cap \mathcal{P}(\gamma)$ also is an indexer and, according to Definition 3.17, the $S$-subspaces of $\mathscr{G}(I)$ are, exactly, the triplets

$$\mathscr{G}(\gamma) = \left( \mathbb{G}\Big(\Gamma(\gamma)\Big), i\Big(\Gamma^2(\gamma)\Big), \Theta\Big(\Delta(\gamma)\Big) \right)$$

one for each $\gamma \in \Gamma(I)$.

A question that naturally arises regarding the $S$-subspaces of $\mathscr{G}(I)$ is the following: If $\mathscr{G}(I)$ has a certain property, would its $S$-subspaces, $\mathscr{G}(\gamma)$, inherit it?

Hence, for example, the property of being a $S$-space, obviously (tautologically) fulfilled by any $S$-space $\mathscr{G}(I)$, is inherited by its $S$-subspaces: every $S$-subspace, $\mathscr{G}(\gamma)$, of a given $S$-space, $\mathscr{G}(I)$, is, itself, a $S$-space, as trivially turn out from the definitions of $S$-space and $S$-subspace.

In what follows, in item 3.30, some other examples related to inheritance of properties are presented, examples that will be relevant for obtaining some results contained in the next chapter.

## 3.30 Inherited Properties: Examples

The reader will find no trouble, resorting to the pertinent definitions in each case, to verify the correctness of the claims (a) to (f) ahead, relative to arbitrarily fixed $S$-spaces $\mathscr{G}(I)$ and $\mathscr{E}(I)$ with the same indexer $\Gamma(I)$.

**(a)** If $\mathscr{G}(I)$ is abelian, surjective or with identity, its $S$-subspaces, $\mathscr{G}(\gamma)$ for each $\gamma \in \Gamma(I)$, also are abelian, surjective or with identity, respectively;

**(b)** If $\mathscr{G}(I)$ is coherent, its $S$-subspaces, $\mathscr{G}(\gamma)$ for each $\gamma \in \Gamma(I)$, also are coherent;

**(c)** If $\mathscr{E}(I)$ is an extension of $\mathscr{G}(I)$, its $S$-subspaces corresponding to each $\gamma \in \Gamma(I)$, $\mathscr{E}(\gamma)$ and $\mathscr{G}(\gamma)$, are such that $\mathscr{E}(\gamma)$ is an extension of $\mathscr{G}(\gamma)$;

**(d)** If $\mathscr{E}(I)$ is a strict (respectively, closed) extension of $\mathscr{G}(I)$, its $S$-subspaces corresponding to each $\gamma \in \Gamma(I)$, $\mathscr{E}(\gamma)$ and $\mathscr{G}(\gamma)$, are such that $\mathscr{E}(\gamma)$ is a strict (respectively, closed) extension of $\mathscr{G}(\gamma)$;

**(e)** If $\mathscr{E}(I)$ is a locally closed extension of $\mathscr{G}(I)$, its $S$-subspaces corresponding to each $\gamma \in \Gamma(I)$, $\mathscr{E}(\gamma)$ and $\mathscr{G}(\gamma)$, are such that $\mathscr{E}(\gamma)$ is a locally closed extension of $\mathscr{G}(\gamma)$;

**(f)** If $\mathscr{G}(I)$ and $\mathscr{E}(I)$ are isomorphic, its $S$-subspaces corresponding to each $\gamma \in \Gamma(I)$, $\mathscr{G}(\gamma)$ and $\mathscr{E}(\gamma)$, also are isomorphic.

Next, closing this chapter, we prove a proposition and a lemma. The proposition establishes an "exclusion relation", regarding the possibility of existence, of certain extensions of an abelian and surjective $S$-space with identity. As for the lemma, it regards a trivial consequence of the definition of coherent $S$-space, it is only presented to facilitate future references.

## 3.31 Proposition

*Let* $\mathscr{G}(I) = (\mathbb{G}(\Gamma(I)), i(\Gamma^2(I)), \Theta(\Delta(I)))$ *be an abelian and surjective $S$-space with identity and* $\widetilde{\mathscr{G}}(I) = (\widetilde{\mathbb{G}}(\Gamma(I)), \widetilde{i}(\Gamma^2(I)), \widetilde{\Theta}(\Delta(I)))$ *be a strict and closed extension* *of $\mathscr{G}(I)$ which, as claimed in 3.26 (and proved in the next chapter), exists and, unless isomorphism, is unique. If there exists an extension $\overline{\mathscr{G}}(I) = (\overline{\mathbb{G}}(\Gamma(I)), \overline{i}(\Gamma^2(I)), \overline{\Theta}(\Delta(I)))$ of $\widetilde{\mathscr{G}}(I)$, $\overline{\mathscr{G}}(I) \neq \widetilde{\mathscr{G}}(I)$, locally closed and coherent, then, $\widetilde{\mathscr{G}}(I)$ is not a coherent $S$-space.*

*Proof.* Suppose $\overline{\mathscr{G}}(I)$ is an extension (of $\widetilde{\mathscr{G}}(I)$) as described above. Since $\overline{\mathscr{G}}(I) \neq \widetilde{\mathscr{G}}(I)$, there exists $\gamma \in \Gamma(I)$ such that $\widetilde{G}(\gamma) \subseteq \overline{G}(\gamma)$ and $\overline{G}(\gamma) \neq \widetilde{G}(\gamma)$. For one such $\gamma$, let $\overline{g} \in \overline{G}(\gamma)$ be such that $\overline{g} \notin \widetilde{G}(\gamma)$ arbitrarily fixed. Let also

$$\overline{\xi}(\gamma) \subseteq \Gamma(\gamma) = \Gamma(I) \cap \mathcal{P}(\gamma)$$

be a cover of $\gamma$ such that

$$\overline{\Theta}_{(\xi,\gamma)}(\overline{g}) \in \widetilde{G}(\xi) \quad \text{for every} \quad \xi \in \overline{\xi}(\gamma).$$

Let us note that the existence of such cover $\overline{\xi}(\gamma)$ is ensured by the fact of the $S$-subspace $\overline{\mathscr{G}}(\gamma)$ of $\overline{\mathscr{G}}(I)$ being a locally closed extension of the $S$-subspace $\widetilde{\mathscr{G}}(\gamma)$ of $\widetilde{\mathscr{G}}(I)$ (recap Definition 3.23(c-2) and item (e) in 3.30). The family

$$\left\{\overline{\Theta}_{(\xi,\gamma)}(\overline{g})\right\}_{\xi \in \overline{\xi}(\gamma)}$$

is, then, clearly, coherent in $\widetilde{\mathscr{G}}(\gamma)$. Let us suppose now that $\widetilde{\mathscr{G}}(I)$ is coherent. In this case, there exists a single $\widetilde{g} \in \widetilde{G}(\gamma)$ such that

$$\widetilde{\Theta}_{(\xi,\gamma)}(\widetilde{g}) = \overline{\Theta}_{(\xi,\gamma)}(\overline{g}) \quad \text{for every} \quad \xi \in \overline{\xi}(\gamma).$$

Since $\overline{\Theta}(\Delta(I))$ is a prolongation of $\widetilde{\Theta}(\Delta(I))$ we have

$$\widetilde{\Theta}_{(\xi,\gamma)}(\widetilde{g}) = \overline{\Theta}_{(\xi,\gamma)}(\widetilde{g}) \quad \text{for every} \quad \xi \in \overline{\xi}(\gamma)$$

and hence we obtain:

$$\overline{\Theta}_{(\xi,\gamma)}(\widetilde{g}) = \overline{\Theta}_{(\xi,\gamma)}(\overline{g}) \quad \text{for every} \quad \xi \in \overline{\xi}(\gamma)$$

from where we conclude (according to Lemma 3.32 ahead) that $\widetilde{g} = \overline{g}$ and, therefore, since $\widetilde{g} \in \widetilde{G}(\gamma)$, that $\overline{g} \in \widetilde{G}(\gamma)$, contradicting our initial hypothesis of $\overline{g} \notin \widetilde{G}(\gamma)$. Therefore, $\widetilde{\mathscr{G}}(I)$ is not coherent. ∎

## 3.32 Lemma

*Let $\mathscr{G}(I) = (\mathbb{G}(\Gamma(I)), i(\Gamma^2(I)), \Theta(\Delta(I)))$ be a coherent $S$-space with indexer $\Gamma(I)$. Let also $\gamma \in \Gamma(I)$ and $\overline{\xi}(\gamma) \subseteq \Gamma(\gamma) = \Gamma(I) \cap \mathcal{P}(\gamma)$ be a cover of $\gamma$, arbitrarily fixed. If $g, g' \in G(\gamma)$ are such that*

$$\Theta_{(\xi,\gamma)}(g) = \Theta_{(\xi,\gamma)}(g') \quad \textit{for every} \quad \xi \in \overline{\xi}(\gamma),$$

*then, $g = g'$.*

*Proof.* Let $\gamma \in \Gamma(I)$, $\overline{\xi}(\gamma) \subseteq \Gamma(I) \cap \mathcal{P}(\gamma)$ a cover of $\gamma$ and $g, g' \in G(\gamma)$ be arbitrarily fixed. Let us consider the following families indexed by $\overline{\xi}(\gamma)$:

$$\left\{h_\xi\right\}_{\xi \in \bar{\xi}(\gamma)} \quad \text{where} \quad h_\xi \coloneqq \Theta_{(\xi,\gamma)}(g)$$

and

$$\left\{h'_\xi\right\}_{\xi \in \bar{\xi}(\gamma)} \quad \text{where} \quad h'_\xi \coloneqq \Theta_{(\xi,\gamma)}(g').$$

It turns out that $\{h_\xi\}_{\xi \in \bar{\xi}(\gamma)}$ and $\left\{h'_\xi\right\}_{\xi \in \bar{\xi}(\gamma)}$ are coherent families in $\mathscr{G}(\gamma)$ (the $S$-subspace of $\mathscr{G}(I)$ corresponding to $\gamma \in \Gamma(I)$). In fact, for the first family, for example, if $\xi, \xi' \in \overline{\xi}(\gamma)$ are such that $\xi \cap \xi' \neq \varnothing$, we get

$$\Theta_{(\xi \cap \xi', \xi)}(h_\xi) = \Theta_{(\xi \cap \xi', \xi)}\Big(\Theta_{(\xi,\gamma)}(g)\Big) = \Theta_{(\xi \cap \xi', \gamma)}(g)$$

and

$$\Theta_{(\xi \cap \xi', \xi')}(h_{\xi'}) = \Theta_{(\xi \cap \xi', \xi')}\Big(\Theta_{(\xi',\gamma)}(g)\Big) = \Theta_{(\xi \cap \xi', \gamma)}(g),$$

that is,

$$\Theta_{(\xi \cap \xi', \xi)}(h_\xi) = \Theta_{(\xi \cap \xi', \xi')}(h_{\xi'}).$$

Thus, as $\mathscr{G}(I)$ is coherent, $g$ and $g'$ are the only elements of $G(\gamma)$ such that

$$\Theta_{(\xi,\gamma)}(g) = h_\xi$$

and

$$\Theta_{(\xi,\gamma)}(g') = h'_\xi$$

for every $\xi \in \overline{\xi}(\gamma)$. Therefore, if $h_\xi = h'_\xi$, that is, $\Theta_{(\xi,\gamma)}(g) = \Theta_{(\xi,\gamma)}(g')$ for every $\xi \in \overline{\xi}(\gamma)$, then $g = g'$. ∎

# 4

# THEOREMS OF EXTENSION OF $S$-SPACES

## Introduction

### 4.1 Goals

In the previous chapter, more specifically in its Introduction (items 3.1 to 3.4), we established our "action plan" to elaborate a definition for the L. Schwartz' distributions, finite order or not, through a "simple" categoric axiomatic. Actually, our plan, in its abstract version (see item 3.4), is more audacious, as it proposes to obtain a categoric axiomatic formulated in terms of an axiomatical structure, the one of a $S$-space, that when specialized (particularized) to the $S$-space of the continuous functions, $\mathscr{C}(\mathbb{R}^n)$, admits as a model the referred L. Schwartz' distributions.

Basically, the referred Chapter 3 provided the definitions which allowed to formulate, in a precise manner, this more audacious version of our purposes, which, essentially, as already highlighted in 3.25, consists of proving the following claims about an abelian and surjective $S$-space $\mathscr{G}(I)$ with identity.

- There exists, unless isomorphism, an unique strict and closed extension of $\mathscr{G}(I)$;
- Let $\widetilde{\mathscr{G}}(I)$ be a strict and closed extension of $\mathscr{G}(I)$. There exists, unless isomorphism, an unique locally closed extension of $\widetilde{\mathscr{G}}(I)$.

This chapter has as goals to prove two theorems; one of them, referred to as the 1st Theorem of Extension of $S$-Spaces (1st TESS), establishes the first claim above, while the other, the 2nd Theorem of Extension of $S$-Spaces (2nd TESS), ensures the following "weaker" (maybe "modified" would be more appropriate) version of the second claim:

- Let $\widetilde{\mathscr{G}}(I)$ be a strict and closed extension of an abelian, surjective, with identity, and also coherent $S$-space $\mathscr{G}(I)$. Then, there exists, unless isomorphism, an unique extension of $\widetilde{\mathscr{G}}(I)$ that is locally closed and coherent.

These two theorems are, however, as we will see over the next chapter, enough for the complete achievement of our axiomatization program.

## 4.2 About the Existence of $\widetilde{\mathscr{G}}(I)$

It is worth remembering here the remarks in item 3.26 regarding the claim that, for a given $S$-space

$$\mathscr{G}(I) = \left(\mathbb{G}\left(\Gamma(I)\right) = \left\{\mathbb{G}(\gamma) = \left(G(\gamma), H(\gamma)\right)\right\}_{\gamma\in\Gamma(I)}, i\left(\Gamma^2(I)\right), \Theta\left(\Delta(I)\right)\right)$$

abelian, surjective, and with identity, there exists, essentially, only one strict and closed extension. As we saw in 3.26,

$$\widetilde{\mathbb{G}}\left(\Gamma(I)\right) = \left\{\widetilde{\mathbb{G}}(\gamma) = \left(\widetilde{G}(\gamma), \widetilde{H}(\gamma)\right)\right\}_{\gamma\in\Gamma(I)},$$

where, for each $\gamma \in \Gamma(I)$, $\widetilde{\mathbb{G}}(\gamma) = (\widetilde{G}(\gamma), \widetilde{H}(\gamma))$ is the strict and closed extension of the $S$-group $\mathbb{G}(\gamma) = (G(\gamma), H(\gamma))$ as described by the Theorem of Extension of $S$-Groups (Theorem 2.16), is a strict and closed extension of the family of $S$-groups,

$$\mathbb{G}\left(\Gamma(I)\right) = \left\{\mathbb{G}(\gamma) = \left(G(\gamma), H(\gamma)\right)\right\}_{\gamma\in\Gamma(I)},$$

of the $S$-space $\mathscr{G}(I)$. Therefore, according to Definition 3.22,

$$\left(\widetilde{\mathbb{G}}\left(\Gamma(I)\right), \widetilde{i}\left(\Gamma^2(I)\right)\right),$$

where $\widetilde{i}(\Gamma^2(I))$ is the extension of the bonding $i(\Gamma^2(I))$ (of $\mathscr{G}(I)$) to the family $\widetilde{\mathbb{G}}(\Gamma(I))$, is an extension of the bonded family (of $\mathscr{G}(I)$)

$$\left(\mathbb{G}\left(\Gamma(I)\right), i\left(\Gamma^2(I)\right)\right).$$

Hence being, according to Definition 3.23, the triplet

$$\widetilde{\mathscr{G}}(I) = \left(\widetilde{\mathbb{G}}\left(\Gamma(I)\right), \widetilde{i}\left(\Gamma^2(I)\right), \widetilde{\Theta}\left(\Delta(I)\right)\right)$$

will be a strict and closed extension of the $\mathscr{G}(I)$, only if $\widetilde{\Theta}(\Delta(I))$ is a restriction for the bonded family

$$\left(\widetilde{\mathbb{G}}\left(\Gamma(I)\right), \widetilde{i}\left(\Gamma^2(I)\right)\right)$$

and a prolongation of $\Theta(\Delta(I))$ to this bonded family.

The question about the existence of such a family $\widetilde{\Theta}(\Delta(I))$ and, therefore, the existence of a strict and closed extension of an abelian and surjective $S$-space with identity is considered ahead in Propositions 4.3 and 4.4; the Proposition 4.3 provides an important result regarding the extension of homomorphisms associated to abelian and surjective $S$-groups with identity and its respective strict and closed extensions, while Proposition 4.4, trivially obtained from Proposition 4.3, establishes the existence and uniqueness of $\widetilde{\Theta}(\Delta(I))$.

# The Restriction $\widetilde{\Theta}(\Delta(I))$

## 4.3 Proposition

*Let $\mathbb{G} = (G, H)$ and $\mathbb{E} = (E, F)$ be abelian and surjective $S$-groups with identity. Let also $\widehat{\mathbb{G}} = (\widehat{G}, \widehat{H})$ and $\widehat{\mathbb{E}} = (\widehat{E}, \widehat{F})$ be strict and closed extensions of $\mathbb{G}$ and $\mathbb{E}$, respectively. Let us suppose that there exist an isomorphism*

$$\begin{aligned} * : H &\longrightarrow F \\ \Phi &\longmapsto \Phi^* \end{aligned}$$

*from the semigroup $H$ onto the semigroup $F$ and a homomorphism*

$$\begin{aligned} h : G &\longrightarrow E \\ g &\longmapsto h(g) \end{aligned}$$

*from the group $G$ into the group $E$, such that:*

**(a)** *for every $\Phi \in H$ and every $g \in G_\Phi$ (the domain of $\Phi$),*

$$h\Big(\Phi(g)\Big) = \Phi^*\Big(h(g)\Big) \quad \textit{and}$$

**(b)** *for every $\Phi \in H$*

$$h\Big(N(\Phi)\Big) \subseteq N(\Phi^*).$$

*Under these conditions, there exists a unique homomorphism*

$$\begin{aligned} \widehat{h} : \widehat{G} &\longrightarrow \widehat{E} \\ \widehat{g} &\longmapsto \widehat{h}(\widehat{g}) \end{aligned}$$

*from the group $\widehat{G}$ into the group $\widehat{E}$ such that*

$$\widehat{h}(g) = h(g) \quad \textit{for every} \quad g \in G$$

*and*

$$\widehat{h}\Big(\widehat{\Phi}(\widehat{g})\Big) = \widehat{\Phi^*}\Big(\widehat{h}(\widehat{g})\Big)$$

*for every $\Phi \in H$ and $\widehat{g} \in \widehat{G}$.*

*Proof.* Let us first prove the following lemma.

**Lemma .** *Let $\mathbb{G} = (G, H)$ be as described in the proposition, that is, an abelian and surjective $S$-group with identity. Let also $\Phi$ and $\Psi$ in $H$, $g$ and $f$ in $G$, $g' \in \Psi^{-1}(g)$ and $f' \in \Phi^{-1}(f)$ be such that*

$$g' - f' \in N(\Phi\Psi).$$

*With these hypotheses we have that:*

$$\Psi^{-1}(g) - \Phi^{-1}(f) \subseteq N(\Phi\Psi).$$

*Proof of the Lemma.* Since $g' \in \Psi^{-1}(g)$ and $f' \in \Phi^{-1}(f)$, then $\Psi^{-1}(g) = g' + N(\Psi)$ and $\Phi^{-1}(f) = f' + N(\Phi)$ and, therefore,

$$\Psi^{-1}(g) - \Phi^{-1}(f) = g' - f' + N(\Psi) + N(\Phi),$$

which, in turn, taking into account the hypothesis of $g' - f' \in N(\Phi\Psi)$, leads to the result.

*Proof of the Proposition.* Let $\widehat{g} \in \widehat{G}$. Since the $S$-group $\widehat{\mathbb{G}}$ is a closed (and strict) extension of $\mathbb{G}$, there exist $g \in G$ and $\Phi \in H$ such that $\widehat{g} = \widehat{\Phi}(g)$. Hence, in the hypothesis of existence of a homomorphism $\widehat{h} \colon \widehat{G} \longrightarrow \widehat{E}$ as described in the proposition we would have:

$$\widehat{h}(\widehat{g}) = \widehat{h}\Big(\widehat{\Phi}(g)\Big) = \widehat{\Phi^*}\Big(\widehat{h}(g)\Big) = \widehat{\Phi^*}\Big(h(g)\Big),$$

that is, $\widehat{h}$ would, necessarily, be given by:

$$\begin{aligned} \widehat{h} : \qquad & \widehat{G} \longrightarrow \widehat{E} \\ \widehat{g} = \widehat{\Phi}(g) & \longmapsto \widehat{h}(\widehat{g}) = \widehat{\Phi^*}\Big(h(g)\Big). \end{aligned}$$

Hence being we must prove, before anything, that $\widehat{h}$ as above is well-defined, that is: for $\Phi, \Psi \in H$ and $g, f \in G$, if $\widehat{\Phi}(g) = \widehat{\Psi}(f)$, then, $\widehat{\Phi^*}(h(g)) = \widehat{\Psi^*}(h(f))$. Let us suppose then that

$$\widehat{\Phi}(g) = \widehat{\Psi}(f)$$

and let $x$ be a generic element of $\Phi^{-1}(f) - \Psi^{-1}(g)$, that is, $x = f' - g'$ where $f' \in \Phi^{-1}(f)$ and $g' \in \Psi^{-1}(g)$. Thus, $f = \Phi(f') = \widehat{\Phi}(f')$ and $g = \Psi(g') = \widehat{\Psi}(g')$ and, consequently, $\widehat{\Phi}(g) = \widehat{\Phi}(\widehat{\Psi}(g')) = \widehat{\Phi\Psi}(g')$ and $\widehat{\Psi}(f) = \widehat{\Psi}(\widehat{\Phi}(f')) = \widehat{\Psi\Phi}(f')$. Since $\widehat{\Phi}(g) = \widehat{\Psi}(f)$, we have $\widehat{\Phi\Psi}(g') = \widehat{\Psi\Phi}(f')$, that is (since, $\widehat{\Psi\Phi} = \widehat{\Psi}\widehat{\Phi} = \widehat{\Phi}\widehat{\Psi} = \widehat{\Phi\Psi}$), $\widehat{\Phi\Psi}(g' - f') = 0$. Thus, $g' - f' \in N(\widehat{\Phi\Psi})$ and hence, taking into account that $N(\widehat{\Phi\Psi}) = N(\Phi\Psi)$ (according to Proposition 1.16(b)),

$$x = f' - g' \in N(\Phi\Psi).$$

Since $x = f' - g'$ is a generic element of $\Phi^{-1}(f) - \Psi^{-1}(g)$, we conclude that

$$\Phi^{-1}(f) - \Psi^{-1}(g) \subseteq N(\Phi\Psi).$$

Now, from the hypotheses of the proposition, specifically that $\;^*\colon H \longrightarrow F$ is an isomorphism and (b), comes:

$$h\Big(N(\Phi\Psi)\Big) \subseteq N(\Phi^*\Psi^*)$$

and hence, since $f' - g' \in \Phi^{-1}(f) - \Psi^{-1}(g) \subseteq N(\Phi\Psi)$, we get that

$$h(f' - g') = h(f') - h(g') \in N(\Phi^*\Psi^*).$$

But $\Psi(g') = g$ and hence $h\Big(\Psi(g')\Big) = h(g)$. On the other hand, from hypothesis (a) we get $h\Big(\Psi(g')\Big) = \Psi^*\Big(h(g')\Big)$ and hence, $\Psi^*\Big(h(g')\Big) = h(g)$ from where we conclude that

$$h(g') \in (\Psi^*)^{-1}\Big(h(g)\Big).$$

Analogously, from $\Phi(f') = f$ one obtains that:

$$h(f') \in (\Phi^*)^{-1}\Big(h(f)\Big).$$

From this, and taking the lemma into account, it results that

$$(\Psi^*)^{-1}\Big(h(g)\Big) - (\Phi^*)^{-1}\Big(h(f)\Big) \subseteq N(\Phi^*\Psi^*)$$

or, resorting again to Proposition 1.16(b),

$$(\Psi^*)^{-1}\Big(h(g)\Big) - (\Phi^*)^{-1}\Big(h(f)\Big) \subseteq N(\widehat{\Phi^*\Psi^*}).$$

Hence, since $h(g') - h(f') \in (\Psi^*)^{-1}(h(g)) - (\Phi^*)^{-1}(h(f))$, it results that

$$\Big(\widehat{\Phi^*\Psi^*}\Big)\Big(h(g') - h(f')\Big) = 0$$

or yet

$$\widehat{\Phi^*}\left(\widehat{\Psi^*}\Big(h(g')\Big)\right) = \widehat{\Psi^*}\left(\widehat{\Phi^*}\Big(h(f')\Big)\right).$$

Since

$$\widehat{\Psi^*}\Big(h(g')\Big) = \Psi^*\Big(h(g')\Big) = h(g)$$

and

$$\widehat{\Phi^*}\Big(h(f')\Big) = \Phi^*\Big(h(f')\Big) = h(f),$$

comes that

$$\widehat{\Phi^*}\Big(h(g)\Big) = \widehat{\Psi^*}\Big(h(f)\Big),$$

which shows us that $\widehat{h}\colon \widehat{G} \longrightarrow \widehat{E}$ such as given in page 131 is well-defined.

Let us prove now that this function $\widehat{h}$ is such that $\widehat{h}(g) = h(g)$ for every $g \in G$. For such, let us remember that the homomorphism $I_G\colon G \longrightarrow G$ (the identity on $G$) belongs to $H$ and that $g = I_G(g) = \widehat{I}_G(g)$ for every $g \in G$. Hence, taking into account the definition of $\widehat{h}$ (in page 131), we have:

$$\widehat{h}(g) = \widehat{h}\Big(\widehat{I}_G(g)\Big) = \widehat{I_G^*}\Big(h(g)\Big).$$

But $I_G^* = I_E$ ($I_E \in F$ is the identity on $E$) for $\;^*\colon H \longrightarrow F$ is an isomorphism. Hence,

$$\widehat{h}(g) = \widehat{I_E}\Big(h(g)\Big) = I_E\Big(h(g)\Big) = h(g).$$

The function $\widehat{h}$ in question is, also, such that

$$\widehat{h}\Big(\widehat{\Phi}(\widehat{g})\Big) = \widehat{\Phi^*}\Big(\widehat{h}(\widehat{g})\Big)$$

for every $\Phi \in H$ and every $\widehat{g} \in \widehat{G}$. In fact, being $\varphi \in H$ and $g \in G$ such that $\widehat{g} = \widehat{\varphi}(g)$, then, by the definition of $\widehat{h}$, we have:

$$\begin{aligned}
\widehat{h}\Big(\widehat{\Phi}(\widehat{g})\Big) &= \widehat{h}\left(\widehat{\Phi}\Big(\widehat{\varphi}(g)\Big)\right) = \widehat{h}\Big(\widehat{\Phi\varphi}(g)\Big) = \\
&= \widehat{\Big(\Phi\varphi\Big)^*}\Big(h(g)\Big) = \widehat{\Phi^*\varphi^*}\Big(h(g)\Big) = \\
&= \Big(\widehat{\Phi^*}\widehat{\varphi^*}\Big)\Big(h(g)\Big) = \widehat{\Phi^*}\left(\widehat{h}\Big(\widehat{\varphi}(g)\Big)\right) = \\
&= \widehat{\Phi^*}\Big(\widehat{h}(\widehat{g})\Big).
\end{aligned}$$

It remains to prove that $\widehat{h}$ is a homomorphism. In order to do so, let $\widehat{g}$ and $\widehat{f}$ be in $\widehat{G}$ and let us take $\Phi$ and $\Psi$ in $H$ and $g$ and $f$ in $G$ such that $\widehat{g} = \widehat{\Phi}(g)$ and $\widehat{f} = \widehat{\Psi}(f)$. Hence,

$$\widehat{g} + \widehat{f} = \widehat{\Phi}(g) + \widehat{\Psi}(f) = \widehat{\Phi\Psi}(g' + f')$$

where

$$g' \in \Psi^{-1}(g) \quad \text{and} \quad f' \in \Phi^{-1}(f).$$

Now let us calculate

$$\begin{aligned}
\widehat{h}(\widehat{g}+\widehat{f}) &= \widehat{h}\left(\widehat{\Phi\Psi}(g'+f')\right) = \left(\widehat{\Phi\Psi}\right)^{*}\left(h(g'+f')\right) = \\
&= \left(\widehat{\Phi\Psi}\right)^{*}\left(h(g')+h(f')\right) = \\
&= \left(\widehat{\Phi\Psi}\right)^{*}\left(h(g')\right) + \left(\widehat{\Phi\Psi}\right)^{*}\left(h(f')\right) = \\
&= \widehat{h}\left(\widehat{\Phi\Psi}(g')\right) + \widehat{h}\left(\widehat{\Phi\Psi}(f')\right) = \\
&= \widehat{h}\left(\widehat{\Phi}\left(\widehat{\Psi}(g')\right)\right) + \widehat{h}\left(\widehat{\Psi}\left(\widehat{\Phi}(f')\right)\right) = \\
&= \widehat{h}\left(\widehat{\Phi}\left(\Psi(g')\right)\right) + \widehat{h}\left(\widehat{\Psi}\left(\Phi(f')\right)\right) = \\
&= \widehat{h}\left(\widehat{\Phi}(g)\right) + \widehat{h}\left(\widehat{\Psi}(f)\right) = \\
&= \widehat{h}(\widehat{g}) + \widehat{h}(\widehat{f}).
\end{aligned}$$

With this, we conclude the proof of the proposition. ■

## 4.4 Proposition

*Let*

$$\mathscr{G}(I) = \left(\mathbb{G}\left(\Gamma(I)\right) = \left\{\mathbb{G}(\gamma) = \left(G(\gamma), H(\gamma)\right)\right\}_{\gamma\in\Gamma(I)}, i\left(\Gamma^2(I)\right), \Theta\left(\Delta(I)\right)\right)$$

*be an abelian and surjective S-space with identity and*

$$\widetilde{\mathbb{G}}\left(\Gamma(I)\right) = \left\{\widetilde{\mathbb{G}}(\gamma) = \left(\widetilde{G}(\gamma), \widetilde{H}(\gamma)\right)\right\}_{\gamma\in\Gamma(I)}$$

*be a family where, for each $\gamma \in \Gamma(I)$, $\widetilde{\mathbb{G}}(\gamma) = (\widetilde{G}(\gamma), \widetilde{H}(\gamma))$ is the strict and closed extension of the S-group $\mathbb{G}(\gamma) = (G(\gamma), H(\gamma))$, as described in the Theorem of Extension of S-Groups (Theorem 2.16). Let also $\widetilde{i}(\Gamma^2(I))$ be the extension to $\widetilde{\mathbb{G}}(\Gamma(I))$ of the bonding $i(\Gamma^2(I))$ of the family of S-groups $\mathbb{G}(\Gamma(I))$ of the S-space $\mathscr{G}(I)$ (see Definition 3.22), and let us considerate the bonded family $(\widetilde{\mathbb{G}}\left(\Gamma(I)\right), \widetilde{i}\left(\Gamma^2(I)\right))$ that, according to Definition 3.22, is an extension of the bonded family $(\mathbb{G}\left(\Gamma(I)\right), i\left(\Gamma^2(I)\right))$ of the S-space $\mathscr{G}(I)$.*

*There exists a unique restriction $\widetilde{\Theta}(\Delta(I))$ for the bonded family $(\widetilde{\mathbb{G}}\left(\Gamma(I)\right), \widetilde{i}\left(\Gamma^2(I)\right))$ (see Definition 3.13), which is a prolongation (see Definition 3.23(a)), to this bonded family, of the restriction $\Theta(\Delta(I))$.*

*Proof.* Let $(\gamma', \gamma) \in \Delta(I)$ be arbitrarily fixed. One easily verifies that all hypotheses of Proposition 4.3 are satisfied with the following choices:

- for the $S$-groups $\mathbb{G} = (G, H)$ and $\mathbb{E} = (E, F)$ we choose $\mathbb{G}(\gamma) = (G(\gamma), H(\gamma))$ and $\mathbb{G}(\gamma') = (G(\gamma'), H(\gamma'))$, respectively;
- for the strict and closed extensions $\widehat{\mathbb{G}}$ and $\widehat{\mathbb{E}}$ of $\mathbb{G}$ and $\mathbb{E}$ we take, respectively, $\widetilde{\mathbb{G}}(\gamma)$ and $\widetilde{\mathbb{G}}(\gamma')$;
- for the isomorphism $\ ^*\colon H \longrightarrow F$ we choose $i_{(\gamma',\gamma)}\colon H(\gamma) \longrightarrow H(\gamma')$;
- for the homomorphism $h\colon G \longrightarrow E$ we take $\Theta_{(\gamma',\gamma)}\colon G(\gamma) \longrightarrow G(\gamma')$.

With these choices, Proposition 4.3 allows us to conclude that there exists an unique homomorphism $\widetilde{\Theta}_{(\gamma',\gamma)}\colon \widetilde{G}(\gamma) \longrightarrow \widetilde{G}(\gamma')$ such that:

$$\widetilde{\Theta}_{(\gamma',\gamma)}(g) = \Theta_{(\gamma',\gamma)}(g) \quad \text{for every} \quad g \in G(\gamma)$$

and

$$\widetilde{\Theta}_{(\gamma',\gamma)}\Big(\widetilde{\Phi}(\widetilde{g})\Big) = \widetilde{i_{(\gamma',\gamma)}(\Phi)}\Big(\widetilde{\Theta}_{(\gamma',\gamma)}(\widetilde{g})\Big) \quad \text{for every} \quad \Phi \in H(\gamma) \quad \text{and} \quad \widetilde{g} \in \widetilde{G}(\gamma).$$

But, taking Definition 3.22 into account,

$$\widetilde{i_{(\gamma',\gamma)}(\Phi)} = \widetilde{\gamma}'\Big(i_{(\gamma',\gamma)}(\Phi)\Big) = \widetilde{\gamma}'\left(i_{(\gamma',\gamma)}\Big((\widetilde{\gamma})^{-1}(\widetilde{\Phi})\Big)\right) = \widetilde{i}_{(\gamma',\gamma)}(\widetilde{\Phi})$$

and hence,

$$\widetilde{\Theta}_{(\gamma',\gamma)}\Big(\widetilde{\Phi}(\widetilde{g})\Big) = \widetilde{i}_{(\gamma',\gamma)}(\widetilde{\Phi})\Big(\widetilde{\Theta}_{(\gamma',\gamma)}(\widetilde{g})\Big)$$

for every $\widetilde{\Phi} \in \widetilde{H}(\gamma)$ and $\widetilde{g} \in \widetilde{G}(\gamma)$.

These results show us that the family

$$\widetilde{\Theta}\Big(\Delta(I)\Big) := \Big\{\widetilde{\Theta}_{(\gamma',\gamma)}\Big\}_{(\gamma',\gamma)\in\Delta(I)}$$

satisfies two of the four requirements, specifically the conditions (a) and (c) of Definition 3.13, to be a restriction for the bonded family $(\widetilde{\mathbb{G}}\,(\Gamma(I))\,, \widetilde{i}\,(\Gamma^2(I)))$. Regarding the requirement (d) of the referred definition, this one is trivially satisfied in this case where the homomorphisms $\widetilde{\Phi} \in \widetilde{H}(\gamma)$, for each $\gamma \in \Gamma(I)$, are, all, endomorphisms on $\widetilde{G}(\gamma)$ and, hence, if $\widetilde{g} \in \widetilde{G}(\gamma)$ then $\widetilde{g} \in \widetilde{G}(\gamma)_{\widetilde{\Phi}}$ (the domain of $\widetilde{\Phi}$) since $\widetilde{G}(\gamma)_{\widetilde{\Phi}} = \widetilde{G}(\gamma)$.

Hence, if the requirement (b) of the definition in question, which in this case requires that

$$\widetilde{\Theta}_{(\gamma'',\gamma')}\Big(\widetilde{\Theta}_{(\gamma',\gamma)}(\widetilde{g})\Big) = \widetilde{\Theta}_{(\gamma'',\gamma)}(\widetilde{g})$$

for every $\widetilde{g} \in \widetilde{G}(\gamma)$ and any $\gamma, \gamma', \gamma'' \in \Gamma(I)$ such that $\gamma'' \subseteq \gamma' \subseteq \gamma$, is also fulfilled, then, the family $\widetilde{\Theta}(\Delta(I))$ will be a restriction for the bonded family $(\widetilde{\mathbb{G}}\,(\Gamma(I))\,, \widetilde{i}\,(\Gamma^2(I)))$ and, also, since, as showed above,

$$\widetilde{\Theta}_{(\gamma',\gamma)}(g) = \Theta_{(\gamma',\gamma)}(g) \quad \text{for every} \quad g \in G(\gamma),$$

a prolongation (to this bonded family) of the restriction $\Theta(\Delta(I))$ (see Definition 3.23(a)).

Hence, in order to complete the proof of this proposition, it remains to show that

$$\widetilde{\Theta}_{(\gamma'',\gamma')}\Big(\widetilde{\Theta}_{(\gamma',\gamma)}(\widetilde{g})\Big) = \widetilde{\Theta}_{(\gamma'',\gamma)}(\widetilde{g})$$

for $\widetilde{g} \in \widetilde{G}(\gamma)$ and $\gamma, \gamma', \gamma'' \in \Gamma(I)$ such that $\gamma'' \subseteq \gamma' \subseteq \gamma$.

Let then $\widetilde{g} \in \widetilde{G}(\gamma)$. Since $(\widetilde{G}(\gamma), \widetilde{H}(\gamma))$ is a closed (and strict) extension of the $S$-group $(G(\gamma), H(\gamma))$, there exist $g \in G(\gamma)$ and $\Phi \in H(\gamma)$ such that $\widetilde{g} = \widetilde{\Phi}(g)$. Now we calculate:

$$\begin{aligned}
\widetilde{\Theta}_{(\gamma'',\gamma)}(\widetilde{g}) &= \widetilde{\Theta}_{(\gamma'',\gamma)}\Big(\widetilde{\Phi}(g)\Big) = \widetilde{i}_{(\gamma'',\gamma)}(\widetilde{\Phi})\Big(\widetilde{\Theta}_{(\gamma'',\gamma)}(g)\Big) = \\
&= \widetilde{i}_{(\gamma'',\gamma)}(\widetilde{\Phi})\Big(\Theta_{(\gamma'',\gamma)}(g)\Big) = \widetilde{i}_{(\gamma'',\gamma)}(\widetilde{\Phi})\Big(\Theta_{(\gamma'',\gamma')}\Big(\Theta_{(\gamma',\gamma)}(g)\Big)\Big) = \\
&= \widetilde{i}_{(\gamma'',\gamma)}(\widetilde{\Phi})\Big(\widetilde{\Theta}_{(\gamma'',\gamma')}\Big(\Theta_{(\gamma',\gamma)}(g)\Big)\Big) = \widetilde{\Theta}_{(\gamma'',\gamma')}\Big(\widetilde{\varphi}\Big(\Theta_{(\gamma',\gamma)}(g)\Big)\Big).
\end{aligned}$$

Hence,

$$\widetilde{\Theta}_{(\gamma'',\gamma)}\left(\widetilde{g}\right) = \widetilde{\Theta}_{(\gamma'',\gamma')}\Big(\widetilde{\varphi}\Big(\Theta_{(\gamma',\gamma)}(g)\Big)\Big) \tag{4.4-1}$$

where $\widetilde{\varphi} \in \widetilde{H}(\gamma')$ is such that

$$\widetilde{i}_{(\gamma'',\gamma')}\left(\widetilde{\varphi}\right) = \widetilde{i}_{(\gamma'',\gamma)}\left(\widetilde{\Phi}\right). \tag{4.4-2}$$

From (4.4-1) we have now that:

$$\begin{aligned}
\widetilde{\Theta}_{(\gamma'',\gamma)}\left(\widetilde{g}\right) &= \widetilde{\Theta}_{(\gamma'',\gamma')}\Big(\widetilde{\varphi}\Big(\widetilde{\Theta}_{(\gamma',\gamma)}(g)\Big)\Big) = \\
&= \widetilde{\Theta}_{(\gamma'',\gamma')}\Big(\widetilde{\Theta}_{(\gamma',\gamma)}\Big(\widetilde{\Psi}(g)\Big)\Big),
\end{aligned} \tag{4.4-3}$$

where $\widetilde{\Psi} \in \widetilde{H}(\gamma)$ is such that

$$\widetilde{i}_{(\gamma',\gamma)}(\widetilde{\Psi}) = \widetilde{\varphi}. \tag{4.4-4}$$

From (4.4-4) we obtain:

$$\widetilde{i}_{(\gamma'',\gamma')}\Big(\widetilde{i}_{(\gamma',\gamma)}\left(\widetilde{\Psi}\right)\Big) = \widetilde{i}_{(\gamma'',\gamma')}\left(\widetilde{\varphi}\right),$$

that is,

$$\widetilde{i}_{(\gamma'',\gamma)}\left(\widetilde{\Psi}\right) = \widetilde{i}_{(\gamma'',\gamma')}\left(\widetilde{\varphi}\right),$$

or yet, taking (4.4-2) into account, that

$$\widetilde{i}_{(\gamma'',\gamma)}\left(\widetilde{\Psi}\right) = \widetilde{i}_{(\gamma'',\gamma)}\left(\widetilde{\Phi}\right).$$

But $\widetilde{i}_{(\gamma'',\gamma)}$ is an isomorphism and hence, from this last expression, we conclude that

$$\widetilde{\Psi} = \widetilde{\Phi}.$$

Taking this result into (4.4-3) we obtain:

$$\widetilde{\Theta}_{(\gamma'',\gamma)}\left(\widetilde{g}\right) = \widetilde{\Theta}_{(\gamma'',\gamma')}\left(\widetilde{\Theta}_{(\gamma',\gamma)}\left(\widetilde{\Phi}(g)\right)\right).$$

But $\widetilde{\Phi}(g) = \widetilde{g}$, and hence,

$$\widetilde{\Theta}_{(\gamma'',\gamma')}\left(\widetilde{\Theta}_{(\gamma',\gamma)}\left(\widetilde{g}\right)\right) = \widetilde{\Theta}_{(\gamma'',\gamma)}\left(\widetilde{g}\right)$$

which completes the proof. ■

## 4.5 Remark

In the Proposition 4.4, the specific family of $S$-groups

$$\widetilde{\mathbb{G}}\Big(\Gamma(I)\Big) = \Big\{\widetilde{\mathbb{G}}(\gamma) = \Big(\widetilde{G}(\gamma), \widetilde{H}(\gamma)\Big)\Big\}_{\gamma\in\Gamma(I)}$$

where, for each $\gamma \in \Gamma(I)$, $\widetilde{\mathbb{G}}(\gamma)$ is the strict and closed extension of the $S$-group $\mathbb{G}(\gamma) = (G(\gamma), H(\gamma))$ as defined in Theorem 2.16 (the Theorem of Extension of $S$-Groups), can be exchanged with any family of $S$-groups, let us say

$$\widehat{\mathbb{G}}\Big(\Gamma(I)\Big) = \Big\{\widehat{\mathbb{G}}(\gamma) = \Big(\widehat{G}(\gamma), \widehat{H}(\gamma)\Big)\Big\}_{\gamma\in\Gamma(I)},$$

such that, for each $\gamma \in \Gamma(I)$, the $S$-group $\widehat{\mathbb{G}}(\gamma)$ is a strict and closed extension of the $S$-group $\mathbb{G}(\gamma)$. Exchanging, in the referred proposition, $\widetilde{\mathbb{G}}(\Gamma(I))$ by $\widehat{\mathbb{G}}(\Gamma(I))$ as above and, also, $\widetilde{i}(\Gamma^2(I))$ by $\widehat{i}(\Gamma^2(I))$ (the extension of $i(\Gamma^2(I))$ to the family $\widehat{\mathbb{G}}(\Gamma(I))$), one obtains a new proposition which admits, *mutadis mutandis*, the same proof presented in 4.4.

# First Theorem of Extension of $S$-Spaces

## 4.6 The Extension $\widetilde{\mathscr{G}}(\mathcal{I})$ of $\mathscr{G}(I)$

Based on Proposition 4.4 we have then established the following result.

Let

$$\mathscr{G}(I) = \bigg(\mathbb{G}\Big(\Gamma(I)\Big), i\Big(\Gamma^2(I)\Big), \Theta\Big(\Delta(I)\Big)\bigg)$$

be as described in Proposition 4.4, that is, an abelian and surjective $S$-space with identity. We define now $\widetilde{\mathscr{G}}(I)$ as follows:

$$\widetilde{\mathscr{G}}(I) \coloneqq \left(\widetilde{\mathbb{G}}\Big(\Gamma(I)\Big), \widetilde{\imath}\Big(\Gamma^2(I)\Big), \widetilde{\Theta}\Big(\Delta(I)\Big)\right)$$

with:

- $\widetilde{\mathbb{G}}(\Gamma(I)) = \{\widetilde{\mathbb{G}}(\gamma) = (\widetilde{G}(\gamma), \widetilde{H}(\gamma))\}_{\gamma\in\Gamma(I)}$ the family of $S$-groups where, for each $\gamma \in \Gamma(I)$, $\widetilde{\mathbb{G}}(\gamma)$ is the strict and closed extension (of the $S$-group $\mathbb{G}(\gamma) \in \mathbb{G}(\Gamma(I))$) described in the Theorem of Extension of $S$-Groups (Theorem 2.16);
- $\widetilde{\imath}(\Gamma^2(I))$ is the extension to the family of $S$-groups $\widetilde{\mathbb{G}}(\Gamma(I))$ of the bonding $i(\Gamma^2(I))$ of the $S$-space $\mathscr{G}(I)$;
- $\widetilde{\Theta}(\Delta(I))$ is the unique restriction for the bonded family $(\widetilde{\mathbb{G}}\,(\Gamma(I))\,, \widetilde{\imath}\,(\Gamma^2(I)))$ that is a prolongation, to this bonded family, of the restriction $\Theta(\Delta(I))$ of the $S$-space $\mathscr{G}(I)$.

Clearly, $\widetilde{\mathscr{G}}(I)$ is a strict and closed extension of the $S$-space $\mathscr{G}(I)$.

It is worth observing that the components of the extension $\widetilde{\mathscr{G}}(I)$ above defined, that is, the family of $S$-groups $\widetilde{\mathbb{G}}(\Gamma(I))$, the bonding $\widetilde{\imath}(\Gamma^2(I))$ and the restriction $\widetilde{\Theta}(\Delta(I))$ are unique in the following sense:

- $\widetilde{\mathbb{G}}(\Gamma(I)) = \{\widetilde{\mathbb{G}}(\gamma) = (\widetilde{G}(\gamma), \widetilde{H}(\gamma))\}_{\gamma\in\Gamma(I)}$ is the only strict and closed extension of the family of $S$-groups $\mathbb{G}(\Gamma(I)) = \{\mathbb{G}(\gamma) = (G(\gamma), H(\gamma))\}_{\gamma\in\Gamma(I)}$ (of $\mathscr{G}(I)$) in the sense of, for each $\gamma \in \Gamma(I)$, by the Theorem of Extension of $S$-Groups (Theorem 2.16), unless isomorphism, $\widetilde{\mathbb{G}}(\gamma) = (\widetilde{G}(\gamma), \widetilde{H}(\gamma))$ is the only strict and closed extension of the $S$-group $\mathbb{G}(\gamma) = (G(\gamma), H(\gamma))$;
- the family $\widetilde{\imath}(\Gamma^2(I))$ — the extension of the family $i(\Gamma^2(I))$ to the family of $S$-groups $\widetilde{\mathbb{G}}(\Gamma(I))$ — is, taking Definition 3.22 into account, uniquely determined by the families $\widetilde{\mathbb{G}}(\Gamma(I))$ and $i(\Gamma^2(I))$ which, in turn, are unique — the first, $\widetilde{\mathbb{G}}(\Gamma(I))$, in the sense above described, and the second, $i(\Gamma^2(I))$, for being **the** specific bonding of the family $\mathbb{G}(\Gamma(I))$ of the $S$-space $\mathscr{G}(I)$;
- finally we have the uniqueness of $\widetilde{\Theta}(\Delta(I))$ such as established by Proposition 4.4.

Given the observations above and our well-defined notion of isomorphism for $S$-spaces (Definition 3.28), we inquire: Are there strict and closed extensions of the $S$-space $\mathscr{G}(I)$ that are not isomorphic to $\widetilde{\mathscr{G}}(I)$? Or, equivalently: $\widetilde{\mathscr{G}}(I)$ is, unless isomorphism, the only strict and closed extension of the abelian, surjective, and with identity $S$-space $\mathscr{G}(I)$?

The answer to this question is provided by the theorem below which, as we already announced in 4.1, will be referred as **1st Theorem of Extension of $S$-Spaces** or, shorter, **1st TESS**.

## 4.7 Theorem (1st TESS)

*Let $\mathscr{G}(I)$ and $\widetilde{\mathscr{G}}(I)$ be the $S$-spaces described in 4.6. $\widetilde{\mathscr{G}}(I)$ is, unless isomorphism, the only strict and closed extension of $\mathscr{G}(I)$.*

*Proof.* Let us suppose that, such as $\widetilde{\mathscr{G}}(I)$,

$$\widehat{\mathscr{G}}(I) = \left(\widehat{\mathbb{G}}\Big(\Gamma(I)\Big) = \left\{\widehat{\mathbb{G}}(\gamma) = \Big(\widehat{G}(\gamma), \widehat{H}(\gamma)\Big)\right\}_{\gamma\in\Gamma(I)}, \widehat{\imath}\Big(\Gamma^2(I)\Big), \widehat{\Theta}\Big(\Delta(I)\Big)\right)$$

is a strict and closed extension of the $S$-space $\mathscr{G}(I)$. Hence being, according to the definition of strict and closed extension of a $S$-space, Definition 3.23, the families of $S$-groups $\widetilde{\mathbb{G}}(\Gamma(I))$ and $\widehat{\mathbb{G}}(\Gamma(I))$ are, both, strict and closed extensions of the family of $S$-groups $\mathbb{G}(\Gamma(I))$, what, in turn, according to Definition 3.8, means that, for each $\gamma \in \Gamma(I)$, the $S$-groups $\widetilde{\mathbb{G}}(\gamma) = (\widetilde{G}(\gamma), \widetilde{H}(\gamma))$ and $\widehat{\mathbb{G}}(\gamma) = (\widehat{G}(\gamma), \widehat{H}(\gamma))$ are strict and closed extensions of the $S$-group $\mathbb{G}(\gamma) = (G(\gamma), H(\gamma))$ and hence, by what establishes the Theorem of Extension of $S$-Groups (Theorem 2.16), $\widetilde{\mathbb{G}}(\gamma)$ and $\widehat{\mathbb{G}}(\gamma)$ are isomorphic in the sense of the Definition 1.4(c). More specifically, the referred theorem tells us that, for each $\gamma \in \Gamma(I)$, the functions $\alpha_\gamma$ and $\beta_\gamma$ defined by (see Theorem 2.16)

$$\begin{aligned} \alpha_\gamma : \quad \widehat{G}(\gamma) &\longrightarrow \widetilde{G}(\gamma) = [H(\gamma) \times G(\gamma)]_\sim \\ \widehat{g} = \widehat{\Phi}(g) &\longmapsto \alpha_\gamma(\widehat{g}) := [\Phi, g] \end{aligned}$$

and

$$\begin{aligned} \beta_\gamma : \widehat{H}(\gamma) &\longrightarrow \widetilde{H}(\gamma) \\ \widehat{\Phi} &\longmapsto \beta_\gamma(\widehat{\Phi}) := \alpha_\gamma \widehat{\Phi} \alpha_\gamma^{-1} = \widetilde{\Phi}, \end{aligned}$$

are isomorphisms between groups ($\widehat{G}(\gamma)$ and $\widetilde{G}(\gamma)$) and semigroups ($\widehat{H}(\gamma)$ and $\widetilde{H}(\gamma)$), respectively, such that

**(a)** $\alpha_\gamma\left(\widehat{\Phi}(\widehat{g})\right) = \beta_\gamma\left(\widehat{\Phi}\right)\left(\alpha_\gamma(\widehat{g})\right)$ for every $\widehat{\Phi} \in \widehat{H}(\gamma)$ and every $\widehat{g} \in \widehat{G}(\gamma)$.

Furthermore, yet by Theorem 2.16, $\alpha_\gamma$ is the only isomorphism from the group $\widehat{G}(\gamma)$ onto the group $\widetilde{G}(\gamma)$ that keeps fixed the elements of $G(\gamma)$ (that is, $\alpha_\gamma(g) = [I_{G(\gamma)}, g] \equiv g$ if $g \in G(\gamma)$).

By the definition of isomorphism for $S$-space, Definition 3.28, in order for the $S$-space $\widehat{\mathscr{G}}(I)$ to be isomorphic to $\widetilde{\mathscr{G}}(I)$, beyond the condition (a) above, which as we saw is satisfied, must also be attended the following requirements:

**(b)** $\alpha_{\gamma'}\left(\widehat{\Theta}_{(\gamma',\gamma)}(\widehat{g})\right) = \widetilde{\Theta}_{(\gamma',\gamma)}\left(\alpha_\gamma(\widehat{g})\right)$ for every $(\gamma',\gamma) \in \Delta(I)$ and every $\widehat{g} \in \widehat{G}(\gamma)$;

**(c)** $\beta_{\gamma'}\left(\widehat{i}_{(\gamma',\gamma)}(\widehat{\Phi})\right) = \widetilde{i}_{(\gamma',\gamma)}\left(\beta_\gamma(\widehat{\Phi})\right)$ for every $(\gamma',\gamma) \in \Gamma^2(I)$ and every $\widehat{\Phi} \in \widehat{H}(\gamma)$.

We will prove (c) first. From the definition of $\beta_\gamma$ we have that $\beta_\gamma\left(\widehat{\Phi}\right) = \widetilde{\Phi}$, that is,

$$\beta_\gamma\left(\widehat{\gamma}(\Phi)\right) = \widetilde{\gamma}\,(\Phi)$$

for every $\gamma \in \Gamma(I)$ and every $\Phi \in H(\gamma)$. Hence,

$$\beta_\gamma \widehat{\gamma} = \widetilde{\gamma} \quad \text{or} \quad \beta_\gamma = \widetilde{\gamma}(\widehat{\gamma})^{-1}.$$

Furthermore, from Definition 3.22 of extension of a bonding, we know that:

$$i_{(\gamma',\gamma)} = \left(\widehat{\gamma}'\right)^{-1} \widehat{i}_{(\gamma',\gamma)}\ \widehat{\gamma}$$

and

$$\widetilde{i}_{(\gamma',\gamma)} = \widetilde{\gamma}'\ i_{(\gamma',\gamma)}(\widetilde{\gamma})^{-1}.$$

Now, with these data, we calculate:

$$\begin{aligned}
\beta_{\gamma'}\left(\widehat{i}_{(\gamma',\gamma)}(\widehat{\Phi})\right) &= \widetilde{\gamma}'\left(\left(\widehat{\gamma}'\right)^{-1}\left(\widehat{i}_{(\gamma',\gamma)}(\widehat{\Phi})\right)\right) = \\
&= \widetilde{\gamma}'\left(\left(\widehat{\gamma}'\right)^{-1}\left(\widehat{i}_{(\gamma',\gamma)}\left(\widehat{\gamma}(\Phi)\right)\right)\right) = \\
&= \widetilde{\gamma}'\left(i_{(\gamma',\gamma)}\,(\Phi)\right)
\end{aligned}$$

and

$$\begin{aligned}
\widetilde{i}_{(\gamma',\gamma)}\left(\beta_\gamma(\widehat{\Phi})\right) &= \widetilde{i}_{(\gamma',\gamma)}(\widetilde{\Phi}) = \widetilde{i}_{(\gamma',\gamma)}\left(\widetilde{\gamma}(\Phi)\right) = \\
&= \widetilde{\gamma}'\left(i_{(\gamma',\gamma)}\left((\widetilde{\gamma})^{-1}\left(\widetilde{\gamma}(\Phi)\right)\right)\right) = \\
&= \widetilde{\gamma}'\left(i_{(\gamma',\gamma)}(\Phi)\right)
\end{aligned}$$

and hence we have demonstrated the requirement (c).

It remains now to prove that the families $\widetilde{\Theta}(\Delta(I))$ and $\widehat{\Theta}(\Delta(I))$ satisfy the condition described in (b). In order to do so, let us first remember, as expressed in Remark 4.5, that Proposition 4.4 is applicable not only to the bonded family

$$\left(\widetilde{\mathbb{G}}\left(\Gamma(I)\right), \widetilde{i}\left(\Gamma^2(I)\right)\right)$$

of $\widetilde{\mathscr{G}}(I)$ but, also, to the bonded family $(\widehat{\mathbb{G}}\,(\Gamma(I))\,, \widehat{i}\,(\Gamma^2(I)))$ of $\widehat{\mathscr{G}}(I)$. Hence, according to the referred proposition, $\widehat{\Theta}(\Delta(I))$ is the only restriction for $(\widehat{\mathbb{G}}\,(\Gamma(I))\,, \widehat{i}\,(\Gamma^2(I)))$ that is also a prolongation, to this bonded family, of the restriction $\Theta(\Delta(I))$. Therefore, due to this uniqueness associated to $\widehat{\Theta}(\Delta(I))$, the proof of (b) or, equivalently, that

$$\widehat{\Theta}_{(\gamma',\gamma)} = \alpha_{\gamma'}^{-1}\widetilde{\Theta}_{(\gamma',\gamma)}\alpha_{\gamma},$$

can be obtained proving that, for $(\gamma',\gamma) \in \Delta(I)$ arbitrarily fixed, the function

$$\chi_{(\gamma',\gamma)} \colon \widehat{G}(\gamma) \longrightarrow \widehat{G}(\gamma')$$

defined by

$$\chi_{(\gamma',\gamma)}(\widehat{g}) := \alpha_{\gamma'}^{-1}\left(\widetilde{\Theta}_{(\gamma',\gamma)}\Big(\alpha_{\gamma}\left(\widehat{g}\right)\Big)\right)$$

for every $\widehat{g} \in \widehat{G}(\gamma)$, is a restriction for the bonded family $(\widehat{\mathbb{G}}\left(\Gamma(I)\right), \widehat{\imath}\left(\Gamma^2(I)\right))$ and, also, a prolongation (to this bonded family) of $\Theta(\Delta(I))$. In fact, if we do so, the Proposition 4.4 would allow us to conclude that

$$\chi_{(\gamma',\gamma)} = \widehat{\Theta}_{(\gamma',\gamma)}$$

that is,

$$\widehat{\Theta}_{(\gamma',\gamma)} = \alpha_{\gamma'}^{-1}\widetilde{\Theta}_{(\gamma',\gamma)}\alpha_{\gamma}$$

Let us then turn back our attention to this function $\chi_{(\gamma',\gamma)}$ and show that it has, among others, the following properties:

**(d)** $\chi_{(\gamma',\gamma)}$ is a homomorphism from the group $\widehat{G}(\gamma)$ into the group $\widehat{G}(\gamma')$, as can be seen directly from its definition;

**(e)** $\chi_{(\gamma',\gamma)}(g) = \Theta_{(\gamma',\gamma)}(g)$ for every $g \in G(\gamma)$;

In fact, since $\alpha_{\gamma} \colon \widehat{G}(\gamma) \longrightarrow \widetilde{G}(\gamma)$ keeps fixed the elements of $G(\gamma)$, we have

$$\begin{aligned}\chi_{(\gamma',\gamma)}(g) &= \alpha_{\gamma'}^{-1}\left(\widetilde{\Theta}_{(\gamma',\gamma)}\Big(\alpha_{\gamma}\left(g\right)\Big)\right) = \alpha_{\gamma'}^{-1}\Big(\widetilde{\Theta}_{(\gamma',\gamma)}\left(g\right)\Big) = \\ &= \alpha_{\gamma'}^{-1}\Big(\Theta_{(\gamma',\gamma)}\left(g\right)\Big) = \Theta_{(\gamma',\gamma)}\left(g\right)\end{aligned}$$

for every $g \in G(\gamma)$;

**(f)** $\chi_{(\gamma',\gamma)}\Big(\widehat{\Phi}\left(\widehat{g}\right)\Big) = \widehat{\imath}_{(\gamma',\gamma)}\left(\widehat{\Phi}\right)\Big(\chi_{(\gamma',\gamma)}\left(\widehat{g}\right)\Big)$ for every $\widehat{\Phi} \in \widehat{H}(\gamma)$ and every $\widehat{g} \in \widehat{G}(\gamma)$;

To prove this property we proceed as follows. Let $\widehat{g} \in \widehat{G}(\gamma)$ and $\widehat{\Phi} \in \widehat{H}(\gamma)$ be arbitrarily fixed. Since the $S$-group $\widehat{\mathbb{G}}(\gamma) = (\widehat{G}(\gamma), \widehat{H}(\gamma))$ is a closed (and strict) extension of the $S$-group $\mathbb{G}(\gamma) = (G(\gamma), H(\gamma))$, then

$$\widehat{g} = \widehat{\Psi}(g)$$

for some $\Psi \in H(\gamma)$ and $g \in G(\gamma)$, and hence

$$\widehat{\Phi}\left(\widehat{g}\right) = \widehat{\Phi}\left(\widehat{\Psi}(g)\right) = \widehat{\Phi\Psi}(g).$$

Thus, taking into account the definition of $\alpha_\gamma$ (at page 139), we have:

$$\alpha_\gamma\left(\widehat{\Phi}\left(\widehat{g}\right)\right) = \alpha_\gamma\left(\widehat{\Phi\Psi}(g)\right) = [\Phi\Psi, g].$$

On the other hand, according to the definition of $\widetilde{\Phi}$ (see Proposition 2.10),

$$\widetilde{\Phi}\left([\Psi, g]\right) = [\Phi\Psi, g].$$

Therefore,

$$\alpha_\gamma\left(\widehat{\Phi}\left(\widehat{g}\right)\right) = \widetilde{\Phi}\left([\Psi, g]\right).$$

Since $\widehat{g} = \widehat{\Psi}(g)$, it results from the definition of $\alpha_\gamma$ that

$$\alpha_\gamma\left(\widehat{g}\right) = [\Psi, g]$$

and, consequently,

$$\alpha_\gamma\left(\widehat{\Phi}\left(\widehat{g}\right)\right) = \widetilde{\Phi}\left(\alpha_\gamma\left(\widehat{g}\right)\right).$$

Thus,

$$\widetilde{\Theta}_{(\gamma',\gamma)}\left(\alpha_\gamma\left(\widehat{\Phi}\left(\widehat{g}\right)\right)\right) = \widetilde{\Theta}_{(\gamma',\gamma)}\left(\widetilde{\Phi}\left(\alpha_\gamma\left(\widehat{g}\right)\right)\right)$$

or yet,

$$\widetilde{\Theta}_{(\gamma',\gamma)}\left(\alpha_\gamma\left(\widehat{\Phi}\left(\widehat{g}\right)\right)\right) = \widetilde{i}_{(\gamma',\gamma)}(\widetilde{\Phi})\left(\widetilde{\Theta}_{(\gamma',\gamma)}\left(\alpha_\gamma\left(\widehat{g}\right)\right)\right),$$

from where we obtain

$$\alpha_{\gamma'}^{-1}\left(\widetilde{\Theta}_{(\gamma',\gamma)}\left(\alpha_\gamma\left(\widehat{\Phi}\left(\widehat{g}\right)\right)\right)\right) = \alpha_{\gamma'}^{-1}\left(\widetilde{i}_{(\gamma',\gamma)}(\widetilde{\Phi})\left(\widetilde{\Theta}_{(\gamma',\gamma)}\left(\alpha_\gamma\left(\widehat{g}\right)\right)\right)\right).$$

Taking now into account the definition of $\chi_{(\gamma',\gamma)}$, the last expression assumes the form

$$\chi_{(\gamma',\gamma)}\left(\widehat{\Phi}\left(\widehat{g}\right)\right) = \alpha_{\gamma'}^{-1}\left(\widetilde{i}_{(\gamma',\gamma)}(\widetilde{\Phi})\left(\widetilde{\Theta}_{(\gamma',\gamma)}\left(\alpha_\gamma\left(\widehat{g}\right)\right)\right)\right)$$

or yet

$$\chi_{(\gamma',\gamma)}\left(\widehat{\Phi}\left(\widehat{g}\right)\right) = \alpha_{\gamma'}^{-1}\left(\widetilde{i}_{(\gamma',\gamma)}(\widetilde{\Phi})\left(\alpha_{\gamma'}\left(\alpha_{\gamma'}^{-1}\left(\widetilde{\Theta}_{(\gamma',\gamma)}\left(\alpha_\gamma\left(\widehat{g}\right)\right)\right)\right)\right)\right)$$

from where, remembering that

$$\alpha_{\gamma'}^{-1}\left(\widetilde{\Theta}_{(\gamma',\gamma)}\left(\alpha_\gamma\left(\widehat{g}\right)\right)\right) = \chi_{(\gamma',\gamma)}\left(\widehat{g}\right),$$

one obtains that

$$\chi_{(\gamma',\gamma)}\left(\widehat{\Phi}\left(\widehat{g}\right)\right) = \alpha_{\gamma'}^{-1}\left(\widetilde{i}_{(\gamma',\gamma)}(\widetilde{\Phi})\left(\alpha_{\gamma'}\left(\chi_{(\gamma',\gamma)}\left(\widehat{g}\right)\right)\right)\right).$$

On the other hand, from the definition of the isomorphism $\beta_{\gamma'}\colon \widehat{H}(\gamma') \longrightarrow \widetilde{H}(\gamma')$ (at page 139) we have

$$\alpha_{\gamma'}^{-1}\left(\widetilde{i}_{(\gamma',\gamma)}(\widetilde{\Phi})\right)\alpha_{\gamma'} = \beta_{\gamma'}^{-1}\left(\widetilde{i}_{(\gamma',\gamma)}(\widetilde{\Phi})\right).$$

Hence being, we can write that

$$\chi_{(\gamma',\gamma)}\left(\widehat{\Phi}\left(\widehat{g}\right)\right) = \beta_{\gamma'}^{-1}\left(\widetilde{i}_{(\gamma',\gamma)}(\widetilde{\Phi})\right)\left(\chi_{(\gamma',\gamma)}\left(\widehat{g}\right)\right).$$

But, from (c) (at page 140) we obtain

$$\widehat{i}_{(\gamma',\gamma)}(\widehat{\Phi}) = \beta_{\gamma'}^{-1}\left(\widetilde{i}_{(\gamma',\gamma)}\left(\beta_{\gamma}(\widehat{\Phi})\right)\right).$$

But $\beta_{\gamma}(\widehat{\Phi}) = \widetilde{\Phi}$, and hence

$$\widehat{i}_{(\gamma',\gamma)}(\widehat{\Phi}) = \beta_{\gamma'}^{-1}\left(\widetilde{i}_{(\gamma',\gamma)}(\widetilde{\Phi})\right),$$

a result that, introduced into the latter expression of $\chi_{(\gamma',\gamma)}$, give us, finally,

$$\chi_{(\gamma',\gamma)}\left(\widehat{\Phi}\left(\widehat{g}\right)\right) = \widehat{i}_{(\gamma',\gamma)}(\widehat{\Phi})\left(\chi_{(\gamma',\gamma)}\left(\widehat{g}\right)\right);$$

**(g)** $\chi_{(\gamma'',\gamma')}(\chi_{(\gamma',\gamma)}(\widehat{g})) = \chi_{(\gamma'',\gamma)}(\widehat{g})$ for every $\widehat{g} \in \widehat{G}(\gamma)$ and any $\gamma, \gamma', \gamma'' \in \Gamma(I)$ such that $\gamma'' \subseteq \gamma' \subseteq \gamma$;

In fact, since by the definition of $\chi_{(\gamma'',\gamma)}$ we have:

$$\chi_{(\gamma'',\gamma)}\left(\widehat{g}\right) = \alpha_{\gamma''}^{-1}\left(\widetilde{\Theta}_{(\gamma'',\gamma)}\left(\alpha_{\gamma}\left(\widehat{g}\right)\right)\right).$$

Since

$$\widetilde{\Theta}_{(\gamma'',\gamma)}\left(\alpha_{\gamma}\left(\widehat{g}\right)\right) = \widetilde{\Theta}_{(\gamma'',\gamma')}\left(\widetilde{\Theta}_{(\gamma',\gamma)}\left(\alpha_{\gamma}\left(\widehat{g}\right)\right)\right),$$

we obtain:

$$\begin{aligned}\chi_{(\gamma'',\gamma)}\left(\widehat{g}\right) &= \alpha_{\gamma''}^{-1}\left(\widetilde{\Theta}_{(\gamma'',\gamma')}\left(\widetilde{\Theta}_{(\gamma',\gamma)}\left(\alpha_{\gamma}\left(\widehat{g}\right)\right)\right)\right) = \\ &= \alpha_{\gamma''}^{-1}\left(\widetilde{\Theta}_{(\gamma'',\gamma')}\left(\alpha_{\gamma'}\left(\alpha_{\gamma'}^{-1}\left(\widetilde{\Theta}_{(\gamma',\gamma)}\left(\alpha_{\gamma}\left(\widehat{g}\right)\right)\right)\right)\right)\right) = \\ &= \chi_{(\gamma'',\gamma')}\left(\chi_{(\gamma',\gamma)}\left(\widehat{g}\right)\right).\end{aligned}$$

The properties (d) to (g) of the functions $\chi_{(\gamma',\gamma)}\colon \widehat{G}_\gamma \longrightarrow \widehat{G}_{\gamma'}$, with $(\gamma',\gamma) \in \Delta(I)$, above proven, allow us to conclude that the family

$$\chi\Big(\Delta(I)\Big) := \Big\{\chi_{(\gamma',\gamma)}\Big\}_{(\gamma',\gamma)\in\Delta(I)}$$

is a restriction for the bonded family $(\widehat{\mathbb{G}}\,(\Gamma(I))\,,\widehat{\imath}\,(\Gamma^2(I)))$ and, also, a prolongation, to this bonded family, of $\Theta(\Delta(I))$ (recap the definitions 3.13 and 3.23(a)). But, as we saw, there exists only one such family, namely, $\widehat{\Theta}(\Delta(I))$. Therefore,

$$\widehat{\Theta}\Big(\Delta(I)\Big) = \chi\Big(\Delta(I)\Big),$$

that is, for every $(\gamma',\gamma) \in \Delta(I)$ we have:

$$\widehat{\Theta}_{(\gamma',\gamma)} = \alpha_{\gamma'}^{-1}\widetilde{\Theta}_{(\gamma',\gamma)}\alpha_\gamma.$$

With this we conclude the proof[17] of the 1st TESS. ■

# The Equivalence Relation $\approx$ and the Set $\overline{G}(\gamma)$

## 4.8 Remarks

The first among the two goals of this chapter, both described in 4.1, was achieved with the proof of the 1st TESS. We start now our journey towards the second goal, namely, to prove a theorem (2nd TESS) that, essentially, establishes the following:

- Being $\widetilde{\mathscr{G}}(I)$ a strict and closed extension of a given $S$-space, $\mathscr{G}(I)$, abelian, surjective, with identity and coherent, then, there exists, unless isomorphism, a unique extension of $\widetilde{\mathscr{G}}(I)$ that is locally closed and coherent.

We will proceed in an analogous manner to that adopted in the proof of the Theorem of Extension of $S$-Groups (Theorem 2.16); more explicitly, guided step-by-step by the hypothesis of the existence of a locally closed and coherent extension of $\widetilde{\mathscr{G}}(I)$, let us say $\widehat{\mathscr{G}}(I)$, we will build a specific extension

$$\overline{\mathscr{G}}(I) = \left(\overline{\mathbb{G}}\Big(\Gamma(I)\Big) = \Big\{\overline{\mathbb{G}}(\gamma) = \Big(\overline{G}(\gamma), \overline{H}(\gamma)\Big)\Big\}_{\gamma\in\Gamma(I)}, \overline{\imath}\Big(\Gamma^2(I)\Big), \overline{\Theta}\Big(\Delta(I)\Big)\right)$$

of $\widetilde{\mathscr{G}}(I)$, with the required properties and isomorphic to the hypothetical extension $\widehat{\mathscr{G}}(I)$. The idea is to obtain from $\widehat{\mathscr{G}}(I)$, more precisely from the properties assumed for this hypothetical extension of $\widetilde{\mathscr{G}}(I)$, indicatives (suggestions) of how to engineer the elements of $\overline{\mathscr{G}}(I)$: the $S$-groups of the family $\overline{\mathbb{G}}(\Gamma(I))$, the members of $\overline{\Theta}(\Delta(I))$, etc.

[17] Observe that the fourth condition of the definition of restriction, the requirement (d) of the Definition 3.13, is trivially satisfied in the context of our proof, since the involved homomorphisms, $\widehat{\Phi}$, are endomorphisms.

Throughout the rest of this chapter,

$$\mathscr{G}(I) = \left(\mathbb{G}\Big(\Gamma(I)\Big) = \Big\{\mathbb{G}(\gamma) = \Big(G(\gamma), H(\gamma)\Big)\Big\}_{\gamma\in\Gamma(I)}, i\Big(\Gamma^2(I)\Big), \Theta\Big(\Delta(I)\Big)\right)$$

and

$$\widetilde{\mathscr{G}}(I) = \left(\widetilde{\mathbb{G}}\Big(\Gamma(I)\Big) = \Big\{\widetilde{\mathbb{G}}(\gamma) = \Big(\widetilde{G}(\gamma), \widetilde{H}(\gamma)\Big)\Big\}_{\gamma\in\Gamma(I)}, \widetilde{i}\Big(\Gamma^2(I)\Big), \widetilde{\Theta}\Big(\Delta(I)\Big)\right)$$

will denote, respectively, an abelian, surjective, with identity and coherent $S$-space, and its strict and closed extension as defined in 4.6. Furthermore,

$$\widehat{\mathscr{G}}(I) = \left(\widehat{\mathbb{G}}\Big(\Gamma(I)\Big) = \Big\{\widehat{\mathbb{G}}(\gamma) = \Big(\widehat{G}(\gamma), \widehat{H}(\gamma)\Big)\Big\}_{\gamma\in\Gamma(I)}, \widehat{i}\Big(\Gamma^2(I)\Big), \widehat{\Theta}\Big(\Delta(I)\Big)\right)$$

will denote a hypothetical locally closed and coherent extension of $\widetilde{\mathscr{G}}(I)$. Any reference to $\mathscr{G}(I)$, $\widetilde{\mathscr{G}}(I)$ or $\widehat{\mathscr{G}}(I)$ henceforth done in the next sections of this chapter, will incorporate the hypotheses above unless otherwise explicitly said.

With some frequency throughout our developments, sometimes without any explicit reference, we will use one or more among the "inherited properties" (a) to (f) listed in 3.30 and, thus, we recommend the reader to revisit them.

## 4.9 Motivation

As we said in the previous item, our goal is to build a $S$-space,

$$\overline{\mathscr{G}}(I) = \left(\overline{\mathbb{G}}\Big(\Gamma(I)\Big) = \Big\{\overline{\mathbb{G}}(\gamma) = \Big(\overline{G}(\gamma), \overline{H}(\gamma)\Big)\Big\}_{\gamma\in\Gamma(I)}, \overline{i}\Big(\Gamma^2(I)\Big), \overline{\Theta}\Big(\Delta(I)\Big)\right),$$

which is a locally closed and coherent extension of $\widetilde{\mathscr{G}}(I)$, isomorphic to $\widehat{\mathscr{G}}(I)$. Hence, our first task must be the one of defining, for each $\gamma \in \Gamma(I)$, the set $\overline{G}(\gamma)$; Which ones could be the members of $\overline{G}(\gamma)$?

Clearly, $\overline{G}(\gamma)$ must be such that there exists between it and the set $\widehat{G}(\gamma)$ a bijection, since for $\overline{\mathscr{G}}(\gamma)$ to be isomorphic to $\widehat{\mathscr{G}}(\gamma)$ it is necessary that, for each $\gamma \in \Gamma(I)$, the $S$-group $\overline{\mathbb{G}}(\gamma) = (\overline{G}(\gamma), \overline{H}(\gamma))$ to be isomorphic to the $S$-group $\widehat{\mathbb{G}}(\gamma) = (\widehat{G}(\gamma), \widehat{H}(\gamma))$ and, therefore, the groups $\overline{G}(\gamma)$ and $\widehat{G}(\gamma)$ must be isomorphic, which requires the same cardinality for the sets $\overline{G}(\gamma)$ and $\widehat{G}(\gamma)$.

On the other hand, it must be clear that the elements of $\overline{G}(\gamma)$ should be "constructed" from what we have, that is, from the objects in $\mathscr{G}(I)$ and $\widetilde{\mathscr{G}}(I)$; $\widehat{\mathscr{G}}(I)$ is only a hypothetical $S$-space that guide the "construction" of $\overline{\mathscr{G}}(\gamma)$.

Let $\gamma \in \Gamma(I)$ be an arbitrarily chosen element such that $\gamma \neq \varnothing$.[18] Let also $\widetilde{F}(\gamma)$ be the set whose members are, exactly, the coherent families in $\widetilde{\mathscr{G}}(\gamma)$ (recap Definition 3.19 and the notation introduced in 3.20). Hence, $x \in \widetilde{F}(\gamma)$ if and only if $x$ is a coherent family in $\widetilde{\mathscr{G}}(\gamma)$, that is,

$$x = \widetilde{g}\left(\overline{\xi}\,(\gamma)\right) = \left\{\widetilde{g}_\xi\right\}_{\xi \in \overline{\xi}(\gamma)}$$

where $\overline{\xi}(\gamma) \subseteq \Gamma(\gamma)$ is a cover of $\gamma$, $\widetilde{g}_\xi \in \widetilde{G}(\xi)$ for each $\xi \in \overline{\xi}(\gamma)$ and

$$\widetilde{\Theta}_{(\xi \cap \xi', \xi)}\left(\widetilde{g}_\xi\right) = \widetilde{\Theta}_{(\xi \cap \xi', \xi')}\left(\widetilde{g}_{\xi'}\right)$$

for every $\xi, \xi' \in \overline{\xi}(\gamma)$ such that $\xi \cap \xi' \neq \varnothing$.

Since $\widehat{\mathscr{G}}(\gamma)$ is an extension of $\widetilde{\mathscr{G}}(\gamma)$ (according to the "inherited property" 3.30(c)), the members of $\widetilde{F}(\gamma)$ also are coherent families in $\widehat{\mathscr{G}}(\gamma)$. Now, $\widehat{\mathscr{G}}(\gamma)$ is a $S$-subspace of the coherent $S$-space $\widehat{\mathscr{G}}(I)$ and hence (according to 3.30(b)) it is also coherent. Therefore, for each $\widetilde{g}(\overline{\xi}(\gamma)) = \{\widetilde{g}_\xi\}_{\xi \in \overline{\xi}(\gamma)}$ in $\widetilde{F}(\gamma)$ there exists a unique $\widehat{g} \in \widehat{G}(\gamma)$ such that

$$\widehat{\Theta}_{(\xi,\gamma)}(\widehat{g}) = \widetilde{g}_\xi \quad \text{for every} \quad \xi \in \overline{\xi}(\gamma),$$

which supports the definition of the function $a_\gamma \colon \widetilde{F}(\gamma) \longrightarrow \widehat{G}(\gamma)$ given as follows:

$$a_\gamma\left(\widetilde{g}\left(\overline{\xi}(\gamma)\right) = \left\{\widetilde{g}_\xi\right\}_{\xi \in \overline{\xi}(\gamma)}\right) \coloneqq \widehat{g}$$

where $\widehat{g}$ is the only element of $\widehat{G}(\gamma)$ such that

$$\widehat{\Theta}_{(\xi,\gamma)}(\widehat{g}) = \widetilde{g}_\xi \quad \text{for every} \quad \xi \in \overline{\xi}(\gamma).$$

The function $a_\gamma \colon \widetilde{F}(\gamma) \longrightarrow \widehat{G}(\gamma)$ is surjective. In fact, let $\widehat{g} \in \widehat{G}(\gamma)$ be arbitrarily fixed. We know that $\widehat{\mathscr{G}}(\gamma)$ is a locally closed extension of the $S$-subspace $\widetilde{\mathscr{G}}(\gamma)$ (according to 3.30(e)) and hence, there exists a cover $\overline{\xi}(\gamma) \subseteq \Gamma(\gamma)$ of $\gamma$[19] such that

$$\widehat{\Theta}_{(\xi,\gamma)}(\widehat{g}) \in \widetilde{G}(\xi) \quad \text{for every} \quad \xi \in \overline{\xi}(\gamma).$$

Now, the family

$$\widetilde{g}\left(\overline{\xi}(\gamma)\right) \coloneqq \left\{\widehat{\Theta}_{(\xi,\gamma)}(\widehat{g})\right\}_{\xi \in \overline{\xi}(\gamma)}$$

is in $\widetilde{F}(\gamma)$, that is, is a coherent family in $\widetilde{\mathscr{G}}(\gamma)$ (as can be easily demonstrated) and yet, with the definition of the function $a_\gamma \colon \widetilde{F}(\gamma) \longrightarrow \widehat{G}(\gamma)$ in mind,

$$a_\gamma\left(\widetilde{g}\left(\overline{\xi}(\gamma)\right)\right) = \widehat{g}.$$

[18] Given an indexer $\Gamma(I)$ and an element $\gamma \in \Gamma(I)$ such that $\gamma \neq \varnothing$, we reserve the symbol "$\Gamma(\gamma)$", henceforth, to denote the indexer $\Gamma(I) \cap \mathcal{P}(\gamma)$, that is, $\Gamma(\gamma) \coloneqq \Gamma(I) \cap \mathcal{P}(\gamma)$.

[19] $\overline{\xi}(\gamma)$ can be obtained as follows: for each $x \in \gamma$, let $\gamma_x \in \Gamma(\gamma)$ be such that $x \in \gamma_x$ and $\widehat{\Theta}_{(\gamma_x,\gamma)}(\widehat{g}) \in \widetilde{G}(\gamma_x)$ (such $\gamma_x$ exists since $\widehat{\mathscr{G}}(\gamma)$ is a locally closed extension of $\widetilde{\mathscr{G}}(\gamma)$); now we take $\overline{\xi}(\gamma) \coloneqq \{\gamma_x \colon x \in \gamma\}$.

Since $\widehat{g} \in \widehat{G}(\gamma)$ is arbitrary, we conclude that $a_\gamma$ is surjective.

If, beyond surjective, $a_\gamma$ were also injective, then, $\widetilde{F}(\gamma)$ and $\widehat{G}(\gamma)$ would have the same cardinality and, hence, the set $\widetilde{F}(\gamma)$, whose definition is entirely independent of the existence of $\widehat{\mathscr{G}}(I)$, would be a candidate for the role of $\overline{G}(\gamma)$. However, $a_\gamma$ is not necessarily injective: nothing, *a priori*, prevents that for a given $\widehat{g} \in \widehat{G}(\gamma)$ there exist $\widetilde{g}(\overline{\xi}(\gamma))$ and $\widetilde{h}(\overline{\eta}(\gamma))$ in $\widetilde{F}(\gamma)$ such that $\widetilde{g}(\overline{\xi}(\gamma)) \neq \widetilde{h}(\overline{\eta}(\gamma))$ and $a_\gamma(\widetilde{g}(\overline{\xi}(\gamma))) = a_\gamma(\widetilde{h}(\overline{\eta}(\gamma)))$.

However, the binary relation $\approx$ in $\widetilde{F}(\gamma)$ defined by,

$$\widetilde{g}\Big(\overline{\xi}(\gamma)\Big) \approx \widetilde{h}\Big(\overline{\eta}(\gamma)\Big) \quad \text{if and only if} \quad a_\gamma\left(\widetilde{g}\Big(\overline{\xi}(\gamma)\Big)\right) = a_\gamma\left(\widetilde{h}\Big(\overline{\eta}(\gamma)\Big)\right)$$

for any $\widetilde{g}(\overline{\xi}(\gamma))$ and $\widetilde{h}(\overline{\eta}(\gamma))$ in $\widetilde{F}(\gamma)$, is, clearly, an equivalence relation in $\widetilde{F}(\gamma)$. Hence, denoting by $[\widetilde{F}(\gamma)]$ the partition of $\widetilde{F}(\gamma)$ determined by $\approx$, that is,

$$\Big[\widetilde{F}(\gamma)\Big] := \left\{\Big[\widetilde{g}\Big(\overline{\xi}(\gamma)\Big)\Big] : \widetilde{g}\Big(\overline{\xi}(\gamma)\Big) \in \widetilde{F}(\gamma)\right\}$$

where

$$\Big[\widetilde{g}\Big(\overline{\xi}(\gamma)\Big)\Big] := \left\{\widetilde{h}\Big(\overline{\eta}(\gamma)\Big) \in \widetilde{F}(\gamma) : \widetilde{h}\Big(\overline{\eta}(\gamma)\Big) \approx \widetilde{g}\Big(\overline{\xi}(\gamma)\Big)\right\},$$

we get the bijective function $A_\gamma$ defined as follows:

$$\begin{aligned} A_\gamma : \quad \Big[\widetilde{F}(\gamma)\Big] &\longrightarrow \widehat{G}(\gamma) \\ \Big[\widetilde{g}\Big(\overline{\xi}(\gamma)\Big)\Big] &\longmapsto A_\gamma\Big(\Big[\widetilde{g}\Big(\overline{\xi}(\gamma)\Big)\Big]\Big) := \widehat{g} \end{aligned}$$

where $\widehat{g} \in \widehat{G}(\gamma)$ is such that

$$a_\gamma\left(\widetilde{g}\Big(\overline{\xi}(\gamma)\Big)\right) = \widehat{g}.$$

We now have a postulant to the role of $\overline{G}(\gamma)$: the partition, $[\widetilde{F}(\gamma)]$, of the set $\widetilde{F}(\gamma)$ of the coherent families in $\widetilde{\mathscr{G}}(\gamma)$, determined by the equivalence relation, $\approx$, above defined. However, the definition of $\approx$ is given through the function $a_\gamma\colon \widetilde{F}(\gamma) \longrightarrow \widehat{G}(\gamma)$ which, in turn, involves the $S$-subspace $\widehat{\mathscr{G}}(\gamma)$ of $\widehat{\mathscr{G}}(I)$; we want to construct $\overline{G}(\gamma)$ only with what we really have, $\mathscr{G}(I)$ and $\widetilde{\mathscr{G}}(I)$. In this sense the Lemma 4.10 ahead is relevant, as it presents an equivalent definition for the relation $\approx$ which does not involve $\widehat{\mathscr{G}}(I)$.

## 4.10 Lemma

*Let $a_\gamma\colon \widetilde{F}(\gamma) \longrightarrow \widehat{G}(\gamma)$ and $\approx \subseteq \widetilde{F}(\gamma) \times \widetilde{F}(\gamma)$ be as defined in 4.9 (page 147). For any*

$$\widetilde{g}\Big(\overline{\xi}(\gamma)\Big) = \Big\{\widetilde{g}_\xi\Big\}_{\xi \in \overline{\xi}(\gamma)} \quad \textit{and} \quad \widetilde{h}\Big(\overline{\eta}(\gamma)\Big) = \Big\{\widetilde{h}_\eta\Big\}_{\eta \in \overline{\eta}(\gamma)}$$

*in* $\widetilde{F}(\gamma)$ *we have:*

$$\widetilde{g}\Big(\overline{\xi}(\gamma)\Big) \approx \widetilde{h}\Big(\overline{\eta}(\gamma)\Big) \quad \textit{if and only if} \quad \widetilde{\Theta}_{(\xi\cap\eta,\xi)}(\widetilde{g}_\xi) = \widetilde{\Theta}_{(\xi\cap\eta,\eta)}(\widetilde{h}_\eta)$$

*for every* $\xi \in \overline{\xi}(\gamma)$ *and* $\eta \in \overline{\eta}(\gamma)$ *such that* $\xi \cap \eta \neq \varnothing$*.*

*Proof.* Let $\widetilde{g}(\overline{\xi}(\gamma)) = \{\widetilde{g}_\xi\}_{\xi\in\overline{\xi}(\gamma)}$ and $\widetilde{h}(\overline{\eta}(\gamma)) = \{\widetilde{h}_\eta\}_{\eta\in\overline{\eta}(\gamma)}$ be elements of $\widetilde{F}(\gamma)$ arbitrarily fixed.

Suppose first that

$$\widetilde{g}\Big(\overline{\xi}(\gamma)\Big) \approx \widetilde{h}\Big(\overline{\eta}(\gamma)\Big).$$

Hence, from the definition of $\approx$, we get

$$a_\gamma\left(\widetilde{g}\Big(\overline{\xi}(\gamma)\Big)\right) = a_\gamma\left(\widetilde{h}\Big(\overline{\eta}(\gamma)\Big)\right)$$

or yet, taking into account the definition of $a_\gamma$, that

$$\widehat{g} := a_\gamma\left(\widetilde{g}\Big(\overline{\xi}(\gamma)\Big)\right) = a_\gamma\left(\widetilde{h}\Big(\overline{\eta}(\gamma)\Big)\right)$$

is the only element of $\widehat{G}(\gamma)$ such that

$$\widehat{\Theta}_{(\xi,\gamma)}(\widehat{g}) = \widetilde{g}_\xi \quad \text{for every} \quad \xi \in \overline{\xi}(\gamma)$$

and

$$\widehat{\Theta}_{(\eta,\gamma)}(\widehat{g}) = \widetilde{h}_\eta \quad \text{for every} \quad \eta \in \overline{\eta}(\gamma).$$

It then results from these two last expressions that, for $\xi$ and $\eta$ such that $\xi \cap \eta \neq \varnothing$,

$$\widehat{\Theta}_{(\xi\cap\eta,\xi)}\Big(\widehat{\Theta}_{(\xi,\gamma)}(\widehat{g})\Big) = \widehat{\Theta}_{(\xi\cap\eta,\xi)}(\widetilde{g}_\xi)$$

and

$$\widehat{\Theta}_{(\xi\cap\eta,\eta)}\Big(\widehat{\Theta}_{(\eta,\gamma)}(\widehat{g})\Big) = \widehat{\Theta}_{(\xi\cap\eta,\eta)}(\widetilde{h}_\eta).$$

Since $\widehat{\Theta}(\Delta(I))$ is a restriction for the bonded family $(\widehat{\mathbb{G}}\,(\Gamma(I))\,, \widehat{\imath}\,(\Gamma^2(I)))$, we get that:

$$\widehat{\Theta}_{(\xi\cap\eta,\xi)}\Big(\widehat{\Theta}_{(\xi,\gamma)}(\widehat{g})\Big) = \widehat{\Theta}_{(\xi\cap\eta,\eta)}\Big(\widehat{\Theta}_{(\eta,\gamma)}(\widehat{g})\Big) = \widehat{\Theta}_{(\xi\cap\eta,\gamma)}\,(\widehat{g})$$

and so we conclude that

$$\widehat{\Theta}_{(\xi\cap\eta,\xi)}(\widetilde{g}_\xi) = \widehat{\Theta}_{(\xi\cap\eta,\eta)}(\widetilde{h}_\eta)$$

for every $\xi \in \overline{\xi}(\gamma)$ and $\eta \in \overline{\eta}(\gamma)$ such that $\xi \cap \eta \neq \varnothing$.

But the restriction $\widehat{\Theta}(\Delta(I))$ for the bonded family $(\widehat{\mathbb{G}}\,(\Gamma(I))\,, \widehat{\imath}\,(\Gamma^2(I)))$ also is a prolongation of $\widetilde{\Theta}(\Delta(I))$ to this bonded family, for $\widehat{\mathscr{G}}(I)$ is an extension of $\widetilde{\mathscr{G}}(I)$. Thus,

$$\widehat{\Theta}_{(\xi\cap\eta,\xi)}(\widetilde{g}_\xi) = \widetilde{\Theta}_{(\xi\cap\eta,\xi)}(\widetilde{g}_\xi) \quad \text{and} \quad \widehat{\Theta}_{(\xi\cap\eta,\eta)}(\widetilde{h}_\eta) = \widetilde{\Theta}_{(\xi\cap\eta,\eta)}(\widetilde{h}_\eta)$$

which allows us to conclude that

$$\widetilde{\Theta}_{(\xi\cap\eta,\xi)}(\widetilde{g}_\xi) = \widetilde{\Theta}_{(\xi\cap\eta,\eta)}(\widetilde{h}_\eta) \tag{4.10-1}$$

for every $\xi \in \overline{\xi}(\gamma)$ and $\eta \in \overline{\eta}(\gamma)$ such that $\xi \cap \eta \neq \varnothing$.

Conversely, let us suppose now that $\widetilde{g}(\overline{\xi}(\gamma)) = \{\widetilde{g}_\xi\}_{\xi\in\overline{\xi}(\gamma)}$ and $\widetilde{h}(\overline{\eta}(\gamma)) = \{\widetilde{h}_\eta\}_{\eta\in\overline{\eta}(\gamma)}$ satisfy this last expression. Since they are coherent families in $\widetilde{\mathscr{G}}(\gamma)$, and that $\widehat{\mathscr{G}}(\gamma)$ is an extension of $\widetilde{\mathscr{G}}(\gamma)$, these families are also coherent in $\widehat{\mathscr{G}}(\gamma)$. Now, $\widehat{\mathscr{G}}(\gamma)$ is coherent (by property 3.30(b) and $\widehat{\mathscr{G}}(I)$ being coherent) and, therefore, there exist (unique) $\widehat{g}$ and $\widehat{h}$ in $\widehat{G}(\gamma)$ such that

$$\widehat{\Theta}_{(\xi,\gamma)}(\widehat{g}) = \widetilde{g}_\xi \quad \text{for every} \quad \xi \in \overline{\xi}(\gamma) \tag{4.10-2}$$

and

$$\widehat{\Theta}_{(\eta,\gamma)}(\widehat{h}) = \widetilde{h}_\eta \quad \text{for every} \quad \eta \in \overline{\eta}(\gamma). \tag{4.10-3}$$

Taking these expressions for $\widetilde{g}_\xi$ and $\widetilde{h}_\eta$ into equation (4.10-1), one obtains:

$$\widetilde{\Theta}_{(\xi\cap\eta,\xi)}\Big(\widehat{\Theta}_{(\xi,\gamma)}(\widehat{g})\Big) = \widetilde{\Theta}_{(\xi\cap\eta,\eta)}\Big(\widehat{\Theta}_{(\eta,\gamma)}(\widehat{h})\Big),$$

or yet, taking again into account that $\widehat{\Theta}(\Delta(I))$ is a prolongation of $\widetilde{\Theta}(\Delta(I))$, that

$$\widehat{\Theta}_{(\xi\cap\eta,\xi)}\Big(\widehat{\Theta}_{(\xi,\gamma)}(\widehat{g})\Big) = \widehat{\Theta}_{(\xi\cap\eta,\eta)}\Big(\widehat{\Theta}_{(\eta,\gamma)}(\widehat{h})\Big),$$

and hence, since $\widehat{\Theta}(\Delta(I))$ is a restriction,

$$\widehat{\Theta}_{(\xi\cap\eta,\gamma)}(\widehat{g}) = \widehat{\Theta}_{(\xi\cap\eta,\gamma)}(\widehat{h})$$

or, equivalently,

$$\widehat{\Theta}_{(\xi\cap\eta,\gamma)}(\widehat{g} - \widehat{h}) = \widehat{0}_{\xi\cap\eta}$$

for every $\xi \in \overline{\xi}(\gamma)$ and $\eta \in \overline{\eta}(\gamma)$ such that $\xi \cap \eta \neq \varnothing$, being $\widehat{0}_{\xi\cap\eta} \in \widehat{G}(\xi \cap \eta)$ the additive neutral of the group $\widehat{G}(\xi \cap \eta)$. This last conclusion can also be written as:

$$\widehat{\Theta}_{(\chi,\gamma)}(\widehat{g} - \widehat{h}) = \widehat{0}_\chi$$

for every $\chi \in \overline{\chi}(\gamma)$, being $\overline{\chi}(\gamma) \subseteq \Gamma(\gamma)$ the cover of $\gamma$ given by

$$\overline{\chi}(\gamma) := \Big\{\chi := \xi \cap \eta \quad : \quad \xi \in \overline{\xi}(\gamma), \quad \eta \in \overline{\eta}(\gamma) \quad \text{and} \quad \xi \cap \eta \neq \varnothing\Big\}.$$

On the other hand, since $\widehat{\Theta}_{(\chi,\gamma)} \colon \widehat{G}(\gamma) \longrightarrow \widehat{G}(\chi)$ is a homomorphism, we have

$$\widehat{\Theta}_{(\chi,\gamma)}(\widehat{0}_\gamma) = \widehat{0}_\chi$$

where $\widehat{0}_\gamma \in \widehat{G}(\gamma)$ is the additive neutral of $\widehat{G}(\gamma)$. Hence, we have

$$\widehat{\Theta}_{(\chi,\gamma)}(\widehat{g} - \widehat{h}) = \widehat{\Theta}_{(\chi,\gamma)}(\widehat{0}_\gamma)$$

for every $\chi \in \overline{\chi}(\gamma)$, from where one concludes, taking the Lemma 3.32 into account, that

$$\widehat{g} - \widehat{h} = \widehat{0}_{\gamma}$$

that is, $\widehat{g} = \widehat{h}$.

From this, the expressions (4.10-2) and (4.10-3) can be rewritten as:

$$\widehat{\Theta}_{(\xi,\gamma)}(\widehat{g}) = \widetilde{g}_{\xi} \quad \text{for every} \quad \xi \in \overline{\xi}(\gamma)$$

and

$$\widehat{\Theta}_{(\eta,\gamma)}(\widehat{g}) = \widetilde{h}_{\eta} \quad \text{for every} \quad \eta \in \overline{\eta}(\gamma),$$

which, in turn, with the definition of the function $a_{\gamma} \colon \widetilde{F}(\gamma) \longrightarrow \widehat{G}(\gamma)$ in mind (see page 146), tells us that

$$a_{\gamma}\left(\widetilde{g}\left(\overline{\xi}(\gamma)\right)\right) = \widehat{g}$$

and

$$a_{\gamma}\left(\widetilde{h}\left(\overline{\eta}(\gamma)\right)\right) = \widehat{g},$$

that is,

$$a_{\gamma}\left(\widetilde{g}\left(\overline{\xi}(\gamma)\right)\right) = a_{\gamma}\left(\widetilde{h}\left(\overline{\eta}(\gamma)\right)\right)$$

and hence, by the definition of the binary relation $\approx$ (see page 147), we finally have that

$$\widetilde{g}\left(\overline{\xi}(\gamma)\right) \approx \widetilde{h}\left(\overline{\eta}(\gamma)\right).$$

■

## 4.11 Remarks About the Lemma 4.10

In item 4.9 we defined (see page 146), based on the hypothesis of existence of $\widehat{\mathscr{G}}(I)$, the function $a_{\gamma} \colon \widetilde{F}(\gamma) \longrightarrow \widehat{G}(\gamma)$ $(\gamma \in \Gamma(I))$ as follows:

$$a_{\gamma}\left(\widetilde{g}\left(\overline{\xi}(\gamma)\right) = \left\{\widetilde{g}_{\xi}\right\}_{\xi \in \overline{\xi}(\gamma)}\right) := \widehat{g}$$

where $\widehat{g} \in \widehat{G}(\gamma)$ is such that

$$\widehat{\Theta}_{(\xi,\gamma)}(\widehat{g}) = \widetilde{g}_{\xi} \quad \text{for every} \quad \xi \in \overline{\xi}(\gamma).$$

With the help of $a_{\gamma} \colon \widetilde{F}(\gamma) \longrightarrow \widehat{G}(\gamma)$ we defined, in the same item 4.9 (see page 147), the binary relation $\approx \subseteq \widetilde{F}(\gamma) \times \widetilde{F}(\gamma)$ as follows: for $\widetilde{g}(\overline{\xi}(\gamma)) = \{\widetilde{g}_{\xi}\}_{\xi \in \overline{\xi}(\gamma)}$ and $\widetilde{h}(\overline{\eta}(\gamma)) = \{\widetilde{h}_{\eta}\}_{\eta \in \overline{\eta}(\gamma)}$ in $\widetilde{F}(\gamma)$,

$$\widetilde{g}\left(\overline{\xi}(\gamma)\right) \approx \widetilde{h}\left(\overline{\eta}(\gamma)\right) \quad \text{if and only if} \quad a_{\gamma}\left(\widetilde{g}\left(\overline{\xi}(\gamma)\right)\right) = a_{\gamma}\left(\widetilde{h}\left(\overline{\eta}(\gamma)\right)\right). \tag{4.11-1}$$

It results immediately from the definition above, as already highlighted in 4.9, that $\approx$ is an equivalence relation in $\widetilde{F}(\gamma)$. On the other hand, Lemma 4.10 tells us that:

$$\widetilde{g}\Big(\overline{\xi}(\gamma)\Big) \approx \widetilde{h}\Big(\overline{\eta}(\gamma)\Big) \quad \text{if and only if} \quad \widetilde{\Theta}_{(\xi\cap\eta,\xi)}\left(\widetilde{g}_\xi\right) = \widetilde{\Theta}_{(\xi\cap\eta,\eta)}(\widetilde{h}_\eta) \tag{4.11-2}$$

for every $\xi \in \overline{\xi}(\gamma)$ and $\eta \in \overline{\eta}(\gamma)$ such that $\xi \cap \eta \neq \varnothing$.

Let us now observe that the latter form, (4.11-2), of expressing the relation $\approx$, provided by Lemma 4.10, only involves elements of the $S$-space $\widetilde{\mathscr{G}}(I)$, in particular it does not make any reference to anything relative to the $S$-space $\widehat{\mathscr{G}}(I)$ and hence, it can be taken as definition of a binary relation in $\widetilde{F}(\gamma)$, independent of any hypothesis relative to $\widehat{\mathscr{G}}(I)$. Clearly, on the hypothesis of existence of $\widehat{\mathscr{G}}(I)$, it results from Lemma 4.10 the equivalence of both characterizations ((4.11-1) and (4.11-2) above) of the relation $\approx$ and the immediate conclusion, obtained from (4.11-1), that $\approx$ is an equivalence relation.

So, in this context, the following question applies: Defining $\approx$ through (4.11-2), is it an equivalence relation in $\widetilde{F}(\gamma)$?

Case affirmative, the partition of $\widetilde{F}(\gamma)$ determined by the relation $\approx$, that is, the set

$$\left[\widetilde{F}(\gamma)\right] = \left\{\left[\widetilde{g}\Big(\overline{\xi}(\gamma)\Big)\right] : \widetilde{g}\Big(\overline{\xi}(\gamma)\Big) \in \widetilde{F}(\gamma)\right\}$$

where

$$\left[\widetilde{g}\Big(\overline{\xi}(\gamma)\Big)\right] = \left\{\widetilde{h}\Big(\overline{\eta}(\gamma)\Big) \in \widetilde{F}(\gamma) : \widetilde{h}\Big(\overline{\eta}(\gamma)\Big) \approx \widetilde{g}\Big(\overline{\xi}(\gamma)\Big)\right\},$$

would be defined without any involvement of $\widehat{\mathscr{G}}(I)$: $[\widetilde{F}(\gamma)]$ built with what we really have, that is, the members of $\widetilde{\mathscr{G}}(I)$.

However, yet assuming an affirmative answer to the question above, in which case we can construct $[\widetilde{F}(\gamma)]$ without any need for $\widehat{\mathscr{G}}(I)$, the considerations of item 4.9 (see pages 146 and 147) show us that, if there exists $\widehat{\mathscr{G}}(I)$, then, the function

$$A_\gamma : \left[\widetilde{F}(\gamma)\right] \longrightarrow \widehat{G}(\gamma)$$

defined by

$$A_\gamma\left(\left[\widetilde{g}\Big(\overline{\xi}(\gamma)\Big)\right]\right) := \widehat{g}$$

where $\widehat{g} \in \widehat{G}(\gamma)$ is such that

$$a_\gamma\left(\widetilde{g}\Big(\overline{\xi}(\gamma)\Big) = \Big\{\widetilde{g}_\xi\Big\}_{\xi\in\overline{\xi}(\gamma)}\right) = \widehat{g}$$

or, equivalently (with the definition of the function $a_\gamma$ in mind), such that

$$\widehat{\Theta}_{(\xi,\gamma)}\left(\widehat{g}\right) = \widetilde{g}_\xi \quad \text{for every} \quad \xi \in \overline{\xi}(\gamma),$$

is a bijection, which turns the set $[\widetilde{F}(\gamma)]$ into a real candidate to the role of $\overline{G}(\gamma)$.

The observations above reveal the importance of the proposition below.

## 4.12 Proposition

*Let $\widetilde{\mathscr{G}}(\gamma)$ be a S-subspace of $\widetilde{\mathscr{G}}(I)$ arbitrarily fixed and $\widetilde{F}(\gamma)$ the class of coherent families in $\widetilde{\mathscr{G}}(\gamma)$. Let also $\approx$ be the binary relation in $\widetilde{F}(\gamma)$ defined as follows:*

*For $\widetilde{g}(\overline{\xi}(\gamma)) = \{\widetilde{g}_\xi\}_{\xi\in\overline{\xi}(\gamma)}$ and $\widetilde{h}(\overline{\eta}(\gamma)) = \{\widetilde{h}_\eta\}_{\eta\in\overline{\eta}(\gamma)}$ in $\widetilde{F}(\gamma)$, we define:*

$$\widetilde{g}\Big(\overline{\xi}(\gamma)\Big) \approx \widetilde{h}\Big(\overline{\eta}(\gamma)\Big) \quad \textit{if and only if} \quad \widetilde{\Theta}_{(\xi\cap\eta,\xi)}(\widetilde{g}_\xi) = \widetilde{\Theta}_{(\xi\cap\eta,\eta)}(\widetilde{h}_\eta)$$

*for every $\xi \in \overline{\xi}(\gamma)$ and $\eta \in \overline{\eta}(\gamma)$ such that $\xi \cap \eta \neq \varnothing$.*

*The relation $\approx$ is an equivalence relation in $\widetilde{F}(\gamma)$. Furthermore, on the hypothesis of existence of $\widehat{\mathscr{G}}(I)$ and being $[\widetilde{F}(\gamma)]$ the partition of $\widetilde{F}(\gamma)$ determined by $\approx$, the function $A_\gamma$ defined as follows is a bijection.*

$$A_\gamma : \quad \Big[\widetilde{F}(\gamma)\Big] \longrightarrow \widehat{G}(\gamma)$$

$$\Big[\widetilde{g}\Big(\overline{\xi}(\gamma)\Big)\Big] \longmapsto A_\gamma\Big(\Big[\widetilde{g}\Big(\overline{\xi}(\gamma)\Big)\Big]\Big) := \widehat{g}$$

*where $\widehat{g} \in \widehat{G}(\gamma)$ is such that*

$$\widehat{\Theta}_{(\xi,\gamma)}(\widehat{g}) = \widetilde{g}_\xi$$

*for every $\xi \in \overline{\xi}(\gamma)$ (the $\widetilde{g}_\xi$ are the members of the family $\widetilde{g}(\overline{\xi}(\gamma)) = \{\widetilde{g}_\xi\}_{\xi\in\overline{\xi}(\gamma)}$).*

*Proof.* There is no doubt that $\approx$, as above defined, is reflexive in $\widetilde{F}(\gamma)$, i.e.,

$$\widetilde{g}\Big(\overline{\xi}(\gamma)\Big) \approx \widetilde{g}\Big(\overline{\xi}(\gamma)\Big) \quad \text{for every} \quad \widetilde{g}\Big(\overline{\xi}(\gamma)\Big) \in \widetilde{F}(\gamma),$$

and symmetrical, that is,

$$\text{if} \quad \widetilde{g}\Big(\overline{\xi}(\gamma)\Big) \approx \widetilde{h}\Big(\overline{\eta}(\gamma)\Big), \quad \text{then} \quad \widetilde{h}\Big(\overline{\eta}(\gamma)\Big) \approx \widetilde{g}\Big(\overline{\xi}(\gamma)\Big).$$

It is then enough to prove that $\approx$ is transitive. Let then

$$\widetilde{f}\Big(\overline{\nu}(\gamma)\Big) = \Big\{\widetilde{f}_\nu\Big\}_{\nu\in\overline{\nu}(\gamma)}, \quad \widetilde{g}\Big(\overline{\xi}(\gamma)\Big) = \Big\{\widetilde{g}_\xi\Big\}_{\xi\in\overline{\xi}(\gamma)} \quad \text{and} \quad \widetilde{h}\Big(\overline{\eta}(\gamma)\Big) = \Big\{\widetilde{h}_\eta\Big\}_{\eta\in\overline{\eta}(\gamma)},$$

be members of $\widetilde{F}(\gamma)$ arbitrarily chosen, such that

$$\widetilde{f}\Big(\overline{\nu}(\gamma)\Big) \approx \widetilde{g}\Big(\overline{\xi}(\gamma)\Big) \quad \text{and} \quad \widetilde{g}\Big(\overline{\xi}(\gamma)\Big) \approx \widetilde{h}\Big(\overline{\eta}(\gamma)\Big),$$

that is,

$$\widetilde{\Theta}_{(\nu\cap\xi,\nu)}(\widetilde{f}_\nu) = \widetilde{\Theta}_{(\nu\cap\xi,\xi)}(\widetilde{g}_\xi) \tag{4.12-1}$$

for every $\nu \in \overline{\nu}(\gamma)$ and $\xi \in \overline{\xi}(\gamma)$ such that $\nu \cap \xi \neq \varnothing$, and

$$\widetilde{\Theta}_{(\xi\cap\eta,\xi)}(\widetilde{g}_\xi) = \widetilde{\Theta}_{(\xi\cap\eta,\eta)}(\widetilde{h}_\eta) \tag{4.12-2}$$

for every $\xi \in \overline{\xi}(\gamma)$ and $\eta \in \overline{\eta}(\gamma)$ such that $\xi \cap \eta \neq \varnothing$.

We want to prove, given (4.12-1) and (4.12-2) above, that

$$\widetilde{f}\Big(\overline{\nu}(\gamma)\Big) \approx \widetilde{h}\Big(\overline{\eta}(\gamma)\Big),$$

i.e., that

$$\widetilde{\Theta}_{(\nu\cap\eta,\nu)}(\widetilde{f}_\nu) = \widetilde{\Theta}_{(\nu\cap\eta,\eta)}(\widetilde{h}_\eta)$$

for every $\nu \in \overline{\nu}(\gamma)$ and $\eta \in \overline{\eta}(\gamma)$ such that $\nu \cap \eta \neq \varnothing$.

From (4.12-1) and (4.12-2), respectively, it results that

$$\widetilde{\Theta}_{(\xi\cap\nu\cap\eta,\nu\cap\xi)}\Big(\widetilde{\Theta}_{(\nu\cap\xi,\nu)}(\widetilde{f}_\nu)\Big) = \widetilde{\Theta}_{(\xi\cap\nu\cap\eta,\nu\cap\xi)}\Big(\widetilde{\Theta}_{(\nu\cap\xi,\xi)}(\widetilde{g}_\xi)\Big)$$

and

$$\widetilde{\Theta}_{(\xi\cap\nu\cap\eta,\xi\cap\eta)}\Big(\widetilde{\Theta}_{(\xi\cap\eta,\xi)}(\widetilde{g}_\xi)\Big) = \widetilde{\Theta}_{(\xi\cap\nu\cap\eta,\xi\cap\eta)}\Big(\widetilde{\Theta}_{(\xi\cap\eta,\eta)}(\widetilde{h}_\eta)\Big),$$

that is,

$$\widetilde{\Theta}_{(\xi\cap\nu\cap\eta,\nu)}(\widetilde{f}_\nu) = \widetilde{\Theta}_{(\xi\cap\nu\cap\eta,\xi)}(\widetilde{g}_\xi)$$

and

$$\widetilde{\Theta}_{(\xi\cap\nu\cap\eta,\xi)}(\widetilde{g}_\xi) = \widetilde{\Theta}_{(\xi\cap\nu\cap\eta,\eta)}(\widetilde{h}_\eta),$$

from where one concludes that

$$\widetilde{\Theta}_{(\xi\cap\nu\cap\eta,\nu)}(\widetilde{f}_\nu) = \widetilde{\Theta}_{(\xi\cap\nu\cap\eta,\eta)}(\widetilde{h}_\eta)$$

for every $\xi \in \overline{\xi}(\gamma)$, $\nu \in \overline{\nu}(\gamma)$ and $\eta \in \overline{\eta}(\gamma)$ such that $\xi \cap \nu \cap \eta \neq \varnothing$.

This last equation can still be written in the following form:

$$\widetilde{\Theta}_{(\xi\cap\nu\cap\eta,\nu\cap\eta)}\Big(\widetilde{\Theta}_{(\nu\cap\eta,\nu)}(\widetilde{f}_\nu)\Big) = \widetilde{\Theta}_{(\xi\cap\nu\cap\eta,\nu\cap\eta)}\Big(\widetilde{\Theta}_{(\nu\cap\eta,\eta)}(\widetilde{h}_\eta)\Big)$$

or yet

$$\widetilde{\Theta}_{(\xi\cap\nu\cap\eta,\nu\cap\eta)}\Big(\widetilde{\Theta}_{(\nu\cap\eta,\nu)}(\widetilde{f}_\nu) - \widetilde{\Theta}_{(\nu\cap\eta,\eta)}(\widetilde{h}_\eta)\Big) = 0$$

for every $\xi \in \overline{\xi}(\gamma)$, $\nu \in \overline{\nu}(\gamma)$ and $\eta \in \overline{\eta}(\gamma)$ such that $\xi \cap \nu \cap \eta \neq \varnothing$.

If now we define $\widetilde{e}_{\nu\cap\eta} \in \widetilde{G}(\nu \cap \eta)$ by

$$\widetilde{e}_{\nu\cap\eta} := \widetilde{\Theta}_{(\nu\cap\eta,\nu)}(\widetilde{f}_\nu) - \widetilde{\Theta}_{(\nu\cap\eta,\eta)}(\widetilde{h}_\eta) \tag{4.12-3}$$

we get that

$$\widetilde{\Theta}_{(\xi\cap\nu\cap\eta,\nu\cap\eta)}(\widetilde{e}_{\nu\cap\eta}) = 0$$

for every $\xi \in \overline{\xi}(\gamma)$, $\nu \in \overline{\nu}(\gamma)$ and $\eta \in \overline{\eta}(\gamma)$ such that $\xi \cap \nu \cap \eta \neq \varnothing$.

Remembering now that the $S$-group $\widetilde{G}(\nu\cap\eta)$ is a strict and closed extension of the $S$-group $G(\nu \cap \eta)$, then, there exists $e_{\nu\cap\eta} \in G(\nu \cap \eta)$ and $\Phi_{\nu\cap\eta} \in H(\nu \cap \eta)$ such that

$$\widetilde{e}_{\nu\cap\eta} = \widetilde{\Phi}_{\nu\cap\eta}(e_{\nu\cap\eta}) \tag{4.12-4}$$

and, thus, we get

$$\widetilde{\Theta}_{(\xi\cap\nu\cap\eta,\nu\cap\eta)}\Big(\widetilde{\Phi}_{\nu\cap\eta}(e_{\nu\cap\eta})\Big) = 0,$$

or yet,

$$\widetilde{i}_{(\xi\cap\nu\cap\eta,\nu\cap\eta)}(\widetilde{\Phi}_{\nu\cap\eta})\Big(\widetilde{\Theta}_{(\xi\cap\nu\cap\eta,\nu\cap\eta)}(e_{\nu\cap\eta})\Big) = 0.$$

Therefore,

$$\widetilde{\Theta}_{(\xi\cap\nu\cap\eta,\nu\cap\eta)}(e_{\nu\cap\eta}) \in N\Big(\widetilde{i}_{(\xi\cap\nu\cap\eta,\nu\cap\eta)}(\widetilde{\Phi}_{\nu\cap\eta})\Big)$$

for every $\xi \in \overline{\xi}(\gamma)$, $\nu \in \overline{\nu}(\gamma)$ and $\eta \in \overline{\eta}(\gamma)$ such that $\xi \cap \nu \cap \eta \neq \varnothing$.

On the other hand, we have (see page 140)

$$\widetilde{i}_{(\chi,\nu\cap\eta)}(\widetilde{\Phi}_{\nu\cap\eta}) = \widetilde{\chi}\Big(i_{(\chi,\nu\cap\eta)}\left(\Phi_{\nu\cap\eta}\right)\Big)$$

being $\chi := \xi \cap \nu \cap \eta$. Hence we see that

$$\widetilde{\Theta}_{(\chi,\nu\cap\eta)}\left(e_{\nu\cap\eta}\right) \in N\left(\widetilde{\chi}\Big(i_{(\chi,\nu\cap\eta)}\left(\Phi_{\nu\cap\eta}\right)\Big)\right)$$

for every $\xi \in \overline{\xi}(\gamma)$, $\nu \in \overline{\nu}(\gamma)$ and $\eta \in \overline{\eta}(\gamma)$ such that $\chi = \xi \cap \nu \cap \eta \neq \varnothing$.

But, by Proposition 1.16(b),

$$N\left(\widetilde{\chi}\Big(i_{(\chi,\nu\cap\eta)}\left(\Phi_{\nu\cap\eta}\right)\Big)\right) = N\Big(i_{(\chi,\nu\cap\eta)}\left(\Phi_{\nu\cap\eta}\right)\Big)$$

and hence

$$\widetilde{\Theta}_{(\chi,\nu\cap\eta)}\left(e_{\nu\cap\eta}\right) = \Theta_{(\chi,\nu\cap\eta)}\left(e_{\nu\cap\eta}\right) \in N\Big(i_{(\chi,\nu\cap\eta)}\left(\Phi_{\nu\cap\eta}\right)\Big),$$

that is,

$$i_{(\chi,\nu\cap\eta)}\left(\Phi_{\nu\cap\eta}\right)\Big(\Theta_{(\chi,\nu\cap\eta)}(e_{\nu\cap\eta})\Big) = 0 \tag{4.12-5}$$

for every $\xi \in \overline{\xi}(\gamma)$, $\nu \in \overline{\nu}(\gamma)$ and $\eta \in \overline{\eta}(\gamma)$ such that $\chi = \xi \cap \nu \cap \eta \neq \varnothing$.

Now, for $\nu_0 \in \overline{\nu}(\gamma)$ and $\eta_0 \in \overline{\eta}(\gamma)$ arbitrarily fixed an such that $\nu_0 \cap \eta_0 \neq \varnothing$, we define:

$$\overline{\alpha}(\nu_0 \cap \eta_0) := \Big\{\alpha = \xi \cap \nu_0 \cap \eta_0 \quad : \quad \xi \in \overline{\xi}(\gamma) \quad \text{and} \quad \xi \cap (\nu_0 \cap \eta_0) \neq \varnothing\Big\}.$$

Clearly $\overline{\alpha}(\nu_0 \cap \eta_0) \subseteq \Gamma(\nu_0 \cap \eta_0)$ is a cover of $\nu_0 \cap \eta_0$. Furthermore, it results from (4.12-5) that

$$i_{(\alpha,\nu_0\cap\eta_0)}(\Phi_{\nu_0\cap\eta_0})\Big(\Theta_{(\alpha,\nu_0\cap\eta_0)}(e_{\nu_0\cap\eta_0})\Big) = 0$$

for every $\alpha \in \overline{\alpha}(\nu_0 \cap \eta_0)$.

If now we resort to Definition 3.13 of the concept of restriction, more specifically its item (d), we will see that the last expression allows us to conclude that $e_{\nu_0\cap\eta_0}$ is an element in the domain of the homomorphism $\Phi_{\nu_0\cap\eta_0}$, that is,

$$e_{\nu_0\cap\eta_0} \in G(\nu_0 \cap \eta_0)_{\Phi_{\nu_0\cap\eta_0}},$$

and hence, taking into account item (c) of the referred definition, we conclude that

$$\Theta_{(\alpha,\nu_0\cap\eta_0)}\Big(\Phi_{\nu_0\cap\eta_0}(e_{\nu_0\cap\eta_0})\Big) = 0 \tag{4.12-6}$$

for every $\alpha \in \overline{\alpha}(\nu_0 \cap \eta_0)$.

To be more precise, the second member of this last equation, 0, is the additive neutral of the group $G(\alpha)$ (since $\Theta_{(\alpha,\nu_0\cap\eta_0)}: G(\nu_0 \cap \eta_0) \longrightarrow G(\alpha)$ is a homomorphism from $G(\nu_0 \cap \eta_0)$ into $G(\alpha)$), henceforth denoted by $0_\alpha$. Hence, and since

$$\Theta_{(\alpha,\nu_0\cap\eta_0)}(0_{\nu_0\cap\eta_0}) = 0_\alpha$$

with $0_{\nu_0\cap\eta_0} \in G(\nu_0 \cap \eta_0)$ the additive neutral of $G(\nu_0 \cap \eta_0)$, our last conclusion, (4.12-6), can be expressed in the following form:

$$\Theta_{(\alpha,\nu_0\cap\eta_0)}\Big(\Phi_{\nu_0\cap\eta_0}(e_{\nu_0\cap\eta_0})\Big) = 0_\alpha = \Theta_{(\alpha,\nu_0\cap\eta_0)}(0_{\nu_0\cap\eta_0})$$

for every $\alpha \in \overline{\alpha}(\nu_0 \cap \eta_0)$.

Taking now into account that $\mathscr{G}(\nu_0 \cap \eta_0)$ is coherent (since $\mathscr{G}(I)$ is), the last equation together with Lemma 3.32 allow us to conclude that

$$\Phi_{\nu_0\cap\eta_0}(e_{\nu_0\cap\eta_0}) = 0_{\nu_0\cap\eta_0}$$

or, once $\widetilde{\Phi}_{\nu_0\cap\eta_0}$ is a prolongation of $\Phi_{\nu_0\cap\eta_0}$,

$$\widetilde{\Phi}_{\nu_0\cap\eta_0}(e_{\nu_0\cap\eta_0}) = 0_{\nu_0\cap\eta_0}$$

from where one obtains, through expressions (4.12-3) and (4.12-4), that

$$\widetilde{\Phi}_{\nu_0\cap\eta_0}(e_{\nu_0\cap\eta_0}) = \widetilde{e}_{\nu_0\cap\eta_0} = \widetilde{\Theta}_{(\nu_0\cap\eta_0,\nu_0)}(\widetilde{f}_{\nu_0}) - \widetilde{\Theta}_{(\nu_0\cap\eta_0,\eta_0)}(\widetilde{h}_{\eta_0}) = 0_{\nu_0\cap\eta_0}$$

and hence

$$\widetilde{\Theta}_{(\nu_0\cap\eta_0,\nu_0)}(\widetilde{f}_{\nu_0}) = \widetilde{\Theta}_{(\nu_0\cap\eta_0,\eta_0)}(\widetilde{h}_{\eta_0}).$$

Since $\nu_0 \in \overline{\nu}(\gamma)$ and $\eta_0 \in \overline{\eta}(\gamma)$ such that $\nu_0 \cap \eta_0 \neq \varnothing$ were arbitrarily fixed, we have

$$\widetilde{\Theta}_{(\nu\cap\eta,\nu)}(\widetilde{f}_\nu) = \widetilde{\Theta}_{(\nu\cap\eta,\eta)}(\widetilde{h}_\eta)$$

for every $\nu \in \overline{\nu}(\gamma)$ and $\eta \in \overline{\eta}(\gamma)$ such that $\nu \cap \eta \neq \varnothing$, which means, according to the definition of the relation $\approx$, that

$$\widetilde{f}\Big(\overline{\nu}(\gamma)\Big) \approx \widetilde{h}\Big(\overline{\eta}(\gamma)\Big)$$

hence concluding the proof of the transitivity of the relation $\approx$. ∎

## 4.13 The Set $\overline{G}(\gamma)$

Regarding the task of building, starting with the elements of the $S$-spaces $\mathscr{G}(I)$ and $\widetilde{\mathscr{G}}(I)$, an extension of $\widetilde{\mathscr{G}}(I)$,

$$\overline{\mathscr{G}}(I) = \left(\overline{\mathbb{G}}\Big(\Gamma(I)\Big) = \Big\{\overline{\mathbb{G}}(\gamma) = \Big(\overline{G}(\gamma), \overline{H}(\gamma)\Big)\Big\}_{\gamma\in\Gamma(I)}, \overline{\imath}\Big(\Gamma^2(I)\Big), \overline{\Theta}\Big(\Delta(I)\Big)\right),$$

locally closed and coherent, isomorphic to $\widehat{\mathscr{G}}(I)$ on the hypothesis of existence of the latter, the Proposition 4.12 allows us to affirm that, for each $\gamma \in \Gamma(I)$, the set $\overline{G}(\gamma)$ defined ahead, not only is well-defined existing $\widehat{\mathscr{G}}(I)$ or not, but also, existing $\widehat{\mathscr{G}}(I)$, it has the same cardinality of $\widehat{G}(\gamma)$ and the function $A_\gamma \colon \overline{G}(\gamma) \longrightarrow \widehat{G}(\gamma)$ defined in the referred proposition is a bijection.

- **Definition of $\overline{G}(\gamma)$:** for $\gamma \in \Gamma(I)$ we define

  $$\overline{G}(\gamma) := \Big[\widetilde{F}(\gamma)\Big],$$

  where $\widetilde{F}(\gamma)$ is the set of coherent families in $\widetilde{\mathscr{G}}(\gamma)$, and $[\widetilde{F}(\gamma)]$ the partition of $\widetilde{F}(\gamma)$ determined by the equivalence relation $\approx$ defined in Proposition 4.12. More explicitly,

  $$\overline{G}(\gamma) = \Big\{\Big[\widetilde{g}\Big(\overline{\xi}(\gamma)\Big)\Big] : \widetilde{g}\Big(\overline{\xi}(\gamma)\Big) \in \widetilde{F}(\gamma)\Big\},$$

  where

  $$\Big[\widetilde{g}\Big(\overline{\xi}(\gamma)\Big)\Big] = \left\{\widetilde{h}\Big(\overline{\eta}(\gamma)\Big) \in \widetilde{F}(\gamma) \colon \widetilde{h}\Big(\overline{\eta}(\gamma)\Big) \approx \widetilde{g}\Big(\overline{\xi}(\gamma)\Big)\right\},$$

  with

  $$\widetilde{h}\Big(\overline{\eta}(\gamma)\Big) = \Big\{\widetilde{h}_\eta\Big\}_{\eta\in\overline{\eta}(\gamma)} \quad \approx \quad \widetilde{g}\Big(\overline{\xi}(\gamma)\Big) = \Big\{\widetilde{g}_\xi\Big\}_{\xi\in\overline{\xi}(\gamma)}$$

  if and only if

  $$\widetilde{\Theta}_{(\xi\cap\eta,\eta)}(\widetilde{h}_\eta) = \widetilde{\Theta}_{(\xi\cap\eta,\xi)}(\widetilde{g}_\xi)$$

  for every $\xi \in \overline{\xi}(\gamma)$ and $\eta \in \overline{\eta}(\gamma)$ such that $\xi \cap \eta \neq \varnothing$.

# The Addition in $\overline{G}(\gamma)$ and the Group $\overline{G}(\gamma)$

## 4.14 Motivation

Our next task, regarding the construction of $\overline{\mathscr{G}}(I)$, consists of defining a binary operation in the set $\overline{G}(\gamma)$, an addition $(+)$, so that $\overline{G}(\gamma)$ equipped with this operation becomes an abelian group that, under the hypothesis of existence of $\widehat{\mathscr{G}}(I)$, is isomorphic to the abelian group $\widehat{G}(\gamma)$. In fact, since we want $\overline{\mathscr{G}}(I)$ isomorphic to $\widehat{\mathscr{G}}(I)$, the definition of isomorphism for $S$-spaces (Definition 3.28) requires the family of $S$-groups of $\overline{\mathscr{G}}(I)$, namely,

$$\overline{\mathbb{G}}\Big(\Gamma(I)\Big) = \Big\{\overline{\mathbb{G}}(\gamma) = \Big(\overline{G}(\gamma), \overline{H}(\gamma)\Big)\Big\}_{\gamma\in\Gamma(I)},$$

to be isomorphic to the family of $S$-groups of $\widehat{\mathscr{G}}(I)$,

$$\widehat{\mathbb{G}}\Big(\Gamma(I)\Big) = \Big\{\widehat{\mathbb{G}}(\gamma) = \Big(\widehat{G}(\gamma), \widehat{H}(\gamma)\Big)\Big\}_{\gamma\in\Gamma(I)},$$

and this means, among other things, that the group $\overline{G}(\gamma)$ is isomorphic to the group $\widehat{G}(\gamma)$, for every $\gamma \in \Gamma(I)$.

Now, we already know that, on the existence of $\widehat{\mathscr{G}}(I)$, the function $A_\gamma : \overline{G}(\gamma) \longrightarrow \widehat{G}(\gamma)$ defined in Proposition 4.12 is a bijection and, hence, the following binary operation, $+$, in $\overline{G}(\gamma)$ is well-defined:

$$\begin{aligned} + : \overline{G}(\gamma) \times \overline{G}(\gamma) &\longrightarrow \overline{G}(\gamma) \\ \big(\overline{g}, \overline{h}\big) &\longmapsto \overline{g} + \overline{h}, \end{aligned}$$

where

$$\overline{g} + \overline{h} \coloneqq A_\gamma^{-1}\Big(A_\gamma\left(\overline{g}\right) + A_\gamma\left(\overline{h}\right)\Big).$$

Thus, for

$$\overline{g} = \Big[\widetilde{g}\Big(\overline{\xi}(\gamma)\Big) = \Big\{\widetilde{g}_\xi\Big\}_{\xi\in\overline{\xi}(\gamma)}\Big] \quad \text{and} \quad \overline{h} = \Big[\widetilde{h}\Big(\overline{\eta}(\gamma)\Big) = \Big\{\widetilde{h}_\eta\Big\}_{\eta\in\overline{\eta}(\gamma)}\Big]$$

arbitrarily chosen in $\overline{G}(\gamma) = [\widetilde{F}(\gamma)]$, and taking into account the definition of $A_\gamma$ given in Proposition 4.12, we have:

**(a)** $A_\gamma(\overline{g}) = \widehat{g}$ and $A_\gamma(\overline{h}) = \widehat{h}$ with $\widehat{g}, \widehat{h} \in \widehat{G}(\gamma)$ such that

$$\widehat{\Theta}_{(\xi,\gamma)}(\widehat{g}) = \widetilde{g}_\xi \quad \text{for every} \quad \xi \in \overline{\xi}(\gamma)$$

and

$$\widehat{\Theta}_{(\eta,\gamma)}(\widehat{h}) = \widetilde{h}_\eta \quad \text{for every} \quad \eta \in \overline{\eta}(\gamma);$$

**(b)** From (a) it results that

$$\widehat{\Theta}_{(\xi\cap\eta,\xi)}\Big(\widehat{\Theta}_{(\xi,\gamma)}(\widehat{g})\Big) = \widehat{\Theta}_{(\xi\cap\eta,\xi)}(\widetilde{g}_\xi)$$

and

$$\widehat{\Theta}_{(\xi\cap\eta,\eta)}\Big(\widehat{\Theta}_{(\eta,\gamma)}(\widehat{h})\Big) = \widehat{\Theta}_{(\xi\cap\eta,\eta)}(\widetilde{h}_\eta),$$

that is,

$$\widehat{\Theta}_{(\xi\cap\eta,\gamma)}(\widehat{g}) = \widehat{\Theta}_{(\xi\cap\eta,\xi)}(\widetilde{g}_\xi) \quad \text{and} \quad \widehat{\Theta}_{(\xi\cap\eta,\gamma)}(\widehat{h}) = \widehat{\Theta}_{(\xi\cap\eta,\eta)}(\widetilde{h}_\eta)$$

for every $\xi \in \overline{\xi}(\gamma)$ and $\eta \in \overline{\eta}(\gamma)$ such that $\xi \cap \eta \neq \varnothing$.

Therefore, since $\widehat{\Theta}_{(\xi\cap\eta,\xi)}(\widetilde{g}_\xi) = \widetilde{\Theta}_{(\xi\cap\eta,\xi)}(\widetilde{g}_\xi)$, $\widehat{\Theta}_{(\xi\cap\eta,\eta)}(\widetilde{h}_\eta) = \widetilde{\Theta}_{(\xi\cap\eta,\eta)}(\widetilde{h}_\eta)$ and $\widehat{\Theta}_{(\xi\cap\eta,\gamma)} : \widehat{G}(\gamma) \longrightarrow \widehat{G}(\xi \cap \eta)$ is a homomorphism, we have:

$$\widehat{\Theta}_{(\xi\cap\eta,\gamma)}(\widehat{g} + \widehat{h}) = \widetilde{\Theta}_{(\xi\cap\eta,\xi)}(\widetilde{g}_\xi) + \widetilde{\Theta}_{(\xi\cap\eta,\eta)}(\widetilde{h}_\eta)$$

for every $\xi \in \overline{\xi}(\gamma)$ and $\eta \in \overline{\eta}(\gamma)$ such that $\xi \cap \eta \neq \varnothing$;

**(c)** From (a) and (b) above we obtain:

$$A_\gamma(\overline{g}) + A_\gamma(\overline{h}) = \widehat{g} + \widehat{h}$$

with $\widehat{g} + \widehat{h} \in \widehat{G}(\gamma)$ such that

$$\widehat{\Theta}_{(\xi\cap\eta,\gamma)}(\widehat{g} + \widehat{h}) = \widetilde{\Theta}_{(\xi\cap\eta,\xi)}(\widetilde{g}_\xi) + \widetilde{\Theta}_{(\xi\cap\eta,\eta)}(\widetilde{h}_\eta)$$

for every $\xi \in \overline{\xi}(\gamma)$ and $\eta \in \overline{\eta}(\gamma)$ such that $\xi \cap \eta \neq \varnothing$;

**(d)** Let now $\overline{\chi}(\gamma) \subseteq \Gamma(\gamma)$ be the cover of $\gamma$ given by

$$\overline{\chi}(\gamma) \coloneqq \Big\{\chi \coloneqq \xi \cap \eta \quad : \quad \xi \in \overline{\xi}(\gamma), \quad \eta \in \overline{\eta}(\gamma), \quad \text{and} \quad \xi \cap \eta \neq \varnothing\Big\}$$

and let also

$$\widetilde{t}\Big(\overline{\chi}(\gamma)\Big) \coloneqq \Big\{\widetilde{t}_\chi \coloneqq \widetilde{\Theta}_{(\xi\cap\eta,\xi)}(\widetilde{g}_\xi) + \widetilde{\Theta}_{(\xi\cap\eta,\eta)}(\widetilde{h}_\eta)\Big\}_{\chi=\xi\cap\eta\in\overline{\chi}(\gamma)}$$

be the family that, as we will see, is coherent in $\widetilde{\mathscr{G}}(\gamma)$, that is, $\widetilde{t}(\overline{\chi}(\gamma)) \in \widetilde{F}(\gamma)$. From (c), we get:

$$\widehat{\Theta}_{(\chi,\gamma)}(\widehat{g} + \widehat{h}) = \widetilde{t}_\chi \quad \text{for every} \quad \chi \in \overline{\chi}(\gamma)$$

which, in turn, by the definition of the function $A_\gamma : \overline{G}(\gamma) \longrightarrow \widehat{G}(\gamma)$, means that

$$\overline{g} + \overline{h} = A_\gamma^{-1}\Big(A_\gamma(\overline{g}) + A_\gamma(\overline{\eta})\Big) = \Big[\widetilde{t}\Big(\overline{\chi}(\gamma)\Big)\Big];$$

**(e)** In short, the operation $+ \colon \overline{G}(\gamma) \times \overline{G}(\gamma) \longrightarrow \overline{G}(\gamma)$ previously defined, can be formulated in the following form:

$$+ \colon \overline{G}(\gamma) \times \overline{G}(\gamma) \longrightarrow \overline{G}(\gamma)$$

is the function that for each pair $(\overline{g}, \overline{h}) \in \overline{G}(\gamma) \times \overline{G}(\gamma)$,

$$\overline{g} = \Big[\widetilde{g}\Big(\overline{\xi}(\gamma)\Big) = \Big\{\widetilde{g}_\xi\Big\}_{\xi\in\overline{\xi}(\gamma)}\Big] \quad \text{and} \quad \overline{h} = \Big[\widetilde{h}\Big(\overline{\eta}(\gamma)\Big) = \Big\{\widetilde{h}_\eta\Big\}_{\eta\in\overline{\eta}(\gamma)}\Big],$$

associates the element $\overline{g} + \overline{h} \in \overline{G}(\gamma)$ defined by

$$\overline{g} + \overline{h} \coloneqq \Big[\widetilde{t}\Big(\overline{\chi}(\gamma)\Big) = \Big\{\widetilde{t}_\chi\Big\}_{\chi\in\overline{\chi}(\gamma)}\Big],$$

with

$$\overline{\chi}(\gamma) \coloneqq \Big\{ \chi \coloneqq \xi \cap \eta \quad : \quad \xi \in \overline{\xi}(\gamma), \quad \eta \in \overline{\eta}(\gamma), \quad \text{and} \quad \xi \cap \eta \neq \varnothing \Big\}$$

and

$$\widetilde{t}_\chi \coloneqq \widetilde{\Theta}_{(\chi,\xi)}(\widetilde{g}_\xi) + \widetilde{\Theta}_{(\chi,\eta)}(\widetilde{h}_\eta)$$

for every $\chi = \xi \cap \eta \in \overline{\chi}(\gamma)$.

Let us note now that the definition of the operation $+ \colon \overline{G}(\gamma) \times \overline{G}(\gamma) \longrightarrow \overline{G}(\gamma)$ as formulated in (e), only involves elements of $\widetilde{\mathscr{G}}(I)$ and, therefore, does not make any mention to $\widehat{\mathscr{G}}(I)$. However, the only assurance we have that the operation $+$ with the formulation given in (e) is well-defined, i.e., that in fact it is a function from $\overline{G}(\gamma) \times \overline{G}(\gamma)$ into $\overline{G}(\gamma)$, is given by its equivalence with the original definition (on page 157) which, in turn, resorts to the function $A_\gamma : \overline{G}(\gamma) \longrightarrow \widehat{G}(\gamma)$ and, therefore, appeals to $\widehat{G}(\gamma)$. So, the following question is appropriate: The addition in $\overline{G}(\gamma)$, $+$, as given in (e), that is, expressed through, only, the elements of $\widetilde{\mathscr{G}}(I)$, without any participation of $\widehat{\mathscr{G}}(I)$, is well-defined? Case affirmative, we will have obtained a candidate for the role of the addition in $\overline{G}(\gamma)$.

The observations above motivate the propositions 4.15 and 4.16 ahead.

## 4.15 Proposition

*Let $\overline{G}(\gamma)$ be the set defined in 4.13 and let also*

$$\overline{g} = \Big[\widetilde{g}\Big(\overline{\xi}(\gamma)\Big) = \Big\{\widetilde{g}_\xi\Big\}_{\xi\in\overline{\xi}(\gamma)}\Big] \quad \textit{and} \quad \overline{h} = \Big[\widetilde{h}\Big(\overline{\eta}(\gamma)\Big) = \Big\{\widetilde{h}_\eta\Big\}_{\eta\in\overline{\eta}(\gamma)}\Big]$$

*be arbitrarily chosen elements of $\overline{G}(\gamma)$. One has a well defined binary operation in $\overline{G}(\gamma)$, $+ \colon \overline{G}(\gamma) \times \overline{G}(\gamma) \longrightarrow \overline{G}(\gamma)$, defining $\overline{g} + \overline{h}$ as follows:*

$$\overline{g} + \overline{h} = \Big[\widetilde{g}\Big(\overline{\xi}(\gamma)\Big)\Big] + \Big[\widetilde{h}\Big(\overline{\eta}(\gamma)\Big)\Big] \coloneqq \Big[\widetilde{t}\Big(\overline{\chi}(\gamma)\Big)\Big],$$

*where $\overline{\chi}(\gamma) \subseteq \Gamma(\gamma)$ is the cover of $\gamma$ defined by*

$$\overline{\chi}(\gamma) \coloneqq \Big\{ \chi \coloneqq \xi \cap \eta \quad : \quad \xi \in \overline{\xi}(\gamma), \quad \eta \in \overline{\eta}(\gamma), \quad \textit{and} \quad \xi \cap \eta \neq \varnothing \Big\}$$

*and $\widetilde{t}(\overline{\chi}(\gamma)) \coloneqq \{\widetilde{t}_\chi\}_{\chi\in\overline{\chi}(\gamma)}$ is the family whose members are*

$$\widetilde{t}_{\chi=\xi\cap\eta} \coloneqq \widetilde{\Theta}_{(\xi\cap\eta,\xi)}(\widetilde{g}_\xi) + \widetilde{\Theta}_{(\xi\cap\eta,\eta)}(\widetilde{h}_\eta)$$

*for every $\chi = \xi \cap \eta \in \overline{\chi}(\gamma)$.*

*Furthermore, under the hypothesis of existence of $\widehat{\mathscr{G}}(I)$, the function $A_\gamma : \overline{G}(\gamma) \longrightarrow \widehat{G}(\gamma)$ as defined in Proposition 4.12, is a bijection and*

$$A_\gamma(\overline{g} + \overline{h}) = A_\gamma(\overline{g}) + A_\gamma(\overline{h}).$$

*Proof.* We must prove first that the family $\widetilde{t}(\overline{\chi}(\gamma))$ is coherent in $\widetilde{\mathscr{G}}(\gamma)$, i.e., that $\widetilde{t}(\overline{\chi}(\gamma)) \in \widetilde{F}(\gamma)$, since for $\overline{g} + \overline{h}$ as above defined to belong to $\overline{G}(\gamma)$ is necessary and sufficient that $\widetilde{t}(\overline{\chi}(\gamma)) \in \widetilde{F}(\gamma)$. Let then $\chi = \xi \cap \eta$ and $\chi' = \xi' \cap \eta'$ be arbitrarily chosen in $\overline{\chi}(\gamma)$ and such that $\chi \cap \chi' \neq \varnothing$. We want to prove that

$$\widetilde{\Theta}_{(\chi\cap\chi',\chi)}(\widetilde{t}_\chi) = \widetilde{\Theta}_{(\chi\cap\chi',\chi')}(\widetilde{t}_{\chi'}).$$

We have

$$\begin{aligned}\widetilde{\Theta}_{(\chi\cap\chi',\chi)}(\widetilde{t}_\chi) &= \widetilde{\Theta}_{(\chi\cap\chi',\chi)}\Big(\widetilde{\Theta}_{(\chi,\xi)}(\widetilde{g}_\xi) + \widetilde{\Theta}_{(\chi,\eta)}(\widetilde{h}_\eta)\Big) = \\ &= \widetilde{\Theta}_{(\chi\cap\chi',\xi)}(\widetilde{g}_\xi) + \widetilde{\Theta}_{(\chi\cap\chi',\eta)}(\widetilde{h}_\eta)\end{aligned}$$

and analogously

$$\widetilde{\Theta}_{(\chi\cap\chi',\chi')}(\widetilde{t}_{\chi'}) = \widetilde{\Theta}_{(\chi\cap\chi',\xi')}(\widetilde{g}_{\xi'}) + \widetilde{\Theta}_{(\chi\cap\chi',\eta')}(\widetilde{h}_{\eta'})$$

We also have

$$\widetilde{\Theta}_{(\chi\cap\chi',\xi)}(\widetilde{g}_\xi) = \widetilde{\Theta}_{(\chi\cap\chi',\xi\cap\xi')}\Big(\widetilde{\Theta}_{(\xi\cap\xi',\xi)}(\widetilde{g}_\xi)\Big)$$

and

$$\widetilde{\Theta}_{(\chi\cap\chi',\xi')}(\widetilde{g}_{\xi'}) = \widetilde{\Theta}_{(\chi\cap\chi',\xi\cap\xi')}\Big(\widetilde{\Theta}_{(\xi\cap\xi',\xi')}(\widetilde{g}_{\xi'})\Big).$$

But, $\widetilde{\Theta}_{(\xi\cap\xi',\xi)}(\widetilde{g}_\xi) = \widetilde{\Theta}_{(\xi\cap\xi',\xi')}(\widetilde{g}_{\xi'})$, since $\widetilde{g}(\overline{\xi}(\gamma)) \in \widetilde{F}(\gamma)$, and hence we obtain:

$$\widetilde{\Theta}_{(\chi\cap\chi',\xi)}(\widetilde{g}_\xi) = \widetilde{\Theta}_{(\chi\cap\chi',\xi')}(\widetilde{g}_{\xi'}).$$

Analogously

$$\widetilde{\Theta}_{(\chi\cap\chi',\eta)}(\widetilde{h}_\eta) = \widetilde{\Theta}_{(\chi\cap\chi',\eta')}(\widetilde{h}_{\eta'}).$$

Introducing the last two results into the expressions obtained above for $\widetilde{\Theta}_{(\chi\cap\chi',\chi)}(\widetilde{t}_\chi)$ and $\widetilde{\Theta}_{(\chi\cap\chi',\chi')}(\widetilde{t}_{\chi'})$, one concludes that

$$\widetilde{\Theta}_{(\chi\cap\chi',\chi)}(\widetilde{t}_\chi) = \widetilde{\Theta}_{(\chi\cap\chi',\chi')}(\widetilde{t}_{\chi'})$$

that is, $\widetilde{t}(\overline{\chi}(\gamma)) \in \widetilde{F}(\gamma)$.

Now, it remains to prove that the equivalence class determined by $\widetilde{t}(\overline{\chi}(\gamma)) \in \widetilde{F}(\gamma)$, that is, $[\widetilde{t}(\overline{\chi}(\gamma))] \in \overline{G}(\gamma)$, is independent of the choice of representatives for the classes $\overline{g} = [\widetilde{g}(\overline{\xi}(\gamma))]$ and $\overline{h} = [\widetilde{h}(\overline{\eta}(\gamma))]$, that is, we must prove that: if

$$\widetilde{u}\Big(\overline{\nu}(\gamma)\Big) = \Big\{\widetilde{u}_\nu\Big\}_{\nu\in\overline{\nu}(\gamma)} \approx \widetilde{g}\Big(\overline{\xi}(\gamma)\Big) \quad \text{and} \quad \widetilde{v}\Big(\overline{\omega}(\gamma)\Big) = \Big\{\widetilde{v}_\omega\Big\}_{\omega\in\overline{\omega}(\gamma)} \approx \widetilde{h}\Big(\overline{\eta}(\gamma)\Big),$$

then, the family

$$\widetilde{s}\Big(\overline{\alpha}(\gamma)\Big) = \Big\{\widetilde{s}_\alpha\Big\}_{\alpha\in\overline{\alpha}(\gamma)},$$

where $\overline{\alpha}(\gamma) \subseteq \Gamma(\gamma)$ is the cover of $\gamma$ given by

$$\overline{\alpha}(\gamma) := \Big\{\alpha := \nu \cap \omega \quad : \quad \nu \in \overline{\nu}(\gamma), \quad \omega \in \overline{\omega}(\gamma) \quad \text{and} \quad \nu \cap \omega \neq \varnothing\Big\}$$

and

$$\widetilde{s}_{\alpha=\nu\cap\omega} \coloneqq \widetilde{\Theta}_{(\nu\cap\omega,\nu)}(\widetilde{u}_\nu) + \widetilde{\Theta}_{(\nu\cap\omega,\omega)}(\widetilde{v}_\omega)$$

for every $\alpha = \nu \cap \omega \in \overline{\alpha}(\gamma)$, is such that

$$\widetilde{s}\Big(\overline{\alpha}(\gamma)\Big) \approx \widetilde{t}\Big(\overline{\chi}(\gamma)\Big).$$

Let then $\widetilde{u}(\overline{\nu}(\gamma))$, $\widetilde{v}(\overline{\omega}(\gamma))$ and $\widetilde{s}(\overline{\alpha}(\gamma))$ be as above. Since, by hypothesis,

$$\widetilde{u}\Big(\overline{\nu}(\gamma)\Big) \approx \widetilde{g}\Big(\overline{\xi}(\gamma)\Big) \quad \text{and} \quad \widetilde{v}\Big(\overline{\omega}(\gamma)\Big) \approx \widetilde{h}\Big(\overline{\eta}(\gamma)\Big),$$

then, taking into account the definition of the relation $\approx$ (see Proposition 4.12), we get:

$$\widetilde{\Theta}_{(\nu\cap\xi,\nu)}(\widetilde{u}_\nu) = \widetilde{\Theta}_{(\nu\cap\xi,\xi)}(\widetilde{g}_\xi) \tag{4.15-1}$$

for every $\nu \in \overline{\nu}(\gamma)$ and $\xi \in \overline{\xi}(\gamma)$ such that $\nu \cap \xi \neq \varnothing$, and

$$\widetilde{\Theta}_{(\omega\cap\eta,\omega)}(\widetilde{v}_\omega) = \widetilde{\Theta}_{(\omega\cap\eta,\eta)}(\widetilde{h}_\eta) \tag{4.15-2}$$

for every $\omega \in \overline{\omega}(\alpha)$ and $\eta \in \overline{\eta}(\gamma)$ such that $\omega \cap \eta \neq \varnothing$.

We want to prove, based on the hypotheses (4.15-1) and (4.15-2) above, that $\widetilde{s}(\overline{\alpha}(\gamma)) \approx \widetilde{t}(\overline{\chi}(\gamma))$, i.e., that

$$\widetilde{\Theta}_{(\alpha\cap\chi,\alpha)}(\widetilde{s}_\alpha) = \widetilde{\Theta}_{(\alpha\cap\chi,\chi)}(\widetilde{t}_\chi)$$

for every $\alpha \in \overline{\alpha}(\gamma)$ and $\chi \in \overline{\chi}(\gamma)$ such that $\alpha \cap \chi \neq \varnothing$.

Let us then calculate $\widetilde{\Theta}_{(\alpha\cap\chi,\alpha)}(\widetilde{s}_\alpha) - \widetilde{\Theta}_{(\alpha\cap\chi,\chi)}(\widetilde{t}_\chi)$ with $\alpha = \nu\cap\omega \in \overline{\alpha}(\gamma)$ and $\chi = \xi\cap\eta \in \overline{\chi}(\gamma)$ arbitrarily chosen and such that $\alpha \cap \chi = \nu \cap \omega \cap \xi \cap \eta \neq \varnothing$. We have:

$$\begin{aligned}
&\widetilde{\Theta}_{(\alpha\cap\chi,\alpha)}(\widetilde{s}_\alpha) - \widetilde{\Theta}_{(\alpha\cap\chi,\chi)}(\widetilde{t}_\chi) = \\
&= \widetilde{\Theta}_{(\alpha\cap\chi,\alpha)}\Big(\widetilde{\Theta}_{(\alpha,\nu)}(\widetilde{u}_\nu) + \widetilde{\Theta}_{(\alpha,\omega)}(\widetilde{v}_\omega)\Big) - \widetilde{\Theta}_{(\alpha\cap\chi,\chi)}\Big(\widetilde{\Theta}_{(\chi,\xi)}(\widetilde{g}_\xi) + \widetilde{\Theta}_{(\chi,\eta)}(\widetilde{h}_\eta)\Big) = \\
&= \widetilde{\Theta}_{(\alpha\cap\chi,\nu)}(\widetilde{u}_\nu) - \widetilde{\Theta}_{(\alpha\cap\chi,\xi)}(\widetilde{g}_\xi) + \widetilde{\Theta}_{(\alpha\cap\chi,\omega)}(\widetilde{v}_\omega) - \widetilde{\Theta}_{(\alpha\cap\chi,\eta)}(\widetilde{h}_\eta) = \\
&= \widetilde{\Theta}_{(\alpha\cap\chi,\nu\cap\xi)}\Big(\widetilde{\Theta}_{(\nu\cap\xi,\nu)}(\widetilde{u}_\nu)\Big) - \widetilde{\Theta}_{(\alpha\cap\chi,\nu\cap\xi)}\Big(\widetilde{\Theta}_{(\nu\cap\xi,\xi)}(\widetilde{g}_\xi)\Big) + \\
&+ \widetilde{\Theta}_{(\alpha\cap\chi,\omega\cap\eta)}\Big(\widetilde{\Theta}_{(\omega\cap\eta,\omega)}(\widetilde{v}_\omega)\Big) - \widetilde{\Theta}_{(\alpha\cap\chi,\omega\cap\eta)}\Big(\widetilde{\Theta}_{(\omega\cap\eta,\eta)}(\widetilde{h}_\eta)\Big) = \\
&= \widetilde{\Theta}_{(\alpha\cap\chi,\nu\cap\xi)}\Big(\widetilde{\Theta}_{(\nu\cap\xi,\nu)}(\widetilde{u}_\nu) - \widetilde{\Theta}_{(\nu\cap\xi,\xi)}(\widetilde{g}_\xi)\Big) + \widetilde{\Theta}_{(\alpha\cap\chi,\omega\cap\eta)}\Big(\widetilde{\Theta}_{(\omega\cap\eta,\omega)}(\widetilde{v}_\omega) - \widetilde{\Theta}_{(\omega\cap\eta,\eta)}(\widetilde{h}_\eta)\Big).
\end{aligned}$$

Finally, considering (4.15-1) and (4.15-2), one obtains that

$$\widetilde{\Theta}_{(\alpha\cap\chi,\alpha)}(\widetilde{s}_\alpha) = \widetilde{\Theta}_{(\alpha\cap\chi,\chi)}(\widetilde{t}_\chi)$$

for $\alpha \in \overline{\alpha}(\gamma)$ and $\chi \in \overline{\chi}(\gamma)$ such that $\alpha \cap \chi \neq \varnothing$. ∎

## 4.16 Proposition

*The set $\overline{G}(\gamma)$ (defined in 4.13) equipped with the addition operation, $+$, defined in Proposition 4.15, is an abelian group. Furthermore, under the hypothesis of existence of $\widehat{\mathscr{G}}(I)$, the function $A_\gamma : \overline{G}(\gamma) \longrightarrow \widehat{G}(\gamma)$ (defined in Proposition 4.12) is an isomorphism from the group $\overline{G}(\gamma)$ onto the group $\widehat{G}(\gamma)$.*

*Proof.*

**(a)** From the definition of $+ : \overline{G}(\gamma) \times \overline{G}(\gamma) \longrightarrow \overline{G}(\gamma)$ trivially follows that

$$\overline{g} + \overline{h} = \overline{h} + \overline{g}$$

for every $\overline{g}, \overline{h} \in \overline{G}(\gamma)$.

**(b)** Let us prove now that the operation $+$ is associative. In order to do so, let

$$\begin{aligned} \overline{f} &= \left[\widetilde{f}\Big(\overline{\nu}(\gamma)\Big) = \left\{\widetilde{f}_\nu\right\}_{\nu\in\overline{\nu}(\gamma)}\right], \\ \overline{g} &= \left[\widetilde{g}\Big(\overline{\xi}(\gamma)\Big) = \left\{\widetilde{g}_\xi\right\}_{\xi\in\overline{\xi}(\gamma)}\right] \quad \text{and} \\ \overline{h} &= \left[\widetilde{h}\Big(\overline{\eta}(\gamma)\Big) = \left\{\widetilde{h}_\eta\right\}_{\eta\in\overline{\eta}(\gamma)}\right] \end{aligned}$$

be arbitrarily chosen in $\overline{G}(\gamma)$, with which we calculate:

**(b-1)** $\overline{f} + \overline{g} = \left[\widetilde{t}\Big(\overline{\chi}(\gamma)\Big) = \left\{\widetilde{t}_\chi\right\}_{\chi\in\overline{\chi}(\gamma)}\right]$ where

$$\overline{\chi}(\gamma) \coloneqq \Big\{\chi \coloneqq \nu \cap \xi \quad : \quad \nu \in \overline{\nu}(\gamma), \quad \xi \in \overline{\xi}(\gamma) \quad \text{and} \quad \nu \cap \xi \neq \varnothing\Big\}$$

and

$$\tilde{t}_{\chi=\nu\cap\xi} \coloneqq \widetilde{\Theta}_{(\chi,\nu)}(\widetilde{f}_\nu) + \widetilde{\Theta}_{(\chi,\xi)}(\widetilde{g}_\xi) \quad \text{for every} \quad \chi = \nu \cap \xi \in \overline{\chi}(\gamma);$$

**(b-2)** $(\overline{f} + \overline{g}) + \overline{h} = \left[\widetilde{t}\Big(\overline{\chi}(\gamma)\Big)\right] + \left[\widetilde{h}\Big(\overline{\eta}(\gamma)\Big)\right] = \left[\widetilde{u}\Big(\overline{\omega}(\gamma)\Big) = \left\{\widetilde{u}_\omega\right\}_{\omega\in\overline{\omega}(\gamma)}\right]$ where

$$\overline{\omega}(\gamma) \coloneqq \Big\{\omega \coloneqq \chi \cap \eta = \nu \cap \xi \cap \eta \quad : \quad \chi = \nu \cap \xi \in \overline{\chi}(\gamma), \quad \eta \in \overline{\eta}(\gamma) \quad \text{and} \quad \chi \cap \eta \neq \varnothing\Big\}$$

and

$$\widetilde{u}_{\omega=\chi\cap\eta} \coloneqq \widetilde{\Theta}_{(\omega,\chi)}(\widetilde{t}_\chi) + \widetilde{\Theta}_{(\omega,\eta)}(\widetilde{h}_\eta)$$

for every $\omega = \chi \cap \eta \in \overline{\omega}(\gamma)$ or, taking into account the expression for $\widetilde{t}_\chi$ in (b-1),

$$\widetilde{u}_{\omega=\chi\cap\eta} = \Big(\widetilde{\Theta}_{(\omega,\nu)}(\widetilde{f}_\nu) + \widetilde{\Theta}_{(\omega,\xi)}(\widetilde{g}_\xi)\Big) + \widetilde{\Theta}_{(\omega,\eta)}(\widetilde{h}_\eta)$$

for every $\omega = \chi \cap \eta = \nu \cap \xi \cap \eta \in \overline{\omega}(\gamma)$;

**(b-3)** $\overline{g} + \overline{h} = \left[\widetilde{v}\Big(\overline{\alpha}(\gamma)\Big) = \Big\{\widetilde{v}_\alpha\Big\}_{\alpha\in\overline{\alpha}(\gamma)}\right]$ where

$$\overline{\alpha}(\gamma) := \Big\{\alpha := \xi \cap \eta \quad : \quad \xi \in \overline{\xi}(\gamma), \quad \eta \in \overline{\eta}(\gamma) \quad \text{and} \quad \xi \cap \eta \neq \varnothing\Big\}$$

and

$$\widetilde{v}_{\alpha=\xi\cap\eta} := \widetilde{\Theta}_{(\alpha,\xi)}(\widetilde{g}_\xi) + \widetilde{\Theta}_{(\alpha,\eta)}(\widetilde{h}_\eta)$$

for every $\alpha = \xi \cap \eta \in \overline{\alpha}(\gamma)$;

**(b-4)** $\overline{f} + (\overline{g} + \overline{h}) = \left[\widetilde{f}\Big(\overline{\nu}(\gamma)\Big)\right] + \left[\widetilde{v}\Big(\overline{\alpha}(\gamma)\Big)\right] = \left[\widetilde{s}\Big(\overline{\pi}(\gamma)\Big) = \Big\{\widetilde{s}_\pi\Big\}_{\pi\in\overline{\pi}(\gamma)}\right]$ where

$$\overline{\pi}(\gamma) := \Big\{\pi := \nu \cap \alpha \quad : \quad \nu \in \overline{\nu}(\gamma), \quad \alpha = \xi \cap \eta \in \overline{\alpha}(\gamma) \quad \text{and} \quad \nu \cap \alpha = \nu \cap \xi \cap \eta \neq \varnothing\Big\}$$

and

$$\widetilde{s}_{\pi=\nu\cap\alpha} := \widetilde{\Theta}_{(\pi,\nu)}(\widetilde{f}_\nu) + \widetilde{\Theta}_{(\pi,\alpha)}(\widetilde{v}_\alpha) \quad \text{for every} \quad \pi = \nu \cap \alpha = \nu \cap \xi \cap \eta \in \overline{\pi}(\gamma),$$

or, considering the expression for $\widetilde{v}_\alpha$ in (b-3),

$$\widetilde{s}_\pi = \widetilde{\Theta}_{(\pi,\nu)}(\widetilde{f}_\nu) + \Big(\widetilde{\Theta}_{(\pi,\xi)}(\widetilde{g}_\xi) + \widetilde{\Theta}_{(\pi,\eta)}(\widetilde{h}_\eta)\Big)$$

or yet,

$$\widetilde{s}_\pi = \Big(\widetilde{\Theta}_{(\pi,\nu)}(\widetilde{f}_\nu) + \widetilde{\Theta}_{(\pi,\xi)}(\widetilde{g}_\xi)\Big) + \widetilde{\Theta}_{(\pi,\eta)}(\widetilde{h}_\eta)$$

for every $\pi \in \overline{\pi}(\gamma)$.

Comparing now the expressions obtained for $\widetilde{u}_\omega$ (with $\omega \in \overline{\omega}(\gamma)$) and $\widetilde{s}_\pi$ (with $\pi \in \overline{\pi}(\gamma)$) in (b-2) and (b-4), respectively, and observing that $\overline{\pi}(\gamma) = \overline{\omega}(\gamma)$, one concludes that

$$\widetilde{u}_\pi = \widetilde{s}_\pi \quad \text{for every} \quad \pi \in \overline{\pi}(\gamma) = \overline{\omega}(\gamma),$$

that is,

$$(\overline{f} + \overline{g}) + \overline{h} = \left[\widetilde{u}\Big(\overline{\omega}(\gamma)\Big)\right] = \left[\widetilde{s}\Big(\overline{\pi}(\gamma)\Big)\right] = \overline{f} + (\overline{g} + \overline{h}).$$

**(c)** Let us prove now the existence of $\overline{0} \in \overline{G}(\gamma)$ such that

$$\overline{g} + \overline{0} = \overline{g} \quad \text{for every} \quad \overline{g} \in \overline{G}(\gamma).$$

In order to do so, let us consider a family

$$\widetilde{0}\Big(\overline{\xi}(\gamma)\Big) = \Big\{\widetilde{0}_\xi\Big\}_{\xi\in\overline{\xi}(\gamma)}$$

where $\widetilde{0}_\xi \in \widetilde{G}(\xi)$ is the additive neutral of the group $\widetilde{G}(\xi)$ ($\overline{\xi}(\gamma) \subseteq \Gamma(\gamma)$ is an arbitrarily fixed cover of $\gamma$). Since $\widetilde{\Theta}_{(\xi\cap\xi',\xi)} : \widetilde{G}(\xi) \longrightarrow \widetilde{G}(\xi\cap\xi')$ as well as $\widetilde{\Theta}_{(\xi\cap\xi',\xi')} : \widetilde{G}(\xi') \longrightarrow \widetilde{G}(\xi\cap\xi')$ are homomorphisms between groups, we have

$$\widetilde{\Theta}_{(\xi\cap\xi',\xi)}(\widetilde{0}_\xi) = \widetilde{0}_{\xi\cap\xi'} = \widetilde{\Theta}_{(\xi\cap\xi',\xi')}(\widetilde{0}_{\xi'})$$

for every $\xi, \xi' \in \overline{\xi}(\gamma)$ such that $\xi \cap \xi' \neq \varnothing$. This proves that the family $\widetilde{0}(\overline{\xi}(\gamma))$ is coherent in $\widetilde{\mathscr{G}}(\gamma)$, that is, $\widetilde{0}(\overline{\xi}(\gamma)) \in \widetilde{F}(\gamma)$. Thus, we define

$$\overline{0} := \left[\widetilde{0}\left(\overline{\xi}(\gamma)\right)\right] \in \overline{G}(\gamma).$$

Let now $\overline{g} = [\widetilde{g}(\overline{\eta}(\gamma)) = \{\widetilde{g}_\eta\}_{\eta \in \overline{\eta}(\gamma)}] \in \overline{G}(\gamma)$ be arbitrarily fixed and let us calculate $\overline{g} + \overline{0}$. We have

$$\overline{g} + \overline{0} = \left[\widetilde{h}\left(\overline{\chi}(\gamma)\right) = \left\{\widetilde{h}_\chi\right\}_{\chi \in \overline{\chi}(\gamma)}\right]$$

where

$$\overline{\chi}(\gamma) := \left\{\chi := \eta \cap \xi \quad : \quad \eta \in \overline{\eta}(\gamma), \quad \xi \in \overline{\xi}(\gamma) \quad \text{and} \quad \eta \cap \xi \neq \varnothing\right\}$$

and

$$\widetilde{h}_{\chi = \eta \cap \xi} := \widetilde{\Theta}_{(\chi,\eta)}(\widetilde{g}_\eta) + \widetilde{\Theta}_{(\chi,\xi)}(\widetilde{0}_\xi) = \widetilde{\Theta}_{(\chi,\eta)}(\widetilde{g}_\eta) + \widetilde{0}_\chi = \widetilde{\Theta}_{(\chi,\eta)}(\widetilde{g}_\eta)$$

for every $\chi = \eta \cap \xi \in \overline{\chi}(\gamma)$.

Hence, we have $\overline{g} + \overline{0} = \overline{g}$ if and if

$$\widetilde{h}\left(\overline{\chi}(\gamma)\right) \approx \widetilde{g}\left(\overline{\eta}(\gamma)\right)$$

that is, if and only if

$$\widetilde{\Theta}_{(\chi \cap \eta',\chi)}(\widetilde{h}_\chi) = \widetilde{\Theta}_{(\chi \cap \eta',\eta')}(\widetilde{g}_{\eta'})$$

for every $\chi = \eta \cap \xi \in \overline{\chi}(\gamma)$ and $\eta' \in \overline{\eta}(\gamma)$ such that $\chi \cap \eta' = \eta \cap \xi \cap \eta' \neq \varnothing$.

Let us then prove this equivalence. From the above expression for $\widetilde{h}_\chi$ comes that, for any $\chi = \eta \cap \xi \in \overline{\chi}(\gamma)$ and $\eta' \in \overline{\eta}(\gamma)$ such that $\chi \cap \eta' = \eta \cap \xi \cap \eta' \neq \varnothing$,

$$\widetilde{\Theta}_{(\chi \cap \eta',\chi)}(\widetilde{h}_\chi) = \widetilde{\Theta}_{(\chi \cap \eta',\chi)}\left(\widetilde{\Theta}_{(\chi,\eta)}(\widetilde{g}_\eta)\right) = \widetilde{\Theta}_{(\chi \cap \eta',\eta)}(\widetilde{g}_\eta) = \widetilde{\Theta}_{(\chi \cap \eta',\eta \cap \eta')}\left(\widetilde{\Theta}_{(\eta \cap \eta',\eta)}(\widetilde{g}_\eta)\right).$$

Now, $\widetilde{g}(\overline{\eta}(\gamma)) \in \widetilde{F}(\gamma)$, that is, $\widetilde{g}(\overline{\eta}(\gamma))$ is a coherent family in $\widetilde{\mathscr{G}}(\gamma)$, and hence, for $\eta \cap \eta' \neq \varnothing$,

$$\widetilde{\Theta}_{(\eta \cap \eta',\eta)}(\widetilde{g}_\eta) = \widetilde{\Theta}_{(\eta \cap \eta',\eta')}(\widetilde{g}_{\eta'}),$$

which allows us to write

$$\widetilde{\Theta}_{(\chi \cap \eta',\chi)}(\widetilde{h}_\chi) = \widetilde{\Theta}_{(\chi \cap \eta',\eta \cap \eta')}\left(\widetilde{\Theta}_{(\eta \cap \eta',\eta')}(\widetilde{g}_{\eta'})\right),$$

that is,

$$\widetilde{\Theta}_{(\chi \cap \eta',\chi)}(\widetilde{h}_\chi) = \widetilde{\Theta}_{(\chi \cap \eta',\eta')}(\widetilde{g}_{\eta'}).$$

**(d)** Finally, we must prove that for each $\overline{g} \in \overline{G}(\gamma)$ there exists $\overline{h} \in \overline{G}(\gamma)$ such that

$$\overline{g} + \overline{h} = \overline{0}.$$

In order to do so, let

$$\overline{g} = \left[\widetilde{g}\left(\overline{\xi}(\gamma)\right) = \left\{\widetilde{g}_\xi\right\}_{\xi \in \overline{\xi}(\gamma)}\right]$$

be arbitrarily chosen in $\overline{G}(\gamma)$ and take

$$\widetilde{h}\Big(\overline{\xi}(\gamma)\Big) := \Big\{\widetilde{h}_\xi := -\widetilde{g}_\xi\Big\}_{\xi\in\overline{\xi}(\gamma)}.$$

First, let us note that, for $\xi,\xi' \in \overline{\xi}(\gamma)$ such that $\xi\cap\xi' \neq \varnothing$,

$$\widetilde{\Theta}_{(\xi\cap\xi',\xi)}(\widetilde{h}_\xi) - \widetilde{\Theta}_{(\xi\cap\xi',\xi')}(\widetilde{h}_{\xi'}) = -\widetilde{\Theta}_{(\xi\cap\xi',\xi)}(\widetilde{g}_\xi) + \widetilde{\Theta}_{(\xi\cap\xi',\xi')}(\widetilde{g}_{\xi'})$$

and since $\widetilde{\Theta}_{(\xi\cap\xi',\xi)}(\widetilde{g}_\xi) = \widetilde{\Theta}_{(\xi\cap\xi',\xi')}(\widetilde{g}_{\xi'})$ once $\widetilde{g}(\overline{\xi}(\gamma)) \in \widetilde{F}(\gamma)$, comes that $\widetilde{\Theta}_{(\xi\cap\xi',\xi)}(\widetilde{h}_\xi) = \widetilde{\Theta}_{(\xi\cap\xi',\xi')}(\widetilde{h}_{\xi'})$, that is, $\widetilde{h}(\overline{\xi}(\gamma)) \in \widetilde{F}(\gamma)$ and, therefore,

$$\overline{h} := \Big[\widetilde{h}\Big(\overline{\xi}(\gamma)\Big)\Big] \in \overline{G}(\gamma).$$

Now, we calculate:

$$\overline{g} + \overline{h} = \Big[\widetilde{g}\Big(\overline{\xi}(\gamma)\Big)\Big] + \Big[\widetilde{h}\Big(\overline{\xi}(\gamma)\Big)\Big] = \Big[\widetilde{f}\Big(\overline{\chi}(\gamma)\Big) = \Big\{\widetilde{f}_\chi\Big\}_{\chi\in\overline{\chi}(\gamma)}\Big]$$

where

$$\overline{\chi}(\gamma) := \Big\{\chi = \xi\cap\xi' \quad : \quad \xi,\xi' \in \overline{\xi}(\gamma) \quad \text{and} \quad \xi\cap\xi' \neq \varnothing\Big\}$$

and

$$\begin{aligned}\widetilde{f}_{\chi=\xi\cap\xi'} &:= \widetilde{\Theta}_{(\xi\cap\xi',\xi)}(\widetilde{g}_\xi) + \widetilde{\Theta}_{(\xi\cap\xi',\xi')}(\widetilde{h}_{\xi'}) = \\ &= \widetilde{\Theta}_{(\xi\cap\xi',\xi)}(\widetilde{g}_\xi) + \widetilde{\Theta}_{(\xi\cap\xi',\xi')}(-\widetilde{g}_{\xi'}) = \\ &= \widetilde{\Theta}_{(\xi\cap\xi',\xi)}(\widetilde{g}_\xi) - \widetilde{\Theta}_{(\xi\cap\xi',\xi')}(\widetilde{g}_{\xi'}).\end{aligned}$$

But, as we saw above, $\widetilde{g}(\overline{\xi}(\gamma)) \in \widetilde{F}(\gamma)$ and, hence, $\widetilde{\Theta}_{(\xi\cap\xi',\xi)}(\widetilde{g}_\xi) - \widetilde{\Theta}_{(\xi\cap\xi',\xi')}(\widetilde{g}_{\xi'}) = \widetilde{0}_{\xi\cap\xi'}$. Thus,

$$\widetilde{f}_\chi = \widetilde{0}_\chi \quad \text{for every} \quad \chi = \xi\cap\xi' \in \overline{\chi}(\gamma).$$

Therefore,

$$\overline{g} + \overline{h} = \Big[\widetilde{f}\Big(\overline{\chi}(\gamma)\Big)\Big] = \Big[\widetilde{0}\Big(\overline{\chi}(\gamma)\Big)\Big]$$

where $\widetilde{0}(\overline{\chi}(\gamma)) := \{\widetilde{0}_\chi\}_{\chi\in\overline{\chi}(\gamma)}$.

Since, being $\overline{\eta}(\gamma) \subseteq \Gamma(\gamma)$ any cover of $\gamma$, the families

$$\widetilde{0}\Big(\overline{\chi}(\gamma)\Big) = \Big\{\widetilde{0}_\chi\Big\}_{\chi\in\overline{\chi}(\gamma)} \quad \text{and} \quad \widetilde{0}\Big(\overline{\eta}(\gamma)\Big) = \Big\{\widetilde{0}_\eta\Big\}_{\eta\in\overline{\eta}(\gamma)}$$

($\widetilde{0}_\chi \in \widetilde{G}(\chi)$ and $\widetilde{0}_\eta \in \widetilde{G}(\eta)$ are, respectively, the additive neutrals of the groups $\widetilde{G}(\chi)$ and $\widetilde{G}(\eta)$) are equivalent to each other, that is,

$$\widetilde{0}\Big(\overline{\chi}(\gamma)\Big) \approx \widetilde{0}\Big(\overline{\eta}(\gamma)\Big)$$

once that

$$\widetilde{\Theta}_{(\chi\cap\eta,\chi)}(\widetilde{0}_\chi) - \widetilde{\Theta}_{(\chi\cap\eta,\eta)}(\widetilde{0}_\eta) = \widetilde{0}_{\chi\cap\eta} - \widetilde{0}_{\chi\cap\eta} = \widetilde{0}_{\chi\cap\eta}$$

for every $\chi \in \overline{\chi}(\gamma)$, $\eta \in \overline{\eta}(\gamma)$ such that $\chi \cap \eta \neq \varnothing$, we then have proved that $[\widetilde{0}(\overline{\chi}(\gamma))]$ is the additive neutral $\overline{0} \in \overline{G}(\gamma)$ defined in (c). Hence,

$$\left[\widetilde{0}\Big(\overline{\chi}(\gamma)\Big)\right] = \overline{0}$$

and

$$\overline{g} + \overline{h} = \overline{0}. \qquad \blacksquare$$

## 4.17 Remarks

Our developments until here provided us, for each $\gamma \in \Gamma(I)$, with an abelian group $\overline{G}(\gamma)$ entirely defined with the elements of the extension $\widetilde{\mathscr{G}}(I)$ of the $S$-space $\mathscr{G}(I)$, such that, under the hypothesis of existence of $\widehat{\mathscr{G}}(I)$, the groups $\overline{G}(\gamma)$ and $\widehat{G}(\gamma)$ are isomorphic and the function $A_\gamma : \overline{G}(\gamma) \longrightarrow \widehat{G}(\gamma)$ defined by

$$A_\gamma\left(\overline{g} = \left[\Big\{\widetilde{g}_\xi\Big\}_{\xi \in \overline{\xi}(\gamma)}\right]\right) = \widehat{g},$$

with $\widehat{g} \in \widehat{G}(\gamma)$ such that

$$\widehat{\Theta}_{(\xi,\gamma)}(\widehat{g}) = \widetilde{g}_\xi \quad \text{for every} \quad \xi \in \overline{\xi}(\gamma),$$

is an isomorphism from the group $\overline{G}(\gamma)$ onto the group $\widehat{G}(\gamma)$. Hence, existing $\widehat{\mathscr{G}}(I)$, since $\widehat{\mathscr{G}}(I)$ is an extension of $\widetilde{\mathscr{G}}(I)$, then, $\widetilde{G}(\gamma) \subseteq \widehat{G}(\gamma)$ is a subgroup of the group $\widehat{G}(\gamma)$, for each $\gamma \in \Gamma(I)$ (more precisely, $\widehat{G}(\gamma)$ admits $\widetilde{G}(\gamma)$ as a subgroup (see Definition 1.6)). Therefore, in this case (of existence of $\widehat{\mathscr{G}}(I)$),

$$A_\gamma^{-1}\Big(\widetilde{G}(\gamma)\Big) = \Big\{\overline{g} \in \overline{G}(\gamma) \colon A_\gamma(\overline{g}) \in \widetilde{G}(\gamma) \subseteq \widehat{G}(\gamma)\Big\} =$$
$$= \Big\{A_\gamma^{-1}(\widetilde{g}) \colon \widetilde{g} \in \widetilde{G}(\gamma) \subseteq \widehat{G}(\gamma)\Big\}$$

is a subgroup of the group $\overline{G}(\gamma)$, isomorphic to the group $\widetilde{G}(\gamma)$.

For $\widetilde{g} \in \widetilde{G}(\gamma) \subseteq \widehat{G}(\gamma)$ arbitrarily fixed, let

$$A_\gamma^{-1}(\widetilde{g}) = \left[\widetilde{h}\Big(\overline{\xi}(\gamma)\Big) = \Big\{\widetilde{h}_\xi\Big\}_{\xi \in \overline{\xi}(\gamma)}\right].$$

Hence, it results from the definition of the function $A_\gamma$ that $\widetilde{g}$ is the only member of $\widetilde{G}(\gamma) \subseteq \widehat{G}(\gamma)$ such that

$$\widehat{\Theta}_{(\xi,\gamma)}(\widetilde{g}) = \widetilde{h}_\xi \quad \text{for every} \quad \xi \in \overline{\xi}(\gamma),$$

or yet, since $\widehat{\Theta}(\Delta(I))$ is a prolongation of $\widetilde{\Theta}(\Delta(I))$, such that

$$\widetilde{\Theta}_{(\xi,\gamma)}(\widetilde{g}) = \widetilde{h}_\xi \quad \text{for every} \quad \xi \in \overline{\xi}(\gamma).$$

Therefore,

$$A_\gamma^{-1}(\widetilde{g}) = \Big[\widetilde{h}\Big(\overline{\xi}(\gamma)\Big) = \Big\{\widetilde{\Theta}_{(\xi,\gamma)}(\widetilde{g})\Big\}_{\xi\in\overline{\xi}(\gamma)}\Big]$$

and hence, the subset $A_\gamma^{-1}(\widetilde{G}(\gamma)) \subseteq \overline{G}(\gamma)$ can be expressed in the following form:

$$A_\gamma^{-1}\Big(\widetilde{G}(\gamma)\Big) = \left\{\Big[\Big\{\widetilde{\Theta}_{(\xi,\gamma)}(\widetilde{g})\Big\}_{\xi\in\overline{\xi}(\gamma)}\Big] : \widetilde{g} \in \widetilde{G}(\gamma)\right\}.$$

Let us observe now that the set

$$\left\{\Big[\Big\{\widetilde{\Theta}_{(\xi,\gamma)}(\widetilde{g})\Big\}_{\xi\in\overline{\xi}(\gamma)}\Big] : \widetilde{g} \in \widetilde{G}(\gamma)\right\} \subseteq \overline{G}(\gamma)$$

is a well-defined subset of $\overline{G}(\gamma)$ either existing $\widehat{\mathscr{G}}(I)$ or not and, if $\widehat{\mathscr{G}}(I)$ does exists, it is a subgroup of $\overline{G}(\gamma)$ isomorphic to the group $\widetilde{G}(\gamma)$, that is, according to Definition 1.6(a), the group $\overline{G}(\gamma)$ admits $\widetilde{G}(\gamma)$ as a subgroup. If now we remember that our purpose is to build an extension $\overline{\mathscr{G}}(I)$ of $\widetilde{\mathscr{G}}(I)$ that is isomorphic to $\widehat{\mathscr{G}}(I)$ if the latter exists, and that, to do so, it is necessary that, for each $\gamma \in \Gamma(I)$, the group $\overline{G}(\gamma)$ admits the group $\widetilde{G}(\gamma)$ as a subgroup, we have no other option than the following proposition.

## 4.18 Proposition

*Let $\widetilde{\mathscr{G}}(\gamma)$ be a S-subspace of $\widetilde{\mathscr{G}}(I)$ arbitrarily fixed and $\overline{G}(\gamma)$ be the abelian group defined in Proposition 4.16. The function*

$$\begin{aligned} b_\gamma : \widetilde{G}(\gamma) &\longrightarrow \overline{G}(\gamma) \\ \widetilde{g} &\longmapsto b_\gamma(\widetilde{g}) \end{aligned}$$

*defined by*

$$b_\gamma(\widetilde{g}) := \Big[\Big\{\widetilde{\Theta}_{(\xi,\gamma)}(\widetilde{g})\Big\}_{\xi\in\overline{\xi}(\gamma)}\Big],$$

*where $\overline{\xi}(\gamma) \subseteq \Gamma(\gamma)$ is a cover of $\gamma$, it is well-defined, injective, and*

$$b_\gamma(\widetilde{g} + \widetilde{h}) = b_\gamma(\widetilde{g}) + b_\gamma(\widetilde{h})$$

*for every $\widetilde{g}, \widetilde{h} \in \widetilde{G}(\gamma)$, that is, $b_\gamma$ is an injective homomorphism from the group $\widetilde{G}(\gamma)$ into the group $\overline{G}(\gamma)$.*

*Furthermore, existing $\widehat{\mathscr{G}}(I)$, the function*

$$\begin{aligned} A_\gamma \circ b_\gamma : \widetilde{G}(\gamma) &\longrightarrow \widehat{G}(\gamma) \\ \widetilde{g} &\longmapsto (A_\gamma \circ b_\gamma)(\widetilde{g}) \end{aligned}$$

*is the identity on $\widetilde{G}(\gamma)$, that is,*

$$(A_\gamma \circ b_\gamma)(\widetilde{g}) = A_\gamma\Big(b_\gamma(\widetilde{g})\Big) = \widetilde{g}$$

*for every $\widetilde{g} \in \widetilde{G}(\gamma)$.*

*Proof.* We must first prove that $b_\gamma$ is well-defined. In other words, that the family $\{\widetilde{\Theta}_{(\xi,\gamma)}(\widetilde{g})\}_{\xi\in\bar{\xi}(\gamma)}$ is a coherent family in $\widetilde{\mathscr{G}}(\gamma)$ whatever the cover $\overline{\xi}(\gamma) \subseteq \Gamma(\gamma)$ of $\gamma$ taken, and yet that the equivalence class determined by the family $\{\widetilde{\Theta}_{(\xi,\gamma)}(\widetilde{g})\}_{\xi\in\bar{\xi}(\gamma)} \in \widetilde{F}(\gamma)$, $[\{\widetilde{\Theta}_{(\xi,\gamma)}(\widetilde{g})\}_{\xi\in\bar{\xi}(\gamma)}]$, does not depend on the cover $\overline{\xi}(\gamma)$ in the following sense: if $\overline{\eta}(\gamma) \subseteq \Gamma(\gamma)$ is another cover of $\gamma$, then,

$$\left\{\widetilde{\Theta}_{(\eta,\gamma)}(\widetilde{g})\right\}_{\eta\in\overline{\eta}(\gamma)} \approx \left\{\widetilde{\Theta}_{(\xi,\gamma)}(\widetilde{g})\right\}_{\xi\in\bar{\xi}(\gamma)}.$$

It is trivial to prove that the family $\{\widetilde{\Theta}_{(\xi,\gamma)}(\widetilde{g})\}_{\xi\in\bar{\xi}(\gamma)}$ belongs to $\widetilde{F}(\gamma)$, whatever the cover $\overline{\xi}(\gamma) \subseteq \Gamma(\gamma)$ of $\gamma$ taken. Let us then prove that, being $\overline{\eta}(\gamma) \subseteq \Gamma(\gamma)$ another cover of $\gamma$ arbitrarily chosen,

$$\left\{\widetilde{\Theta}_{(\eta,\gamma)}(\widetilde{g})\right\}_{\eta\in\overline{\eta}(\gamma)} \approx \left\{\widetilde{\Theta}_{(\xi,\gamma)}(\widetilde{g})\right\}_{\xi\in\bar{\xi}(\gamma)}.$$

But this is trivial as well, since for any $\xi \in \overline{\xi}(\gamma)$ and $\eta \in \overline{\eta}(\gamma)$ such that $\xi \cap \eta \neq \varnothing$, we have

$$\widetilde{\Theta}_{(\xi\cap\eta,\eta)}\Big(\widetilde{\Theta}_{(\eta,\gamma)}(\widetilde{g})\Big) - \widetilde{\Theta}_{(\xi\cap\eta,\xi)}\Big(\widetilde{\Theta}_{(\xi,\gamma)}(\widetilde{g})\Big) = \widetilde{\Theta}_{(\xi\cap\eta,\gamma)}(\widetilde{g}) - \widetilde{\Theta}_{(\xi\cap\eta,\gamma)}(\widetilde{g}) = 0_{\xi\cap\eta}$$

that is,

$$\widetilde{\Theta}_{(\xi\cap\eta,\eta)}\Big(\widetilde{\Theta}_{(\eta,\gamma)}(\widetilde{g})\Big) = \widetilde{\Theta}_{(\xi\cap\eta,\xi)}\Big(\widetilde{\Theta}_{(\xi,\gamma)}(\widetilde{g})\Big)$$

which, taking into account the definition of the relation $\approx$ (see Proposition 4.12), establishes the equivalence between the families $\{\widetilde{\Theta}_{(\eta,\gamma)}(\widetilde{g})\}_{\eta\in\overline{\eta}(\gamma)}$ and $\{\widetilde{\Theta}_{(\xi,\gamma)}(\widetilde{g})\}_{\xi\in\bar{\xi}(\gamma)}$.

Next we prove that, for $\widetilde{g}$ and $\widetilde{h}$ arbitrarily chosen in $\widetilde{G}(\gamma)$,

$$b_\gamma(\widetilde{g} + \widetilde{h}) = b_\gamma(\widetilde{g}) + b_\gamma(\widetilde{g}).$$

We have:

$$b_\gamma(\widetilde{g}) = \left[\left\{\widetilde{\Theta}_{(\xi,\gamma)}(\widetilde{g})\right\}_{\xi\in\bar{\xi}(\gamma)}\right] \quad \text{and} \quad b_\gamma(\widetilde{h}) = \left[\left\{\widetilde{\Theta}_{(\eta,\gamma)}(\widetilde{h})\right\}_{\eta\in\overline{\eta}(\gamma)}\right]$$

being $\overline{\xi}(\gamma) \subseteq \Gamma(\gamma)$ and $\overline{\eta}(\gamma) \subseteq \Gamma(\gamma)$ covers of $\gamma$. Hence, by the definition of addition in $\overline{G}(\gamma)$ (see Proposition 4.15),

$$b_\gamma(\widetilde{g}) + b_\gamma(\widetilde{h}) = \left[\widetilde{t}\Big(\overline{\chi}(\gamma)\Big) = \left\{\widetilde{t}_\chi\right\}_{\chi\in\overline{\chi}(\gamma)}\right]$$

where

$$\overline{\chi}(\gamma) \coloneqq \left\{\chi \coloneqq \xi \cap \eta \quad : \quad \xi \in \overline{\xi}(\gamma), \quad \eta \in \overline{\eta}(\gamma) \quad \text{and} \quad \xi \cap \eta \neq \varnothing\right\}$$

and

$$\begin{aligned}
\widetilde{t}_{\chi=\xi\cap\eta} &= \widetilde{\Theta}_{(\xi\cap\eta,\xi)}\Big(\widetilde{\Theta}_{(\xi,\gamma)}(\widetilde{g})\Big) + \widetilde{\Theta}_{(\xi\cap\eta,\eta)}\Big(\widetilde{\Theta}_{(\eta,\gamma)}(\widetilde{h})\Big) = \\
&= \widetilde{\Theta}_{(\xi\cap\eta,\gamma)}(\widetilde{g}) + \widetilde{\Theta}_{(\xi\cap\eta,\gamma)}(\widetilde{h}) = \\
&= \widetilde{\Theta}_{(\xi\cap\eta,\gamma)}(\widetilde{g} + \widetilde{h})
\end{aligned}$$

for every $\chi = \xi \cap \eta \in \overline{\chi}(\gamma)$. That is,

$$b_\gamma(\widetilde{g}) + b_\gamma(\widetilde{h}) = \left[\left\{\widetilde{\Theta}_{(\chi,\gamma)}(\widetilde{g}+\widetilde{h})\right\}_{\chi\in\overline{\chi}(\gamma)}\right].$$

But, from the definition of the function $b_\gamma : \widetilde{G}(\gamma) \longrightarrow \overline{G}(\gamma)$ we have that

$$\left[\left\{\widetilde{\Theta}_{(\chi,\gamma)}(\widetilde{g}+\widetilde{h})\right\}_{\chi\in\overline{\chi}(\gamma)}\right] = b_\gamma(\widetilde{g}+\widetilde{h})$$

and, hence,

$$b_\gamma(\widetilde{g}) + b_\gamma(\widetilde{h}) = b_\gamma(\widetilde{g}+\widetilde{h}).$$

Finally, we must prove that $b_\gamma$ is injective. In order to do so, let $\widetilde{g}, \widetilde{h} \in \widetilde{G}(\gamma)$ be arbitrarily chosen and such that

$$b_\gamma(\widetilde{g}) = b_\gamma(\widetilde{h}),$$

that is,

$$b_\gamma(\widetilde{g}-\widetilde{h}) = \overline{0}$$

where, according to Proposition 4.16 (more specifically, item (c) of its proof),

$$\overline{0} = \left[\widetilde{0}\left(\overline{\xi}(\gamma)\right) = \left\{\widetilde{0}_\xi\right\}_{\xi\in\overline{\xi}(\gamma)}\right],$$

where $\widetilde{0}_\xi \in \widetilde{G}(\xi)$ is the additive neutral of the group $\widetilde{G}(\xi)$, is the additive neutral of the group $\overline{G}(\gamma)$.

From the definition of the function $b_\gamma$ we have:

$$b_\gamma(\widetilde{g}-\widetilde{h}) = \left[\widetilde{f}\left(\overline{\eta}(\gamma)\right) = \left\{\widetilde{f}_\eta = \widetilde{\Theta}_{(\eta,\gamma)}(\widetilde{g}-\widetilde{h})\right\}_{\eta\in\overline{\eta}(\gamma)}\right]$$

and, hence,

$$\left[\widetilde{f}\left(\overline{\eta}(\gamma)\right)\right] = \left[\widetilde{0}\left(\overline{\xi}(\gamma)\right)\right]$$

that is,

$$\widetilde{f}\left(\overline{\eta}(\gamma)\right) \approx \widetilde{0}\left(\overline{\xi}(\gamma)\right),$$

which means, taking into account the definition of the relation $\approx$, that

$$\widetilde{\Theta}_{(\xi\cap\eta,\eta)}(\widetilde{f}_\eta) = \widetilde{\Theta}_{(\xi\cap\eta,\xi)}(\widetilde{0}_\xi)$$

for every $\xi \in \overline{\xi}(\gamma)$ and $\eta \in \overline{\eta}(\gamma)$ such that $\xi \cap \eta \neq \varnothing$ or yet, since $\widetilde{\Theta}_{(\xi\cap\eta,\xi)}(\widetilde{0}_\xi) = \widetilde{0}_{\xi\cap\eta}$ (the additive neutral of the group $\widetilde{G}(\xi\cap\eta)$),

$$\widetilde{\Theta}_{(\xi\cap\eta,\eta)}(\widetilde{f}_\eta) = \widetilde{0}_{\xi\cap\eta}$$

for every $\xi \in \overline{\xi}(\gamma)$ and $\eta \in \overline{\eta}(\gamma)$ such that $\xi \cap \eta \neq \varnothing$.

On the other hand, the $S$-group $\widetilde{\mathbb{G}}(\eta) = (\widetilde{G}(\eta), \widetilde{H}(\eta))$ is a (strict and) closed extension of the $S$-group $\mathbb{G}(\eta) = (G(\eta), H(\eta))$ and, hence, there exist $\Phi_\eta \in H(\eta)$ and $f_\eta \in G(\eta)$ such that

$$\widetilde{f}_\eta = \widetilde{\Phi}_\eta(f_\eta)$$

and, with this, we obtain that

$$\widetilde{\Theta}_{(\xi\cap\eta,\eta)}(\widetilde{f}_\eta) = \widetilde{\Theta}_{(\xi\cap\eta,\eta)}\Big(\widetilde{\Phi}_\eta(f_\eta)\Big) = \widetilde{0}_{\xi\cap\eta},$$

from where we conclude that

$$\widetilde{i}_{(\xi\cap\eta,\eta)}(\widetilde{\Phi}_\eta)\Big(\widetilde{\Theta}_{(\xi\cap\eta,\eta)}(f_\eta)\Big) = \widetilde{0}_{\xi\cap\eta},$$

or yet, since $\widetilde{\Theta}_{(\xi\cap\eta,\eta)}(f_\eta) = \Theta_{(\xi\cap\eta,\eta)}(f_\eta)$,

$$\widetilde{i}_{(\xi\cap\eta,\eta)}(\widetilde{\Phi}_\eta)\Big(\Theta_{(\xi\cap\eta,\eta)}(f_\eta)\Big) = \widetilde{0}_{\xi\cap\eta}$$

for every $\xi \in \overline{\xi}(\gamma)$ and $\eta \in \overline{\eta}(\gamma)$ such that $\xi \cap \eta \neq \varnothing$. Therefore,

$$\Theta_{(\xi\cap\eta,\eta)}(f_\eta) \in N\Big(\widetilde{i}_{(\xi\cap\eta,\eta)}(\widetilde{\Phi}_\eta)\Big).$$

Resorting now to the definition of extension of a bonding, Definition 3.22, we obtain that

$$\widetilde{i}_{(\xi\cap\eta,\eta)}(\widetilde{\Phi}_\eta) = \widetilde{(\xi \cap \eta)}\Big(i_{(\xi\cap\eta,\eta)}(\Phi_\eta)\Big)$$

and, hence,

$$\Theta_{(\xi\cap\eta,\eta)}(f_\eta) \in N\bigg(\widetilde{(\xi \cap \eta)}\Big(i_{(\xi\cap\eta,\eta)}(\Phi_\eta)\Big)\bigg)$$

for every $\xi \in \overline{\xi}(\gamma)$ and $\eta \in \overline{\eta}(\gamma)$ such that $\xi \cap \eta \neq \varnothing$.

But, from Proposition 1.16(b),

$$N\bigg(\widetilde{(\xi \cap \eta)}\Big(i_{(\xi\cap\eta,\eta)}(\Phi_\eta)\Big)\bigg) = N\Big(i_{(\xi\cap\eta,\eta)}(\Phi_\eta)\Big)$$

and, hence,

$$\Theta_{(\xi\cap\eta,\eta)}(f_\eta) \in N\Big(i_{(\xi\cap\eta,\eta)}(\Phi_\eta)\Big),$$

which means that

$$i_{(\xi\cap\eta,\eta)}(\Phi_\eta)\Big(\Theta_{(\xi\cap\eta,\eta)}(f_\eta)\Big) = 0_{\xi\cap\eta}$$

for every $\xi \in \overline{\xi}(\gamma)$ and $\eta \in \overline{\eta}(\gamma)$ such that $\xi \cap \eta \neq \varnothing$.

Let us take now any $\eta_0 \in \overline{\eta}(\gamma)$ and fix it. Clearly,

$$\overline{\nu}(\eta_0) := \Big\{\nu := \xi \cap \eta_0 \quad : \quad \xi \in \overline{\xi}(\gamma) \quad \text{and} \quad \xi \cap \eta_0 \neq \varnothing\Big\} \subseteq \Gamma(\eta_0)$$

is a cover of $\eta_0$, and from the last equation above results that

$$i_{(\nu,\eta_0)}(\Phi_{\eta_0})\Big(\Theta_{(\nu,\eta_0)}(f_{\eta_0})\Big) = 0_\nu$$

for every $\nu \in \overline{\nu}(\eta_0)$, from where one concludes, with Definition 3.13 in mind (more specifically, items (d) and (c) of this definition), that $f_{\eta_0}$ is in the domain of $\Phi_{\eta_0}$ ($f_{\eta_0} \in G(\eta_0)_{\Phi_{\eta_0}}$) and that

$$\Theta_{(\nu,\eta_0)}\Big(\Phi_{\eta_0}(f_{\eta_0})\Big) = 0_\nu$$

for every $\nu \in \overline{\nu}(\eta_0)$.

Since

$$\Theta_{(\nu,\eta_0)}(0_{\eta_0}) = 0_\nu \quad \text{for every} \quad \nu \in \overline{\nu}(\eta_0),$$

comes that

$$\Theta_{(\nu,\eta_0)}\Big(\Phi_{\eta_0}(f_{\eta_0})\Big) = \Theta_{(\nu,\eta_0)}(0_{\eta_0}) \quad \text{for every} \quad \nu \in \overline{\nu}(\eta_0)$$

and, hence, remembering that $\mathscr{G}(\eta_0)$ is coherent (since $\mathscr{G}(I)$ is) and taking Lemma 3.32 into account, we can conclude that

$$\Phi_{\eta_0}(f_{\eta_0}) = 0_{\eta_0}.$$

Once $\eta_0 \in \overline{\eta}(\gamma)$ was arbitrarily fixed, we have

$$\Phi_\eta(f_\eta) = 0_\eta \quad \text{for every} \quad \eta \in \overline{\eta}(\gamma)$$

and, hence, since $\widetilde{\Phi}_\eta$ is a prolongation of $\Phi_\eta$ to $\widetilde{G}(\eta)$,

$$\widetilde{\Phi}_\eta(f_\eta) = 0_\eta \quad \text{for every} \quad \eta \in \overline{\eta}(\gamma).$$

Now, remembering ourselves that

$$\widetilde{f}_\eta = \widetilde{\Theta}_{(\eta,\gamma)}(\widetilde{g} - \widetilde{h})$$

and that

$$\widetilde{f}_\eta = \widetilde{\Phi}_\eta(f_\eta),$$

we obtain,

$$\widetilde{\Theta}_{(\eta,\gamma)}(\widetilde{g} - \widetilde{h}) = 0_\eta \quad \text{for every} \quad \eta \in \overline{\eta}(\gamma).$$

Let now $e \in G(\gamma)$ and $\Psi \in H(\gamma)$ be such that

$$\widetilde{\Psi}(e) = \widetilde{g} - \widetilde{h}$$

($\Psi \in H(\gamma)$ and $e \in G(\gamma)$ exist since the $S$-group $\widetilde{\mathbb{G}}(\gamma) = (\widetilde{G}(\gamma), \widetilde{H}(\gamma))$ is a closed extension of the $S$-group $\mathbb{G}(\gamma) = (G(\gamma), H(\gamma))$). Then, we have

$$\widetilde{\Theta}_{(\eta,\gamma)}\Big(\widetilde{\Psi}(e)\Big) = 0_\eta \quad \text{for every} \quad \eta \in \overline{\eta}(\gamma)$$

and, hence,

$$\widetilde{i}_{(\eta,\gamma)}(\widetilde{\Psi})\Big(\widetilde{\Theta}_{(\eta,\gamma)}(e)\Big) = 0_\eta \quad \text{for every} \quad \eta \in \overline{\eta}(\gamma).$$

Therefore,

$$\widetilde{\Theta}_{(\eta,\gamma)}(e) = \Theta_{(\eta,\gamma)}(e) \in N\Big(\widetilde{i}_{(\eta,\gamma)}(\widetilde{\Psi})\Big),$$

and since (resorting again to Definition 3.22 of bonding extension, as well as Proposition 1.16(b))

$$N\Big(\widetilde{i}_{(\eta,\gamma)}(\widetilde{\Psi})\Big) = N\bigg(\widetilde{\eta}\Big(i_{(\eta,\gamma)}(\Psi)\Big)\bigg) = N\Big(i_{(\eta,\gamma)}(\Psi)\Big),$$

we have

$$\Theta_{(\eta,\gamma)}(e) \in N\Big(i_{(\eta,\gamma)}(\Psi)\Big) \quad \text{for every} \quad \eta \in \overline{\eta}(\gamma).$$

Hence,

$$i_{(\eta,\gamma)}(\Psi)\Big(\Theta_{(\eta,\gamma)}(e)\Big) = 0_\eta \quad \text{for every} \quad \eta \in \overline{\eta}(\gamma)$$

and, with this, appealing once more to items (d) and (c) of Definition 3.13, we get

$$\Theta_{(\eta,\gamma)}\Big(\Psi(e)\Big) = 0_\eta = \Theta_{(\eta,\gamma)}(0_\gamma) \quad \text{for every} \quad \eta \in \overline{\eta}(\gamma)$$

which, with a new appeal to Lemma 3.32, leads us to

$$\Psi(e) = \widetilde{\Psi}(e) = 0_\gamma.$$

But, $\widetilde{\Psi}(e) = \widetilde{g} - \widetilde{h}$ and, hence, we conclude that

$$\widetilde{g} = \widetilde{h},$$

which completes the proof of injectivity of the function $b_\gamma : \widetilde{G}(\gamma) \longrightarrow \overline{G}(\gamma)$. ∎

# The Prolongation $\overline{\Phi}$ of $\widetilde{\Phi}$ and the Semigroup $\overline{H}(\gamma)$

## 4.19 Motivation

It is convenient to remember at this point that when referring to the $S$-spaces

$$\mathscr{G}(I) = \bigg(\mathbb{G}\Big(\Gamma(I)\Big) = \Big\{\mathbb{G}(\gamma) = \Big(G(\gamma), H(\gamma)\Big)\Big\}_{\gamma\in\Gamma(I)}, i\Big(\Gamma^2(I)\Big), \Theta\Big(\Delta(I)\Big)\bigg),$$

$$\widetilde{\mathscr{G}}(I) = \bigg(\widetilde{\mathbb{G}}\Big(\Gamma(I)\Big) = \Big\{\widetilde{\mathbb{G}}(\gamma) = \Big(\widetilde{G}(\gamma), \widetilde{H}(\gamma)\Big)\Big\}_{\gamma\in\Gamma(I)}, \widetilde{i}\Big(\Gamma^2(I)\Big), \widetilde{\Theta}\Big(\Delta(I)\Big)\bigg)$$

and

$$\widehat{\mathscr{G}}(I) = \bigg(\widehat{\mathbb{G}}\Big(\Gamma(I)\Big) = \Big\{\widehat{\mathbb{G}}(\gamma) = \Big(\widehat{G}(\gamma), \widehat{H}(\gamma)\Big)\Big\}_{\gamma\in\Gamma(I)}, \widehat{i}\Big(\Gamma^2(I)\Big), \widehat{\Theta}\Big(\Delta(I)\Big)\bigg),$$

as established in Remarks 4.8, the following hypotheses are implicit: $\mathscr{G}(I)$ is abelian, surjective, coherent, and with identity; $\widetilde{\mathscr{G}}(I)$ is the strict and closed extension of $\mathscr{G}(I)$, unique unless isomorphism, described in Theorem 4.7 (1st TESS); $\widehat{\mathscr{G}}(I)$, a hypothetical locally closed and coherent extension of $\widetilde{\mathscr{G}}(I)$. Furthermore, yet in the Remarks 4.8, we established our goals: to build an extension of $\widetilde{\mathscr{G}}(I)$,

$$\overline{\mathscr{G}}(I) = \left(\overline{\mathbb{G}}\Big(\Gamma(I)\Big) = \left\{\overline{\mathbb{G}}(\gamma) = \Big(\overline{G}(\gamma), \overline{H}(\gamma)\Big)\right\}_{\gamma\in\Gamma(I)}, \bar{i}\Big(\Gamma^2(I)\Big), \overline{\Theta}\Big(\Delta(I)\Big)\right),$$

isomorphic to $\widehat{\mathscr{G}}(I)$. Regarding this goal, a review of what we got until now is presented ahead:

**(a)** We already have a family of abelian groups (Proposition 4.16), $\{\overline{G}(\gamma)\}_{\gamma\in\Gamma(I)}$, such that the group $\overline{G}(\gamma)$, for each $\gamma \in \Gamma(I)$, admits the group $\widetilde{G}(\gamma)$ as a subgroup, that is, the function $b_\gamma : \widetilde{G}(\gamma) \longrightarrow \overline{G}(\gamma)$, defined in Proposition 4.18, is an injective homomorphism and, hence, $b_\gamma(\widetilde{G}(\gamma)) \subseteq \overline{G}(\gamma)$ is a subgroup of $\overline{G}(\gamma)$ isomorphic to the group $\widetilde{G}(\gamma)$;

**(b)** Under the hypothesis of existence of $\widehat{\mathscr{G}}(I)$, the group $\overline{G}(\gamma)$, for each $\gamma \in \Gamma(I)$, is isomorphic to the group $\widehat{G}(\gamma)$ and the function $A_\gamma : \overline{G}(\gamma) \longrightarrow \widehat{G}(\gamma)$, defined in the Proposition 4.12, is an isomorphism that keeps fixed the elements of $\widetilde{G}(\gamma)$ in the sense of the function $A_\gamma \circ b_\gamma : \widetilde{G}(\gamma) \longrightarrow \widehat{G}(\gamma)$ be such that

$$(A_\gamma \circ b_\gamma)(\tilde{g}) = A_\gamma\Big(b_\gamma(\tilde{g})\Big) = \tilde{g}$$

for every $\tilde{g} \in \widetilde{G}(\gamma)$.

The review above makes clear our next task to obtain the extension $\overline{\mathscr{G}}(I)$ for $\widetilde{\mathscr{G}}(I)$, namely: for each $\gamma \in \Gamma(I)$, to equip the abelian group $\overline{G}(\gamma)$ with a semigroup of endomorphisms on $\overline{G}(\gamma)$ that are, exactly and precisely, prolongations of the homomorphisms $\widetilde{\Phi} \in \widetilde{H}(\gamma)$ to $\overline{G}(\gamma)$ (recap Definition 1.6). In other terms, we must obtain, for each $\gamma \in \Gamma(I)$, a semigroup $\overline{H}(\gamma)$ of endomorphisms on $\overline{G}(\gamma)$ in such a way that the $S$-group $(\overline{G}(\gamma), \overline{H}(\gamma))$ becomes an extension of the $S$-group $(\widetilde{G}(\gamma), \widetilde{H}(\gamma))$, beyond being isomorphic to the $S$-group $(\widehat{G}(\gamma), \widehat{H}(\gamma))$ under the hypothesis of existence of $\widehat{\mathscr{G}}(I)$.

Hence, for each $\widetilde{\Phi} : \widetilde{G}(\gamma) \longrightarrow \widetilde{G}(\gamma)$ in $\widetilde{H}(\gamma)$, we must define $\overline{\Phi} : \overline{G}(\gamma) \longrightarrow \overline{G}(\gamma)$ in such a way that $\overline{\Phi}$ becomes an endomorphism on $\overline{G}(\gamma)$ and

$$\overline{\Phi}(\tilde{g}) = \widetilde{\Phi}(\tilde{g}) \quad \text{for every} \quad \tilde{g} \in \widetilde{G}(\gamma),$$

or, more precisely, without any abuse of language (recap Definition 1.6, in particular its part (b), where the abuse of language above referred is explained),

$$\overline{\Phi}\Big(b_\gamma(\tilde{g})\Big) = b_\gamma\Big(\widetilde{\Phi}(\tilde{g})\Big) \quad \text{for every} \quad \tilde{g} \in \widetilde{G}(\gamma). \tag{4.19-1}$$

On the other hand, on the hypothesis of existence of $\widehat{\mathscr{G}}(I)$, $\widehat{\Phi} : \widehat{G}(\gamma) \longrightarrow \widehat{G}(\gamma)$ is a prolongation of $\widetilde{\Phi}$ to $\widehat{G}(\gamma)$, that is, $\widehat{\Phi} \in \widehat{H}(\gamma)$ is given by

$$\widehat{\Phi} = \widehat{\gamma}(\widetilde{\Phi}),$$

and, hence,

$$\widehat{\Phi}(\widetilde{g}) = \widetilde{\Phi}(\widetilde{g}) \quad \text{for every} \quad \widetilde{g} \in \widetilde{G}(\gamma), \tag{4.19-2}$$

From (4.19-1) and (4.19-2) follows

$$\overline{\Phi}\Big(b_\gamma(\widetilde{g})\Big) = b_\gamma\Big(\widehat{\Phi}(\widetilde{g})\Big) \quad \text{for every} \quad \widetilde{g} \in \widetilde{G}(\gamma). \tag{4.19-3}$$

But,

$$A_\gamma\Big(b_\gamma(\widetilde{g})\Big) = \widetilde{g}$$

or, equivalently,

$$b_\gamma(\widetilde{g}) = A_\gamma^{-1}(\widetilde{g}) \quad \text{for every} \quad \widetilde{g} \in \widetilde{G}(\gamma),$$

and since $\widehat{\Phi}(\widetilde{g}) = \widetilde{\Phi}(\widetilde{g}) \in \widetilde{G}(\gamma)$ for every $\widetilde{g} \in \widetilde{G}(\gamma)$, then

$$b_\gamma\Big(\widehat{\Phi}(\widetilde{g})\Big) = A_\gamma^{-1}\Big(\widehat{\Phi}(\widetilde{g})\Big) \quad \text{for every} \quad \widetilde{g} \in \widetilde{G}(\gamma).$$

Taking this results into (4.19-3), one obtains:

$$\overline{\Phi}\Big(A_\gamma^{-1}(\widetilde{g})\Big) = A_\gamma^{-1}\Big(\widehat{\Phi}(\widetilde{g})\Big)$$

or yet

$$\overline{\Phi}\Big(A_\gamma^{-1}(\widetilde{g})\Big) = A_\gamma^{-1}\left(\widehat{\Phi}\left(A_\gamma\Big(A_\gamma^{-1}(\widetilde{g})\Big)\right)\right) \quad \text{for every} \quad \widetilde{g} \in \widetilde{G}(\gamma).$$

Making now

$$\overline{g} \coloneqq A_\gamma^{-1}(\widetilde{g}) = b_\gamma(\widetilde{g})$$

one gets that

$$\overline{\Phi}(\overline{g}) = A_\gamma^{-1}\left(\widehat{\Phi}\Big(A_\gamma(\overline{g})\Big)\right)$$

for every $\overline{g} \in b_\gamma(\widetilde{G}(\gamma)) = \{b_\gamma(\widetilde{g}) : \widetilde{g} \in \widetilde{G}(\gamma)\}$.

We can summarize the considerations above as follows: if $\widehat{\mathscr{G}}(I)$ does exists, then, an endomorphism on $\overline{G}(\gamma)$, $\overline{\Phi} : \overline{G}(\gamma) \longrightarrow \overline{G}(\gamma)$, is a prolongation of $\widetilde{\Phi} \in \widetilde{H}(\gamma)$ (to $\overline{G}(\gamma)$) if and only if,

$$\overline{\Phi}(\overline{g}) = A_\gamma^{-1}\left(\widehat{\Phi}\Big(A_\gamma(\overline{g})\Big)\right) \quad \text{for every} \quad \overline{g} \in b_\gamma\Big(\widetilde{G}(\gamma)\Big), \tag{4.19-4}$$

where $\widehat{\Phi} \in \widehat{H}(\gamma)$ is the prolongation of $\widetilde{\Phi} \in \widetilde{H}(\gamma)$ to $\widehat{G}(\gamma)$, that is, $\widehat{\Phi} = \widehat{\gamma}(\widetilde{\Phi})$.

Hence, under the existence of $\widehat{\mathscr{G}}(I)$, defining $\overline{\Phi} : \overline{G}(\gamma) \longrightarrow \overline{G}(\gamma)$ by

$$\overline{\Phi}(\overline{g}) \coloneqq A_\gamma^{-1}\left(\widehat{\Phi}\Big(A_\gamma(\overline{g})\Big)\right) \quad \text{for every} \quad \overline{g} \in \overline{G}(\gamma), \tag{4.19-5}$$

where $\widehat{\Phi} = \widehat{\gamma}(\widetilde{\Phi})$, it results immediately that $\overline{\Phi} : \overline{G}(\gamma) \longrightarrow \overline{G}(\gamma)$ is an endomorphism on $\overline{G}(\gamma)$ (since $A_\gamma : \overline{G}(\gamma) \longrightarrow \widehat{G}(\gamma)$ is an isomorphism between the groups $\overline{G}(\gamma)$ and $\widehat{G}(\gamma)$, and $\widehat{\Phi} : \widehat{G}(\gamma) \longrightarrow \widehat{G}(\gamma)$ is an endomorphism on $\widehat{G}(\gamma)$) that satisfies the condition (4.19-4) above and, therefore, it is a prolongation of $\widetilde{\Phi} \in \widetilde{H}(\gamma)$ to $\overline{G}(\gamma)$.

**Remark.** To reach the conclusion above, we made use of the following hypothesis: the group $\widehat{G}(\gamma)$ not only admits the group $\widetilde{G}(\gamma)$ as its subgroup (in the sense of the Definition 1.6(a)), but "trivially" does so, that is, $\widetilde{G}(\gamma) \subseteq \widehat{G}(\gamma)$ is, in fact, a subgroup of $\widehat{G}(\gamma)$ and not only a group that is isomorphic to a subgroup of $\widehat{G}(\gamma)$. This hypothesis is explicitly put on (4.19-2) as well as, right after the expression (4.19-3), when we take $A_\gamma(b_\gamma(\widetilde{g})) = \widetilde{g}$ for every $\widetilde{g} \in \widetilde{G}(\gamma)$. However, this hypothesis does not imply any loss of generality in the following sense: if $\widetilde{G}(\gamma)$ is not a subset of $\widehat{G}(\gamma)$ and, therefore, $A_\gamma(b_\gamma(\widetilde{g})) \neq \widetilde{g}$, even so $A_\gamma \circ b_\gamma : \widetilde{G}(\gamma) \longrightarrow \widehat{G}(\gamma)$ is an injective homomorphism and thus

$$A_\gamma\left(b_\gamma\Big(\widetilde{G}(\gamma)\Big)\right) = \Big\{A_\gamma\Big(b_\gamma(\widetilde{g})\Big) : \widetilde{g} \in \widetilde{G}(\gamma)\Big\} \subseteq \widehat{G}(\gamma)$$

is a subgroup of $\widehat{G}(\gamma)$ isomorphic to $\widetilde{G}(\gamma)$, that is, $\widehat{G}(\gamma)$ admits (although not trivially any more) the group $\widetilde{G}(\gamma)$ as one of its subgroups; with this, and by this, as we will prove ahead, one obtains the same conclusion reached above on the hypothesis of $\widetilde{G}(\gamma) \subseteq \widehat{G}(\gamma)$, that is, the same expression (4.19-4) is obtained for $\overline{\Phi}(\overline{g})$. In fact, in the case where $\widehat{G}(\gamma)$ does not contain $\widetilde{G}(\gamma)$ but contains $A_\gamma(b_\gamma(\widetilde{G}(\gamma)))$, which is a subgroup of $\widehat{G}(\gamma)$ isomorphic to $\widetilde{G}(\gamma)$, the expression (4.19-2) involves the previously referred abuse of language and must be understood as an abbreviated form of the following equation:

$$\widehat{\Phi}\Big((A_\gamma \circ b_\gamma)(\widetilde{g})\Big) = (A_\gamma \circ b_\gamma)\Big(\widetilde{\Phi}(\widetilde{g})\Big) \quad \text{for every} \quad \widetilde{g} \in \widetilde{G}(\gamma)$$

from where one obtains that

$$A_\gamma^{-1}\left(\widehat{\Phi}\left(A_\gamma\Big(b_\gamma(\widetilde{g})\Big)\right)\right) = b_\gamma\Big(\widetilde{\Phi}(\widetilde{g})\Big)$$

that when taken into 4.19-1 gives us

$$\overline{\Phi}\Big(b_\gamma(\widetilde{g})\Big) = A_\gamma^{-1}\left(\widehat{\Phi}\left(A_\gamma\Big(b_\gamma(\widetilde{g})\Big)\right)\right) \quad \text{for every} \quad \widetilde{g} \in \widetilde{G}(\gamma).$$

Taking now

$$\overline{g} := b_\gamma(\widetilde{g})$$

one obtains

$$\overline{\Phi}(\overline{g}) = A_\gamma^{-1}\left(\widehat{\Phi}\Big(A_\gamma(\overline{g})\Big)\right),$$

for every $\overline{g} \in b_\gamma(\widetilde{G}(\gamma)) = \{b_\gamma(\widetilde{g}) : \widetilde{g} \in \widetilde{G}(\gamma)\}$, which is the same expression (4.19-4) previously obtained.

Let us resume, after this remark, to the definition of $\overline{\Phi} : \overline{G}(\gamma) \longrightarrow \overline{G}(\gamma)$, given by (4.19-5), backed by the hypothesis of existence of $\widehat{\mathscr{G}}(I)$, namely:

$$\begin{aligned} \overline{\Phi} : \overline{G}(\gamma) &\longrightarrow \overline{G}(\gamma) \\ \overline{g} &\longmapsto \overline{\Phi}(\overline{g}) := A_\gamma^{-1}\left(\widehat{\Phi}\Big(A_\gamma(\overline{g})\Big)\right) \end{aligned} \tag{4.19-5}$$

where $\widehat{\Phi} \in \widehat{H}(\gamma)$ is the prolongation of $\widetilde{\Phi} \in \widetilde{H}(\gamma)$ to $\widehat{G}(\gamma)$, that is, $\widehat{\Phi} = \widehat{\gamma}(\widetilde{\Phi})$. As we had concluded, $\overline{\Phi}$ is an endomorphism on $\overline{G}(\gamma)$ which is, also, a prolongation of $\widetilde{\Phi} \in \widetilde{H}(\gamma)$ to $\overline{G}(\gamma)$.

Let now

$$\overline{g} = \left[\widetilde{g}\Big(\overline{\xi}(\gamma)\Big) = \Big\{\widetilde{g}_\xi\Big\}_{\xi \in \overline{\xi}(\gamma)}\right] \in \overline{G}(\gamma)$$

be arbitrarily fixed and let us calculate $\overline{\Phi}(\overline{g})$ such as defined above by (4.19-5). We have, by the definition of $A_\gamma : \overline{G}(\gamma) \longrightarrow \widehat{G}(\gamma)$, that $A_\gamma(\overline{g})$ is the only element of $\widehat{G}(\gamma)$ such that

$$\widehat{\Theta}_{(\xi,\gamma)}\Big(A_\gamma(\overline{g})\Big) = \widetilde{g}_\xi \quad \text{for every} \quad \xi \in \overline{\xi}(\gamma).$$

Thus, for any $\xi \in \overline{\xi}(\gamma)$, we have:

$$\widehat{\Theta}_{(\xi,\gamma)}\left(\widehat{\Phi}\Big(A_\gamma(\overline{g})\Big)\right) = \widehat{i}_{(\xi,\gamma)}(\widehat{\Phi})\left(\widehat{\Theta}_{(\xi,\gamma)}\Big(A_\gamma(\overline{g})\Big)\right) = \widehat{i}_{(\xi,\gamma)}(\widehat{\Phi})(\widetilde{g}_\xi)$$

and since

$$\widehat{i}_{(\xi,\gamma)}(\widehat{\Phi}) = \widehat{\xi}\left(\widetilde{i}_{(\xi,\gamma)}\Big((\widehat{\gamma})^{-1}(\widehat{\Phi})\Big)\right) = \widehat{\xi}\Big(\widetilde{i}_{(\xi,\gamma)}(\widetilde{\Phi})\Big)$$

then

$$\widehat{\Theta}_{(\xi,\gamma)}\left(\widehat{\Phi}\Big(A_\gamma(\overline{g})\Big)\right) = \left(\widehat{\xi}\Big(\widetilde{i}_{(\xi,\gamma)}(\widetilde{\Phi})\Big)\right)(\widetilde{g}_\xi)$$

But $\widehat{\xi}(\widetilde{i}_{(\xi,\gamma)}(\widetilde{\Phi}))$ is the prolongation to $\widehat{G}(\xi)$ of $\widetilde{i}_{(\xi,\gamma)}(\widetilde{\Phi})$ and, hence,

$$\left(\widehat{\xi}\Big(\widetilde{i}_{(\xi,\gamma)}(\widetilde{\Phi})\Big)\right)(\widetilde{g}_\xi) = \widetilde{i}_{(\xi,\gamma)}(\widetilde{\Phi})(\widetilde{g}_\xi).$$

Therefore,

$$\widehat{\Theta}_{(\xi,\gamma)}\left(\widehat{\Phi}\Big(A_\gamma(\overline{g})\Big)\right) = \widetilde{i}_{(\xi,\gamma)}(\widetilde{\Phi})(\widetilde{g}_\xi) \quad \text{for every} \quad \xi \in \overline{\xi}(\gamma).$$

From the expression above, taking into account the definition of $A_\gamma : \overline{G}(\gamma) \longrightarrow \widehat{G}(\gamma)$, it results that

$$A_\gamma^{-1}\left(\widehat{\Phi}\Big(A_\gamma(\overline{g})\Big)\right) = \Big[\Big\{\widetilde{i}_{(\xi,\gamma)}(\widetilde{\Phi})(\widetilde{g}_\xi)\Big\}_{\xi\in\overline{\xi}(\gamma)}\Big]$$

and, hence, by the definition of $\overline{\Phi}$ given at (4.19-5), we have:

$$\overline{\Phi}\left(\overline{g} = \Big[\widetilde{g}\Big(\overline{\xi}(\gamma)\Big)\Big]\right) = \Big[\Big\{\widetilde{i}_{(\xi,\gamma)}(\widetilde{\Phi})(\widetilde{g}_\xi)\Big\}_{\xi\in\overline{\xi}(\gamma)}\Big]. \tag{4.19-6}$$

Observe now that the expression obtained for $\overline{\Phi}(\overline{g})$, equation (4.19-6) above, does not make any mention to $\widehat{\mathscr{G}}(I)$; it only involves objects of $\widetilde{\mathscr{G}}(I)$, namely, the endomorphisms $\widetilde{\Phi} \in \widetilde{H}(\gamma)$ and $\widetilde{i}_{(\xi,\gamma)}(\widetilde{\Phi}) \in \widetilde{H}(\xi)$ for each $\xi \in \overline{\xi}(\gamma)$, and $\overline{g} = [\widetilde{g}(\overline{\xi}(\gamma))] \in \overline{G}(\gamma)$ (worth remember that $\overline{G}(\gamma)$ was "built" with only what we have in $\widetilde{\mathscr{G}}(I)$).

This observation raises the following question: The expression (4.19-6) defines, by itself, that is, independently of the existence of $\widehat{\mathscr{G}}(I)$, an endomorphism on $\overline{G}(\gamma)$ which is a prolongation of $\widetilde{\Phi} \in \widetilde{H}(\gamma)$ to $\overline{G}(\gamma)$? Clearly, if so, and $\widehat{\mathscr{G}}(I)$ (also) does exist, we would have: $\overline{\Phi} = A_\gamma^{-1}\widehat{\Phi}A_\gamma$.

Our next proposition intends to answer the question above.

## 4.20 Proposition

*Let*

$$\widetilde{\mathscr{G}}(I) = \left(\widetilde{\mathbb{G}}\Big(\Gamma(I)\Big) = \Big\{\widetilde{\mathbb{G}}(\gamma) = \Big(\widetilde{G}(\gamma), \widetilde{H}(\gamma)\Big)\Big\}_{\gamma\in\Gamma(I)}, \widetilde{i}\Big(\Gamma^2(I)\Big), \widetilde{\Theta}\Big(\Delta(I)\Big)\right)$$

*be as described in 4.6, $\gamma \in \Gamma(I)$ be arbitrarily fixed and $\overline{G}(\gamma)$ be the abelian group defined in Proposition 4.16. Then, for each $\widetilde{\Phi} \in \widetilde{H}(\gamma)$, the function*

$$\overline{\Phi} : \overline{G}(\gamma) \longrightarrow \overline{G}(\gamma)$$

*where*

$$\overline{\Phi}\left(\Big[\widetilde{g}\Big(\overline{\xi}(\gamma)\Big) = \Big\{\widetilde{g}_\xi\Big\}_{\xi\in\overline{\xi}(\gamma)}\Big]\right) := \Big[\Big\{\widetilde{i}_{(\xi,\gamma)}(\widetilde{\Phi})(\widetilde{g}_\xi)\Big\}_{\xi\in\overline{\xi}(\gamma)}\Big],$$

*is well-defined, it is an endomorphism on $\overline{G}(\gamma)$ and a prolongation of $\widetilde{\Phi}$ to $\overline{G}(\gamma)$. Furthermore, if $\widehat{\mathscr{G}}(I)$ does exist, then*

$$\overline{\Phi} = A_\gamma^{-1}\ \widehat{\Phi}\ A_\gamma$$

*where $\widehat{\Phi} = \widehat{\gamma}(\widetilde{\Phi})$ (the prolongation of $\widetilde{\Phi} \in \widetilde{H}(\gamma)$ to $\widehat{G}(\gamma)$) and $A_\gamma : \overline{G}(\gamma) \longrightarrow \widehat{G}(\gamma)$ is the function defined in Proposition 4.12.*

*Proof.* First, we must prove that

$$\left[\left\{\widetilde{i}_{(\xi,\gamma)}(\widetilde{\Phi})(\widetilde{g}_\xi)\right\}_{\xi\in\overline{\xi}(\gamma)}\right] \in \overline{G}(\gamma),$$

i.e., that the family

$$\left\{\widetilde{i}_{(\xi,\gamma)}(\widetilde{\Phi})(\widetilde{g}_\xi)\right\}_{\xi\in\overline{\xi}(\gamma)}$$

is coherent in $\widetilde{\mathscr{G}}(\gamma)$ (that is, it belongs to $\widetilde{F}(\gamma)$). Hence, for $\xi, \xi' \in \overline{\xi}(\gamma)$ such that $\xi \cap \xi' \neq \varnothing$, we calculate:

$$\begin{aligned}\widetilde{\Theta}_{(\xi\cap\xi',\xi)}\left(\widetilde{i}_{(\xi,\gamma)}(\widetilde{\Phi})(\widetilde{g}_\xi)\right) &= \widetilde{i}_{(\xi\cap\xi',\xi)}\left(\widetilde{i}_{(\xi,\gamma)}(\widetilde{\Phi})\right)\left(\widetilde{\Theta}_{(\xi\cap\xi',\xi)}(\widetilde{g}_\xi)\right) = \\ &= \widetilde{i}_{(\xi\cap\xi',\gamma)}(\widetilde{\Phi})\left(\widetilde{\Theta}_{(\xi\cap\xi',\xi)}(\widetilde{g}_\xi)\right)\end{aligned}$$

and

$$\begin{aligned}\widetilde{\Theta}_{(\xi\cap\xi',\xi')}\left(\widetilde{i}_{(\xi',\gamma)}(\widetilde{\Phi})(\widetilde{g}_{\xi'})\right) &= \widetilde{i}_{(\xi\cap\xi',\xi')}\left(\widetilde{i}_{(\xi',\gamma)}(\widetilde{\Phi})\right)\left(\widetilde{\Theta}_{(\xi\cap\xi',\xi')}(\widetilde{g}_{\xi'})\right) = \\ &= \widetilde{i}_{(\xi\cap\xi',\gamma)}(\widetilde{\Phi})\left(\widetilde{\Theta}_{(\xi\cap\xi',\xi')}(\widetilde{g}_{\xi'})\right).\end{aligned}$$

Once $\widetilde{\Theta}_{(\xi\cap\xi',\xi)}(\widetilde{g}_\xi) = \widetilde{\Theta}_{(\xi\cap\xi',\xi')}(\widetilde{g}_{\xi'})$ since the family $\widetilde{g}(\overline{\xi}(\gamma)) = \{\widetilde{g}_\xi\}_{\xi\in\overline{\xi}(\gamma)}$ is coherent in $\widetilde{\mathscr{G}}(\gamma)$, we conclude that

$$\widetilde{\Theta}_{(\xi\cap\xi',\xi)}\left(\widetilde{i}_{(\xi,\gamma)}(\widetilde{\Phi})(\widetilde{g}_\xi)\right) = \widetilde{\Theta}_{(\xi\cap\xi',\xi')}\left(\widetilde{i}_{(\xi',\gamma)}(\widetilde{\Phi})(\widetilde{g}_{\xi'})\right),$$

that is, $\{\widetilde{i}_{(\xi,\gamma)}(\widetilde{\Phi})(\widetilde{g}_\xi)\}_{\xi\in\overline{\xi}(\gamma)}$ is a coherent family in $\widetilde{\mathscr{G}}(\gamma)$.

Now it is necessary to prove that the function $\overline{\Phi}$ is well-defined, that is: if

$$\widetilde{g}\left(\overline{\xi}(\gamma)\right) = \left\{\widetilde{g}_\xi\right\}_{\xi\in\overline{\xi}(\gamma)} \approx \left\{\widetilde{h}_\eta\right\}_{\eta\in\overline{\eta}(\gamma)} = \widetilde{h}\left(\overline{\eta}(\gamma)\right),$$

then

$$\left\{\widetilde{i}_{(\xi,\gamma)}(\widetilde{\Phi})(\widetilde{g}_\xi)\right\}_{\xi\in\overline{\xi}(\gamma)} \approx \left\{\widetilde{i}_{(\eta,\gamma)}(\widetilde{\Phi})(\widetilde{h}_\eta)\right\}_{\eta\in\overline{\eta}(\gamma)}.$$

Hence, let us suppose that $\widetilde{g}(\overline{\xi}(\gamma)) \approx \widetilde{h}(\overline{\eta}(\gamma))$, that is,

$$\widetilde{\Theta}_{(\xi\cap\eta,\xi)}(\widetilde{g}_\xi) = \widetilde{\Theta}_{(\xi\cap\eta,\eta)}(\widetilde{h}_\eta)$$

for every $\xi \in \overline{\xi}(\gamma)$ and $\eta \in \overline{\eta}(\gamma)$ such that $\xi \cap \eta \neq \varnothing$.

We then have

$$\widetilde{i}_{(\xi\cap\eta,\xi)}\left(\widetilde{i}_{(\xi,\gamma)}(\widetilde{\Phi})\right)\left(\widetilde{\Theta}_{(\xi\cap\eta,\xi)}(\widetilde{g}_\xi)\right) = \widetilde{i}_{(\xi\cap\eta,\eta)}\left(\widetilde{i}_{(\eta,\gamma)}(\widetilde{\Phi})\right)\left(\widetilde{\Theta}_{(\xi\cap\eta,\eta)}(\widetilde{h}_\eta)\right)$$

or yet,

$$\widetilde{\Theta}_{(\xi\cap\eta,\xi)}\left(\widetilde{i}_{(\xi,\gamma)}(\widetilde{\Phi})(\widetilde{g}_\xi)\right) = \widetilde{\Theta}_{(\xi\cap\eta,\eta)}\left(\widetilde{i}_{(\eta,\gamma)}(\widetilde{\Phi})(\widetilde{h}_\eta)\right)$$

for every $\xi \in \overline{\xi}(\gamma)$ and $\eta \in \overline{\eta}(\gamma)$ such that $\xi \cap \eta \neq \varnothing$, that is,

$$\left\{\widetilde{i}_{(\xi,\gamma)}(\widetilde{\Phi})(\widetilde{g}_\xi)\right\}_{\xi\in\overline{\xi}(\gamma)} \approx \left\{\widetilde{i}_{(\eta,\gamma)}(\widetilde{\Phi})(\widetilde{h}_\eta)\right\}_{\eta\in\overline{\eta}(\gamma)}.$$

Let us now prove that $\overline{\Phi}$ is an endomorphism on $\overline{G}(\gamma)$. In fact, since for

$$\overline{g} = \left[\left\{\widetilde{g}_\xi\right\}_{\xi\in\overline{\xi}(\gamma)}\right] \quad \text{and} \quad \overline{h} = \left[\left\{\widetilde{h}_\eta\right\}_{\eta\in\overline{\eta}(\gamma)}\right]$$

arbitrarily chosen in $\overline{G}(\gamma)$, we have, by the definition of addition of the group $\overline{G}(\gamma)$ (see Proposition 4.15), that

$$\overline{g} + \overline{h} = \left[\left\{\widetilde{\Theta}_{(\xi\cap\eta,\xi)}(\widetilde{g}_\xi) + \widetilde{\Theta}_{(\xi\cap\eta,\eta)}(\widetilde{h}_\eta)\right\}_{\xi\cap\eta\in\overline{\chi}(\gamma)}\right]$$

where $\overline{\chi}(\gamma) \subseteq \Gamma(\gamma)$ is the cover of $\gamma$ given by

$$\overline{\chi}(\gamma) := \left\{\chi := \xi \cap \eta \quad : \quad \xi \in \overline{\xi}(\gamma), \quad \eta \in \overline{\eta}(\gamma) \quad \text{and} \quad \xi \cap \eta \neq \varnothing\right\}.$$

Therefore, by the definition of $\overline{\Phi}$ follows that

$$\overline{\Phi}(\overline{g} + \overline{h}) = \left[\left\{\widetilde{i}_{(\xi\cap\eta,\gamma)}(\widetilde{\Phi})\Big(\widetilde{\Theta}_{(\xi\cap\eta,\xi)}(\widetilde{g}_\xi) + \widetilde{\Theta}_{(\xi\cap\eta,\eta)}(\widetilde{h}_\eta)\Big)\right\}_{\xi\cap\eta\in\overline{\chi}(\gamma)}\right].$$

Now,

$$\begin{aligned}
&\widetilde{i}_{(\xi\cap\eta,\gamma)}(\widetilde{\Phi})\Big(\widetilde{\Theta}_{(\xi\cap\eta,\xi)}(\widetilde{g}_\xi) + \widetilde{\Theta}_{(\xi\cap\eta,\eta)}(\widetilde{h}_\eta)\Big) = \\
&= \widetilde{i}_{(\xi\cap\eta,\gamma)}(\widetilde{\Phi})\Big(\widetilde{\Theta}_{(\xi\cap\eta,\xi)}(\widetilde{g}_\xi)\Big) + \widetilde{i}_{(\xi\cap\eta,\gamma)}(\widetilde{\Phi})\Big(\widetilde{\Theta}_{(\xi\cap\eta,\eta)}(\widetilde{h}_\eta)\Big) = \\
&= \widetilde{i}_{(\xi\cap\eta,\xi)}\Big(\widetilde{i}_{(\xi,\gamma)}(\widetilde{\Phi})\Big)\Big(\widetilde{\Theta}_{(\xi\cap\eta,\xi)}(\widetilde{g}_\xi)\Big) + \widetilde{i}_{(\xi\cap\eta,\eta)}\Big(\widetilde{i}_{(\eta,\gamma)}(\widetilde{\Phi})\Big)\Big(\widetilde{\Theta}_{(\xi\cap\eta,\eta)}(\widetilde{h}_\eta)\Big) = \\
&= \widetilde{\Theta}_{(\xi\cap\eta,\xi)}\Big(\widetilde{i}_{(\xi,\gamma)}(\widetilde{\Phi})(\widetilde{g}_\xi)\Big) + \widetilde{\Theta}_{(\xi\cap\eta,\eta)}\Big(\widetilde{i}_{(\eta,\gamma)}(\widetilde{\Phi})(\widetilde{h}_\eta)\Big).
\end{aligned}$$

Thus,

$$\overline{\Phi}(\overline{g} + \overline{h}) = \left[\left\{\widetilde{\Theta}_{(\xi\cap\eta,\xi)}\Big(\widetilde{i}_{(\xi,\gamma)}(\widetilde{\Phi})(\widetilde{g}_\xi)\Big) + \widetilde{\Theta}_{(\xi\cap\eta,\eta)}\Big(\widetilde{i}_{(\eta,\gamma)}(\widetilde{\Phi})(\widetilde{h}_\eta)\Big)\right\}_{\xi\cap\eta\in\overline{\chi}(\gamma)}\right].$$

But, from the definition of addition in $\overline{G}(\gamma)$, we have

$$\begin{aligned}
&\left[\left\{\widetilde{i}_{(\xi,\gamma)}(\widetilde{\Phi})(\widetilde{g}_\xi)\right\}_{\xi\in\overline{\xi}(\gamma)}\right] + \left[\left\{\widetilde{i}_{(\eta,\gamma)}(\widetilde{\Phi})(\widetilde{h}_\eta)\right\}_{\eta\in\overline{\eta}(\gamma)}\right] = \\
&= \left[\left\{\widetilde{\Theta}_{(\xi\cap\eta,\xi)}\Big(\widetilde{i}_{(\xi,\gamma)}(\widetilde{\Phi})(\widetilde{g}_\xi)\Big) + \widetilde{\Theta}_{(\xi\cap\eta,\eta)}\Big(\widetilde{i}_{(\eta,\gamma)}(\widetilde{\Phi})(\widetilde{h}_\eta)\Big)\right\}_{\xi\cap\eta\in\overline{\chi}(\gamma)}\right],
\end{aligned}$$

and, hence,

$$\begin{aligned}
\overline{\Phi}(\overline{g} + \overline{h}) &= \left[\left\{\widetilde{i}_{(\xi,\gamma)}(\widetilde{\Phi})(\widetilde{g}_\xi)\right\}_{\xi\in\overline{\xi}(\gamma)}\right] + \left[\left\{\widetilde{i}_{(\eta,\gamma)}(\widetilde{\Phi})(\widetilde{h}_\eta)\right\}_{\eta\in\overline{\eta}(\gamma)}\right] = \\
&= \overline{\Phi}(\overline{g}) + \overline{\Phi}(\overline{h}).
\end{aligned}$$

Finally, we must prove that $\overline{\Phi}$ is a prolongation to $\overline{G}(\gamma)$ of $\widetilde{\Phi}$, that is,

$$\overline{\Phi}\Big(b_\gamma(\widetilde{g})\Big) = b_\gamma\Big(\widetilde{\Phi}(\widetilde{g})\Big)$$

for every $\widetilde{g} \in \widetilde{G}(\gamma)$, being $b_\gamma : \widetilde{G}(\gamma) \longrightarrow \overline{G}(\gamma)$ the injective homomorphism defined in Proposition 4.18 by

$$b_\gamma(\widetilde{g}) = \Big[\Big\{\widetilde{\Theta}_{(\xi,\gamma)}(\widetilde{g})\Big\}_{\xi\in\overline{\xi}(\gamma)}\Big],$$

where $\overline{\xi}(\gamma) \subseteq \Gamma(\gamma)$ is an arbitrary cover of $\gamma$.

We have, for $\widetilde{g} \in \widetilde{G}(\gamma)$ arbitrarily chosen, that

$$\overline{\Phi}\Big(b_\gamma(\widetilde{g})\Big) = \Big[\Big\{\widetilde{i}_{(\xi,\gamma)}(\widetilde{\Phi})\Big(\widetilde{\Theta}_{(\xi,\gamma)}(\widetilde{g})\Big)\Big\}_{\xi\in\overline{\xi}(\gamma)}\Big].$$

Now,

$$\widetilde{i}_{(\xi,\gamma)}(\widetilde{\Phi})\Big(\widetilde{\Theta}_{(\xi,\gamma)}(\widetilde{g})\Big) = \widetilde{\Theta}_{(\xi,\gamma)}\Big(\widetilde{\Phi}(\widetilde{g})\Big),$$

and, hence,

$$\overline{\Phi}\Big(b_\gamma(\widetilde{g})\Big) = \Big[\Big\{\widetilde{\Theta}_{(\xi,\gamma)}\Big(\widetilde{\Phi}(\widetilde{g})\Big)\Big\}_{\xi\in\overline{\xi}(\gamma)}\Big] = b_\gamma\Big(\widetilde{\Phi}(\widetilde{g})\Big). \qquad \blacksquare$$

## 4.21 The Semigroup $\overline{H}(\gamma)$

Associating to each $\widetilde{\Phi} \in \widetilde{H}(\gamma)$ the endomorphism $\overline{\Phi}$ on $\overline{G}(\gamma)$ defined in Proposition 4.20, one obtains the set

$$\overline{H}(\gamma) := \Big\{\overline{\Phi} : \widetilde{\Phi} \in \widetilde{H}(\gamma)\Big\}$$

whose members are, therefore, exact and precisely, the prolongations to $\overline{G}(\gamma)$, as defined in the referred proposition, of the endomorphisms $\widetilde{\Phi} \in \widetilde{H}(\gamma)$.

Given $\overline{\Phi}$ and $\overline{\Psi}$ in $\overline{H}(\gamma)$, we define $\overline{\Phi}\ \overline{\Psi}$ as the following endomorphism on $\overline{G}(\gamma)$:

$$\begin{aligned} \overline{\Phi}\ \overline{\Psi} : \overline{G}(\gamma) &\longrightarrow \overline{G}(\gamma) \\ \overline{g} &\longmapsto (\overline{\Phi}\ \overline{\Psi})(\overline{g}) \end{aligned}$$

where

$$(\overline{\Phi}\ \overline{\Psi})(\overline{g}) := \overline{\Phi}\Big(\overline{\Psi}(\overline{g})\Big).$$

It results that $\overline{\Phi}\ \overline{\Psi} \in \overline{H}(\gamma)$. In fact, let $\overline{g} = [\{\widetilde{g}_\xi\}_{\xi\in\overline{\xi}(\gamma)}] \in \overline{G}(\gamma)$ be arbitrarily chosen and let us calculate $(\overline{\Phi}\ \overline{\Psi})(\overline{g})$. We have:

$$\begin{aligned} (\overline{\Phi}\ \overline{\Psi})(\overline{g}) = \overline{\Phi}\Big(\overline{\Psi}(\overline{g})\Big) &= \overline{\Phi}\left(\Big[\Big\{\widetilde{i}_{(\xi,\gamma)}(\widetilde{\Psi})(\widetilde{g}_\xi)\Big\}_{\xi\in\overline{\xi}(\gamma)}\Big]\right) = \\ &= \Big[\Big\{\widetilde{i}_{(\xi,\gamma)}(\widetilde{\Phi})\Big(\widetilde{i}_{(\xi,\gamma)}(\widetilde{\Psi})(\widetilde{g}_\xi)\Big)\Big\}_{\xi\in\overline{\xi}(\gamma)}\Big]. \end{aligned}$$

On the other hand (see Definition 3.11)

$$\widetilde{i}_{(\xi,\gamma)}(\widetilde{\Phi})\widetilde{i}_{(\xi,\gamma)}(\widetilde{\Psi}) = \widetilde{i}_{(\xi,\gamma)}(\widetilde{\Phi}\widetilde{\Psi}) = \widetilde{i}_{(\xi,\gamma)}(\widetilde{\Phi\Psi}),$$

and, hence,

$$(\overline{\Phi}\ \overline{\Psi})(\overline{g}) = \left[\left\{\widetilde{i}_{(\xi,\gamma)}(\widetilde{\Phi\Psi})(\widetilde{g}_\xi)\right\}_{\xi\in\bar{\xi}(\gamma)}\right] = (\overline{\Phi\Psi})(\overline{g})$$

for every $\overline{g} \in \overline{G}(\gamma)$, that is,

$$\overline{\Phi}\ \overline{\Psi} = \overline{\Phi\Psi},$$

i.e., $\overline{\Phi}\ \overline{\Psi}$ is the endomorphism on $\overline{G}(\gamma)$ associated (through the definition given in Proposition 4.20) with the endomorphism $\widetilde{\Phi\Psi} \in \widetilde{H}(\gamma)$.

Therefore, since the members of $\overline{H}(\gamma)$ are, exact and precisely, the prolongations (as defined in Proposition 4.20) of the members of $\widetilde{H}(\gamma)$,

$$\overline{\Phi}\ \overline{\Psi} = \overline{\Phi\Psi} \in \overline{H}(\gamma).$$

This result allows us to affirm that the operation $\bullet$ defined as

$$\begin{aligned} \bullet : \overline{H}(\gamma) \times \overline{H}(\gamma) &\longrightarrow \overline{H}(\gamma) \\ (\overline{\Phi}, \overline{\Psi}) &\longmapsto \overline{\Phi} \bullet \overline{\Psi} := \overline{\Phi}\ \overline{\Psi}, \end{aligned}$$

where

$$\begin{aligned} \overline{\Phi}\ \overline{\Psi} : \overline{G}(\gamma) &\longrightarrow \overline{G}(\gamma) \\ \overline{g} &\longmapsto (\overline{\Phi}\ \overline{\Psi})(\overline{g}) := \overline{\Phi}\Big(\overline{\Psi}(\overline{g})\Big), \end{aligned}$$

is a well-defined binary operation in $\overline{H}(\gamma)$ and, yet, that $\overline{H}(\gamma)$ with this operation is a semigroup of endomorphisms on $\overline{G}(\gamma)$. Since $\widetilde{\Phi\Psi} = \widetilde{\Psi\Phi}$, then, $\overline{\Phi\Psi} = \overline{\Psi\Phi}$ and, hence, $\overline{\Phi}\ \overline{\Psi} = \overline{\Psi}\ \overline{\Phi}$, that is, the semigroup $\overline{H}(\gamma)$ is abelian.

The considerations above compose a proof for the proposition stated below.

## 4.22 Proposition

*The set*

$$\overline{H}(\gamma) := \left\{\overline{\Phi} : \widetilde{\Phi} \in \widetilde{H}(\gamma)\right\}$$

*equipped with the operation*

$$\begin{aligned} \bullet : \overline{H}(\gamma) \times \overline{H}(\gamma) &\longrightarrow \overline{H}(\gamma) \\ (\overline{\Phi}, \overline{\Psi}) &\longmapsto \overline{\Phi} \bullet \overline{\Psi} := \overline{\Phi}\ \overline{\Psi}, \end{aligned}$$

*both defined in 4.21, is an abelian semigroup isomorphic to the semigroup $\widetilde{H}(\gamma)$, being the function*

$$\begin{aligned} \overline{\gamma} : \widetilde{H}(\gamma) &\longrightarrow \overline{H}(\gamma) \\ \widetilde{\Phi} &\longmapsto \overline{\gamma}(\widetilde{\Phi}) := \overline{\Phi} \end{aligned}$$

*an isomorphism. Furthermore, if $\widehat{\mathscr{G}}(I)$ does exist, we have (according to Proposition 4.20) that*

$$\overline{\Phi} = A_\gamma^{-1}\widehat{\Phi}A_\gamma,$$

*where*

$$\widehat{\Phi} = \widehat{\gamma}(\widetilde{\Phi})$$

*and $A_\gamma : \overline{G}(\gamma) \longrightarrow \widehat{G}(\gamma)$ is the function defined in Proposition 4.12.*

## 4.23 The Bonded Family $\left(\overline{\mathbb{G}}(\Gamma(I)), \overline{i}(\Gamma^2(I))\right)$

The propositions 4.16, 4.18 and 4.22, related, respectively, with the group $\overline{G}(\gamma)$, the injective homomorphism $b_\gamma : \widetilde{G}(\gamma) \longrightarrow \overline{G}(\gamma)$ and the semigroup $\overline{H}(\gamma)$, allow us to affirm, regarding the extension $\widetilde{\mathscr{G}}(I)$ of the $S$-space $\mathscr{G}(I)$, for $\gamma \in \Gamma(I)$ arbitrarily fixed, that:

**(a)** the group $\overline{G}(\gamma)$ referred to in Proposition 4.16 admits, as follows from Proposition 4.18, the group $\widetilde{G}(\gamma)$ as one of its subgroups;

**(b)** the semigroup $\overline{H}(\gamma)$, defined in Proposition 4.22, is a prolongation (or, more precisely, a $b_\gamma$-prolongation, according to Definition 1.6(c)) of the semigroup $\widetilde{H}(\gamma)$ to $\overline{G}(\gamma)$;

**(c)** with $\overline{\mathbb{G}}(\gamma)$ defined by

$$\overline{\mathbb{G}}(\gamma) := \left(\overline{G}(\gamma), \overline{H}(\gamma)\right),$$

it results from (a) and (b) above that $\overline{\mathbb{G}}(\gamma)$ is a $S$-group which is an extension (more precisely, a $b_\gamma$-extension, according to Definition 1.6(d)) of the $S$-group $\widetilde{\mathbb{G}}(\gamma) = (\widetilde{G}(\gamma), \widetilde{H}(\gamma))$;

**(d)** with $\overline{\mathbb{G}}(\Gamma(I))$ defined by

$$\overline{\mathbb{G}}\Big(\Gamma(I)\Big) := \Big\{\overline{\mathbb{G}}(\gamma) = \Big(\overline{G}(\gamma), \overline{H}(\gamma)\Big)\Big\}_{\gamma\in\Gamma(I)},$$

it results from (c) that $\overline{\mathbb{G}}(\Gamma(I))$ is a family of $S$-groups which is an extension of the family

$$\widetilde{\mathbb{G}}\Big(\Gamma(I)\Big) = \Big\{\widetilde{\mathbb{G}}(\gamma) = \Big(\widetilde{G}(\gamma), \widetilde{H}(\gamma)\Big)\Big\}_{\gamma\in\Gamma(I)}$$

(recap Definition 3.8);

**(e)** since $\overline{\mathbb{G}}(\Gamma(I))$ is an extension of the family $\widetilde{\mathbb{G}}(\Gamma(I))$, then, by Definition 3.22, the extension to $\overline{\mathbb{G}}(\Gamma(I))$ of the bonding $\widetilde{i}(\Gamma^2(I))$ (of the bonded family $(\widetilde{G}(\Gamma(I)), \widetilde{i}(\Gamma^2(I)))$ of the $S$-space $\widetilde{\mathscr{G}}(I))$ is the family

$$\overline{i}\Big(\Gamma^2(I)\Big) = \Big\{\overline{i}_{(\gamma',\gamma)}\Big\}_{(\gamma',\gamma)\in\Gamma^2(I)}$$

where, for each $(\gamma', \gamma) \in \Gamma^2(I)$, $\overline{i}_{(\gamma',\gamma)}$ is the isomorphism from $\overline{H}(\gamma)$ onto $\overline{H}(\gamma')$ defined by:

$$\overline{i}_{(\gamma',\gamma)} : \overline{H}(\gamma) \longrightarrow \overline{H}(\gamma')$$
$$\overline{\Phi} \longmapsto \overline{i}_{(\gamma',\gamma)}(\overline{\Phi}) := \overline{\gamma}'\left(\widetilde{i}_{(\gamma',\gamma)}\left((\overline{\gamma})^{-1}(\overline{\Phi})\right)\right).$$

As we know (recap item 3.21), $\overline{i}(\Gamma^2(I))$ is a bonding of the family $\overline{\mathbb{G}}(\Gamma(I))$;

**(f)** finally, from (e) and the definitions of bonded family (Definition 3.11) and extension of bonded family (Definition 3.22), we conclude that

$$\left(\overline{\mathbb{G}}\Big(\Gamma(I)\Big) = \Big\{\overline{\mathbb{G}}(\gamma) = \Big(\overline{G}(\gamma), \overline{H}(\gamma)\Big)\Big\}_{\gamma\in\Gamma(I)}, \overline{i}\Big(\Gamma^2(I)\Big) = \Big\{\overline{i}_{(\gamma',\gamma)}\Big\}_{(\gamma',\gamma)\in\Gamma^2(I)}\right)$$

is a bonded family which is an extension of the bonded family $(\widetilde{\mathbb{G}}(\Gamma(I)), \widetilde{i}(\Gamma^2(I)))$ of the $S$-space $\widetilde{\mathscr{G}}(I)$.

# The Restriction $\overline{\Theta}\,(\Delta(I))$

## 4.24 Preliminaries

Now, we would like to define a family

$$\overline{\Theta}\Big(\Delta(I)\Big) = \Big\{\overline{\Theta}_{(\gamma',\gamma)}\Big\}_{(\gamma',\gamma)\in\Delta(I)}$$

whose members, $\overline{\Theta}_{(\gamma',\gamma)}$, are functions from $\overline{G}(\gamma)$ into $\overline{G}(\gamma')$ ($\overline{\Theta}_{(\gamma',\gamma)} : \overline{G}(\gamma) \longrightarrow \overline{G}(\gamma')$), satisfying the conditions required (by Definition 3.13) in order for $\overline{\Theta}(\Delta(I))$ to be a restriction for the bonded family $(\overline{\mathbb{G}}(\Gamma(I)), \overline{i}(\Gamma^2(I)))$ defined in 4.23.

Furthermore, we would like $\overline{\Theta}(\Delta(I))$ to be a prolongation of $\widetilde{\Theta}(\Delta(I))$ (the restriction of the $S$-space $\widetilde{\mathscr{G}}(I)$) to $(\overline{\mathbb{G}}(\Gamma(I)), \overline{i}(\Gamma^2(I)))$, that is,

$$\overline{\Theta}_{(\gamma',\gamma)}(\widetilde{g}) = \widetilde{\Theta}_{(\gamma',\gamma)}(\widetilde{g})$$

or, more precisely (recap the "abuse of language" mentioned in item (b) of Definition 1.6), for each $(\gamma', \gamma) \in \Delta(I)$,

$$\overline{\Theta}_{(\gamma',\gamma)}\Big(b_\gamma(\widetilde{g})\Big) = b_{\gamma'}\Big(\widetilde{\Theta}_{(\gamma',\gamma)}(\widetilde{g})\Big) \quad \text{for every} \quad \widetilde{g} \in \widetilde{G}(\gamma), \tag{4.24-1}$$

with $b_\gamma : \widetilde{G}(\gamma) \longrightarrow \overline{G}(\gamma)$ ($b_{\gamma'} : \widetilde{G}(\gamma') \longrightarrow \overline{G}(\gamma')$) as in Proposition 4.18.

With such a family $\overline{\Theta}(\Delta(I))$, defining

$$\overline{\mathscr{G}}(I) := \left(\overline{\mathbb{G}}\Big(\Gamma(I)\Big), \overline{i}\Big(\Gamma^2(I)\Big), \overline{\Theta}\Big(\Delta(I)\Big)\right),$$

$\overline{\mathscr{G}}(I)$ would be an extension of $\widetilde{\mathscr{G}}(I)$ (recap Definition 3.23). However, the extension $\overline{\mathscr{G}}(I)$ that we are willing to obtain for $\widetilde{\mathscr{G}}(I)$ must be isomorphic to $\widehat{\mathscr{G}}(I)$ on the hypothesis of existence of the latter. Hence, in order for our odds of success in this pursuit for $\overline{\Theta}(\Delta(I))$ to be higher, we will (such as we did to obtain the family of $S$-groups $\overline{\mathbb{G}}(\Gamma(I)) = \{\overline{\mathbb{G}}(\gamma) = (\overline{G}(\gamma), \overline{H}(\gamma))\}_{\gamma\in\Gamma(I)}$) resort to the hypothesis of existence of $\widehat{\mathscr{G}}(I)$, looking for suggestions on how to define the functions $\overline{\Theta}_{(\gamma',\gamma)}$.

We assume then the existence of $\widehat{\mathscr{G}}(I)$. In this case, the homomorphisms $\widehat{\Theta}_{(\gamma',\gamma)} : \widehat{G}(\gamma) \longrightarrow \widehat{G}(\gamma')$, for each $(\gamma',\gamma) \in \Delta(I)$, of the family $\widehat{\Theta}(\Delta(I))$ (the restriction of the $S$-space $\widehat{\mathscr{G}}(I)$) are, as well as the members $\overline{\Theta}_{(\gamma',\gamma)} : \overline{G}(\gamma) \longrightarrow \overline{G}(\gamma')$ of the family $\overline{\Theta}(\Delta(I))$ also should be, extensions (to $\widehat{G}(\gamma) \supseteq \widetilde{G}(\gamma)$) of the homomorphisms $\widetilde{\Theta}_{(\gamma',\gamma)} : \widetilde{G}(\gamma) \longrightarrow \widetilde{G}(\gamma')$. Thus, it is necessary that, for each $(\gamma',\gamma) \in \Delta(I)$, $\overline{\Theta}_{(\gamma',\gamma)}$ and $\widetilde{\Theta}_{(\gamma',\gamma)}$ meet the condition (4.24-1), whereas $\widehat{\Theta}_{(\gamma',\gamma)}$ and $\widetilde{\Theta}_{(\gamma',\gamma)}$ satisfy the following requirement:

$$\widehat{\Theta}_{(\gamma',\gamma)}(\widetilde{g}) = \widetilde{\Theta}_{(\gamma',\gamma)}(\widetilde{g}) \quad \text{for every} \quad \widetilde{g} \in \widetilde{G}(\gamma). \tag{4.24-2}$$

It results from (4.24-1) and (4.24-2) that

$$\overline{\Theta}_{(\gamma',\gamma)}\Big(b_\gamma(\widetilde{g})\Big) = b_{\gamma'}\Big(\widehat{\Theta}_{(\gamma',\gamma)}(\widetilde{g})\Big) \quad \text{for every} \quad \widetilde{g} \in \widetilde{G}(\gamma).$$

Taking

$$\overline{g} := b_\gamma(\widetilde{g}) \quad \text{and, therefore,} \quad \widetilde{g} = b_\gamma^{-1}(\overline{g})$$

one obtains

$$\overline{\Theta}_{(\gamma',\gamma)}(\overline{g}) = b_{\gamma'}\left(\widehat{\Theta}_{(\gamma',\gamma)}\Big(b_\gamma^{-1}(\overline{g})\Big)\right) \quad \text{for every} \quad \overline{g} \in b_\gamma\Big(\widetilde{G}(\gamma)\Big),$$

where $b_\gamma(\widetilde{G}(\gamma)) = \{b_\gamma(\widetilde{g}) : \widetilde{g} \in \widetilde{G}(\gamma)\}$.

But, from Proposition 4.18 we know that, for every $\gamma \in \Gamma(I)$,

$$A_\gamma\Big(b_\gamma(\widetilde{g})\Big) = \widetilde{g} \quad \text{for every} \quad \widetilde{g} \in \widetilde{G}(\gamma),$$

that is,

$$A_\gamma(\overline{g}) = b_\gamma^{-1}(\overline{g}) \quad \text{for every} \quad \overline{g} \in b_\gamma(\widetilde{G}(\gamma)),$$

and, therefore,

$$b_\gamma^{-1} = A_\gamma \quad \text{in} \quad b_\gamma\Big(\widetilde{G}(\gamma)\Big)$$

or

$$b_\gamma = A_\gamma^{-1} \quad \text{in} \quad \widetilde{G}(\gamma).$$

Thus, the last expression obtained for $\overline{\Theta}_{(\gamma',\gamma)}(\overline{g})$ can be written in the following form:

$$\overline{\Theta}_{(\gamma',\gamma)}(\overline{g}) = b_{\gamma'}\left(\widehat{\Theta}_{(\gamma',\gamma)}\Big(A_\gamma(\overline{g})\Big)\right) = A_{\gamma'}^{-1}\left(\widehat{\Theta}_{(\gamma',\gamma)}\Big(A_\gamma(\overline{g})\Big)\right) \tag{4.24-3}$$

for every $\overline{g} \in b_\gamma(\widetilde{G}(\gamma)) \subseteq \overline{G}(\gamma)$.

We can summarize the considerations above as follows: if $\widehat{\mathscr{G}}(I)$ does exist, then, for $(\gamma', \gamma) \in \Delta(I)$, a function $\overline{\Theta}_{(\gamma',\gamma)} : \overline{G}(\gamma) \longrightarrow \overline{G}(\gamma')$ is an "extension" (to $\overline{G}(\gamma) \supseteq b_\gamma(\widetilde{G}(\gamma)) \equiv \widetilde{G}(\gamma)$) of $\widetilde{\Theta}_{(\gamma',\gamma)} : \widetilde{G}(\gamma) \longrightarrow \widetilde{G}(\gamma')$, if an only if it satisfies the condition (4.24-3) above.

Now, if $\widehat{\mathscr{G}}(I)$ does exist, then, for each $(\gamma', \gamma) \in \Delta(I)$, the function

$$\begin{aligned} \overline{\Theta}_{(\gamma',\gamma)} : \overline{G}(\gamma) &\longrightarrow \overline{G}(\gamma') \\ \overline{g} &\longmapsto \overline{\Theta}_{(\gamma',\gamma)}(\overline{g}) \end{aligned}$$

defined by

$$\overline{\Theta}_{(\gamma',\gamma)}(\overline{g}) := A_{\gamma'}^{-1}\left(\widehat{\Theta}_{(\gamma',\gamma)}\Big(A_\gamma(\overline{g})\Big)\right) \quad \text{for every} \quad \overline{g} \in \overline{G}(\gamma) \tag{4.24-4}$$

beyond well-defined, is a homomorphism (since $A_\gamma : \overline{G}(\gamma) \longrightarrow \widehat{G}(\gamma)$ and $A_{\gamma'} : \overline{G}(\gamma') \longrightarrow \widehat{G}(\gamma')$ are isomorphisms and $\widehat{\Theta}_{(\gamma',\gamma)} : \widehat{G}(\gamma) \longrightarrow \widehat{G}(\gamma')$ is a homomorphism) that satisfies, due to its definition, the condition (4.24-3) and, therefore, also is an "extension" of $\widetilde{\Theta}_{(\gamma',\gamma)} : \widetilde{G}(\gamma) \longrightarrow \widetilde{G}(\gamma')$ (to $\overline{G}(\gamma) \supseteq b_\gamma(\widetilde{G}(\gamma)) \equiv \widetilde{G}(\gamma)$). Let us calculate $\overline{\Theta}_{(\gamma',\gamma)}(\overline{g})$ for

$$\overline{g} = \left[\widetilde{g}\Big(\overline{\xi}(\gamma)\Big) = \Big\{\widetilde{g}_\xi\Big\}_{\xi \in \overline{\xi}(\gamma)}\right]$$

arbitrarily chosen in $\overline{G}(\gamma)$. As we know, $A_\gamma(\overline{g})$ is the only member of $\widehat{G}(\gamma)$ such that

$$\widehat{\Theta}_{(\xi,\gamma)}\Big(A_\gamma(\overline{g})\Big) = \widetilde{g}_\xi \quad \text{for every} \quad \xi \in \overline{\xi}(\gamma).$$

Let us take:

$$\overline{\eta}(\gamma') := \Big\{\eta := \xi \cap \gamma' \quad : \quad \xi \in \overline{\xi}(\gamma) \quad \text{and} \quad \xi \cap \gamma' \neq \varnothing\Big\}.$$

Hence, $\overline{\eta}(\gamma') \subseteq \Gamma(\gamma')$ is a cover of $\gamma'$. Now, from the expression above for $\widehat{\Theta}_{(\xi,\gamma)}(A_\gamma(\overline{g}))$ we obtain

$$\widehat{\Theta}_{(\xi\cap\gamma',\xi)}\left(\widehat{\Theta}_{(\xi,\gamma)}\Big(A_\gamma(\overline{g})\Big)\right) = \widehat{\Theta}_{(\xi\cap\gamma',\xi)}(\widetilde{g}_\xi)$$

that is,

$$\widehat{\Theta}_{(\xi\cap\gamma',\gamma)}\Big(A_{(\gamma)}(\overline{g})\Big) = \widehat{\Theta}_{(\xi\cap\gamma',\xi)}(\widetilde{g}_\xi)$$

or yet,

$$\widehat{\Theta}_{(\xi\cap\gamma',\gamma')}\left(\widehat{\Theta}_{(\gamma',\gamma)}\Big(A_\gamma(\overline{g})\Big)\right) = \widehat{\Theta}_{(\xi\cap\gamma',\xi)}(\widetilde{g}_\xi)$$

for every $\xi \cap \gamma' = \eta \in \overline{\eta}(\gamma')$.

But, $\widehat{\Theta}_{(\xi\cap\gamma',\xi)}$ is an extension of $\widetilde{\Theta}_{(\xi\cap\gamma',\xi)}$, that is,

$$\widehat{\Theta}_{(\xi\cap\gamma',\xi)}(\widetilde{g}_\xi) = \widetilde{\Theta}_{(\xi\cap\gamma',\xi)}(\widetilde{g}_\xi)$$

and, hence, we obtain that

$$\widehat{\Theta}_{(\xi\cap\gamma',\gamma')}\left(\widehat{\Theta}_{(\gamma',\gamma)}\left(A_{\gamma}(\overline{g})\right)\right) = \widetilde{\Theta}_{(\xi\cap\gamma',\xi)}(\widetilde{g}_{\xi})$$

for every $\xi \cap \gamma' = \eta \in \overline{\eta}(\gamma')$.

From this last expression and taking into account the definition of the function $A_{\gamma'}$ one concludes that

$$A_{\gamma'}^{-1}\left(\widehat{\Theta}_{(\gamma',\gamma)}\left(A_{\gamma}(\overline{g})\right)\right) = \left[\left\{\widetilde{\Theta}_{(\xi\cap\gamma',\xi)}(\widetilde{g}_{\xi})\right\}_{\xi\cap\gamma'\in\overline{\eta}(\gamma')}\right],$$

that is, according to the definition of $\overline{\Theta}_{(\gamma',\gamma)}$,

$$\overline{\Theta}_{(\gamma',\gamma)}\left(\overline{g} = \left[\left\{\widetilde{g}_{\xi}\right\}_{\xi\in\overline{\xi}(\gamma)}\right]\right) = \left[\left\{\widetilde{\Theta}_{(\xi\cap\gamma',\xi)}(\widetilde{g}_{\xi})\right\}_{\xi\cap\gamma'\in\overline{\eta}(\gamma')}\right]. \tag{4.24-5}$$

Let us remark now that, although the expression above for $\overline{\Theta}_{(\gamma',\gamma)}(\overline{g})$ has been obtained under hypothesis of existence of $\widehat{\mathscr{G}}(I)$, the referred expression does not involve elements of $\widehat{\mathscr{G}}(I)$ but, only, objects of $\widetilde{\mathscr{G}}(I)$. Hence, we inquire: If we define $\overline{\Theta}_{(\gamma',\gamma)}$ : $\overline{G}(\gamma) \longrightarrow \overline{G}(\gamma')$ with the expression above, would we have a well-defined function that is an "extension" of $\widetilde{\Theta}_{(\gamma',\gamma)} : \widetilde{G}(\gamma) \longrightarrow \widetilde{G}(\gamma')$?

The answer to this question is provided by Proposition 4.26. Before, however, in item 4.25 ahead, we clarify some points of the considerations which led us to the above formula for $\overline{\Theta}_{(\gamma',\gamma)}(\overline{g})$.

## 4.25 Remark

Let us, again, consider the abuse of language introduced by Definition 1.6(b), "officialized" at Remark 1.7, retaken at Warning 3.24 and analysed in the context of Motivation 4.19 (see the Remark in 4.19). Our goal here is analogous to the one considered in Motivation 4.19, namely, to show how "harmless" the referred abuse is regarding the conclusions in the previous item (4.24); more precise and specifically, we are going to prove that the formula for $\overline{\Theta}_{(\gamma',\gamma)}(\overline{g})$, given in (4.24-5), obtained from expressions involving the referred abuse of language, does not change when deduced using "formal language".

Resuming item 4.24, we see that the obtainment of the expression (4.24-5) for $\overline{\Theta}_{(\gamma',\gamma)}(\overline{g})$ has as starting point the equations (4.24-1) and (4.24-2), that is,

$$\overline{\Theta}_{(\gamma',\gamma)}\left(b_{\gamma}(\widetilde{g})\right) = b_{\gamma'}\left(\widetilde{\Theta}_{(\gamma',\gamma)}(\widetilde{g})\right) \quad \text{for every} \quad \widetilde{g} \in \widetilde{G}(\gamma)$$

and

$$\widehat{\Theta}_{(\gamma',\gamma)}(\widetilde{g}) = \widetilde{\Theta}_{(\gamma',\gamma)}(\widetilde{g}) \quad \text{for every} \quad \widetilde{g} \in \widetilde{G}(\gamma),$$

from which we obtained that

$$\overline{\Theta}_{(\gamma',\gamma)}(\overline{g}) = b_{\gamma'}\left(\widehat{\Theta}_{(\gamma',\gamma)}\left(b_\gamma^{-1}(\overline{g})\right)\right) \quad \text{for every} \quad \overline{g} \in b_\gamma\left(\widetilde{G}(\gamma)\right).$$

Regarding the first of the last three equations above there is nothing to remark, however, the second, and consequently the latter one (once it its obtained with the first two), assume that the group $\widehat{G}(\gamma)$ "trivially" admits the group $\widetilde{G}(\gamma)$ as one of its subgroups (see Warning 3.24, as well as the Remark in Motivation 4.19). But, what we know about $\widehat{G}(\gamma)$ regarding $\widetilde{G}(\gamma)$ is that $\widehat{G}(\gamma)$ $\delta_\gamma$-admits the group $\widetilde{G}(\gamma)$ as a subgroup (see Definition 1.6(b)), being $\delta_\gamma : \widetilde{G}(\gamma) \longrightarrow \widehat{G}(\gamma)$ an injective homomorphism. Then, all we can say is that $\delta_\gamma(\widetilde{G}(\gamma)) \subseteq \widehat{G}(\gamma)$ is a subgroup of $\widehat{G}(\gamma)$, not that $\widetilde{G}(\gamma) \subseteq \widehat{G}(\gamma)$, and much less that $\widetilde{G}(\gamma)$ is a subgroup of $\widehat{G}(\gamma)$. Clearly, one possibility is that $\delta_\gamma = I_{\widetilde{G}(\gamma)}$ (see Warning 3.24), in which case we would have

$$\delta_\gamma\left(\widetilde{G}(\gamma)\right) = I_{\widetilde{G}(\gamma)}\left(\widetilde{G}(\gamma)\right) = \widetilde{G}(\gamma) \subseteq \widehat{G}(\gamma)$$

and, hence, and only at this case (of $\delta_\gamma = I_{\widetilde{G}(\gamma)}$), $\widetilde{G}(\gamma)$ would be a subgroup of $\widehat{G}(\gamma)$.

Thus, the second one of the equations we are taking as our "starting point" involves the abuse of language under consideration, therefore being an abbreviated form of the equation

$$\widehat{\Theta}_{(\gamma',\gamma)}\left(\delta_\gamma(\widetilde{g})\right) = \delta_{\gamma'}\left(\widetilde{\Theta}_{\gamma',\gamma}(\widetilde{g})\right) \quad \text{for every} \quad \widetilde{g} \in \widetilde{G}(\gamma)$$

where (for each $\gamma \in \Gamma(I)$)

$$\begin{aligned} \delta_\gamma : \widetilde{G}(\gamma) &\longrightarrow \widehat{G}(\gamma) \\ \widetilde{g} &\longmapsto \delta_\gamma(\widetilde{g}), \end{aligned}$$

as previously remarked, is an injective homomorphism (which allows us to identify the group $\widetilde{G}(\gamma)$ with the subgroup $\delta_\gamma(\widetilde{G}(\gamma)) \subseteq \widehat{G}(\gamma)$ of the group $\widehat{G}(\gamma)$). Our starting point for the obtainment of $\overline{\Theta}_{(\gamma',\gamma)}(\overline{g})$ then becomes

$$\overline{\Theta}_{(\gamma',\gamma)}\left(b_\gamma(\widetilde{g})\right) = b_{\gamma'}\left(\widetilde{\Theta}_{(\gamma',\gamma)}(\widetilde{g})\right) \quad \text{for every} \quad \widetilde{g} \in \widetilde{G}(\gamma)$$

and

$$\widehat{\Theta}_{(\gamma',\gamma)}\left(\delta_\gamma(\widetilde{g})\right) = \delta_{\gamma'}\left(\widetilde{\Theta}_{(\gamma',\gamma)}(\widetilde{g})\right) \quad \text{for every} \quad \widetilde{g} \in \widetilde{G}(\gamma).$$

Another element used at the obtainment of expression (4.24-5) for $\overline{\Theta}_{(\gamma',\gamma)}(\overline{g})$ and which also involves the abuse of language in question, is the function $A_\gamma : \overline{G}(\gamma) \longrightarrow \widehat{G}(\gamma)$. In fact, as defined in Proposition 4.12, we have (for each $\gamma \in \Gamma(I)$):

$$\begin{aligned} A_\gamma : \overline{G}(\gamma) = \left[\widetilde{F}(\gamma)\right] &\longrightarrow \widehat{G}(\gamma) \\ \overline{g} = \left[\{\widetilde{g}_\xi\}_{\xi \in \bar{\xi}(\gamma)}\right] &\longmapsto A_\gamma(\overline{g}) = \widehat{g} \end{aligned}$$

where $\widehat{g} \in \widehat{G}(\gamma)$ is such that

$$\widehat{\Theta}_{(\xi,\gamma)}(\widehat{g}) = \widetilde{g}_\xi \quad \text{for every} \quad \xi \in \overline{\xi}(\gamma).$$

But, $\widehat{\Theta}_{(\xi,\gamma)}(\widehat{g}) \in \widehat{G}(\xi)$ while $\widetilde{g}_\xi \in \widetilde{G}(\xi)$ and, hence, we have $\widetilde{g}_\xi \in \widehat{G}(\xi)$ for every $\xi \in \overline{\xi}(\gamma)$. However, what we know, as we saw above, is that $\delta_\xi(\widetilde{g}_\xi) \in \widehat{G}(\xi)$ for every $\xi \in \overline{\xi}(\gamma)$ and not that $\widetilde{g}_\xi \in \widehat{G}(\xi)$; the formula $\widehat{\Theta}_{\xi,\gamma}(\widehat{g}) = \widetilde{g}_\xi$ in the definition of $A_\gamma$ requires $\widetilde{G}(\xi)$ (for each $\xi \in \overline{\xi}(\gamma)$) to be, in fact, a subgroup of $\widehat{G}(\xi)$ and not only that $\widehat{G}(\xi)$ $\delta_\xi$-admits $\widetilde{G}(\xi)$ as a subgroup, unless the abuse of language in question had already been previously "officialized".

The above considerations about the function $A_\gamma$ suggest that this function, undressed of the camouflage in which $\widetilde{G}(\gamma)$ denotes $\delta_\gamma(\widetilde{G}(\gamma))$ and $\widetilde{g} \in \widetilde{G}(\gamma)$ represents the element $\delta_\gamma(\widetilde{g}) \in \widehat{G}(\gamma)$, has the following definition:

$$\begin{aligned} A_\gamma : \qquad & \overline{G}(\gamma) \longrightarrow \widehat{G}(\gamma) \\ \overline{g} = \Big[\{\widetilde{g}_\xi\}_{\xi \in \overline{\xi}(\gamma)}\Big] & \longmapsto A_\gamma(\overline{g}) = \widehat{g} \end{aligned}$$

where $\widehat{g} \in \widehat{G}(\gamma)$ is such that

$$\widehat{\Theta}_{(\xi,\gamma)}(\widehat{g}) = \delta_\xi(\widetilde{g}_\xi) \quad \text{for every} \quad \xi \in \overline{\xi}(\gamma).$$

Although quite reasonable, the modification suggested in the definition of $A_\gamma$ requires justification. In fact, the definition of the "camouflaged" version of $A_\gamma$, in which

$$A_\gamma\left(\overline{g} = \left[\left\{\widetilde{g}_\xi\right\}_{\xi \in \overline{\xi}(\gamma)}\right]\right) = \widehat{g}$$

being $\widehat{g} \in \widehat{G}(\gamma)$ such that $\widehat{\Theta}_{(\xi,\gamma)}(\widehat{g}) = \widetilde{g}_\xi$ for every $\xi \in \overline{\xi}(\gamma)$, was justified as follows: since $\overline{g} = [\{\widetilde{g}_\xi\}_{\xi \in \overline{\xi}(\gamma)}] \in \overline{G}(\gamma)$, we have $\{\widetilde{g}_\xi\}_{\xi \in \overline{\xi}(\gamma)} \in \widetilde{F}(\gamma)$, that is, $\{\widetilde{g}_\xi\}_{\xi \in \overline{\xi}(\gamma)}$ is a coherent family in $\widetilde{\mathscr{G}}(\gamma)$ and, hence, once $\widetilde{G}(\gamma) \subseteq \widehat{G}(\gamma)$ for $\gamma \in \Gamma(I)$ and $\widehat{\Theta}_{(\gamma',\gamma)}(\widetilde{g}) = \widetilde{\Theta}_{(\gamma',\gamma)}(\widetilde{g})$ for every $\widetilde{g} \in \widetilde{G}(\gamma)$ since $\widehat{\Theta}_{(\gamma',\gamma)}$ is an extension of $\widetilde{\Theta}_{(\gamma',\gamma)}$, $\{\widetilde{g}_\xi\}_{\xi \in \overline{\xi}(\gamma)}$ also is a coherent family in $\widehat{\mathscr{G}}(\gamma)$ and, therefore, since $\widehat{\mathscr{G}}(\gamma)$ is coherent, there exists a single $\widehat{g} \in \widehat{G}(\gamma)$ such that $\widehat{\Theta}_{(\xi,\gamma)}(\widehat{g}) = \widetilde{g}_\xi$ for every $\xi \in \overline{\xi}(\gamma)$, which validates the definition (in the camouflaged version) of $A_\gamma$.

Now, the modified version of the $A_\gamma$ definition requires, so that it is justified, the existence, for each $\overline{g} = [\{\widetilde{g}_\xi\}_{\xi \in \overline{\xi}(\gamma)}] \in \overline{G}(\gamma)$, of a single $\widehat{g} \in \widehat{G}(\gamma)$ such that

$$\widehat{\Theta}_{(\xi,\gamma)}(\widehat{g}) = \delta_\xi(\widetilde{g}_\xi) \quad \text{for every} \quad \xi \in \overline{\xi}(\gamma).$$

In other words, to justify the suggested version for the definition of $A_\gamma$, we must prove that the family $\{\delta_\xi(\widetilde{g}_\xi)\}_{\xi \in \overline{\xi}(\gamma)}$ is coherent in $\widehat{\mathscr{G}}(\gamma)$ whenever the family $\{\widetilde{g}_\xi\}_{\xi \in \overline{\xi}(\gamma)}$ is coherent in $\widetilde{\mathscr{G}}(\gamma)$.

Another ingredient that played an important role when getting formula (4.24-5) for $\overline{\Theta}_{(\gamma',\gamma)}(\overline{g})$, and which also involves the abuse of language in question, consists of the following property of the functions $A_\gamma : \overline{G}(\gamma) \longrightarrow \widehat{G}(\gamma)$ and $b_\gamma : \widetilde{G}(\gamma) \longrightarrow \overline{G}(\gamma)$ (described in Proposition 4.18):

$$A_\gamma\Big(b_\gamma(\widetilde{g})\Big) = \widetilde{g} \quad \text{for every} \quad \widetilde{g} \in \widetilde{G}(\gamma).$$

Clearly, $b_\gamma(\widetilde{g}) \in \overline{G}(\gamma)$ and, hence, $A_\gamma(b_\gamma(\widetilde{g})) \in \widehat{G}(\gamma)$, which, taking into account the above property, leads us to conclude that $\widetilde{g} \in \widehat{G}(\gamma)$, hence revealing the use of the alluded abuse of language. Here, also seems reasonable that the "right" version of the referred property is

$$A_\gamma\Big(b_\gamma(\widetilde{g})\Big) = \delta_\gamma(\widetilde{g}) \quad \text{for every} \quad \widetilde{g} \in \widetilde{G}(\gamma),$$

which also needs justification.

At the end of this item, we will provide the due justifications for the versions given above for the definition of the function $A_\gamma$ and the property mentioned, among others, such as, for instance, that the new version of $A_\gamma$ is an isomorphism. In what follows, assuming them justified, we proceed to the derivation of expression (4.24-5) for $\overline{\Theta}_{\gamma',\gamma}(\overline{g})$. Let us then resume to our new ("uncloaked") "starting point", namely, the expressions

$$\overline{\Theta}_{(\gamma',\gamma)}\Big(b_\gamma(\widetilde{g})\Big) = b_{\gamma'}\Big(\widetilde{\Theta}_{(\gamma',\gamma)}(\widetilde{g})\Big) \quad \text{for every} \quad \widetilde{g} \in \widetilde{G}(\gamma)$$

and

$$\widehat{\Theta}_{(\gamma',\gamma)}\Big(\delta_\gamma(\widetilde{g})\Big) = \delta_{\gamma'}\Big(\widetilde{\Theta}_{(\gamma',\gamma)}(\widetilde{g})\Big) \quad \text{for every} \quad \widetilde{g} \in \widetilde{G}(\gamma).$$

From the last expression one obtains

$$\widetilde{\Theta}_{(\gamma',\gamma)}(\widetilde{g}) = \delta_{\gamma'}^{-1}\left(\widehat{\Theta}_{(\gamma',\gamma)}\Big(\delta_\gamma(\widetilde{g})\Big)\right)$$

which, when taken into the first equation, provides:

$$\overline{\Theta}_{(\gamma',\gamma)}\Big(b_\gamma(\widetilde{g})\Big) = b_{\gamma'}\left(\delta_{\gamma'}^{-1}\left(\widehat{\Theta}_{(\gamma',\gamma)}\Big(\delta_\gamma(\widetilde{g})\Big)\right)\right)$$

for every $\widetilde{g} \in \widetilde{G}(\gamma)$.

Taking

$$\overline{g} := b_\gamma(\widetilde{g}) \quad \text{and, therefore,} \quad \widetilde{g} = b_\gamma^{-1}(\overline{g})$$

one gets

$$\overline{\Theta}_{(\gamma',\gamma)}(\overline{g}) = b_{\gamma'}\left(\delta_{\gamma'}^{-1}\left(\widehat{\Theta}_{(\gamma',\gamma)}\left(\delta_\gamma\Big(b_\gamma^{-1}(\overline{g})\Big)\right)\right)\right)$$

for every

$$\overline{g} \in b_\gamma\Big(\widetilde{G}(\gamma)\Big) = \Big\{b_\gamma(\widetilde{g}) : \widetilde{g} \in \widetilde{G}(\gamma)\Big\} \subseteq \overline{G}(\gamma).$$

But, as assumed above, for every $\gamma \in \Gamma(I)$,

$$A_\gamma\Big(b_\gamma(\widetilde{g})\Big) = \delta_\gamma(\widetilde{g}) \quad \text{for every} \quad \widetilde{g} \in \widetilde{G}(\gamma),$$

that is,

$$A_\gamma(\overline{g}) = \delta_\gamma\Big(b_\gamma^{-1}(\overline{g})\Big) \quad \text{for every} \quad \overline{g} \in b_\gamma\Big(\widetilde{G}(\gamma)\Big),$$

and, therefore,

$$A_\gamma = \delta_\gamma b_\gamma^{-1} \quad \text{in} \quad b_\gamma\Big(\widetilde{G}(\gamma)\Big)$$

or

$$A_\gamma^{-1} = b_\gamma \delta_\gamma^{-1} \quad \text{in} \quad \delta_\gamma\Big(\widetilde{G}(\gamma)\Big) = \Big\{\delta_\gamma(\widetilde{g}) : \widetilde{g} \in \widetilde{G}(\gamma)\Big\}$$

for every $\gamma \in \Gamma(I)$.

With this, the last expression above for $\overline{\Theta}_{\gamma',\gamma}(\overline{g})$ assumes the following form:

$$\overline{\Theta}_{(\gamma',\gamma)}(\overline{g}) = b_{\gamma'}\left(\delta_{\gamma'}^{-1}\left(\widehat{\Theta}_{(\gamma',\gamma)}\Big(A_\gamma(\overline{g})\Big)\right)\right) = A_{\gamma'}^{-1}\left(\widehat{\Theta}_{(\gamma',\gamma)}\Big(A_\gamma(\overline{g})\Big)\right)$$

for every $\overline{g} \in b_\gamma(\widetilde{G}(\gamma))$.

In short, such as in the previous item (review the conclusion formulated in the last paragraph of page 184), the development above allows us to affirm that if $\widehat{\mathscr{G}}(I)$ does exist, then, for each $(\gamma', \gamma) \in \Delta(I)$, a function $\overline{\Theta}_{(\gamma',\gamma)} : \overline{G}(\gamma) \longrightarrow \overline{G}(\gamma')$ is an "extension" of $\widetilde{\Theta}_{(\gamma',\gamma)} : \widetilde{G}(\gamma) \longrightarrow \widetilde{G}(\gamma')$, that is, it attends the condition expressed by the first of the two expressions that compose our (new) "starting point", if and only if,

$$\overline{\Theta}_{(\gamma',\gamma)}(\overline{g}) = A_{\gamma'}^{-1}\left(\widehat{\Theta}_{(\gamma',\gamma)}\Big(A_\gamma(\overline{g})\Big)\right) \quad \text{for every} \quad \overline{g} \in b_\gamma\Big(\widetilde{G}(\gamma)\Big). \tag{4.25-1}$$

Now, again as in the previous item, if $\widehat{\mathscr{G}}(I)$ does exist, then, for $(\gamma', \gamma) \in \Delta(I)$, the function $\overline{\Theta}_{(\gamma',\gamma)}$ defined in **every** $\overline{G}(\gamma)$ by

$$\begin{aligned} \overline{\Theta}_{(\gamma',\gamma)} : \overline{G}(\gamma) &\longrightarrow \overline{G}(\gamma') \\ \overline{g} &\longmapsto \overline{\Theta}_{(\gamma',\gamma)}(\overline{g}) := A_{\gamma'}^{-1}\left(\widehat{\Theta}_{(\gamma',\gamma)}\Big(A_\gamma(\overline{g})\Big)\right) \end{aligned}$$

is well-defined, is a homomorphism (since $A_\gamma$ and $A_{\gamma'}$ are isomorphisms and $\widehat{\Theta}_{(\gamma',\gamma)}$ a homomorphism), and it satisfies, due to its definition, the condition (4.25-1) above and, therefore, it is also an "extension" of $\widetilde{\Theta}_{(\gamma',\gamma)} : \widetilde{G}(\gamma) \longrightarrow \widetilde{G}(\gamma')$.

Let now $\overline{g} = [\{\widetilde{g}_\xi\}_{\xi \in \overline{\xi}(\gamma)}]$ be an arbitrarily chosen element in the domain $\overline{G}(\gamma)$ of the function $\overline{\Theta}_{(\gamma',\gamma)}$ above defined, and let us calculate $\overline{\Theta}_{(\gamma',\gamma)}(\overline{g})$. In order to do so, we remember that $A_\gamma(\overline{g})$, according to the non-camouflaged version of the definition of $A_\gamma$, is the only element of $\widehat{G}(\gamma)$ such that

$$\widehat{\Theta}_{(\xi,\gamma)}\Big(A_\gamma(\overline{g})\Big) = \delta_\xi(\widetilde{g}_\xi) \quad \text{for every} \quad \xi \in \overline{\xi}(\gamma).$$

Let us take

$$\overline{\eta}(\gamma') := \Big\{ \eta := \xi \cap \gamma' \quad : \quad \xi \in \overline{\xi}(\gamma) \quad \text{and} \quad \xi \cap \gamma' \neq \varnothing \Big\}.$$

Clearly, $\overline{\eta}(\gamma') \subseteq \Gamma(\gamma')$ is a cover of $\gamma'$. It then results from the above expression of $\widehat{\Theta}_{(\xi,\gamma)}(A_\gamma(\overline{g}))$ that

$$\widehat{\Theta}_{(\xi\cap\gamma',\xi)}\left(\widehat{\Theta}_{(\xi,\gamma)}\Big(A_\gamma(\overline{g})\Big)\right) = \widehat{\Theta}_{(\xi\cap\gamma',\xi)}\Big(\delta_\xi(\widetilde{g}_\xi)\Big),$$

that is,

$$\widehat{\Theta}_{(\xi\cap\gamma',\gamma)}\Big(A_\gamma(\overline{g})\Big) = \widehat{\Theta}_{(\xi\cap\gamma',\xi)}\Big(\delta_\xi(\widetilde{g}_\xi)\Big),$$

or, yet,

$$\widehat{\Theta}_{(\xi\cap\gamma',\gamma')}\left(\widehat{\Theta}_{(\gamma',\gamma)}\Big(A_\gamma(\overline{g})\Big)\right) = \widehat{\Theta}_{(\xi\cap\gamma',\xi)}\Big(\delta_\xi(\widetilde{g}_\xi)\Big)$$

for every $\xi \cap \gamma' = \eta \in \overline{\eta}(\gamma')$. But, $\widehat{\Theta}_{(\xi\cap\gamma',\xi)}$ is an "extension" of $\widetilde{\Theta}_{(\xi\cap\gamma',\xi)}$, which means, using the "formal language" (without abuse of language), that

$$\widehat{\Theta}_{(\xi\cap\gamma',\xi)}\Big(\delta_\xi(\widetilde{g}_\xi)\Big) = \delta_{\xi\cap\gamma'}\Big(\widetilde{\Theta}_{(\xi\cap\gamma',\xi)}(\widetilde{g}_\xi)\Big)$$

and, hence, we can conclude that

$$\widehat{\Theta}_{(\xi\cap\gamma',\gamma')}\left(\widehat{\Theta}_{(\gamma',\gamma)}\Big(A_\gamma(\overline{g})\Big)\right) = \delta_{\xi\cap\gamma'}\Big(\widetilde{\Theta}_{(\xi\cap\gamma',\xi)}(\widetilde{g}_\xi)\Big)$$

for every $\xi \cap \gamma' = \eta \in \overline{\eta}(\gamma')$.

This last equation and the definition of the function $A_{\gamma'}$ allow us to conclude that

$$A_{\gamma'}^{-1}\left(\widehat{\Theta}_{(\gamma',\gamma)}\Big(A_\gamma(\overline{g})\Big)\right) = \left[\Big\{\widetilde{\Theta}_{(\xi\cap\gamma',\xi)}(\widetilde{g}_\xi)\Big\}_{\xi\cap\gamma'\in\overline{\eta}(\gamma')}\right]$$

and, hence, by the definition given for $\overline{\Theta}_{(\gamma',\gamma)}(\overline{g})$, we obtain:

$$\overline{\Theta}_{(\gamma',\gamma)}\left(\overline{g} = \left[\Big\{\widetilde{g}_\xi\Big\}_{\xi\in\overline{\xi}(\gamma)}\right]\right) = \left[\Big\{\widetilde{\Theta}_{(\xi\cap\gamma',\xi)}(\widetilde{g}_\xi)\Big\}_{\xi\cap\gamma'\in\overline{\eta}(\gamma')}\right]$$

which is the same expression (4.24-5) obtained in the item 4.24.

It remains now, in order for us to complete our goal of proving that the abuse of language under analysis is, on the considered circumstances, healthy, that is, harmless in the sense of not leading us to errors, to justify the following claims that were assumed and used in the obtainment of this purpose:

**(a)** the function

$$\begin{aligned} A_\gamma : \qquad \overline{G}(\gamma) &\longrightarrow \widehat{G}(\gamma) \\ \overline{g} = \left[\{\widetilde{g}_\xi\}_{\xi\in\overline{\xi}(\gamma)}\right] &\longmapsto A_\gamma(\overline{g}) \end{aligned}$$

where $A_\gamma(\overline{g}) \in \widehat{G}(\gamma)$ is such that

$$\widehat{\Theta}_{(\xi,\gamma)}\Big(A_\gamma(\overline{g})\Big) = \delta_\xi(\widetilde{g}_\xi) \quad \text{for every} \quad \xi \in \overline{\xi}(\gamma),$$

is well-defined;

**(b)** $A_\gamma$ as defined in (a) is such that

$$A_\gamma\Big(b_\gamma(\widetilde{g})\Big) = \delta_\gamma(\widetilde{g}) \quad \text{for every} \quad \widetilde{g} \in \widetilde{G}(\gamma);$$

**(c)** $A_\gamma$ as defined in (a) is an isomorphism.

We verify below each one of the claims (a) to (c).

**(a)** Here it is enough (and necessary) to prove that: if $\{\widetilde{g}_\xi\}_{\xi\in\overline{\xi}(\gamma)}$ is a coherent family in $\widetilde{\mathscr{G}}(\gamma)$, then, the family $\{\delta_\xi(\widetilde{g}_\xi)\}_{\xi\in\overline{\xi}(\gamma)}$ also is coherent in $\widehat{\mathscr{G}}(\gamma)$. Let then $\{\widetilde{g}_\xi\}_{\xi\in\overline{\xi}(\gamma)}$ be coherent in $\widetilde{\mathscr{G}}(\gamma)$ and take $\xi, \xi' \in \overline{\xi}(\gamma)$ such that $\xi \cap \xi' \neq \varnothing$. Since $\widehat{\mathscr{G}}(I)$ is an extension of $\widetilde{\mathscr{G}}(I)$, then, for any $(\gamma', \gamma) \in \Delta(I)$,

$$\widehat{\Theta}_{(\gamma',\gamma)}\Big(\delta_\gamma(\widetilde{g})\Big) = \delta_{\gamma'}\Big(\widetilde{\Theta}_{(\gamma',\gamma)}(\widetilde{g})\Big) \quad \text{for every} \quad \widetilde{g} \in \widetilde{G}(\gamma)$$

and, hence, taking, first, $\gamma' = \xi \cap \xi'$, $\gamma = \xi$ and $\widetilde{g} = \widetilde{g}_\xi$ and, then, $\gamma' = \xi \cap \xi'$, $\gamma = \xi'$ and $\widetilde{g} = \widetilde{g}_{\xi'}$, one obtains that

$$\widehat{\Theta}_{(\xi\cap\xi',\xi)}\Big(\delta_\xi(\widetilde{g}_\xi)\Big) = \delta_{\xi\cap\xi'}\Big(\widetilde{\Theta}_{(\xi\cap\xi',\xi)}(\widetilde{g}_\xi)\Big)$$

and

$$\widehat{\Theta}_{(\xi\cap\xi',\xi')}\Big(\delta_{\xi'}(\widetilde{g}_{\xi'})\Big) = \delta_{\xi\cap\xi'}\Big(\widetilde{\Theta}_{(\xi\cap\xi',\xi')}(\widetilde{g}_{\xi'})\Big).$$

Now, since $\{\widetilde{g}_\xi\}_{\xi\in\overline{\xi}(\gamma)}$ is coherent in $\widetilde{\mathscr{G}}(\gamma)$, we have

$$\widetilde{\Theta}_{(\xi\cap\xi',\xi)}(\widetilde{g}_\xi) = \widetilde{\Theta}_{(\xi\cap\xi',\xi')}(\widetilde{g}_{\xi'})$$

and, with this, we conclude that

$$\widehat{\Theta}_{(\xi\cap\xi',\xi)}\Big(\delta_\xi(\widetilde{g}_\xi)\Big) = \widehat{\Theta}_{(\xi\cap\xi',\xi')}\Big(\delta_{\xi'}(\widetilde{g}_{\xi'})\Big),$$

which, in turn, tells us that the family $\{\delta_\xi(\widetilde{g}_\xi)\}_{\xi\in\overline{\xi}(\gamma)}$ is coherent in $\widehat{\mathscr{G}}(\gamma)$.

**(b)** We want to prove that

$$A_\gamma\Big(b_\gamma(\widetilde{g})\Big) = \delta_\gamma(\widetilde{g}) \quad \text{for every} \quad \widetilde{g} \in \widetilde{G}(\gamma).$$

Let then $\widetilde{g} \in \widetilde{G}(\gamma)$ arbitrarily fixed. Keeping in mind the definition of the function $b_\gamma : \widetilde{G}(\gamma) \longrightarrow \overline{G}(\gamma)$ (review Proposition 4.18), we have that

$$b_\gamma(\widetilde{g}) = \Big[\Big\{\widetilde{\Theta}_{(\xi,\gamma)}(\widetilde{g})\Big\}_{\xi\in\overline{\xi}(\gamma)}\Big].$$

Now, by the definition of $A_\gamma$ (as given in (a)),

$$A_\gamma\Big(b_\gamma(\widetilde{g})\Big) = \widehat{h}$$

where $\widehat{h} \in \widehat{G}(\gamma)$ is such that

$$\widehat{\Theta}_{(\xi,\gamma)}(\widehat{h}) = \delta_\xi\Big(\widetilde{\Theta}_{(\xi,\gamma)}(\widetilde{g})\Big) \quad \text{for every} \quad \xi \in \overline{\xi}(\gamma).$$

But, on the other hand, since $\widehat{\mathscr{G}}(I)$ is an extension of $\widetilde{\mathscr{G}}(I)$, we have

$$\widehat{\Theta}_{(\xi,\gamma)}\Big(\delta_\gamma(\widetilde{g})\Big) = \delta_\xi\Big(\widetilde{\Theta}_{(\xi,\gamma)}(\widetilde{g})\Big) \quad \text{for every} \quad \xi \in \overline{\xi}(\gamma)$$

and, hence, from the last two equations, we conclude that

$$\widehat{\Theta}_{(\xi,\gamma)}(\widehat{h}) = \widehat{\Theta}_{(\xi,\gamma)}\Big(\delta_\gamma(\widetilde{g})\Big) \quad \text{for every} \quad \xi \in \overline{\xi}(\gamma),$$

from where one obtains, taking Lemma 3.32 into account, that

$$\widehat{h} = \delta_\gamma(\widetilde{g}),$$

that is,

$$A_\gamma\Big(b_\gamma(\widetilde{g})\Big) = \delta_\gamma(\widetilde{g}).$$

**(c)** Let us first prove that the function $A_\gamma$ defined in (a) is surjective.

We know that $\widehat{\mathscr{G}}(I)$ is a locally closed extension of $\widetilde{\mathscr{G}}(I)$, which, by Definition 3.23(c-2) expressed without the abuse of language (see Warning 3.24), means that: for each $\gamma \in \Gamma(I)$, if $\widehat{g} \in \widehat{G}(\gamma)$ and $x \in \gamma$, there exists $\gamma' \in \Gamma(\gamma)$ such that $x \in \gamma'$ and

$$\widehat{\Theta}_{(\gamma',\gamma)}(\widehat{g}) \in \delta_{\gamma'}\Big(\widetilde{G}(\gamma')\Big).$$

Let now $\widehat{g} \in \widehat{G}(\gamma)$ be arbitrarily fixed and $\overline{\xi}(\gamma) \subseteq \Gamma(\gamma)$ be a cover of $\gamma$ such that[20]:

$$\widehat{\Theta}_{(\xi,\gamma)}(\widehat{g}) \in \delta_\xi\Big(\widetilde{G}(\xi)\Big) \quad \text{for every} \quad \xi \in \overline{\xi}(\gamma).$$

Hence, for each $\xi \in \overline{\xi}(\gamma)$ there exists $\widetilde{g}_\xi \in \widetilde{G}(\xi)$ such that

$$\widehat{\Theta}_{(\xi,\gamma)}(\widehat{g}) = \delta_\xi(\widetilde{g}_\xi).$$

Let us consider the family $\widetilde{g}(\overline{\xi}(\gamma))$ defined as follows:

$$\widetilde{g}\Big(\overline{\xi}(\gamma)\Big) := \Big\{\widetilde{g}_\xi = \delta_\xi^{-1}\Big(\widehat{\Theta}_{(\xi,\gamma)}(\widehat{g})\Big)\Big\}_{\xi\in\overline{\xi}(\gamma)}.$$

The family $\widetilde{g}(\overline{\xi}(\gamma))$ is coherent in $\widetilde{\mathscr{G}}(\gamma)$, that is, $\widetilde{g}(\overline{\xi}(\gamma)) \in \widetilde{F}(\gamma)$. In fact, remembering that (see our "starting point" at page 189), for every $(\gamma',\gamma) \in \Delta(I)$,

$$\widehat{\Theta}_{(\gamma',\gamma)}\Big(\delta_\gamma(\widetilde{g})\Big) = \delta_{\gamma'}\Big(\widetilde{\Theta}_{(\gamma',\gamma)}(\widetilde{g})\Big) \quad \text{for every} \quad \widetilde{g} \in \widetilde{G}(\gamma),$$

[20] One such cover $\overline{\xi}(\gamma)$ does exist, according to the footnote 19 at page 146.

we have, for $\xi, \xi' \in \overline{\xi}(\gamma)$ such that $\xi \cap \xi' \neq \emptyset$,

$$\widehat{\Theta}_{(\xi\cap\xi',\xi)}\left(\widehat{\Theta}_{(\xi,\gamma)}(\widehat{g})\right) = \widehat{\Theta}_{(\xi\cap\xi',\xi)}\left(\delta_{\xi}(\widetilde{g}_{\xi})\right) = \delta_{\xi\cap\xi'}\left(\widetilde{\Theta}_{(\xi\cap\xi',\xi)}(\widetilde{g}_{\xi})\right)$$

and

$$\widehat{\Theta}_{(\xi\cap\xi',\xi')}\left(\widehat{\Theta}_{(\xi',\gamma)}(\widehat{g})\right) = \widehat{\Theta}_{(\xi\cap\xi',\xi')}\left(\delta_{\xi'}(\widetilde{g}_{\xi'})\right) = \delta_{\xi\cap\xi'}\left(\widetilde{\Theta}_{(\xi\cap\xi',\xi')}(\widetilde{g}_{\xi'})\right).$$

But,

$$\widehat{\Theta}_{(\xi\cap\xi',\xi)}\left(\widehat{\Theta}_{(\xi,\gamma)}(\widehat{g})\right) = \widehat{\Theta}_{(\xi\cap\xi',\xi')}\left(\widehat{\Theta}_{(\xi',\gamma)}(\widehat{g})\right) = \widehat{\Theta}_{(\xi\cap\xi',\gamma)}(\widehat{g}),$$

and, hence, we have

$$\delta_{\xi\cap\xi'}\left(\widetilde{\Theta}_{(\xi\cap\xi',\xi)}(\widetilde{g}_{\xi})\right) = \delta_{\xi\cap\xi'}\left(\widetilde{\Theta}_{(\xi\cap\xi',\xi')}(\widetilde{g}_{\xi'})\right)$$

from where we conclude, taking into account that $\delta_{\xi\cap\xi'} : \widetilde{G}(\xi\cap\xi') \longrightarrow \widehat{G}(\xi\cap\xi')$ is an injective homomorphism, that

$$\widetilde{\Theta}_{(\xi\cap\xi',\xi)}(\widetilde{g}_{\xi}) = \widetilde{\Theta}_{(\xi\cap\xi',\xi')}(\widetilde{g}_{\xi'}),$$

that is,

$$\widetilde{g}\left(\overline{\xi}(\gamma)\right) = \left\{\widetilde{g}_{\xi} = \delta_{\xi}^{-1}\left(\widehat{\Theta}_{(\xi,\gamma)}(\widehat{g})\right)\right\}_{\xi\in\overline{\xi}(\gamma)} \in \widetilde{F}(\gamma).$$

From this, we can define $\overline{g} \in \overline{G}(\gamma)$ as

$$\overline{g} := \left[\widetilde{g}\left(\overline{\xi}(\gamma)\right)\right].$$

Hence, by the definition of $A_{\gamma}$, we get that $A_{\gamma}(\overline{g}) \in \widehat{G}(\gamma)$ is such that

$$\widehat{\Theta}_{(\xi,\gamma)}\left(A_{\gamma}(\overline{g})\right) = \delta_{\xi}(\widetilde{g}_{\xi}) = \delta_{\xi}\left(\delta_{\xi}^{-1}\left(\widehat{\Theta}_{(\xi,\gamma)}(\widehat{g})\right)\right) = \widehat{\Theta}_{(\xi,\gamma)}(\widehat{g})$$

for every $\xi \in \overline{\xi}(\gamma)$. From this, by Lemma 3.32, follows that

$$A_{\gamma}(\overline{g}) = \widehat{g},$$

which concludes the proof that $A_{\gamma}$ is a surjective function.

Let us now prove that $A_{\gamma} : \overline{G}(\gamma) \longrightarrow \widehat{G}(\gamma)$ is injective. In order to do so, let

$$\overline{g} = \left[\left\{\widetilde{g}_{\xi}\right\}_{\xi\in\overline{\xi}(\gamma)}\right] \in \overline{G}(\gamma) \quad \text{and} \quad \overline{h} = \left[\left\{\widetilde{h}_{\eta}\right\}_{\eta\in\overline{\eta}(\gamma)}\right] \in \overline{G}(\gamma)$$

be arbitrarily chosen in $\overline{G}(\gamma)$. Hence, $A_{\gamma}(\overline{g}), A_{\gamma}(\overline{h}) \in \widehat{G}(\gamma)$ are such that

$$\widehat{\Theta}_{(\xi,\gamma)}\left(A_{\gamma}(\overline{g})\right) = \delta_{\xi}(\widetilde{g}_{\xi}) \quad \text{for every} \quad \xi \in \overline{\xi}(\gamma)$$

and

$$\widehat{\Theta}_{(\eta,\gamma)}\left(A_{\gamma}(\overline{h})\right) = \delta_{\eta}(\widetilde{h}_{\eta}) \quad \text{for every} \quad \eta \in \overline{\eta}(\gamma).$$

Let now $\overline{\chi}(\gamma) \subseteq \Gamma(\gamma)$ be the following cover of $\gamma$:

$$\overline{\chi}(\gamma) := \Big\{\chi := \xi \cap \eta \quad : \quad \xi \in \overline{\xi}(\gamma), \quad \eta \in \overline{\eta}(\gamma) \quad \text{and} \quad \xi \cap \eta \neq \varnothing\Big\}.$$

For an arbitrary $\xi \cap \eta \in \overline{\chi}(\gamma)$ we calculate:

$$\widehat{\Theta}_{(\xi\cap\eta,\gamma)}\Big(A_\gamma(\overline{g})\Big) = \widehat{\Theta}_{(\xi\cap\eta,\xi)}\left(\widehat{\Theta}_{(\xi,\gamma)}\Big(A_\gamma(\overline{g})\Big)\right) = \widehat{\Theta}_{(\xi\cap\eta,\xi)}\Big(\delta_\xi(\widetilde{g}_\xi)\Big) =$$
$$= \delta_{\xi\cap\eta}\Big(\widetilde{\Theta}_{(\xi\cap\eta,\xi)}(\widetilde{g}_\xi)\Big)$$

and

$$\widehat{\Theta}_{(\xi\cap\eta,\gamma)}\Big(A_\gamma(\overline{h})\Big) = \widehat{\Theta}_{(\xi\cap\eta,\eta)}\left(\widehat{\Theta}_{(\eta,\gamma)}\Big(A_\gamma(\overline{h})\Big)\right) = \widehat{\Theta}_{(\xi\cap\eta,\eta)}\Big(\delta_\eta(\widetilde{h}_\eta)\Big) =$$
$$= \delta_{\xi\cap\eta}\Big(\widetilde{\Theta}_{(\xi\cap\eta,\eta)}(\widetilde{h}_\eta)\Big).$$

If $A_\gamma(\overline{g}) = A_\gamma(\overline{h})$, it results from the expressions above that

$$\delta_{\xi\cap\eta}\Big(\widetilde{\Theta}_{(\xi\cap\eta,\xi)}(\widetilde{g}_\xi)\Big) = \delta_{\xi\cap\eta}\Big(\widetilde{\Theta}_{(\xi\cap\eta,\eta)}(\widetilde{h}_\eta)\Big)$$

that is,

$$\widetilde{\Theta}_{(\xi\cap\eta,\xi)}(\widetilde{g}_\xi) = \widetilde{\Theta}_{(\xi\cap\eta,\eta)}(\widetilde{h}_\eta) \quad \text{for every} \quad \xi \cap \eta \in \overline{\chi}(\gamma).$$

Thus,

$$\Big\{\widetilde{g}_\xi\Big\}_{\xi\in\overline{\xi}(\gamma)} \approx \Big\{\widetilde{h}_\eta\Big\}_{\eta\in\overline{\eta}(\gamma)}$$

and, therefore,

$$\overline{g} = \left[\Big\{\widetilde{g}_\xi\Big\}_{\xi\in\overline{\xi}(\gamma)}\right] = \left[\Big\{\widetilde{h}_\eta\Big\}_{\eta\in\overline{\eta}(\gamma)}\right] = \overline{h},$$

that is, $A_\gamma$ is injective.

It remains now to prove that $A_\gamma$ is a homomorphism, that is,

$$A_\gamma(\overline{g} + \overline{h}) = A_\gamma(\overline{g}) + A_\gamma(\overline{h})$$

for every $\overline{g}, \overline{h} \in \overline{G}(\gamma)$.

Let then

$$\overline{g} = \left[\Big\{\widetilde{g}_\xi\Big\}_{\xi\in\overline{\xi}(\gamma)}\right] \in \overline{G}(\gamma) \quad \text{and} \quad \overline{h} = \left[\Big\{\widetilde{h}_\eta\Big\}_{\eta\in\overline{\eta}(\gamma)}\right] \in \overline{G}(\gamma)$$

be arbitrarily fixed. Hence, taking into account the definition of the function $A_\gamma : \overline{G}(\gamma) \longrightarrow \widehat{G}(\gamma)$ as well as the addition in $\overline{G}(\gamma)$ (see Proposition 4.15), we have:

$$\overline{g} + \overline{h} = \left[\Big\{\widetilde{\Theta}_{(\xi\cap\eta,\xi)}(\widetilde{g}_\xi) + \widetilde{\Theta}_{(\xi\cap\eta,\eta)}(\widetilde{h}_\eta)\Big\}_{\xi\cap\eta\in\overline{\chi}(\gamma)}\right]$$

where $\overline{\chi}(\gamma) := \{\chi := \xi \cap \eta \quad : \quad \xi \in \overline{\xi}(\gamma), \quad \eta \in \overline{\eta}(\gamma) \quad \text{and} \quad \xi \cap \eta \neq \varnothing\}$,

$$\widehat{\Theta}_{(\xi,\gamma)}\Big(A_\gamma(\overline{g})\Big) = \delta_\xi(\widetilde{g}_\xi) \quad \text{for every} \quad \xi \in \overline{\xi}(\gamma)$$

and

$$\widehat{\Theta}_{(\eta,\gamma)}\Big(A_\gamma(\overline{h})\Big) = \delta_\eta(\widetilde{h}_\eta) \quad \text{for every} \quad \eta \in \overline{\eta}(\gamma).$$

Now, for an arbitrary $\xi \cap \eta \in \overline{\chi}(\gamma)$, we calculate:

$$\begin{aligned}
\widehat{\Theta}_{(\xi\cap\eta,\gamma)}\Big(A_\gamma(\overline{g}+\overline{h})\Big) &= \delta_{\xi\cap\eta}\Big(\widetilde{\Theta}_{(\xi\cap\eta,\xi)}(\widetilde{g}_\xi) + \widetilde{\Theta}_{(\xi\cap\eta,\eta)}(\widetilde{h}_\eta)\Big) = \\
&= \delta_{\xi\cap\eta}\Big(\widetilde{\Theta}_{(\xi\cap\eta,\xi)}(\widetilde{g}_\xi)\Big) + \delta_{\xi\cap\eta}\Big(\widetilde{\Theta}_{(\xi\cap\eta,\eta)}(\widetilde{h}_\eta)\Big) = \\
&= \widehat{\Theta}_{(\xi\cap\eta,\xi)}\Big(\delta_\xi(\widetilde{g}_\xi)\Big) + \widehat{\Theta}_{(\xi\cap\eta,\eta)}\Big(\delta_\eta(\widetilde{h}_\eta)\Big) = \\
&= \widehat{\Theta}_{(\xi\cap\eta,\xi)}\bigg(\widehat{\Theta}_{(\xi,\gamma)}\Big(A_\gamma(\overline{g})\Big)\bigg) + \widehat{\Theta}_{(\xi\cap\eta,\eta)}\bigg(\widehat{\Theta}_{(\eta,\gamma)}\Big(A_\gamma(\overline{h})\Big)\bigg).
\end{aligned}$$

Hence, we have that:

$$\widehat{\Theta}_{(\xi\cap\eta,\gamma)}\Big(A_\gamma(\overline{g}+\overline{h})\Big) = \widehat{\Theta}_{(\xi\cap\eta,\gamma)}\Big(A_\gamma(\overline{g})\Big) + \widehat{\Theta}_{(\xi\cap\eta,\gamma)}\Big(A_\gamma(\overline{h})\Big) = \widehat{\Theta}_{(\xi\cap\eta,\gamma)}\Big(A_\gamma(\overline{g}) + A_\gamma(\overline{h})\Big)$$

for every $\xi \cap \eta \in \overline{\chi}(\gamma)$. From here one concludes, once more resorting to Lemma 3.32, that

$$A_\gamma(\overline{g}+\overline{h}) = A_\gamma(\overline{g}) + A_\gamma(\overline{h}).$$

We return now, after this long item regarding the harmlessness of the abuse of language used in the obtainment (in item 4.24) of the expression (4.24-5) for $\overline{\Theta}_{(\gamma',\gamma)}(\overline{g})$, to the question formulated at the end of item 4.24, namely: If we define $\overline{\Theta}_{(\gamma',\gamma)} : \overline{G}(\gamma) \longrightarrow \overline{G}(\gamma')$ through the expression (4.24-5), would we have a well-defined function which is an "extension" of $\widetilde{\Theta}_{(\gamma',\gamma)} : \widetilde{G}(\gamma) \longrightarrow \widetilde{G}(\gamma')$? As we promised there, the answer to this question is provided by the following proposition.

## 4.26 Proposition

*Let*

$$\bigg(\overline{\mathbb{G}}\Big(\Gamma(I)\Big) = \Big\{\overline{G}(\gamma) = \Big(\overline{G}(\gamma), \overline{H}(\gamma)\Big)\Big\}_{\gamma\in\Gamma(I)}, \overline{i}\Big(\Gamma^2(I)\Big)\bigg)$$

*be the extension, described in 4.23, of the bonded family*

$$\bigg(\widetilde{\mathbb{G}}\Big(\Gamma(I)\Big) = \Big\{\widetilde{G}(\gamma) = \Big(\widetilde{G}(\gamma), \widetilde{H}(\gamma)\Big)\Big\}_{\gamma\in\Gamma(I)}, \widetilde{i}\Big(\Gamma^2(I)\Big)\bigg)$$

*of the S-space*

$$\widetilde{\mathscr{G}}(I) = \bigg(\widetilde{\mathbb{G}}\Big(\Gamma(I)\Big), \widetilde{i}\Big(\Gamma^2(I)\Big), \widetilde{\Theta}\Big(\Delta(I)\Big)\bigg)$$

*described in 4.6.*

*Let also, for each* $(\gamma', \gamma) \in \Delta(I)$,

$$\begin{aligned}\overline{\Theta}_{(\gamma',\gamma)} : \overline{G}(\gamma) &\longrightarrow \overline{G}(\gamma') \\ \overline{g} &\longmapsto \overline{\Theta}_{(\gamma',\gamma)}(\overline{g})\end{aligned}$$

*defined by*

$$\overline{\Theta}_{(\gamma',\gamma)}\left(\overline{g} = \left[\left\{\widetilde{g}_\xi\right\}_{\xi\in\overline{\xi}(\gamma)}\right]\right) := \left[\left\{\widetilde{\Theta}_{(\xi\cap\gamma',\xi)}(\widetilde{g}_\xi)\right\}_{\xi\cap\gamma'\in\overline{\eta}(\gamma')}\right],$$

*where* $\overline{\eta}(\gamma') \subseteq \Gamma(\gamma')$ *is the cover of* $\gamma'$ *defined by*

$$\overline{\eta}(\gamma') := \left\{\eta := \xi \cap \gamma' \quad : \quad \xi \in \overline{\xi}(\gamma) \quad \textit{and} \quad \xi \cap \gamma' \neq \varnothing\right\}.$$

*Then,* $\overline{\Theta}_{(\gamma',\gamma)}$ *is a well-defined function and, also, a homomorphism from the group* $\overline{G}(\gamma)$ *into the group* $\overline{G}(\gamma')$ *such that*

$$\overline{\Theta}_{(\gamma',\gamma)}\Big(b_\gamma(\widetilde{g})\Big) = b_{\gamma'}\Big(\widetilde{\Theta}_{(\gamma',\gamma)}(\widetilde{g})\Big)$$

*for every* $\widetilde{g} \in \widetilde{G}(\gamma)$*, being* $b_\gamma$ $(b_{\gamma'})$ *the injective homomorphism defined in Proposition 4.18. Furthermore, under the hypothesis of existence of* $\widehat{\mathscr{G}}(I)$*, one has*

$$\overline{\Theta}_{(\gamma',\gamma)} = A_{\gamma'}^{-1}\ \widehat{\Theta}_{(\gamma',\gamma)}\ A_\gamma,$$

*with* $A_\gamma$ $(A_{\gamma'})$ *as defined in Proposition 4.12 (or, without the abuse of language, in (a) on page 191).*

*Proof.* Let us first prove that $\overline{\Theta}_{(\gamma',\gamma)}$ is well-defined, for what it is enough to show that: if the families

$$\left\{\widetilde{g}_\xi\right\}_{\xi\in\overline{\xi}(\gamma)} \quad \text{and} \quad \left\{\widetilde{h}_\chi\right\}_{\chi\in\overline{\chi}(\gamma)}$$

are coherent in $\widetilde{\mathscr{G}}(\gamma)$ and such that

$$\left\{\widetilde{g}_\xi\right\}_{\xi\in\overline{\xi}(\gamma)} \approx \left\{\widetilde{h}_\chi\right\}_{\chi\in\overline{\chi}(\gamma)}, \tag{4.26-1}$$

then, the families

$$\left\{\widetilde{\Theta}_{(\xi\cap\gamma',\xi)}(\widetilde{g}_\xi)\right\}_{\xi\cap\gamma'\in\overline{\eta}(\gamma')} \quad \text{and} \quad \left\{\widetilde{\Theta}_{(\chi\cap\gamma',\chi)}(\widetilde{h}_\chi)\right\}_{\chi\cap\gamma'\in\overline{\nu}(\gamma')},$$

where $\overline{\eta}(\gamma') \subseteq \Gamma(\gamma')$ and $\overline{\nu}(\gamma') \subseteq \Gamma(\gamma')$ are the covers of $\gamma'$ defined by

$$\overline{\eta}(\gamma') := \left\{\xi \cap \gamma' \quad : \quad \xi \in \overline{\xi}(\gamma) \quad \text{and} \quad \xi \cap \gamma' \neq \varnothing\right\}$$

and

$$\overline{\nu}(\gamma') := \left\{\chi \cap \gamma' \quad : \quad \chi \in \overline{\chi}(\gamma) \quad \text{and} \quad \chi \cap \gamma' \neq \varnothing\right\},$$

are coherent in $\widetilde{\mathscr{G}}(\gamma')$ and

$$\left\{\widetilde{\Theta}_{(\xi\cap\gamma',\xi)}(\widetilde{g}_\xi)\right\}_{\xi\cap\gamma'\in\overline{\eta}(\gamma')} \approx \left\{\widetilde{\Theta}_{(\chi\cap\gamma',\chi)}(\widetilde{h}_\chi)\right\}_{\chi\cap\gamma'\in\overline{\nu}(\gamma')}. \tag{4.26-2}$$

Let us then suppose, first, that, for instance, $\{\widetilde{g}_\xi\}_{\xi\in\overline{\xi}(\gamma)}$ is coherent in $\widetilde{\mathscr{G}}(\gamma)$, i.e., that

$$\widetilde{\Theta}_{(\xi\cap\xi',\xi)}(\widetilde{g}_\xi) = \widetilde{\Theta}_{(\xi\cap\xi',\xi')}(\widetilde{g}_{\xi'}) \tag{4.26-3}$$

for every $\xi,\xi'\in\overline{\xi}(\gamma)$ such that $\xi\cap\xi'\neq\varnothing$, and let us prove that $\{\widetilde{\Theta}_{(\xi\cap\gamma',\xi)}(\widetilde{g}_\xi)\}_{\xi\cap\gamma'\in\overline{\eta}(\gamma')}$ is coherent in $\widetilde{\mathscr{G}}(\gamma')$, i.e., that for every $\xi\cap\gamma',\xi'\cap\gamma'\in\overline{\eta}(\gamma')$ such that $\xi\cap\gamma'\cap\xi'\cap\gamma' = \xi\cap\xi'\cap\gamma'\neq\varnothing$,

$$\widetilde{\Theta}_{(\xi\cap\xi'\cap\gamma',\xi\cap\gamma')}\left(\widetilde{\Theta}_{(\xi\cap\gamma',\xi)}(\widetilde{g}_\xi)\right) = \widetilde{\Theta}_{(\xi\cap\xi'\cap\gamma',\xi'\cap\gamma')}\left(\widetilde{\Theta}_{(\xi'\cap\gamma',\xi')}(\widetilde{g}_{\xi'})\right),$$

or, equivalently,

$$\widetilde{\Theta}_{(\xi\cap\xi'\cap\gamma',\xi)}(\widetilde{g}_\xi) = \widetilde{\Theta}_{(\xi\cap\xi'\cap\gamma',\xi')}(\widetilde{g}_{\xi'}).$$

However, this trivially follows from (4.26-3) and the identity

$$\widetilde{\Theta}_{(\xi\cap\xi'\cap\gamma',\xi)}(\widetilde{g}_\xi) = \widetilde{\Theta}_{(\xi\cap\xi'\cap\gamma',\xi\cap\xi')}\left(\widetilde{\Theta}_{(\xi\cap\xi',\xi)}(\widetilde{g}_\xi)\right).$$

Now, supposing that (4.26-1) occurs, i.e., that

$$\widetilde{\Theta}_{(\xi\cap\chi,\xi)}(\widetilde{g}_\xi) = \widetilde{\Theta}_{(\xi\cap\chi,\chi)}(\widetilde{h}_\chi) \tag{4.26-4}$$

for every $\xi\in\overline{\xi}(\gamma)$ and $\chi\in\overline{\chi}(\gamma)$ such that $\xi\cap\chi\neq\varnothing$, we must prove that (4.26-2) is also true, i.e., that

$$\widetilde{\Theta}_{(\xi\cap\gamma'\cap\chi\cap\gamma',\xi\cap\gamma')}\left(\widetilde{\Theta}_{(\xi\cap\gamma',\xi)}(\widetilde{g}_\xi)\right) = \widetilde{\Theta}_{(\xi\cap\gamma'\cap\chi\cap\gamma',\chi\cap\gamma')}\left(\widetilde{\Theta}_{(\chi\cap\gamma',\chi)}(\widetilde{h}_\chi)\right),$$

that is,

$$\widetilde{\Theta}_{(\xi\cap\chi\cap\gamma',\xi)}(\widetilde{g}_\xi) = \widetilde{\Theta}_{(\xi\cap\chi\cap\gamma',\chi)}(\widetilde{h}_\chi)$$

for every $\xi\cap\gamma'\in\overline{\eta}(\gamma')$ and $\chi\cap\gamma'\in\overline{\nu}(\gamma')$ such that $\xi\cap\gamma'\cap\chi\cap\gamma' = \xi\cap\chi\cap\gamma'\neq\varnothing$.

Well, since

$$\widetilde{\Theta}_{(\xi\cap\chi\cap\gamma',\chi)}(\widetilde{h}_\chi) = \widetilde{\Theta}_{(\xi\cap\chi\cap\gamma',\xi\cap\chi)}\left(\widetilde{\Theta}_{(\xi\cap\chi,\chi)}(\widetilde{h}_\chi)\right)$$

and taking (4.26-4) into account, we get

$$\widetilde{\Theta}_{(\xi\cap\chi\cap\gamma',\chi)}(\widetilde{h}_\chi) = \widetilde{\Theta}_{(\xi\cap\chi\cap\gamma',\xi\cap\chi)}\left(\widetilde{\Theta}_{(\xi\cap\chi,\xi)}(\widetilde{g}_\xi)\right)$$

from where one obtains that

$$\widetilde{\Theta}_{(\xi\cap\chi\cap\gamma',\chi)}(\widetilde{h}_\chi) = \widetilde{\Theta}_{(\xi\cap\chi\cap\gamma',\xi)}(\widetilde{g}_\xi)$$

which concludes the proof of $\overline{\Theta}_{(\gamma',\gamma)}$ being well-defined.

Let us prove now that $\overline{\Theta}_{(\gamma',\gamma)} : \overline{G}(\gamma) \longrightarrow \overline{G}(\gamma')$ is a homomorphism. Let then

$$\overline{g} = \left[\left\{\widetilde{g}_\xi\right\}_{\xi\in\overline{\xi}(\gamma)}\right] \quad \text{and} \quad \overline{h} = \left[\left\{\widetilde{h}_\chi\right\}_{\chi\in\overline{\chi}(\gamma)}\right]$$

be arbitrarily chosen elements in $\overline{G}(\gamma)$. From the definition of addition in $\overline{G}(\gamma)$ follows that

$$\overline{g} + \overline{h} = \left[\left\{\widetilde{\Theta}_{(\xi\cap\chi,\xi)}(\widetilde{g}_\xi) + \widetilde{\Theta}_{(\xi\cap\chi,\chi)}(\widetilde{h}_\chi)\right\}_{\xi\cap\chi\in\overline{\nu}_1(\gamma)}\right],$$

where $\overline{\nu}_1(\gamma) \subseteq \Gamma(\gamma)$ is the cover of $\gamma$ given by

$$\overline{\nu}_1(\gamma) \coloneqq \left\{\xi\cap\chi \quad : \quad \xi\in\overline{\xi}(\gamma), \quad \chi\in\overline{\chi}(\gamma) \quad \text{and} \quad \xi\cap\chi\neq\varnothing\right\}.$$

Hence, from the definition of $\overline{\Theta}_{(\gamma',\gamma)}$, we have

$$\overline{\Theta}_{(\gamma',\gamma)}(\overline{g}+\overline{h}) = \left[\left\{\widetilde{\Theta}_{(\xi\cap\chi\cap\gamma',\xi\cap\chi)}\Big(\widetilde{\Theta}_{(\xi\cap\chi,\xi)}(\widetilde{g}_\xi) + \widetilde{\Theta}_{(\xi\cap\chi,\chi)}(\widetilde{h}_\chi)\Big)\right\}_{\xi\cap\chi\cap\gamma'\in\overline{\eta}_1(\gamma')}\right]$$

where $\overline{\eta}_1(\gamma') \subseteq \Gamma(\gamma')$ is the cover of $\gamma'$ defined as

$$\overline{\eta}_1(\gamma') \coloneqq \left\{\xi\cap\chi\cap\gamma' \quad : \quad \xi\cap\chi\in\overline{\nu}_1(\gamma) \quad \text{and} \quad \xi\cap\chi\cap\gamma'\neq\varnothing\right\}.$$

Since

$$\begin{aligned}
&\widetilde{\Theta}_{(\xi\cap\chi\cap\gamma',\xi\cap\chi)}\Big(\widetilde{\Theta}_{(\xi\cap\chi,\xi)}(\widetilde{g}_\xi) + \widetilde{\Theta}_{(\xi\cap\chi,\chi)}(\widetilde{h}_\chi)\Big) = \\
&= \widetilde{\Theta}_{(\xi\cap\chi\cap\gamma',\xi)}(\widetilde{g}_\xi) + \widetilde{\Theta}_{(\xi\cap\chi\cap\gamma',\chi)}(\widetilde{h}_\chi) = \\
&= \widetilde{\Theta}_{(\xi\cap\chi\cap\gamma',\xi\cap\gamma')}\Big(\widetilde{\Theta}_{(\xi\cap\gamma',\xi)}(\widetilde{g}_\xi)\Big) + \widetilde{\Theta}_{(\xi\cap\chi\cap\gamma',\chi\cap\gamma')}\Big(\widetilde{\Theta}_{(\chi\cap\gamma',\chi)}(\widetilde{h}_\chi)\Big),
\end{aligned}$$

then,

$$\begin{aligned}
&\overline{\Theta}_{(\gamma',\gamma)}(\overline{g}+\overline{h}) = \\
&= \left[\left\{\widetilde{\Theta}_{(\xi\cap\chi\cap\gamma',\xi\cap\gamma')}\Big(\widetilde{\Theta}_{(\xi\cap\gamma',\xi)}(\widetilde{g}_\xi)\Big) + \widetilde{\Theta}_{(\xi\cap\chi\cap\gamma',\chi\cap\gamma')}\Big(\widetilde{\Theta}_{(\chi\cap\gamma',\chi)}(\widetilde{h}_\chi)\Big)\right\}_{\xi\cap\chi\cap\gamma'\in\overline{\eta}_1(\gamma')}\right]
\end{aligned}$$

On the other hand, being $\overline{\alpha}(\gamma') \subseteq \Gamma(\gamma')$ and $\overline{\beta}(\gamma') \subseteq \Gamma(\gamma')$ the covers of $\gamma'$ defined by

$$\overline{\alpha}(\gamma') \coloneqq \left\{\alpha \coloneqq \xi\cap\gamma' \quad : \quad \xi\in\overline{\xi}(\gamma) \quad \text{and} \quad \xi\cap\gamma'\neq\varnothing\right\}$$

and

$$\overline{\beta}(\gamma') \coloneqq \left\{\beta \coloneqq \chi\cap\gamma' \quad : \quad \chi\in\overline{\chi}(\gamma) \quad \text{and} \quad \chi\cap\gamma'\neq\varnothing\right\},$$

then, from the definitions of $\overline{\Theta}_{(\gamma',\gamma)}$ and addition of the group $\overline{G}(\gamma)$, we have:

$$\begin{aligned}
&\overline{\Theta}_{(\gamma',\gamma)}(\overline{g}) + \overline{\Theta}_{(\gamma',\gamma)}(\overline{h}) = \\
&= \left[\left\{\widetilde{\Theta}_{(\xi\cap\gamma',\xi)}(\widetilde{g}_\xi)\right\}_{\xi\cap\gamma'\in\overline{\alpha}(\gamma')}\right] + \left[\left\{\widetilde{\Theta}_{(\chi\cap\gamma',\chi)}(\widetilde{h}_\chi)\right\}_{\chi\cap\gamma'\in\overline{\beta}(\gamma')}\right] = \\
&= \left[\left\{\widetilde{\Theta}_{(\xi\cap\chi\cap\gamma',\xi\cap\gamma')}\Big(\widetilde{\Theta}_{(\xi\cap\gamma',\xi)}(\widetilde{g}_\xi)\Big) + \widetilde{\Theta}_{(\xi\cap\chi\cap\gamma',\chi\cap\gamma')}\Big(\widetilde{\Theta}_{(\chi\cap\gamma',\chi)}(\widetilde{h}_\chi)\Big)\right\}_{\xi\cap\chi\cap\gamma'\in\overline{\eta}_1(\gamma')}\right].
\end{aligned}$$

From the expressions above obtained for $\overline{\Theta}_{(\gamma',\gamma)}(\overline{g}+\overline{h})$ and $\overline{\Theta}_{(\gamma',\gamma)}(\overline{g})+\overline{\Theta}_{(\gamma',\gamma)}(\overline{h})$, it results that

$$\overline{\Theta}_{(\gamma',\gamma)}(\overline{g}+\overline{h}) = \overline{\Theta}_{(\gamma',\gamma)}(\overline{g}) + \overline{\Theta}_{(\gamma',\gamma)}(\overline{h}),$$

that is, $\overline{\Theta}_{(\gamma',\gamma)} : \overline{G}(\gamma) \longrightarrow \overline{G}(\gamma')$ is a homomorphism.

For us to complete the proof of our proposition, it remains to demonstrate that $\overline{\Theta}_{(\gamma',\gamma)}$ is an "extension" of $\widetilde{\Theta}_{(\gamma',\gamma)}$, that is,

$$\overline{\Theta}_{(\gamma',\gamma)}\Big(b_\gamma(\widetilde{g})\Big) = b_{\gamma'}\Big(\widetilde{\Theta}_{(\gamma',\gamma)}(\widetilde{g})\Big)$$

for every $\widetilde{g} \in \widetilde{G}(\gamma)$.

We know, by the definition of $b_\gamma : \widetilde{G}(\gamma) \longrightarrow \overline{G}(\gamma)$, that

$$b_\gamma(\widetilde{g}) = \Big[\Big\{\widetilde{\Theta}_{(\xi,\gamma)}(\widetilde{g})\Big\}_{\xi\in\overline{\xi}(\gamma)}\Big]$$

for every $\widetilde{g} \in \widetilde{G}(\gamma)$, being $\overline{\xi}(\gamma) \subseteq \Gamma(\gamma)$ a cover of $\gamma$. Hence, from the definition of $\overline{\Theta}_{(\gamma',\gamma)}$, we have:

$$\overline{\Theta}_{(\gamma',\gamma)}\Big(b_\gamma(\widetilde{g})\Big) = \Big[\Big\{\widetilde{\Theta}_{(\xi\cap\gamma',\xi)}\Big(\widetilde{\Theta}_{(\xi,\gamma)}(\widetilde{g})\Big)\Big\}_{\xi\cap\gamma'\in\overline{\eta}(\gamma')}\Big] = \Big[\Big\{\widetilde{\Theta}_{(\xi\cap\gamma',\gamma)}(\widetilde{g})\Big\}_{\xi\cap\gamma'\in\overline{\eta}(\gamma')}\Big]$$

being $\overline{\eta}(\gamma') \subseteq \Gamma(\gamma')$ the cover of $\gamma'$ given by

$$\overline{\eta}(\gamma') = \Big\{\xi \cap \gamma' \quad : \quad \xi \in \overline{\xi}(\gamma) \quad \text{and} \quad \xi \cap \gamma' \neq \varnothing\Big\}.$$

On the other hand, resorting again to the definition of, now, $b_{\gamma'} : \widetilde{G}(\gamma') \longrightarrow \overline{G}(\gamma')$, we have:

$$b_{\gamma'}\Big(\widetilde{\Theta}_{(\gamma',\gamma)}(\widetilde{g})\Big) = \Big[\Big\{\widetilde{\Theta}_{(\nu,\gamma')}\Big(\widetilde{\Theta}_{(\gamma',\gamma)}(\widetilde{g})\Big)\Big\}_{\nu\in\overline{\nu}_2(\gamma')}\Big] = \Big[\Big\{\widetilde{\Theta}_{(\nu,\gamma)}(\widetilde{g})\Big\}_{\nu\in\overline{\nu}_2(\gamma')}\Big]$$

where $\overline{\nu}_2(\gamma') \subseteq \Gamma(\gamma')$ is a cover of $\gamma'$.

Now,

$$\Big\{\widetilde{\Theta}_{(\xi\cap\gamma',\gamma)}(\widetilde{g})\Big\}_{\xi\cap\gamma'\in\overline{\eta}(\gamma')} \approx \Big\{\widetilde{\Theta}_{(\nu,\gamma)}(\widetilde{g})\Big\}_{\nu\in\overline{\nu}_2(\gamma')},$$

that is,

$$\widetilde{\Theta}_{(\xi\cap\gamma'\cap\nu,\xi\cap\gamma')}\Big(\widetilde{\Theta}_{(\xi\cap\gamma',\gamma)}(\widetilde{g})\Big) = \widetilde{\Theta}_{(\xi\cap\gamma'\cap\nu,\nu)}\Big(\widetilde{\Theta}_{(\nu,\gamma)}(\widetilde{g})\Big)$$

for every $\xi \cap \gamma' \in \overline{\eta}(\gamma')$ and $\nu \in \overline{\nu}_2(\gamma')$ such that $\xi \cap \gamma' \cap \nu \neq \varnothing$, since $\widetilde{\Theta}_{(\alpha,\beta)}\widetilde{\Theta}_{(\beta,\delta)} = \widetilde{\Theta}_{(\alpha,\delta)}$ whatever $\alpha, \beta, \delta \in \Gamma(I)$ such that $\alpha \subseteq \beta \subseteq \delta$.

Therefore,

$$\Big[\Big\{\widetilde{\Theta}_{(\xi\cap\gamma',\gamma)}(\widetilde{g})\Big\}_{\xi\cap\gamma'\in\overline{\eta}(\gamma')}\Big] = \Big[\Big\{\widetilde{\Theta}_{(\nu,\gamma)}(\widetilde{g})\Big\}_{\nu\in\overline{\nu}_2(\gamma')}\Big]$$

and, hence,

$$\overline{\Theta}_{(\gamma',\gamma)}\Big(b_\gamma(\widetilde{g})\Big) = b_{\gamma'}\Big(\widetilde{\Theta}_{(\gamma',\gamma)}(\widetilde{g})\Big). \qquad \blacksquare$$

## 4.27 Remark

Based on what establishes Proposition 4.26, the following family of functions is well-defined:

$$\overline{\Theta}\Big(\Delta(I)\Big) := \Big\{\overline{\Theta}_{(\gamma',\gamma)}\Big\}_{(\gamma',\gamma)\in\Delta(I)}$$

where, for each $(\gamma',\gamma) \in \Delta(I)$, $\overline{\Theta}_{(\gamma',\gamma)} : \overline{G}(\gamma) \longrightarrow \overline{G}(\gamma')$ is the function defined in the referred proposition. As we saw, for each $(\gamma',\gamma) \in \Delta(I)$, the corresponding member $\overline{\Theta}_{(\gamma',\gamma)} : \overline{G}(\gamma) \longrightarrow \overline{G}(\gamma')$ of $\overline{\Theta}(\Delta(I))$ is a homomorphism of the group $\overline{G}(\gamma)$ into the group $\overline{G}(\gamma')$. Hence, remembering ourselves of the definition of restriction for a bonded family (Definition 3.13), we can say that: the family $\overline{\Theta}(\Delta(I))$ above defined is a restriction for the bonded family

$$\left(\overline{\mathbb{G}}\Big(\Gamma(I)\Big) = \Big\{\overline{\mathbb{G}}(\gamma) = \Big(\overline{G}(\gamma), \overline{H}(\gamma)\Big)\Big\}_{\gamma\in\Gamma(I)}, \overline{i}\Big(\Gamma^2(I)\Big)\right)$$

described in 4.23, if and only if, for every $\gamma, \gamma', \gamma'' \in \Gamma(I)$ such that $\gamma'' \subseteq \gamma' \subseteq \gamma$,

- $\overline{\Theta}_{(\gamma'',\gamma')}\Big(\overline{\Theta}_{(\gamma',\gamma)}(\overline{g})\Big) = \overline{\Theta}_{(\gamma'',\gamma)}(\overline{g})$ for every $\overline{g} \in \overline{G}(\gamma)$,
- $\overline{\Theta}_{(\gamma',\gamma)}\Big(\overline{\Phi}(\overline{g})\Big) = \overline{i}_{(\gamma',\gamma)}(\overline{\Phi})\Big(\overline{\Theta}_{(\gamma',\gamma)}(\overline{g})\Big)$ for every $\overline{\Phi} \in \overline{H}(\gamma)$ and every $\overline{g} \in \overline{G}(\gamma)$ and
- if $\overline{\Phi} \in \overline{H}(\gamma)$, $\overline{g} \in \overline{G}(\gamma)$ and $\overline{i}_{(\gamma',\gamma)}(\overline{\Phi})\Big(\overline{\Theta}_{(\gamma',\gamma)}(\overline{g})\Big) = 0$ for every $\gamma' \in \overline{\Gamma}(\gamma)$, where $\overline{\Gamma}(\gamma) \subseteq \Gamma(\gamma)$ is a cover of $\gamma$, then, $\overline{g} \in \overline{G}(\gamma)_{\overline{\Phi}}$ (the domain of $\overline{\Phi}$).

The last condition above does not impose any restriction in this case, since $\overline{\Phi} \in \overline{H}(\gamma)$ is an endomorphism on $\overline{G}(\gamma)$ and, therefore, any $\overline{g} \in \overline{G}(\gamma)$ belongs to the domain of $\overline{\Phi}$. Hence, the family of homomorphisms $\overline{\Theta}(\Delta(I))$ described above is a restriction for the bonded family $(\overline{G}(\Gamma(I)), \overline{i}(\Gamma^2(I)))$ above referred, if and only if, the first two conditions above are satisfied. This is the content of our next proposition.

## 4.28 Proposition

*Let*

$$\left(\overline{\mathbb{G}}\Big(\Gamma(I)\Big) = \Big\{\overline{\mathbb{G}}(\gamma) = \Big(\overline{G}(\gamma), \overline{H}(\gamma)\Big)\Big\}_{\gamma\in\Gamma(I)}, \overline{i}\Big(\Gamma^2(I)\Big) = \Big\{\overline{i}_{(\gamma',\gamma)}\Big\}_{(\gamma',\gamma)\in\Gamma^2(I)}\right)$$

*be the bonded family described in 4.23 and*

$$\overline{\Theta}\Big(\Delta(I)\Big) = \Big\{\overline{\Theta}_{(\gamma',\gamma)}\Big\}_{(\gamma',\gamma)\in\Delta(I)}$$

*be the family of homomorphisms defined in 4.27.* $\overline{\Theta}(\Delta(I))$ *satisfies, for every* $\gamma, \gamma', \gamma'' \in \Gamma(I)$ *such that* $\gamma'' \subseteq \gamma' \subseteq \gamma$*, the following conditions and, thus, is a restriction for the bonded family* $(\overline{\mathbb{G}}(\Gamma(I)), \overline{i}(\Gamma^2(I)))$.

**(a)** $\overline{\Theta}_{(\gamma'',\gamma')}\Big(\overline{\Theta}_{(\gamma',\gamma)}(\overline{g})\Big) = \overline{\Theta}_{(\gamma'',\gamma)}(\overline{g})$ *for every* $\overline{g} \in \overline{G}(\gamma)$*;*

**(b)** $\overline{\Theta}_{(\gamma',\gamma)}\Big(\overline{\Phi}(\overline{g})\Big) = \overline{i}_{(\gamma',\gamma)}(\overline{\Phi})\Big(\overline{\Theta}_{(\gamma',\gamma)}(\overline{g})\Big)$ *for every* $\overline{\Phi} \in \overline{H}(\gamma)$ *and every* $\overline{g} \in \overline{G}(\gamma)$*.*

*Proof.* **(a)** Let $\overline{g} = [\{\widetilde{g}_\xi\}_{\xi\in\overline{\xi}(\gamma)}] \in \overline{G}(\gamma)$ be arbitrarily chosen. By the definition of $\overline{\Theta}_{(\gamma',\gamma)}$ we have:

$$\overline{\Theta}_{(\gamma',\gamma)}(\overline{g}) = \Big[\Big\{\widetilde{h}_{\xi\cap\gamma'}\Big\}_{\xi\cap\gamma'\in\overline{\eta}(\gamma')}\Big]$$

being $\overline{\eta}(\gamma') \subseteq \Gamma(\gamma')$ the cover of $\gamma'$ given by

$$\overline{\eta}(\gamma') = \Big\{\xi \cap \gamma' \quad : \quad \xi \in \overline{\xi}(\gamma) \quad \text{and} \quad \xi \cap \gamma' \neq \varnothing\Big\}$$

and

$$\widetilde{h}_{\xi\cap\gamma'} = \widetilde{\Theta}_{(\xi\cap\gamma',\xi)}(\widetilde{g}_\xi) \quad \text{for every} \quad \xi \cap \gamma' \in \overline{\eta}(\gamma').$$

Hence, resorting now to the definition of $\overline{\Theta}_{(\gamma'',\gamma')}$, we have:

$$\overline{\Theta}_{(\gamma'',\gamma')}\Big(\overline{\Theta}_{(\gamma',\gamma)}(\overline{g})\Big) = \Big[\Big\{\widetilde{\Theta}_{(\xi\cap\gamma'\cap\gamma'',\xi\cap\gamma')}(\widetilde{h}_{\xi\cap\gamma'})\Big\}_{\xi\cap\gamma'\cap\gamma''\in\overline{\nu}(\gamma'')}\Big]$$

where $\overline{\nu}(\gamma'') \subseteq \Gamma(\gamma'')$ is the cover of $\gamma''$ given by

$$\overline{\nu}(\gamma'') := \Big\{\xi \cap \gamma' \cap \gamma'' = \xi \cap \gamma'' \quad : \quad \xi \in \overline{\xi}(\gamma) \quad \text{and} \quad \xi \cap \gamma'' \neq \varnothing\Big\}.$$

Taking now into account the expression of $\widetilde{h}_{\xi\cap\gamma'}$, we get:

$$\widetilde{\Theta}_{(\xi\cap\gamma'\cap\gamma'',\xi\cap\gamma')}(\widetilde{h}_{\xi\cap\gamma'}) = \widetilde{\Theta}_{(\xi\cap\gamma'',\xi\cap\gamma')}(\widetilde{h}_{\xi\cap\gamma'}) = \widetilde{\Theta}_{(\xi\cap\gamma'',\xi\cap\gamma')}\Big(\widetilde{\Theta}_{(\xi\cap\gamma',\xi)}(\widetilde{g}_\xi)\Big) = \widetilde{\Theta}_{(\xi\cap\gamma'',\xi)}(\widetilde{g}_\xi).$$

Hence, we obtain that:

$$\overline{\Theta}_{(\gamma'',\gamma')}\Big(\overline{\Theta}_{(\gamma',\gamma)}(\overline{g})\Big) = \Big[\Big\{\widetilde{\Theta}_{(\xi\cap\gamma'',\xi)}(\widetilde{g}_\xi)\Big\}_{\xi\cap\gamma''\in\overline{\nu}(\gamma'')}\Big].$$

On the other hand, directly from the definition of $\overline{\Theta}_{(\gamma'',\gamma)}$ follows that

$$\overline{\Theta}_{(\gamma'',\gamma)}(\overline{g}) = \Big[\Big\{\widetilde{\Theta}_{(\xi\cap\gamma'',\xi)}(\widetilde{g}_\xi)\Big\}_{\xi\cap\gamma''\in\overline{\nu}(\gamma'')}\Big]$$

and, hence, we conclude that

$$\overline{\Theta}_{(\gamma'',\gamma')}\Big(\overline{\Theta}_{(\gamma',\gamma)}(\overline{g})\Big) = \overline{\Theta}_{(\gamma'',\gamma)}(\overline{g}).$$

**(b)** Let, again, $\overline{g} = [\{\widetilde{g}_\xi\}_{\xi\in\overline{\xi}(\gamma)}]$ be arbitrarily chosen in $\overline{G}(\gamma)$. Let also $\overline{\Phi} \in \overline{H}(\gamma)$. From the definition of $\overline{\Phi}$ (see Proposition 4.20) we have:

$$\overline{\Phi}(\overline{g}) = \Big[\Big\{\widetilde{i}_{(\xi,\gamma)}(\widetilde{\Phi})(\widetilde{g}_\xi)\Big\}_{\xi\in\overline{\xi}(\gamma)}\Big]$$

where

$$\widetilde{\Phi} = (\overline{\gamma})^{-1}(\overline{\Phi}),$$

being $\overline{\gamma} : \widetilde{H}(\gamma) \longrightarrow \overline{H}(\gamma)$ the isomorphism defined in Proposition 4.22.

Keeping now in mind the definition of $\overline{\Theta}_{(\gamma',\gamma)}$, one obtains:

$$\overline{\Theta}_{(\gamma',\gamma)}\Big(\overline{\Phi}(\overline{g})\Big) = \Big[\Big\{\widetilde{\Theta}_{(\xi\cap\gamma',\xi)}\Big(\widetilde{i}_{(\xi,\gamma)}(\widetilde{\Phi})(\widetilde{g}_\xi)\Big)\Big\}_{\xi\cap\gamma'\in\overline{\eta}(\gamma')}\Big]$$

where $\overline{\eta}(\gamma') \subseteq \Gamma(\gamma')$ is the cover of $\gamma'$ given by

$$\overline{\eta}(\gamma') = \Big\{\xi\cap\gamma' \quad : \quad \xi\in\overline{\xi}(\gamma) \quad \text{and} \quad \xi\cap\gamma' \neq \varnothing\Big\}.$$

But, as we know,

$$\widetilde{\Theta}_{(\xi\cap\gamma',\xi)}\Big(\widetilde{i}_{(\xi,\gamma)}(\widetilde{\Phi})(\widetilde{g}_\xi)\Big) = \widetilde{i}_{(\xi\cap\gamma',\xi)}\Big(\widetilde{i}_{(\xi,\gamma)}(\widetilde{\Phi})\Big)\Big(\widetilde{\Theta}_{(\xi\cap\gamma',\xi)}(\widetilde{g}_\xi)\Big) = \widetilde{i}_{(\xi\cap\gamma',\gamma)}(\widetilde{\Phi})\Big(\widetilde{\Theta}_{(\xi\cap\gamma',\xi)}(\widetilde{g}_\xi)\Big).$$

Thus,

$$\overline{\Theta}_{(\gamma',\gamma)}\Big(\overline{\Phi}(\overline{g})\Big) = \Big[\Big\{\widetilde{i}_{(\xi\cap\gamma',\gamma)}(\widetilde{\Phi})\Big(\widetilde{\Theta}_{(\xi\cap\gamma',\xi)}(\widetilde{g}_\xi)\Big)\Big\}_{\xi\cap\gamma'\in\overline{\eta}(\gamma')}\Big].$$

We also know (review Definition 3.22) that

$$\widetilde{i}_{(\gamma',\gamma)}(\widetilde{\Phi}) = (\overline{\gamma'})^{-1}\left(\overline{i}_{(\gamma',\gamma)}\Big(\overline{\gamma}(\widetilde{\Phi})\Big)\right) = (\overline{\gamma'})^{-1}\Big(\overline{i}_{(\gamma',\gamma)}(\overline{\Phi})\Big)$$

and, hence, resorting first to the definition of $\overline{\Theta}_{(\gamma',\gamma)}$ in 4.26 (to calculate $\overline{\Theta}_{(\gamma',\gamma)}(\overline{g})$) and, next, to the definition given in Proposition 4.20 (to calculate $\overline{i}_{(\gamma',\gamma)}(\overline{\Phi})(\overline{\Theta}_{(\gamma',\gamma)}(\overline{g}))$ where, then, we will employ the information above highlighted, $\widetilde{i}_{(\gamma',\gamma)}(\widetilde{\Phi}) = (\overline{\gamma'})^{-1}(\overline{i}_{(\gamma',\gamma)}(\overline{\Phi}))$), we get:

$$\begin{aligned}
\overline{i}_{(\gamma',\gamma)}(\overline{\Phi})\Big(\overline{\Theta}_{(\gamma',\gamma)}(\overline{g})\Big) &= \overline{i}_{(\gamma',\gamma)}(\overline{\Phi})\left(\Big[\Big\{\widetilde{\Theta}_{(\xi\cap\gamma',\xi)}(\widetilde{g}_\xi)\Big\}_{\xi\cap\gamma'\in\overline{\eta}(\gamma')}\Big]\right) = \\
&= \Big[\Big\{\widetilde{i}_{(\xi\cap\gamma',\gamma')}\Big(\widetilde{i}_{(\gamma',\gamma)}(\widetilde{\Phi})\Big)\Big(\widetilde{\Theta}_{(\xi\cap\gamma',\xi)}(\widetilde{g}_\xi)\Big)\Big\}_{\xi\cap\gamma'\in\overline{\eta}(\gamma')}\Big] = \\
&= \Big[\Big\{\widetilde{i}_{(\xi\cap\gamma',\gamma)}(\widetilde{\Phi})\Big(\widetilde{\Theta}_{(\xi\cap\gamma',\xi)}(\widetilde{g}_\xi)\Big)\Big\}_{\xi\cap\gamma'\in\overline{\eta}(\gamma')}\Big].
\end{aligned}$$

Taking now this result into the expression above obtained for $\overline{\Theta}_{(\gamma',\gamma)}(\overline{\Phi}(\overline{g}))$, we conclude that

$$\overline{\Theta}_{(\gamma',\gamma)}\Big(\overline{\Phi}(\overline{g})\Big) = \overline{i}_{(\gamma',\gamma)}(\overline{\Phi})\Big(\overline{\Theta}_{(\gamma',\gamma)}(\overline{g})\Big),$$

which completes the proof of the proposition. ■

# The Extension $\overline{\mathscr{G}}(I)$ of $\widetilde{\mathscr{G}}(I)$

## 4.29 The $S$-Space $\overline{\mathscr{G}}(I)$

Let

$$\left(\overline{\mathbb{G}}\Big(\Gamma(I)\Big) = \Big\{\overline{\mathbb{G}}(\gamma) = \Big(\overline{G}(\gamma), \overline{H}(\gamma)\Big)\Big\}_{\gamma\in\Gamma(I)}, \overline{i}\Big(\Gamma^2(I)\Big) = \Big\{\overline{i}_{(\gamma',\gamma)}\Big\}_{(\gamma',\gamma)\in\Gamma^2(I)}\right)$$

be as described in 4.23. As we saw in the referred item, the ordered pair above is a bonded family (Definition 3.11). The Proposition 4.28 establishes that the family of homomorphisms

$$\overline{\Theta}\Big(\Delta(I)\Big) = \Big\{\overline{\Theta}_{(\gamma',\gamma)}\Big\}_{(\gamma',\gamma)\in\Delta(I)},$$

defined in 4.27, is a restriction for the bonded family above (Definition 3.13). Hence being, and taking into account the definition of $S$-space (Definition 3.15), the triplet $\overline{\mathscr{G}}(I)$ defined below is a $S$-space:

$$\overline{\mathscr{G}}(I) \coloneqq \Bigg(\overline{\mathbb{G}}\Big(\Gamma(I)\Big) = \Big\{\overline{\mathbb{G}}(\gamma) = \Big(\overline{G}(\gamma), \overline{H}(\gamma)\Big)\Big\}_{\gamma\in\Gamma(I)},$$
$$\overline{i}\Big(\Gamma^2(I)\Big) = \Big\{\overline{i}_{(\gamma',\gamma)}\Big\}_{(\gamma',\gamma)\in\Gamma^2(I)},$$
$$\overline{\Theta}\Big(\Delta(I)\Big) = \Big\{\overline{\Theta}_{(\gamma',\gamma)}\Big\}_{(\gamma',\gamma)\in\Delta(I)}\Bigg).$$

We also saw in 4.23 that the bonded family of the $S$-space $\overline{\mathscr{G}}(I)$, namely,

$$\Bigg(\overline{\mathbb{G}}\Big(\Gamma(I)\Big), \overline{i}\Big(\Gamma^2(I)\Big)\Bigg),$$

is an extension of the bonded family (Definition 3.22)

$$\Bigg(\widetilde{\mathbb{G}}\Big(\Gamma(I)\Big), \widetilde{i}\Big(\Gamma^2(I)\Big)\Bigg)$$

of the $S$-space

$$\widetilde{\mathscr{G}}(I) = \Bigg(\widetilde{\mathbb{G}}\Big(\Gamma(I)\Big), \widetilde{i}\Big(\Gamma^2(I)\Big), \widetilde{\Theta}\Big(\Delta(I)\Big)\Bigg)$$

defined in 4.6.

Since, from Proposition 4.26, the members of the families $\overline{\Theta}(\Delta(I))$ and $\widetilde{\Theta}(\Delta(I))$ are such that, for each $(\gamma',\gamma) \in \Delta(I)$,

$$\overline{\Theta}_{(\gamma',\gamma)}\Big(b_\gamma(\tilde{g})\Big) = b_{\gamma'}\Big(\widetilde{\Theta}_{(\gamma',\gamma)}(\tilde{g})\Big) \quad \text{for every} \quad \tilde{g} \in \widetilde{G}(\gamma),$$

we conclude, based on Definition 3.23(a), that $\overline{\Theta}(\Delta(I))$ is a prolongation of $\widetilde{\Theta}(\Delta(I))$ to $(\overline{\mathbb{G}}(\Gamma(I)), \overline{i}(\Gamma^2(I)))$ and, therefore, now backed by Definition 3.23(b), that the $S$-space $\overline{\mathscr{G}}(I)$ is an extension of $\widetilde{\mathscr{G}}(I)$.

A particularly important aspect of $\overline{\mathscr{G}}(I)$, regarding our "axiomatic purposes" (described in 4.1 or, more detailed, at the Introduction of Chapter 3), consists of, if there exists an extension

$$\widehat{\mathscr{G}}(I) = \Bigg(\widehat{\mathbb{G}}\Big(\Gamma(I)\Big), \widehat{i}\Big(\Gamma^2(I)\Big), \widehat{\Theta}\Big(\Delta(I)\Big)\Bigg)$$

of $\widetilde{\mathscr{G}}(I)$ that is locally closed and coherent, then, one has, for any $\gamma, \gamma' \in \Gamma(I)$ such that $\gamma' \subseteq \gamma$, that:

**(a')** the function $A_\gamma : \overline{G}(\gamma) \longrightarrow \widehat{G}(\gamma)$ (defined in Proposition 4.12), according to Proposition 4.16, is an isomorphism from the group $\overline{G}(\gamma)$ onto the group $\widetilde{G}(\gamma)$;

**(b')** for every $\overline{\Phi} \in \overline{H}(\gamma)$, according to Proposition 4.22,

$$\overline{\Phi} = A_\gamma^{-1}\ \widehat{\Phi}\ A_\gamma;$$

**(c')** for $\overline{\Theta}_{(\gamma',\gamma)} \in \overline{\Theta}(\Delta(I))$, by Proposition 4.26,

$$\overline{\Theta}_{(\gamma',\gamma)} = A_{\gamma'}^{-1}\ \widehat{\Theta}_{(\gamma',\gamma)}\ A_\gamma.$$

Let us observe that, due to (b'), we can define the following function from the semigroup $\overline{H}(\gamma)$ into the semigroup $\widehat{H}(\gamma)$:

$$\begin{aligned} B_\gamma : \overline{H}(\gamma) &\longrightarrow \widehat{H}(\gamma) \\ \overline{\Phi} &\longmapsto B_\gamma(\overline{\Phi}) \coloneqq A_\gamma \overline{\Phi} A_\gamma^{-1} = \widehat{\Phi}. \end{aligned}$$

It immediately results, as a consequence of $A_\gamma : \overline{G}(\gamma) \longrightarrow \widehat{G}(\gamma)$ being an isomorphism, that $B_\gamma : \overline{H}(\gamma) \longrightarrow \widehat{H}(\gamma)$ is an isomorphism from the semigroup $\overline{H}(\gamma)$ onto the semigroup $\widehat{H}(\gamma)$.

In summary, $\overline{\mathscr{G}}(I)$ is an extension of $\widetilde{\mathscr{G}}(I)$ that, on the existence of $\widehat{\mathscr{G}}(I)$, for each $\gamma \in \Gamma(I)$, the functions $A_\gamma : \overline{G}(\gamma) \longrightarrow \widehat{G}(\gamma)$ and $B_\gamma : \overline{H}(\gamma) \longrightarrow \widehat{H}(\gamma)$ above referred are isomorphisms, between groups and semigroups, respectively, such that:

**(a)** taking (b') and the definition of $B_\gamma$ into account,

$$A_\gamma\Big(\overline{\Phi}(\overline{g})\Big) = B_\gamma(\overline{\Phi})\Big(A_\gamma(\overline{g})\Big)$$

for every $\overline{\Phi} \in \overline{H}(\gamma)$ and every $\overline{g} \in \overline{G}(\gamma)$;

**(b)** by (c'),

$$A_{\gamma'}\Big(\overline{\Theta}_{(\gamma',\gamma)}(\overline{g})\Big) = \widehat{\Theta}_{(\gamma',\gamma)}\Big(A_\gamma(\overline{g})\Big)$$

for every $(\gamma', \gamma) \in \Delta(I)$ and every $\overline{g} \in \overline{G}(\gamma)$.

Returning now to the definition of isomorphism for $S$-spaces, Definition 3.28, we see that the conditions (a) and (b) of the referred definition, translated to the case of the $S$-spaces $\overline{\mathscr{G}}(I)$ and $\widehat{\mathscr{G}}(I)$ are, exactly, the conditions (a) and (b) above which, as we saw, are satisfied by the $S$-spaces $\overline{\mathscr{G}}(I)$ and $\widehat{\mathscr{G}}(I)$. The condition (c) of the definition in question also is satisfied by $\overline{\mathscr{G}}(I)$ and $\widehat{\mathscr{G}}(I)$, that is,

$$B_{\gamma'}\Big(\overline{\imath}_{(\gamma',\gamma)}(\overline{\Phi})\Big) = \widehat{\imath}_{(\gamma',\gamma)}\Big(B_\gamma(\overline{\Phi})\Big)$$

for every $(\gamma', \gamma) \in \Gamma^2(I)$ and every $\overline{\Phi} \in \overline{H}(\gamma)$.

In fact, as we know, for every $\gamma \in \Gamma(I)$,

$$\overline{\gamma}(\widetilde{\Phi}) = \overline{\Phi} \quad \text{for every} \quad \widetilde{\Phi} \in \widetilde{H}(\gamma),$$

$$\widehat{\gamma}(\widetilde{\Phi}) = \widehat{\Phi} \quad \text{for every} \quad \widetilde{\Phi} \in \widetilde{H}(\gamma)$$

and, hence, since (by the definition of $B_\gamma$)

$$B_\gamma(\overline{\Phi}) = \widehat{\Phi},$$

we have

$$B_\gamma\Big(\overline{\gamma}(\widetilde{\Phi})\Big) = \widehat{\gamma}(\widetilde{\Phi}) \quad \text{for every} \quad \widetilde{\Phi} \in \widetilde{H}(\gamma),$$

that is,

$$B_\gamma \overline{\gamma} = \widehat{\gamma}$$

and, therefore,

$$B_\gamma = \widehat{\gamma}(\overline{\gamma})^{-1}.$$

Furthermore, from the definition of bonding extension, Definition 3.22, we obtain:

$$\widetilde{i}_{(\gamma',\gamma)} = (\overline{\gamma'})^{-1}\, \overline{i}_{(\gamma',\gamma)}\, \overline{\gamma}$$

and

$$\widehat{i}_{(\gamma',\gamma)} = \widehat{\gamma'}\, \widetilde{i}_{(\gamma',\gamma)}(\widehat{\gamma})^{-1}.$$

Now, with these data, we calculate:

$$\begin{aligned} B_{\gamma'}\Big(\overline{i}_{(\gamma',\gamma)}(\overline{\Phi})\Big) &= \widehat{\gamma'}\bigg((\overline{\gamma'})^{-1}\Big(\overline{i}_{(\gamma',\gamma)}(\overline{\Phi})\Big)\bigg) = \\ &= \widehat{\gamma'}\bigg((\overline{\gamma'})^{-1}\Big(\overline{i}_{(\gamma',\gamma)}\Big(\overline{\gamma}(\widetilde{\Phi})\Big)\Big)\bigg) = \\ &= \widehat{\gamma'}\Big(\widetilde{i}_{(\gamma',\gamma)}(\widetilde{\Phi})\Big) \end{aligned}$$

and

$$\begin{aligned} \widehat{i}_{(\gamma',\gamma)}\Big(B_\gamma(\overline{\Phi})\Big) &= \widehat{i}_{(\gamma',\gamma)}(\widehat{\Phi}) = \widehat{i}_{(\gamma',\gamma)}\Big(\widehat{\gamma}(\widetilde{\Phi})\Big) = \\ &= \widehat{\gamma'}\bigg(\widetilde{i}_{(\gamma',\gamma)}\bigg((\widehat{\gamma})^{-1}\Big(\widehat{\gamma}(\widetilde{\Phi})\Big)\bigg)\bigg) = \\ &= \widehat{\gamma'}\Big(\widetilde{i}_{(\gamma',\gamma)}(\widetilde{\Phi})\Big), \end{aligned}$$

from where we conclude that

$$B_{\gamma'}\Big(\overline{i}_{(\gamma',\gamma)}(\overline{\Phi})\Big) = \widehat{i}_{(\gamma',\gamma)}\Big(B_\gamma(\overline{\Phi})\Big),$$

that is, the condition (c) of the definition of isomorphism for $S$-spaces, as the two others, (a) and (b), also is satisfied regarding the $S$-spaces $\overline{\mathscr{G}}(I)$ and $\widehat{\mathscr{G}}(I)$. It is then established the following result:

- If there exists an extension $\widehat{\mathscr{G}}(I)$ of $\widetilde{\mathscr{G}}(I)$ that is locally closed and coherent, then, $\widehat{\mathscr{G}}(I)$ is isomorphic to $\overline{\mathscr{G}}(I)$.

The following question then arises: Is there any locally closed and coherent extension of $\widetilde{\mathscr{G}}(I)$?

Clearly, given the result above, such an extension of $\widetilde{\mathscr{G}}(I)$ exists if and only if $\overline{\mathscr{G}}(I)$ is a locally closed and coherent extension of $\widetilde{\mathscr{G}}(I)$.

The Propositions 4.30 and 4.31 below answer our question.

## 4.30 Proposition

*The S-space $\overline{\mathscr{G}}(I)$ defined in 4.29 is coherent.*

*Proof.* Let $\overline{\mathscr{G}}(\gamma)$ be a $S$-subspace of $\overline{\mathscr{G}}(I)$ and

$$\overline{g}\Big(\overline{\xi}(\gamma)\Big) = \Big\{\overline{g}_{\xi}\Big\}_{\xi\in\overline{\xi}(\gamma)}$$

be a coherent family in $\overline{\mathscr{G}}(\gamma)$, arbitrarily fixed. Hence, since $\overline{g}(\overline{\xi}(\gamma))$ is coherent in $\overline{\mathscr{G}}(\gamma)$, we have

$$\overline{\Theta}_{(\xi\cap\xi',\xi)}(\overline{g}_{\xi}) = \overline{\Theta}_{(\xi\cap\xi',\xi')}(\overline{g}_{\xi'}) \tag{4.30-1}$$

for every $\xi, \xi' \in \overline{\xi}(\gamma)$ such that $\xi \cap \xi' \neq \varnothing$.

Furthermore, since $\overline{g}_{\xi} \in \overline{G}(\xi)$ for each $\xi \in \overline{\xi}(\gamma)$, we also have

$$\overline{g}_{\xi} = \Big[\widetilde{g}\Big(\overline{\nu}(\xi)\Big)\Big]$$

where

$$\widetilde{g}\Big(\overline{\nu}(\xi)\Big) = \Big\{\widetilde{g}_{\nu_{\xi}}\Big\}_{\nu_{\xi}\in\overline{\nu}(\xi)}$$

is a coherent family in $\widetilde{\mathscr{G}}(\xi)$, that is, such that:

$$\widetilde{\Theta}_{(\nu_{\xi}\cap\nu'_{\xi},\nu_{\xi})}(\widetilde{g}_{\nu_{\xi}}) = \widetilde{\Theta}_{(\nu_{\xi}\cap\nu'_{\xi},\nu'_{\xi})}(\widetilde{g}_{\nu'_{\xi}}) \tag{4.30-2}$$

for every $\nu_{\xi}, \nu'_{\xi} \in \overline{\nu}(\xi)$ such that $\nu_{\xi} \cap \nu'_{\xi} \neq \varnothing$.

Now, with the definition of $\overline{\Theta}_{(\gamma'',\gamma')}$ in mind, we get:

$$\begin{aligned} \overline{\Theta}_{(\xi\cap\xi',\xi)}(\overline{g}_{\xi}) &= \overline{\Theta}_{(\xi\cap\xi',\xi)}\left(\Big[\Big\{\widetilde{g}_{\nu_{\xi}}\Big\}_{\nu_{\xi}\in\overline{\nu}(\xi)}\Big]\right) = \\ &= \Big[\Big\{\widetilde{\Theta}_{(\nu_{\xi}\cap\xi\cap\xi',\nu_{\xi})}(\widetilde{g}_{\nu_{\xi}})\Big\}_{\nu_{\xi}\cap\xi\cap\xi'\in\overline{\eta}(\xi\cap\xi')}\Big] = \\ &= \Big[\Big\{\widetilde{\Theta}_{(\nu_{\xi}\cap\xi',\nu_{\xi})}(\widetilde{g}_{\nu_{\xi}})\Big\}_{\nu_{\xi}\cap\xi'\in\overline{\eta}(\xi\cap\xi')}\Big] \end{aligned} \tag{4.30-3}$$

where $\overline{\eta}(\xi \cap \xi') \subseteq \Gamma(\xi \cap \xi')$ is the cover of $\xi \cap \xi'$ given by

$$\overline{\eta}(\xi \cap \xi') := \Big\{\eta := \nu_\xi \cap \xi \cap \xi' = \nu_\xi \cap \xi' \quad : \quad \nu_\xi \in \overline{\nu}(\xi) \quad \text{and} \quad \xi' \cap \nu_\xi \neq \varnothing\Big\}.$$

Analogously we obtain that

$$\begin{aligned}\overline{\Theta}_{(\xi\cap\xi',\xi')}(\overline{g}_{\xi'}) &= \overline{\Theta}_{(\xi\cap\xi',\xi')}\left(\left[\left\{\widetilde{g}_{\nu_{\xi'}}\right\}_{\nu_{\xi'}\in\overline{\nu}(\xi')}\right]\right) = \\ &= \left[\left\{\widetilde{\Theta}_{(\nu_{\xi'}\cap\xi\cap\xi',\nu_{\xi'})}(\widetilde{g}_{\nu_{\xi'}})\right\}_{\nu_{\xi'}\cap\xi\cap\xi'\in\overline{\alpha}(\xi\cap\xi')}\right] = \\ &= \left[\left\{\widetilde{\Theta}_{(\nu_{\xi'}\cap\xi,\nu_{\xi'})}(\widetilde{g}_{\nu_{\xi'}})\right\}_{\nu_{\xi'}\cap\xi\in\overline{\alpha}(\xi\cap\xi')}\right]\end{aligned} \tag{4.30-4}$$

where $\overline{\alpha}(\xi \cap \xi') \subseteq \Gamma(\xi \cap \xi')$ is the cover of $\xi \cap \xi'$ defined by

$$\overline{\alpha}(\xi \cap \xi') := \Big\{\alpha := \nu_{\xi'} \cap \xi \cap \xi' = \nu_{\xi'} \cap \xi \quad : \quad \nu_{\xi'} \in \overline{\nu}(\xi') \quad \text{and} \quad \xi \cap \nu_{\xi'} \neq \varnothing\Big\}.$$

From (4.30-1), (4.30-3) and (4.30-4) one concludes that

$$\left\{\widetilde{\Theta}_{(\nu_\xi\cap\xi',\nu_\xi)}(\widetilde{g}_{\nu_\xi})\right\}_{\nu_\xi\cap\xi'\in\overline{\eta}(\xi\cap\xi')} \approx \left\{\widetilde{\Theta}_{(\nu_{\xi'}\cap\xi,\nu_{\xi'})}(\widetilde{g}_{\nu_{\xi'}})\right\}_{\nu_{\xi'}\cap\xi\in\overline{\alpha}(\xi\cap\xi')}$$

and, therefore,

$$\widetilde{\Theta}_{(\nu_\xi\cap\xi'\cap\nu_{\xi'}\cap\xi,\nu_\xi\cap\xi')}\Big(\widetilde{\Theta}_{(\nu_\xi\cap\xi',\nu_\xi)}(\widetilde{g}_{\nu_\xi})\Big) = \widetilde{\Theta}_{(\nu_\xi\cap\xi'\cap\nu_{\xi'}\cap\xi,\nu_{\xi'}\cap\xi)}\Big(\widetilde{\Theta}_{(\nu_{\xi'}\cap\xi,\nu_{\xi'})}(\widetilde{g}_{\nu_{\xi'}})\Big),$$

that is,

$$\widetilde{\Theta}_{(\nu_\xi\cap\nu_{\xi'},\nu_\xi)}(\widetilde{g}_{\nu_\xi}) = \widetilde{\Theta}_{(\nu_\xi\cap\nu_{\xi'},\nu_{\xi'})}(\widetilde{g}_{\nu_{\xi'}}) \tag{4.30-5}$$

for every $\nu_\xi \in \overline{\nu}(\xi)$ and $\nu_{\xi'} \in \overline{\nu}(\xi')$ such that $\nu_\xi \cap \nu_{\xi'} \neq \varnothing$.

Let now $\overline{\chi}(\gamma) \subseteq \Gamma(\gamma)$ be the cover of $\gamma$ defined as follows: $\chi \in \overline{\chi}(\gamma)$ if and only if $\chi = \nu_\xi \in \overline{\nu}(\xi)$ for some $\xi \in \overline{\xi}(\gamma)$, that is,

$$\overline{\chi}(\gamma) := \bigcup_{\xi\in\overline{\xi}(\gamma)} \overline{\nu}(\xi).$$

We also define $\widetilde{g}(\overline{\chi}(\gamma))$ as the following family:

$$\widetilde{g}\Big(\overline{\chi}(\gamma)\Big) := \Big\{\widetilde{g}_\chi\Big\}_{\chi\in\overline{\chi}(\gamma)}.$$

It results from 4.30-5 that $\widetilde{g}(\overline{\chi}(\gamma))$ is a coherent family in $\widetilde{\mathscr{G}}(\gamma)$ and, hence,

$$\overline{g} := \Big[\widetilde{g}\Big(\overline{\chi}(\gamma)\Big)\Big]$$

is a well-defined element of the group $\overline{G}(\gamma)$. Below, we prove that this member $\overline{g} \in \overline{G}(\gamma)$ is such that, for each $\xi \in \overline{\xi}(\gamma)$,

$$\overline{\Theta}_{(\xi,\gamma)}(\overline{g}) = \overline{g}_\xi.$$

In order to do so, with $\xi \in \overline{\xi}(\gamma)$ arbitrarily fixed we calculate:

$$\overline{\Theta}_{(\xi,\gamma)}(\overline{g}) = \overline{\Theta}_{(\xi,\gamma)}\left(\left[\left\{\widetilde{g}_\chi\right\}_{\chi\in\overline{\chi}(\gamma)}\right]\right) = \left[\left\{\widetilde{\Theta}_{(\chi\cap\xi,\chi)}(\widetilde{g}_\chi)\right\}_{\chi\cap\xi\in\overline{\beta}(\xi)}\right]$$

where $\overline{\beta}(\xi) \subseteq \Gamma(\xi)$ is the cover of $\xi$ given by

$$\overline{\beta}(\xi) := \left\{\beta := \chi \cap \xi \quad : \quad \chi \in \overline{\chi}(\gamma) \quad \text{and} \quad \chi \cap \xi \neq \varnothing\right\}.$$

Now, as we know,

$$\overline{\Theta}_{(\xi,\gamma)}(\overline{g}) = \left[\left\{\widetilde{\Theta}_{(\chi\cap\xi,\chi)}(\widetilde{g}_\chi)\right\}_{\chi\cap\xi\in\overline{\beta}(\xi)}\right] = \overline{g}_\xi = \left[\left\{\widetilde{g}_{\nu_\xi}\right\}_{\nu_\xi\in\overline{\nu}(\xi)}\right]$$

if and only if

$$\left\{\widetilde{\Theta}_{(\chi\cap\xi,\chi)}(\widetilde{g}_\chi)\right\}_{\chi\cap\xi\in\overline{\beta}(\xi)} \approx \left\{\widetilde{g}_{\nu_\xi}\right\}_{\nu_\xi\in\overline{\nu}(\xi)}$$

that is, if and only if

$$\widetilde{\Theta}_{(\chi\cap\xi\cap\nu_\xi,\chi\cap\xi)}\left(\widetilde{\Theta}_{(\chi\cap\xi,\chi)}(\widetilde{g}_\chi)\right) = \widetilde{\Theta}_{(\chi\cap\xi\cap\nu_\xi,\nu_\xi)}(\widetilde{g}_{\nu_\xi})$$

for every $\chi \cap \xi \in \overline{\beta}(\xi)$ and $\nu_\xi \in \overline{\nu}(\xi)$ such that $\chi \cap \xi \cap \nu_\xi = \chi \cap \nu_\xi \neq \varnothing$, or yet,

$$\widetilde{\Theta}_{(\chi\cap\nu_\xi,\chi)}(\widetilde{g}_\chi) = \widetilde{\Theta}_{(\chi\cap\nu_\xi,\nu_\xi)}(\widetilde{g}_{\nu_\xi}) \tag{4.30-6}$$

for every $\chi \in \overline{\chi}(\gamma)$ and $\nu_\xi \in \overline{\nu}(\xi)$ such that $\chi \cap \nu_\xi \neq \varnothing$.

In fact, (4.30-6) is true since, if $\chi \in \overline{\chi}(\gamma) = \bigcup_{\lambda\in\overline{\xi}(\gamma)} \overline{\nu}(\lambda)$ is such that $\chi \in \overline{\nu}(\xi)$, i.e., $\chi = \nu'_\xi \in \overline{\nu}(\xi)$, (4.30-2) shows us that (4.30-6) is true. If, on the other hand, $\chi \notin \overline{\nu}(\xi)$, then, $\chi = \nu_{\xi'} \in \overline{\nu}(\xi')$ for some $\xi' \in \overline{\xi}(\gamma)$ such that $\xi' \neq \xi$ and, in this case, the truthness of (4.30-6) is ensured by (4.30-5). Therefore,

$$\overline{\Theta}_{(\xi,\gamma)}(\overline{g}) = \overline{g}_\xi \quad \text{for every} \quad \xi \in \overline{\xi}(\gamma).$$

Finally, we must prove that $\overline{g}$, as we defined it, is the only element of $\overline{G}(\gamma)$ that satisfies the condition above. Let us suppose, then, that

$$\overline{h} = \left[\left\{\widetilde{h}_\sigma\right\}_{\sigma\in\overline{\sigma}(\gamma)}\right] \in \overline{G}(\gamma)$$

is such that

$$\overline{\Theta}_{(\xi,\gamma)}(\overline{h}) = \overline{g}_\xi \quad \text{for every} \quad \xi \in \overline{\xi}(\gamma).$$

Hence, taking into account the definition of $\overline{\Theta}_{(\xi,\gamma)}$, we are supposing that, for $\xi \in \overline{\xi}(\gamma)$ arbitrarily fixed,

$$\overline{\Theta}_{(\xi,\gamma)}(\overline{h}) = \left[\left\{\widetilde{\Theta}_{(\sigma\cap\xi,\sigma)}(\widetilde{h}_\sigma)\right\}_{\sigma\cap\xi\in\overline{\Pi}(\xi)}\right] = \overline{g}_\xi = \left[\left\{\widetilde{g}_{\nu_\xi}\right\}_{\nu_\xi\in\overline{\nu}(\xi)}\right]$$

where

$$\overline{\Pi}(\xi) := \left\{\Pi := \sigma \cap \xi \quad : \quad \sigma \in \overline{\sigma}(\gamma) \quad \text{and} \quad \sigma \cap \xi \neq \varnothing\right\} \subseteq \Gamma(\xi)$$

is a cover of $\xi$.

From this hypothesis then follows that

$$\left\{\widetilde{\Theta}_{(\sigma\cap\xi,\sigma)}(\widetilde{h}_\sigma)\right\}_{\sigma\cap\xi\in\overline{\Pi}(\xi)} \approx \left\{\widetilde{g}_{\nu_\xi}\right\}_{\nu_\xi\in\overline{\nu}(\xi)}$$

that is,

$$\widetilde{\Theta}_{(\sigma\cap\xi\cap\nu_\xi,\sigma)}(\widetilde{h}_\sigma) = \widetilde{\Theta}_{(\sigma\cap\xi\cap\nu_\xi,\nu_\xi)}(\widetilde{g}_{\nu_\xi})$$

or yet, since $\sigma \cap \xi \cap \nu_\xi = \sigma \cap \nu_\xi$,

$$\widetilde{\Theta}_{(\sigma\cap\nu_\xi,\sigma)}(\widetilde{h}_\sigma) = \widetilde{\Theta}_{(\sigma\cap\nu_\xi,\nu_\xi)}(\widetilde{g}_{\nu_\xi}) \tag{4.30-7}$$

for every $\sigma \in \overline{\sigma}(\gamma)$ and $\nu_\xi \in \overline{\nu}(\xi)$ such that $\sigma \cap \nu_\xi \neq \varnothing$. Since $\xi \in \overline{\xi}(\gamma)$ is arbitrary, this last equation is valid for every $\xi \in \overline{\xi}(\gamma)$, that is, (4.30-7) is true for every $\sigma \in \overline{\sigma}(\gamma)$ and every $\nu_\xi \in \overline{\nu}(\xi)$ whatever $\xi \in \overline{\xi}(\gamma)$, such that $\sigma \cap \nu_\xi \neq \varnothing$ and, hence, since $\overline{\chi}(\gamma) = \bigcup_{\xi\in\overline{\xi}(\gamma)} \overline{\nu}(\xi)$, we have

$$\widetilde{\Theta}_{(\sigma\cap\chi,\sigma)}(\widetilde{h}_\sigma) = \widetilde{\Theta}_{(\sigma\cap\chi,\chi)}(\widetilde{g}_\chi)$$

for every $\sigma \in \overline{\sigma}(\gamma)$ and $\chi \in \overline{\chi}(\gamma)$ such that $\sigma \cap \chi \neq \varnothing$. But this means that

$$\left\{\widetilde{h}_\sigma\right\}_{\sigma\in\overline{\sigma}(\gamma)} \approx \left\{\widetilde{g}_\chi\right\}_{\chi\in\overline{\chi}(\gamma)},$$

that is,

$$\overline{h} = \overline{g}. \qquad \blacksquare$$

## 4.31 Proposition

*The S-space $\overline{\mathscr{G}}(I)$ defined in 4.29 is a locally closed extension of the S-space $\widetilde{\mathscr{G}}(I)$ described in 4.6.*

*Proof.* Let $\gamma \in \Gamma(I)$, $x \in \gamma$ and

$$\overline{g} = \left[\left\{\widetilde{g}_\xi\right\}_{\xi\in\overline{\xi}(\gamma)}\right] \in \overline{G}(\gamma)$$

be arbitrarily fixed. We must prove (according to Definition 3.23(c-2)) that there exists $\gamma' \in \Gamma(\gamma)$ such that $x \in \gamma'$ and

$$\overline{\Theta}_{(\gamma',\gamma)}(\overline{g}) \in b_{\gamma'}\left(\widetilde{G}(\gamma')\right).$$

Since $\overline{\xi}(\gamma) \subseteq \Gamma(\gamma)$ is a cover of $\gamma$ and $x \in \gamma$, then, there exists $\xi_x \in \overline{\xi}(\gamma)$ such that $x \in \xi_x$; let us take

$$\gamma' := \xi_x.$$

Hence, $\gamma' \in \Gamma(\gamma)$ and $x \in \gamma'$. We select now among the members of the family $\{\widetilde{g}_\xi\}_{\xi\in\overline{\xi}(\gamma)}$, the one that corresponds to $\gamma' \in \overline{\xi}(\gamma)$, that is,

$$\widetilde{g}_{\gamma'} \in \widetilde{G}(\gamma').$$

Taking now into account the definition of the function $b_{\gamma'} : \widetilde{G}(\gamma') \longrightarrow \overline{G}(\gamma')$ (see Proposition 4.18), we get

$$b_{\gamma'}(\widetilde{g}_{\gamma'}) = \left[\left\{\widetilde{\Theta}_{(\chi,\gamma')}(\widetilde{g}_{\gamma'})\right\}_{\chi\in\overline{\chi}(\gamma')}\right] \tag{4.31-1}$$

where $\overline{\chi}(\gamma') \subseteq \Gamma(\gamma')$ is a cover of $\gamma'$.

On the other hand, by the definition of $\overline{\Theta}_{(\gamma',\gamma)}$ (see Proposition 4.26),

$$\overline{\Theta}_{(\gamma',\gamma)}(\overline{g}) = \overline{\Theta}_{(\gamma',\gamma)}\left(\left[\left\{\widetilde{g}_\xi\right\}_{\xi\in\overline{\xi}(\gamma)}\right]\right) = \left[\left\{\widetilde{\Theta}_{(\xi\cap\gamma',\xi)}(\widetilde{g}_\xi)\right\}_{\xi\cap\gamma'\in\overline{\eta}(\gamma')}\right]$$

where

$$\overline{\eta}(\gamma') := \left\{\xi \cap \gamma' \quad : \quad \xi \in \overline{\xi}(\gamma) \quad \text{and} \quad \xi \cap \gamma' \neq \varnothing\right\}.$$

Since the family $\{\widetilde{g}_\xi\}_{\xi\in\overline{\xi}(\gamma)}$ is coherent in $\widetilde{\mathscr{G}}(\gamma)$, we have

$$\widetilde{\Theta}_{(\xi\cap\gamma',\xi)}(\widetilde{g}_\xi) = \widetilde{\Theta}_{(\xi\cap\gamma',\gamma')}(\widetilde{g}_{\gamma'})$$

for every $\xi \in \overline{\xi}(\gamma)$ such that $\xi \cap \gamma' \neq \varnothing$, allowing us to write $\overline{\Theta}_{(\gamma',\gamma)}(\overline{g})$ as follows:

$$\overline{\Theta}_{(\gamma',\gamma)}(\overline{g}) = \left[\left\{\widetilde{\Theta}_{(\xi\cap\gamma',\gamma')}(\widetilde{g}_{\gamma'})\right\}_{\xi\cap\gamma'\in\overline{\eta}(\gamma')}\right]. \tag{4.31-2}$$

Let us note now that

$$\left\{\widetilde{\Theta}_{(\xi\cap\gamma',\gamma')}(\widetilde{g}_{\gamma'})\right\}_{\xi\cap\gamma'\in\overline{\eta}(\gamma')} \approx \left\{\widetilde{\Theta}_{(\chi,\gamma')}(\widetilde{g}_{\gamma'})\right\}_{\chi\in\overline{\chi}(\gamma')}.$$

In fact, since

$$\widetilde{\Theta}_{(\xi\cap\gamma'\cap\chi,\xi\cap\gamma')}\left(\widetilde{\Theta}_{(\xi\cap\gamma',\gamma')}(\widetilde{g}_{\gamma'})\right) = \widetilde{\Theta}_{(\xi\cap\gamma'\cap\chi,\chi)}\left(\widetilde{\Theta}_{(\chi,\gamma')}(\widetilde{g}_{\gamma'})\right)$$

for every $\xi \cap \gamma' \in \overline{\eta}(\gamma')$ and $\chi \in \overline{\chi}(\gamma')$ such that $\xi \cap \gamma' \cap \chi \neq \varnothing$.

Therefore,

$$\left[\left\{\widetilde{\Theta}_{(\xi\cap\gamma',\gamma')}(\widetilde{g}_{\gamma'})\right\}_{\xi\cap\gamma'\in\overline{\eta}(\gamma')}\right] = \left[\left\{\widetilde{\Theta}_{(\chi,\gamma')}(\widetilde{g}_{\gamma'})\right\}_{\chi\in\overline{\chi}(\gamma')}\right]$$

and, hence, taking (4.31-1) and (4.31-2) into account, we conclude that

$$\overline{\Theta}_{(\gamma',\gamma)}(\overline{g}) = b_{\gamma'}(\widetilde{g}_{\gamma'}),$$

which shows us that

$$\overline{\Theta}_{(\gamma',\gamma)}(\overline{g}) \in b_{\gamma'}\left(\widetilde{G}(\gamma')\right).$$

■

# Second Theorem of Extension of $S$-Spaces

## 4.32 Remark

Throughout this chapter we obtained two extensions, $\widetilde{\mathscr{G}}(I)$ and $\overline{\mathscr{G}}(I)$, for a given abelian, surjective, with identity and corehent $S$-space, $\mathscr{G}(I)$:

**(a)** $\widetilde{\mathscr{G}}(I)$ — a strict and closed extension of $\mathscr{G}(I)$ (in this case, $\mathscr{G}(I)$ does not need to be coherent);

**(b)** $\overline{\mathscr{G}}(I)$ — a coherent and locally closed extension of $\widetilde{\mathscr{G}}(I)$ (and, therefore, an extension of $\mathscr{G}(I)$, according to the conclusion obtained at the end of item 3.26 (page 120)).

For $\widetilde{\mathscr{G}}(I)$ we proved the 1st TESS (Theorem 4.7), which establishes its uniqueness unless isomorphism. Regarding $\overline{\mathscr{G}}(I)$, our considerations in item 4.29 together with Propositions 4.30 and 4.31 prove that $\overline{\mathscr{G}}(I)$ also is unique unless isomorphism. This result, as well as the 1st TESS, as we will see in the next chapter, is fundamental for our goals (described in itens 3.3 and 3.4) and, therefore, it is explicitly stated below; we will refer to it as the **Second Theorem of Extension of $S$-Spaces** or, shortly, **2nd TESS**.

## 4.33 Theorem (2nd TESS)

*Let*

$$\mathscr{G}(I) = \Bigg( \mathbb{G}\Big(\Gamma(I)\Big) = \Big\{ \mathbb{G}(\gamma) = \Big(G(\gamma), H(\gamma)\Big) \Big\}_{\gamma\in\Gamma(I)}, \quad i\Big(\Gamma^2(I)\Big) = \Big\{ i_{(\gamma',\gamma)} \Big\}_{(\gamma',\gamma)\in\Gamma^2(I)}, \quad \Theta\Big(\Delta(I)\Big) = \Big\{ \Theta_{(\gamma',\gamma)} \Big\}_{(\gamma',\gamma)\in\Delta(I)} \Bigg)$$

*be an abelian, surjective, with identity and coherent $S$-space, and*

$$\widetilde{\mathscr{G}}(I) = \Bigg( \widetilde{\mathbb{G}}\Big(\Gamma(I)\Big) = \Big\{ \widetilde{\mathbb{G}}(\gamma) = \Big(\widetilde{G}(\gamma), \widetilde{H}(\gamma)\Big) \Big\}_{\gamma\in\Gamma(I)}, \quad \tilde{i}\Big(\Gamma^2(I)\Big) = \Big\{ \tilde{i}_{(\gamma',\gamma)} \Big\}_{(\gamma',\gamma)\in\Gamma^2(I)}, \quad \widetilde{\Theta}\Big(\Delta(I)\Big) = \Big\{ \widetilde{\Theta}_{(\gamma',\gamma)} \Big\}_{(\gamma',\gamma)\in\Delta(I)} \Bigg)$$

*its strict and closed extension described in 4.6, unique unless isomorphism. Let also*

$$\overline{\mathscr{G}}(I) = \left(\overline{\mathbb{G}}\Big(\Gamma(I)\Big) = \left\{\overline{\mathbb{G}}(\gamma) = \Big(\overline{G}(\gamma), \overline{H}(\gamma)\Big)\right\}_{\gamma\in\Gamma(I)}, \quad \overline{i}\Big(\Gamma^2(I)\Big) = \left\{\overline{i}_{(\gamma',\gamma)}\right\}_{(\gamma',\gamma)\in\Gamma^2(I)}, \quad \overline{\Theta}\Big(\Delta(I)\Big) = \left\{\overline{\Theta}_{(\gamma',\gamma)}\right\}_{(\gamma',\gamma)\in\Delta(I)}\right)$$

*be the extension of* $\widetilde{\mathscr{G}}(I)$ *described in 4.29.*

*Then,* $\overline{\mathscr{G}}(I)$ *is a coherent and locally closed extension of* $\widetilde{\mathscr{G}}(I)$*. Furthermore,* $\overline{\mathscr{G}}(I)$ *is unique unless isomorphism, that is, if*

$$\widehat{\mathscr{G}}(I) = \left(\widehat{\mathbb{G}}\Big(\Gamma(I)\Big) = \left\{\widehat{\mathbb{G}}(\gamma) = \Big(\widehat{G}(\gamma), \widehat{H}(\gamma)\Big)\right\}_{\gamma\in\Gamma(I)}, \quad \widehat{i}\Big(\Gamma^2(I)\Big) = \left\{\widehat{i}_{(\gamma',\gamma)}\right\}_{(\gamma',\gamma)\in\Gamma^2(I)}, \quad \widehat{\Theta}\Big(\Delta(I)\Big) = \left\{\widehat{\Theta}_{(\gamma',\gamma)}\right\}_{(\gamma',\gamma)\in\Delta(I)}\right)$$

*is another coherent and locally closed extension of* $\widetilde{\mathscr{G}}(I)$*, then, for each* $\gamma \in \Gamma(I)$*, the functions* $A_\gamma : \overline{G}(\gamma) \longrightarrow \widehat{G}(\gamma)$ *and* $B_\gamma : \overline{H}(\gamma) \longrightarrow \widehat{H}(\gamma)$ *defined, respectively, in Proposition 4.12 (or, without abuse of language, in (a) on page 191) and item 4.29 (page 205), are isomorphisms such that:*

**(a)** $A_\gamma\Big(\overline{\Phi}(\overline{g})\Big) = B_\gamma(\overline{\Phi})\Big(A_\gamma(\overline{g})\Big)$ *for every* $\overline{\Phi} \in \overline{H}(\gamma)$ *and every* $\overline{g} \in \overline{G}(\gamma)$*;*

**(b)** $A_{\gamma'}\Big(\overline{\Theta}_{(\gamma',\gamma)}(\overline{g})\Big) = \widehat{\Theta}_{(\gamma',\gamma)}\Big(A_\gamma(\overline{g})\Big)$ *for every* $(\gamma',\gamma) \in \Delta(I)$ *and every* $\overline{g} \in \overline{G}(\gamma)$*;*

**(c)** $B_{\gamma'}\Big(\overline{i}_{(\gamma',\gamma)}(\overline{\Phi})\Big) = \widehat{i}_{(\gamma',\gamma)}\Big(B_\gamma(\overline{\Phi})\Big)$ *for every* $(\gamma',\gamma) \in \Gamma^2(I)$ *and every* $\overline{\Phi} \in \overline{H}(\gamma)$*.*

*Evenmore,* $A_\gamma : \overline{G}(\gamma) \longrightarrow \widehat{G}(\gamma)$ *is the only isomorphism that keeps fixed the elements of* $\widetilde{G}(\gamma)$*, that is,*

$$A_\gamma(\widetilde{g}) = \widetilde{g} \quad \textit{for every} \quad \widetilde{g} \in \widetilde{G}(\gamma)$$

*or, more precisely (that is, without abuse of language),*

$$A_\gamma\Big(b_\gamma(\widetilde{g})\Big) = \delta_\gamma(\widetilde{g}) \quad \textit{for every} \quad \widetilde{g} \in \widetilde{G}(\gamma),$$

*being* $b_\gamma : \widetilde{G}(\gamma) \longrightarrow \overline{G}(\gamma)$ *the injective homomorphism defined in Proposition 4.18 and* $\delta_\gamma : \widetilde{G}(\gamma) \longrightarrow \widehat{G}(\gamma)$ *the injective homomorphism that backs the fact of* $\widehat{G}(\gamma)$ *admiting the group* $\widetilde{G}(\gamma)$ *as a subgroup (see Definition 1.6(a) and Remark 4.25).* ■

We finish this chapter with the proof of a proposition that will show itself useful at the final part of Chapter 5, more precisely at the considerations of items 5.27 and 5.38. To state it, we will use the notation below, where $\overset{*}{\mathscr{G}}(I)$ and $\overset{\circ}{\mathscr{G}}(I)$ represent arbitrary $S$-spaces.

**Notation.**

$\overset{\circ}{\mathscr{G}}(I) \leqslant \overset{*}{\mathscr{G}}(I)$ if and only if $\overset{*}{\mathscr{G}}(I)$ is an extension of $\overset{\circ}{\mathscr{G}}(I)$;

$\overset{\circ}{\mathscr{G}}(I) \simeq \overset{*}{\mathscr{G}}(I)$ if and only if $\overset{\circ}{\mathscr{G}}(I)$ is isomorphic to $\overset{*}{\mathscr{G}}(I)$.

## 4.34 Proposition

*Let the $S$-spaces $\mathscr{G}(I)$, $\widetilde{\mathscr{G}}(I)$ and $\overline{\mathscr{G}}(I)$ be as described in Theorem 4.33 (2nd TESS). Let also*

$$\overset{\circ}{\mathscr{G}}(I) = \left( \begin{aligned} \overset{\circ}{\mathbb{G}}\Big(\Gamma(I)\Big) &= \Big\{\overset{\circ}{\mathbb{G}}(\gamma) = \Big(\overset{\circ}{G}(\gamma), \overset{\circ}{H}(\gamma)\Big)\Big\}_{\gamma\in\Gamma(I)}, \\ \overset{\circ}{i}\Big(\Gamma^2(I)\Big) &= \Big\{\overset{\circ}{i}_{(\gamma',\gamma)}\Big\}_{(\gamma',\gamma)\in\Gamma^2(I)}, \\ \overset{\circ}{\Theta}\Big(\Delta(I)\Big) &= \Big\{\overset{\circ}{\Theta}_{(\gamma',\gamma)}\Big\}_{(\gamma',\gamma)\in\Delta(I)} \end{aligned} \right)$$

*and*

$$\overset{*}{\mathscr{G}}(I) = \left( \begin{aligned} \overset{*}{\mathbb{G}}\Big(\Gamma(I)\Big) &= \Big\{\overset{*}{\mathbb{G}}(\gamma) = \Big(\overset{*}{G}(\gamma), \overset{*}{H}(\gamma)\Big)\Big\}_{\gamma\in\Gamma(I)}, \\ \overset{*}{i}\Big(\Gamma^2(I)\Big) &= \Big\{\overset{*}{i}_{(\gamma',\gamma)}\Big\}_{(\gamma',\gamma)\in\Gamma^2(I)}, \\ \overset{*}{\Theta}\Big(\Delta(I)\Big) &= \Big\{\overset{*}{\Theta}_{(\gamma',\gamma)}\Big\}_{(\gamma',\gamma)\in\Delta(I)} \end{aligned} \right)$$

*be $S$-spaces such that:*

$$\mathscr{G}(I) \leqslant \overset{\circ}{\mathscr{G}}(I), \quad \overset{\circ}{\mathscr{G}}(I) \simeq \widetilde{\mathscr{G}}(I)$$

*and $\overset{*}{\mathscr{G}}(I)$ is a coherent and locally closed extension of $\overset{\circ}{\mathscr{G}}(I)$. Under these conditions we have*

$$\overset{*}{\mathscr{G}}(I) \simeq \overline{\mathscr{G}}(I).$$

*Proof.* We have

$$\mathscr{G}(I) \leqslant \overset{\circ}{\mathscr{G}}(I) \quad \text{(by hypothesis)},$$

$$\mathscr{G}(I) \leqslant \widetilde{\mathscr{G}}(I) \quad \text{(by 2nd TESS)} \quad \text{and}$$

$$\overset{\circ}{\mathscr{G}}(I) \simeq \widetilde{\mathscr{G}}(I) \quad \text{(by hypothesis)}.$$

Now, it results from the definitions of $S$-space extension (Definition 3.23) and $S$-space isomorphism (Definition 3.28) that: if each one among two isomorphic $S$-spaces ($\mathring{\mathscr{G}}(I)$ and $\widetilde{\mathscr{G}}(I)$ in the present case) is an extension of a same $S$-space ($\mathscr{G}(I)$, in our case, since $\mathscr{G}(I) \leqslant \mathring{\mathscr{G}}(I)$ and $\mathscr{G}(I) \leqslant \widetilde{\mathscr{G}}(I)$), then, any one of these two ($\mathring{\mathscr{G}}(I)$ or $\widetilde{\mathscr{G}}(I)$) is an extension of the other. Hence, we have that:

$$\widetilde{\mathscr{G}}(I) \leqslant \mathring{\mathscr{G}}(I) \quad (\text{as well as } \mathring{\mathscr{G}}(I) \leqslant \widetilde{\mathscr{G}}(I)).$$

Since (by hypothesis)

$$\mathring{\mathscr{G}}(I) \leqslant \overset{*}{\mathscr{G}}(I),$$

and taking into account the last result of item 3.26 (presented at the foot of page 120), we conclude that

$$\widetilde{\mathscr{G}}(I) \leqslant \overset{*}{\mathscr{G}}(I).$$

We prove below that $\overset{*}{\mathscr{G}}(I)$ is a coherent and locally closed extension of $\widetilde{\mathscr{G}}(I)$. In fact, by hypothesis, $\overset{*}{\mathscr{G}}(I)$ is a coherent and locally closed extension of $\mathring{\mathscr{G}}(I)$ and, hence, the $S$-space $\overset{*}{\mathscr{G}}(I)$ is coherent and yet, for $\gamma \in \Gamma(I)$, $\overset{*}{g} \in \overset{*}{G}(\gamma)$ and $x \in \gamma$ arbitrarily fixed, there exists $\gamma' \in \Gamma(\gamma)$ such that $x \in \gamma'$ and

$$\overset{*}{\Theta}_{(\gamma',\gamma)}(\overset{*}{g}) \in \mathring{G}(\gamma'),$$

that is,

$$\overset{*}{\Theta}_{(\gamma',\gamma)}(\overset{*}{g}) = \mathring{g}$$

for some $\mathring{g} \in \mathring{G}(\gamma')$. But $\mathring{\mathscr{G}}(I) \simeq \widetilde{\mathscr{G}}(I)$ (by hypothesis) and, hence, with the definition of isomorphism for $S$-spaces in mind, there exists an isomorphism from the group $\widetilde{G}(\gamma')$ onto the group $\mathring{G}(\gamma')$, $\alpha_{\gamma'} : \widetilde{G}(\gamma') \longrightarrow \mathring{G}(\gamma')$ and, therefore, a (unique) $\widetilde{g} \in \widetilde{G}(\gamma')$ such that

$$\overset{*}{\Theta}_{(\gamma',\gamma)}(\overset{*}{g}) = \alpha_{\gamma'}(\widetilde{g}) = \mathring{g}.$$

Thus,

$$\overset{*}{\Theta}_{(\gamma',\gamma)}(\overset{*}{g}) \in \alpha_{\gamma'}\Big(\widetilde{G}(\gamma')\Big),$$

which show us that $\overset{*}{\mathscr{G}}(I)$ is a locally closed extension of $\widetilde{\mathscr{G}}(I)$ (consider Warning 3.24). Since $\overset{*}{\mathscr{G}}(I)$ is a coherent $S$-space, we then have that $\overset{*}{\mathscr{G}}(I)$ is a coherent and locally closed extension of $\widetilde{\mathscr{G}}(I)$ and, hence, by the 2nd TESS, we conclude that

$$\overset{*}{\mathscr{G}}(I) \simeq \overline{\mathscr{G}}(I). \qquad \blacksquare$$

# 5

# DISTRIBUTIONS AND ITS AXIOMATICS

## Introduction

### 5.1 Retrospective

**(a)** The starting point.

The concepts of distribution and derivative of a distribution were informal and intuitively introduced in the first items of this monograph (itens 1.1 and 1.3), as generalizations (extensions) of the notions of continuous function and derivative of a function, respectively, to circumvent the difficulties and limitations imposed by the existence of continuous but not differentiable functions; these difficulties led us to conceive a class of objects, the distributions, which properly includes the collection of continuous functions, and, for the members of this new class, a notion of differentiability which extends the corresponding classical notion in such a way that every distribution, and therefore every continuous function, is infinitely differentiable.

Guided by a set of conditions which naturally impose itself as a requirement to such a generalization of the notions of continuous function and classical differentiability (see 1.3(a)-(d)), we identified a basic domain in relation to which the questions relative to our goal could be formulated, namely, the ordered pair

$$\Big(C(\Omega), \partial(\Omega)\Big),$$

where $C(\Omega)$ is the abelian group of the functions defined and continuous on the open set $\Omega \subseteq \mathbb{R}^n$, and $\partial(\Omega)$ is the semigroup of the partial derivatives (see 1.2 and 1.5). Hence, we were led to a first version, still informal, of the problem posed by the intended generalization:

- To obtain an extension $(\widetilde{C}(\Omega), \widetilde{\partial}(\Omega))$ of the domain $(C(\Omega), \partial(\Omega))$, that satisfies the "naturally imposed conditions" (above referred to).

$(\widetilde{C}(\Omega), \widetilde{\partial}(\Omega))$ then starts to represent the generalization we want to obtain, the "habitat of the distributions", being the distributions the elements of $\widetilde{C}(\Omega)$ and its derivatives the members of $\widetilde{\partial}(\Omega)$.

**(b)** The precise characterization of the goal.

After a formalization, through precise definitions (Definitions 1.4, 1.6, 1.10 and 1.14), of several concepts suggested by an informal and intuitive approach that, partially, we remembered above (in (a)), it was possible to characterize our goal, rigorously, in the form of the following problem (see item 1.15):

- Strict and Closed Extension Problem (or, in short, SCEP):
  Given the abelian, surjective, and with identity $S$-group

$$\mathbb{C}(\Omega) = \Big(C(\Omega), \partial(\Omega)\Big),$$

  to construct an extension,

$$\widetilde{\mathbb{C}}(\Omega) = \Big(\widetilde{C}(\Omega), \widetilde{\partial}(\Omega)\Big),$$

  of $\mathbb{C}(\Omega)$ that is strict and closed.

Furthermore, since the definitions of the concepts involved in the SCEP could be given in an abstract (axiomatic) way, it was possible to state it in the following abstract version:

- SCEP (abstract version):
  Given the abelian, surjective, and with identity $S$-group

$$\mathbb{G} = \Big(G, H\Big),$$

  to construct an extension,

$$\widetilde{\mathbb{G}} = \Big(\widetilde{G}, \widetilde{H}\Big),$$

  of $\mathbb{G}$ that is strict and closed.

From this, the search for our generalization of the continuous functions and its derivatives, the "habitat of the distributions", became to be the search for a solution $(\widetilde{C}(\Omega), \widetilde{\partial}(\Omega))$ to the SCEP. It was clear that one such "universe" would also be obtained by specializing (particularizing) a solution $(\widetilde{G}, \widetilde{H})$ of the abstract version of the SCEP, if we could find any, to the case where $G = C(\Omega)$ and $H = \partial(\Omega)$.

**(c)** The solution to the SCEP.

With the proof of the Theorem of Extension of $S$-Groups (Theorem 2.16), we not only obtained a specific and well-defined solution $(\widetilde{G}, \widetilde{H})$ to the SCEP in its abstract version, but also proved that the obtained solution is, essentially, the only solution to the problem; translating this solution to the case where $G = C(\Omega)$ and $H = \partial(\Omega)$, we determined the (essentially unique) "habitat of the distributions", namely, the solution $(\widetilde{C}(\Omega), \widetilde{\partial}(\Omega))$ of the SCEP in its original, non-abstract, version.

**(d)** The $\widetilde{\mathbb{G}}$-distributions and its derivatives.

Once obtained, through the proof of the Theorem of Extension of $S$-Groups, the essentially unique solution to the SCEP in its abstract version, and, consequently, the one for the referred problem in its formulation associated to the $S$-group $(C(\Omega), \partial(\Omega))$, we were left with no doubt about how we should precisely define the notions of distribution and its derivatives that, until then, only had an informal and intuitive conceptualization; formal and precisely, a distribution is any element of the class $\widetilde{C}(\Omega)$, and derivative of distribution is any of the members of the collection $\widetilde{\partial}(\Omega)$, being $(\widetilde{C}(\Omega), \widetilde{\partial}(\Omega))$ the solution to the SCEP in its original version.

Since the solution $(\widetilde{C}(\Omega), \widetilde{\partial}(\Omega))$ is obtained from the solution $(\widetilde{G}, \widetilde{H})$ to the SCEP in its abstract version, specializing the latter to the case where $G = C(\Omega)$ and $H = \partial(\Omega)$, we went naturally induced to introduce abstract (axiomatical) notions of distribution and derivative of distribution; hence we defined $\widetilde{\mathbb{G}}$-distribution as any element of the set $\widetilde{G}$ and $\widetilde{\mathbb{G}}$-distribution derivative as any member of the class $\widetilde{H}$ (Definition 2.18) and, from this, the distributions and its derivatives as the members of $\widetilde{C}(\Omega)$ and $\widetilde{\partial}(\Omega)$, respectively, turned out to be exactly the $\widetilde{\mathbb{C}}(\Omega)$-distributions and the $\widetilde{\mathbb{C}}(\Omega)$-distributions derivatives, respectively.

**(e)** The $\widetilde{\mathbb{G}}$-distributions axiomatics.

Inspired and motivated by the Theorem of Extension of $S$-Groups (see Remark 2.19), we formulated two equivalent axiomatics (in 2.20 and 2.25) that categorically define the $\widetilde{\mathbb{G}}$-distributions and its derivatives. From this axiomatics, specializing them to the model associated to the $S$-group of continuous functions $\mathbb{C}(\Omega) = (C(\Omega), \partial(\Omega))$, it then resulted two axiomatics (categoric and equivalents) to our, now, $\widetilde{\mathbb{C}}(\Omega)$-distributions and its derivatives.

The comparison of these results regarding $\widetilde{\mathbb{C}}(\Omega)$-distributions and its derivatives, with the corresponding concepts defined in the Distributions Theory elaborated by the French mathematician L. Schwartz, revealed that the generalization (extension) of the notions of continuous function and differentiability promoted by Schwartz' theory, is broader than the one suggested by the Theorem of Extension of $S$-Groups and crystallized in the definition of $\widetilde{\mathbb{C}}(\Omega)$-distributions and its derivatives; as we proved in our schematical presentation of the Schwartz' theory (items 2.27 to 2.40), our axiomatic definition of the $\widetilde{\mathbb{C}}(\Omega)$-distributions and its derivatives (in any of its versions) does **not** admit as a model the class $\mathscr{D}'(\Omega)$ of all the Schwartz's distributions in $\Omega$, but only the proper subclass $\mathscr{D}'_f(\Omega) \subseteq \mathscr{D}'(\Omega)$ of the finite order distributions in $\Omega$.

**(f)** The $S$-spaces.

What structure should replace the one of $S$-group, and what type of extension we should require of this new structure in the formulation of an extension problem that, like the SCEP for $S$-groups, would lead us to a theorem that, such as the Theorem of Extension of $S$-Groups, would ensure the existence of one, essentially unique, solution that could capture not only the finite order distributions but the whole class of the Schwartz' distributions?

Based on the concept of local equality of distributions and one result, both belonging to Schwartz' theory, establishing that every distribution is, locally, the derivative of a continuous function, the considerations made at the Introduction of Chapter 3 (items 3.1 and 3.4) led us, naturally, to an informal notion of a "space", $\mathscr{C}(\mathbb{R}^n)$, constituted of a family of $S$-groups (one $S$-group $\mathbb{C}(\Omega) = (C(\Omega), \partial(\Omega))$ for each open set $\Omega \subseteq \mathbb{R}^n$) and a family of functions bonding the groups of these $S$-groups (the restriction functions $\rho_{(\Omega',\Omega)} : C(\Omega) \longrightarrow C(\Omega')$); even more, these same considerations suggested, also naturally, two types of extension that we could require of the space $\mathscr{C}(\mathbb{R}^n)$ that could promote a generalization, to the notions of continuous function and differentiability, broader than that provided by the strict and closed extension, $\widetilde{\mathbb{C}}(\Omega)$, of the $S$-group $\mathbb{C}(\Omega)$, possibly as broad as the corresponding generalization promoted by Schwartz' theory. Hence, equipped with these informally obtained notions, the one of "space" $\mathscr{C}(\mathbb{R}^n)$ and its "promising" extensions, we formulated, informally, two extension problems associated to the space $\mathscr{C}(\mathbb{R}^n)$, and dedicated all Chapter 3 to, through a network of rigorously defined concepts, characterize in a precise way the space $\mathscr{C}(\mathbb{R}^n)$ and the two extension problems associated to it; in reality, the referred Chapter 3 went well beyond that, insofar as it provided the definition of an abstract (axiomatic) structure, the one of $S$-space, and the precise formulation of two extension problems associated to it, which admit as models, respectively, the space $\mathscr{C}(\mathbb{R}^n)$ and its corresponding extension problems.

In summary, Chapter 3 was intended to motivate and establish in a precise and abstract (axiomatic) way the concepts that allowed us to rigorously formulate the following extension problems:

- Given an abelian, surjective, with identity, and coherent $S$-space $\mathscr{G}(I)$, one asks:
  - **(i)** To obtain a strict and closed extension, $\widetilde{\mathscr{G}}(I)$, of $\mathscr{G}(I)$;
  - **(ii)** Being $\widetilde{\mathscr{G}}(I)$ as in (i), to obtain an extension $\overline{\mathscr{G}}(I)$ of $\widetilde{\mathscr{G}}(I)$, that is locally closed and coherent.

**(g)** The theorems of extension of $S$-spaces.

Finally, in Chapter 4, we were able to prove two theorems, the 1st and 2nd TESS (Theorems 4.7 and 4.33), that solved, respectively, the problems (i) and (ii) above formulated, establishing not only the existence of one, essentially unique, solution to each one of them, but also exhibiting them.

## 5.2 The Current Stage

We find ourselves now, regarding an abelian, surjective, with identity, and coherent $S$-space $\mathscr{G}(I)$, at an analogous point where we found ourselves regarding an abelian, surjective, and with identity $S$-group $\mathbb{G} = (G, H)$ when we proved the Theorem of Extension of $S$-Groups. At that stage, backed by the referred theorem, we could formulate precise definitions for the notions, until then informal, of distributions and its derivatives, that came to formally characterize these concepts not only for the situation associated with the continuous functions and the classic notion of differentiability (the case corresponding to the $S$-group $\mathbb{C}(\Omega) = (C(\Omega), \partial(\Omega))$), but also, and mainly, defining them in a context associated to any $S$-group $\mathbb{G} = (G, H)$ that, as $\mathbb{C}(\Omega)$, is abelian, surjective, and with identity; in this way originated, at that stage, the $\widetilde{\mathbb{G}}$-distributions an its derivatives and, as a particular case (a model), the $\widetilde{\mathbb{C}}(\Omega)$-distributions and its respective derivatives — the latter ones attending our basic goal of freeing differential calculus of the limitations imposed by the classic notion of differentiability, whereas the first, the $\widetilde{\mathbb{G}}$-distributions, provides a corresponding generalization to abstract $S$-groups.

Why not now, at the current stage of our development, at which we have the 1st and 2nd TESS, to proceed in an analogous manner to that above described, formulating precise definitions of two types of distributions, the $\widetilde{\mathscr{G}}(I)$-distributions and the $\overline{\mathscr{G}}(I)$-distributions, and its corresponding derivatives?

We will not only carry this out, that is, define the referred concepts, but we will also repeat, *mutatis mutandis*, what we did after defining (item 2.18) the $\widetilde{\mathbb{G}}$-distributions an its derivatives. More explicit and specifically, we indicate and comment our next steps below.

## 5.3 Our Next Steps

In what follows, we highlight the basic stages of our future considerations.

**(a)** First of all we will prove that every $S$-group $\mathbb{G}$ can be considered as, or identified to, a special $S$-space, denoted by $\mathscr{G}[\mathbb{G}]$ and denominated the $S$-space generated by $\mathbb{G}$. We will see that the two Theorems of Extension of $S$-Spaces, the 1st and 2nd TESS, corresponding to a $S$-space $\mathscr{G}(I) = \mathscr{G}[\mathbb{G}]$ generated by an abelian, surjective, and with identity $S$-group, $\mathbb{G}$, collide into a single theorem which is, essentially, the Theorem of Extension of $S$-Groups for the $S$-group $\mathbb{G}$;

**(b)** Next, taking (a) into account, i.e., that the $S$-group $\mathbb{G}$ and the $S$-space generated by it, $\mathscr{G}[\mathbb{G}]$, are, essentially, the same object, we will consider a domain of distributions, $(\mathbb{G}, \widehat{\mathbb{G}})$, associated to an abelian, surjective, and with identity $S$-group, $\mathbb{G}$, introduced by Definition 2.18(a), as a particular case of the more general notion of domain of distributions defined as an ordered pair of $S$-spaces,

$$\Big(\mathscr{G}(I), \widehat{\mathscr{G}}(I)\Big),$$

not necessarily generated by $S$-groups; as we will see (in Motivation 5.8 and Definition 5.10), one such domain can be of two types or species, as $\widehat{\mathscr{G}}(I)$ is characterized through the 1st or 2nd TESS, from which the $\widetilde{\mathscr{G}}(I)$ and $\overline{\mathscr{G}}(I)$-distributions (and derivatives) will be defined;

**(c)** Basing ourselves on considerations analogous to those made in Remark 2.19, about the possibility of obtaining a categoric axiomatic definition for the $\widetilde{\mathbb{G}}$-distributions and its derivatives, through axioms suggested by the very own theorem (of extension of $S$-groups) that allowed and backed the definition of the referred concepts, we will formulate two axiomatics, inspiring ourselves now into the 1st and 2nd TESS to choose the axioms, one for the $\widetilde{\mathscr{G}}(I)$-distributions and another for the $\overline{\mathscr{G}}(I)$-distributions, as well as we will prove the consistency and categoricity of each one of them;

**(d)** Such as was possible for us to accomplish for the $\widetilde{\mathbb{G}}$-distributions, we will simplify the axiomatics formulated for the $\widetilde{\mathscr{G}}(I)$ and $\overline{\mathscr{G}}(I)$-distributions. More specifically, we will prove that the primitive term "addition", present on both axiomatics in question, can be taken as a defined term, which, in turn, allows the formulation of two simpler equivalent axiomatic versions for the $\widetilde{\mathscr{G}}(I)$ and $\overline{\mathscr{G}}(I)$-distributions;

**(e)** We will revisit the Schwartz' distributions theory and prove that it allows the construction of two $S$-spaces: one of them, that we will denote by $\mathscr{D}'_f(\mathbb{R}^n)$, associated to the $S$-groups of the finite order distributions $(\mathbb{D}'_f(\Omega) = (D'_f(\Omega), D_f(\Omega)), \Omega \in \Gamma(\mathbb{R}^n)$, according to item 2.36), and the other, $\mathscr{D}'(\mathbb{R}^n)$, involving the $S$-groups of all the distributions in $\Omega$, of finite or infinite order $(\mathbb{D}'(\Omega) = (D'(\Omega), D(\Omega)), \Omega \in \mathbb{R}^n$, according to item 2.36); furthermore, and mainly, we will see that the $S$-space $\mathscr{D}'_f(\mathbb{R}^n)$ is a strict and closed extension of the $S$-space $\mathscr{C}(\mathbb{R}^n)$ of the continuous functions and, therefore, by the 1st TESS, it is isomorphic to the extension $\widetilde{\mathscr{C}}(\mathbb{R}^n)$ of $\mathscr{C}(\mathbb{R}^n)$ as described in this theorem, and also that the $S$-space $\mathscr{D}'(\mathbb{R}^n)$ is a locally closed and coherent extension of $\widetilde{\mathscr{C}}(\mathbb{R}^n)$ and, hence, by the 2nd TESS, it is isomorphic to the extension $\overline{\mathscr{C}}(\mathbb{R}^n)$ of $\widetilde{\mathscr{C}}(\mathbb{R}^n)$ defined in the 2nd TESS.

From this, we will then have proved that our definitions of $\widetilde{\mathscr{C}}(\mathbb{R}^n)$-distributions and its derivatives and of $\overline{\mathscr{C}}(\mathbb{R}^n)$-distributions and its corresponding derivatives (that is, the $\widetilde{\mathscr{G}}(I)$ and $\overline{\mathscr{G}}(I)$-distributions translated to the case where $\mathscr{G}(I) = \mathscr{C}(\mathbb{R}^n)$), in any of its versions (axiomatic or via the proper domain of distributions), are equivalent to the ones of finite order distribution and distribution (of all orders), respectively, given by the Schwartz' formulation;

**(f)** Finally, translating the axiomatics formulated for the $\widetilde{\mathscr{G}}(I)$ and $\overline{\mathscr{G}}(I)$-distributions to the particular case associated to the $S$-space $\mathscr{G}(I) = \mathscr{C}(\mathbb{R}^n)$ of the continuous functions, we will obtain simple axiomatics defining, categorically, extensions (generalizations) of the notions of continuous function and differentiability which are, taking into account what is established in (e), equivalent to those provided by Schwartz' distributions theory.

# $\widetilde{\mathscr{G}}(I)$ and $\overline{\mathscr{G}}(I)$-Distributions

## 5.4 Preliminaries

Right after obtaining the Theorem of Extension of $S$-Groups (Theorem 2.16), the considerations in 2.17 suggested that the notions, at that time yet informal, of distribution and derivative of distribution, that was being employed by us to refer to the desired generalizations of the concepts of continuous function and differentiability, were formalized through what the theorem establishes. More explicitly, since the $S$-group $\widetilde{\mathbb{C}}(\Omega) = (\widetilde{C}(\Omega), \widetilde{\partial}(\Omega))$ described in the theorem in question is, unless isomorphism, the only strict and closed extension of the $S$-group $\mathbb{C}(\Omega) = (C(\Omega), \partial(\Omega))$ of the continuous functions, that is, roughly speaking, $\widetilde{C}(\Omega)$ and $\widetilde{\partial}(\Omega)$ are the only extensions of the continuous functions ($C(\Omega)$) and its derivatives ($\partial(\Omega)$) that fulfill conditions which any such extension, in order to be useful, must attend (according to 1.3(a)-(d)), it then becomes natural the formalization of the notions of distribution and its derivatives, informally employed to refer to these extensions, defining them as being, precisely, the elements of the group $\widetilde{C}(\Omega)$ and the semigroup $\widetilde{\partial}(\Omega)$, respectively. Furthermore, and taking into account that the Theorem of Extension of $S$-Groups refers not only to the $S$-group $\mathbb{C}(\Omega) = (C(\Omega), \partial(\Omega))$ of the continuous functions, but to any $S$-group $\mathbb{G} = (G, H)$ that, as $\mathbb{C}(\Omega)$, is abelian, surjective, and with identity, it also results natural that the essentially unique strict and closed extension $\widetilde{\mathbb{G}} = (\widetilde{G}, \widetilde{H})$, of such a $S$-group $\mathbb{G} = (G, H)$, to be taken as a "habitat of distributions" where the distributions or, more adequately, the $\widetilde{\mathbb{G}}$-distributions and its derivatives, seen as extensions (generalizations) not of the continuous functions and its derivatives anymore, but of the elements of the abstract group $G$ and the semigroup H, respectively, would be defined, precise and respectively, as the elements of the group $\widetilde{\mathbb{G}}$ and the endomorphisms of the semigroup $\widetilde{H}$.

These were, basically, the considerations of item 2.17 (above referred) that led us to the definitions of the concept of domain of distribution (Definition 2.18(a)) and the notions, associated with one such domain, of $\widetilde{\mathbb{G}}$-distribution and $\widetilde{\mathbb{G}}$-distribution derivative (Definition 2.18(b)).

Backed now by the defined concepts and obtained results throughout Chapters 3 and 4, particularly by the notions of $S$-space and its extensions and the Theorems of Extension of $S$-Spaces, the 1st and 2nd TESS, we can say that we have enough material for a generalization of the concepts of domain of distribution, $\widetilde{\mathbb{G}}$-distribution and $\widetilde{\mathbb{G}}$-distribution derivative. The next items of this section are then destined to elaborate this claim giving it precise content.

## 5.5 $S$-Spaces Generated by $S$-Groups

Let us observe that the singleton set $\Gamma(I) \coloneqq \{I\}$, with $I$ being any set, trivially attends the defining conditions of an indexer (see Definition 3.6). For this trivial indexer, the associated sets $\Gamma^2(I)$ and $\Delta(I)$ are:

$$\Gamma^2(I) \coloneqq \Big\{(\gamma', \gamma) : \gamma', \gamma \in \Gamma(I)\Big\} = \Big\{(I, I)\Big\}$$

and

$$\Delta(I) \coloneqq \Big\{(\gamma', \gamma) \in \Gamma^2(I) : \gamma' \subseteq \gamma\Big\} = \Big\{(I, I)\Big\} = \Gamma^2(I).$$

Let now $\mathbb{G} = (G, H)$ be an arbitrarily fixed $S$-group and let us consider the following families indexed, respectively, by the sets $\Gamma(I)$, $\Gamma^2(I)$ and $\Delta(I)$:

$$\begin{aligned} \mathbb{G}\Big(\Gamma(I)\Big) &\coloneqq \Big\{\mathbb{G}(\gamma)\Big\}_{\gamma\in\Gamma(I)} = \Big\{\mathbb{G}(I)\Big\}, \\ i\Big(\Gamma^2(I)\Big) &\coloneqq \Big\{i_{(\gamma',\gamma)}\Big\}_{(\gamma',\gamma)\in\Gamma^2(I)} = \Big\{i_{(I,I)}\Big\} \quad \text{and} \\ \Theta\Big(\Delta(I)\Big) &\coloneqq \Big\{\Theta_{(\gamma',\gamma)}\Big\}_{(\gamma',\gamma)\in\Delta(I)} = \Big\{\Theta_{(I,I)}\Big\}, \end{aligned}$$

where $\mathbb{G}(I) \coloneqq \mathbb{G} = (G, H)$, and $i_{(I,I)}$ and $\Theta_{(I,I)}$ are, respectively, the identity functions on $H$ and $G$, that is, $i_{(I,I)} \coloneqq I_H$ and $\Theta_{(I,I)} \coloneqq I_G$, where

$$\begin{aligned} I_H : H &\longrightarrow H \\ \Phi &\longmapsto I_H(\Phi) \coloneqq \Phi \end{aligned} \qquad \text{and} \qquad \begin{aligned} I_G : G &\longrightarrow G \\ g &\longmapsto I_G(g) \coloneqq g, \end{aligned}$$

It immediately results from the pertinent definitions that:

- $\mathbb{G}(\Gamma(I)) = \{\mathbb{G} = (G, H)\}$ is a family of $S$-groups (Definition 3.7);
- $i(\Gamma^2(I)) = \{I_H\}$ is a bonding of the family of $S$-groups $\mathbb{G}(\Gamma(I)) = \{\mathbb{G} = (G, H)\}$ (Definition 3.11);
- $(\mathbb{G}(\Gamma(I)) = \{\mathbb{G} = (G, H)\}, i(\Gamma^2(I)) = \{I_H\})$ is a bonded family (Definition 3.11);

- $\Theta(\Delta(I)) = \{I_G\}$ is a restriction for the bonded family $(\mathbb{G}(\Gamma(I)) = \{\mathbb{G} = (G,H)\}, \{I_H\})$ (Definition 3.13).

Therefore, taking into account the definition of $S$-space (Definition 3.15), one concludes that the triplet

$$\left(\mathbb{G}\Big(\Gamma(I)\Big) = \Big\{\mathbb{G} = (G,H)\Big\}, i\Big(\Gamma^2(I)\Big) = \Big\{I_H\Big\}, \Theta\Big(\Delta(I)\Big) = \Big\{I_G\Big\}\right)$$

is a $S$-space; let us denote this $S$-space by $\mathscr{G}[\mathbb{G}]$, that is,

$$\mathscr{G}[\mathbb{G}] := \left(\Big\{\mathbb{G} = (G,H)\Big\}, \Big\{I_H\Big\}, \Big\{I_G\Big\}\right),$$

and say that $\mathscr{G}[\mathbb{G}]$ is the **$S$-space generated by the $S$-group** $\mathbb{G} = (G,H)$. For obvious reasons, we will make no distinction between a $S$-group $\mathbb{G}$ and the $S$-space $\mathscr{G}[\mathbb{G}]$.

## 5.6 1st TESS Regarding $S$-Spaces Generated by $S$-Groups

Let us now suppose that $\mathbb{G} = (G,H)$ is an abelian, surjective, and with identity $S$-group. In this case, the $S$-space generated by $\mathbb{G}$,

$$\mathscr{G}[\mathbb{G}] := \left(\Big\{\mathbb{G} = (G,H)\Big\}, \Big\{I_H\Big\}, \Big\{I_G\Big\}\right),$$

is, according to Definition 3.15, an abelian, surjective, and with identity $S$-space. By the Theorem of Extension of $S$-Groups (Theorem 2.16), there exists a single (unless isomorphism) strict and closed extension of the $S$-group $\mathbb{G} = (G,H)$, namely, $\widetilde{\mathbb{G}} = (\widetilde{G}, \widetilde{H})$ as defined in the referred theorem. Hence, the family of $S$-groups $\{\widetilde{\mathbb{G}} = (\widetilde{G}, \widetilde{H})\}$ is a strict and closed extension of the family of $S$-groups $\{\mathbb{G} = (G,H)\}$ of the $S$-space $\mathscr{G}[\mathbb{G}]$ (according to Definition 3.8). It immediately results (see Definition 3.22) that the extension of the bonding $\{i_{(I,I)} = I_H\}$ of the family of $S$-groups $\{\mathbb{G} = (G,H)\}$ to the family $\{\widetilde{\mathbb{G}} = (\widetilde{G}, \widetilde{H})\}$, is $\{\widetilde{i}_{(I,I)} = I_{\widetilde{H}}\}$, where $I_{\widetilde{H}} : \widetilde{H} \longrightarrow \widetilde{H}$ is the identity on $\widetilde{H}$ ($I_{\widetilde{H}}(\widetilde{\Phi}) = \widetilde{\Phi}$ for every $\widetilde{\Phi} \in \widetilde{H}$). Hence, the ordered pair

$$\left(\Big\{\widetilde{\mathbb{G}} = (\widetilde{G}, \widetilde{H})\Big\}, \Big\{\widetilde{i}_{(I,I)} = I_{\widetilde{H}}\Big\}\right)$$

is a bonded family that is an extension of the bonded family

$$\left(\Big\{\mathbb{G} = (G,H)\Big\}, \Big\{i_{(I,I)} = I_H\Big\}\right)$$

of the $S$-space $\mathscr{G}[\mathbb{G}]$ (according to Definition 3.22).

Let us suppose now that $\{\widetilde{\Theta}_{(I,I)}\}$ is a restriction for the bonded family

$$\left(\left\{\widetilde{\mathbb{G}} = (\widetilde{G}, \widetilde{H})\right\}, \left\{\widetilde{i}_{(I,I)} = I_{\widetilde{H}}\right\}\right)$$

(see Definition 3.13). In this case, since the $S$-group $\widetilde{\mathbb{G}} = (\widetilde{G}, \widetilde{H})$ is a closed (and strict) extension of $\mathbb{G} = (G, H)$, for $\widetilde{g} \in \widetilde{G}$ arbitrarily fixed we have

$$\widetilde{g} = \widetilde{\Phi}(g),$$

for some $\widetilde{\Phi} \in \widetilde{H}$ and $g \in G$, and hence

$$\widetilde{\Theta}_{(I,I)}(\widetilde{g}) = \widetilde{\Theta}_{(I,I)}\left(\widetilde{\Phi}(g)\right) = \widetilde{i}_{(I,I)}(\widetilde{\Phi})\left(\widetilde{\Theta}_{(I,I)}(g)\right) = \widetilde{\Phi}\left(\widetilde{\Theta}_{(I,I)}(g)\right).$$

Thus, if we suppose that the family $\{\widetilde{\Theta}_{(I,I)}\}$ is, not only a restriction for the bonded family $(\{\widetilde{\mathbb{G}}\}, \{\widetilde{i}_{(I,I)}\})$, but also a prolongation to this bonded family of the restriction $\{\Theta_{(I,I)} = I_G\}$ of the $S$-space $\mathscr{G}[\mathbb{G}]$ (see Definition 3.23(a)), we have

$$\widetilde{\Theta}_{(I,I)}(g) = \Theta_{(I,I)}(g) = g$$

and, therefore,

$$\widetilde{\Theta}_{(I,I)}(\widetilde{g}) = \widetilde{\Phi}\left(\widetilde{\Theta}_{(I,I)}(g)\right) = \widetilde{\Phi}(g) = \widetilde{g}.$$

Since $\widetilde{g} \in \widetilde{G}$ was arbitrarily fixed, we conclude that

$$\widetilde{\Theta}_{(I,I)} = I_{\widetilde{G}},$$

where $I_{\widetilde{G}} : \widetilde{G} \longrightarrow \widetilde{G}$ is the identity function (on $\widetilde{G}$).

Reciprocally, it is trivial to prove that the family $\{\widetilde{\Theta}_{(I,I)} := I_{\widetilde{G}}\}$ is a restriction for the bonded family $(\{\widetilde{\mathbb{G}}\} = (\widetilde{G}, \widetilde{H}), \{\widetilde{i}_{(I,I)} = I_{\widetilde{H}}\})$ and also a prolongation of $\{\Theta_{(I,I)} = I_G\}$ to this bonded family.

Taking now into account the definition of extension of a $S$-space (Definition 3.23(b)), as well as 3.23(c-1), the conclusions above allow us to say that the triplet

$$\left(\left\{\widetilde{\mathbb{G}} = (\widetilde{G}, \widetilde{H})\right\}, \left\{\widetilde{i}_{(I,I)} = I_{\widetilde{H}}\right\}, \left\{\widetilde{\Theta}_{(I,I)} = I_{\widetilde{G}}\right\}\right)$$

is a strict and closed extension of the $S$-space $\mathscr{G}[\mathbb{G}]$ generated by $\mathbb{G}$. But, the triplet above is the $S$-space generated by the $S$-group $\widetilde{\mathbb{G}} = (\widetilde{G}, \widetilde{H})$ (according to 5.5), that is,

$$\mathscr{G}[\widetilde{\mathbb{G}}] = \left(\left\{\widetilde{\mathbb{G}} = (\widetilde{G}, \widetilde{H})\right\}, \left\{\widetilde{i}_{(I,I)} = I_{\widetilde{H}}\right\}, \left\{\widetilde{\Theta}_{(I,I)} = I_{\widetilde{G}}\right\}\right).$$

Since $\mathscr{G}[\mathbb{G}]$ is an abelian, surjective, and with identity $S$-space, then, by the 1st TESS (Theorem 4.7), $\mathscr{G}[\widetilde{\mathbb{G}}]$ is the only (unless isomorphism) strict and closed extension of $\mathscr{G}[\mathbb{G}]$ and, hence, based on the natural identification of the $S$-groups with the respective $S$-spaces generated by them, one concludes that the 1st TESS regarding a $S$-space generated by an abelian, surjective, and with identity $S$-group is, essentially, the Theorem of Extension of $S$-Groups.

## 5.7 2nd TESS Regarding $S$-Spaces Generated by $S$-Groups

Let us note now that the $S$-spaces generated by $S$-groups are, all, trivially coherent, that is, they trivially fulfil the coherency condition (see Definition 3.19). Furthermore, if $\widehat{\mathbb{G}} = (\widehat{G}, \widehat{H})$ is an extension of a $S$-group $\mathbb{G} = (G, H)$ arbitrarily fixed, then, $\mathscr{G}[\widehat{\mathbb{G}}]$, the $S$-space generated by $\widehat{\mathbb{G}}$, is, also trivially, a locally closed extension of itself (see Definition 3.23(c-2)).

Thus, the $S$-space $\mathscr{G}[\widetilde{\mathbb{G}}]$ generated by the (essentially unique) strict and closed extension, $\widetilde{\mathbb{G}}$, of an abelian, surjective, and with identity $S$-group, $\mathbb{G}$, is a coherent and locally closed extension of itself. But, as we saw, $\mathscr{G}[\widetilde{\mathbb{G}}]$ is the strict and closed extension of the abelian, surjective, with identity, and (trivially) coherent $S$-space $\mathscr{G}[\mathbb{G}]$, and since, by the 2nd TESS, there exists a single (unless isomorphism) extension $\overline{\mathscr{G}}[\mathbb{G}]$ of $\mathscr{G}[\widetilde{\mathbb{G}}]$ which is coherent and locally closed, we conclude that $\overline{\mathscr{G}}[\mathbb{G}] = \mathscr{G}[\widetilde{\mathbb{G}}]$ or, more precisely, $\overline{\mathscr{G}}[\mathbb{G}]$ is isomorphic to $\mathscr{G}[\widetilde{\mathbb{G}}]$.

In short, regarding $S$-spaces generated by $S$-groups, the Theorems of Extension of $S$-spaces, that is, the 1st and 2nd TESS (Theorems 4.7 and 4.33), collide into a single theorem which is, essentially, the Theorem of Extension of $S$-Groups (Theorem 2.16).

## 5.8 Motivation

A domain of distribution was defined in Chapter 2 as an ordered pair

$$(\mathbb{G}, \widehat{\mathbb{G}}),$$

where $\mathbb{G} = (G, H)$ is an abelian, surjective, and with identity $S$-group and $\widehat{\mathbb{G}} = (\widehat{G}, \widehat{H})$ is a $S$-group isomorphic to the strict and closed extension of $\mathbb{G}$, $\widetilde{\mathbb{G}} = (\widetilde{G}, \widetilde{H})$, as described in the Theorem of Extension of $S$-Groups (Theorem 2.16); associated to such a domain $(\mathbb{G}, \widehat{\mathbb{G}})$ were also defined, in the referred chapter, the $\widetilde{\mathbb{G}}$-distributions and its derivatives as being, respectively, the elements of the group $\widehat{G}$ and the endomorphisms of the semigroup $\widehat{H}$ (Definition 2.18(b)). Taking now into account that the $S$-groups $\mathbb{G} = (G, H)$ and $\widehat{\mathbb{G}} = (\widehat{G}, \widehat{H})$, as we highlighted in 5.5, do not differ, essentially, from its respective generated $S$-spaces, $\mathscr{G}[\mathbb{G}]$ and $\mathscr{G}[\widehat{\mathbb{G}}]$, also the domain of distribution $(\mathbb{G}, \widehat{\mathbb{G}})$ does not differ, to the same extend, from the ordered pair $(\mathscr{G}[\mathbb{G}], \mathscr{G}[\widehat{\mathbb{G}}])$.

The observations above suggest a natural generalization of the notion of domain of distribution (and of its corresponding concepts of distribution and derivative of distribution), where one such domain would be defined as an ordered pair of $S$-spaces (generated or not by $S$-groups),

$$\Big(\mathscr{G}(I), \widehat{\mathscr{G}}(I)\Big),$$

being $\mathscr{G}(I)$ abelian, surjective, and with identity (in particular, $\mathscr{G}(I) = \mathscr{G}[\mathbb{G}] \equiv \mathbb{G}$, where $\mathbb{G}$ is an abelian, surjective, and with identity $S$-group) and $\widehat{\mathscr{G}}(I)$ an extension of $\mathscr{G}(I)$ such that:

**(a)** is isomorphic to the extension $\widetilde{\mathscr{G}}(I)$ of $\mathscr{G}(I)$ defined in the 1st TESS or,

**(b)** is isomorphic to the extension $\overline{\mathscr{G}}(I)$ of $\mathscr{G}(I)$ defined in the 2nd TESS[21], on the hypothesis of $\mathscr{G}(I)$ be coherent (beyond abelian, surjective, and with identity).

Hence, in the context of this generalization, we have not only one, but two "types" of domain of distribution: the ones of "1st species", $(\mathscr{G}(I), \widehat{\mathscr{G}}(I) \equiv \widetilde{\mathscr{G}}(I))$, with $\widehat{\mathscr{G}}(I)$ as described in (a), and the ones of "2nd species", $(\mathscr{G}(I), \widehat{\mathscr{G}}(I) \equiv \overline{\mathscr{G}}(I))$, with $\widehat{\mathscr{G}}(I)$ as specified in (b). Clearly, for $S$-spaces $\mathscr{G}(I)$ generated by abelian, surjective, and with identity $S$-groups, for which, as we saw, $\widetilde{\mathscr{G}}(I)$ is isomorphic to $\overline{\mathscr{G}}(I)$, the two species of domain of distribution merge into a single notion which is, essentially, the one originally introduced in Chapter 2 through Definition 2.18(a).

In what follows, in item 5.10, motivated and guided by the considerations above, we formulate the due definitions. Before, however, we present, in 5.9, a simple lemma that, as we will see, will back one of the concepts defined in 5.10.

## 5.9 Lemma

*Let*

$$\mathbb{G}\Big(\Gamma(I)\Big) = \Big\{\mathbb{G}(\gamma) = \Big(G(\gamma), H(\gamma)\Big)\Big\}_{\gamma\in\Gamma(I)}$$

*be a family of $S$-groups such that*

$$G(\gamma) \cap G(\gamma') = \varnothing$$

*for every $\gamma, \gamma' \in \Gamma(I)$ such that $\gamma \neq \gamma'$. If*

$$\widehat{\mathbb{G}}\Big(\Gamma(I)\Big) = \Big\{\widehat{\mathbb{G}}(\gamma) = \Big(\widehat{G}(\gamma), \widehat{H}(\gamma)\Big)\Big\}_{\gamma\in\Gamma(I)}$$

*is an extension of the family $\mathbb{G}(\Gamma(I))$, then,*

$$\widehat{G}(\gamma) \cap \widehat{G}(\gamma') = \varnothing$$

*for every $\gamma, \gamma' \in \Gamma(I)$ such that $\gamma \neq \gamma'$.*

[21] Worth remember here that $\overline{\mathscr{G}}(I)$ is an extension of the extension $\widetilde{\mathscr{G}}(I)$ of $\mathscr{G}(I)$ and, therefore, it also is an extension of $\mathscr{G}(I)$, as established at page 120.

*Proof.* In fact, for if $\gamma, \gamma' \in \Gamma(I)$ such that $\gamma \neq \gamma'$ and $\widehat{G}(\gamma) \cap \widehat{G}(\gamma') \neq \varnothing$ do exist, then, for $\widehat{g} \in \widehat{G}(\gamma) \cap \widehat{G}(\gamma')$ we would have:

$$\widehat{g} - \widehat{g} = \widehat{0}_\gamma \in \widehat{G}(\gamma) \quad \text{and} \quad \widehat{g} - \widehat{g} = \widehat{0}_{\gamma'} \in \widehat{G}(\gamma')$$

where $\widehat{0}_\gamma$ and $\widehat{0}_{\gamma'}$ are the additive neutrals of the groups $\widehat{G}(\gamma)$ and $\widehat{G}(\gamma')$, respectively. On the other hand, $G(\gamma)$ and $G(\gamma')$ are subgroups of $\widehat{G}(\gamma)$ and $\widehat{G}(\gamma')$, respectively, once that the family of $S$-groups $\widehat{\mathbb{G}}(\Gamma(I))$ is an extension of the family $\mathbb{G}(\Gamma(I))$. Therefore, $G(\gamma)$ and $\widehat{G}(\gamma)$ have the same additive neutral, as well as $G(\gamma')$ and $\widehat{G}(\gamma')$, that is, $\widehat{g} - \widehat{g} \in G(\gamma)$ and, also, $\widehat{g} - \widehat{g} \in G(\gamma')$. Consequently,

$$\widehat{g} - \widehat{g} \in G(\gamma) \cap G(\gamma')$$

and, hence, $G(\gamma) \cap G(\gamma') \neq \varnothing$. ■

## 5.10 Definition

**(a)** Let

$$\mathscr{G}(I) = \left(\mathbb{G}\Big(\Gamma(I)\Big) = \Big\{\mathbb{G}(\gamma) = \Big(G(\gamma), H(\gamma)\Big)\Big\}_{\gamma\in\Gamma(I)}, i\Big(\Gamma^2(I)\Big), \Theta\Big(\Delta(I)\Big)\right)$$

be an arbitrarily fixed $S$-space such that

$$G(\gamma) \cap G(\gamma') = \varnothing$$

for every $\gamma, \gamma' \in \Gamma(I)$ such that $\gamma \neq \gamma'$.

**$\mathscr{G}(I)$-individual** or, simply, **individual** if there is not risk of ambiguity, is, by definition, any element of the groups of the $S$-space $\mathscr{G}(I)$; hence, $g$ is a $\mathscr{G}(I)$-individual, if only if, there exists $\gamma \in \Gamma(I)$ such that $g \in G(\gamma)$. Since the groups of the $S$-space $\mathscr{G}(I)$ are two by two disjoints, to each individual $g$ corresponds a single $\gamma \in \Gamma(I)$ such that $g \in G(\gamma)$, that will be denominated the **domain of (the individual)** $g$.

**(b)** Let now

$$\mathscr{G}(I) = \left(\mathbb{G}\Big(\Gamma(I)\Big) = \Big\{\mathbb{G}(\gamma) = \Big(G(\gamma), H(\gamma)\Big)\Big\}_{\gamma\in\Gamma(I)}, i\Big(\Gamma^2(I)\Big), \Theta\Big(\Delta(I)\Big)\right)$$

be an abelian, surjective, and with identity $S$-space arbitrarily fixed and such that

$$G(\gamma) \cap G(\gamma') = \varnothing$$

for every $\gamma, \gamma' \in \Gamma(I)$ such that $\gamma \neq \gamma'$. We say that the ordered pair

$$\Big(\mathscr{G}(I), \widehat{\mathscr{G}}(I)\Big)$$

is

**(b-1)** a **1st species domain of distribution**, if and only if, $\widehat{\mathscr{G}}(I)$ is isomorphic to the strict and closed extension $\widetilde{\mathscr{G}}(I)$ of $\mathscr{G}(I)$, as defined in the 1st TESS;

**(b-2)** a **2nd species domain of distribution**, if and only if, $\mathscr{G}(I)$ as above described also is coherent and $\widehat{\mathscr{G}}(I)$ is isomorphic to the locally closed and coherent extension $\overline{\mathscr{G}}(I)$ of $\widetilde{\mathscr{G}}(I)$, as defined in the 2nd TESS.

**(c)** Regarding a 1st species domain of distribution, $(\mathscr{G}(I), \widehat{\mathscr{G}}(I))$, we define:

**(c-1)** $\widetilde{\mathscr{G}}(I)$**-distribution** is any element $\widehat{g}$ of anyone of the groups $\widehat{G}(\gamma)$ of the $S$-space $\widehat{\mathscr{G}}(I)$, i.e., a $\widehat{\mathscr{G}}(I)$-individual;

**(c-2)** if $\widehat{g}$ is a $\widetilde{\mathscr{G}}(I)$-distribution, it results, taking Lemma 5.9 into account, that there exists a single $\gamma \in \Gamma(I)$ such that $\widehat{g} \in \widehat{G}(\gamma)$, the domain of the individual $\widehat{g}$ (as defined in (a)), that will be called (too) the **domain of the** $\widetilde{\mathscr{G}}(I)$**-distribution** $\widehat{g}$. Hence, $\widehat{G}(\gamma)$ is the class of all $\widetilde{\mathscr{G}}(I)$-distributions with domain $\gamma$;

**(c-3)** **derivative of** $\widetilde{\mathscr{G}}(I)$**-distribution** is any endomorphism $\widehat{\Phi}$ of anyone of the semigroups $\widehat{H}(\gamma)$ of the $S$-space $\widehat{\mathscr{G}}(I)$; we will also say that the members $\widehat{\Phi}, \widehat{\Psi}, \ldots$ of $\widehat{H}(\gamma)$ are derivatives of the $\widetilde{\mathscr{G}}(I)$-distributions with domain $\gamma$ (since $\widehat{\Phi}, \widehat{\Psi}, \ldots \in \widehat{H}(\gamma)$ are endomorphisms on the group $\widehat{G}(\gamma)$, the class of all $\widetilde{\mathscr{G}}(I)$-distributions with domain $\gamma$);

**(c-4)** **restriction of** $\widetilde{\mathscr{G}}(I)$**-distribution** is any homomorphism $\widehat{\Theta}_{(\gamma',\gamma)} \in \widehat{\Theta}(\Delta(I))$ of the family $\widehat{\Theta}(\Delta(I))$ of the $S$-space $\widehat{\mathscr{G}}(I)$.

**(d)** Regarding a 2nd species domain of distribution, $(\mathscr{G}(I), \widehat{\mathscr{G}}(I))$, we define as in (c), *mutatis mutandis*, the notions of $\overline{\mathscr{G}}(I)$**-distribution**, **domain**, **derivative** and **restriction of** $\overline{\mathscr{G}}(I)$**-distribution**.

## 5.11 Remarks

To each $S$-Space $\mathscr{G}(I)$ that is abelian, surjective, with identity, coherent and whose constituting groups are two by two disjoints, corresponds, through the Definitions 5.10(c) and (d), two classes of concepts: one of them composed by the notions of $\widetilde{\mathscr{G}}(I)$-distribution, domain, derivative and restriction of $\widetilde{\mathscr{G}}(I)$-distribution regarding a 1st species domain of distribution, $(\mathscr{G}(I), \widehat{\mathscr{G}}(I))$, and the other one composed of the corresponding notions regarding a 2nd species domain of distribution, that is, the concepts of $\overline{\mathscr{G}}(I)$-distribution, domain, derivative and restriction of $\overline{\mathscr{G}}(I)$-distribution.

As we know, if $\mathbb{G} = (G, H)$ is an abelian, surjective, and with identity $S$-group, then, the $S$-space generated by $\mathbb{G}$,

$$\mathscr{G}[\mathbb{G}] := \Big(\Big\{\mathbb{G} = (G, H)\Big\}, \Big\{i_{(I,I)} = I_H\Big\}, \Big\{\Theta_{(I,I)} = I_G\Big\}\Big),$$

is an abelian, surjective, with identity, (trivially) coherent, and whose constituting groups (only one, $G$) are (trivially) two by two disjoints. Hence, we can take $\mathscr{G}(I) = \mathscr{G}[\mathbb{G}]$ in the definition under analysis (Definition 5.10) and, for such a $\mathscr{G}(I)$, as we saw in 5.7, the strict and closed extension, $\widetilde{\mathscr{G}}(I)$, of $\mathscr{G}(I)$, given by the 1st TESS, is equal (isomorphic) to the locally closed and coherent extension, $\overline{\mathscr{G}}(I)$, of $\widetilde{\mathscr{G}}(I)$, provided by the 2nd TESS, and both are equal (isomorphic) to the $S$-space $\mathscr{G}[\widetilde{\mathbb{G}}]$ generated by the strict and closed extension, $\widetilde{\mathbb{G}}$, of the $S$-group $\mathbb{G}$:

$$\widetilde{\mathscr{G}}(I) = \overline{\mathscr{G}}(I) = \mathscr{G}[\widetilde{\mathbb{G}}] = \Big(\Big\{\widetilde{\mathbb{G}} = (\widetilde{G}, \widetilde{H})\Big\}, \Big\{\widetilde{i}_{(I,I)} = I_{\widetilde{H}}\Big\}, \Big\{\widetilde{\Theta}_{(I,I)} = I_{\widetilde{G}}\Big\}\Big).$$

Hence, in this case, where $\mathscr{G}(I) = \mathscr{G}[\mathbb{G}]$, corresponds a single domain of distribution, namely,

$$\Big(\mathscr{G}(I) = \mathscr{G}[\mathbb{G}], \widehat{\mathscr{G}}(I)\Big)$$

where $\widehat{\mathscr{G}}(I)$ is isomorphic to the $S$-space $\widetilde{\mathscr{G}}(I) = \mathscr{G}[\widetilde{\mathbb{G}}]$, that is, $\widehat{\mathscr{G}}(I) = \mathscr{G}[\widehat{\mathbb{G}}]$ where $\widehat{\mathbb{G}} = (\widehat{G}, \widehat{H})$ is a $S$-group isomorphic to $\widetilde{\mathbb{G}} = (\widetilde{G}, \widetilde{H})$ (the strict and closed extension of $\mathbb{G} = (G, H)$ defined in the Theorem of Extension of $S$-Groups). Consequently, the two classes of concepts above mentioned, associated with $\mathscr{G}(I)$, collide into a single class if $\mathscr{G}(I) = \mathscr{G}[\mathbb{G}]$, composed of the concepts of $\mathscr{G}[\widetilde{\mathbb{G}}]$-distribution, domain, derivative and restriction of $\mathscr{G}[\widetilde{\mathbb{G}}]$-distribution, whose definitions given in 5.10(c) assume the following form when translated (particularized) to this case.

- Definition 5.10(c): translated for $\mathscr{G}(I) = \mathscr{G}[\mathbb{G}]$.

  Regarding a domain of distribution

  $$\Big(\mathscr{G}[\mathbb{G}], \mathscr{G}[\widehat{\mathbb{G}}]\Big)$$

  where $\mathbb{G} = (G, H)$ is an abelian, surjective, and with identity $S$-group and $\widehat{\mathbb{G}} = (\widehat{G}, \widehat{H})$ is a $S$-group isomorphic to the strict and closed extension of $\mathbb{G}$, $\widetilde{\mathbb{G}} = (\widetilde{G}, \widetilde{H})$, described in the Theorem of Extension of $S$-Groups, we define:

  **(c-1)** $\widetilde{\mathscr{G}}[\mathbb{G}]$-distribution is any element $\widehat{g}$ of the group $\widehat{G}$;

**(c-2)** if $\widehat{g}$ is a $\widetilde{\mathscr{G}}[\mathbb{G}]$-distribution, the domain of $\widehat{g}$ is the only $\gamma \in \Gamma(I)$ such that $\widehat{g} \in \widehat{G}(\gamma)$.

Since the indexer of the $S$-spaces $\mathscr{G}[\mathbb{G}]$ and $\mathscr{G}[\widehat{\mathbb{G}}]$ is $\Gamma(I) = \{I\}$ and $\widehat{G}(I) = \widehat{G}$, it results that all the $\widetilde{\mathscr{G}}[\mathbb{G}]$-distributions, that is, every $\widehat{g} \in \widehat{G}$, have the same domain, namely, the set $I$ which, in turn, is an arbitrary set; we then see that the concept of domain of a distribution is, in this case where $\mathscr{G}(I) = \mathscr{G}[\mathbb{G}]$, artificial and unnecessary;

**(c-3)** derivative of $\widetilde{\mathscr{G}}[\mathbb{G}]$-distribution is any endomorphism $\widehat{\Phi} \in \widehat{H}$;

**(c-4)** restriction of $\widetilde{\mathscr{G}}[\mathbb{G}]$-distribution is any member of the family $\widehat{\Theta}(\Delta(I))$ of the $S$-space $\mathscr{G}[\widehat{\mathbb{G}}]$.

Since $\widehat{\Theta}(\Delta(I)) = \{I_{\widehat{G}}\}$, we see that "restriction of $\widetilde{\mathscr{G}}[\mathbb{G}]$-distribution" and "identity function on $\widehat{G}$" are perfect synonyms and, hence, such as the concept of domain of distribution, the one of restriction of $\widetilde{\mathscr{G}}[\mathbb{G}]$-distribution also is, in this case where $\mathscr{G}(I) = \mathscr{G}[\mathbb{G}]$, artificial and unnecessary.

In short, of the concepts associated to a domain of distribution introduced by Definition 5.10(c) and (d), when referred to a domain $(\mathscr{G}[\mathbb{G}], \mathscr{G}[\widehat{\mathbb{G}}])$ determined by an $S$-space $\mathscr{G}[\mathbb{G}]$ generated by an abelian, surjective, and with identity $S$-group $\mathbb{G} = (G, H)$, only two are relevant: the ones of $\widetilde{\mathscr{G}}[\mathbb{G}]$-distribution and $\widetilde{\mathscr{G}}[\mathbb{G}]$-distribution derivative. Furthermore, the definitions of these two concepts are exactly and respectively equal to those of $\widetilde{\mathbb{G}}$-distribution and $\widetilde{\mathbb{G}}$-distribution derivative formulated in Definition 2.18(b). This "coincidence" is consistent with (beyond reinforcing it) the identification of the $S$-groups $\mathbb{G}$ and $\widehat{\mathbb{G}}$ with the respective $S$-spaces generated by them, $\mathscr{G}[\mathbb{G}]$ and $\mathscr{G}[\widehat{\mathbb{G}}]$, through which we identify the domain of distribution $(\mathscr{G}[\mathbb{G}], \mathscr{G}[\widehat{\mathbb{G}}])$ with the domain $(\mathbb{G}, \widehat{\mathbb{G}})$ introduced by Definition 2.18(a).

# $\widetilde{\mathscr{G}}(I)$-Distributions Axiomatics

## 5.12 $\widetilde{\mathscr{G}}(I)$-Distributions Axioms: 1st Version

Just like the $\widetilde{\mathbb{G}}$-distributions and its derivatives, defined through the concept of domain of distribution (Definition 2.18), admit a categoric axiomatic definition, as was seen in Chapter 2 (items 2.20 to 2.22), so do the $\widetilde{\mathscr{G}}(I)$-distributions and its derivatives, defined as they were, in 5.10(c), via a generalization of the referred concept of domain of distribution, can be introduced through a categoric axiomatic as we will see next (in this item and the three next ones, 5.13 to 5.15).[22]

For the $\widetilde{\mathbb{G}}$-distributions, guided by the considerations presented in Remark 2.19 (review this item), we "extracted" the axioms from the Theorem of Extension of $S$-Groups. Now, for the $\widetilde{\mathscr{G}}(I)$-distributions, based on analogous observations to those in 2.19 translated to this more general case, the choice of the axioms will be driven by the 1st TESS that, as we saw, is a generalization of the Theorem of Extension of $S$-Groups; roughly speaking, assigning the "status" of axioms to the defining predicates of the concepts of extension, strict extension, and closed extension regarding, not an $S$-group anymore, but now, a $S$-space, we will get an axiomatic whose primitive terms implicitly defined by it are $\widetilde{\mathscr{G}}(I)$-distribution, domain, addition, derivative and restriction of $\widetilde{\mathscr{G}}(I)$-distributions. Once obtained, this axiomatic, as we will see, is consistent, categoric, and the meaning assigned by the axioms to its primitive terms, coincide with those imputed to these same terms by Definition 5.10(c) (regarding a 1st species domain of distribution).

With the considerations above, we hope to have provided a heuristic indication that allows one to consider reasonable the following axiomatic for the $\widetilde{\mathscr{G}}(I)$-distributions:

- **Primitive Terms:** $\widetilde{\mathscr{G}}(I)$-distribution, domain, addition, derivative and restriction of $\widetilde{\mathscr{G}}(I)$-distribution;
- **Precedent Theories:** Classical Logic, Set Theory and $S$-Spaces Theory;
- **Axioms:** the ones formulated ahead, regarding an abelian, surjective, and with identity $S$-space

$$\mathscr{G}(I) = \Bigg(\mathbb{G}\Big(\Gamma(I)\Big) = \Big\{\mathbb{G}(\gamma) = \Big(G(\gamma), H(\gamma)\Big)\Big\}_{\gamma\in\Gamma(I)},$$
$$i\Big(\Gamma^2(I)\Big) = \Big\{i_{(\gamma',\gamma)}\Big\}_{(\gamma',\gamma)\in\Gamma^2(I)},$$
$$\Theta\Big(\Delta(I)\Big) = \Big\{\Theta_{(\gamma',\gamma)}\Big\}_{(\gamma',\gamma)\in\Delta(I)}\Bigg),$$

arbitrarily fixed, such that $G(\gamma) \cap G(\gamma') = \varnothing$ if $\gamma \neq \gamma'$.

**Axiom 1** Every $\mathscr{G}(I)$-individual is a $\widetilde{\mathscr{G}}(I)$-distribution.

**Axiom 2** To each $\widetilde{\mathscr{G}}(I)$-distribution, $\widehat{g}$, corresponds a single $\gamma \in \Gamma(I)$, denominated the domain of $\widehat{g}$, in such a way that if $\widehat{g}$ is a $\mathscr{G}(I)$-individual, the domain $\gamma$ is the domain of the $\mathscr{G}(I)$-individual $\widehat{g}$.

[22] The same occurs regarding the $\overline{\mathscr{G}}(I)$-distributions, as we will see in the section "$\overline{\mathscr{G}}(I)$-Distributions Axiomatics".

**Axiom 3** The addition is a function (operation) that to each pair $(\widehat{g}, \widehat{h})$ of $\widetilde{\mathscr{G}}(I)$-distributions with the same domain $\gamma$, associates a (single) $\widetilde{\mathscr{G}}(I)$-distribution with the same domain $\gamma$, denoted by $\widehat{g} + \widehat{h}$ and denominated the sum of $\widehat{g}$ with $\widehat{h}$, in such a way that if $\widehat{g}$ and $\widehat{h}$ are $\mathscr{G}(I)$-individuals, then, the sum $\widehat{g} + \widehat{h}$ is the same as the one obtained through the addition of the group $G(\gamma)$ to which belong the $\mathscr{G}(I)$-individuals $\widehat{g}$ and $\widehat{h}$.

**Notation.** We will denote by $\widehat{G}(\gamma)$ the class of all the $\widetilde{\mathscr{G}}(I)$-distributions with domain $\gamma \in \Gamma(I)$. It results from Axiom 1 and Axiom 2 that $G(\gamma)$, the class of all the $\mathscr{G}(I)$-individuals with domain $\gamma$, is such that $G(\gamma) \subseteq \widehat{G}(\gamma)$.

It is convenient to remember here that, for $\Phi \in H(\gamma)$, the domain of the homomorphism $\Phi$, denoted by $G(\gamma)_\Phi$ (according to the notation introduced in Definition 1.4), is a subgroup of the group $G(\gamma)$ and, therefore, it is included in $\widehat{G}(\gamma)$. As always, we will say that a function $\widehat{\Phi} : \widehat{G}(\gamma) \longrightarrow \widehat{G}(\gamma)$ is an extension of $\Phi : G(\gamma)_\Phi \longrightarrow G(\gamma)$ to $\widehat{G}(\gamma)$, if and only if,

$$\widehat{\Phi}(g) = \Phi(g)$$

for every $g \in G(\gamma)_\Phi \subseteq G(\gamma) \subseteq \widehat{G}(\gamma)$.

**Axiom 4** The derivatives of the $\widetilde{\mathscr{G}}(I)$-distributions with domain $\gamma$, whose class is denoted by $\widehat{H}(\gamma)$, are functions with domain and codomain equal to $\widehat{G}(\gamma)$, such that:

**(4-1)** $\widehat{\Phi}(\widehat{g}+\widehat{h}) = \widehat{\Phi}(\widehat{g})+\widehat{\Phi}(\widehat{h})$ for every $\widehat{g}, \widehat{h} \in \widehat{G}(\gamma)$ and every derivative $\widehat{\Phi} \in \widehat{H}(\gamma)$;

**(4-2)** each derivative $\widehat{\Phi} \in \widehat{H}(\gamma)$ is an extension to $\widehat{G}(\gamma)$ of a single $\Phi \in H(\gamma)$, that is, $\widehat{\Phi}(g) = \Phi(g)$ for every $g \in G(\gamma)_\Phi \subseteq \widehat{G}(\gamma)$; conversely, for each $\Phi \in H(\gamma)$ there exists a single derivative $\widehat{\Phi} \in \widehat{H}(\gamma)$ such that $\Phi(g) = \widehat{\Phi}(g)$ for every $g \in G(\gamma)_\Phi \subseteq \widehat{G}(\gamma)$;

**(4-3)** the class $\widehat{H}(\gamma)$, equipped with the usual operation of composition of functions, is a semigroup isomorphic to the semigroup $H(\gamma)$, with the function

$$\begin{aligned} \widehat{\gamma} : H(\gamma) &\longrightarrow \widehat{H}(\gamma) \\ \Phi &\longmapsto \widehat{\gamma}(\Phi) := \text{the extension of } \Phi \text{ to } \widehat{G}(\gamma) \end{aligned}$$

as an isomorphism.

**Notation.** We will denote the image of $\Phi \in H(\gamma)$ by the isomorphism $\widehat{\gamma} : H(\gamma) \longrightarrow \widehat{H}(\gamma)$ as $\widehat{\Phi}$, that is,

$$\widehat{\Phi} := \widehat{\gamma}(\Phi).$$

**Axiom 5** Let $\gamma \in \Gamma(I)$ and $\Phi \in H(\gamma)$ arbitrarily fixed. If $g \in G(\gamma)$ is such that $g \notin G(\gamma)_\Phi$, then, $\widehat{\Phi}(g) \notin G(\gamma)$.

**Axiom 6** Let $\gamma \in \Gamma(I)$ be arbitrarily fixed. For each $\widehat{g} \in \widehat{G}(\gamma)$ there exist $\Phi \in H(\gamma)$ and $g \in G(\gamma)$ such that

$$\widehat{g} = \widehat{\Phi}(g).$$

In other terms, every $\widetilde{\mathscr{G}}(I)$-distribution is a derivative of some $\mathscr{G}(I)$-individual.

The next axiom makes reference to a term that is not a primitive term, namely, $\widehat{i}_{(\gamma',\gamma)}$. Therefore, it is necessary to define it through primitive ones.

**Definition.** For each $(\gamma', \gamma) \in \Gamma^2(I)$ we define $\widehat{i}_{(\gamma',\gamma)}$ as the following function:

$$\widehat{i}_{(\gamma',\gamma)} : \widehat{H}(\gamma) \longrightarrow \widehat{H}(\gamma')$$

$$\widehat{\Phi} \longmapsto \widehat{i}_{(\gamma',\gamma)}(\widehat{\Phi}) := \widehat{\gamma'}\left( i_{(\gamma',\gamma)}\Big( (\widehat{\gamma})^{-1}(\widehat{\Phi}) \Big) \right)$$

It immediately results that $\widehat{i}_{(\gamma',\gamma)}$ is an isomorphism from the semigroup $\widehat{H}(\gamma)$ onto the semigroup $\widehat{H}(\gamma')$.

**Axiom 7** The restrictions of $\widetilde{\mathscr{G}}(I)$-distributions are functions $\widehat{\Theta}_{(\gamma',\gamma)} : \widehat{G}(\gamma) \longrightarrow \widehat{G}(\gamma')$, a single one for each $(\gamma', \gamma) \in \Delta(I)$, that associate to each $\widetilde{\mathscr{G}}(I)$-distribution with domain $\gamma$, a single $\widetilde{\mathscr{G}}(I)$-distribution with domain $\gamma' \subseteq \gamma$, in such a way that:

**(7-1)** $\widehat{\Theta}_{(\gamma',\gamma)}(\widehat{g} + \widehat{h}) = \widehat{\Theta}_{(\gamma',\gamma)}(\widehat{g}) + \widehat{\Theta}_{(\gamma',\gamma)}(\widehat{h})$ for every $\widehat{g}, \widehat{h} \in \widehat{G}(\gamma)$;

**(7-2)** $\widehat{\Theta}_{(\gamma',\gamma)}(g) = \Theta_{(\gamma',\gamma)}(g)$ for every $g \in G(\gamma) \subseteq \widehat{G}(\gamma)$;

**(7-3)** $\widehat{\Theta}_{(\gamma'',\gamma')}\Big(\widehat{\Theta}_{(\gamma',\gamma)}(\widehat{g})\Big) = \widehat{\Theta}_{(\gamma'',\gamma)}(\widehat{g})$ for every $\widehat{g} \in \widehat{G}(\gamma)$ and $\gamma, \gamma', \gamma'' \in \Gamma(I)$ such that $\gamma'' \subseteq \gamma' \subseteq \gamma$;

**(7-4)** $\widehat{\Theta}_{(\gamma',\gamma)}\Big(\widehat{\Phi}(\widehat{g})\Big) = \widehat{i}_{(\gamma',\gamma)}(\widehat{\Phi})\Big(\widehat{\Theta}_{(\gamma',\gamma)}(\widehat{g})\Big)$ for every $\widehat{g} \in \widehat{G}(\gamma)$ and every $\widehat{\Phi} \in \widehat{H}(\gamma)$.

## 5.13 The Consistency of the Axioms

The $S$-space $\mathscr{G}(I)$, in relation to which the axioms in item 5.12 are formulated, is abelian, surjective, and with identity. Therefore, the $S$-space

$$\widetilde{\mathscr{G}}(I) = \left( \widetilde{\mathbb{G}}\Big(\Gamma(I)\Big) = \Big\{ \widetilde{\mathbb{G}}(\gamma) = \Big(\widetilde{G}(\gamma), \widetilde{H}(\gamma)\Big) \Big\}_{\gamma \in \Gamma(I)}, \right.$$

$$\widetilde{i}\Big(\Gamma^2(I)\Big) = \Big\{ \widetilde{i}_{(\gamma',\gamma)} \Big\}_{(\gamma',\gamma) \in \Gamma^2(I)},$$

$$\left. \widetilde{\Theta}\Big(\Delta(I)\Big) = \Big\{ \widetilde{\Theta}_{(\gamma',\gamma)} \Big\}_{(\gamma',\gamma) \in \Delta(I)} \right),$$

as defined in the 1st TESS (Theorem 4.7), is a strict and closed extension of the $S$-space $\mathscr{G}(I)$. Furthermore, since the groups $G(\gamma)$ that compose $\mathscr{G}(I)$ are such that $G(\gamma) \cap G(\gamma') = \varnothing$ for $\gamma, \gamma' \in \Gamma(I)$ such that $\gamma \neq \gamma'$, the same occurs, taking Lemma 5.9 into account, with the groups $\widetilde{G}(\gamma)$ of $\widetilde{\mathscr{G}}(I)$ and, hence, not only the $\mathscr{G}(I)$-individuals have well-defined the concept of domain (see Definition 5.10 (a)), but also the $\widetilde{\mathscr{G}}(I)$-individuals do have it.

Keeping now in mind the definitions, regarding a $S$-space, of the concepts of domain of individuals, extension, strict extension, and closed extension, one easily verifies that the interpretation given ahead for the primitive terms of the axiomatic formulated in 5.12, transforms each one of its axioms into a true statement about the objects that compose the extension $\widetilde{\mathscr{G}}(I)$ described in the 1st TESS.

| **Primitive terms** | **Interpretation** |
|---|---|
| $\widetilde{\mathscr{G}}(I)$-distribution | $\widetilde{\mathscr{G}}(I)$-individual |
| $\widetilde{\mathscr{G}}(I)$-distribution domain | $\widetilde{\mathscr{G}}(I)$-individual domain |
| $\widetilde{\mathscr{G}}(I)$-distributions addition | $\widetilde{\mathscr{G}}(I)$-individuals addition |
| $\widetilde{\mathscr{G}}(I)$-distribution derivative | endomorphism $\widetilde{\Phi}$ of the semigroups $\widetilde{H}(\gamma)$ of the $S$-space $\widetilde{\mathscr{G}}(I)$ |
| $\widetilde{\mathscr{G}}(I)$-distribution restriction | homomorphism $\widetilde{\Theta}_{(\gamma',\gamma)}$ of the family $\widetilde{\Theta}(\Delta(I))$ of the $S$-space $\widetilde{\mathscr{G}}(I)$ |

In other terms, the interpretation above composes a model of the $\widetilde{\mathscr{G}}(I)$-distributions axiomatic, with which the consistency of its axioms is established.

Furthermore, as can be easily verified, an interpretation such as the one above, formulated from the particular $S$-space $\widetilde{\mathscr{G}}(I)$ described in the 1st TESS, if defined with another $S$-space which is, such as $\widetilde{\mathscr{G}}(I)$, a strict and closed extension of $\mathscr{G}(I)$, it would also satisfy the axioms and, therefore, would be another model of our axiomatic that, due to the 1st TESS, is isomorphic to the first one obtained from the $S$-space $\widetilde{\mathscr{G}}(I)$.

Now we inquire: Are all the models of this axiomatic obtained from the "interpretation schema" given in the table above, where $\widetilde{\mathscr{G}}(I)$ must be interpreted as an arbitrary $S$-space, unless the requirement of being a strict and closed extension of $\mathscr{G}(I)$? Or, in other terms: The axiomatic in question is categoric?

In what follows, in item 5.14, we derive a set of logical consequences from the axioms under analysis, which will allow us, in item 5.15, to conclude that we are dealing with a categoric axiomatic.

## 5.14 Elementary Consequences of the Axioms

Among the logical consequences of the axioms formulated in 5.12 we highlight the following:

**(a)** $\widehat{g}+\widehat{h}=\widehat{h}+\widehat{g}$ for every $\widehat{g},\widehat{h}\in\widehat{G}(\gamma)$.

In fact, if $\widehat{g}$ and $\widehat{h}$ belongs to $\widehat{G}(\gamma)$, then, by Axiom 6, there exist $g,h\in G(\gamma)$ and $\Phi,\Psi\in H(\gamma)$ such that

$$\widehat{g}=\widehat{\gamma}(\Phi)(g)=\widehat{\Phi}(g)$$

and

$$\widehat{h}=\widehat{\gamma}(\Psi)(h)=\widehat{\Psi}(h).$$

But, $\Phi,\Psi\in H(\gamma)$ are surjective homomorphisms and, hence, there exist $\overset{*}{g}\in G(\gamma)_{\Psi}$ and $\overset{*}{h}\in G(\gamma)_{\Phi}$ such that

$$g=\Psi(\overset{*}{g})$$

and

$$h=\Phi(\overset{*}{h})$$

or yet, taking the Axioms (4-2) and (4-3) into account,

$$g=\Psi(\overset{*}{g})=\widehat{\gamma}(\Psi)(\overset{*}{g})=\widehat{\Psi}(\overset{*}{g})$$

and

$$h=\Phi(\overset{*}{h})=\widehat{\gamma}(\Phi)(\overset{*}{h})=\widehat{\Phi}(\overset{*}{h}).$$

Thus, we have

$$\widehat{g}=\widehat{\Phi}(g)=\widehat{\Phi}\Big(\widehat{\Psi}(\overset{*}{g})\Big)$$

and

$$\widehat{h}=\widehat{\Psi}(h)=\widehat{\Psi}\Big(\widehat{\Phi}(\overset{*}{h})\Big).$$

On the other hand, since $\widehat{\gamma}:H(\gamma)\longrightarrow\widehat{H}(\gamma)$ is an isomorphism (Axiom (4-3)), we have:

$$\widehat{\Phi}\widehat{\Psi}=\widehat{\gamma}(\Phi)\widehat{\gamma}(\Psi)=\widehat{\gamma}(\Phi\Psi)=\widehat{\gamma}(\Psi\Phi)=\widehat{\gamma}(\Psi)\widehat{\gamma}(\Phi)=\widehat{\Psi}\widehat{\Phi}.$$

From this, we can then write

$$\widehat{g}+\widehat{h}=\widehat{\Phi}\Big(\widehat{\Psi}(\overset{*}{g})\Big)+\widehat{\Psi}\Big(\widehat{\Phi}(\overset{*}{h})\Big)=(\widehat{\Phi}\widehat{\Psi})(\overset{*}{g})+(\widehat{\Phi}\widehat{\Psi})(\overset{*}{h}).$$

Finally, by Axiom (4-1), we have:

$$\begin{aligned}\widehat{g}+\widehat{h}&=(\widehat{\Phi}\widehat{\Psi})(\overset{*}{g}+\overset{*}{h})=(\widehat{\Phi}\widehat{\Psi})(\overset{*}{h}+\overset{*}{g})=\\&=(\widehat{\Phi}\widehat{\Psi})(\overset{*}{h})+(\widehat{\Phi}\widehat{\Psi})(\overset{*}{g})=\\&=\widehat{\Psi}\Big(\widehat{\Phi}(\overset{*}{h})\Big)+\widehat{\Phi}\Big(\widehat{\Psi}(\overset{*}{g})\Big)=\\&=\widehat{h}+\widehat{g}.\end{aligned}$$

**(b)** $\widehat{f} + (\widehat{g} + \widehat{h}) = (\widehat{f} + \widehat{g}) + \widehat{h}$ for every $\widehat{f}, \widehat{g}, \widehat{h} \in \widehat{G}(\gamma)$.

With a procedure analogous to the one presented above for (a), one can prove (b).

**(c)** There exists $\widehat{0} \in \widehat{G}(\gamma)$ such that, for every $\widehat{g} \in \widehat{G}(\gamma)$, $\widehat{g} + \widehat{0} = \widehat{g}$.

Let $\widehat{g} \in \widehat{G}(\gamma)$ be arbitrarily chosen. By Axiom 6 we have

$$\widehat{g} = \widehat{\gamma}(\Phi)(g) = \widehat{\Phi}(g)$$

for some $\Phi \in H(\gamma)$ and some $g \in G(\gamma)$. Since $\Phi : G(\gamma)_\Phi \subseteq G(\gamma) \longrightarrow G(\gamma)$ is a homomorphism, we have

$$\Phi(0) = 0$$

being $0 \in G(\gamma)_\Phi \subseteq G(\gamma)$ the additive neutral of the group $G(\gamma)$. Now, by Axioms (4-2) and (4-3), we have

$$\Phi(0) = \widehat{\gamma}(\Phi)(0) = \widehat{\Phi}(0),$$

and, hence,

$$\widehat{\Phi}(0) = 0.$$

From this, and taking Axiom (4-1) and Axiom 3 into account, comes

$$\widehat{g} + 0 = \widehat{\Phi}(g) + \widehat{\Phi}(0) = \widehat{\Phi}(g + 0) = \widehat{\Phi}(g) = \widehat{g},$$

which proves (c) once that $0 \in \widehat{G}(\gamma)$ (since $0 \in G(\gamma) \subseteq \widehat{G}(\gamma)$).

**(d)** For each $\widehat{g} \in \widehat{G}(\gamma)$, there exists $\widehat{h} \in \widehat{G}(\gamma)$ such that

$$\widehat{g} + \widehat{h} = \widehat{0} = 0,$$

with $0 \in G(\gamma)$ as in (c), that is, the additive neutral of the group $G(\gamma)$.

In fact, let $\widehat{g} \in \widehat{G}(\gamma)$ be arbitrarily fixed and let also $g \in G(\gamma)$ and $\Phi \in H(\gamma)$, as Axiom 6 establishes, such that

$$\widehat{g} = \widehat{\gamma}(\Phi)(g) = \widehat{\Phi}(g).$$

Now, we define

$$\widehat{h} := \widehat{\Phi}(-g),$$

and, hence, taking Axiom (4-1) and Axiom 3 into account, we obtain that

$$\widehat{g} + \widehat{h} = \widehat{\Phi}(g) + \widehat{\Phi}(-g) = \widehat{\Phi}\left(g + (-g)\right) = \widehat{\Phi}(0).$$

But, as we saw in the proof of (c), $\widehat{\Phi}(0) = 0$. Therefore,

$$\widehat{g} + \widehat{h} = 0.$$

**(e)** The class $\widehat{G}(\gamma)$, of all the $\widetilde{\mathscr{G}}(I)$-distributions with domain $\gamma$, with the addition operation is an abelian group that has $G(\gamma)$ as its subgroup.

In fact, from (a) to (d) above one concludes that $\widehat{G}(\gamma)$ with the addition operation is an abelian group, while from Axiom 1 and Axiom 2 one has that the group $G(\gamma) \subseteq \widehat{G}(\gamma)$ is a subgroup of the group $\widehat{G}(\gamma)$.

**(f)** The derivatives $\widehat{\Phi} : \widehat{G}(\gamma) \longrightarrow \widehat{G}(\gamma)$ of the $\widetilde{\mathscr{G}}(I)$-distributions with domain $\gamma$ are endomorphisms on the group $\widehat{G}(\gamma)$; the class $\widehat{H}(\gamma)$ of these derivatives, with the usual composition of functions, is a semigroup (of endomorphisms on $\widehat{G}(\gamma)$) that is a prolongation of the semigroup $H(\gamma)$ to $\widehat{G}(\gamma)$ (according to Definition 1.6(c)).

In fact, from (e) and Axiom (4-1) one obtains that the derivatives $\widehat{\Phi} : \widehat{G}(\gamma) \longrightarrow \widehat{G}(\gamma)$ are endomorphisms on the group $\widehat{G}(\gamma)$, while from the Axioms (4-2) and (4-3) one deducts that $\widehat{H}(\gamma)$, with the operation of composition of functions, is a semigroup and that this semigroup is a prolongation of the semigroup $H(\gamma)$ to $\widehat{G}(\gamma)$.

**(g)** The ordered pair $\widehat{\mathbb{G}}(\gamma)$, $\gamma \in \Gamma(I)$, defined by

$$\widehat{\mathbb{G}}(\gamma) := \left(\widehat{G}(\gamma), \widehat{H}(\gamma)\right),$$

is a $S$-group which is an extension of the $S$-group $\mathbb{G}(\gamma) = (G(\gamma), H(\gamma))$.

This trivially follows from (e) and (f) above and the definitions of $S$-group and extension of $S$-group.

**(h)** $\widehat{\mathbb{G}}(\gamma) = (\widehat{G}(\gamma), \widehat{H}(\gamma))$ is a strict and closed extension of the $S$-group $\mathbb{G}(\gamma) = (G(\gamma), H(\gamma))$.

As one deducts from (g), Axiom 5 and Axiom 6.

**(i)** The family $\widehat{\mathbb{G}}(\Gamma(I))$ defined by

$$\widehat{\mathbb{G}}\left(\Gamma(I)\right) := \left\{\widehat{\mathbb{G}}(\gamma) = \left(\widehat{G}(\gamma), \widehat{H}(\gamma)\right)\right\}_{\gamma \in \Gamma(I)},$$

is a strict and closed extension of the family of $S$-groups,

$$\mathbb{G}\left(\Gamma(I)\right) := \left\{\mathbb{G}(\gamma) = \left(G(\gamma), H(\gamma)\right)\right\}_{\gamma \in \Gamma(I)},$$

of the $S$-space

$$\mathscr{G}(I) = \left(\mathbb{G}\left(\Gamma(I)\right), i\left(\Gamma^2(I)\right), \Theta\left(\Delta(I)\right)\right).$$

Immediate consequence of (h) and of the definition of strict and closed extension of a family of $S$-groups (Definition 3.8).

**(j)** The family $\widehat{i}(\Gamma^2(I))$ defined by

$$\widehat{i}\Big(\Gamma^2(I)\Big) \coloneqq \Big\{\widehat{i}_{(\gamma',\gamma)}\Big\}_{(\gamma',\gamma)\in\Gamma^2(I)},$$

where $\widehat{i}_{(\gamma',\gamma)}$ are the functions defined right after Axiom 6, is the extension of the bonding

$$i\Big(\Gamma^2(I)\Big) = \Big\{i_{(\gamma',\gamma)}\Big\}_{(\gamma',\gamma)\in\Gamma^2(I)}$$

(of the $S$-space $\mathscr{G}(I) = (\mathbb{G}(\Gamma(I)), i(\Gamma^2(I)), \Theta(\Delta(I)))$ to the family of $S$-groups $\widehat{\mathbb{G}}(\Gamma(I))$.

This immediately follows from (i) and the definition of extension of a bonding (Definition 3.22).

**(k)** The ordered pair

$$\left(\widehat{\mathbb{G}}\Big(\Gamma(I)\Big), \widehat{i}\Big(\Gamma^2(I)\Big)\right)$$

is a bonded family that is an extension of the bonded family

$$\left(\mathbb{G}\Big(\Gamma(I)\Big), i\Big(\Gamma^2(I)\Big)\right)$$

of the $S$-space

$$\mathscr{G}(I) = \left(\mathbb{G}\Big(\Gamma(I)\Big), i\Big(\Gamma^2(I)\Big), \Theta\Big(\Delta(I)\Big)\right).$$

It follows trivially from (j) and the definitions of bonded family and extension of bonded family (Definitions 3.11 and 3.22, respectively).

**(l)** The family of the restrictions of the $\widetilde{\mathscr{G}}(I)$-distributions, that is, the family $\widehat{\Theta}(\Delta(I))$ defined by

$$\widehat{\Theta}\Big(\Delta(I)\Big) \coloneqq \Big\{\widehat{\Theta}_{(\gamma',\gamma)}\Big\}_{(\gamma',\gamma)\in\Delta(I)}$$

is a restriction for the bonded family

$$\left(\widehat{\mathbb{G}}\Big(\Gamma(I)\Big), \widehat{i}\Big(\Gamma^2(I)\Big)\right).$$

In fact, reporting ourselves to Definition 3.13, which introduces the concept of restriction, we observe that all the conditions, (a) to (d) in 3.13, required by it are fulfilled: the requirement (d) is, in this case where $\widehat{\Phi} \in \widehat{H}(\gamma)$ is an endomorphism, trivially attended, while (a),(b) and (c) are fulfilled due to Axioms (7-1),(7-3) and (7-4), respectively.

**(m)** The triplet $\widehat{\mathscr{G}}(I)$ defined by

$$\widehat{\mathscr{G}}(I) := \Big(\widehat{\mathbb{G}}\Big(\Gamma(I)\Big), \widehat{\imath}\Big(\Gamma^2(I)\Big), \widehat{\Theta}\Big(\Delta(I)\Big)\Big)$$

is a $S$-space.

Immediate consequence of (l) and Definition 3.15 of $S$-space.

**(n)** The family $\widehat{\Theta}(\Delta(I))$ is a prolongation of $\Theta(\Delta(I))$ (the restriction of the $S$-space $\mathscr{G}(I) = (\mathbb{G}(\Gamma(I)), \imath(\Gamma^2(I)), \Theta(\Delta(I))))$ to the bonded family $(\widehat{\mathbb{G}}(\Gamma(I)), \widehat{\imath}(\Gamma^2(I)))$.

Consequence directly obtained from (l), Axiom (7-2) and the definition of prolongation (Definition 3.23(a)).

**(o)** The $S$-space

$$\widehat{\mathscr{G}}(I) = \Big(\widehat{\mathbb{G}}\Big(\Gamma(I)\Big), \widehat{\imath}\Big(\Gamma^2(I)\Big), \widehat{\Theta}\Big(\Delta(I)\Big)\Big)$$

is a strict and closed extension of the $S$-space

$$\mathscr{G}(I) = \Big(\mathbb{G}\Big(\Gamma(I)\Big), \imath\Big(\Gamma^2(I)\Big), \Theta\Big(\Delta(I)\Big)\Big).$$

In fact, we have:

- $\widehat{\mathscr{G}}(I)$ is a $S$-space (by (m)),
- the bonded family $(\widehat{\mathbb{G}}(\Gamma(I)), \widehat{\imath}(\Gamma^2(I)))$ is an extension of the bonded family $(\mathbb{G}(\Gamma(I)), \imath(\Gamma^2(I)))$ (by (k)),
- $\widehat{\Theta}(\Delta(I))$ is a prolongation of the restriction $\Theta(\Delta(I))$ to the bonded family $(\widehat{\mathbb{G}}(\Gamma(I)), \widehat{\imath}(\Gamma^2(I)))$ (by (k)) and
- $\widehat{\mathbb{G}}(\Gamma(I))$ is a strict and closed extension of the family $\mathbb{G}(\Gamma(I))$ (by (i)).

Hence, and taking into account the definitions of extension of a $S$-space (Definition 3.23(b)) and of strict and closed extension of a $S$-space (Definition 3.23(c-1)), we can conclude that $\widehat{\mathscr{G}}(I)$ is a strict and closed extension of $\mathscr{G}(I)$.

## 5.15 The Categoricity of the Axioms

We already know that the $\widetilde{\mathscr{G}}(I)$-distributions axioms are consistent: as we saw in item 5.13, the objects that compose the strict and closed extension,

$$\widetilde{\mathscr{G}}(I) = \left(\widetilde{\mathbb{G}}\Big(\Gamma(I)\Big) = \Big\{\widetilde{\mathbb{G}}(\gamma) = \Big(\widetilde{G}(\gamma), \widetilde{H}(\gamma)\Big)\Big\}_{\gamma\in\Gamma(I)}, \widetilde{i}\Big(\Gamma^2(I)\Big), \widetilde{\Theta}\Big(\Delta(I)\Big)\right),$$

of the $S$-space $\mathscr{G}(I)$ (relative to which the axioms refer to), defined in the 1st TESS, that is, the $\widetilde{\mathscr{G}}(I)$-individuals (elements of the groups $\widetilde{G}(\gamma)$), the endomorphisms of the semigroups $\widetilde{H}(\gamma)$ and the homomorphisms of the family $\widetilde{\Theta}(\Delta(I))$, provide a model of the referred axiomatic.

On the other hand, taking into account the logical consequences above obtained (in 5.14), of which the one described in (o) is a synthesis, one concludes that: the constituent elements of a (any) model of the $\widetilde{\mathscr{G}}(I)$-distributions axiomatic, define a $S$-space, $\widehat{\mathscr{G}}(I)$ (in (o)), which is a strict and closed extension of $\mathscr{G}(I)$. Now, by the 1st TESS, $\widehat{\mathscr{G}}(I)$ is isomorphic to $\widetilde{\mathscr{G}}(I)$, that is, the axiomatic under analysis is categoric.

## 5.16 Remark

As we saw in item 5.14, to each model of the $\widetilde{\mathscr{G}}(I)$-distributions axiomatic corresponds a $S$-space,

$$\widehat{\mathscr{G}}(I) = \left(\widehat{\mathbb{G}}\Big(\Gamma(I)\Big) = \Big\{\widehat{\mathbb{G}}(\gamma) = \Big(\widehat{G}(\gamma), \widehat{H}(\gamma)\Big)\Big\}_{\gamma\in\Gamma(I)}, \widehat{i}\Big(\Gamma^2(I)\Big), \widehat{\Theta}\Big(\Delta(I)\Big)\right),$$

described in 5.14(o), "composed" by the elements of the model, which is a strict and closed extension of the underlying $S$-space

$$\mathscr{G}(I) = \left(\mathbb{G}\Big(\Gamma(I)\Big) = \Big\{\mathbb{G}(\gamma) = \Big(G(\gamma), H(\gamma)\Big)\Big\}_{\gamma\in\Gamma(I)}, i\Big(\Gamma^2(I)\Big), \Theta\Big(\Delta(I)\Big)\right)$$

of the referred axiomatic (that $S$-space to which the axioms refer to). Thus, for each $\gamma \in \Gamma(I)$, the group $\widehat{G}(\gamma)$ and the endomorphisms $\widehat{\Phi} : \widehat{G}(\gamma) \longrightarrow \widehat{G}(\gamma)$ of the set $\widehat{H}(\gamma)$ are such that:

**(a)** $G(\gamma) \subseteq \widehat{G}(\gamma)$;

**(b)** $\widehat{\Phi}(g) = \Phi(g)$ for every $g \in G(\gamma)_\Phi$ and the set

$$\widehat{H}(\gamma) = \Big\{\widehat{\Phi} : \Phi \in H(\gamma)\Big\}$$

with the usual operation of composition of functions, is a semigroup isomorphic to the semigroup $H(\gamma)$, with the function

$$\begin{aligned} \widehat{\gamma} : H(\gamma) &\longrightarrow \widehat{H}(\gamma) \\ \Phi &\longmapsto \widehat{\gamma}(\Phi) := \widehat{\Phi} \end{aligned}$$

as an isomorphism;

**(c)** for each $\widehat{g} \in \widehat{G}(\gamma)$, there exist $\Phi \in H(\gamma)$ and $g \in G(\gamma)$ such that

$$\widehat{g} = \widehat{\Phi}(g).$$

Furthermore, taking into account Proposition 1.16(b) (which applies here once $\widehat{\mathbb{G}} = (\widehat{G}(\gamma), \widehat{H}(\gamma))$ is a strict and closed extension of $\mathbb{G} = (G(\gamma), H(\gamma))$), we also have:

**(d)** if $g, h \in G(\gamma)$ and $\Phi \in H(\gamma)$, then,

$$\widehat{\Phi}(g) = \widehat{\Phi}(h) \quad \text{if and only if} \quad g - h \in N(\Phi).$$

On the other hand, according to Proposition 2.24, the properties (a) to (d) above determine the existence of a single extension of the addition of the group $G(\gamma)$ to the set $\widehat{G}(\gamma) \supseteq G(\gamma)$, in such a way that with this operation the set $\widehat{G}(\gamma)$ becomes an abelian group and the elements $\widehat{\Phi} \in \widehat{H}(\gamma)$ endomorphisms on this group $\widehat{G}(\gamma)$. Well, since $\widehat{G}(\gamma)$ and $\widehat{\Phi} \in \widehat{H}(\gamma)$, as members of the $S$-space $\widehat{\mathscr{G}}(I)$, are, respectively, an abelian group that has $G(\gamma)$ as subgroup and an endomorphism on $\widehat{G}(\gamma)$, the only extension of the addition of the group $G(\gamma)$ to the set $\widehat{G}(\gamma)$ determined by conditions (a) to (d), is the addition of the group $\widehat{G}(\gamma)$ (of the $S$-space $\widehat{\mathscr{G}}(I)$) itself.

Now, remembering ourselves that for the $\widetilde{\mathbb{G}}$-distributions axiomatic formulated in item 2.20, by the same reason that its models also attend the conditions (a) to (d) of Proposition 2.24, it was possible to obtain an equivalent and simpler axiomatic, in which the term "addition" does not figure among the primitive ones anymore, we are naturally lead to conjecture that, also for the $\widetilde{\mathscr{G}}(I)$-distributions, it is possible to formulate an axiomatic equivalent to that in 5.12 containing a smaller number of primitive terms.

In what follows, guided by the trail that conducted us to the simplified version of the $\widetilde{\mathbb{G}}$-distributions axiomatic (formulated in item 2.25), we will prove that the conjecture above is true: in the next item, 5.17, we formulate the axioms and, next, in item 5.18, we prove the equivalence of this new axiomatic with the one in item 5.12.

## 5.17 $\widetilde{\mathscr{G}}(I)$-Distributions Axioms: Simplified Version

Such as for the formulation of $\widetilde{\mathscr{G}}(I)$-distributions axiomatic in item 5.12, here also, the axioms of this simplified version are stated regarding an abelian, surjective, and with

identity $S$-space

$$\mathscr{G}(I) = \Bigg( \mathbb{G}\Big(\Gamma(I)\Big) = \Big\{ \mathbb{G}(\gamma) = \Big( G(\gamma), H(\gamma) \Big) \Big\}_{\gamma \in \Gamma(I)}, \\ i\Big(\Gamma^2(I)\Big) = \Big\{ i_{(\gamma',\gamma)} \Big\}_{(\gamma',\gamma) \in \Gamma^2(I)}, \\ \Theta\Big(\Delta(I)\Big) = \Big\{ \Theta_{(\gamma',\gamma)} \Big\}_{(\gamma',\gamma) \in \Delta(I)} \Bigg),$$

arbitrarily fixed, such that $G(\gamma) \cap G(\gamma') = \varnothing$ for every $\gamma, \gamma' \in \Gamma(I)$ such that $\gamma \neq \gamma'$.

The **primitive terms** that the axioms ahead implicitly define are: $\widetilde{\mathscr{G}}(I)$-distribution, domain, derivative and restriction of $\widetilde{\mathscr{G}}(I)$-distribution.

The **precedent theories**, the same ones assumed for the previous axiomatic (in 5.12): Classical Logic, Set Theory, and Theory of $S$-Spaces.

The **axioms** are the ones presented below.

**Axiom 1** Every $\mathscr{G}(I)$-individual is a $\widetilde{\mathscr{G}}(I)$-distribution.

**Axiom 2** To each $\widetilde{\mathscr{G}}(I)$-distribution, $\widehat{g}$, corresponds a single $\gamma \in \Gamma(I)$, denominated the domain of $\widehat{g}$, in such a way that if $\widehat{g}$ is a $\mathscr{G}(I)$-individual, the domain $\gamma$ is the domain of the $\mathscr{G}(I)$-individual $\widehat{g}$.

**Notation.** We will denote by $\widehat{G}(\gamma)$ the class of all the $\widetilde{\mathscr{G}}(I)$-distributions with domain $\gamma \in \Gamma(I)$. It results from Axiom 1 and Axiom 2 that the class of $\mathscr{G}(I)$-individuals with domain $\gamma$, $G(\gamma)$, is such that $G(\gamma) \subseteq \widehat{G}(\gamma)$.

**Axiom 3** The derivatives of the $\widetilde{\mathscr{G}}(I)$-distributions with domain $\gamma$, whose class is denoted by $\widehat{H}(\gamma)$, are functions with domain and codomain equal to $\widehat{G}(\gamma)$, such that:

**(3-1)** each derivative $\widehat{\Phi} \in \widehat{H}(\gamma)$ is an extension to $\widehat{G}(\gamma)$ of a single $\Phi \in H(\gamma)$, that is, $\widehat{\Phi}(g) = \Phi(g)$ for every $g \in G(\gamma)_\Phi \subseteq \widehat{G}(\gamma)$; conversely, for each $\Phi \in H(\gamma)$ there exists a single derivative $\widehat{\Phi} \in \widehat{H}(\gamma)$ such that $\Phi(g) = \widehat{\Phi}(g)$ for every $g \in G(\gamma)_\Phi \subseteq \widehat{G}(\gamma)$;

**(3-2)** the class $\widehat{H}(\gamma)$, equipped with the usual operation of composition of functions, is a semigroup isomorphic to the semigroup $H(\gamma)$ with the function

$$\begin{aligned} \widehat{\gamma} :\ H(\gamma) &\longrightarrow \widehat{H}(\gamma) \\ \Phi &\longmapsto \widehat{\gamma}(\Phi) := \text{the extension of } \Phi \text{ to } \widehat{G}(\gamma) \end{aligned}$$

as an isomorphism.

**Notation.** For $\Phi \in H(\gamma)$, we will denote $\widehat{\gamma}(\Phi)$ by $\widehat{\Phi}$, that is, $\widehat{\Phi} := \widehat{\gamma}(\Phi)$.

**Axiom 4** Let $\gamma \in \Gamma(I)$ be arbitrarily fixed. For each $\widehat{g} \in \widehat{G}(\gamma)$ there exist $\Phi \in H(\gamma)$ and $g \in G(\gamma)$ such that

$$\widehat{g} = \widehat{\Phi}(g).$$

In other terms, every $\widetilde{\mathscr{G}}(I)$-distribution is a derivative of some $\mathscr{G}(I)$-individual.

**Axiom 5** For $\gamma \in \Gamma(I)$, $g, h \in G(\gamma)$ and $\Phi \in H(\gamma)$ arbitrarily fixed, one has:

$$\widehat{\Phi}(g) = \widehat{\Phi}(h) \quad \text{if and only if} \quad g - h \in N(\Phi).$$

**Definition.** For each $(\gamma', \gamma) \in \Gamma^2(I)$, we define $\widehat{i}_{(\gamma',\gamma)}$ as the following function:

$$\begin{aligned} \widehat{i}_{(\gamma',\gamma)} : \widehat{H}(\gamma) &\longrightarrow \widehat{H}(\gamma') \\ \widehat{\Phi} &\longmapsto \widehat{i}_{(\gamma',\gamma)}(\widehat{\Phi}) := \widehat{\gamma'}\left( i_{(\gamma',\gamma)}\left( (\widehat{\gamma})^{-1}(\widehat{\Phi}) \right) \right) \end{aligned}$$

It results that $\widehat{i}_{(\gamma',\gamma)}$ is an isomorphism from the semigroup $\widehat{H}(\gamma)$ onto the semigroup $\widehat{H}(\gamma')$.

**Axiom 6** The restrictions of $\widetilde{\mathscr{G}}(I)$-distributions are functions $\widehat{\Theta}_{(\gamma',\gamma)} : \widehat{G}(\gamma) \longrightarrow \widehat{G}(\gamma')$, a single one for each $(\gamma', \gamma) \in \Delta(I)$, that associates to each $\widetilde{\mathscr{G}}(I)$-distribution with domain $\gamma$, a single $\widetilde{\mathscr{G}}(I)$-distribution with domain $\gamma' \subseteq \gamma$, in such a way that:

**(6-1)** $\widehat{\Theta}_{(\gamma',\gamma)}(g) = \Theta_{(\gamma',\gamma)}(g)$ for every $g \in G(\gamma) \subseteq \widehat{G}(\gamma)$;

**(6-2)** $\widehat{\Theta}_{(\gamma'',\gamma')}\left(\widehat{\Theta}_{(\gamma',\gamma)}(\widehat{g})\right) = \widehat{\Theta}_{(\gamma'',\gamma)}(\widehat{g})$ for every $\widehat{g} \in \widehat{G}(\gamma)$ and every $\gamma, \gamma', \gamma'' \in \Gamma(I)$ such that $\gamma'' \subseteq \gamma' \subseteq \gamma$;

**(6-3)** $\widehat{\Theta}_{(\gamma',\gamma)}\left(\widehat{\Phi}(\widehat{g})\right) = \widehat{i}_{(\gamma',\gamma)}(\widehat{\Phi})\left(\widehat{\Theta}_{(\gamma',\gamma)}(\widehat{g})\right)$ for every $\widehat{g} \in \widehat{G}(\gamma)$ and every $\widehat{\Phi} \in \widehat{H}(\gamma)$.

## 5.18 Equivalence of the Axiomatics Formulated in Items 5.12 and 5.17

As we know, two axiomatics are equivalent if and only if the primitive terms of each one of them are primitive or defined terms of the other, and the axioms of each one of them are axioms or theorems of the other. In order for us to prove that this occurs regarding the axioms formulated in items 5.12 and 5.17, we will consider again the Proposition 2.24 and proceed in an analogous manner to that in item 2.26, where we proved the equivalence of the two axiomatics formulated for the $\widetilde{\mathbb{G}}$-distributions. Hence, in what follows, with the adaptations and complements, required for this case, to the procedure in 2.26, the equivalence between the two $\widetilde{\mathscr{G}}(I)$-distributions axiomatics is proved.

Let $\gamma \in \Gamma(I)$ be arbitrarily fixed and $\mathbb{G}(\gamma) = (G(\gamma), H(\gamma))$ be the corresponding $S$-group of the underlying $S$-space $\mathscr{G}(I)$ of the axiomatics in question (the $S$-space to which the axioms refer to). Since $\mathscr{G}(I)$ is abelian, surjective, and with identity, its $S$-groups also are abelian, surjective, and with identity. Hence, we can take the $S$-group $\mathbb{G} = (G, H)$ of Proposition 2.24 as $\mathbb{G} = (G(\gamma), H(\gamma))$, which allows us to conclude that the Axioms 1 and 2, 3, 4, and 5 of the axiomatic in 5.17 are, exact and respectively, the conditions (a) to (d) of Proposition 2.24 with $\widehat{G} = \widehat{G}(\gamma)$ and $\widehat{H} = \widehat{H}(\gamma)$.

Taking now into account that the proof of the referred proposition, composed of the proof of Lemmas 1 to 6 there contained, leans, exclusively, on the hypotheses (a) to (d) and of the $S$-group $\mathbb{G} = (G, H)$ being abelian, surjective, and with identity, we concluded that these lemmas translated to the case where $\mathbb{G} = \mathbb{G}(\gamma)$ (the $S$-group of the $S$-space $\mathscr{G}(I)$ corresponding to $\gamma \in \Gamma(I)$), i.e., $\widehat{G} = \widehat{G}(\gamma)$ (the class of the $\widetilde{\mathscr{G}}(I)$-distributions with domain $\gamma$) and $\widehat{H} = \widehat{H}(\gamma)$ (the set of derivatives of $\widetilde{\mathscr{G}}(I)$-distributions with domain $\gamma$), they are logical consequences of Axioms 1, 2, 3, 4 and 5 above referred and, therefore, are theorems of the axiomatic formulated in 5.17. We highlight, for convenience, the following composition of logical consequences of these axioms, established in Lemmas 3 to 6 (contained in the proof of Proposition 2.24):

- For each $\gamma \in \Gamma(I)$, let $+ : \widehat{G}(\gamma) \times \widehat{G}(\gamma) \longrightarrow \widehat{G}(\gamma)$ be the binary operation in $\widehat{G}(\gamma)$, denominated addition, defined as follows: for $\widehat{g}, \widehat{h} \in \widehat{G}(\gamma)$ arbitrarily chosen, let $\Phi, \Psi \in H(\gamma)$ and $g, h \in G(\gamma)$ such that

  $$\widehat{g} = \widehat{\Phi}(g) \quad \text{and} \quad \widehat{h} = \widehat{\Psi}(h)$$

  and, hence, we define

  $$\widehat{g} + \widehat{h} = \widehat{\Phi}(g) + \widehat{\Psi}(h) := \widehat{\Phi\Psi}(g^* + h^*)$$

  with

  $$g^* \in \Psi^{-1}(g) \quad \text{and} \quad h^* \in \Phi^{-1}(h)$$

  arbitrarily chosen.

  Regarding this operation, one has that:

  **(i)** The $+$ operation (addition) above is well-defined;

  **(ii)** $\widehat{G}(\gamma)$ with the addition above is an abelian group;

  **(iii)** The group $G(\gamma)$ is a subgroup of the group $\widehat{G}(\gamma)$ (given in (ii));

**(iv)** For each $\Phi \in H(\gamma)$, $\widehat{\Phi} \in \widehat{H}(\gamma)$ is an endomorphism on the group $\widehat{G}(\gamma)$ (of (ii));

**(v)** The $+$ operation above defined is the only binary operation in $\widehat{G}(\gamma)$ with the properties (ii) to (iv).

The consequence (i) above allows to "enrich" the axiomatic described in 5.17, henceforth axiomatic 2, with a new term, "addition", defining it as the name of the binary operation, $+$, above described. It immediately results, taking into account the consequences (ii) to (v), that among the theorems of axiomatic 2 enriched with the addition as a defined concept, are Axiom 3 and Axiom (4-1) of the axiomatic in item 5.12, henceforth referred to as axiomatic 1.

And what about Axiom 5 and Axiom (7-1) of axiomatic 1, which does not figure as axioms in the axiomatic 2 (the only ones if we disconsider Axiom 3 and Axiom (4-1) that, as we saw, are theorems of axiomatic 2), would they also be theorems of axiomatic 2? This is what in fact occurs as we will prove below.

**Remark.** At the proofs that the Axiom 5 and Axiom (7-1) are theorems of axiomatic 2, presented below in (vii) and (vi), respectively, any reference to an axiom without mention to what axiomatic (1 or 2) it belongs to, means that the referred axiom belongs to axiomatic 2.

**(vi)** Axiom (7-1) of axiomatic 1:

$$\widehat{\Theta}_{(\gamma',\gamma)}(\widehat{g}+\widehat{h}) = \widehat{\Theta}_{(\gamma',\gamma)}(\widehat{g}) + \widehat{\Theta}_{(\gamma',\gamma)}(\widehat{h})$$

for every $(\gamma',\gamma) \in \Delta(I)$ and any $\widehat{g},\widehat{h} \in \widehat{G}(\gamma)$.

In fact, for $\widehat{g},\widehat{h} \in \widehat{G}(\gamma)$ arbitrarily fixed, Axiom 3 and Axiom 4 allow us to say that there exist $\Phi, \Psi \in H(\gamma)$ and $g, h \in G(\gamma)$ such that

$$\widehat{g} = \widehat{\gamma}(\Phi)(g) = \widehat{\Phi}(g) \quad \text{and} \quad \widehat{h} = \widehat{\gamma}(\Psi)(h) = \widehat{\Psi}(h).$$

Since $\Phi : G(\gamma)_\Phi \longrightarrow G(\gamma)$ and $\Psi : G(\gamma)_\Psi \longrightarrow G(\gamma)$ are surjective homomorphisms, there exist $g_1 \in G(\gamma)_\Psi$ and $h_1 \in G(\gamma)_\Phi$ such that

$$g = \Psi(g_1) \quad \text{and} \quad h = \Phi(h_1)$$

and, hence, taking Axiom 3 into account, we get

$$g = \Psi(g_1) = \widehat{\Psi}(g_1) \quad \text{and} \quad h = \Phi(h_1) = \widehat{\Phi}(h_1),$$

which lead us to

$$\widehat{g} = \widehat{\Phi}(g) = \widehat{\Phi}\Big(\widehat{\Psi}(g_1)\Big) \quad \text{and} \quad \widehat{h} = \widehat{\Psi}(h) = \widehat{\Psi}\Big(\widehat{\Phi}(h_1)\Big).$$

Resorting again to Axiom 3, as well as the fact of $H(\gamma)$ being an abelian semigroup, we get that:

$$\widehat{\Phi}\widehat{\Psi} = \widehat{\gamma}(\Phi)\widehat{\gamma}(\Psi) = \widehat{\gamma}(\Phi\Psi) = \widehat{\gamma}(\Psi\Phi) = \widehat{\gamma}(\Psi)\widehat{\gamma}(\Phi) = \widehat{\Psi}\widehat{\Phi}.$$

With the results above, we obtain:

$$\begin{aligned}\widehat{g} + \widehat{h} &= \widehat{\Phi}\Big(\widehat{\Psi}(g_1)\Big) + \widehat{\Psi}\Big(\widehat{\Phi}(h_1)\Big) = \\ &= (\widehat{\Phi}\widehat{\Psi})(g_1) + (\widehat{\Phi}\widehat{\Psi})(h_1) = \\ &= \widehat{\Phi\Psi}(g_1) + \widehat{\Phi\Psi}(h_1)\end{aligned}$$

where $\widehat{\Phi\Psi} = \widehat{\gamma}(\Phi\Psi) \in \widehat{H}(\gamma)$.

But, by consequence (iv) (page 246), $\widehat{\Phi\Psi} \in \widehat{H}(\gamma)$ is an endomorphism on the group $\widehat{G}(\gamma)$ and, hence,

$$\widehat{g} + \widehat{h} = \widehat{\Phi\Psi}(g_1 + h_1).$$

Consequently,

$$\widehat{\Theta}_{(\gamma',\gamma)}(\widehat{g} + \widehat{h}) = \widehat{\Theta}_{(\gamma',\gamma)}\Big(\widehat{\Phi\Psi}(g_1 + h_1)\Big)$$

from where one obtains, applying the Axiom (6-3), that

$$\widehat{\Theta}_{(\gamma',\gamma)}(\widehat{g} + \widehat{h}) = \widehat{i}_{(\gamma',\gamma)}(\widehat{\Phi\Psi})\Big(\widehat{\Theta}_{(\gamma',\gamma)}(g_1 + h_1)\Big).$$

Now, since $g_1 \in G(\gamma)_\Psi \subseteq G(\gamma)$ and $h_1 \in G(\gamma)_\Phi \subseteq G(\gamma)$, it comes that $g_1 + h_1 \in G(\gamma)$ and, hence, by Axiom (6-1),

$$\widehat{\Theta}_{(\gamma',\gamma)}(g_1 + h_1) = \Theta_{(\gamma',\gamma)}(g_1 + h_1),$$

and, since $\Theta_{(\gamma',\gamma)} \in \Theta(\Delta(I))$ is a homomorphism, we obtain

$$\widehat{\Theta}_{(\gamma',\gamma)}(g_1 + h_1) = \Theta_{(\gamma',\gamma)}(g_1) + \Theta_{(\gamma',\gamma)}(h_1) = \widehat{\Theta}_{(\gamma',\gamma)}(g_1) + \widehat{\Theta}_{(\gamma',\gamma)}(h_1).$$

Thus,

$$\widehat{\Theta}_{(\gamma',\gamma)}(\widehat{g} + \widehat{h}) = \widehat{i}_{(\gamma',\gamma)}(\widehat{\Phi\Psi})\Big(\widehat{\Theta}_{(\gamma',\gamma)}(g_1) + \widehat{\Theta}_{(\gamma',\gamma)}(h_1)\Big).$$

Reporting ourselves to the definition of $\widehat{i}_{(\gamma',\gamma)}$ we get that

$$\widehat{i}_{(\gamma',\gamma)}(\widehat{\Phi\Psi}) \in \widehat{H}(\gamma')$$

and, hence, by consequence (iv) (page 246), $\widehat{i}_{(\gamma',\gamma)}(\widehat{\Phi\Psi})$ is an endomorphism on the group $\widehat{G}(\gamma')$, which allows us to write

$$\widehat{\Theta}_{(\gamma',\gamma)}(\widehat{g} + \widehat{h}) = \widehat{i}_{(\gamma',\gamma)}(\widehat{\Phi\Psi})\Big(\widehat{\Theta}_{(\gamma',\gamma)}(g_1)\Big) + \widehat{i}_{(\gamma',\gamma)}(\widehat{\Phi\Psi})\Big(\widehat{\Theta}_{(\gamma',\gamma)}(h_1)\Big).$$

From this, with a new application of Axiom (6-3), it results that

$$\widehat{\Theta}_{(\gamma',\gamma)}(\widehat{g} + \widehat{h}) = \widehat{\Theta}_{(\gamma',\gamma)}\Big(\widehat{\Phi\Psi}(g_1)\Big) + \widehat{\Theta}_{(\gamma',\gamma)}\Big(\widehat{\Phi\Psi}(h_1)\Big).$$

But, as we saw above,

$$\widehat{\Phi\Psi}(g_1) = \widehat{g} \quad \text{and} \quad \widehat{\Phi\Psi}(h_1) = \widehat{h}$$

and, hence

$$\widehat{\Theta}_{(\gamma',\gamma)}(\widehat{g} + \widehat{h}) = \widehat{\Theta}_{(\gamma',\gamma)}(\widehat{g}) + \widehat{\Theta}_{(\gamma',\gamma)}(\widehat{h}).$$

**(vii)** Axiom 5 of axiomatic 1.

For $\gamma \in \Gamma(I)$ and $\Phi \in H(\gamma)$ arbitrarily fixed, if $g \in G(\gamma)$ is such that $g \notin G(\gamma)_\Phi$, then, $\widehat{\Phi}(\gamma) \notin G(\gamma)$.

The proof of (vii) is the same, with the due adaptations, as the one at page 65 upon its establishment, in Chapter 2, of the equivalence between the two axiomatics formulated for the $\widetilde{\mathbb{G}}$-distributions. However, for the reader's convenience, we re-present it below.

It is easy to see that (vii) is equivalent to:

$$\text{if} \quad g \in G(\gamma) \quad \text{and} \quad \widehat{\Phi}(g) \in G(\gamma), \quad \text{then} \quad g \in G(\gamma)_\Phi.$$

Let us then prove this implication and, in order to do so, let $\Phi \in H(\gamma)$ and $g \in G(\gamma)$ be arbitrarily fixed and such that

$$\widehat{\Phi}(g) \in G(\gamma).$$

Once that $\Phi : G(\gamma)_\Phi \longrightarrow G(\gamma)$ is surjective, there exists $g_1 \in G(\gamma)_\Phi$ such that

$$\Phi(g_1) = \widehat{\Phi}(g).$$

But, by Axiom 3,

$$\Phi(g_1) = \widehat{\Phi}(g_1)$$

and, hence, we have

$$\widehat{\Phi}(g_1) = \widehat{\Phi}(g),$$

from where one concludes, resorting to Axiom 5, that

$$g_1 - g \in N(\Phi).$$

Therefore, since $N(\Phi) \subseteq G(\gamma)_\Phi$, $g_1 - g \in G(\gamma)_\Phi$ and, once that $G(\gamma)_\Phi$ is a group (subgroup of $G(\gamma)$) and $g_1 \in G(\gamma)_\Phi$, $g_1 - (g_1 - g) = g \in G(\gamma)_\Phi$.

Finally, for us to complete the proof of the equivalence of axiomatics 1 and 2, we shall proof that Axiom 5 of axiomatic 2, the only one that does not figure among the ones of axiomatic 1, is a theorem of the latter.

We also observe here that the proof provided below is, up to some few notational changes, the same one at page 65 presented when studying the equivalence of the $\widetilde{\mathbb{G}}$-distributions axiomatics.

Let us then proof that the following statement is a theorem of axiomatic 1:

- for $g, h \in G(\gamma)$ and $\Phi \in H(\gamma)$, one has $\widehat{\Phi}(g) = \widehat{\Phi}(h)$ if and only if $g - h \in N(\Phi)$.

**Remark.** In the proof ahead, different from those for (vi) and (vii) above, the reference to an axiom or a consequence of the axioms, without mention to the axiomatic (1 or 2) it belongs to, means that the referred axiom or consequence belongs to axiomatic 1.

Let $g, h \in G(\gamma)$ and $\Phi \in H(\gamma)$ be arbitrarily fixed. If $g - h \in N(\Phi)$, then, $g - h \in G(\gamma)_\Phi$ and, hence, by Axiom 4,

$$\widehat{\Phi}(g - h) = \Phi(g - h) = 0.$$

Now, among the logical consequences of the axioms obtained in item 5.14, we have that $\widehat{G}(\gamma)$ is an abelian group that has $G(\gamma) \subseteq \widehat{G}(\gamma)$ as a subgroup (5.14(e)) and $\widehat{H}(\gamma)$ is a semigroup of endomorphisms on $\widehat{G}(\gamma)$ (5.14(f)). Hence, from $\widehat{\Phi}(g - h) = 0$ results that

$$\widehat{\Phi}(g) = \widehat{\Phi}(h).$$

Conversely, let us suppose now that $\widehat{\Phi}(g) = \widehat{\Phi}(h)$, that is, taking into account the logical consequences above referred (5.14(e) and 5.14(f)), that $g - h \in N(\widehat{\Phi})$. It is also a logical consequence of the axioms that $(\widehat{G}(\gamma), \widehat{H}(\gamma))$ is a strict and closed extension of the $S$-group (abelian, surjective, and with identity) $\mathbb{G}(\gamma) = (G(\gamma), H(\gamma))$ (5.14(h)). Hence being, and now with Proposition 1.16(b) in mind (which also is a theorem of axiomatic 1, once it is derived from hypotheses fulfilled by this axiomatic), by which $N(\widehat{\Phi}) = N(\Phi)$, we conclude that $g - h \in N(\Phi)$.

It is then established the equivalence of the two axiomatics formulated in 5.12 and 5.17.

## 5.19 $\widetilde{\mathbb{G}}$-Distributions Obtained from $\widetilde{\mathscr{G}}(I)$-Distributions

Let us consider now the $\widetilde{\mathscr{G}}(I)$-distributions axiomatic (retake the axioms in item 5.12) particularized to the case where the underlying $S$-space, $\mathscr{G}(I)$ (the one that the axioms refer to), is the $S$-space generated by a $S$-group, $\mathbb{G} = (G, H)$, abelian, surjective, and with identity. In this case, as we saw in 5.5, $\Gamma(I) = \{I\}$, $\Gamma^2(I) = \Delta(I) = \{(I, I)\}$ being I an arbitrarily fixed set, and

$$\mathscr{G}(I) = \mathscr{G}[\mathbb{G}] = \left(\mathbb{G}\Big(\Gamma(I)\Big), i\Big(\Gamma^2(I)\Big), \Theta\Big(\Delta(I)\Big)\right)$$

where

$$\mathbb{G}\Big(\Gamma(I)\Big) = \Big\{\mathbb{G}(\gamma) = \Big(G(\gamma), H(\gamma)\Big) := \mathbb{G} = (G, H)\Big\}_{\gamma \in \{I\}} = \Big\{\mathbb{G}(I) = \mathbb{G}\Big\},$$

$$i\Big(\Gamma^2(I)\Big) = \Big\{i_{(\gamma', \gamma)} := I_H\Big\}_{(\gamma', \gamma) \in \{(I, I)\}} = \Big\{i_{(I,I)} := I_H\Big\}$$

and

$$\Theta\Big(\Delta(I)\Big) = \Big\{\Theta_{(\gamma',\gamma)} \coloneqq I_G\Big\}_{(\gamma',\gamma)\in\{(I,I)\}} = \Big\{\Theta_{(I,I)} \coloneqq I_G\Big\},$$

with $I_H$ and $I_G$ the identity functions on $H$ and $G$, respectively.

Clearly, in this case, since $\Gamma(I) = \{I\}$, the primitive notion of domain of a $(\widetilde{\mathscr{G}}(I) =)$ $\widetilde{\mathscr{G}}[\mathbb{G}]$-distribution, referred to in Axiom 2 (in 5.12), trivializes: every $\widetilde{\mathscr{G}}[\mathbb{G}]$-distribution has the same set $I$ (arbitrarily fixed) as its domain. Furthermore, the reader will have no problem noticing that Axioms 1, 3, 4, 5 and 6 (the first six ones in 5.12, discarding Axiom 2) are, with an obvious adaptation in the notation ($G$ instead of $G(\gamma)$, $H$ instead of $H(\gamma)$, etc.), exactly, the $\widetilde{\mathbb{G}}$-distributions axioms formulated in item 2.20. Regarding Axiom 7 ((7-1) to (7-4)), for which one has no correspondent in the $\widetilde{\mathbb{G}}$-distributions axiomatic (in 2.20), we observe, yet in the case of $\mathscr{G}(I) = \mathscr{G}[\mathbb{G}]$, that:

- the function $\widehat{i}_{(I,I)} : \widehat{H} \longrightarrow \widehat{H}$ that figures in Axiom (7-4), as one easily concludes from its definition (at page 234), is the identity function on $\widehat{H}$:

$$\widehat{i}_{(I,I)} = I_{\widehat{H}};$$

- $\widehat{\Theta}_{(I,I)}\widehat{G} \longrightarrow \widehat{G}$ satisfies Axiom 7 if and only if $\widehat{\Theta}_{(I,I)}$ is the identity function on $\widehat{G}$:

$$\widehat{\Theta}_{(I,I)} = I_{\widehat{G}}.$$

  In fact, on the hypothesis of a function $\widehat{\Theta}_{(I,I)} : \widehat{G} \longrightarrow \widehat{G}$ to attend the requirements of Axiom 7, and being $\widehat{g} \in \widehat{G}$ arbitrarily fixed, we have: by Axiom 6, $\widehat{g} = \widehat{\Phi}(g)$ with $\Phi \in H$ and $g \in G$, hence, taking into account Axioms (7-2) and (7-4),

$$\widehat{\Theta}_{(I,I)}(\widehat{g}) = \widehat{\Theta}_{(I,I)}\Big(\widehat{\Phi}(g)\Big) = \widehat{i}_{(I,I)}(\widehat{\Phi})\Big(\widehat{\Theta}_{(I,I)}(g)\Big) = \widehat{\Phi}\Big(\Theta_{(I,I)}(g)\Big) = \widehat{\Phi}(g) = \widehat{g}.$$

  Conversely, if $\widehat{\Theta}_{(I,I)} = I_{\widehat{G}}$, all the requirements ((7-1) to (7-4)) of Axiom 7 are trivially fullfiled.

Therefore, of the primitive terms that figure in the $\widetilde{\mathscr{G}}[\mathbb{G}]$-distributions axioms, assume a trivial role not only the term "domain", but also the one of "restriction": all $\widetilde{\mathscr{G}}[\mathbb{G}]$-distributions have the same domain, $I$, and the restriction of any $\widetilde{\mathscr{G}}[\mathbb{G}]$-distribution is the own $\widetilde{\mathscr{G}}[\mathbb{G}]$-distribution.

Discarding of the axiomatic in question these two notions, as well as the axioms corresponding to them (Axiom 2 and Axiom 7), one obtains an axiomatic composed by Axioms 1, 3, 4, 5 and 6 that involve and, therefore, implicitly define, only the relevant primitive terms, namely: $\widetilde{\mathscr{G}}[\mathbb{G}]$-distribution, addition and $\widetilde{\mathscr{G}}[\mathbb{G}]$-distribution derivative. But, as above remarked, these axioms are exactly the $\widetilde{\mathbb{G}}$-distribution ones formulated in 2.20 and, hence, taking into account that this axiomatic (in 2.20) is categoric, we conclude that the primitive concepts defined by it, namely, $\widetilde{\mathbb{G}}$-distribution, addition, and $\widetilde{\mathbb{G}}$-distribution derivative, are, respectively, equal to the ones of $\widetilde{\mathscr{G}}[\mathbb{G}]$-distribution, addition and $\widetilde{\mathscr{G}}[\mathbb{G}]$-distribution derivative.

Worth highlighting here that the conclusion above, i.e., that the defining axioms of the $\widetilde{\mathbb{G}}$-distributions are the same axioms that define the $\widetilde{\mathscr{G}}[\mathbb{G}]$-distributions in its non-trivial aspects, is in perfect harmony with the identification (proposed in item 5.5) of a $S$-group with the $S$-space generated by it, that is,

$$\mathbb{G} = (G, H) \equiv \mathscr{G}[\mathbb{G}] = \left(\left\{\mathbb{G} = (G, H)\right\}, \left\{I_H\right\}, \left\{I_G\right\}\right),$$

along with the result (proved in 5.6) about $\widehat{\mathbb{G}} = (\widehat{G}, \widehat{H})$ being a strict and closed extension of an abelian, surjective, and with identity $S$-group $\mathbb{G} = (G, H)$, then, the $S$-space

$$\mathscr{G}[\widehat{\mathbb{G}}] \equiv \widehat{\mathbb{G}} = (\widehat{G}, \widehat{H})$$

is isomorphic to the (unique unless isomorphism) strict and closed extension of $\mathscr{G}[\mathbb{G}] \equiv \mathbb{G}$,

$$\widetilde{\mathscr{G}}[\mathbb{G}] = \mathscr{G}[\widetilde{\mathbb{G}}] \equiv \widetilde{\mathbb{G}},$$

described in the 1st TESS, which, in this case (where $\mathscr{G}(I) = \mathscr{G}[\mathbb{G}] \equiv \mathbb{G}$) is essentially the Theorem of Extension of $S$-Groups.

# $\overline{\mathscr{G}}(I)$-Distributions Axiomatics

## 5.20 Remarks

Let now

$$\mathscr{G}(I) = \Bigg(\mathbb{G}\Big(\Gamma(I)\Big) = \Big\{\mathbb{G}(\gamma) = \Big(G(\gamma), H(\gamma)\Big)\Big\}_{\gamma\in\Gamma(I)},$$
$$i\Big(\Gamma^2(I)\Big) = \Big\{i_{(\gamma',\gamma)}\Big\}_{(\gamma',\gamma)\in\Gamma^2(I)},$$
$$\Theta\Big(\Delta(I)\Big) = \Big\{\Theta_{(\gamma',\gamma)}\Big\}_{(\gamma',\gamma)\in\Delta(I)}\Bigg),$$

be an abelian, surjective, with identity, and coherent $S$-space, arbitrarily fixed, such that

$$G(\gamma) \cap G(\gamma') = \varnothing$$

for every $\gamma, \gamma' \in \Gamma(I)$ such that $\gamma \neq \gamma'$, and let us consider the 2nd species domain of distribution (Definition 5.10(b)) determined by $\mathscr{G}(I)$, i.e., the ordered pair

$$\left(\mathscr{G}(I), \widehat{\mathscr{G}}(I)\right)$$

where $\widehat{\mathscr{G}}(I)$ is a $S$-space isomorphic to the extension $\overline{\mathscr{G}}(I)$ of $\mathscr{G}(I)$ defined in the 2nd TESS (Theorem 4.33). As we know, regarding one such domain of distribution, unless isomorphism, are uniquely characterized, via Definition 5.10(d), the concepts of $\overline{\mathscr{G}}(I)$-distribution, domain, addition, derivative and restriction of $\overline{\mathscr{G}}(I)$-distributions.

Taking into account that $\widetilde{\mathscr{G}}(I)$-distribution and the associated notions of domain, addition, derivative, and restriction of $\widetilde{\mathscr{G}}(I)$-distributions are, in an analogous manner, defined regarding a 1st species domain of distribution and that, as we saw, these concepts can be characterized via categoric axiomatics, it is natural to inquire if the same would not happen for the $\overline{\mathscr{G}}(I)$-distributions.

It is the purpose of this section to prove that, just as for the $\widetilde{\mathscr{G}}(I)$-distributions, so the $\overline{\mathscr{G}}(I)$-distributions (and its associated notions) can be defined through a categoric axiomatic; such as was done for the $\widetilde{\mathscr{G}}(I)$-distributions, in the next items, we will formulate two axiomatics (categoric and equivalents to each other) for the $\overline{\mathscr{G}}(I)$-distributions.

## 5.21 $\overline{\mathscr{G}}(I)$-Distributions Axioms: 1st Version

It is worth remembering at this point what guided us in the choice of the $\widetilde{\mathscr{G}}(I)$-distributions axioms: at the initial considerations of item 5.12 we referred to this question in the following terms, "the choice will be driven by the 1st TESS"…"roughly speaking, assigning the "status" of axiom to the defining predicates of the concepts of extension, strict extension, and closed extension regarding…a $S$-space".

In view of the similarity between the definitions of $\widetilde{\mathscr{G}}(I)$- and $\overline{\mathscr{G}}(I)$-distributions (and its associated notions) given through the notion of domain of distribution (Definitions 5.10(c) and (d)), it is certainly reasonable to follow an analogous direction for the obtainment of axioms that could characterize the $\overline{\mathscr{G}}(I)$-distributions. Thus, now guided by the 2nd TESS (Theorem 4.33), and assigning, roughly speaking, the "status" of axiom to the defining predicates of the concepts referred to in this theorem, we formulate the following axiomatic:

- **Primitive Terms:** $\overline{\mathscr{G}}(I)$-distribution, domain, addition, derivative, and restriction of $\overline{\mathscr{G}}(I)$-distribution;
- **Precedent Theories:** Classical Logic, Set Theory, and $S$-Spaces Theory;
- **Axioms:** The ones stated ahead, regarding an abelian, surjective, with identity, and

coherent $S$-space,

$$\mathscr{G}(I) = \Bigg(\mathbb{G}\Big(\Gamma(I)\Big) = \Big\{\mathbb{G}(\gamma) = \Big(G(\gamma), H(\gamma)\Big)\Big\}_{\gamma\in\Gamma(I)},$$
$$i\Big(\Gamma^2(I)\Big) = \Big\{i_{(\gamma',\gamma)}\Big\}_{(\gamma',\gamma)\in\Gamma^2(I)},$$
$$\Theta\Big(\Delta(I)\Big) = \Big\{\Theta_{(\gamma',\gamma)}\Big\}_{(\gamma',\gamma)\in\Delta(I)}\Bigg),$$

arbitrarily fixed, with $G(\gamma) \cap G(\gamma') = \varnothing$ for $\gamma, \gamma' \in \Gamma(I)$ such that $\gamma \neq \gamma'$.

**Axiom 1** Every $\mathscr{G}(I)$-individual is a $\overline{\mathscr{G}}(I)$-distribution.

**Axiom 2** To each $\overline{\mathscr{G}}(I)$-distribution, $\overset{*}{g}$, corresponds a single $\gamma \in \Gamma(I)$, denominated the domain of $\overset{*}{g}$, in such a way that if $\overset{*}{g}$ is an $\mathscr{G}(I)$-individual, the domain $\gamma$ is the domain of the $\mathscr{G}(I)$-individual $\overset{*}{g}$.

**Notation.** We will denote by $\overset{*}{G}(\gamma)$ the class of $\overline{\mathscr{G}}(I)$-distributions with domain $\gamma$. It results from Axiom 1 and Axiom 2 that the class $G(\gamma)$ of $\mathscr{G}(I)$-individuals with domain $\gamma$ is a subclass of $\overset{*}{G}(\gamma)$: $G(\gamma) \subseteq \overset{*}{G}(\gamma)$.

**Axiom 3** The addition is a function (operation) that to each pair $(\overset{*}{g}, \overset{*}{h})$ of $\overline{\mathscr{G}}(I)$-distributions with the same domain $\gamma$, associates a (single) $\overline{\mathscr{G}}(I)$-distribution with the same domain $\gamma$, denoted by $\overset{*}{g} + \overset{*}{h}$ and denominated the sum of $\overset{*}{g}$ and $\overset{*}{h}$, in such a way that if $\overset{*}{g}$ and $\overset{*}{h}$ are $\mathscr{G}(I)$-individuals, then, the sum $\overset{*}{g} + \overset{*}{h}$ is the same as the one obtained through the addition of the group $G(\gamma)$ to which belong the $\mathscr{G}(I)$-individuals $\overset{*}{g}$ and $\overset{*}{h}$.

**Axiom 4** The derivatives of the $\overline{\mathscr{G}}(I)$-distributions with domain $\gamma$, whose class is denoted by $\overset{*}{H}(\gamma)$, are functions with domain and codomain equal to $\overset{*}{G}(\gamma)$, such that:

**(4-1)** $\overset{*}{\Phi}(\overset{*}{g}+\overset{*}{h}) = \overset{*}{\Phi}(\overset{*}{g})+\overset{*}{\Phi}(\overset{*}{h})$ for every $\overset{*}{g}, \overset{*}{h} \in \overset{*}{G}(\gamma)$ and every derivative $\overset{*}{\Phi} \in \overset{*}{H}(\gamma)$;

**(4-2)** each derivative $\overset{*}{\Phi} \in \overset{*}{H}(\gamma)$ is an extension to $\overset{*}{G}(\gamma)$ of a single $\Phi \in H(\gamma)$, that is, $\overset{*}{\Phi}(g) = \Phi(g)$ for every $g \in G(\gamma)_\Phi \subseteq \overset{*}{G}(\gamma)$; conversely, for each $\Phi \in H(\gamma)$ there exists a single derivative $\overset{*}{\Phi} \in \overset{*}{H}(\gamma)$ such that $= \Phi(g) = \overset{*}{\Phi}(g)$ for every $g \in G(\gamma)_\Phi \subseteq \overset{*}{G}(\gamma)$;

**(4-3)** the class $\overset{*}{H}(\gamma)$, equipped with the usual operation of composition of functions, is a semigroup isomorphic to the semigroup $H(\gamma)$, with the function

$$\overset{*}{\gamma} : H(\gamma) \longrightarrow \overset{*}{H}(\gamma)$$
$$\Phi \longmapsto \overset{*}{\gamma}(\Phi) := \text{the extension of } \Phi \text{ to } \overset{*}{G}(\gamma)$$

as an isomorphism.

**Notation.** The image of $\Phi \in H(\gamma)$ by the isomorphism $\overset{*}{\gamma} : H(\gamma) \longrightarrow \overset{*}{H}(\gamma)$ will be denoted by $\overset{*}{\Phi}$, that is,

$$\overset{*}{\Phi} := \overset{*}{\gamma}(\Phi).$$

**Axiom 5** Let $\gamma \in \Gamma(I)$ and $\Phi \in H(\gamma)$ be arbitrarily fixed. If $g \in G(\gamma)$ is such that $g \notin G(\gamma)_\Phi$, then

$$\overset{*}{\Phi}(g) \notin G(\gamma).$$

The next axiom makes reference to a term that is not primitive, namely, $\overset{*}{i}_{(\gamma',\gamma)}$. Therefore, it is necessary to define it through primitive terms.

**Definition.** For each $(\gamma', \gamma) \in \Gamma^2(I)$ we define $\overset{*}{i}_{(\gamma',\gamma)}$ as the following function:

$$\overset{*}{i}_{(\gamma',\gamma)} : \overset{*}{H}(\gamma) \longrightarrow \overset{*}{H}(\gamma')$$
$$\overset{*}{\Phi} \longmapsto \overset{*}{i}_{(\gamma',\gamma)}(\overset{*}{\Phi}) := \overset{*}{\gamma'}\left( i_{(\gamma',\gamma)}\left( (\overset{*}{\gamma})^{-1}(\overset{*}{\Phi}) \right) \right).$$

It results immediately that $\overset{*}{i}_{(\gamma',\gamma)}$ is an isomorphism from the semigroup $\overset{*}{H}(\gamma)$ onto the semigroup $\overset{*}{H}(\gamma')$.

**Axiom 6** The restrictions of $\overline{\mathscr{G}}(I)$-distributions are functions $\overset{*}{\Theta}_{(\gamma',\gamma)} : \overset{*}{G}(\gamma) \longrightarrow \overset{*}{G}(\gamma')$, a single one for each $(\gamma', \gamma) \in \Delta(I)$, that associate to each $\overline{\mathscr{G}}(I)$-distribution with domain $\gamma$, a single $\overline{\mathscr{G}}(I)$-distribution with domain $\gamma' \subseteq \gamma$, in such a way that:

**(6-1)** $\overset{*}{\Theta}_{(\gamma',\gamma)}(\overset{*}{g} + \overset{*}{h}) = \overset{*}{\Theta}_{(\gamma',\gamma)}(\overset{*}{g}) + \overset{*}{\Theta}_{(\gamma',\gamma)}(\overset{*}{h})$ for every $\overset{*}{g}, \overset{*}{h} \in \overset{*}{G}(\gamma)$;

**(6-2)** $\overset{*}{\Theta}_{(\gamma',\gamma)}(g) = \Theta_{(\gamma',\gamma)}(g)$ for every $g \in G(\gamma) \subseteq \overset{*}{G}(\gamma)$;

**(6-3)** $\overset{*}{\Theta}_{(\gamma'',\gamma')}\left(\overset{*}{\Theta}_{(\gamma',\gamma)}(\overset{*}{g})\right) = \overset{*}{\Theta}_{(\gamma'',\gamma)}(\overset{*}{g})$ for every $\overset{*}{g} \in \overset{*}{G}(\gamma)$ and $\gamma, \gamma', \gamma'' \in \Gamma(I)$ such that $\gamma'' \subseteq \gamma' \subseteq \gamma$;

**(6-4)** $\overset{*}{\Theta}_{(\gamma',\gamma)}\left(\overset{*}{\Phi}(\overset{*}{g})\right) = \overset{*}{i}_{(\gamma',\gamma)}(\overset{*}{\Phi})\left(\overset{*}{\Theta}_{(\gamma',\gamma)}(\overset{*}{g})\right)$ for every $\overset{*}{g} \in \overset{*}{G}(\gamma)$ and every $\overset{*}{\Phi} \in \overset{*}{H}(\gamma)$.

**Axiom 7** Let $\gamma \in \Gamma(I)$ be arbitrarily fixed. If $\overset{*}{g} \in \overset{*}{G}(\gamma)$ and $x \in \gamma$, there exists $\gamma' \in \Gamma(\gamma)$[23] such that $x \in \gamma'$ and

$$\overset{*}{\Theta}_{(\gamma',\gamma)}(\overset{*}{g}) = \overset{*}{\Phi}(g)$$

for some $\Phi \in H(\gamma')$ and $g \in G(\gamma')$.

[23] $\Gamma(\gamma) = \Gamma(I) \cap \mathcal{P}(\gamma)$, according to the footnote 18 at page 146; see also Definition 3.6.

In order for us to state the last axiom, a new term is introduced by the following definition.

**Definition.** Let $\gamma \in \Gamma(I)$ and $\overline{\xi}(\gamma) \subseteq \Gamma(\gamma)$ be a cover of $\gamma$. Let also

$$\left\{ \overset{*}{g}_\xi \right\}_{\xi \in \overline{\xi}(\gamma)},$$

be a family indexed by $\overline{\xi}(\gamma)$, whose members, $\overset{*}{g}_\xi$, are $\overline{\mathscr{G}}(I)$-distributions with domain $\xi$, that is, $\overset{*}{g}_\xi \in \overset{*}{G}(\xi)$ for every $\xi \in \overline{\xi}(\gamma)$. We sill say that the family $\{\overset{*}{g}_\xi\}_{\xi \in \overline{\xi}(\gamma)}$ is **coherent** if and only if

$$\overset{*}{\Theta}_{(\xi \cap \xi', \xi)}(\overset{*}{g}_\xi) = \overset{*}{\Theta}_{(\xi \cap \xi', \xi')}(\overset{*}{g}_{\xi'})$$

for every $\xi, \xi' \in \overline{\xi}(\gamma)$ such that $\xi \cap \xi' \neq \varnothing$.

**Axiom 8** Let $\gamma \in \Gamma(I)$ be arbitrarily fixed. If $\{\overset{*}{g}_\xi\}_{\xi \in \overline{\xi}(\gamma)}$ is a coherent family of $\overline{\mathscr{G}}(I)$-distributions, then, there exists a single $\overline{\mathscr{G}}(I)$-distribution $\overset{*}{g} \in \overset{*}{G}(\gamma)$ such that

$$\overset{*}{\Theta}_{(\xi, \gamma)}(\overset{*}{g}) = \overset{*}{g}_\xi$$

for every $\xi \in \overline{\xi}(\gamma)$.

## 5.22 The Consistency of the Axioms

The consistency of the $\overline{\mathscr{G}}(I)$-distributions axiomatic above formulated, can be established in an analogous manner to the corresponding proof provided in item 5.13 for the $\widetilde{\mathscr{G}}(I)$-distributions axioms, resorting now, not to the 1st TESS as we did in that case, but to the 2nd TESS, which is the theorem that applies in this case. In fact, keeping in mind that the $S$-space

$$\overline{\mathscr{G}}(I) = \left( \overline{\mathbb{G}}\Big(\Gamma(I)\Big) = \left\{ \overline{\mathbb{G}}(\gamma) = \Big(\overline{G}(\gamma), \overline{H}(\gamma)\Big) \right\}_{\gamma \in \Gamma(I)}, \right.$$
$$\overline{i}\Big(\Gamma^2(I)\Big) = \left\{ \overline{i}_{(\gamma', \gamma)} \right\}_{(\gamma', \gamma) \in \Gamma^2(I)},$$
$$\left. \overline{\Theta}\Big(\Delta(I)\Big) = \left\{ \overline{\Theta}_{(\gamma', \gamma)} \right\}_{(\gamma', \gamma) \in \Delta(I)} \right)$$

such as defined in the 2nd TESS (Theorem 4.33), is a locally closed and coherent extension of the extension $\widetilde{\mathscr{G}}(I)$ (described in the 1st TESS (Theorem 4.7)) of the underlying $S$-space $\mathscr{G}(I)$ to which the axioms in 5.21 refers to, the reader will find no difficulty to verify, resorting to the pertinent definitions, that the interpretation described below for the primitive terms constitutes a model of the $\overline{\mathscr{G}}(I)$-distributions axiomatic, i.e., that this axiomatic is consistent.

| **Primitive terms** | **Interpretation** |
|---|---|
| $\overline{\mathscr{G}}(I)$-distribution | $\overline{\mathscr{G}}(I)$-individual |
| $\overline{\mathscr{G}}(I)$-distribution domain | $\overline{\mathscr{G}}(I)$-individual domain |
| $\overline{\mathscr{G}}(I)$-distributions addition | $\overline{\mathscr{G}}(I)$-individuals addition |
| $\overline{\mathscr{G}}(I)$-distribution derivative | endomorphism $\overline{\Phi}$ of the semigroups $\overline{H}(\gamma)$ of the $S$-space $\overline{\mathscr{G}}(I)$ |
| $\overline{\mathscr{G}}(I)$-distribution restriction | homomorphism $\overline{\Theta}_{(\gamma',\gamma)}$ of the family $\overline{\Theta}(\Delta(I))$ of the $S$-space $\overline{\mathscr{G}}(I)$ |

Such as for the $\widetilde{\mathscr{G}}(I)$-distributions axiomatic formulated in 5.12, for which the interpretation described in 5.13 (table at page 234) can be taken as an "interpretation schema" where $\widetilde{\mathscr{G}}(I)$ represents any strict and closed extension of the $S$-space $\mathscr{G}(I)$ (recap item 5.13), also here the table above can be so interpreted. More precise and specifically, replacing, in the referred interpretation, the $S$-space $\overline{\mathscr{G}}(I)$ by another one isomorphic to it, one obtains a new model for the axiomatic, which, by what establishes the 2nd TESS, is isomorphic to the that corresponding to the $S$-space $\overline{\mathscr{G}}(I)$.

In the next items, 5.23 to 5.26, several logical consequences of the axioms formulated in 5.21 are obtained, as well as new concepts, such as, for example, the one of $\mathring{\mathscr{G}}(I)$-distribution, are introduced via definition. With these results, we will immediately see that any model of these axioms is, essentially, the one constituted by the interpretation given in the table above, and yet that the $\mathring{\mathscr{G}}(I)$-distributions and associated concepts (domain, addition, derivative, and restriction of $\mathring{\mathscr{G}}(I)$-distributions) compose a new model of the $\widetilde{\mathscr{G}}(I)$-distributions axioms formulated in 5.12.

In what follows, we prove, for some axioms, as examples, that the interpretation given in the table above satisfies the axiomatic requirements.

- Axiom 1.

  Translates, according to the given interpretation, to the following statement:

  "Every $\mathscr{G}(I)$-individual is a $\overline{\mathscr{G}}(I)$-individual"

  This is a true statement once that, taking into account the definition of individual (Definition 5.10(a)), we have:

  $$g \text{ is a } \mathscr{G}(I)\text{-individual if and only if } g \in G(\gamma) \text{ for some } \gamma \in \Gamma(I)$$

  and

  $$\overline{g} \text{ is a } \overline{\mathscr{G}}(I)\text{-individual if and only if } \overline{g} \in \overline{G}(\gamma) \text{ for some } \gamma \in \Gamma(I).$$

But, if $g \in G(\gamma)$, then, $g \in \overline{G}(\gamma)$, since $G(\gamma) \subseteq \overline{G}(\gamma)$.

- Axiom 2.

On the suggested interpretation, this axiom states that:

> "To each $\overline{\mathscr{G}}(I)$-individual, $\overline{g}$, corresponds a single $\gamma \in \Gamma(I)$, denominated the domain of $\overline{g}$, in such a way that if $\overline{g}$ is a $\mathscr{G}(I)$-individual, the domain $\gamma$ is the domain of the $\mathscr{G}(I)$-individual $\overline{g}$."

Also here we obtained a true statement. In fact, since the underlying $S$-space $\mathscr{G}(I)$ of the axiomatic in question is such that its groups, $G(\gamma)$, are two-by-two disjoints ($G(\gamma) \cap G(\gamma') = \varnothing$ if $\gamma \neq \gamma'$), the same occurs, due to Lemma 5.9, with the groups $\overline{G}(\gamma)$ of the extension $\overline{\mathscr{G}}(I)$ and, hence, if $\overline{g}$ is a $\overline{\mathscr{G}}(I)$-individual, then, by definition of individual, $\overline{g} \in \overline{G}(\gamma)$ for some $\gamma \in \Gamma(I)$ that, in this case is unique and it is, now taking into account the definition of domain of an individual (Definition 5.10(a)), the domain of $\overline{g}$ being $\overline{g}$ a $\mathscr{G}(I)$-individual (i.e., $\overline{g} \in G(\gamma) \subseteq \overline{G}(\gamma)$) or not (that is, $\overline{g} \in \overline{G}(\gamma)$ and $\overline{g} \notin G(\gamma)$).

- Axiom 7.

Translated, this axiom states that:

> "Let $\gamma \in \Gamma(I)$ be arbitrarily fixed. If $\overline{g} \in \overline{G}(\gamma)$ and $x \in \gamma$, there exists $\gamma' \in \Gamma(\gamma)$ such that $x \in \gamma'$ and
>
> $$\overline{\Theta}_{(\gamma',\gamma)}(\overline{g}) = \overline{\Phi}(g)$$
>
> for some $\Phi \in H(\gamma')$ and $g \in G(\gamma')$."

The truthiness of this statement follows, basically, from the fact of the extension $\overline{\mathscr{G}}(I)$ described in the 2nd TESS be a locally closed extension of $\widetilde{\mathscr{G}}(I)$ (defined in the 1st TESS). In fact, by the definition of locally closed extension (Definition 3.23(c-2)), for $\gamma \in \Gamma(I)$ and $\overline{g} \in \overline{G}(\gamma)$, if $x \in \gamma$, there exists $\gamma' \in \Gamma(\gamma)$ such that $x \in \gamma'$ and

$$\overline{\Theta}_{(\gamma',\gamma)}(\overline{g}) \in \widetilde{G}(\gamma').$$

Thus,

$$\overline{\Theta}_{(\gamma',\gamma)}(\overline{g}) = \widetilde{g}$$

for some $\widetilde{g} \in \widetilde{G}(\gamma')$. But $\widetilde{\mathscr{G}}(I)$ is a (strict and) closed extension of $\mathscr{G}(I)$ and, hence, according to the definition of closed extension of a $S$-space (Definition 3.23(c-1)),

$$\widetilde{g} = \widetilde{\Phi}(g)$$

for some $\Phi \in H(\gamma')$ and $g \in G(\gamma')$.

However,

$$\overline{\Phi}(g) = \widetilde{\Phi}(g),$$

since $\overline{\mathscr{G}}(I)$ is an extension of $\widetilde{\mathscr{G}}(I)$, and with this we have that

$$\overline{\Theta}_{(\gamma',\gamma)}(\overline{g}) = \overline{\Phi}(g)$$

for some $\Phi \in H(\gamma')$ and $g \in G(\gamma')$.

## 5.23 The $S$-Group $\mathring{\mathbb{G}}(\gamma) = (\mathring{G}(\gamma), \mathring{H}(\gamma))$

As previously remarked, it results from Axiom 1 and Axiom 2[24] that the class of $\mathscr{G}(I)$-individuals with domain $\gamma$, $G(\gamma)$, is a subclass of the class $\overset{*}{G}(\gamma)$ of $\overline{\mathscr{G}}(I)$-distributions with same domain $\gamma$: $G(\gamma) \subseteq \overset{*}{G}(\gamma)$ for each $\gamma \in \Gamma(I)$. Since, by Axiom 4, each derivative $\overset{*}{\Phi} \in \overset{*}{H}(\gamma)$ is a function with domain $\overset{*}{G}(\gamma)$ assuming values in $\overset{*}{G}(\gamma)$, are well-defined, as $\overline{\mathscr{G}}(I)$-distributions (with domain $\gamma$), the derivatives of $g \in G(\gamma)$, that is, $\overset{*}{\Phi}(g)$ where $\overset{*}{\Phi}(= \overset{*}{\gamma}(\Phi)) \in \overset{*}{H}(\gamma)$; we will denominate these derivatives of the $\mathscr{G}(I)$-individuals as $\mathring{\mathscr{G}}(I)$-distributions. More formally, the term $\mathring{\mathscr{G}}(I)$-distribution is introduced into our axiomatic through the following definition.

**Definition.** Let $\gamma \in \Gamma(I)$ and $\mathring{G}(\gamma)$ be the set defined by:

$$\mathring{G}(\gamma) := \left\{ \overset{*}{g} \in \overset{*}{G}(\gamma) \quad : \quad \overset{*}{g} = \overset{*}{\Phi}(g) \quad \text{for some} \quad \Phi \in H(\gamma) \quad \text{and} \quad g \in G(\gamma) \right\}.$$

The elements of $\mathring{G}(\gamma)$ are denominated $\mathring{\mathscr{G}}(I)$**-distributions with domain** $\gamma$.

In what follows, we highlight some properties of the $\mathring{\mathscr{G}}(I)$-distributions, implicit in our axioms. In each of these properties, $\gamma$ is an arbitrary element of $\Gamma(I)$.

**(a)** $G(\gamma) \subseteq \mathring{G}(\gamma)$

$\mathscr{G}(I)$ is a $S$-space with identity, that is, $I_{G(\gamma)} \in H(\gamma)$. Thus, by Axiom 4 ((4-2) and (4-3)) we have that

$$\overset{*}{\gamma}(I_{G(\gamma)}) = \overset{*}{I}_{G(\gamma)} \in \overset{*}{H}(\gamma)$$

and

$$\overset{*}{I}_{G(\gamma)}(g) = I_{G(\gamma)}(g) \quad \text{for every} \quad g \in G(\gamma).$$

But $I_{G(\gamma)}(g) = g$ for every $g \in G(\gamma)$ and, hence,

$$g = \overset{*}{I}_{G(\gamma)}(g) \quad \text{for every} \quad g \in G(\gamma),$$

which show us, taking into account the definition of $\mathring{G}(\gamma)$, that $g \in \mathring{G}(\gamma)$ for every $g \in G(\gamma)$, i.e., that $G(\gamma) \subseteq \mathring{G}(\gamma)$.

[24] In this item, as well as in the next three ones (5.24 to 5.26), any mention to an axiom or to an axiomatic without any other reference, must be interpreted as regarding the $\overline{\mathscr{G}}(I)$-distributions axiomatic in 5.21.

**(b)** For every $\mathring{g}, \mathring{h} \in \mathring{G}(\gamma)$ one has $\mathring{g} + \mathring{h} \in \mathring{G}(\gamma)$ and $\mathring{g} + \mathring{h} = \mathring{h} + \mathring{g}$.

In fact, let $\mathring{g}, \mathring{h} \in \mathring{G}(\gamma)$ be arbitrarily chosen. From the definition of $\mathring{G}(\gamma)$, there exist $\Phi, \Psi \in H(\gamma)$ amd $g, h \in G(\gamma)$ such that:

$$\mathring{g} = \overset{*}{\Phi}(g) \quad \text{and} \quad \mathring{h} = \overset{*}{\Psi}(h).$$

Since $\Phi : G(\gamma)_\Phi \longrightarrow G(\gamma)$ and $\Psi : G(\gamma)_\Psi \longrightarrow G(\gamma)$ are surjective homomorphisms (once that $\mathscr{G}(I)$ is a surjective $S$-space), there exist

$$g_1 \in G(\gamma)_\Psi \quad \text{and} \quad h_1 \in G(\gamma)_\Phi$$

such that

$$g = \Psi(g_1) \quad \text{and} \quad h = \Phi(h_1),$$

or, by what establishes Axiom 4 ((4-2)),

$$g = \overset{*}{\Psi}(g_1) \quad \text{and} \quad h = \overset{*}{\Phi}(h_1).$$

Thus, we have:

$$\mathring{g} + \mathring{h} = \overset{*}{\Phi}(g) + \overset{*}{\Psi}(h) = \overset{*}{\Phi}\Big(\overset{*}{\Psi}(g_1)\Big) + \overset{*}{\Psi}\Big(\overset{*}{\Phi}(h_1)\Big).$$

Resorting again to Axiom 4 ((4-3)), $\overset{*}{\gamma} : H(\gamma) \longrightarrow \overset{*}{H}(\gamma)$ is an isomorphism and, therefore,

$$\overset{*}{\gamma}(\Phi)\overset{*}{\gamma}(\Psi) = \overset{*}{\gamma}(\Phi\Psi) = \overset{*}{\gamma}(\Psi\Phi) = \overset{*}{\gamma}(\Psi)\overset{*}{\gamma}(\Phi)$$

that is,

$$\overset{*}{\Phi}\overset{*}{\Psi} = \overset{*}{\widehat{\Phi\Psi}} = \overset{*}{\widehat{\Psi\Phi}} = \overset{*}{\Psi}\overset{*}{\Phi}.$$

With this, and a new resort to Axiom 4 ((4-1)), we obtain:

$$\mathring{g} + \mathring{h} = \overset{*}{\Phi}\Big(\overset{*}{\Psi}(g_1)\Big) + \overset{*}{\Psi}\Big(\overset{*}{\Phi}(h_1)\Big) = \overset{*}{\widehat{\Phi\Psi}}(g_1 + h_1)$$

with $\Phi\Psi \in H(\gamma)$ and, by Axiom 3, with $g_1 + h_1 \in G(\gamma)$. From this, results that $\mathring{g} + \mathring{h} \in \mathring{G}(\gamma)$ and, also, taking into account that $G(\gamma)$ is an abelian group, that

$$\begin{aligned} \mathring{g} + \mathring{h} &= \overset{*}{\widehat{\Phi\Psi}}(g_1 + h_1) = \overset{*}{\widehat{\Phi\Psi}}(h_1 + g_1) = \overset{*}{\Phi}\overset{*}{\Psi}(h_1 + g_1) = \\ &= (\overset{*}{\Phi}\overset{*}{\Psi})(h_1) + (\overset{*}{\Phi}\overset{*}{\Psi})(g_1) = \overset{*}{\Psi}\Big(\overset{*}{\Phi}(h_1)\Big) + \overset{*}{\Phi}\Big(\overset{*}{\Psi}(g_1)\Big) = \\ &= \overset{*}{\Psi}(h) + \overset{*}{\Phi}(g) = \mathring{h} + \mathring{g}. \end{aligned}$$

**(c)** There exists $\mathring{0} \in \mathring{G}(\gamma)$ such that $\mathring{g} + \mathring{0} = \mathring{g}$ for every $\mathring{g} \in \mathring{G}(\gamma)$.

Let us take $\mathring{0} := 0 \in G(\gamma) \subseteq \mathring{G}(\gamma)$, being 0 the additive neutral of the group $G(\gamma)$. Let $\mathring{g} \in \mathring{G}(\gamma)$ be arbitrarily fixed. Hence, there exist $\Phi \in H(\gamma)$ and $g \in G(\gamma)$ such that

$$\mathring{g} = \overset{*}{\Phi}(g).$$

Since $\Phi$ is a homomorphism, $\Phi(0) = 0$ and, hence, by Axiom 4 ((4-2)), $\overset{*}{\Phi}(0) = \Phi(0) = 0$ and, with this, we have that

$$\mathring{g} + 0 = \overset{*}{\Phi}(g) + \overset{*}{\Phi}(0).$$

Resorting now to Axiom 4 ((4-1)) and Axiom 3 comes that

$$\mathring{g} + 0 = \overset{*}{\Phi}(g + 0) = \overset{*}{\Phi}(g) = \overset{*}{g}.$$

**(d)** For each $\mathring{g} \in \mathring{G}(\gamma)$ there exists $\mathring{g}_{(-)} \in \mathring{G}(\gamma)$ such that

$$\mathring{g} + \mathring{g}_{(-)} = \mathring{0}(= 0 \in G(\gamma)).$$

Let $\mathring{g} \in \mathring{G}(\gamma)$. From the definition of $\mathring{G}(\gamma)$ comes that

$$\mathring{g} = \overset{*}{\Phi}(g)$$

for some $\Phi \in H(\gamma)$ and $g \in G(\gamma)$. Let us take:

$$\mathring{g}_{(-)} := \overset{*}{\Phi}(-g)$$

where $-g \in G(\gamma)$ is the opposite element, in the group $G(\gamma)$, of $g \in G(\gamma)$, i.e., such that

$$g + (-g) = 0.$$

Clearly, $\mathring{g}_{(-)} \in \mathring{G}(\gamma)$. Furthermore, taking Axiom 4 ((4-1) and (4-2)) and Axiom 3 into account, we have:

$$\mathring{g} + \mathring{g}_{(-)} = \overset{*}{\Phi}(g) + \overset{*}{\Phi}(-g) = \overset{*}{\Phi}\Big(g + (-g)\Big) = \overset{*}{\Phi}(0) = \Phi(0)$$

But, $\Phi : G(\gamma)_\Phi \longrightarrow G(\gamma)$ is a homomorphism and, hence, $\Phi(0) = 0$. Thus, $\mathring{g} + \mathring{g}_{(-)} = 0$.

**(e)** $\mathring{f} + (\mathring{g} + \mathring{h}) = (\mathring{f} + \mathring{g}) + \mathring{h}$ for every $\mathring{f}, \mathring{g}, \mathring{h} \in \mathring{G}(\gamma)$.

Let $\mathring{f}, \mathring{g}, \mathring{h} \in \mathring{G}(\gamma)$ be arbitrarily chosen. We then have

$$\mathring{f} = \overset{*}{\Phi}(f), \quad \mathring{g} = \overset{*}{\Psi}(g) \quad \text{and} \quad \mathring{h} = \overset{*}{\varphi}(h)$$

with $\Phi, \Psi, \varphi \in H(\gamma)$ and $f, g, h \in G(\gamma)$.

Since $\Psi\varphi$, $\Phi\varphi$ and $\Phi\Psi$ are also members of $H(\gamma)$, then, they are surjective homomorphisms (onto $G(\gamma)$) and, hence, there exist

$$f_1 \in G(\gamma)_{\Psi\varphi}, \quad g_1 \in G(\gamma)_{\Phi\varphi} \quad \text{and} \quad h_1 \in G(\gamma)_{\Phi\Psi}$$

such that

$$f = (\Psi\varphi)(f_1), \quad g = (\Phi\varphi)(g_1), \quad \text{and} \quad h = (\Phi\Psi)(h_1)$$

or, by Axiom 4 ((4-2)),

$$f = \overset{*}{\widehat{\Psi\varphi}}(f_1), \quad g = \overset{*}{\widehat{\Phi\varphi}}(g_1), \quad \text{and} \quad h = \overset{*}{\widehat{\Phi\Psi}}(h_1).$$

From this, we can write

$$\mathring{f} + (\mathring{g} + \mathring{h}) = \overset{*}{\Phi}(f) + \Big(\overset{*}{\Psi}(g) + \overset{*}{\varphi}(h)\Big) = \overset{*}{\Phi}\Big(\overset{*}{\widehat{\Psi\varphi}}(f_1)\Big) + \\ + \Big[\overset{*}{\Psi}\Big(\overset{*}{\widehat{\Phi\varphi}}(g_1)\Big) + \overset{*}{\varphi}\Big(\overset{*}{\widehat{\Phi\Psi}}(h_1)\Big)\Big].$$

On the other hand, Axiom 4 ((4-3)) tells us that

$$\overset{*}{\Phi}(\overset{*}{\widehat{\Psi\varphi}}) = \overset{*}{\Psi}(\overset{*}{\widehat{\Phi\varphi}}) = \overset{*}{\varphi}(\overset{*}{\widehat{\Phi\Psi}}) = \overset{*}{\Phi}\overset{*}{\Psi}\overset{*}{\varphi}.$$

Thus, comes that

$$\mathring{f} + (\mathring{g} + \mathring{h}) = (\overset{*}{\Phi}\overset{*}{\Psi}\overset{*}{\varphi})(f_1) + \Big[(\overset{*}{\Phi}\overset{*}{\Psi}\overset{*}{\varphi})(g_1) + (\overset{*}{\Phi}\overset{*}{\Psi}\overset{*}{\varphi})(h_1)\Big].$$

Finally, resorting again to Axiom 4 ((4-1) and (4-3)), as well as Axiom 3 and the fact of $f_1$, $g_1$ and $h_1$ being members of the group $G(\gamma)$, we obtain:

$$\begin{aligned}
\mathring{f} + (\mathring{g} + \mathring{h}) &= (\overset{*}{\Phi}\overset{*}{\Psi}\overset{*}{\varphi})(f_1) + (\overset{*}{\Phi}\overset{*}{\Psi}\overset{*}{\varphi})(g_1 + h_1) = \\
&= (\overset{*}{\Phi}\overset{*}{\Psi}\overset{*}{\varphi})\Big(f_1 + (g_1 + h_1)\Big) = \\
&= (\overset{*}{\Phi}\overset{*}{\Psi}\overset{*}{\varphi})\Big((f_1 + g_1) + h_1\Big) = \\
&= (\overset{*}{\Phi}\overset{*}{\Psi}\overset{*}{\varphi})(f_1 + g_1) + (\overset{*}{\Phi}\overset{*}{\Psi}\overset{*}{\varphi})(h_1) = \\
&= \Big[\overset{*}{\Phi}\Big((\overset{*}{\Psi}\overset{*}{\varphi})(f_1)\Big) + \overset{*}{\Psi}\Big((\overset{*}{\Phi}\overset{*}{\varphi})(g_1)\Big)\Big] + \overset{*}{\varphi}\Big((\overset{*}{\Phi}\overset{*}{\Psi})(h_1)\Big) = \\
&= \Big[\overset{*}{\Phi}(f) + \overset{*}{\Psi}(g)\Big] + \overset{*}{\varphi}(h) = \\
&= (\mathring{f} + \mathring{g}) + \mathring{h}.
\end{aligned}$$

**(f)** The set $\mathring{G}(\gamma)$ with the addition of the $\overline{\mathscr{G}}(I)$-distributions (with domain $\gamma$) restricted to it, is an abelian group that has the group $G(\gamma)$ as a subgroup.

This follows immediately from (a) to (e) above.

**(g)** The derivatives of a $\mathring{\mathscr{G}}(I)$-distribution also are $\mathring{\mathscr{G}}(I)$-distributions, that is, if $\mathring{g} \in \mathring{G}(\gamma)$ and $\overset{*}{\Phi} \in \overset{*}{H}(\gamma)$, then

$$\overset{*}{\Phi}(\mathring{g}) \in \mathring{G}(\gamma).$$

In fact, if $\overset{\circ}{g} \in \overset{\circ}{G}(\gamma)$, then,

$$\overset{\circ}{g} = \overset{*}{\Psi}(g)$$

for some $\Psi \in H(\gamma)$ and $g \in G(\gamma)$. Thus, if $\overset{*}{\Phi} \in \overset{*}{H}(\gamma)$, we have that

$$\overset{*}{\Phi}(\overset{\circ}{g}) = \overset{*}{\Phi}\Big(\overset{*}{\Psi}(g)\Big)$$

or, taking into account Axiom 4 ((4-3)) and that $\Phi\Psi \in H(\gamma)$,

$$\overset{*}{\Phi}(\overset{\circ}{g}) = \overset{*}{\widehat{\Phi\Psi}}(g),$$

from where one concludes that $\overset{*}{\Phi}(g_0) \in \overset{\circ}{G}(\gamma)$.

Let now $\overset{*}{\Phi} : \overset{*}{G}(\gamma) \longrightarrow \overset{*}{G}(\gamma)$ be an arbitrarily chosen element in $\overset{*}{H}(\gamma)$ and let us consider the (usual) restriction function of $\overset{*}{\Phi}$ to the subset $\overset{\circ}{G}(\gamma)$ of its domain $\overset{*}{G}(\gamma)$, that is, the function $\overset{*}{\Phi}|_{\overset{\circ}{G}(\gamma)}$ defined by:

$$\begin{aligned} \overset{*}{\Phi}|_{\overset{\circ}{G}(\gamma)} : \overset{\circ}{G}(\gamma) &\longrightarrow \overset{*}{G}(\gamma) \\ \overset{\circ}{g} &\longmapsto \overset{*}{\Phi}|_{\overset{\circ}{G}(\gamma)}(\overset{\circ}{g}) := \overset{*}{\Phi}(\overset{\circ}{g}). \end{aligned}$$

The result (g) above shows us that the function $\overset{*}{\Phi}|_{\overset{\circ}{G}(\gamma)}$ assumes values in $\overset{\circ}{G}(\gamma)$, that is, it is a function from $\overset{\circ}{G}(\gamma)$ into $\overset{\circ}{G}(\gamma)$.

It is worth, for reasons that will become clear ahead, to highlight the classes of these functions through a definition.

**Definition.** Let $\gamma \in \Gamma(I)$ be arbitrarily fixed. $\overset{\circ}{H}(\gamma)$ is the set defined by

$$\overset{\circ}{H}(\gamma) := \Big\{\overset{*}{\Phi}|_{\overset{\circ}{G}(\gamma)} : \overset{*}{\Phi} \in \overset{*}{H}(\gamma)\Big\}.$$

We will say that $\overset{\circ}{H}(\gamma)$ is the class of the derivatives of the $\overset{\circ}{\mathscr{G}}(I)$-distributions with domain $\gamma$.

**(h)** The members of $\overset{\circ}{H}(\gamma)$ are endomorphisms on the group $\overset{\circ}{G}(\gamma)$.

In fact, let $\overset{\circ}{\Phi} \in \overset{\circ}{H}(\gamma)$ be arbitrarily chosen. By the definition of $\overset{\circ}{H}(\gamma)$, we have that

$$\begin{aligned} \overset{\circ}{\Phi} : \overset{\circ}{G}(\gamma) &\longrightarrow \overset{\circ}{G}(\gamma) \\ \overset{\circ}{g} &\longmapsto \overset{\circ}{\Phi}(\overset{\circ}{g}) = \overset{*}{\Phi}|_{\overset{\circ}{G}(\gamma)}(\overset{\circ}{g}) \end{aligned}$$

for some $\overset{*}{\Phi} \in \overset{*}{H}(\gamma)$.

As we saw in (f), $\overset{\circ}{G}(\gamma)$ is an abelian group; for any $\overset{\circ}{g}, \overset{\circ}{h} \in \overset{\circ}{G}(\gamma)$ we have $\overset{\circ}{g} + \overset{\circ}{h} \in \overset{\circ}{G}(\gamma)$ and

$$\overset{\circ}{\Phi}(\overset{\circ}{g} + \overset{\circ}{h}) = \overset{*}{\Phi}|_{\overset{\circ}{G}(\gamma)}(\overset{\circ}{g} + \overset{\circ}{h}) = \overset{*}{\Phi}(\overset{\circ}{g} + \overset{\circ}{h}).$$

On the other hand, by Axiom 4 ((4-1)),

$$\overset{*}{\Phi}(\mathring{g}+\mathring{h}) = \overset{*}{\Phi}(\mathring{g}) + \overset{*}{\Phi}(\mathring{h}) = \overset{*}{\Phi}|_{\mathring{G}(\gamma)}(\mathring{g}) + \overset{*}{\Phi}|_{\mathring{G}(\gamma)}(\mathring{h}) = \mathring{\Phi}(\mathring{g}) + \mathring{\Phi}(\mathring{h}).$$

Thus,

$$\mathring{\Phi}(\mathring{g}+\mathring{h}) = \mathring{\Phi}(\mathring{g}) + \mathring{\Phi}(\mathring{h})$$

that is, $\mathring{\Phi} : \mathring{G}(\gamma) \longrightarrow \mathring{G}(\gamma)$ is an endomorphisms on the group $\mathring{G}(\gamma)$.

**(i)** The function $\mathring{\gamma}$ defined ahead is bijective.

$$\begin{aligned} \mathring{\gamma} :\ & H(\gamma) \longrightarrow \mathring{H}(\gamma) \\ & \Phi \longmapsto \mathring{\gamma}(\Phi) := \overset{*}{\gamma}(\Phi)|_{\mathring{G}(\gamma)} \end{aligned}$$

(according to Axiom 4, $\overset{*}{\gamma}(\Phi) = \overset{*}{\Phi}$ is the extension of $\Phi$ to $\overset{*}{G}(\gamma)$).

Clearly, $\mathring{\gamma}$ is surjective since, by the definition of $\mathring{H}(\gamma)$, if $\mathring{\Phi} \in \mathring{H}(\gamma)$, then,

$$\mathring{\Phi} = \overset{*}{\Phi}|_{\mathring{G}(\gamma)} = \overset{*}{\gamma}(\Phi)|_{\mathring{G}(\gamma)}$$

for some $\Phi \in H(\gamma)$, that is, $\mathring{\Phi} = \mathring{\gamma}(\Phi)$. Let now $\Phi, \Psi \in H(\gamma)$ be such that $\Phi \neq \Psi$.

Taking into account that (according to what Axiom 4 ((4-2)) establishes) $\overset{*}{\Phi}$ and $\overset{*}{\Psi}$ are, respectively, extensions of $\Phi$ and $\Psi$ to $\overset{*}{G}(\gamma)$ and, yet, that

$$G(\gamma)_\Phi \subseteq G(\gamma) \subseteq \mathring{G}(\gamma) \subseteq \overset{*}{G}(\gamma)$$

and

$$G(\gamma)_\Psi \subseteq G(\gamma) \subseteq \mathring{G}(\gamma) \subseteq \overset{*}{G}(\gamma),$$

it results that

$$\mathring{\gamma}(\Phi) = \overset{*}{\Phi}|_{\mathring{G}(\gamma)} \quad \text{and} \quad \mathring{\gamma}(\Psi) = \overset{*}{\Psi}|_{\mathring{G}(\gamma)}$$

are, respectively, extensions of $\Phi$ and $\Psi$ to $\mathring{G}(\gamma)$. Therefore, since $\Phi \neq \Psi$, its extensions to $\mathring{G}(\gamma)$ can not be equal, i.e.,

$$\mathring{\gamma}(\Phi) \neq \mathring{\gamma}(\Psi)$$

and, hence, we conclude that $\mathring{\gamma}$ is injective.

**Notation.** The image of $\Phi \in H(\gamma)$ by the function $\mathring{\gamma} : H(\gamma) \longrightarrow \mathring{H}(\gamma)$ will be denoted by $\mathring{\Phi}$, that is,

$$\mathring{\Phi} := \mathring{\gamma}(\Phi).$$

**(j)** $\mathring{H}(\gamma)$, with the usual operation of composition of functions, is a semigroup isomorphic to the semigroup $H(\gamma)$ and the function $\mathring{\gamma} : H(\gamma) \longrightarrow \mathring{H}(\gamma)$ defined in (i) is an isomorphism.

We must prove that, for $\mathring{\Phi} = \mathring{\gamma}(\Phi)$ and $\mathring{\Psi} = \mathring{\gamma}(\Psi)$ arbitrarily chosen in $\mathring{H}(\gamma)$,

$$\mathring{\gamma}(\Phi)\mathring{\gamma}(\Psi) \in \mathring{H}(\gamma)$$

and that

$$\mathring{\gamma}(\Phi\Psi) = \mathring{\gamma}(\Phi)\mathring{\gamma}(\Psi).$$

In order to do so, let $\mathring{g} \in \mathring{G}(\gamma)$ be arbitrarily chosen and let us calculate $(\mathring{\gamma}(\Phi)\mathring{\gamma}(\Psi))(\mathring{g})$. We have:

$$\begin{aligned}\Big(\mathring{\gamma}(\Phi)\mathring{\gamma}(\Psi)\Big)(\mathring{g}) &= \mathring{\gamma}(\Phi)\Big(\mathring{\gamma}(\Psi)(\mathring{g})\Big) = \mathring{\gamma}(\Phi)\Big(\overset{*}{\Psi}|_{\mathring{G}(\gamma)}(\mathring{g})\Big) = \\ &= \mathring{\gamma}(\Phi)\Big(\overset{*}{\Psi}(\mathring{g})\Big) = \overset{*}{\Phi}|_{\mathring{G}(\gamma)}\Big(\overset{*}{\Psi}(\mathring{g})\Big) = \overset{*}{\Phi}\Big(\overset{*}{\Psi}(\mathring{g})\Big) = \\ &= \overset{*}{\widehat{\Phi\Psi}}(\mathring{g}) = \overset{*}{\widehat{\Phi\Psi}}|_{\mathring{G}(\gamma)}(\mathring{g}) = \mathring{\gamma}(\Phi\Psi)(\mathring{g}).\end{aligned}$$

Once that $\mathring{g} \in \mathring{G}(\gamma)$ is arbitrary, we obtain

$$\mathring{\gamma}(\Phi)\mathring{\gamma}(\Psi) = \mathring{\gamma}(\Phi\Psi),$$

and since $\Phi\Psi \in H(\gamma)$ (for $H(\gamma)$ is a semigroup), it results from this last equation not only that $\mathring{\gamma}(\Phi)\mathring{\gamma}(\Psi)$ belongs to $\mathring{H}(\gamma)$ (being the image of $\Phi\Psi \in H(\gamma)$ by $\mathring{\gamma}$), i.e., that $\mathring{H}(\gamma)$ is closed regarding the operation of composition and, therefore, it is a semigroup, but also that the bijective function $\mathring{\gamma} : H(\gamma) \longrightarrow \mathring{H}(\gamma)$ defined in (k) is an isomorphism.

**(k)** $\mathring{H}(\gamma)$ is a prolongation of $H(\gamma)$ to $\mathring{G}(\gamma)$.

In fact, by (h) and (j) above we have that $\mathring{H}(\gamma)$ is a semigroup of endomorphisms on $\mathring{G}(\gamma)$ isomorphic to the semigroup $H(\gamma)$, being the function

$$\begin{aligned}\mathring{\gamma} : H(\gamma) &\longrightarrow \mathring{H}(\gamma) \\ \Phi &\longmapsto \mathring{\gamma}(\Phi) := \overset{*}{\gamma}(\Phi)|_{\mathring{G}(\gamma)}\end{aligned}$$

an isomorphism such that

$$\mathring{\Phi}(g) = \Phi(g) \quad \text{for every} \quad \Phi \in H(\gamma) \quad \text{and} \quad g \in G(\gamma)_\Phi.$$

But, this means, taking Definition 1.6(c) into account, that $\mathring{H}(\gamma)$ is a prolongation to $\mathring{G}(\gamma)$ of $H(\gamma)$.

**(l)** For each $\gamma \in \Gamma(I)$, the ordered pair

$$\mathring{\mathbb{G}}(\gamma) := \Big(\mathring{G}(\gamma), \mathring{H}(\gamma)\Big)$$

is a $S$-group that is a strict and closed extension of the $S$-group

$$\mathbb{G}(\gamma) = \Big(G(\gamma), H(\gamma)\Big).$$

In fact, by (f) we know that $\mathring{G}(\gamma)$ is an abelian group that has the group $G(\gamma)$ as subgroup, and by (k) we have that $\mathring{H}(\gamma)$ is a prolongation of $H(\gamma)$ to $\mathring{G}(\gamma)$. Thus being, and taking into account the definitions of $S$-group (Definition 1.4(a)) and of $S$-group extension (Definition 1.6(d)), we conclude that $\mathring{\mathbb{G}}(\gamma)$ is a $S$-group and also an extension of the $S$-group $\mathbb{G}(\gamma)$. Furthermore, trivially follows from Axiom 5 that $\mathring{\mathbb{G}}(\gamma)$ is a strict extension (Definition 1.10(b)) of $\mathbb{G}(\gamma)$. Finally, due to the own definition of the set $\mathring{G}(\gamma)$ (at page 258), we have that for each $\mathring{g} \in \mathring{G}(\gamma)$, there exist $g \in G(\gamma)$ and $\Phi \in H(\gamma)$ such that

$$\mathring{g} = \overset{*}{\Phi}(g) = \overset{*}{\Phi}|_{\mathring{G}(\gamma)}(g) = \mathring{\Phi}(g)$$

and, hence, according to Definition 1.14, $\mathring{\mathbb{G}}(\gamma)$ is a closed extension of the $S$-group $\mathbb{G}(\gamma)$.

## 5.24 The $S$-Space $\mathring{\mathscr{G}}(I) = (\mathring{\mathbb{G}}(\Gamma(I)), \mathring{\imath}(\Gamma^2(I)), \mathring{\Theta}(\Delta(I)))$

In what follows, via definition, we introduce new terms into our axiomatic: $\mathring{\mathbb{G}}(\Gamma(I))$ and $\mathring{\imath}(\Gamma^2(I))$.

**Definition.** $\mathring{\mathbb{G}}(\Gamma(I))$ and $\mathring{\imath}(\Gamma^2(I))$ are families indexed, respectively, by $\Gamma(I)$ and $\Gamma^2(I)$, defined as follows:

$$\mathring{\mathbb{G}}\Big(\Gamma(I)\Big) := \Big\{\mathring{\mathbb{G}}(\gamma) = \Big(\mathring{G}(\gamma), \mathring{H}(\gamma)\Big)\Big\}_{\gamma\in\Gamma(I)}$$

with $\mathring{\mathbb{G}}(\gamma)$ as defined in 5.23, and

$$\mathring{\imath}\Big(\Gamma^2(I)\Big) := \Big\{\mathring{\imath}_{(\gamma',\gamma)}\Big\}_{(\gamma',\gamma)\in\Gamma^2(I)}$$

where, for each $(\gamma', \gamma) \in \Gamma^2(I)$, $\mathring{\imath}_{(\gamma',\gamma)}$ is the following function:

$$\begin{aligned} \mathring{\imath}_{(\gamma',\gamma)} : \mathring{H}(\gamma) &\longrightarrow \mathring{H}(\gamma') \\ \mathring{\Phi} &\longmapsto \mathring{\imath}_{(\gamma',\gamma)}(\mathring{\Phi}) := \mathring{\gamma}'\bigg(i_{(\gamma',\gamma)}\Big((\mathring{\gamma})^{-1}(\mathring{\Phi})\Big)\bigg). \end{aligned}$$

Regarding the families $\mathring{\mathbb{G}}(\Gamma(I))$ and $\mathring{\imath}(\Gamma^2(I))$ above defined we have that:

**(a)** $\mathring{\mathbb{G}}(\Gamma(I))$ is a family of $S$-groups that is a strict and closed extension of the family of $S$-groups $\mathbb{G}(\Gamma(I))$ (of the $S$-space $\mathscr{G}(I)$).

Immediate consequence of the definitions of family of $S$-groups (Definition 3.7), of strict and closed extension of a family of $S$-groups (Definition 3.8) and of 5.23(l).

**(b)** $\mathring{\imath}(\Gamma^2(I))$ is the extension of the bonding $i(\Gamma^2(I))$ (of the $S$-space $\mathscr{G}(I)$) to the family $\mathring{\mathbb{G}}(\Gamma(I))$.

This follows directly from (a) and the definition of extension of a bonding (Definition 3.22).

**(c)** The ordered pair

$$\left(\mathring{\mathbb{G}}\Big(\Gamma(I)\Big), \mathring{i}\Big(\Gamma^2(I)\Big)\right)$$

is a bonded family that is an extension of the bonded family (of the $S$-space $\mathscr{G}(I)$) $(\mathbb{G}(\Gamma(I)), i(\Gamma^2(I)))$.

It results immediately from (a) and (b) above and from the definitions of bonded family (Defnition 3.11) and extension of a bonded family (Definition 3.22).

Let us turn now our attention to the restrictions of $\overline{\mathscr{G}}(I)$-distributions. The header of Axiom 6 states that:

> The restrictions of $\overline{\mathscr{G}}(I)$-distributions are functions $\overset{*}{\Theta}_{(\gamma',\gamma)} : \overset{*}{G}(\gamma) \longrightarrow \overset{*}{G}(\gamma')$, a **single** one for each $(\gamma', \gamma) \in \Delta(I)$, that associate to each $\overline{\mathscr{G}}(I)$-distribution with domain $\gamma$, a single $\overline{\mathscr{G}}(I)$-distribution with domain $\gamma' \subseteq \gamma$, ...

In other words, the referred axiom, in its header, states the existence of a function, let us say $\lambda$, with domain $\Delta(I)$, whose image is the set of **all** the restrictions of $\overline{\mathscr{G}}(I)$-distributions, such that, for each $(\gamma', \gamma) \in \Delta(I)$,

> $\lambda(\gamma', \gamma) =$ the unique restriction of $\overline{\mathscr{G}}(I)$-distribution with domain $\overset{*}{G}(\gamma)$ and codomain $\overset{*}{G}(\gamma')$;

furthermore, the axiom introduce the symbol "$\overset{*}{\Theta}_{(\gamma',\gamma)}$" to denote $\lambda(\gamma', \gamma)$, that is,

$$\overset{*}{\Theta}_{(\gamma',\gamma)} \coloneqq \lambda(\gamma', \gamma).$$

An alternative manner of referring to the function $\lambda$, which will be used here, consists of treating it as being the family

$$\left\{\overset{*}{\Theta}_{(\gamma',\gamma)}\right\}_{(\gamma',\gamma)\in\Delta(I)}.$$

In short, Axiom 6 ensures the definition below.[25]

**Definition.** $\overset{*}{\Theta}(\Delta(I))$ is a family indexed by $\Delta(I)$ and defined by:

$$\overset{*}{\Theta}\Big(\Delta(I)\Big) \coloneqq \left\{\overset{*}{\Theta}_{(\gamma',\gamma)}\right\}_{(\gamma',\gamma)\in\Delta(I)}.$$

A property of the members of the family $\overset{*}{\Theta}(\Delta(I))$ that worths to highlight at this point is presented below.

[25] The definition, in 5.14(l), of the family $\widehat{\Theta}(\Delta(I))$ of the restrictions of the $\widetilde{\mathscr{G}}(I)$-distributions, is equally backed by the header of Axiom 7 of the axiomatic in 5.12.

**(d)** $\overset{*}{\Theta}_{(\gamma',\gamma)}(\mathring{g}) \in \mathring{G}(\gamma')$ for every $(\gamma',\gamma) \in \Delta(I)$ and $\mathring{g} \in \mathring{G}(\gamma)$.

In fact, for $(\gamma',\gamma) \in \Delta(I)$ and $\mathring{g} \in \mathring{G}(\gamma)$ arbitrarily chosen we have, by the definition of $\mathring{G}(\gamma)$ (at page 258), that $\mathring{g} = \overset{*}{\Phi}(g)$ for some $\Phi \in H(\gamma)$ and $g \in G(\gamma)$. Hence, resorting to Axiom 6 ((6-4) and (6-2)),

$$\overset{*}{\Theta}_{(\gamma',\gamma)}(\mathring{g}) = \overset{*}{\Theta}_{(\gamma',\gamma)}\Big(\overset{*}{\Phi}(g)\Big) = \overset{*}{\imath}_{(\gamma',\gamma)}(\overset{*}{\Phi})\Big(\overset{*}{\Theta}_{(\gamma',\gamma)}(g)\Big) = \overset{*}{\imath}_{(\gamma',\gamma)}(\overset{*}{\Phi})\Big(\Theta_{(\gamma',\gamma)}(g)\Big).$$

Since $\Theta_{(\gamma',\gamma)}(g) \in G(\gamma')$ and $\overset{*}{\imath}_{(\gamma',\gamma)}(\overset{*}{\Phi}) \in \overset{*}{H}(\gamma')$, the expression above obtained for $\overset{*}{\Theta}_{(\gamma',\gamma)}(\mathring{g})$ and the definition of $\mathring{G}(\gamma')$ allow us to conclude that

$$\overset{*}{\Theta}_{(\gamma',\gamma)}(\mathring{g}) \in \mathring{G}(\gamma').$$

Let now $\overset{*}{\Theta}_{(\gamma',\gamma)} : \overset{*}{G}(\gamma) \longrightarrow \overset{*}{G}(\gamma')$ be arbitrarily chosen in $\overset{*}{\Theta}(\Delta(I))$ and let us consider its restriction (in the usual sense of restriction of a function) to the subset $\mathring{G}(\gamma)$ of its domain $\overset{*}{G}(\gamma)$, that is, the function $\overset{*}{\Theta}_{(\gamma',\gamma)}|_{\mathring{G}(\gamma)}$ with domain $\mathring{G}(\gamma)$ defined by

$$\overset{*}{\Theta}_{(\gamma',\gamma)}|_{\mathring{G}(\gamma)}(\mathring{g}) := \overset{*}{\Theta}_{(\gamma',\gamma)}(\mathring{g}) \quad \text{for every} \quad \mathring{g} \in \mathring{G}(\gamma).$$

The result (d) shows us that this function assume values in $\mathring{G}(\gamma')$, that is, it regards a function from $\mathring{G}(\gamma)$ into $\mathring{G}(\gamma')$:

$$\begin{aligned} \overset{*}{\Theta}_{(\gamma',\gamma)}|_{\mathring{G}(\gamma)} : \mathring{G}(\gamma) &\longrightarrow \mathring{G}(\gamma') \\ \mathring{g} &\longmapsto \overset{*}{\Theta}_{(\gamma',\gamma)}|_{\mathring{G}(\gamma)}(\mathring{g}) := \overset{*}{\Theta}_{(\gamma',\gamma)}(\mathring{g}). \end{aligned}$$

These functions, $\overset{*}{\Theta}_{(\gamma',\gamma)}\Big|_{\mathring{G}(\gamma)}$, can be taken as being members of a well-defined family indexed by $\Delta(I)$, namely, the family $\mathring{\Theta}(\Delta(I))$, defined ahead.

**Definition.** $\mathring{\Theta}(\Delta(I))$ is the family indexed by $\Delta(I)$ defined by:

$$\mathring{\Theta}\Big(\Delta(I)\Big) := \Big\{\mathring{\Theta}_{(\gamma',\gamma)}\Big\}_{(\gamma',\gamma)\in\Delta(I)}$$

where

$$\mathring{\Theta}_{(\gamma',\gamma)} := \overset{*}{\Theta}_{(\gamma',\gamma)}|_{\mathring{G}(\gamma)} \quad \text{for each} \quad (\gamma',\gamma) \in \Delta(I),$$

that is, $\mathring{\Theta}_{(\gamma',\gamma)}$ is the function

$$\begin{aligned} \mathring{\Theta}_{(\gamma',\gamma)} : \mathring{G}(\gamma) &\longrightarrow \mathring{G}(\gamma') \\ \mathring{g} &\longmapsto \mathring{\Theta}_{(\gamma',\gamma)}(\mathring{g}) := \overset{*}{\Theta}_{(\gamma',\gamma)}(\mathring{g}). \end{aligned}$$

The families $\Theta(\Delta(I))$, $\mathring{\Theta}(\Delta(I))$ and $\overset{*}{\Theta}(\Delta(I))$ are "injective families", that is, they are such that:

**(e)** Let $(\gamma', \gamma), (\beta', \beta) \in \Delta(I)$. One has that:

**(e-1)** $\Theta_{(\gamma',\gamma)} = \Theta_{(\beta',\beta)}$ (if and) only if $\gamma' = \beta'$ and $\gamma = \beta$;

**(e-2)** $\overset{\circ}{\Theta}_{(\gamma',\gamma)} = \overset{\circ}{\Theta}_{(\beta',\beta)}$ (if and) only if $\gamma' = \beta'$ and $\gamma = \beta$;

**(e-3)** $\overset{*}{\Theta}_{(\gamma',\gamma)} = \overset{*}{\Theta}_{(\beta',\beta)}$ (if and) only if $\gamma' = \beta'$ and $\gamma = \beta$;

In order for us to prove this, let us remember ourselves that the groups, $G(\gamma)$, of the $S$-space $\mathscr{G}(I)$ are such that

$$G(\gamma) \cap G(\gamma') = \varnothing \quad \text{if} \quad \gamma \neq \gamma'.$$

We also know that, for each $\gamma \in \Gamma(I)$, $\overset{\circ}{G}(\gamma)$ is an abelian group that has $G(\gamma)$ as subgroup and, hence, taking Lemma 5.9 into account, we conclude that the groups $\overset{\circ}{G}(\gamma)$ are, as the groups $G(\gamma)$, such that

$$\overset{\circ}{G}(\gamma) \cap \overset{\circ}{G}(\gamma') = \varnothing \quad \text{if} \quad \gamma \neq \gamma'.$$

Regarding the sets $\overset{*}{G}(\gamma)$, we will prove in item 5.25 that they also are abelian groups with $\overset{\circ}{G}(\gamma)$ and $G(\gamma)$ as subgroups. Thus, again by Lemma 5.9, we have

$$\overset{*}{G}(\gamma) \cap \overset{*}{G}(\gamma') = \varnothing \quad \text{if} \quad \gamma \neq \gamma'.$$

With these data we will prove, for instance, (e-3). Hence, let us suppose that the functions $\overset{*}{\Theta}_{(\gamma',\gamma)}$ and $\overset{*}{\Theta}_{(\beta',\beta)}$ are equal:

$$\overset{*}{\Theta}_{(\gamma',\gamma)} = \overset{*}{\Theta}_{(\beta',\beta)}.$$

Thus, its domains are equal, i.e.,

$$\overset{*}{G}(\gamma) = \overset{*}{G}(\beta)$$

and, hence, $\gamma = \beta$ for if $\gamma \neq \beta$ we would have $\overset{*}{G}(\gamma) \cap \overset{*}{G}(\beta) = \varnothing$. Now, $\overset{*}{\Theta}_{(\gamma',\gamma)}$ and $\overset{*}{\Theta}_{(\beta',\beta=\gamma)}$ assume values in $\overset{*}{G}(\gamma')$ and $\overset{*}{G}(\beta')$, respectively. Therefore, if $\overset{*}{\Theta}_{(\gamma',\gamma)} = \overset{*}{\Theta}_{(\beta',\beta=\gamma)}$, it is necessary that $\gamma' = \beta'$ once that, again, if $\gamma' \neq \beta'$, then $\overset{*}{G}(\gamma') \cap \overset{*}{G}(\beta') = \varnothing$ and, consequently, the values of the functions in question would not be equal.

Regarding the families $\overset{\circ}{\mathbb{G}}(\Gamma(I))$, $\overset{\circ}{\imath}(\Gamma^2(I))$, and $\overset{\circ}{\Theta}(\Delta(I))$, introduced into our axiomatic through the definitions above formulated, our axioms allow to state, among other things, that:

**(f)** $\overset{\circ}{\Theta}(\Delta(I))$ is a family of homomorphisms, i.e., for each $(\gamma', \gamma) \in \Delta(I)$, $\overset{\circ}{\Theta}_{(\gamma',\gamma)} : \overset{\circ}{G}(\gamma) \longrightarrow \overset{\circ}{G}(\gamma')$ is a homomorphism from the group $\overset{\circ}{G}(\gamma)$ into the group $\overset{\circ}{G}(\gamma')$.

In fact, let $\overset{\circ}{\Theta}_{(\gamma',\gamma)} : \overset{\circ}{G}(\gamma) \longrightarrow \overset{\circ}{G}(\gamma')$ be an arbitrary member of the family $\overset{\circ}{\Theta}(\Delta(I))$. We know, by 5.23(f), that $\overset{\circ}{G}(\gamma)$ and $\overset{\circ}{G}(\gamma')$ are (abelian) groups. Furthermore, for $\overset{\circ}{g}, \overset{\circ}{h} \in \overset{\circ}{G}(\gamma)$ arbitrarily chosen, we have, resorting to the definition of $\overset{\circ}{\Theta}_{(\gamma',\gamma)}$ and Axiom 6 ((6-1)), that:

$$\overset{\circ}{\Theta}_{(\gamma',\gamma)}(\overset{\circ}{g} + \overset{\circ}{h}) = \overset{*}{\Theta}_{(\gamma',\gamma)}(\overset{\circ}{g} + \overset{\circ}{h}) = \overset{*}{\Theta}_{(\gamma',\gamma)}(\overset{\circ}{g}) + \overset{*}{\Theta}_{(\gamma',\gamma)}(\overset{\circ}{h}) = \overset{\circ}{\Theta}_{(\gamma',\gamma)}(\overset{\circ}{g}) + \overset{\circ}{\Theta}_{(\gamma',\gamma)}(\overset{\circ}{h}).$$

**(g)** For every $\gamma, \gamma', \gamma'' \in \Gamma(I)$ such that $\gamma'' \subseteq \gamma' \subseteq \gamma$ and every $\mathring{g} \in \mathring{G}(\gamma)$,

$$\mathring{\Theta}_{(\gamma'',\gamma')}\Big(\mathring{\Theta}_{(\gamma',\gamma)}(\mathring{g})\Big) = \mathring{\Theta}_{(\gamma'',\gamma)}(\mathring{g}).$$

An immediate consequence of the definition given for the members of the family $\mathring{\Theta}(\Delta(I))$ and Axiom 6 ((6-3)).

**(h)** For each $(\gamma', \gamma) \in \Delta(I)$, every $\mathring{\Phi} \in \mathring{H}(\gamma)$ and every $\mathring{g} \in \mathring{G}(\gamma)$,

$$\mathring{\Theta}_{(\gamma',\gamma)}\Big(\mathring{\Phi}(\mathring{g})\Big) = \mathring{i}_{(\gamma',\gamma)}(\mathring{\Phi})\Big(\mathring{\Theta}_{(\gamma',\gamma)}(\mathring{g})\Big).$$

By the definitions of $\mathring{\Theta}_{(\gamma',\gamma)}$ and $\mathring{\Phi}$ we have

$$\mathring{\Theta}_{(\gamma',\gamma)}\Big(\mathring{\Phi}(\mathring{g})\Big) = \overset{*}{\Theta}_{(\gamma',\gamma)}\Big(\overset{*}{\Phi}(\mathring{g})\Big),$$

and, by Axiom 6((6-4)), that

$$\overset{*}{\Theta}_{(\gamma',\gamma)}\Big(\overset{*}{\Phi}(\mathring{g})\Big) = \overset{*}{i}_{(\gamma',\gamma)}(\overset{*}{\Phi})\Big(\overset{*}{\Theta}_{(\gamma',\gamma)}(\mathring{g})\Big) = \overset{*}{i}_{(\gamma',\gamma)}(\overset{*}{\Phi})\Big(\mathring{\Theta}_{(\gamma',\gamma)}(\mathring{g})\Big).$$

Thus,

$$\mathring{\Theta}_{(\gamma',\gamma)}\Big(\mathring{\Phi}(\mathring{g})\Big) = \overset{*}{i}_{(\gamma',\gamma)}(\overset{*}{\Phi})\Big(\mathring{\Theta}_{(\gamma',\gamma)}(\mathring{g})\Big). \tag{h-1}$$

On the other hand, the definitions of $\overset{*}{i}_{(\gamma',\gamma)}$ (at page 254) and $\mathring{i}_{(\gamma',\gamma)}$ (at page 265) give us

$$\overset{*}{i}_{(\gamma',\gamma)}(\overset{*}{\Phi}) = \overset{*}{\gamma'}\left(i_{(\gamma',\gamma)}\Big((\overset{*}{\gamma})^{-1}(\overset{*}{\Phi})\Big)\right)$$

and

$$\mathring{i}_{(\gamma',\gamma)}(\mathring{\Phi}) = \mathring{\gamma}'\left(i_{(\gamma',\gamma)}\Big((\mathring{\gamma})^{-1}(\mathring{\Phi})\Big)\right),$$

and, since $(\overset{*}{\gamma})^{-1}(\overset{*}{\Phi}) = \Phi = (\mathring{\gamma})^{-1}(\mathring{\Phi})$, we have:

$$\overset{*}{i}_{(\gamma',\gamma)}(\overset{*}{\Phi}) = \overset{*}{\gamma'}\Big(i_{(\gamma',\gamma)}(\Phi)\Big)$$

and

$$\mathring{i}_{(\gamma',\gamma)}(\mathring{\Phi}) = \mathring{\gamma}'\Big(i_{(\gamma',\gamma)}(\Phi)\Big).$$

Taking now

$$\Psi := i_{(\gamma',\gamma)}(\Phi) \in H(\gamma')$$

one obtains

$$\overset{*}{i}_{(\gamma',\gamma)}(\overset{*}{\Phi}) = \overset{*}{\gamma'}(\Psi) = \overset{*}{\Psi}$$

and

$$\mathring{i}_{(\gamma',\gamma)}(\mathring{\Phi}) = \mathring{\gamma}'(\Psi) = \mathring{\Psi}.$$

Thus, in possession of these data and remembering ourselves that $\mathring{\Theta}_{(\gamma',\gamma)}(\mathring{g}) \in \mathring{G}(\gamma')$ and $\overset{*}{\Psi}|_{\mathring{G}(\gamma')} = \mathring{\Psi}$, we obtain that:

$$\overset{*}{\imath}_{(\gamma',\gamma)}(\overset{*}{\Phi})\left(\mathring{\Theta}_{(\gamma',\gamma)}(\mathring{g})\right) = \overset{*}{\Psi}\left(\mathring{\Theta}_{(\gamma',\gamma)}(\mathring{g})\right) = \mathring{\Psi}\left(\mathring{\Theta}_{(\gamma',\gamma)}(\mathring{g})\right) = \mathring{\imath}_{(\gamma',\gamma)}(\mathring{\Phi})\left(\mathring{\Theta}_{(\gamma',\gamma)}(\mathring{g})\right).$$

Taking this last result into (h-1) we get that

$$\mathring{\Theta}_{(\gamma',\gamma)}\left(\mathring{\Phi}(\mathring{g})\right) = \mathring{\imath}_{(\gamma',\gamma)}(\mathring{\Phi})\left(\mathring{\Theta}_{(\gamma',\gamma)}(\mathring{g})\right).$$

**(i)** $\mathring{\Theta}(\Delta(I))$ is a restriction for the bonded family $(\mathring{\mathbb{G}}(\Gamma(I)), \mathring{\imath}(\Gamma^2(I)))$ (in (c)).

In fact, of the four requirements of Definition 3.13 for $\mathring{\Theta}(\Delta(I))$ to be a restriction for the bonded family $(\mathring{\mathbb{G}}(\Gamma(I)), \mathring{\imath}(\Gamma^2(I)))$, the results (f), (g) and (h) above show that the three first ones are fullfiled. Regarding the fourth requirement, this one is trivially fullfiled at the present case once that $\mathring{\Phi}$ is an endomorphism on $\mathring{G}(\gamma)$.

**(j)** $\mathring{\Theta}(\Delta(I))$ is a prolongation of the family $\Theta(\Delta(I))$ (of the $S$-space $\mathscr{G}(I)$) to the bonded family $(\mathring{\mathbb{G}}(\Gamma(I)), \mathring{\imath}(\Gamma^2(I)))$ (in (c)).

By (c), $(\mathring{\mathbb{G}}(\Gamma(I)), \mathring{\imath}(\Gamma^2(I)))$ is an extension of the bonded family $(\mathbb{G}(\Gamma(I)), \imath(\Gamma^2(I)))$ (of the $S$-space $\mathscr{G}(I)$) and, by (i), $\mathring{\Theta}(\Delta(I))$ is a restriction for $(\mathring{\mathbb{G}}(\Gamma(I)), \mathring{\imath}(\Gamma^2(I)))$. Furthermore, for every $(\gamma',\gamma) \in \Delta(I)$ and every $g \in G(\gamma) \subseteq \mathring{G}(\gamma)$, since $\overset{*}{\Theta}_{(\gamma',\gamma)}(g) = \mathring{\Theta}_{(\gamma',\gamma)}(g)$, follows from Axiom 6 ((6-2)) that

$$\mathring{\Theta}_{(\gamma',\gamma)}(g) = \Theta_{(\gamma',\gamma)}(g),$$

and, hence, according to Definition 3.23(a), we have established the result (j).

**(k)** The triplet $\mathring{\mathscr{G}}(I)$ defined below is a $S$-space that is a strict and closed extension of the $S$-space $\mathscr{G}(I)$.

$$\mathring{\mathscr{G}}(I) := \left(\mathring{\mathbb{G}}\left(\Gamma(I)\right), \mathring{\imath}\left(\Gamma^2(I)\right), \mathring{\Theta}\left(\Delta(I)\right)\right).$$

An immediate consequence from the results above and the definitions of $S$-space (Definition 3.15), $S$-space extension (Definition 3.23(b)) and strict and closed extension of a $S$-space (Definition 3.23(c-1)).

**(l)** $\mathring{\mathscr{G}}(I)$ is a model of the $\widetilde{\mathscr{G}}(I)$-distributions axiomatic. More precisely, the interpretation given to the primitive terms that figure into the axioms formulated in 5.12, through the elements that integrate the $S$-space $\mathring{\mathscr{G}}(I)$, such as prescribe the table ahead, define a model of that axiomatic.

| **Primitive terms** | **Interpretation** |
|---|---|
| $\widetilde{\mathscr{G}}(I)$-distribution | $\overset{\circ}{\mathscr{G}}(I)$-individual |
| $\widetilde{\mathscr{G}}(I)$-distribution domain | $\overset{\circ}{\mathscr{G}}(I)$-individual domain |
| $\widetilde{\mathscr{G}}(I)$-distributions addition | $\overset{\circ}{\mathscr{G}}(I)$-individuals addition |
| $\widetilde{\mathscr{G}}(I)$-distribution derivative | endomorphism $\overset{\circ}{\Phi}$ of the semigroups $\overset{\circ}{H}(\gamma)$ of the $S$-space $\overset{\circ}{\mathscr{G}}(I)$ |
| $\widetilde{\mathscr{G}}(I)$-distribution restriction | homomorphism $\overset{\circ}{\Theta}_{(\gamma',\gamma)}$ of the family $\overset{\circ}{\Theta}(\Delta(I))$ of the $S$-space $\overset{\circ}{\mathscr{G}}(I)$ |

In fact, as we saw in item 5.13, the interpretation there provided (at page 236) for the $\widetilde{\mathscr{G}}(I)$-distributions axiomatic (formulated in 5.12), can be seem as an "interpretation schema" in which the $S$-space $\widetilde{\mathscr{G}}(I)$ (defined in the 1st TESS) can be replaced by any other $S$-space isomorphic to it: for each such $S$-space, the corresponding interpretation constitutes a model of that axiomatic. Hence, since $\overset{\circ}{\mathscr{G}}(I)$ as above defined in (k) is a strict and closed extension of the $S$-space $\mathscr{G}(I)$ we have, taking the 1st TESS into account, that $\overset{\circ}{\mathscr{G}}(I)$ is isomorphic to $\widetilde{\mathscr{G}}(I)$ and, therefore, the interpretation given above, which is exactly that of item 5.13 with the $S$-space $\widetilde{\mathscr{G}}(I)$ replaced by $\overset{\circ}{\mathscr{G}}(I)$, determines a model of the axiomatic in 5.12.

## 5.25 The $S$-Group $\overset{*}{\mathbb{G}}(\gamma) = (\overset{*}{G}(\gamma), \overset{*}{H}(\gamma))$

In this item (and the next one), we will continue the task of exploring the $\widetilde{\mathscr{G}}(I)$-distributions axiomatic, in the sense of obtaining new logical consequences from its axioms. In what follows, $\gamma$ is an arbitrary element of $\Gamma(I)$ and, moreover, the reader must keep in mind the footnote 24 (at page 258).

**(a)** For any $\overset{*}{g}, \overset{*}{h} \in \overset{*}{G}(\gamma)$ one has

$$\overset{*}{g} + \overset{*}{h} = \overset{*}{h} + \overset{*}{g}.$$

Let $\overset{*}{g}, \overset{*}{h} \in \overset{*}{G}(\gamma)$ be arbitrarily fixed and take:

$$\overset{*}{e} := \overset{*}{g} + \overset{*}{h} \quad \text{and} \quad \overset{*}{f} := \overset{*}{h} + \overset{*}{g}$$

By Axiom 3 we have $\overset{*}{e}, \overset{*}{f} \in \overset{*}{G}(\gamma)$.

Let also $\overline{\xi}(\gamma) \subseteq \Gamma(\gamma)$ be a cover of $\gamma$ such that:

$$\overset{*}{\Theta}_{(\xi,\gamma)}(\overset{*}{g}), \overset{*}{\Theta}_{(\xi,\gamma)}(\overset{*}{h}) \in \overset{\circ}{G}(\xi) \tag{a-1}$$

for every $\xi \in \overline{\xi}(\gamma)$.

The existence of one such cover of $\gamma$ is ensured by Axiom 7. In fact, taking into account the definition of $\mathring{G}(\delta)$, $\delta \in \Gamma(I)$, in 5.23, the referred axiom allows us to state that: for each $x \in \gamma$ there exist $\gamma_x, \gamma'_x \in \Gamma(\gamma)$ such that $x \in \gamma_x$, $x \in \gamma'_x$ and

$$\overset{*}{\Theta}_{(\gamma_x,\gamma)}(\overset{*}{g}) \in \mathring{G}(\gamma_x) \quad \text{and} \quad \overset{*}{\Theta}_{(\gamma'_x,\gamma)}(\overset{*}{h}) \in \mathring{G}(\gamma'_x).$$

Now, we define

$$\overline{\xi}(\gamma) := \Big\{\gamma_x \cap \gamma'_x : x \in \gamma\Big\}.$$

Clearly, $\overline{\xi}(\gamma)$ is a cover of $\gamma$ and $\overline{\xi}(\gamma) \subseteq \Gamma(\gamma)$.

Since $\overset{*}{\Theta}_{(\gamma_x,\gamma)}(\overset{*}{g})$ and $\overset{*}{\Theta}_{(\gamma'_x,\gamma)}(\overset{*}{h})$ belongs, respectively, to the domains of the homomorphisms

$$\mathring{\Theta}_{(\gamma_x\cap\gamma'_x,\gamma_x)} : \mathring{G}(\gamma_x) \longrightarrow \mathring{G}(\gamma_x \cap \gamma'_x) \quad \text{and} \quad \mathring{\Theta}_{(\gamma_x\cap\gamma'_x,\gamma'_x)} : \mathring{G}(\gamma'_x) \longrightarrow \mathring{G}(\gamma_x \cap \gamma'_x)$$

which (both), in turn, assume values in $\mathring{G}(\gamma_x \cap \gamma'_x)$ for each $x \in \gamma$, we obtain that:

$$\mathring{\Theta}_{(\gamma_x\cap\gamma'_x,\gamma_x)}\left(\overset{*}{\Theta}_{(\gamma_x,\gamma)}(\overset{*}{g})\right) \in \mathring{G}(\gamma_x \cap \gamma'_x)$$

and

$$\mathring{\Theta}_{(\gamma_x\cap\gamma'_x,\gamma'_x)}\left(\overset{*}{\Theta}_{(\gamma'_x,\gamma)}(\overset{*}{h})\right) \in \mathring{G}(\gamma_x \cap \gamma'_x)$$

for every $x \in \gamma$.

On the other hand, resorting to the definition of $\mathring{\Theta}_{(\delta',\delta)}$, $(\delta', \delta) \in \Delta(I)$, in 5.24 (at page 267), as well as Axiom 6((6-3)), we obtain that:

$$\mathring{\Theta}_{(\gamma_x\cap\gamma'_x,\gamma_x)}\left(\overset{*}{\Theta}_{(\gamma_x,\gamma)}(\overset{*}{g})\right) = \overset{*}{\Theta}_{(\gamma_x\cap\gamma'_x,\gamma_x)}\left(\overset{*}{\Theta}_{(\gamma_x,\gamma)}(\overset{*}{g})\right) = \overset{*}{\Theta}_{(\gamma_x\cap\gamma'_x,\gamma)}(\overset{*}{g})$$

and

$$\mathring{\Theta}_{(\gamma_x\cap\gamma'_x,\gamma'_x)}\left(\overset{*}{\Theta}_{(\gamma'_x,\gamma)}(\overset{*}{h})\right) = \overset{*}{\Theta}_{(\gamma_x\cap\gamma'_x,\gamma'_x)}\left(\overset{*}{\Theta}_{(\gamma'_x,\gamma)}(\overset{*}{h})\right) = \overset{*}{\Theta}_{(\gamma_x\cap\gamma'_x,\gamma)}(\overset{*}{h}).$$

Taking now these results into the last two pertinence relations, we conclude that

$$\overset{*}{\Theta}_{(\gamma_x\cap\gamma'_x,\gamma)}(\overset{*}{g}) \in \mathring{G}(\gamma_x \cap \gamma'_x) \quad \text{and} \quad \overset{*}{\Theta}_{(\gamma_x\cap\gamma'_x,\gamma)}(\overset{*}{h}) \in \mathring{G}(\gamma_x \cap \gamma'_x)$$

for every $x \in \gamma$, or, equivalently, that

$$\overset{*}{\Theta}_{(\xi,\gamma)}(\overset{*}{g}) \in \mathring{G}(\xi) \quad \text{and} \quad \overset{*}{\Theta}_{(\xi,\gamma)}(\overset{*}{h}) \in \mathring{G}(\xi)$$

for every $\xi \in \overline{\xi}(\gamma) = \{\gamma_x \cap \gamma'_x : x \in \gamma\}$.

It is thus ensured the existence of a cover $\overline{\xi}(\gamma) \subseteq \Gamma(\gamma)$ of $\gamma$ that attends condition (a-1).

Now, remembering ourselves that $\mathring{G}(\xi)$, for each $\xi \in \overline{\xi}(\gamma)$, is an abelian group and, therefore, that

$$\overset{*}{\Theta}_{(\xi,\gamma)}(\overset{*}{g}) + \overset{*}{\Theta}_{(\xi,\gamma)}(\overset{*}{h}) = \overset{*}{\Theta}_{(\xi,\gamma)}(\overset{*}{h}) + \overset{*}{\Theta}_{(\xi,\gamma)}(\overset{*}{g})$$

for every $\xi \in \overline{\xi}(\gamma)$, and taking yet into account what establishes Axiom 6((6-1)), we obtain for $\overset{*}{e} = \overset{*}{g} + \overset{*}{h}$ and $\overset{*}{f} = \overset{*}{h} + \overset{*}{g}$, that:

$$\begin{aligned}\overset{*}{\Theta}_{(\xi,\gamma)}(\overset{*}{e}) &= \overset{*}{\Theta}_{(\xi,\gamma)}(\overset{*}{g} + \overset{*}{h}) = \overset{*}{\Theta}_{(\xi,\gamma)}(\overset{*}{g}) + \overset{*}{\Theta}_{(\xi,\gamma)}(\overset{*}{h}) = \\ &= \overset{*}{\Theta}_{(\xi,\gamma)}(\overset{*}{h}) + \overset{*}{\Theta}_{(\xi,\gamma)}(\overset{*}{g}) = \overset{*}{\Theta}_{(\xi,\gamma)}(\overset{*}{h} + \overset{*}{g}) = \\ &= \overset{*}{\Theta}_{(\xi,\gamma)}(\overset{*}{f})\end{aligned}$$

for every $\xi \in \overline{\xi}(\gamma)$. Now, the family $\{\overset{*}{\Theta}_{(\xi,\gamma)}(\overset{*}{e})\}_{\xi\in\overline{\xi}(\gamma)}$ is, clearly, a coherent family. Thus, taking Axiom 8 into account, the last expression allows us to conclude that

$$\overset{*}{e} = \overset{*}{f},$$

i.e., that $\overset{*}{g} + \overset{*}{h} = \overset{*}{h} + \overset{*}{g}$.

**(b)** There exists $\overset{*}{0} \in \overset{*}{G}(\gamma)$ such that $\overset{*}{g} + \overset{*}{0} = \overset{*}{g}$ for every $\overset{*}{g} \in \overset{*}{G}(\gamma)$.

Let $\overset{*}{g} \in \overset{*}{G}(\gamma)$ be arbitrarily chosen and take

$$\overset{*}{0} := \overset{\circ}{0},$$

where $\overset{\circ}{0} \in \overset{\circ}{G}(\gamma) \subseteq \overset{*}{G}(\gamma)$ is the additive neutral of the group $\overset{\circ}{G}(\gamma)$. Let now $\overline{\xi}(\gamma) \subseteq \Gamma(\gamma)$ be a cover of $\gamma$ such that

$$\overset{*}{\Theta}_{(\xi,\gamma)}(\overset{*}{g}) \in \overset{\circ}{G}(\xi) \quad \text{for every} \quad \xi \in \overline{\xi}(\gamma).$$

The existence of one such cover was proved in (a) above.

Now, by Axiom 6((6-1)) and taking into account the definition of $\overset{\circ}{\Theta}_{(\xi,\gamma)}$ ($\overset{\circ}{\Theta}_{(\xi,\gamma)} := \overset{*}{\Theta}_{(\xi,\gamma)}|_{\overset{\circ}{G}(\gamma)}$), as well as the fact that $\overset{\circ}{\Theta}_{(\xi,\gamma)} : \overset{\circ}{G}(\gamma) \longrightarrow \overset{\circ}{G}(\xi)$ being a homomorphism, we have:

$$\begin{aligned}\overset{*}{\Theta}_{(\xi,\gamma)}(\overset{*}{g} + \overset{*}{0}) &= \overset{*}{\Theta}_{(\xi,\gamma)}(\overset{*}{g}) + \overset{*}{\Theta}_{(\xi,\gamma)}(\overset{*}{0}) = \overset{*}{\Theta}_{(\xi,\gamma)}(\overset{*}{g}) + \overset{*}{\Theta}_{(\xi,\gamma)}(\overset{\circ}{0}) = \\ &= \overset{*}{\Theta}_{(\xi,\gamma)}(\overset{*}{g}) + \overset{\circ}{\Theta}_{(\xi,\gamma)}(\overset{\circ}{0}) = \overset{*}{\Theta}_{(\xi,\gamma)}(\overset{*}{g}) + \overset{\circ}{0}_{\xi}\end{aligned}$$

where $\overset{\circ}{0}_{\xi} \in \overset{\circ}{G}(\xi)$ is the additive neutral of the group $\overset{\circ}{G}(\xi)$. But, $\overset{*}{\Theta}_{(\xi,\gamma)}(\overset{*}{g}) \in \overset{\circ}{G}(\xi)$ for every $\xi \in \overline{\xi}(\gamma)$, and since $\overset{\circ}{G}(\xi)$, with the addition operation of $\overset{*}{G}(\xi)$ restricted to it, is a group (as we saw in 5.23(f)), we have that $\overset{*}{\Theta}_{(\xi,\gamma)}(\overset{*}{g}) + \overset{\circ}{0}_{\xi} = \overset{*}{\Theta}_{(\xi,\gamma)}(\overset{*}{g})$ for every $\xi \in \overline{\xi}(\gamma)$. Thus,

$$\overset{*}{\Theta}_{(\xi,\gamma)}(\overset{*}{g} + \overset{*}{0}) = \overset{*}{\Theta}_{(\xi,\gamma)}(\overset{*}{g}) \quad \text{for every} \quad \xi \in \overline{\xi}(\gamma).$$

On the other hand, the family $\{\overset{*}{\Theta}_{(\xi,\gamma)}(\overset{*}{g})\}_{\xi\in\overline{\xi}(\gamma)}$ is, clearly, a coherent familly and, thus, taking Axiom 8 into account, from the last equation above one concludes that

$$\overset{*}{g} + \overset{*}{0} = \overset{*}{g}.$$

**(c)** For each $\overset{*}{g} \in \overset{*}{G}(\gamma)$, there exists $\overset{*}{g}_{(-)} \in \overset{*}{G}(\gamma)$ such that

$$\overset{*}{g} + \overset{*}{g}_{(-)} = \overset{*}{0} (= \overset{\circ}{0} \in \overset{\circ}{G}(\gamma)).$$

In fact, let $\overline{\xi} \subseteq \Gamma(\gamma)$ be a cover of $\gamma$ such that

$$\overset{*}{\Theta}_{(\xi,\gamma)}(\overset{*}{g}) \in \overset{\circ}{G}(\xi) \quad \text{for every} \quad \xi \in \overline{\xi}(\gamma).$$

As we know, Axiom 7 ensures the existence of such a cover.

Let now, for each $\xi \in \overline{\xi}(\gamma)$, $-\overset{*}{\Theta}_{(\xi,\gamma)}(\overset{*}{g}) \in \overset{\circ}{G}(\xi)$ be the opposite element of $\overset{*}{\Theta}_{(\xi,\gamma)}(\overset{*}{g}) \in \overset{\circ}{G}(\xi)$ in the group $\overset{\circ}{G}(\xi)$, that is, such that

$$\overset{*}{\Theta}_{(\xi,\gamma)}(\overset{*}{g}) + \left(-\overset{*}{\Theta}_{(\xi,\gamma)}(\overset{*}{g})\right) = \overset{\circ}{0}_{\xi},$$

where $\overset{\circ}{0}_{\xi} \in \overset{\circ}{G}(\xi)$ is the additive neutral of the group $\overset{\circ}{G}(\xi)$.

Let us now consider the family $\{-\overset{*}{\Theta}_{(\xi,\gamma)}(\overset{*}{g})\}_{\xi\in\bar{\xi}(\gamma)}$. This family is coherent, since for $\xi, \xi' \in \overline{\xi}(\gamma)$ such that $\xi \cap \xi' \neq \varnothing$ we have (taking into account that $\overset{\circ}{\Theta}_{(\delta',\delta)} = \overset{*}{\Theta}_{(\delta',\delta)}|_{\overset{\circ}{G}(\gamma)}$ and yet that $\overset{\circ}{\Theta}_{(\delta',\delta)}$ is a homomorphism) that:

$$\begin{aligned}\overset{*}{\Theta}_{(\xi\cap\xi',\xi)}\left(-\overset{*}{\Theta}_{(\xi,\gamma)}(\overset{*}{g})\right) &= \overset{\circ}{\Theta}_{(\xi\cap\xi',\xi)}\left(-\overset{*}{\Theta}_{(\xi,\gamma)}(\overset{*}{g})\right) = \\ &= -\overset{\circ}{\Theta}_{(\xi\cap\xi',\xi)}\left(\overset{*}{\Theta}_{(\xi,\gamma)}(\overset{*}{g})\right) = -\overset{*}{\Theta}_{(\xi\cap\xi',\xi)}\left(\overset{*}{\Theta}_{(\xi,\gamma)}(\overset{*}{g})\right)\end{aligned}$$

and, hence, resorting to Axiom 6((6-3)),

$$\overset{*}{\Theta}_{(\xi\cap\xi',\xi)}\left(-\overset{*}{\Theta}_{(\xi,\gamma)}(\overset{*}{g})\right) = -\overset{*}{\Theta}_{(\xi\cap\xi',\gamma)}(\overset{*}{g}).$$

Analogously, one obtains that

$$\overset{*}{\Theta}_{(\xi\cap\xi',\xi')}\left(-\overset{*}{\Theta}_{(\xi',\gamma)}(\overset{*}{g})\right) = -\overset{*}{\Theta}_{(\xi\cap\xi',\gamma)}(\overset{*}{g}),$$

which, along with the last equation, show that the family $\{-\overset{*}{\Theta}_{(\xi,\gamma)}(\overset{*}{g})\}_{\xi\in\bar{\xi}(\gamma)}$ is coherent. Thus, by what Axiom 8 establishes, there exists a single element in $\overset{*}{G}(\gamma)$, let us say $\overset{*}{g}_{(-)}$, such that

$$\overset{*}{\Theta}_{(\xi,\gamma)}(\overset{*}{g}_{(-)}) = -\overset{*}{\Theta}_{(\xi,\gamma)}(\overset{*}{g})$$

for every $\xi \in \overline{\xi}(\gamma)$. With this, and now taking Axiom 6((6-1)) into account, we have:

$$\overset{*}{\Theta}_{(\xi,\gamma)}(\overset{*}{g} + \overset{*}{g}_{(-)}) = \overset{*}{\Theta}_{(\xi,\gamma)}(\overset{*}{g}) + \overset{*}{\Theta}_{(\xi,\gamma)}(\overset{*}{g}_{(-)}) = \overset{*}{\Theta}_{(\xi,\gamma)}(\overset{*}{g}) + \left(-\overset{*}{\Theta}_{(\xi,\gamma)}(\overset{*}{g})\right) = \overset{\circ}{0}_{\xi}$$

for every $\xi \in \overline{\xi}(\gamma)$.

On the other hand, for $\overset{*}{0} = \overset{\circ}{0} \in \overset{\circ}{G}(\gamma)$, the additive neutral of the group $\overset{\circ}{G}(\gamma)$, we have (once that, again, $\overset{\circ}{\Theta}_{(\xi,\gamma)} = \overset{*}{\Theta}_{(\xi,\gamma)}|_{\overset{\circ}{G}(\gamma)}$ is a homomorphism from the group $\overset{\circ}{G}(\gamma)$ into the group $\overset{\circ}{G}(\xi)$) that

$$\overset{*}{\Theta}_{(\xi,\gamma)}(\overset{*}{0}) = \overset{\circ}{\Theta}_{(\xi,\gamma)}(\overset{\circ}{0}) = \overset{\circ}{0}_{\xi} \quad \text{for every} \quad \xi \in \overline{\xi}(\gamma),$$

Therefore, we obtain that

$$\overset{*}{\Theta}_{(\xi,\gamma)}(\overset{*}{g} + \overset{*}{g}_{(-)}) = \overset{*}{\Theta}_{(\xi,\gamma)}(\overset{*}{0}) \quad \text{for every} \quad \xi \in \overline{\xi}(\gamma),$$

and, since, clearly, $\{\overset{*}{\Theta}_{(\xi,\gamma)}(\overset{*}{0})\}_{\xi\in\bar{\xi}(\gamma)}$ is a coherent family, this result together with Axiom 8 allow us to conclude that

$$\overset{*}{g} + \overset{*}{g}_{(-)} = \overset{*}{0}.$$

**(d)** $\overset{*}{f} + (\overset{*}{g} + \overset{*}{h}) = (\overset{*}{f} + \overset{*}{g}) + \overset{*}{h}$ for every $\overset{*}{f}, \overset{*}{g}, \overset{*}{h} \in \overset{*}{G}(\gamma)$.

Again endorsed by Axiom 7, let $\overline{\xi}(\gamma) \subseteq \Gamma(\gamma)$ be a cover of $\gamma$ such that

$$\overset{*}{\Theta}_{(\xi,\gamma)}(\overset{*}{f}), \overset{*}{\Theta}_{(\xi,\gamma)}(\overset{*}{g}), \overset{*}{\Theta}_{(\xi,\gamma)}(\overset{*}{h}) \in \overset{\circ}{G}(\xi)$$

for every $\xi \in \overline{\xi}(\gamma)$, being $\overset{*}{f}$, $\overset{*}{g}$ and $\overset{*}{h}$ arbitrarily chosen in $\overset{*}{G}(\gamma)$.

Resorting to Axiom 6((6-1)), Axiom 3 and the fact of $\overset{\circ}{G}(\xi)$ being a group, we calculate, for each $\xi \in \overline{\xi}(\gamma)$:

$$\begin{aligned}\overset{*}{\Theta}_{(\xi,\gamma)}\Big(\overset{*}{f} + (\overset{*}{g} + \overset{*}{h})\Big) &= \overset{*}{\Theta}_{(\xi,\gamma)}(\overset{*}{f}) + \overset{*}{\Theta}_{(\xi,\gamma)}(\overset{*}{g} + \overset{*}{h}) = \overset{*}{\Theta}_{(\xi,\gamma)}(\overset{*}{f}) + \Big(\overset{*}{\Theta}_{(\xi,\gamma)}(\overset{*}{g}) + \overset{*}{\Theta}_{(\xi,\gamma)}(\overset{*}{h})\Big) = \\ &= \Big(\overset{*}{\Theta}_{(\xi,\gamma)}(\overset{*}{f}) + \overset{*}{\Theta}_{(\xi,\gamma)}(\overset{*}{g})\Big) + \overset{*}{\Theta}_{(\xi,\gamma)}(\overset{*}{h}) = \overset{*}{\Theta}_{(\xi,\gamma)}(\overset{*}{f} + \overset{*}{g}) + \overset{*}{\Theta}_{(\xi,\gamma)}(\overset{*}{h}) = \\ &= \overset{*}{\Theta}_{(\xi,\gamma)}\Big((\overset{*}{f} + \overset{*}{g}) + \overset{*}{h}\Big)\end{aligned}$$

for every $\xi \in \overline{\xi}(\gamma)$.

From this last equation, together with Axiom 8 (taking into account that $\{\overset{*}{\Theta}_{(\xi,\gamma)}(\overset{*}{f} + (\overset{*}{g} + \overset{*}{h}))\}_{\xi\in\bar{\xi}(\gamma)}$ is a coherent family) we obtain that

$$\overset{*}{f} + (\overset{*}{g} + \overset{*}{h}) = (\overset{*}{f} + \overset{*}{g}) + \overset{*}{h}.$$

**(e)** The set $\overset{*}{G}(\gamma)$ with the $\overline{\mathscr{G}}(I)$-distributions (with domain $\gamma$) addition is an abelian group that has the groups $G(\gamma)$ and $\overset{\circ}{G}(\gamma)$ as subgroups.

In fact, (a) to (d) above tell us that $\overset{*}{G}(\gamma)$ with the addition mentioned in Axiom 3 is an abelian group. Furthermore, by the definition of $\overset{\circ}{G}(\gamma)$ in 5.23 (page 258),

$$G(\gamma) \subseteq \overset{\circ}{G}(\gamma) \subseteq \overset{*}{G}(\gamma),$$

and, according to what we saw in 5.23(f), $\overset{\circ}{G}(\gamma)$ with the addition above referred is an abelian group that has the group $G(\gamma)$ as subgroup.

**(f)** $\overset{*}{H}(\gamma)$ is a prolongation of $H(\gamma)$ to $\overset{*}{G}(\gamma)$.

As we saw above, $\overset{*}{G}(\gamma)$ is a group. For $\overset{*}{\Phi} : \overset{*}{G}(\gamma) \longrightarrow \overset{*}{G}(\gamma)$ arbitrarily chosen in $\overset{*}{H}(\gamma)$ we have, by Axiom 4((4-1) and (4-2)), that

$$\overset{*}{\Phi}(\overset{*}{g} + \overset{*}{h}) = \overset{*}{\Phi}(\overset{*}{g}) + \overset{*}{\Phi}(\overset{*}{h}) \quad \text{for every} \quad \overset{*}{g}, \overset{*}{h} \in \overset{*}{G}(\gamma)$$

and, for a single $\Phi \in H(\gamma)$,

$$\overset{*}{\Phi}(g) = \Phi(g) \quad \text{for every} \quad g \in G(\gamma)_{\Phi} \subseteq \overset{*}{G}(\gamma),$$

that is, $\overset{*}{\Phi}$ is an endomorphism on the group $\overset{*}{G}(\gamma)$ and also an extension of $\Phi$ to $\overset{*}{G}(\gamma)$. In short, according to Definition 1.6(b), $\overset{*}{\Phi}$ is a prolongation of $\Phi$ to $\overset{*}{G}(\gamma)$.

Furthermore, yet by Axiom 4((4-3)), $\overset{*}{H}(\gamma)$, with the usual operation of composition of functions, is a semigroup isomorphic to the semigroup $H(\gamma)$ and the function

$$\begin{aligned} \overset{*}{\gamma} : H(\gamma) &\longrightarrow \overset{*}{H}(\gamma) \\ \Phi &\longmapsto \overset{*}{\gamma}(\Phi) = \overset{*}{\Phi} := \text{the extension of } \Phi \text{ to } \overset{*}{G}(\gamma) \end{aligned}$$

is an isomorphism.

This, however, according to Definition 1.6(c), means that $\overset{*}{H}(\gamma)$ is a prolongation of $H(\gamma)$ to $\overset{*}{G}(\gamma)$, that is, (f).

**(g)** $\overset{*}{H}(\gamma)$ is a prolongation of $\overset{\circ}{H}(\gamma)$ to $\overset{*}{G}(\gamma)$.

We already know (by (f)) that the members of $\overset{*}{H}(\gamma)$ are endomorphisms on the group $\overset{*}{G}(\gamma)$ and that (by (e)) $\overset{\circ}{G}(\gamma) \subseteq \overset{*}{G}(\gamma)$ is a subgroup of $\overset{*}{G}(\gamma)$. We also know, according to 5.23(j), that $\overset{\circ}{H}(\gamma)$ (defined at page 262), with the usual operation of composition of functions, is a semigroup isomorphic to the semigroup $H(\gamma)$ and that the function (defined in 5.23(i))

$$\begin{aligned} \overset{\circ}{\gamma} : H(\gamma) &\longrightarrow \overset{\circ}{H}(\gamma) \\ \Phi &\longmapsto \overset{\circ}{\gamma}(\Phi) = \overset{\circ}{\Phi} := \overset{*}{\gamma}(\Phi)|_{\overset{\circ}{G}(\gamma)} \end{aligned}$$

is an isomorphism.

Since, by Axiom 4((4-3)),

$$\begin{aligned} \overset{*}{\gamma} : H(\gamma) &\longrightarrow \overset{*}{H}(\gamma) \\ \Phi &\longmapsto \overset{*}{\gamma}(\Phi) = \overset{*}{\Phi} := \text{the extension of } \Phi \text{ to } \overset{*}{G}(\gamma) \end{aligned}$$

is an isomorphism, it immediately results that the function

$$\begin{aligned} \overset{*}{\gamma}(\overset{\circ}{\gamma})^{-1} : \overset{\circ}{H}(\gamma) &\longrightarrow \overset{*}{H}(\gamma) \\ \overset{\circ}{\Phi} &\longmapsto \overset{*}{\gamma}\Big((\overset{\circ}{\gamma})^{-1}(\overset{\circ}{\Phi})\Big) = \overset{*}{\Phi} \end{aligned}$$

is an isomorphism such that

$$\left(\left(\overset{*}{\gamma}(\overset{\circ}{\gamma})^{-1}\right)(\overset{\circ}{\Phi})\right)(\overset{\circ}{g}) = \overset{*}{\Phi}(\overset{\circ}{g}) = \overset{\circ}{\Phi}(\overset{\circ}{g}) \quad \text{for every} \quad \overset{\circ}{g} \in \overset{\circ}{G}(\gamma) \subseteq \overset{*}{G}(\gamma),$$

that is, $(\overset{*}{\gamma}(\overset{\circ}{\gamma})^{-1})(\overset{\circ}{\Phi})$ is an extension of $\overset{\circ}{\Phi}$ to $\overset{*}{G}(\gamma)$.

In short, $\overset{*}{H}(\gamma)$ is a semigroup of endomorphisms on $\overset{*}{G}(\gamma)$, isomorphic to the semigroup $\overset{\circ}{H}(\gamma)$, being $\overset{*}{\gamma}(\overset{\circ}{\gamma})^{-1} : \overset{\circ}{H}(\gamma) \longrightarrow \overset{*}{H}(\gamma)$ an isomorphism that to each $\overset{\circ}{\Phi} \in \overset{\circ}{H}(\gamma)$ associates (in $\overset{*}{H}(\gamma)$) a prolongation of $\overset{\circ}{\Phi}$ to $\overset{*}{G}(\gamma)$. Again, according to Definition 1.6(c), this means that $\overset{*}{H}(\gamma)$ is a prolongation of $\overset{\circ}{H}(\gamma)$ to $\overset{*}{G}(\gamma)$, as stated in (g).

**(h)** For each $\gamma \in \Gamma(I)$, the ordered pair

$$\overset{*}{\mathbb{G}}(\gamma) := \left(\overset{*}{G}(\gamma), \overset{*}{H}(\gamma)\right)$$

is a $S$-group. Evenmore, $\overset{*}{\mathbb{G}}(\gamma)$ is an extension of the $S$-groups $\mathbb{G}(\gamma) = (G(\gamma), H(\gamma))$ and $\overset{\circ}{\mathbb{G}}(\gamma) = (\overset{\circ}{G}(\gamma), \overset{\circ}{H}(\gamma))$.

This trivially follows from (e), (f) and (g) above and the definitions of $S$-group and $S$-group extension.

## 5.26 The $S$-Space $\overset{*}{\mathscr{G}}(I) = (\overset{*}{\mathbb{G}}(\Gamma(I)), \overset{*}{\imath}(\Gamma^2(I)), \overset{*}{\Theta}(\Delta(I)))$

We start with the definition of a new term: $\overset{*}{\mathbb{G}}(\Gamma(I))$.

**Definition.** $\overset{*}{\mathbb{G}}(\Gamma(I))$ is the family indexed by $\Gamma(I)$ and defined by:

$$\overset{*}{\mathbb{G}}\left(\Gamma(I)\right) := \left\{\overset{*}{\mathbb{G}}(\gamma) = \left(\overset{*}{G}(\gamma), \overset{*}{H}(\gamma)\right)\right\}_{\gamma\in\Gamma(I)}$$

with $\overset{*}{\mathbb{G}}(\gamma)$ as defined in 5.25.

It is worth remembering here the definition of the family $\overset{*}{\imath}(\Gamma^2(I))$, formulated right after Axiom 5, at page 254, namely:

$$\overset{*}{\imath}\left(\Gamma^2(I)\right) := \left\{\overset{*}{\imath}_{(\gamma',\gamma)}\right\}_{(\gamma',\gamma)\in\Gamma^2(I)}$$

being

$$\begin{aligned} \overset{*}{\imath}_{(\gamma',\gamma)} : \overset{*}{H}(\gamma) &\longrightarrow \overset{*}{H}(\gamma') \\ \overset{*}{\Phi} &\longmapsto \overset{*}{\imath}_{(\gamma',\gamma)}(\overset{*}{\Phi}) := \overset{*}{\gamma}{}'\left(i_{(\gamma',\gamma)}\left((\overset{*}{\gamma})^{-1}(\overset{*}{\Phi})\right)\right) \end{aligned}$$

for each $(\gamma', \gamma) \in \Gamma^2(I)$.

For these families we highlight that:

**(a)** $\overset{*}{\mathbb{G}}(\Gamma(I))$ is a family of $S$-groups that is an extension of the family of $S$-groups $\overset{\circ}{\mathbb{G}}(\Gamma(I))$.

This directly results from 5.25(h) and the definitions of family of $S$-groups (Definition 3.7) and extension of a family of $S$-groups (Definition 3.8).

**(b)** $\overset{*}{\imath}(\Gamma^2(I))$ is the extension of the bonding $\overset{\circ}{\imath}(\Gamma^2(I))$ (of the $S$-space $\overset{\circ}{\mathscr{G}}(I)$) to the family $\overset{*}{\mathbb{G}}(\Gamma(I))$.

In fact, from 5.23(j) we know that the function

$$\begin{aligned} \overset{\circ}{\gamma} : H(\gamma) &\longrightarrow \overset{\circ}{H}(\gamma) \\ \Phi &\longmapsto \overset{\circ}{\gamma}(\Phi) = \overset{\circ}{\Phi} := \overset{*}{\gamma}(\Phi)|_{\overset{\circ}{G}(\gamma)} \end{aligned}$$

is an isomorphism from the semigroup $H(\gamma)$ onto the semigroup $\overset{\circ}{H}(\gamma)$. By Axiom 4((4-3)),

$$\begin{aligned} \overset{*}{\gamma} : H(\gamma) &\longrightarrow \overset{*}{H}(\gamma) \\ \Phi &\longmapsto \overset{*}{\gamma}(\Phi) = \overset{*}{\Phi} := \text{the extension of } \Phi \text{ to } \overset{*}{G}(\gamma) \end{aligned}$$

is an isomorphism from the semigroup $H(\gamma)$ onto the semigroup $\overset{*}{H}(\gamma)$. Defining now, for each $\gamma \in \Gamma(I)$, the function $\overline{\gamma}$ by

$$\overline{\gamma} := \overset{*}{\gamma}(\overset{\circ}{\gamma})^{-1},$$

that is,

$$\begin{aligned} \overline{\gamma} : \overset{\circ}{H}(\gamma) &\longrightarrow \overset{*}{H}(\gamma) \\ \overset{\circ}{\Phi} &\longmapsto \overline{\gamma}(\overset{\circ}{\Phi}) := \overset{*}{\gamma}\Big((\overset{\circ}{\gamma})^{-1}(\overset{\circ}{\Phi})\Big) = \overset{*}{\Phi}, \end{aligned}$$

it immediately results that $\overline{\gamma}$ is an isomorphism from the semigroup $\overset{\circ}{H}(\gamma)$ onto the semigroup $\overset{*}{H}(\gamma)$, such that, for each $\overset{\circ}{\Phi} \in \overset{\circ}{H}(\gamma)$,

$$\overline{\gamma}(\overset{\circ}{\Phi}) = \text{the prolongation } (\in \overset{*}{H}(\gamma)) \text{ of } \overset{\circ}{\Phi} \text{ to } \overset{*}{G}(\gamma).$$

Thus being, and taking into account the Definition 3.22 of extension of a bonding, we have that the extension of the bonding $\overset{\circ}{\imath}(\Gamma^2(I))$ (of the $S$-space $\overset{\circ}{\mathscr{G}}(I)$) to the family $\overset{*}{\mathbb{G}}(\Gamma(I))$, is the family $\overline{\imath}(\Gamma^2(I))$ defined by:

$$\overline{\imath}\Big(\Gamma^2(I)\Big) := \Big\{\overline{\imath}_{(\gamma',\gamma)}\Big\}_{(\gamma',\gamma)\in\Gamma^2(I)},$$

being $\overline{\imath}_{(\gamma',\gamma)}$ (for each $(\gamma',\gamma) \in \Gamma^2(I)$) the following function,

$$\begin{aligned} \overline{\imath}_{(\gamma',\gamma)} : \overset{*}{H}(\gamma) &\longrightarrow \overset{*}{H}(\gamma') \\ \overset{*}{\Phi} &\longmapsto \overline{\imath}_{(\gamma',\gamma)}(\overset{*}{\Phi}) := \overline{\gamma'}\Bigg(\overset{\circ}{\imath}_{(\gamma',\gamma)}\Big((\overline{\gamma})^{-1}(\overset{*}{\Phi})\Big)\Bigg). \end{aligned}$$

In terms of the definitions above the statement (b) translates to:

$$\overset{*}{\imath}\Big(\Gamma^2(\gamma)\Big) = \overline{\imath}\Big(\Gamma^2(\gamma)\Big).$$

Let then $(\gamma', \gamma) \in \Gamma^2(I)$ be arbitrarily chosen. By the definition of the family $\mathring{\imath}(\Gamma^2(I))$ (in 5.24) we have that $\mathring{\imath}_{(\gamma',\gamma)} : \mathring{H}(\gamma) \longrightarrow \mathring{H}(\gamma')$ is the function given by:

$$\mathring{\imath}_{(\gamma',\gamma)} = \mathring{\gamma}' \; i_{(\gamma',\gamma)}(\mathring{\gamma})^{-1}.$$

Thus, for $\overset{*}{\Phi} \in \overset{*}{H}(\gamma)$ arbitrarily fixed,

$$\overline{\imath}_{(\gamma',\gamma)}(\overset{*}{\Phi}) = \overline{\gamma'}\Bigg(\mathring{\imath}_{(\gamma',\gamma)}\Big((\overline{\gamma})^{-1}(\overset{*}{\Phi})\Big)\Bigg) = \overline{\gamma'}\Bigg(\mathring{\gamma}'\Bigg(i_{(\gamma',\gamma)}\Bigg((\mathring{\gamma})^{-1}\Big((\overline{\gamma})^{-1}(\overset{*}{\Phi})\Big)\Bigg)\Bigg)\Bigg)$$

and, hence, taking into account that $\overline{\gamma'} = \overset{*}{\gamma}'(\mathring{\gamma}')^{-1}$ and $(\overline{\gamma})^{-1} = \mathring{\gamma}(\overset{*}{\gamma})^{-1}$, we obtain that

$$\overline{\imath}_{(\gamma',\gamma)}(\overset{*}{\Phi}) = \overset{*}{\gamma}'\Bigg(i_{(\gamma',\gamma)}\Big((\overset{*}{\gamma})^{-1}(\overset{*}{\Phi})\Big)\Bigg) = \overset{*}{\imath}_{(\gamma',\gamma)}(\overset{*}{\Phi}).$$

Since $\overset{*}{\Phi} \in \overset{*}{H}(\gamma)$ is arbitrary, comes that

$$\overline{\imath}_{(\gamma',\gamma)} = \overset{*}{\imath}_{(\gamma',\gamma)}$$

from where one concludes, taking into account that $(\gamma', \gamma) \in \Gamma^2(I)$ is arbitrary, that

$$\overline{\imath}\Big(\Gamma^2(I)\Big) = \overset{*}{\imath}\Big(\Gamma^2(I)\Big).$$

**(c)** The ordered pair

$$\Bigg(\overset{*}{\mathbb{G}}\Big(\Gamma(I)\Big), \overset{*}{\imath}\Big(\Gamma^2(I)\Big)\Bigg)$$

is a bonded family that is an extension of the bonded family $(\mathring{\mathbb{G}}(\Gamma(I)), \mathring{\imath}(\Gamma^2(I)))$ (of the $S$-space $\mathring{\mathscr{G}}(I)$).

This directly follows from (a) and (b) above and the definition of extension of a bonded family (Definition 3.22).

Regarding the family of restrictions of the $\overline{\mathscr{G}}(I)$-distributions, that is, the family

$$\overset{*}{\Theta}\Big(\Delta(I)\Big) = \Big\{\overset{*}{\Theta}_{(\gamma',\gamma)}\Big\}_{(\gamma',\gamma)\in\Delta(I)}$$

defined in 5.24 (at page 267), we highlight the following properties:

**(d)** $\overset{*}{\Theta}(\Delta(I))$ is a restriction for the bonded family $(\overset{*}{\mathbb{G}}(\Gamma(I)), \overset{*}{\imath}(\Gamma^2(I)))$.

In fact, as we saw in 5.25(e), $\overset{*}{G}(\gamma)$, for each $\gamma \in \Gamma(I)$, is a (abelian) group and, hence, taking Axiom 6((6-1)) into account,

**(d-1)** for every $(\gamma', \gamma) \in \Delta(I)$, $\overset{*}{\Theta}_{(\gamma',\gamma)} \in \overset{*}{\Theta}(\Delta(I))$ is a homomorphism.

Furthermore, yet from Axiom 6((6-3) and (6-4)), the elements of $\overset{*}{\Theta}(\Delta(I))$ are such that:

**(d-2)** $\overset{*}{\Theta}_{(\gamma'',\gamma')}(\overset{*}{\Theta}_{(\gamma',\gamma)}(\overset{*}{g})) = \overset{*}{\Theta}_{(\gamma'',\gamma)}(\overset{*}{g})$ for every $\gamma, \gamma', \gamma'' \in \Gamma(I)$ such that $\gamma'' \subseteq \gamma' \subseteq \gamma$ and every $\overset{*}{g} \in \overset{*}{G}(\gamma)$;

**(d-3)** $\overset{*}{\Theta}_{(\gamma',\gamma)}(\overset{*}{\Phi}(\overset{*}{g})) = \overset{*}{\imath}_{(\gamma',\gamma)}(\overset{*}{\Phi})(\overset{*}{\Theta}_{(\gamma',\gamma)}(\overset{*}{g}))$ for every $\overset{*}{\Phi} \in \overset{*}{H}(\gamma)$, $\overset{*}{g} \in \overset{*}{G}(\gamma)$ and $(\gamma', \gamma) \in \Delta(I)$.

Reporting ourselves now to the definition of restriction for a bonded family, Definition 3.13, we see that (d-1), (d-2) and (d-3) above are, exact and respectively, the conditions (a), (b) and (c) of the referred definition. Regarding the fourth and last requirement, (d), it is, at the present case, trivially attended since the elements $\overset{*}{\Phi} \in \overset{*}{H}(\gamma)$ are endomorphisms.

**(e)** $\overset{*}{\Theta}(\Delta(I))$ is a prolongation of the family $\overset{\circ}{\Theta}(\Delta(I))$ (of the $S$-space $\overset{\circ}{\mathscr{G}}(I)$) to the bonded family $(\overset{*}{\mathbb{G}}(\Gamma(I)), \overset{*}{\imath}(\Gamma^2(I)))$ (in (c)).

In fact, for as we saw in (c), the ordered pair $(\overset{*}{\mathbb{G}}(\Gamma(I)), \overset{*}{\imath}(\Gamma^2(I)))$ is an extension of the bonded family $(\overset{\circ}{\mathbb{G}}(\Gamma(I)), \overset{\circ}{\imath}(\Gamma^2(I)))$ (of the $S$-space $\overset{\circ}{\mathscr{G}}(I)$) and, by (d), $\overset{*}{\Theta}(\Delta(I))$ is a restriction for $(\overset{*}{\mathbb{G}}(\Gamma(I)), \overset{*}{\imath}(\Gamma^2(I)))$. Furthermore, remembering ourselves the definition of the members $\overset{\circ}{\Theta}_{(\gamma',\gamma)}$ of the family $\overset{\circ}{\Theta}(\Delta(I))$, namely,

$$\overset{\circ}{\Theta}_{(\gamma',\gamma)} := \overset{*}{\Theta}_{(\gamma',\gamma)}|_{\overset{\circ}{G}(\gamma)},$$

we get that

$$\overset{*}{\Theta}_{(\gamma',\gamma)}(\overset{\circ}{g}) = \overset{\circ}{\Theta}_{(\gamma',\gamma)}(\overset{\circ}{g})$$

for any $(\gamma', \gamma) \in \Delta(I)$ and every $\overset{\circ}{g} \in \overset{\circ}{G}(\gamma)$. Therefore, taking Definition 3.23(a) into account, the results above ensure the truthiness of (e).

**(f)** The triplet $\overset{*}{\mathscr{G}}(I)$ defined below is a $S$-space that is a locally closed and coherent extension of the $S$-space $\overset{\circ}{\mathscr{G}}(I)$.

$$\overset{*}{\mathscr{G}}(I) := \left(\overset{*}{\mathbb{G}}\Big(\Gamma(I)\Big), \overset{*}{\imath}\Big(\Gamma^2(I)\Big), \overset{*}{\Theta}\Big(\Delta(I)\Big)\right).$$

It results from (c), (d) and (e) above, taking into account the Definitions 3.15 (of $S$-space) and 3.23(b) (of $S$-space extension), that $\overset{*}{\mathscr{G}}(I)$ is a $S$-space and also an extension of the $S$-space $\overset{\circ}{\mathscr{G}}(I)$. Finally, resorting to the definitions of locally closed extension (Definition 3.23(c-2)) and of coherent extension (Definition 3.23(c-3)), we see that Axiom 7 and Axiom 8 establish, respectively, that the $S$-space $\overset{*}{\mathscr{G}}(I)$ is a locally closed and coherent extension of the $S$-space $\overset{\circ}{\mathscr{G}}(I)$.

## 5.27 The Categoricity of the $\overline{\mathscr{G}}(I)$-Distributions Axiomatic

As we saw in item 5.22, for each $S$-space

$$\check{\mathscr{G}}(I) = \Bigg(\check{\mathbb{G}}\Big(\Gamma(I)\Big) = \Big\{\check{\mathbb{G}}(\gamma) = \Big(\check{G}(\gamma), \check{H}(\gamma)\Big)\Big\}_{\gamma\in\Gamma(I)}, \check{\imath}\Big(\Gamma^2(I)\Big), \check{\Theta}\Big(\Delta(I)\Big)\Bigg)$$

isomorphic to $\overline{\mathscr{G}}(I)$, this last one defined in the 2nd TESS, the table below defines a model of the $\overline{\mathscr{G}}(I)$-distributions axiomatic.

Now, given a model of the referred axiomatic, the logical consequences obtained in the items 5.23 to 5.26 show us that with the elements that constitute the model we obtain two $S$-spaces, namely, $\overset{\circ}{\mathscr{G}}(I)$ as in 5.24(k) and $\overset{*}{\mathscr{G}}(I)$ as defined in 5.26(f): the first, $\overset{\circ}{\mathscr{G}}(I)$, is a strict and closed extension of $\mathscr{G}(I)$ (the underlying $S$-space of the axiomatic), and the second, $\overset{*}{\mathscr{G}}(I)$, is a locally closed and coherent extension of $\overset{\circ}{\mathscr{G}}(I)$. But, taking Proposition 4.34 into account, $\overset{*}{\mathscr{G}}(I)$ is isomorphic to $\overline{\mathscr{G}}(I)$ ($\overset{*}{\mathscr{G}}(I) \simeq \overline{\mathscr{G}}(I)$) and, thus, we can conclude that **every** model of the $\overline{\mathscr{G}}(I)$-distributions axiomatic is obtained through the table below for some $S$-space $\check{\mathscr{G}}(I)$ isomorphic to $\overline{\mathscr{G}}(I)$. Hence, any two models of this axiomatic are isomorphic to each other and isomorphic to $\overline{\mathscr{G}}(I)$; it is a categoric axiomatic.

| **Primitive terms** | **Interpretation** |
|---|---|
| $\overline{\mathscr{G}}(I)$-distribution | $\check{\mathscr{G}}(I)$-individual |
| $\overline{\mathscr{G}}(I)$-distribution domain | $\check{\mathscr{G}}(I)$-individual domain |
| $\overline{\mathscr{G}}(I)$-distributions addition | $\check{\mathscr{G}}(I)$-individuals addition |
| $\overline{\mathscr{G}}(I)$-distribution derivative | endomorphism $\check{\Phi}$ of the semigroups $\check{H}(\gamma)$ of the $S$-space $\check{\mathscr{G}}(I)$ |
| $\overline{\mathscr{G}}(I)$-distribution restriction | homomorphism $\check{\Theta}_{(\gamma',\gamma)}$ of the family $\check{\Theta}(\Delta(I))$ of the $S$-space $\check{\mathscr{G}}(I)$ |

## 5.28 $\overline{\mathscr{G}}(I)$-Distributions Axioms: Simplified Version

Such as for the $\tilde{\mathscr{G}}(I)$-distributions (see 5.17), the $\overline{\mathscr{G}}(I)$-distributions can also be defined through an axiomatic where the term "addition" does not figure as a primitive term. In what follows, we formulate the axioms of this simplified axiomatic version, leaving for the next items the proof of its equivalence with the axiomatic given in 5.21.

The primitive terms and precedent theories of this axiomatic version are:

- **Primitive Terms:** $\overline{\mathscr{G}}(I)$-distribution, domain, derivative, and restriction of $\overline{\mathscr{G}}(I)$-distribution;
- **Precedent Theories:** Classical Logic, Set Theory, and $S$-Spaces Theory;
- **Axioms**: The ones stated ahead, regarding an abelian, surjective, with identity, and coherent $S$-space,

$$\mathscr{G}(I) = \Bigg(\mathbb{G}\Big(\Gamma(I)\Big) = \Big\{\mathbb{G}(\gamma) = \Big(G(\gamma), H(\gamma)\Big)\Big\}_{\gamma\in\Gamma(I)}, \quad i\Big(\Gamma^2(I)\Big) = \Big\{i_{(\gamma',\gamma)}\Big\}_{(\gamma',\gamma)\in\Gamma^2(I)}, \quad \Theta\Big(\Delta(I)\Big) = \Big\{\Theta_{(\gamma',\gamma)}\Big\}_{(\gamma',\gamma)\in\Delta(I)}\Bigg),$$

arbitrarily fixed, with

$$G(\gamma) \cap G(\gamma') = \varnothing$$

for $\gamma, \gamma' \in \Gamma(I)$ such that $\gamma \neq \gamma'$.

**Axiom 1** Every $\mathscr{G}(I)$-individual is a $\overline{\mathscr{G}}(I)$-distribution.

**Axiom 2** To each $\overline{\mathscr{G}}(I)$-distribution, $\overset{*}{g}$, corresponds a single $\gamma \in \Gamma(I)$, denominated the domain of $\overset{*}{g}$, in such a way that if $\overset{*}{g}$ is a $\mathscr{G}(I)$-individual, the domain $\gamma$ is the domain of the $\mathscr{G}(I)$-individual $\overset{*}{g}$.

**Notation.** We will denote by $\overset{*}{G}(\gamma)$ the class of $\overline{\mathscr{G}}(I)$-distributions with domain $\gamma$. It results from Axiom 1 and Axiom 2 that the class $G(\gamma)$ of $\mathscr{G}(I)$-individuals with domain $\gamma$ is a subclass of $\overset{*}{G}(\gamma)$: $G(\gamma) \subseteq \overset{*}{G}(\gamma)$.

**Axiom 3** The derivatives of the $\overline{\mathscr{G}}(I)$-distributions with domain $\gamma$, whose class is denoted by $\overset{*}{H}(\gamma)$, are functions with domain and codomain equal to $\overset{*}{G}(\gamma)$, such that:

**(3-1)** each derivative $\overset{*}{\Phi} \in \overset{*}{H}(\gamma)$ is an extension to $\overset{*}{G}(\gamma)$ of a single $\Phi \in H(\gamma)$, that is, $\overset{*}{\Phi}(g) = \Phi(g)$ for every $g \in G(\gamma)_\Phi \subseteq \overset{*}{G}(\gamma)$; conversely, for each $\Phi \in H(\gamma)$ there exists a single derivative $\overset{*}{\Phi} \in \overset{*}{H}(\gamma)$ such that $\Phi(g) = \overset{*}{\Phi}(g)$ for every $g \in G(\gamma)_\Phi \subseteq \overset{*}{G}(\gamma)$;

**(3-2)** the class $\overset{*}{H}(\gamma)$, equipped with the usual operation of composition of functions, is a semigroup isomorphic to the semigroup $H(\gamma)$, with the function

$$\begin{aligned} \overset{*}{\gamma} : H(\gamma) &\longrightarrow \overset{*}{H}(\gamma) \\ \Phi &\longmapsto \overset{*}{\gamma}(\Phi) := \text{the extension of } \Phi \text{ to } \overset{*}{G}(\gamma) \end{aligned}$$

as an isomorphism.

**Notation.** The image of $\Phi \in H(\gamma)$ by the isomorphism $\overset{*}{\gamma} : H(\gamma) \longrightarrow \overset{*}{H}(\gamma)$ will be denoted by $\overset{*}{\Phi}$, that is,

$$\overset{*}{\Phi} := \overset{*}{\gamma}(\Phi).$$

**Axiom 4** For $\gamma \in \Gamma(I)$, $g, h \in G(\gamma)$ and $\Phi \in H(\gamma)$ arbitrarily fixed, one has:

$$\overset{*}{\Phi}(g) = \overset{*}{\Phi}(h) \quad \text{if and only if} \quad g - h \in N(\Phi).$$

**Definition.** For each $(\gamma', \gamma) \in \Gamma^2(I)$ we define $\overset{*}{i}_{(\gamma',\gamma)}$ as the following function:

$$\begin{aligned} \overset{*}{i}_{(\gamma',\gamma)} : \overset{*}{H}(\gamma) &\longrightarrow \overset{*}{H}(\gamma') \\ \overset{*}{\Phi} &\longmapsto \overset{*}{i}_{(\gamma',\gamma)}(\overset{*}{\Phi}) := \overset{*}{\gamma}{}' \left( i_{(\gamma',\gamma)} \left( (\overset{*}{\gamma})^{-1}(\overset{*}{\Phi}) \right) \right). \end{aligned}$$

It results that $\overset{*}{i}_{(\gamma',\gamma)}$ is an isomorphism from the semigroup $\overset{*}{H}(\gamma)$ onto the semigroup $\overset{*}{H}(\gamma')$.

**Axiom 5** The restrictions of $\overline{\mathscr{G}}(I)$-distributions are functions $\overset{*}{\Theta}_{(\gamma',\gamma)} : \overset{*}{G}(\gamma) \longrightarrow \overset{*}{G}(\gamma')$, a single one for each $(\gamma', \gamma) \in \Delta(I)$, that associate to each $\overline{\mathscr{G}}(I)$-distribution with domain $\gamma$, a single $\overline{\mathscr{G}}(I)$-distribution with domain $\gamma' \subseteq \gamma$, in such a way that:

**(5-1)** $\overset{*}{\Theta}_{(\gamma',\gamma)}(g) = \Theta_{(\gamma',\gamma)}(g)$ for every $g \in G(\gamma) \subseteq \overset{*}{G}(\gamma)$;

**(5-2)** $\overset{*}{\Theta}_{(\gamma'',\gamma')}\left(\overset{*}{\Theta}_{(\gamma',\gamma)}(\overset{*}{g})\right) = \overset{*}{\Theta}_{(\gamma'',\gamma)}(\overset{*}{g})$ for every $\overset{*}{g} \in \overset{*}{G}(\gamma)$ and $\gamma, \gamma', \gamma'' \in \Gamma(I)$ such that $\gamma'' \subseteq \gamma' \subseteq \gamma$;

**(5-3)** $\overset{*}{\Theta}_{(\gamma',\gamma)}\left(\overset{*}{\Phi}(\overset{*}{g})\right) = \overset{*}{i}_{(\gamma',\gamma)}(\overset{*}{\Phi})\left(\overset{*}{\Theta}_{(\gamma',\gamma)}(\overset{*}{g})\right)$ for every $\overset{*}{g} \in \overset{*}{G}(\gamma)$ and every $\overset{*}{\Phi} \in \overset{*}{H}(\gamma)$.

**Axiom 6** Let $\gamma \in \Gamma(I)$ be arbitrarily fixed. If $\overset{*}{g} \in \overset{*}{G}(\gamma)$ and $x \in \gamma$, there exists $\gamma' \in \Gamma(\gamma)$ such that $x \in \gamma'$ and

$$\overset{*}{\Theta}_{(\gamma',\gamma)}(\overset{*}{g}) = \overset{*}{\Phi}(g)$$

for some $\Phi \in H(\gamma')$ and $g \in G(\gamma')$.

**Definition.** Let $\gamma \in \Gamma(I)$ and $\overline{\xi}(\gamma) \subseteq \Gamma(\gamma)$ be a cover of $\gamma$. Let also

$$\left\{ \overset{*}{g}_\xi \right\}_{\xi \in \overline{\xi}(\gamma)}$$

be a family indexed by $\overline{\xi}(\gamma)$, whose members, $\overset{*}{g}_\xi$, are $\overline{\mathscr{G}}(I)$-distributions with domain $\xi$, that is, $\overset{*}{g}_\xi \in \overset{*}{G}(\xi)$ for every $\xi \in \overline{\xi}(\gamma)$. We sill say that the family $\{\overset{*}{g}_\xi\}_{\xi \in \overline{\xi}(\gamma)}$ is **coherent** if and only if

$$\overset{*}{\Theta}_{(\xi \cap \xi', \xi)}(\overset{*}{g}_\xi) = \overset{*}{\Theta}_{(\xi \cap \xi', \xi')}(\overset{*}{g}_{\xi'})$$

for every $\xi, \xi' \in \overline{\xi}(\gamma)$ such that $\xi \cap \xi' \neq \emptyset$.

**Axiom 7** Let $\gamma \in \Gamma(I)$ be arbitrarily fixed. If $\{\overset{*}{g}_\xi\}_{\xi \in \overline{\xi}(\gamma)}$ is a coherent family of $\overline{\mathscr{G}}(I)$-distributions, then, there exists a single $\overline{\mathscr{G}}(I)$-distribution $\overset{*}{g} \in \overset{*}{G}(\gamma)$ such that

$$\overset{*}{\Theta}_{(\xi,\gamma)}(\overset{*}{g}) = \overset{*}{g}_\xi$$

for every $\xi \in \overline{\xi}(\gamma)$.

A subclass of $\overset{*}{G}(\gamma)$, that will prove itself important in our considerations regarding the axioms above, is defined below.

**Definition.** $\overset{\circ}{G}(\gamma)$ is the subclass of the $\overline{\mathscr{G}}(I)$-distributions with domain $\gamma$ defined by:

$$\overset{\circ}{G}(\gamma) := \left\{ \overset{*}{g} \in \overset{*}{G}(\gamma) \quad : \quad \overset{*}{g} = \overset{*}{\Phi}(g) \quad \text{for some} \quad \Phi \in H(\gamma) \quad \text{and} \quad g \in G(\gamma) \right\}.$$

The elements of $\overset{\circ}{G}(\gamma)$ are denominated $\overset{\circ}{\mathscr{G}}(I)$-distributions with domain $\gamma$.

Before we start an analysis of the axioms above, in the sense of obtaining some of its logical consequences, done in the next items (5.29 and 5.30), it is opportune to observe that, although the primitive terms of this axiomatic — henceforth "axiomatic B" — have the same denomination of the ones in the axiomatic formulated in 5.21 — henceforth "axiomatic A" — except for the "addition" term that does not figure among the axioms above, once they are different axiomatics, yet with a considerable number of axioms present in both, these terms with the same denomination in both axiomatics may not share the same properties.

The same observation is pertinent regarding the terms introduced into both axiomatics through definitions formally identical and for which one attributes the same symbol or name; for instance, it is part of axiomatic A a class introduced by a definition (at page 258) formally identical to the definition above for $\overset{\circ}{G}(\gamma)$, that, in the first axiomatic, is also denoted by $\overset{\circ}{G}(\gamma)$. The same occurs with the definitions of $\overset{*}{i}_{(\gamma',\gamma)}$ (at page 254 for axiomatic A and page 283 for axiomatic B) and coherent family (axiomatic A: page 254; axiomatic B: page 284). One must keep in mind that theses objects, even with identical formal definitions and represented by the same symbols, they may differ.

Clearly, the properties that can be proved for these primitive or defined terms that appear with the same denomination and symbology in both axiomatics, resorting only to the axioms present in both axiomatics, are, in fact, its properties either considered as terms of one or the other between these axiomatics; hence, in what follows, specially in the item 5.29, whenever one such property is needed, we will simply present it and refer the reader to the proof provided in the context of axiomatic A.

However, we will see that these axiomatics, A and B, are equivalent and, hence, we will conclude that, the above referred terms and symbols, are equal not only in its form, but also in its content; but, until we prove this equivalence, it is worth keeping the observations above in mind.

## 5.29 Elementary Consequences of the Axioms

In the developments presented here and in the next item, the reference to an axiom or an axiomatic without any other information must be interpreted as regarding axiomatic B, that is, the one constituted by the axioms presented in 5.28.

Are relevant for the establishment of the equivalence between axiomatics A and B, as we will see in item 5.31, the following logical consequences of the axioms.

**(a)** $G(\gamma) \subseteq \overset{\circ}{G}(\gamma)$

Result proved in 5.23(a).[26]

**(b)** If $\overset{\circ}{g} \in \overset{\circ}{G}(\gamma)$ and $\overset{*}{\Phi} \in \overset{*}{H}(\gamma)$, then, $\overset{*}{\Phi}(\overset{\circ}{g}) \in \overset{\circ}{G}(\gamma)$.

Result proved in 5.23(g).

In what follows, we introduce into our axiomatic B the set $\overset{\circ}{H}(\gamma)$ through a definition formally identical to the one in axiomatic A, of the set there denoted by the same symbol.

**Definition.** Let $\gamma \in \Gamma(I)$ be arbitrarily fixed. $\overset{\circ}{H}(\gamma)$ is the set defined by

$$\overset{\circ}{H}(\gamma) := \left\{ \overset{*}{\Phi}|_{\overset{\circ}{G}(\gamma)} : \overset{*}{\Phi} \in \overset{*}{H}(\gamma) \right\},$$

where

$$\begin{aligned} \overset{*}{\Phi}|_{\overset{\circ}{G}(\gamma)} : \overset{\circ}{G}(\gamma) &\longrightarrow \overset{*}{G}(\gamma) \\ \overset{\circ}{g} &\longmapsto \overset{*}{\Phi}|_{\overset{\circ}{G}(\gamma)}(\overset{\circ}{g}) := \overset{*}{\Phi}(\overset{\circ}{g}). \end{aligned}$$

(Observe that (b) ensures that the functions belonging to $\overset{\circ}{H}(\gamma)$ assume values in $\overset{\circ}{G}(\gamma)$, that is, they are function from $\overset{\circ}{G}(\gamma)$ into $\overset{\circ}{G}(\gamma)$).

Analogously, through a definition given ahead, formally identical to the one provided for the homonymous element, $\overset{\circ}{\gamma}$, in axiomatic A, we introduce into our axiomatic the function $\overset{\circ}{\gamma}$.

[26] In order for the reader to certify itself that (a) is a logical consequence of the axiomatic B, one must verify that the proof provided in 5.23(a), even elaborated in the context of axiomatic A, it only uses axioms of A that are also present in axiomatic B and, thus, it is also a proof of (a) in the axiomatic B. In an analogous manner, one must proceed in similar cases.

**Definition.** $\overset{\circ}{\gamma}$ is the function defined by:

$$\begin{aligned} \overset{\circ}{\gamma} : H(\gamma) &\longrightarrow \overset{\circ}{H}(\gamma) \\ \Phi &\longmapsto \overset{\circ}{\gamma}(\Phi) := \overset{*}{\gamma}(\Phi)|_{\overset{\circ}{G}(\gamma)} \end{aligned}$$

(according to Axiom 3, $\overset{*}{\gamma}(\Phi) = \overset{*}{\Phi}$ is the extension of $\Phi$ to $\overset{*}{G}(\gamma)$).

Here, also, such as we did in axiomatic A, we introduce the following notation:

**Notation.** $\overset{\circ}{\Phi} := \overset{\circ}{\gamma}(\Phi) = \overset{*}{\Phi}|_{\overset{\circ}{G}(\gamma)}$.

Regarding the class $\overset{\circ}{H}(\gamma)$ and the function $\overset{\circ}{\gamma} : H(\gamma) \longrightarrow \overset{\circ}{H}(\gamma)$ above defined we have that:

**(c)** $\overset{\circ}{H}(\gamma)$, with the usual operation of composition of functions, is a semigroup isomorphic to the semigroup $H(\gamma)$ and the function $\overset{\circ}{\gamma} : H(\gamma) \longrightarrow \overset{\circ}{H}(\gamma)$ defined above is an isomorphism.

Result proved in 5.23(i) and 5.23(j).

Other three logical consequences of our axiomatic are:

**(d)** For each $\overset{\circ}{g} \in \overset{\circ}{G}(\gamma)$, there exist $\Phi \in H(\gamma)$ and $g \in G(\gamma)$ such that

$$\overset{\circ}{g} = \overset{\circ}{\Phi}(g).$$

In fact, if $\overset{\circ}{g} \in \overset{\circ}{G}(\gamma)$, then, from the definition of $\overset{\circ}{G}(\gamma)$ we can write that

$$\overset{\circ}{g} = \overset{*}{\Phi}(g)$$

for some $\Phi \in H(\gamma)$ and $g \in G(\gamma)$. Since $\overset{\circ}{\Phi} = \overset{*}{\Phi}|_{\overset{\circ}{G}(\gamma)}$ and, by (a), $G(\gamma) \subseteq \overset{\circ}{G}(\gamma)$, we have $\overset{*}{\Phi}(g) = \overset{\circ}{\Phi}(g)$ and, hence, $\overset{\circ}{g} = \overset{\circ}{\Phi}(g)$.

**(e)** $\overset{\circ}{\gamma} : H(\gamma) \longrightarrow \overset{\circ}{H}(\gamma)$ is an isomorphism such that, for each $\Phi \in H(\gamma)$ one has that:

$$\Big(\overset{\circ}{\gamma}(\Phi)\Big)(g) = \overset{\circ}{\Phi}(g) = \Phi(g)$$

for every $g \in G(\gamma)_\Phi$

In fact, as we know,

$$\overset{\circ}{\Phi} = \overset{\circ}{\gamma}(\Phi) = \overset{*}{\gamma}(\Phi)|_{\overset{\circ}{G}(\gamma)} = \overset{*}{\Phi}|_{\overset{\circ}{G}(\gamma)}.$$

Hence, if $g \in G(\gamma)_\Phi \subseteq G(\gamma) \subseteq \overset{\circ}{G}(\gamma)$, we have

$$\overset{\circ}{\Phi}(g) = \overset{*}{\gamma}(\Phi)|_{\overset{\circ}{G}(\gamma)}(g) = \overset{*}{\Phi}|_{\overset{\circ}{G}(\gamma)}(g) = \overset{*}{\Phi}(g).$$

But, taking Axiom 3((3-1)) into account,

$$\overset{*}{\Phi}(g) = \Phi(g) \quad \text{for} \quad g \in G(\gamma)_{\Phi}.$$

Therefore,

$$\overset{\circ}{\Phi}(g) = \Phi(g).$$

**(f)** For $g, h \in G(\gamma)$ and $\Phi \in H(\gamma)$ arbitrarily fixed, one has that:

$$\overset{\circ}{\Phi}(g) = \overset{\circ}{\Phi}(h) \quad \text{if and only if} \quad g - h \in N(\Phi).$$

An immediate consequence of Axiom 4 and the definition of $\overset{\circ}{\Phi}$.

Another logical consequence of the axioms can now be obtained resorting to Proposition 2.24 and the results (a), (c) and (e), (d), and (f) above. In fact, for each $\gamma \in \Gamma(I)$, $\mathbb{G}(\gamma) = (G(\gamma), H(\gamma))$ is an abelian, surjective, and with identity $S$-group, and then can assume the role of the $S$-group $\mathbb{G} = (G, H)$ of the referred proposition; furthermore, the results (a), (c) and (e), (d), and (f) above referred correspond, for $\overset{\circ}{G}(\gamma)$, $\overset{\circ}{H}(\gamma)$ and $\overset{\circ}{\gamma}$, respectively, to the conditions (a) to (d) of the Proposition 2.24 required for $\widehat{G}$, $\widehat{H}$ and $\widehat{\ }$. It then results:

**(g)** For each $\gamma \in \Gamma(I)$, there exists a single binary operation in $\overset{\circ}{G}(\gamma)$, let us say $+ : \overset{\circ}{G}(\gamma) \times \overset{\circ}{G}(\gamma) \longrightarrow \overset{\circ}{G}(\gamma)$, which we will denominate "addition", that makes of $\overset{\circ}{G}(\gamma)$ an abelian group that has the group $G(\gamma)$ as a subgroup and of $\overset{\circ}{\Phi} : \overset{\circ}{G}(\gamma) \longrightarrow \overset{\circ}{G}(\gamma)$, for each $\overset{\circ}{\Phi} \in \overset{\circ}{H}(\gamma)$, an endomorphism on the group $\overset{\circ}{G}(\gamma)$.

In what follows, any mention to addition in $\overset{\circ}{G}(\gamma)$ or to $\overset{\circ}{G}(\gamma)$ as a group, refers to the (unique) addition described in (g).

Several other consequences can now be obtained taking into account the addition characterized in (g). Hence, we have:

**(h)** $\overset{\circ}{H}(\gamma)$ is a prolongation of $H(\gamma)$ to $\overset{\circ}{G}(\gamma)$.

In fact, by (c) and (e) we have that $\overset{\circ}{H}(\gamma)$ is a semigroup isomorphic to the semigroup $H(\gamma)$ and

$$\begin{aligned} \overset{\circ}{\gamma} : H(\gamma) &\longrightarrow \overset{\circ}{H}(\gamma) \\ \Phi &\longmapsto \overset{\circ}{\gamma}(\Phi) := \overset{*}{\gamma}(\Phi)|_{\overset{\circ}{G}(\gamma)} \end{aligned}$$

is an isomorphism such that

$$\overset{\circ}{\gamma}(\Phi)(g) = \overset{\circ}{\Phi}(g) = \Phi(g)$$

for every $g \in G(\gamma)_{\Phi}$.

Furthermore, by (g), each $\mathring{\Phi} \in \mathring{H}(\gamma)$ is an endomorphism on the group $\mathring{G}(\gamma)$. Thus, keeping Definition 1.6(c) in mind, we can conclude that $\mathring{H}(\gamma)$ is a prolongation of $H(\gamma)$ to $\mathring{G}(\gamma)$.

**(i)** For each $\gamma \in \Gamma(I)$, the ordered pair

$$\mathring{\mathbb{G}}(\gamma) \coloneqq \Big(\mathring{G}(\gamma), \mathring{H}(\gamma)\Big)$$

is a $S$-group and, also, an extension of the $S$-group $\mathbb{G}(\gamma) = (G(\gamma), H(\gamma))$.

Direct consequence of (g), (h) and the definitions of $S$-group and $S$-group extension.

**(j)** The $S$-group $\mathring{\mathbb{G}}(\gamma)$ defined in (i) is a strict and closed extension of the $S$-group $\mathbb{G}(\gamma)$.

It directly follows from (d) that $\mathring{\mathbb{G}}(\gamma)$ is a closed extension of the $S$-group $\mathbb{G}(\gamma)$. Let us them prove that $\mathring{\mathbb{G}}(\gamma)$ also is a strict extension of $\mathbb{G}(\gamma)$, i.e., that: for each $\Phi \in H(\gamma)$, if $g \in G(\gamma)$ and $g \notin G(\gamma)_\Phi$, then $\mathring{\Phi}(g) \notin G(\gamma)$.

In order to do so, let $\Phi \in H(\gamma)$ and $g \in G(\gamma)$ be such that $g \notin G(\gamma)_\Phi$. Suppose that

$$h \coloneqq \mathring{\Phi}(g) \in G(\gamma).$$

Now, $h = \Phi(h_1)$ for some $h_1 \in G(\gamma)_\Phi$, once that $\Phi : G(\gamma)_\Phi \longrightarrow G(\gamma)$ is a surjective homomorphism and, thus,

$$\mathring{\Phi}(g) = \Phi(h_1).$$

However, by (e), $\Phi(h_1) = \mathring{\Phi}(h_1)$. Therefore,

$$\mathring{\Phi}(g) = \mathring{\Phi}(h_1)$$

which, by (f), means that

$$g - h_1 \in N(\Phi) \subseteq G(\gamma)_\Phi.$$

Since $h_1 \in G(\gamma)_\Phi$ and $G(\gamma)_\Phi$ is a group (subgroup of the group $G(\gamma)$), then

$$g = (g - h_1) + h_1 \in G(\gamma)_\Phi,$$

which contradicts our hypothesis of $g \notin G(\gamma)_\Phi$. Therefore, $\mathring{\Phi}(\gamma) \notin G(\gamma)$.

In what follows, we are going to introduce two families into our axiomatic, $\mathring{\mathbb{G}}(\Gamma(I))$ and $\mathring{\imath}(\Gamma^2(I))$, through definitions formally identical to the homonymous ones in axiomatic A ($\mathring{\mathbb{G}}(\Gamma(I))$ and $\mathring{\imath}(\Gamma^2(I))$ defined in 5.24 at page 265).

**Definition.** $\mathring{\mathbb{G}}(\Gamma(I))$ and $\mathring{\imath}(\Gamma^2(I))$ are families indexed, respectively, by $\Gamma(I)$ and $\Gamma^2(I)$, defined as follows:

$$\mathring{\mathbb{G}}\Big(\Gamma(I)\Big) \coloneqq \Big\{\mathring{\mathbb{G}}(\gamma) = \Big(\mathring{G}(\gamma), \mathring{H}(\gamma)\Big)\Big\}_{\gamma \in \Gamma(I)}$$

with $\mathring{\mathbb{G}}(\gamma)$ as in (i) above, and

$$\mathring{i}\Big(\Gamma^2(I)\Big) := \Big\{\mathring{i}_{(\gamma',\gamma)}\Big\}_{(\gamma',\gamma)\in\Gamma^2(I)}$$

where, for each $(\gamma',\gamma) \in \Gamma^2(I)$, $\mathring{i}_{(\gamma',\gamma)}$ is the following function,

$$\mathring{i}_{(\gamma',\gamma)} : \mathring{H}(\gamma) \longrightarrow \mathring{H}(\gamma')$$

$$\mathring{\Phi} \longmapsto \mathring{i}_{(\gamma',\gamma)}(\mathring{\Phi}) := \mathring{\gamma}'\left(i_{(\gamma',\gamma)}\Big((\mathring{\gamma})^{-1}(\mathring{\Phi})\Big)\right).$$

Regarding these families, the properties 5.24(a), 5.24(b) and 5.24(c) proved in item 5.24 for its homonymous (in axiomatic A), also are valid here (in axiomatic B) and by the same arguments (proofs) there presented. Hence, we have that:

**(k)** $\mathring{\mathbb{G}}(\Gamma(I))$ is a family of $S$-groups that is a strict and closed extension of the family of $S$-groups $\mathbb{G}(\Gamma(I))$ (of the $S$-space $\mathscr{G}(I)$).

**(l)** $\mathring{i}(\Gamma^2(I))$ is the extension of the bonding $i(\Gamma^2(I))$ (of the $S$-space $\mathscr{G}(I)$) to the family $\mathring{\mathbb{G}}(\Gamma(I))$.

**(m)** The ordered pair

$$\left(\mathring{\mathbb{G}}\Big(\Gamma(I)\Big), \mathring{i}\Big(\Gamma^2(I)\Big)\right)$$

is a bonded family that is an extension of the bonded family (of the $S$-space $\mathscr{G}(I)$) $(\mathbb{G}(\Gamma(I)), i(\Gamma^2(I)))$.

Next, in the same way, i.e., that of to enrich our axiomatic, whenever possible, with definitions formally identical to those present in axiomatic A, we are going to define two other families, $\overset{*}{\Theta}(\Delta(I))$ and $\mathring{\Theta}(\Delta(I))$, both indexed by the set $\Delta(I)$.

**Definition.** $\overset{*}{\Theta}(\Delta(I))$ and $\mathring{\Theta}(\Delta(I))$ are the families defined as follows:

$$\overset{*}{\Theta}\Big(\Delta(I)\Big) := \Big\{\overset{*}{\Theta}_{(\gamma',\gamma)}\Big\}_{(\gamma',\gamma)\in\Delta(I)}$$

and

$$\mathring{\Theta}\Big(\Delta(I)\Big) := \Big\{\mathring{\Theta}_{(\gamma',\gamma)}\Big\}_{(\gamma',\gamma)\in\Delta(I)},$$

where, for each $(\gamma',\gamma) \in \Delta(I)$,

$$\mathring{\Theta}_{(\gamma',\gamma)} := \overset{*}{\Theta}_{(\gamma',\gamma)}|_{\mathring{G}(\gamma)},$$

that is, $\mathring{\Theta}_{(\gamma',\gamma)}$ is the function

$$\mathring{\Theta}_{(\gamma',\gamma)} : \mathring{G}(\gamma) \longrightarrow \mathring{G}(\gamma')$$

$$\mathring{g} \longmapsto \mathring{\Theta}_{(\gamma',\gamma)}(\mathring{g}) := \overset{*}{\Theta}_{(\gamma',\gamma)}(\mathring{g})$$

(by (n) ahead, one has that $\mathring{\Theta}_{(\gamma',\gamma)}(\mathring{g}) \in \mathring{G}(\gamma')$).

Also for the families $\overset{*}{\Theta}(\Delta(I))$ and $\overset{\circ}{\Theta}(\Delta(I))$ above defined, many of the properties obtained in item 5.24 for its homonymous ones in axiomatic A are, by the same arguments there presented, true here. The reader will find no problem verifying, reviewing the respective proofs provided in 5.24, that the results 5.24(d), (e), (g) and (h) rewritten below as (n), (o), (p) and (q), respectively, are valid in the context of axiomatic B.

**(n)** $\overset{*}{\Theta}_{(\gamma',\gamma)}(\mathring{g}) \in \mathring{G}(\gamma')$ for every $(\gamma',\gamma) \in \Delta(I)$ and $\mathring{g} \in \mathring{G}(\gamma)$.

**(o)** Let $(\gamma',\gamma), (\beta',\beta) \in \Delta(I)$. One has that:

> **(o-1)** $\overset{\circ}{\Theta}_{(\gamma',\gamma)} = \overset{\circ}{\Theta}_{(\beta',\beta)}$ (if and) only if $\gamma' = \beta'$ and $\gamma = \beta$;
>
> **(o-2)** $\overset{*}{\Theta}_{(\gamma',\gamma)} = \overset{*}{\Theta}_{(\beta',\beta)}$ (if and) only if $\gamma' = \beta'$ and $\gamma = \beta$.[27]

**(p)** For every $\gamma, \gamma', \gamma'' \in \Gamma(I)$ such that $\gamma'' \subseteq \gamma' \subseteq \gamma$ and every $\mathring{g} \in \mathring{G}(\gamma)$,

$$\overset{\circ}{\Theta}_{(\gamma'',\gamma')}\left(\overset{\circ}{\Theta}_{(\gamma',\gamma)}(\mathring{g})\right) = \overset{\circ}{\Theta}_{(\gamma'',\gamma)}(\mathring{g}).$$

**(q)** For each $(\gamma',\gamma) \in \Delta(I)$, every $\mathring{\Phi} \in \mathring{H}(\gamma)$ and every $\mathring{g} \in \mathring{G}(\gamma)$,

$$\overset{\circ}{\Theta}_{(\gamma',\gamma)}\left(\mathring{\Phi}(\mathring{g})\right) = \mathring{\imath}_{(\gamma',\gamma)}(\mathring{\Phi})\left(\mathring{\Phi}_{(\gamma',\gamma)}(\mathring{g})\right).$$

In axiomatic A, from Axiom 6((6-1)) and the definition of the functions $\overset{\circ}{\Theta}_{(\gamma',\gamma)} \in \overset{\circ}{\Theta}(\Delta(I))$, it immediately results that the members of the family $\overset{\circ}{\Theta}(\Delta(I))$ are homomorphisms (of group). In axiomatic B, we does not have the referred axiom at our disposal. However, also here (in B), the members of $\overset{\circ}{\Theta}(\Delta(I))$ are homomorphisms, as proved ahead.

**(r)** For each $(\gamma',\gamma) \in \Delta(I)$, $\overset{\circ}{\Theta}_{(\gamma',\gamma)} \in \overset{\circ}{\Theta}(\Delta(I))$ is a homomorphism.

In fact, let $\mathring{g}, \mathring{h} \in \mathring{G}(\gamma)$ be arbitrarily fixed. By (d), there exist $\Phi, \Psi \in H(\gamma)$ and $g, h \in G(\gamma)$ such that

$$\mathring{g} = \mathring{\Phi}(g) \quad \text{and} \quad \mathring{h} = \mathring{\Psi}(h),$$

and, since $\Phi : G(\gamma)_\Phi \longrightarrow G(\gamma)$ and $\Psi : G(\gamma)_\Psi \longrightarrow G(\gamma)$ are surjective homomorphisms, there exist $g_1 \in G(\gamma)_\Psi$ and $h_1 \in G(\gamma)_\Phi$ such that

$$g = \Psi(g_1) \quad \text{and} \quad h = \Phi(h_1),$$

[27] The reader that reexamine the proof of this proposition in the context of axiomatic A, given in 5.24(e), will conclude that it can only be taken as a proof in axiomatic B, if $\overset{*}{G}(\gamma)$ is, also in this axiomatic, an abelian group that has the groups $G(\gamma)$ and $\mathring{G}(\gamma)$ as its subgroups; we will see in item 5.30 that this in fact occurs in axiomatic B.

from where one obtains, taking (e) into account, that

$$g = \mathring{\Psi}(g_1) \quad \text{and} \quad h = \mathring{\Phi}(h_1)$$

and, consequently, that

$$\mathring{g} = \mathring{\Phi}\Big(\mathring{\Psi}(g_1)\Big) \quad \text{and} \quad \mathring{h} = \mathring{\Psi}\Big(\mathring{\Phi}(h_1)\Big).$$

Since, by (e), $\mathring{\gamma} : H(\gamma) \longrightarrow \mathring{H}(\gamma)$ is an isomorphism, then,

$$\mathring{\Phi}\mathring{\Psi} = \mathring{\gamma}(\Phi)\mathring{\gamma}(\Psi) = \mathring{\gamma}(\Phi\Psi) = \mathring{\gamma}(\Psi\Phi) = \mathring{\gamma}(\Psi)\mathring{\gamma}(\Phi) = \mathring{\Psi}\mathring{\Phi},$$

and, hence, we have that

$$\mathring{g} = (\mathring{\Phi}\mathring{\Psi})(g_1) \quad \text{and} \quad \mathring{h} = (\mathring{\Phi}\mathring{\Psi})(h_1).$$

Since $\mathring{\Phi}\mathring{\Psi} \in \mathring{H}(\gamma)$, then, by (g), $\mathring{\Phi}\mathring{\Psi}$ is an endomorphism on $\mathring{G}(\gamma)$ and, thus,

$$\mathring{g} + \mathring{h} = (\mathring{\Phi}\mathring{\Psi})(g_1) + (\mathring{\Phi}\mathring{\Psi})(h_1) = (\mathring{\Phi}\mathring{\Psi})(g_1 + h_1).$$

Now, we calculate:

$$\begin{aligned}
\mathring{\Theta}_{(\gamma',\gamma)}(\mathring{g} + \mathring{h}) &= \mathring{\Theta}_{(\gamma',\gamma)}\Big((\mathring{\Phi}\mathring{\Psi})(g_1 + h_1)\Big) = \\
&= \mathring{\imath}_{(\gamma',\gamma)}(\mathring{\Phi}\mathring{\Psi})\Big(\mathring{\Theta}_{(\gamma',\gamma)}(g_1 + h_1)\Big) = && \text{(by (q))} \\
&= \mathring{\imath}_{(\gamma',\gamma)}(\mathring{\Phi}\mathring{\Psi})\Big(\overset{*}{\Theta}_{(\gamma',\gamma)}(g_1 + h_1)\Big) = && \text{(by the def. of } \mathring{\Theta}_{(\gamma',\gamma)}) \\
&= \mathring{\imath}_{(\gamma',\gamma)}(\mathring{\Phi}\mathring{\Psi})\Big(\Theta_{(\gamma',\gamma)}(g_1 + h_1)\Big) = && \text{(by Axiom (5-1))} \\
&= \mathring{\imath}_{(\gamma',\gamma)}(\mathring{\Phi}\mathring{\Psi})\Big(\Theta_{(\gamma',\gamma)}(g_1) + \Theta_{(\gamma',\gamma)}(h_1)\Big).
\end{aligned}$$

In this last passage we resorted to the fact of the functions $\Theta_{(\gamma',\gamma)} \in \Theta(\Delta(I))$ being homomorphisms.

Now, from the definition of $\mathring{\imath}_{(\gamma',\gamma)}$ we have that $\mathring{\imath}_{(\gamma',\gamma)}(\mathring{\Phi}\mathring{\Psi}) \in \mathring{H}(\gamma')$ being it, therefore, an endomorphism on the group $\mathring{G}(\gamma')$. With this, we have that (remembering that $\Theta_{(\gamma',\gamma)}(g_1)$ and $\Theta_{(\gamma',\gamma)}(h_1)$ belongs to $G(\gamma')$ and that, by (a), $G(\gamma') \subseteq \mathring{G}(\gamma')$):

$$\begin{aligned}
\mathring{\Theta}_{(\gamma',\gamma)}(\mathring{g} + \mathring{h}) &= \mathring{\imath}_{(\gamma',\gamma)}(\mathring{\Phi}\mathring{\Psi})\Big(\Theta_{(\gamma',\gamma)}(g_1) + \Theta_{(\gamma',\gamma)}(h_1)\Big) = \\
&= \mathring{\imath}_{(\gamma',\gamma)}(\mathring{\Phi}\mathring{\Psi})\Big(\Theta_{(\gamma',\gamma)}(g_1)\Big) + \mathring{\imath}_{(\gamma',\gamma)}(\mathring{\Phi}\mathring{\Psi})\Big(\Theta_{(\gamma',\gamma)}(h_1)\Big) = \\
&= \mathring{\imath}_{(\gamma',\gamma)}(\mathring{\Phi}\mathring{\Psi})\Big(\overset{*}{\Theta}_{(\gamma',\gamma)}(g_1)\Big) + \mathring{\imath}_{(\gamma',\gamma)}(\mathring{\Phi}\mathring{\Psi})\Big(\overset{*}{\Theta}_{(\gamma',\gamma)}(h_1)\Big) = && \text{(by Axiom (5-1))} \\
&= \mathring{\imath}_{(\gamma',\gamma)}(\mathring{\Phi}\mathring{\Psi})\Big(\mathring{\Theta}_{(\gamma',\gamma)}(g_1)\Big) + \mathring{\imath}_{(\gamma',\gamma)}(\mathring{\Phi}\mathring{\Psi})\Big(\mathring{\Theta}_{(\gamma',\gamma)}(h_1)\Big) = && \text{(def. of } \mathring{\Theta}_{(\gamma',\gamma)}) \\
&= \mathring{\Theta}_{(\gamma',\gamma)}\Big((\mathring{\Phi}\mathring{\Psi})(g_1)\Big) + \mathring{\Theta}_{(\gamma',\gamma)}\Big((\mathring{\Phi}\mathring{\Psi})(h_1)\Big). && \text{(by (q))}
\end{aligned}$$

Since $\mathring{g} = (\mathring{\Phi}\mathring{\Psi})(g_1)$ and $\mathring{h} = (\mathring{\Phi}\mathring{\Psi})(h_1)$, we conclude that

$$\mathring{\Theta}_{(\gamma',\gamma)}(\mathring{g} + \mathring{h}) = \mathring{\Theta}_{(\gamma',\gamma)}(\mathring{g}) + \mathring{\Theta}_{(\gamma',\gamma)}(\mathring{h}).$$

Based on the results above obtained ((a) to (r)) and taking into account the pertinent definitions, one can easily conclude that the family

$$\mathring{\Theta}\Big(\Delta(I)\Big) = \Big\{\mathring{\Theta}_{(\gamma',\gamma)}\Big\}_{(\gamma',\gamma)\in\Delta(I)}$$

is a restriction (Definition 3.13) for the bonded family (in (m))

$$\bigg(\mathring{\mathbb{G}}\Big(\Gamma(I)\Big), \mathring{\imath}\Big(\Gamma^2(I)\Big)\bigg),$$

and, also, a prolongation (Definition 3.23), to this bonded family, of the family

$$\Theta\Big(\Delta(I)\Big) = \Big\{\Theta_{(\gamma',\gamma)}\Big\}_{(\gamma',\gamma)\in\Delta(I)}.$$

In short, we have that:

**(s)** The triplet $\mathring{\mathscr{G}}(I)$ defined by

$$\mathring{\mathscr{G}}(I) := \bigg(\mathring{\mathbb{G}}\Big(\Gamma(I)\Big), \mathring{\imath}\Big(\Gamma^2(I)\Big), \mathring{\Theta}\Big(\Delta(I)\Big)\bigg)$$

is a $S$-space and also a strict and closed extension of the $S$-space $\mathscr{G}(I)$.

This last result, that, in a certain way, encapsulate the logical consequences of axiomatic B obtained until here, shows us that the $\widetilde{\mathscr{G}}(I)$-distributions defined by the (categoric) axiomatic formulated in item 5.12 (or by its equivalent simplified version given in 5.17), are present into axiomatic B. That is, the axioms of the axiomatic B allow (support) the definitions of the families $\mathring{\mathbb{G}}(\Gamma(I))$, $\mathring{\imath}(\Gamma^2(I))$ and $\mathring{\Theta}(\Delta(I))$ such as above formulated and, evenmore, they ensure, through its logical consequences, that these defined terms have the properties synthesized in (s) which, taking the 1st TESS into account, permit us to conclude that the triplet

$$\mathring{\mathscr{G}}(I) := \bigg(\mathring{\mathbb{G}}\Big(\Gamma(I)\Big), \mathring{\imath}\Big(\Gamma^2(I)\Big), \mathring{\Theta}\Big(\Delta(I)\Big)\bigg)$$

is a $S$-space isomorphic to the $S$-space $\widetilde{\mathscr{G}}(I)$ such as defined in the 1st TESS that, as we know, unless isomorphism, is the only model of the $\widetilde{\mathscr{G}}(I)$-distributions axiomatic formulated in 5.12;

What about the $\overline{\mathscr{G}}(I)$-distributions defined categorically by axiomatic A? Are they present into axiomatic B? The next proposition provides some other consequences of our axioms that will allow us to answer this question.

## 5.30 Proposition

*Let $\gamma \in \Gamma(I)$ be arbitrarily fixed. There exists a unique binary operation in $\overset{*}{G}(\gamma)$, $\oplus : \overset{*}{G}(\gamma) \times \overset{*}{G}(\gamma) \longrightarrow \overset{*}{G}(\gamma)$, denominated addition, such that:*

**(a)** *$\overset{*}{G}(\gamma)$ with the $\oplus$ operation is an abelian group that has the group $\overset{\circ}{G}(\gamma)$ (described in 5.29(g)) as a subgroup;*

**(b)** *the derivatives $\overset{*}{\Phi} \in \overset{*}{H}(\gamma)$ are endomorphisms on the group $\overset{*}{G}(\gamma)$ (described in (a));*

**(c)** *the restrictions $\overset{*}{\Theta}_{(\gamma',\gamma)} \in \overset{*}{\Theta}(\Delta(I))$ are homomorphisms between the groups $\overset{*}{G}(\gamma)$ and $\overset{*}{G}(\gamma')$.*

*Proof.*

**Part 1**: Definition of $\oplus$.

Let $\overset{*}{g}, \overset{*}{h} \in \overset{*}{G}(\gamma)$ be arbitrarily chosen. By Axiom 6 (and taking into account the definition of $\overset{\circ}{G}(\delta)$, with $\delta \in \Gamma(I)$), for each $x \in \gamma$ there exist $\gamma_x, \gamma'_x \in \Gamma(\gamma)$ such that $x \in \gamma_x \cap \gamma'_x$ and

$$\overset{*}{\Theta}_{(\gamma_x,\gamma)}(\overset{*}{g}) \in \overset{\circ}{G}(\gamma_x) \quad \text{and} \quad \overset{*}{\Theta}_{(\gamma'_x,\gamma)}(\overset{*}{h}) \in \overset{\circ}{G}(\gamma'_x).$$

Let us take

$$\overline{\xi}(\gamma) := \Big\{\gamma_x \cap \gamma'_x : x \in \gamma\Big\}.$$

Clearly, $\overline{\xi}(\gamma) \subseteq \Gamma(\gamma)$ and it is a cover of $\gamma$.

Taking now into account the Axiom 5((5-2)), as well as the definition of the functions $\overset{\circ}{\Theta}_{(\gamma',\gamma)} \in \overset{\circ}{\Theta}(\Delta(I))$, one has that:

$$\overset{*}{\Theta}_{(\gamma_x\cap\gamma'_x,\gamma)}(\overset{*}{g}) = \overset{*}{\Theta}_{(\gamma_x\cap\gamma'_x,\gamma_x)}\Big(\overset{*}{\Theta}_{(\gamma_x,\gamma)}(\overset{*}{g})\Big) = \overset{\circ}{\Theta}_{(\gamma_x\cap\gamma'_x,\gamma_x)}\Big(\overset{*}{\Theta}_{(\gamma_x,\gamma)}(\overset{*}{g})\Big) \in \overset{\circ}{G}(\gamma_x \cap \gamma'_x)$$

and

$$\overset{*}{\Theta}_{(\gamma_x\cap\gamma'_x,\gamma)}(\overset{*}{h}) = \overset{*}{\Theta}_{(\gamma_x\cap\gamma'_x,\gamma'_x)}\Big(\overset{*}{\Theta}_{(\gamma'_x,\gamma)}(\overset{*}{h})\Big) = \overset{\circ}{\Theta}_{(\gamma_x\cap\gamma'_x,\gamma'_x)}\Big(\overset{*}{\Theta}_{(\gamma'_x,\gamma)}(\overset{*}{h})\Big) \in \overset{\circ}{G}(\gamma_x \cap \gamma'_x),$$

that is, $\overline{\xi}(\gamma)$ is such that

$$\overset{*}{\Theta}_{(\xi,\gamma)}(\overset{*}{g}) \in \overset{\circ}{G}(\xi) \quad \text{and} \quad \overset{*}{\Theta}_{(\xi,\gamma)}(\overset{*}{h}) \in \overset{\circ}{G}(\xi)$$

for every $\xi = \gamma_x \cap \gamma'_x \in \overline{\xi}(\gamma)$.

Let now the family $\overset{\circ}{f}(\overline{\xi}(\gamma))$ be defined as follows:

$$\overset{\circ}{f}\Big(\overline{\xi}(\gamma)\Big) := \Big\{\overset{\circ}{f}_\xi := \overset{*}{\Theta}_{(\xi,\gamma)}(\overset{*}{g}) + \overset{*}{\Theta}_{(\xi,\gamma)}(\overset{*}{h})\Big\}_{\xi\in\overline{\xi}(\gamma)}.$$

Observe that, for each $\xi \in \overline{\xi}(\gamma)$, we have $\overset{\circ}{f}_\xi \in \overset{\circ}{G}(\xi)$ once that, according to 5.29(g), $\overset{\circ}{G}(\xi)$ with the $+$ operation is a group and the terms $\overset{*}{\Theta}_{(\xi,\gamma)}(\overset{*}{g})$ and $\overset{*}{\Theta}_{(\xi,\gamma)}(\overset{*}{h})$ of the sum $\overset{\circ}{f}_\xi$ belongs to $\mathring{G}(\xi)$. Furthermore, the family $\mathring{f}(\overline{\xi}(\gamma))$ is coherent (Definition at page 284). In fact, let $\xi, \xi' \in \overline{\xi}(\gamma)$ be such that $\xi \cap \xi' \neq \varnothing$. Resorting to the definition of the functions $\mathring{\Theta}_{(\gamma',\gamma)} \in \mathring{\Theta}(\Delta(I))$, to the fact of them being homomorphisms (according to 5.29(r)) and Axiom 5((5-2)), we get:

$$\begin{aligned}\overset{*}{\Theta}_{(\xi\cap\xi',\xi)}(\mathring{f}_\xi) = \mathring{\Theta}_{(\xi\cap\xi',\xi)}(\mathring{f}_\xi) = \mathring{\Theta}_{(\xi\cap\xi',\xi)}\left(\overset{*}{\Theta}_{(\xi,\gamma)}(\overset{*}{g}) + \overset{*}{\Theta}_{(\xi,\gamma)}(\overset{*}{h})\right) = \\ = \mathring{\Theta}_{(\xi\cap\xi',\xi)}\left(\overset{*}{\Theta}_{(\xi,\gamma)}(\overset{*}{g})\right) + \mathring{\Theta}_{(\xi\cap\xi',\xi)}\left(\overset{*}{\Theta}_{(\xi,\gamma)}(\overset{*}{h})\right) = \\ = \overset{*}{\Theta}_{(\xi\cap\xi',\xi)}\left(\overset{*}{\Theta}_{(\xi,\gamma)}(\overset{*}{g})\right) + \overset{*}{\Theta}_{(\xi\cap\xi',\xi)}\left(\overset{*}{\Theta}_{(\xi,\gamma)}(\overset{*}{h})\right) = \\ = \overset{*}{\Theta}_{(\xi\cap\xi',\gamma)}(\overset{*}{g}) + \overset{*}{\Theta}_{(\xi\cap\xi',\gamma)}(\overset{*}{h})\end{aligned}$$

and, analogously,

$$\overset{*}{\Theta}_{(\xi\cap\xi',\xi')}(\mathring{f}_{\xi'}) = \overset{*}{\Theta}_{(\xi\cap\xi',\gamma)}(\overset{*}{g}) + \overset{*}{\Theta}_{(\xi\cap\xi',\gamma)}(\overset{*}{h}),$$

that is,

$$\overset{*}{\Theta}_{(\xi\cap\xi',\xi)}(\mathring{f}_\xi) = \overset{*}{\Theta}_{(\xi\cap\xi',\xi')}(\mathring{f}_{\xi'}).$$

Thus being, by Axiom 7, there exists a single $\overset{*}{f} \in \overset{*}{G}(\gamma)$ such that

$$\overset{*}{\Theta}_{(\xi,\gamma)}(\overset{*}{f}) = \mathring{f}_\xi = \overset{*}{\Theta}_{(\xi,\gamma)}(\overset{*}{g}) + \overset{*}{\Theta}_{(\xi,\gamma)}(\overset{*}{h}) \tag{5.30-1}$$

for every $\xi \in \overline{\xi}(\gamma)$. Furthermore, this $\overset{*}{f}$ does not depend on the cover $\overline{\xi}(\gamma)$ of $\gamma$ in the following sense: being $\overline{\eta} \subseteq \Gamma(\gamma)$ another cover of $\gamma$ such that

$$\overset{*}{\Theta}_{(\eta,\gamma)}(\overset{*}{g}) \in \mathring{G}(\eta) \quad \text{and} \quad \overset{*}{\Theta}_{(\eta,\gamma)}(\overset{*}{h}) \in \mathring{G}(\eta)$$

for every $\eta \in \overline{\eta}(\gamma)$, and $\overset{*}{e} \in \overset{*}{G}(\gamma)$ (the only one) such that

$$\overset{*}{\Theta}_{(\eta,\gamma)}(\overset{*}{e}) = \overset{*}{\Theta}_{(\eta,\gamma)}(\overset{*}{g}) + \overset{*}{\Theta}_{(\eta,\gamma)}(\overset{*}{h}) \tag{5.30-2}$$

for every $\eta \in \overline{\eta}(\gamma)$, then, $\overset{*}{e} = \overset{*}{f}$.

In fact, let, for each $x \in \gamma$, $\xi_x \in \overline{\xi}(\gamma)$ and $\eta_x \in \overline{\eta}(\gamma)$ be such that $x \in \xi_x \cap \eta_x$ and take

$$\overline{\nu}(\gamma) := \left\{\nu_x := \xi_x \cap \eta_x : x \in \gamma\right\}.$$

It results that $\overline{\nu}(\gamma) \subseteq \Gamma(\gamma)$ and it is a cover of $\gamma$. Furthermore, taking (5.30-1) into account and using Axiom 5((5-2)), from the definition of $\mathring{\Theta}_{(\gamma',\gamma)} \in \mathring{\Theta}(\Delta(I))$ and the fact of them being homomorphisms, we have, for each $x \in \gamma$, that

$$\begin{aligned}\overset{*}{\Theta}_{(\nu_x,\gamma)}(\overset{*}{f}) = \overset{*}{\Theta}_{(\nu_x,\xi_x)}\left(\overset{*}{\Theta}_{(\xi_x,\gamma)}(\overset{*}{f})\right) = \overset{*}{\Theta}_{(\nu_x,\xi_x)}\left(\overset{*}{\Theta}_{(\xi_x,\gamma)}(\overset{*}{g}) + \overset{*}{\Theta}_{(\xi_x,\gamma)}(\overset{*}{h})\right) = \\ = \mathring{\Theta}_{(\nu_x,\xi_x)}\left(\overset{*}{\Theta}_{(\xi_x,\gamma)}(\overset{*}{g}) + \overset{*}{\Theta}_{(\xi_x,\gamma)}(\overset{*}{h})\right) = \\ = \mathring{\Theta}_{(\nu_x,\xi_x)}\left(\overset{*}{\Theta}_{(\xi_x,\gamma)}(\overset{*}{g})\right) + \mathring{\Theta}_{(\nu_x,\xi_x)}\left(\overset{*}{\Theta}_{(\xi_x,\gamma)}(\overset{*}{h})\right) = \\ = \overset{*}{\Theta}_{(\nu_x,\xi_x)}\left(\overset{*}{\Theta}_{(\xi_x,\gamma)}(\overset{*}{g})\right) + \overset{*}{\Theta}_{(\nu_x,\xi_x)}\left(\overset{*}{\Theta}_{(\xi_x,\gamma)}(\overset{*}{h})\right) = \\ = \overset{*}{\Theta}_{(\nu_x,\gamma)}(\overset{*}{g}) + \overset{*}{\Theta}_{(\nu_x,\gamma)}(\overset{*}{h}),\end{aligned}$$

that is,

$$\overset{*}{\Theta}_{(\nu,\gamma)}(\overset{*}{f}) = \overset{*}{\Theta}_{(\nu,\gamma)}(\overset{*}{g}) + \overset{*}{\Theta}_{(\nu,\gamma)}(\overset{*}{h}) \quad \text{for every} \quad \nu \in \overline{\nu}(\gamma).$$

Analogously, one obtains from (5.30-2) that

$$\overset{*}{\Theta}_{(\nu,\gamma)}(\overset{*}{e}) = \overset{*}{\Theta}_{(\nu,\gamma)}(\overset{*}{g}) + \overset{*}{\Theta}_{(\nu,\gamma)}(\overset{*}{h}) \quad \text{for every} \quad \nu \in \overline{\nu}(\gamma).$$

Therefore,

$$\overset{*}{\Theta}_{(\nu,\gamma)}(\overset{*}{f}) = \overset{*}{\Theta}_{(\nu,\gamma)}(\overset{*}{e}) \quad \text{for every} \quad \nu \in \overline{\nu}(\gamma),$$

from where one concludes, due to Axiom 7, that

$$\overset{*}{f} = \overset{*}{e}.$$

The results above backs the definition given ahead for the binary operation $\oplus$ in $\overset{*}{G}(\gamma)$.

**Definition.** The addition, $\oplus$, in $\overset{*}{G}(\gamma)$ is the function defined as follows:

$$\begin{aligned} \oplus : \overset{*}{G}(\gamma) \times \overset{*}{G}(\gamma) &\longrightarrow \overset{*}{G}(\gamma) \\ (\overset{*}{g}, \overset{*}{h}) &\longmapsto \overset{*}{g} \oplus \overset{*}{h} \end{aligned}$$

where $\overset{*}{g} \oplus \overset{*}{h} \in \overset{*}{G}(\gamma)$ is such that

$$\overset{*}{\Theta}_{(\xi,\gamma)}(\overset{*}{g} \oplus \overset{*}{h}) = \overset{*}{\Theta}_{(\xi,\gamma)}(\overset{*}{g}) + \overset{*}{\Theta}_{(\xi,\gamma)}(\overset{*}{h})$$

for every $\xi \in \overline{\xi}(\gamma)$, being $\overline{\xi}(\gamma) \subseteq \Gamma(\gamma)$ a cover of $\gamma$ such that

$$\overset{*}{\Theta}_{(\xi,\gamma)}(\overset{*}{g}) \in \overset{\circ}{G}(\xi) \quad \text{and} \quad \overset{*}{\Theta}_{(\xi,\gamma)}(\overset{*}{h}) \in \overset{\circ}{G}(\xi)$$

for every $\xi \in \overline{\xi}(\gamma)$.

Let us remark that, for $\overset{\circ}{g}, \overset{\circ}{h} \in \overset{\circ}{G}(\gamma) \subseteq \overset{*}{G}(\gamma)$,

$$\overset{\circ}{g} \oplus \overset{\circ}{h} = \overset{\circ}{g} + \overset{\circ}{h},$$

where $+$ (that figures at the definition above for $\oplus$) is the addition of the group $\overset{\circ}{G}(\gamma)$ (according to 5.29(g)); in other terms, $\oplus$ is an extension to $\overset{*}{G}(\gamma)$ of the addition, $+$, of the group $\overset{\circ}{G}(\gamma)$. In fact, for any cover $\overline{\xi}(\gamma) \subseteq \Gamma(\gamma)$ of $\gamma$ we have, for $\overset{\circ}{g}$ and $\overset{\circ}{h}$ in $\overset{\circ}{G}(\gamma)$, that

$$\overset{*}{\Theta}_{(\xi,\gamma)}(\overset{\circ}{g}) \Big( = \overset{\circ}{\Theta}_{(\xi,\gamma)}(\overset{\circ}{g}) \Big) \in \overset{\circ}{G}(\xi)$$

and

$$\overset{*}{\Theta}_{(\xi,\gamma)}(\overset{\circ}{h}) \Big( = \overset{\circ}{\Theta}_{(\xi,\gamma)}(\overset{\circ}{h}) \Big) \in \overset{\circ}{G}(\xi)$$

for every $\xi \in \overline{\xi}(\gamma)$. Furthermore, since $\overset{\circ}{g} + \overset{\circ}{h} \in \overset{\circ}{G}(\gamma)$, once that $\overset{\circ}{G}(\gamma)$ with the $+$ operation is a group, comes that (remembering ourselves that the functions $\overset{\circ}{\Theta}_{(\gamma',\gamma)} \in \overset{\circ}{\Theta}(\Delta(I))$ are homomorphisms):

$$\overset{*}{\Theta}_{(\xi,\gamma)}(\overset{\circ}{g} + \overset{\circ}{h}) = \overset{\circ}{\Theta}_{(\xi,\gamma)}(\overset{\circ}{g} + \overset{\circ}{h}) = \overset{\circ}{\Theta}_{(\xi,\gamma)}(\overset{\circ}{g}) + \overset{\circ}{\Theta}_{(\xi,\gamma)}(\overset{\circ}{h}) = \overset{*}{\Theta}_{(\xi,\gamma)}(\overset{\circ}{g}) + \overset{*}{\Theta}_{(\xi,\gamma)}(\overset{\circ}{h})$$

for every $\xi \in \overline{\xi}(\gamma)$. Therefore, taking into account the definition of $\oplus$, we get that

$$\overset{\circ}{g} \oplus \overset{\circ}{h} = \overset{\circ}{g} + \overset{\circ}{h}.$$

It trivially follows from the definition of $\oplus$, taking into account the commutativity and associativity of the addition $+$ in $\overset{\circ}{G}(\delta)$, $\delta \in \Gamma(I)$, that $\oplus$ also is a commutative and associative operation.

**Part 2**: $(\overset{*}{G}(\gamma), \oplus)$ is an abelian group.

We already saw in Part 1 that $\oplus$ is a commutative and associative operation. It remains then to prove the existence of an additive neutral and, for each $\overset{*}{g} \in \overset{*}{G}(\gamma)$, the existence of an opposite element.

Let us start proving that the additive neutral $\overset{\circ}{0} \in \overset{\circ}{G}(\gamma)$ of the abelian group $(\overset{\circ}{G}(\gamma), +)$ also is a neutral element for the $\oplus$ operation, i.e., that

$$\overset{*}{g} \oplus \overset{\circ}{0} = \overset{*}{g} \quad \text{for every} \quad \overset{*}{g} \in \overset{*}{G}(\gamma).$$

In fact, let $\overline{\xi}(\gamma) \subseteq \Gamma(\gamma)$ be a cover of $\gamma$ such that

$$\overset{*}{\Theta}_{(\xi,\gamma)}(\overset{*}{g}) \in \overset{\circ}{G}(\xi) \quad \text{for every} \quad \xi \in \overline{\xi}(\gamma).$$

(We already know that Axiom 6 ensures the existence of one such cover $\overline{\xi}(\gamma)$).

Since, by the definition of $\overset{\circ}{\Theta}_{(\xi,\gamma)}$,

$$\overset{*}{\Theta}_{(\xi,\gamma)}(\overset{\circ}{0}) = \overset{\circ}{\Theta}_{(\xi,\gamma)}(\overset{\circ}{0}) \in \overset{\circ}{G}(\xi)$$

for every $\xi \in \overline{\xi}(\gamma)$, we have, by the definition of the $\oplus$ operation, that $\overset{*}{g} \oplus \overset{\circ}{0}$ is the element of $\overset{*}{G}(\gamma)$ such that

$$\overset{*}{\Theta}_{(\xi,\gamma)}(\overset{*}{g} \oplus \overset{\circ}{0}) = \overset{*}{\Theta}_{(\xi,\gamma)}(\overset{*}{g}) + \overset{*}{\Theta}_{(\xi,\gamma)}(\overset{\circ}{0})$$

for every $\xi \in \overline{\xi}(\gamma)$.

But, $\overset{\circ}{\Theta}_{(\xi,\gamma)} : \overset{\circ}{G}(\gamma) \longrightarrow \overset{\circ}{G}(\xi)$ is a homomorphism and, thus, $\overset{\circ}{\Theta}_{(\xi,\gamma)}(\overset{\circ}{0})$ is the additive neutral of the group $\overset{\circ}{G}(\xi)$. Therefore,

$$\overset{*}{\Theta}_{(\xi,\gamma)}(\overset{*}{g} \oplus \overset{\circ}{0}) = \overset{*}{\Theta}_{(\xi,\gamma)}(\overset{*}{g}) \quad \text{for every} \quad \xi \in \overline{\xi}(\gamma).$$

Now, since $\{\overset{*}{\Theta}_{\xi,\gamma}(\overset{*}{g})\}_{\xi\in\overline{\xi}(\gamma)}$ is, clearly, a coherent family, it results from the last equation, along with Axiom 7, that

$$\overset{*}{g} \oplus \overset{\circ}{0} = \overset{*}{g}.$$

Let us prove now that, for each $\overset{*}{g} \in \overset{*}{G}(\gamma)$, there exists $\overset{*}{g}_{(-)} \in \overset{*}{G}(\gamma)$ such that

$$\overset{*}{g} \oplus \overset{*}{g}_{(-)} = \overset{\circ}{0}.$$

In order to do so, we resort again to a cover $\overline{\xi}(\gamma) \subseteq \Gamma(\gamma)$ of $\gamma$ such that

$$\overset{*}{\Theta}_{(\xi,\gamma)}(\overset{*}{g}) \in \overset{\circ}{G}(\xi) \quad \text{for every} \quad \xi \in \overline{\xi}(\gamma),$$

with $\overset{*}{g} \in \overset{*}{G}(\gamma)$ arbitrarily fixed.

Let now

$$\left\{-\overset{*}{\Theta}_{(\xi,\gamma)}(\overset{*}{g})\right\}_{\xi \in \overline{\xi}(\gamma)}$$

be the family where $-\overset{*}{\Theta}_{(\xi,\gamma)}(\overset{*}{g}) \in \overset{\circ}{G}(\xi)$ is the opposite element of $\overset{*}{\Theta}_{(\xi,\gamma)}(\overset{*}{g})$ in the group $\overset{\circ}{G}(\xi)$, that is, the only element in $\overset{\circ}{G}(\xi)$ such that

$$\overset{*}{\Theta}_{(\xi,\gamma)}(\overset{*}{g}) + \left(-\overset{*}{\Theta}_{(\xi,\gamma)}(\overset{*}{g})\right) = \overset{\circ}{0}_\xi, \tag{5.30-3}$$

being $\overset{\circ}{0}_\xi \in \overset{\circ}{G}(\xi)$ the additive neutral of the group $\overset{\circ}{G}(\xi)$. This family is coherent, once that for $\xi, \xi' \in \overline{\xi}(\gamma)$ such that $\xi \cap \xi' \neq \varnothing$, we have (employing the definition of $\overset{\circ}{\Theta}_{(\gamma',\gamma)} \in \overset{\circ}{\Theta}(\Delta(I))$, that they are homomorphisms and Axiom 5((5-2))):

$$\begin{aligned}
\overset{*}{\Theta}_{(\xi\cap\xi',\xi)}\left(-\overset{*}{\Theta}_{(\xi,\gamma)}(\overset{*}{g})\right) &= \overset{\circ}{\Theta}_{(\xi\cap\xi',\xi)}\left(-\overset{*}{\Theta}_{(\xi,\gamma)}(\overset{*}{g})\right) = \\
&= -\overset{\circ}{\Theta}_{(\xi\cap\xi',\xi)}\left(\overset{*}{\Theta}_{(\xi,\gamma)}(\overset{*}{g})\right) = \\
&= -\overset{*}{\Theta}_{(\xi\cap\xi',\xi)}\left(\overset{*}{\Theta}_{(\xi,\gamma)}(\overset{*}{g})\right) = \\
&= -\overset{*}{\Theta}_{(\xi\cap\xi',\gamma)}(\overset{*}{g})
\end{aligned}$$

and, analogously,

$$\overset{*}{\Theta}_{(\xi\cap\xi',\xi')}\left(-\overset{*}{\Theta}_{(\xi',\gamma)}(\overset{*}{g})\right) = -\overset{*}{\Theta}_{(\xi\cap\xi',\gamma)}(\overset{*}{g}),$$

that is,

$$\overset{*}{\Theta}_{(\xi\cap\xi',\xi)}\left(-\overset{*}{\Theta}_{(\xi,\gamma)}(\overset{*}{g})\right) = \overset{*}{\Theta}_{(\xi\cap\xi',\xi')}\left(-\overset{*}{\Theta}_{(\xi',\gamma)}(\overset{*}{g})\right).$$

Now, since $\{-\overset{*}{\Theta}_{(\xi,\gamma)}(\overset{*}{g})\}_{\xi\in\overline{\xi}(\gamma)}$ is coherent, there exists, as established by Axiom 7, a single element in $\overset{*}{G}(\gamma)$, let us say $\overset{*}{g}_{(-)}$, such that

$$\overset{*}{\Theta}_{(\xi,\gamma)}(\overset{*}{g}_{(-)}) = -\overset{*}{\Theta}_{(\xi,\gamma)}(\overset{*}{g})\left(\in \overset{\circ}{G}(\xi)\right)$$

for every $\xi \in \overline{\xi}(\gamma)$.

Taking this result into (5.30-3) we obtain that

$$\overset{*}{\Theta}_{(\xi,\gamma)}(\overset{*}{g}) + \overset{*}{\Theta}_{(\xi,\gamma)}(\overset{*}{g}_{(-)}) = \overset{\circ}{0}_\xi$$

for every $\xi \in \overline{\xi}(\gamma)$, an expression that, if we remember ourselves that $\overset{\circ}{\Theta}_{(\xi,\gamma)} : \overset{\circ}{G}(\gamma) \longrightarrow \overset{\circ}{G}(\xi)$ is a homomorphism and, therefore, that

$$\overset{*}{\Theta}_{(\xi,\gamma)}(\overset{\circ}{0}) = \overset{\circ}{\Theta}_{(\xi,\gamma)}(\overset{\circ}{0}) = \overset{\circ}{0}_\xi,$$

can be written as follows:

$$\overset{*}{\Theta}_{(\xi,\gamma)}(\overset{\circ}{0}) = \overset{*}{\Theta}_{(\xi,\gamma)}(\overset{*}{g}) + \overset{*}{\Theta}_{(\xi,\gamma)}(\overset{*}{g}_{(-)})$$

for every $\xi \in \overline{\xi}(\gamma)$, with, as seen above,

$$\overset{*}{\Theta}_{(\xi,\gamma)}(\overset{*}{g}) \in \overset{\circ}{G}(\xi) \quad \text{and} \quad \overset{*}{\Theta}_{(\xi,\gamma)}(\overset{*}{g}_{(-)}) \in \overset{\circ}{G}(\xi)$$

for every $\xi \in \overline{\xi}(\gamma)$.

This conclusion, in turn, in light of the definition of $\oplus$, tells us exactly that

$$\overset{*}{g} \oplus \overset{*}{g}_{(-)} = \overset{\circ}{0}.$$

Thus we conclude the proof of $(\overset{*}{G}(\gamma), \oplus)$ being an abelian group. This result, along with the conclusion obtained in Part 1 about $\oplus$ being an extension to $\overset{*}{G}(\gamma) \supseteq \overset{\circ}{G}(\gamma)$ of the addition $+$ of the group $(\overset{\circ}{G}(\gamma), +)$, show us that $(\overset{\circ}{G}(\gamma), +)$ is a subgroup of the group $(\overset{*}{G}(\gamma), \oplus)$. Henceforth, we will use the same symbol, "+", to denote any of these two operations, $\oplus$ or $+$, leaving for the context to eliminate the ambiguity.

**Part 3**: $\overset{*}{\Phi} \in \overset{*}{H}(\gamma)$ is an endomorphism.

Let $\overset{*}{g}, \overset{*}{h} \in \overset{*}{G}(\gamma)$ be arbitrarily chosen and $\overset{*}{\Phi}$ be any element of $\overset{*}{H}(\gamma)$. Let also $\overline{\xi}(\gamma) \subseteq \Gamma(\gamma)$ be a cover of $\gamma$ such that, for each $\xi \in \overline{\xi}(\gamma)$, $\overset{*}{\Theta}_{(\xi,\gamma)}(\overset{*}{g})$, $\overset{*}{\Theta}_{(\xi,\gamma)}(\overset{*}{h})$, $\overset{*}{\Theta}_{(\xi,\gamma)}(\overset{*}{\Phi}(\overset{*}{g}))$ and $\overset{*}{\Theta}_{(\xi,\gamma)}(\overset{*}{\Phi}(\overset{*}{h}))$ belongs to $\overset{\circ}{G}(\xi)$.[28]

Hence, by the definition of addition given in Part 1, $\overset{*}{g} + \overset{*}{h}$ is the only element in $\overset{*}{G}(\gamma)$ such that

$$\overset{*}{\Theta}_{(\xi,\gamma)}(\overset{*}{g} + \overset{*}{h}) = \overset{*}{\Theta}_{(\xi,\gamma)}(\overset{*}{g}) + \overset{*}{\Theta}_{(\xi,\gamma)}(\overset{*}{h}) \quad \text{for every} \quad \xi \in \overline{\xi}(\gamma).$$

With this, and taking Axiom 5((5-3)) into account, we obtain, for each $\xi \in \overline{\xi}(\gamma)$, that

$$\overset{*}{\Theta}_{(\xi,\gamma)}\Big(\overset{*}{\Phi}(\overset{*}{g} + \overset{*}{h})\Big) = \overset{*}{\imath}_{(\xi,\gamma)}(\overset{*}{\Phi})\Big(\overset{*}{\Theta}_{(\xi,\gamma)}(\overset{*}{g} + \overset{*}{h})\Big) = \overset{*}{\imath}_{(\xi,\gamma)}(\overset{*}{\Phi})\Big(\overset{*}{\Theta}_{(\xi,\gamma)}(\overset{*}{g}) + \overset{*}{\Theta}_{(\xi,\gamma)}(\overset{*}{h})\Big) \quad (5.30\text{-}4)$$

On the other hand, from the definitions of $\overset{*}{\imath}_{(\xi,\gamma)}$ and $\overset{\circ}{\imath}_{(\xi,\gamma)}$ (at pages 283 and 289, respectively), we know that

$$\overset{*}{\imath}_{(\xi,\gamma)}(\overset{*}{\Phi}) = \overset{*}{\xi}\left(i_{(\xi,\gamma)}\Big((\overset{*}{\gamma})^{-1}(\overset{*}{\Phi})\Big)\right) = \overset{*}{\xi}\Big(i_{(\xi,\gamma)}(\Phi)\Big)$$

and

$$\overset{\circ}{\imath}_{(\xi,\gamma)}(\overset{\circ}{\Phi}) = \overset{\circ}{\xi}\left(i_{(\xi,\gamma)}\Big((\overset{\circ}{\gamma})^{-1}(\overset{\circ}{\Phi})\Big)\right) = \overset{\circ}{\xi}\Big(i_{(\xi,\gamma)}(\Phi)\Big),$$

and, hence, taking

$$\Psi := i_{(\xi,\gamma)}(\Phi) \in H(\xi),$$

[28] The existence of one such cover is ensured by Axiom 6, which can be verified through a construction intirely analogous to that, in Part 1 (when defining $\oplus$), of the cover $\overline{\xi}(\gamma)$.

one gets that

$$\overset{*}{\imath}_{(\xi,\gamma)}(\overset{*}{\Phi}) = \overset{*}{\xi}(\Psi) = \overset{*}{\Psi} \in \overset{*}{H}(\xi)$$

and

$$\overset{\circ}{\imath}_{(\xi,\gamma)}(\overset{\circ}{\Phi}) = \overset{\circ}{\xi}(\Psi) = \overset{\circ}{\Psi} \in \overset{\circ}{H}(\xi).$$

Since $\overset{\circ}{\Psi} = \overset{*}{\Psi}|_{\overset{\circ}{G}(\xi)}$, we have that

$$\overset{*}{\imath}_{(\xi,\gamma)}(\overset{*}{\Phi})(\overset{\circ}{g}) = \overset{\circ}{\imath}_{(\xi,\gamma)}(\overset{\circ}{\Phi})(\overset{\circ}{g}) \quad \text{for every} \quad \overset{\circ}{g} \in \overset{\circ}{G}(\xi).$$

Since $\overset{*}{\Theta}_{(\xi,\gamma)}(\overset{*}{g}) + \overset{*}{\Theta}_{(\xi,\gamma)}(\overset{*}{h}) \in \overset{\circ}{G}(\xi)$ and taking the last result into account, the expression (5.30-4) can be written as follows:

$$\overset{*}{\Theta}_{(\xi,\gamma)}\Big(\overset{*}{\Phi}(\overset{*}{g}+\overset{*}{h})\Big) = \overset{\circ}{\imath}_{(\xi,\gamma)}(\overset{\circ}{\Phi})\Big(\overset{*}{\Theta}_{(\xi,\gamma)}(\overset{*}{g}) + \overset{*}{\Theta}_{(\xi,\gamma)}(\overset{*}{h})\Big)$$

for every $\xi \in \overline{\xi}(\gamma)$.

Now, as we know, $\overset{\circ}{\imath}_{(\xi,\gamma)}(\overset{\circ}{\Phi}) \in \overset{\circ}{H}(\xi)$ is an endomorphism and, hence,

$$\begin{aligned}\overset{*}{\Theta}_{(\xi,\gamma)}\Big(\overset{*}{\Phi}(\overset{*}{g}+\overset{*}{h})\Big) &= \overset{\circ}{\imath}_{(\xi,\gamma)}(\overset{\circ}{\Phi})\Big(\overset{*}{\Theta}_{(\xi,\gamma)}(\overset{*}{g})\Big) + \overset{\circ}{\imath}_{(\xi,\gamma)}(\overset{\circ}{\Phi})\Big(\overset{*}{\Theta}_{(\xi,\gamma)}(\overset{*}{h})\Big) = \\ &= \overset{*}{\imath}_{(\xi,\gamma)}(\overset{*}{\Phi})\Big(\overset{*}{\Theta}_{(\xi,\gamma)}(\overset{*}{g})\Big) + \overset{*}{\imath}_{(\xi,\gamma)}(\overset{*}{\Phi})\Big(\overset{*}{\Theta}_{(\xi,\gamma)}(\overset{*}{h})\Big),\end{aligned}$$

that, considering Axiom 5((5-3)), can be written as

$$\overset{*}{\Theta}_{(\xi,\gamma)}\Big(\overset{*}{\Phi}(\overset{*}{g}+\overset{*}{h})\Big) = \overset{*}{\Theta}_{(\xi,\gamma)}\Big(\overset{*}{\Phi}(\overset{*}{g})\Big) + \overset{*}{\Theta}_{(\xi,\gamma)}\Big(\overset{*}{\Phi}(\overset{*}{h})\Big)$$

for every $\xi \in \overline{\xi}(\gamma)$. Taking now into account the definition of addition given in Part 1 and remembering ourselves that the cover $\overline{\xi}(\gamma)$ of $\gamma$ is such that

$$\overset{*}{\Theta}_{(\xi,\gamma)}\Big(\overset{*}{\Phi}(\overset{*}{g})\Big) \in \overset{\circ}{G}(\xi) \quad \text{and} \quad \overset{*}{\Theta}_{(\xi,\gamma)}\Big(\overset{*}{\Phi}(\overset{*}{h})\Big) \in \overset{\circ}{G}(\xi)$$

for every $\xi \in \overline{\xi}(\gamma)$, we finally conclude that

$$\overset{*}{\Phi}(\overset{*}{g}+\overset{*}{h}) = \overset{*}{\Phi}(\overset{*}{g}) + \overset{*}{\Phi}(\overset{*}{h}).$$

**Part 4**: $\overset{*}{\Theta}_{(\gamma',\gamma)} \in \overset{*}{\Theta}(\Delta(I))$ is a homomorphism.

Let $\overset{*}{g}, \overset{*}{h} \in \overset{*}{G}(\gamma)$ and $\overset{*}{\Theta}_{(\gamma',\gamma)} \in \overset{*}{\Theta}(\Delta(I))$ be arbitrarily chosen and let us take:

$$\overset{*}{g}_1 := \overset{*}{\Theta}_{(\gamma',\gamma)}(\overset{*}{g}), \quad \overset{*}{h}_1 := \overset{*}{\Theta}_{(\gamma',\gamma)}(\overset{*}{h}), \quad \text{and} \quad \overset{*}{f}_1 := \overset{*}{\Theta}_{(\gamma',\gamma)}(\overset{*}{g}+\overset{*}{h}).$$

Let also (backed by Axiom 6) $\overline{\eta}(\gamma) \subseteq \Gamma(\gamma)$ and $\overline{\chi}(\gamma') \subseteq \Gamma(\gamma')$ covers of $\gamma$ and $\gamma'$, respectively, such that:

$$\overset{*}{\Theta}_{(\eta,\gamma)}(\overset{*}{g}) \in \overset{\circ}{G}(\eta) \quad \text{and} \quad \overset{*}{\Theta}_{(\eta,\gamma)}(\overset{*}{h}) \in \overset{\circ}{G}(\eta) \quad \text{for every} \quad \eta \in \overline{\eta}(\gamma) \tag{5.30-5}$$

and

$$\overset{*}{\Theta}_{(\chi,\gamma')}(\overset{*}{g}_1) \in \overset{\circ}{G}(\chi) \quad \text{and} \quad \overset{*}{\Theta}_{(\chi,\gamma')}(\overset{*}{h}_1) \in \overset{\circ}{G}(\chi) \quad \text{for every} \quad \chi \in \overline{\chi}(\gamma').$$

Now, we define

$$\overline{\xi}(\gamma') := \Big\{ \xi := \chi \cap \eta \quad : \quad \chi \in \overline{\chi}(\gamma'), \quad \eta \in \overline{\eta}(\gamma) \quad \text{and} \quad \chi \cap \eta \neq \varnothing \Big\}.$$

Clearly, $\overline{\xi}(\gamma') \subseteq \Gamma(\gamma')$ is a cover of $\gamma'$ such that, for every $\xi = \chi \cap \eta \in \overline{\xi}(\gamma')$,

$$\overset{*}{\Theta}_{(\xi,\gamma')}(\overset{*}{g}_1) \in \overset{\circ}{G}(\xi) \quad \text{and} \quad \overset{*}{\Theta}_{(\xi,\gamma')}(\overset{*}{h}_1) \in \overset{\circ}{G}(\xi), \tag{5.30-6}$$

since

$$\overset{*}{\Theta}_{(\xi=\chi\cap\eta,\gamma')}(\overset{*}{g}_1) = \overset{*}{\Theta}_{(\chi\cap\eta,\chi)}\Big(\overset{*}{\Theta}_{(\chi,\gamma')}(\overset{*}{g}_1)\Big),$$

$\overset{*}{\Theta}_{(\chi,\gamma')}(\overset{*}{g}_1) \in \overset{\circ}{G}(\chi)$ and $\overset{*}{\Theta}_{(\chi\cap\eta,\chi)}|_{\overset{\circ}{G}(\chi)} = \overset{\circ}{\Theta}_{(\chi\cap\eta,\chi)}$ is a homomorphism from $\overset{\circ}{G}(\chi)$ into $\overset{\circ}{G}(\xi = \chi \cap \eta)$ (analogously for $\overset{*}{\Theta}_{(\xi,\gamma')}(\overset{*}{h}_1)$).

Now, the pertinence relations in (5.30-5) and (5.30-6), along with the definitions (given in Part 1) of the additions in $\overset{*}{G}(\gamma)$ and $\overset{*}{G}(\gamma')$, respectively, allow us to conclude that:

**(a)** $\overset{*}{g} + \overset{*}{h}$ is the only element in $\overset{*}{G}(\gamma)$ such that

$$\overset{*}{\Theta}_{(\eta,\gamma)}(\overset{*}{g} + \overset{*}{h}) = \overset{*}{\Theta}_{(\eta,\gamma)}(\overset{*}{g}) + \overset{*}{\Theta}_{(\eta,\gamma)}(\overset{*}{h}) \tag{5.30-7}$$

for every $\eta \in \overline{\eta}(\gamma)$;

**(b)** $\overset{*}{g}_1 + \overset{*}{h}_1$ is the only element in $\overset{*}{G}(\gamma')$ such that

$$\overset{*}{\Theta}_{(\xi,\gamma')}(\overset{*}{g}_1 + \overset{*}{h}_1) = \overset{*}{\Theta}_{(\xi,\gamma')}(\overset{*}{g}_1) + \overset{*}{\Theta}_{(\xi,\gamma')}(\overset{*}{h}_1),$$

or, taking into account the definitions of $\overset{*}{g}_1$ and $\overset{*}{h}_1$ as well as Axiom 5((5-2)), such that

$$\begin{aligned}\overset{*}{\Theta}_{(\xi,\gamma')}(\overset{*}{g}_1 + \overset{*}{h}_1) &= \overset{*}{\Theta}_{(\xi,\gamma')}\Big(\overset{*}{\Theta}_{(\gamma',\gamma)}(\overset{*}{g})\Big) + \overset{*}{\Theta}_{(\xi,\gamma')}\Big(\overset{*}{\Theta}_{(\gamma',\gamma)}(\overset{*}{h})\Big) = \\ &= \overset{*}{\Theta}_{(\xi,\gamma)}(\overset{*}{g}) + \overset{*}{\Theta}_{(\xi,\gamma)}(\overset{*}{h})\end{aligned} \tag{5.30-8}$$

for every $\xi \in \overline{\xi}(\gamma')$.

From the expression (5.30-7), taking into account the pertinence relations in (5.30-5) and also that $\overset{\circ}{\Theta}_{(\chi\cap\eta,\eta)} = \overset{*}{\Theta}_{(\chi\cap\eta,\eta)}|_{\overset{\circ}{G}(\eta)}$ is a homomorphism, we obtain:

$$\begin{aligned}\overset{*}{\Theta}_{(\chi\cap\eta,\eta)}\Big(\overset{*}{\Theta}_{(\eta,\gamma)}(\overset{*}{g} + \overset{*}{h})\Big) &= \overset{*}{\Theta}_{(\chi\cap\eta,\eta)}\Big(\overset{*}{\Theta}_{(\eta,\gamma)}(\overset{*}{g}) + \overset{*}{\Theta}_{(\eta,\gamma)}(\overset{*}{h})\Big) = \\ &= \overset{\circ}{\Theta}_{(\chi\cap\eta,\eta)}\Big(\overset{*}{\Theta}_{(\eta,\gamma)}(\overset{*}{g}) + \overset{*}{\Theta}_{(\eta,\gamma)}(\overset{*}{h})\Big) = \\ &= \overset{\circ}{\Theta}_{(\chi\cap\eta,\eta)}\Big(\overset{*}{\Theta}_{(\eta,\gamma)}(\overset{*}{g})\Big) + \overset{\circ}{\Theta}_{(\chi\cap\eta,\eta)}\Big(\overset{*}{\Theta}_{(\eta,\gamma)}(\overset{*}{h})\Big) = \\ &= \overset{*}{\Theta}_{(\chi\cap\eta,\eta)}\Big(\overset{*}{\Theta}_{(\eta,\gamma)}(\overset{*}{g})\Big) + \overset{*}{\Theta}_{(\chi\cap\eta,\eta)}\Big(\overset{*}{\Theta}_{(\eta,\gamma)}(\overset{*}{h})\Big),\end{aligned}$$

that is, using Axiom 5((5-2)),

$$\overset{*}{\Theta}_{(\xi,\gamma)}(\overset{*}{g}+\overset{*}{h}) = \overset{*}{\Theta}_{(\xi,\gamma)}(\overset{*}{g}) + \overset{*}{\Theta}_{(\xi,\gamma)}(\overset{*}{h}) \tag{5.30-9}$$

for every $\xi \in \overline{\xi}(\gamma')$.

On the other hand, since

$$\overset{*}{f}_1 := \overset{*}{\Theta}_{(\gamma',\gamma)}(\overset{*}{g}+\overset{*}{h}),$$

we have, for each $\xi \in \overline{\xi}(\gamma')$, that

$$\overset{*}{\Theta}_{(\xi,\gamma')}(\overset{*}{f}_1) = \overset{*}{\Theta}_{(\xi,\gamma')}\Big(\overset{*}{\Theta}_{(\gamma',\gamma)}(\overset{*}{g}+\overset{*}{h})\Big)$$

that is, again resorting to axiom Axiom 5((5-2)),

$$\overset{*}{\Theta}_{(\xi,\gamma')}(\overset{*}{f}_1) = \overset{*}{\Theta}_{(\xi,\gamma)}(\overset{*}{g}+\overset{*}{h})$$

and, hence, taking (5.30-9) into account, we obtain:

$$\overset{*}{\Theta}_{(\xi,\gamma')}(\overset{*}{f}_1) = \overset{*}{\Theta}_{(\xi,\gamma)}(\overset{*}{g}) + \overset{*}{\Theta}_{(\xi,\gamma)}(\overset{*}{h})$$

for every $\xi \in \overline{\xi}(\gamma')$.

But, as we saw in (b), $\overset{*}{g}_1 + \overset{*}{h}_1$ is the only element in $\overset{*}{G}(\gamma')$ that satisfies (5.30-8). Since the last equation shows us that $\overset{*}{f}_1$ also "satisfies" (5.30-8), we conclude that

$$\overset{*}{f}_1 = \overset{*}{g}_1 + \overset{*}{h}_1,$$

i.e., that

$$\overset{*}{\Theta}_{(\gamma',\gamma)}(\overset{*}{g}+\overset{*}{h}) = \overset{*}{\Theta}_{(\gamma',\gamma)}(\overset{*}{g}) + \overset{*}{\Theta}_{(\gamma',\gamma)}(\overset{*}{h}).$$

**Part 5**: The uniqueness of the addition.

We must prove now and, hence, complete the proof of our proposition, that the addition, $+$, in $\overset{*}{G}(\gamma)$ as defined in Part 1, is the only possible extension to $\overset{*}{G}(\gamma)$ of the addition of the group $\overset{\circ}{G}(\gamma)$, that attends the following conditions:

- $(\overset{*}{G}(\gamma), +)$ is an abelian group;
- $\overset{*}{\Phi} \in \overset{*}{H}(\gamma)$ is an endomorphism on the group $(\overset{*}{G}(\gamma), +)$;
- $\overset{*}{\Theta}_{(\gamma',\gamma)} \in \overset{*}{\Theta}(\Delta(I))$ is a homomorphism from the group $(\overset{*}{G}(\gamma), +)$ into the group $(\overset{*}{G}(\gamma'), +)$.

Let us then suppose that there exists another binary operation in $\overset{*}{G}(\gamma)$, let us say $\boxplus$, which is an extension to $\overset{*}{G}(\gamma)$ of the addition of the group $\overset{\circ}{G}(\gamma)$ and that attends the conditions above (with "$\boxplus$" in the place of "+"). Let $\overset{*}{g}, \overset{*}{h} \in \overset{*}{G}(\gamma)$ be arbitrarily fixed and $\overline{\xi}(\gamma) \subseteq \Gamma(\gamma)$ be a cover of $\gamma$ such that

$$\overset{*}{\Theta}_{(\xi,\gamma)}(\overset{*}{g}) \in \overset{\circ}{G}(\xi) \quad \text{and} \quad \overset{*}{\Theta}_{(\xi,\gamma)}(\overset{*}{h}) \in \overset{\circ}{G}(\xi)$$

for every $\xi \in \overline{\xi}(\gamma)$.

Hence, we have that:

$$\overset{*}{\Theta}_{(\xi,\gamma)}(\overset{*}{g} \boxplus \overset{*}{h}) = \overset{*}{\Theta}_{(\xi,\gamma)}(\overset{*}{g}) \boxplus \overset{*}{\Theta}_{(\xi,\gamma)}(\overset{*}{h}) = \overset{*}{\Theta}_{(\xi,\gamma)}(\overset{*}{g}) + \overset{*}{\Theta}_{(\xi,\gamma)}(\overset{*}{h})$$

for every $\xi \in \overline{\xi}(\gamma)$.

On the other hand, according to the definition of the addition, +, in $\overset{*}{G}(\gamma)$ given in Part 1, $\overset{*}{g} + \overset{*}{h}$ is the only element in $\overset{*}{G}(\gamma)$ such that

$$\overset{*}{\Theta}_{(\xi,\gamma)}(\overset{*}{g} + \overset{*}{h}) = \overset{*}{\Theta}_{(\xi,\gamma)}(\overset{*}{g}) + \overset{*}{\Theta}_{(\xi,\gamma)}(\overset{*}{h})$$

for every $\xi \in \overline{\xi}(\gamma)$. Therefore,

$$\overset{*}{g} \boxplus \overset{*}{h} = \overset{*}{g} + \overset{*}{h}.$$

■

## 5.31 The Equivalence of Axiomatics A and B

As we know, two axiomatics are said to be equivalent when the primitive terms of any of them are primitive or defined terms of the other and, also, the axioms of any of them are axioms or theorems of the other. Regarding axiomatics A and B, formulated, respectively, in items 5.21 and 5.28, an inspection in the list of its primitive and defined terms reveals that:

**(a)** all primitive terms of axiomatic B, namely, $\overline{\mathscr{G}}(I)$-distribution, domain, derivative and restriction of $\overline{\mathscr{G}}(I)$-distributions, also figure as primitive terms in axiomatic A;

**(b)** all primitive terms of axiomatic A, except for the term "addition of $\overline{\mathscr{G}}(I)$-distributions", also are primitive terms in axiomatic B; regarding the term "addition" of axiomatic A, this one, as we saw in the proof of Proposition 5.30, can be introduced in axiomatic B through the definition formulated in Part 1 of the referred proof, then being a defined term of axiomatic B.

Regarding the axioms of the axiomatics in question, also by inspection, we have:

**(c)** except for Axiom 5 and the ones related with the primitive term "addition" of axiomatic A, namely, Axiom 3, Axiom 4((4-1)) and Axiom 6((6-1)), all others also are axioms of axiomatic B.

With respect to Axiom 3, Axiom 4((4-1)) and Axiom 6((6-1)) of axiomatic A, these, as Proposition 5.30 shows us, are logical consequences, that is, theorems of axiomatic B equipped with the term "addition" (introduced into it through the definition formulated in Part 1 of the proof of the referred proposition).

Also Axiom 5 of axiomatic A, namely,

- If $g \in G(\gamma)$ and $g \notin G(\gamma)_\Phi$, then, $\overset{*}{\Phi}(g) \notin G(\gamma)$,

is a theorem of axiomatic B. In fact, one of the logical consequences of axiomatic B, that in 5.29(j), establishes that:

- The $S$-group $\mathring{\mathbb{G}}(\gamma) = (\mathring{G}(\gamma), \mathring{H}(\gamma))$ is a strict and closed extension of the $S$-group $\mathbb{G}(\gamma) = (G(\gamma), H(\gamma))$.

Thus, with $\mathring{\mathbb{G}}(\gamma)$ as strict extension of $\mathbb{G}(\gamma)$, we have that:

- If $g \in G(\gamma)$ and $g \notin G(\gamma)_\Phi$, then, $\mathring{\Phi}(g) \notin G(\gamma)$.

On the other hand, $G(\gamma) \subseteq \mathring{G}(\gamma)$ (5.29(a)) and $\mathring{\Phi} := \overset{*}{\Phi}|_{\mathring{G}(\gamma)}$, which allow us to conclude that $\mathring{\Phi}(g) = \overset{*}{\Phi}(g)$ for every $g \in G(\gamma)$. Therefore, we have as a theorem of axiomatic B that:

- If $g \in G(\gamma)$ and $g \notin G(\gamma)_\Phi$, then, $\overset{*}{\Phi}(g) \notin G(\gamma)$,

which is, exactly, the Axiom 5 of axiomatic A.

In short, all axioms of axiomatic A are axioms or theorems of axiomatic B.

**(d)** Finally, we remark that the only axiom of axiomatic B that does not figure among the axioms of axiomatic A is the fourth one, namely:

- For $\gamma \in \Gamma(I)$, $g, h \in G(\gamma)$ and $\Phi \in H(\gamma)$ arbitrarily fixed, one has:

$$\overset{*}{\Phi}(g) = \overset{*}{\Phi}(h) \quad \text{if and only if} \quad g - h \in N(\Phi).$$

However, as proved ahead, the statement above is a theorem of axiomatic A. In fact, among the logical consequences of axiomatic A obtained in item 5.23, we highlight that in 5.23(l), i.e.:

- For each $\gamma \in \Gamma(I)$, the ordered pair

$$\mathring{\mathbb{G}}(\gamma) := \left(\mathring{G}(\gamma), \mathring{H}(\gamma)\right)$$

  is a $S$-group that is a strict and closed extension of the $S$-group

$$\mathbb{G}(\gamma) = \left\{G(\gamma), H(\gamma)\right\}.$$

Since $\mathbb{G}(\gamma) = (G(\gamma), H(\gamma))$ is a surjective $S$-group (remember that the underlying $S$-space $\mathscr{G}(I)$ of the axiomatic A is abelian, **surjective**, coherent and with identity), then, by Proposition 1.16(b),

$$N(\mathring{\Phi}) = N(\Phi),$$

from where it results, for $g, h \in G(\gamma)$, that

$$g - h \in N(\mathring{\Phi}) \quad \text{if and only if} \quad g - h \in N(\Phi),$$

that is,

$$\mathring{\Phi}(g) = \mathring{\Phi}(h) \quad \text{if and only if} \quad g - h \in N(\Phi).$$

But, $\mathring{\Phi}(g) = \overset{*}{\Phi}(g)$ and $\mathring{\Phi}(h) = \overset{*}{\Phi}(h)$ since $\mathring{\Phi} := \overset{*}{\Phi}|_{\mathring{G}(\gamma)}$ and $G(\gamma) \subseteq \mathring{G}(\gamma)$. Hence, we have

$$\overset{*}{\Phi}(g) = \overset{*}{\Phi}(h) \quad \text{if and only if} \quad g - h \in N(\Phi),$$

which is what Axiom 4 of axiomatic A establishes.

In short, the considerations made in (a) to (d) above allow us to conclude that axiomatics A and B are equivalent.

# Schwartz' Distributions Revisited

## 5.32 Preliminaries

The items composing this section have as goals to prove that:

**(a)** the finite order distributions[29] and its derivatives are, essentially, the $\widetilde{\mathscr{C}}(\mathbb{R}^n)$-distributions and its corresponding derivatives, respectively;

**(b)** the distributions (of all orders, finite or infinite) and its derivatives are, essentially, the $\overline{\mathscr{C}}(\mathbb{R}^n)$-distributions and its corresponding derivatives, respectively.

In this endeavor, we will widely use the notions, concepts, and results present in the section "Schwartz' Distributions" in Chapter 2, more precisely the content of items 2.27 to 2.37, which, probably, will demand the reader a review of this material.

Also, the $S$-space of continuous functions,

$$\mathscr{C}(\mathbb{R}^n) = \left(\mathbb{C}\Big(\Gamma(\mathbb{R}^n)\Big), j\Big(\Gamma^2(\mathbb{R}^n)\Big), \rho\Big(\Delta(\mathbb{R}^n)\Big)\right),$$

[29] In the remaining of this chapter, the term "distribution" refers to the concept of Schwartz' distribution (Definition 2.30(a)).

will be heavily present in what follows, and so we also suggest a review of its definition in item 3.16.

Our procedure to reach the goals (a) and (b) above described consists of, roughly speaking, the execution of the "plan" described in 5.3(e), also recommended to be reviewed by the reader.

## 5.33 The Families $\mathbb{D}'(\Gamma(\mathbb{R}^n))$ and $\mathbb{D}'_f(\Gamma(\mathbb{R}^n))$

As we know, for each $\Omega \in \Gamma(\mathbb{R}^n)$ ($\Gamma(\mathbb{R}^n) = \{\Omega \subseteq \mathbb{R}^n : \Omega \text{ is an open set of } \mathbb{R}^n\}$), the $S$-group of the distributions of domain $\Omega$ (2.36),

$$\mathbb{D}'(\Omega) = \Big(D'(\Omega), D(\Omega) = \Big\{D^\alpha_\Omega : \alpha \in \mathbb{N}^n\Big\}\Big),$$

where $D'(\Omega)$ (2.31) is the abelian group of the distributions of domain $\Omega$ and $D(\Omega)$ (2.34) the semigroup of its derivatives, is an extension of the $S$-group

$$\mathbb{D}'_f(\Omega) = \Big(D'_f(\Omega), D_f(\Omega) = \Big\{D^\alpha_{\Omega-f} : \alpha \in \mathbb{N}^n\Big\}\Big),$$

of the finite order distributions of domain $\Omega$ (2.36), where $D'_f(\Omega)$ (2.31) is the subgroup, of the group $D'(\Omega)$, of the finite order distributions of domain $\Omega$, and $D_f(\Omega)$ (2.35) is the semigroup of its derivatives. We also know (2.37) that $\mathbb{D}'_f(\Omega)$ is a strict and closed extension of the $S$-group

$$\mathbb{C}(\Omega) = \Big(C(\Omega), \partial(\Omega) = \Big\{\partial^\alpha_\Omega : \alpha \in \mathbb{N}^n\Big\}\Big)$$

of the continuous functions on $\Omega$.

Hence, the families $\mathbb{D}'(\Gamma(\mathbb{R}^n))$ and $\mathbb{D}'_f(\Gamma(\mathbb{R}^n))$, indexed by $\Gamma(\mathbb{R}^n)$, defined by

$$\mathbb{D}'\Big(\Gamma(\mathbb{R}^n)\Big) := \Big\{\mathbb{D}'(\Omega) = \Big(D'(\Omega), D(\Omega)\Big)\Big\}_{\Omega\in\Gamma(\mathbb{R}^n)}$$

and

$$\mathbb{D}'_f\Big(\Gamma(\mathbb{R}^n)\Big) := \Big\{\mathbb{D}'_f(\Omega) = \Big(D'_f(\Omega), D_f(\Omega)\Big)\Big\}_{\Omega\in\Gamma(\mathbb{R}^n)}$$

are families of $S$-groups (Definition 3.7) such that: $\mathbb{D}'(\Gamma(\mathbb{R}^n))$ is an extension (Definition 3.8) of the family $\mathbb{D}'_f(\Gamma(\mathbb{R}^n))$ that, in turn, is a strict and closed extension (Definition 3.8) of the family of $S$-groups

$$\mathbb{C}\Big(\Gamma(\mathbb{R}^n)\Big) = \Big\{\mathbb{C}(\Omega) = \Big(C(\Omega), \partial(\Omega)\Big)\Big\}_{\Omega\in\Gamma(\mathbb{R}^n)}$$

of the $S$-space $\mathscr{C}(\mathbb{R}^n)$.

## 5.34 The Bonded Families $(\mathbb{D}'(\Gamma(\mathbb{R}^n)), j'(\Gamma^2(\mathbb{R}^n)))$ and $(\mathbb{D}'_f(\Gamma(\mathbb{R}^n)), j^{(f)}(\Gamma^2(\mathbb{R}^n)))$

The statement that, for each $\Omega \in \Gamma(\mathbb{R}^n)$, the $S$-group

$$\mathbb{D}'_f(\Omega) = \Big(D'_f(\Omega), D_f(\Omega) = \Big\{D^\alpha_{\Omega-f} : \alpha \in \mathbb{N}^n\Big\}\Big),$$

is an extension of the $S$-group

$$\mathbb{C}(\Omega) = \Big(C(\Omega), \partial(\Omega) = \Big\{\partial^\alpha_\Omega : \alpha \in \mathbb{N}^n\Big\}\Big)$$

means, among other things (review items 2.35 and 2.36), that the semigroup $D_f(\Omega)$ is isomorphic to the semigroup $\partial(\Omega)$ and that the function ($\mathrm{d}_f$ in 2.35(d), here denoted by $\mathrm{d}_{\Omega-f}$)

$$\begin{aligned} \mathrm{d}_{\Omega-f} : \partial(\Omega) &\longrightarrow D_f(\Omega) \\ \partial^\alpha_\Omega &\longmapsto \mathrm{d}_{\Omega-f}(\partial^\alpha_\Omega) := D^\alpha_{\Omega-f} \end{aligned}$$

is an isomorphism.

Analogously, the statement that, for each $\Omega \in \Gamma(\mathbb{R}^n)$, the $S$-group

$$\mathbb{D}'(\Omega) = \Big(D'(\Omega), D(\Omega) = \Big\{D^\alpha_\Omega : \alpha \in \mathbb{N}^n\Big\}\Big),$$

is an extension of the $S$-group

$$\mathbb{D}'_f(\Omega) = \Big(D'_f(\Omega), D_f(\Omega) = \Big\{D^\alpha_{\Omega-f} : \alpha \in \mathbb{N}^n\Big\}\Big),$$

tells us, among other things, that the semigroups $D_f(\Omega)$ and $D(\Omega)$ are isomorphic and that the function ($e$ in 2.36(d), here denoted by $e_\Omega$)

$$\begin{aligned} e_\Omega : D_f(\Omega) &\longrightarrow D(\Omega) \\ D^\alpha_{\Omega-f} &\longmapsto e_\Omega(D^\alpha_{\Omega-f}) := D^\alpha_\Omega \end{aligned}$$

is an isomorphism.

Now, since the family of $S$-groups

$$\mathbb{D}'_f\Big(\Gamma(\mathbb{R}^n)\Big) := \Big\{\mathbb{D}'_f(\Omega) = \Big(D'_f(\Omega), D_f(\Omega)\Big)\Big\}_{\Omega\in\Gamma(\mathbb{R}^n)}$$

is an extension of the family of $S$-groups

$$\mathbb{C}\Big(\Gamma(\mathbb{R}^n)\Big) = \Big\{\mathbb{C}(\Omega) = \Big(C(\Omega), \partial(\Omega)\Big)\Big\}_{\Omega\in\Gamma(\mathbb{R}^n)},$$

one then has well-defined (according to Definition 3.22) the extension, to $\mathbb{D}'_f(\Gamma(\mathbb{R}^n))$, of the bonding $j(\Gamma^2(\mathbb{R}^n))$ of the bonded family

$$\Big(\mathbb{C}\Big(\Gamma(\mathbb{R}^n)\Big), j\Big(\Gamma^2(\mathbb{R}^n)\Big)\Big)$$

of the $S$-space

$$\mathscr{C}(\mathbb{R}^n) = \Big(\mathbb{C}\Big(\Gamma(\mathbb{R}^n)\Big), j\Big(\Gamma^2(\mathbb{R}^n)\Big), \rho\Big(\Delta(\mathbb{R}^n)\Big)\Big).$$

Denoting by

$$j^{(f)}\Big(\Gamma^2(\mathbb{R}^n)\Big) = \Big\{j^{(f)}_{(\Omega',\Omega)}\Big\}_{(\Omega',\Omega)\in\Gamma^2(\mathbb{R}^n)}$$

the referred extension of the bonding $j(\Gamma^2(\mathbb{R}^n))$, one has, taking Definition 3.22 into account, that the members $j^{(f)}_{(\Omega',\gamma)}$ of this family are the following functions:

$$\begin{aligned} j^{(f)}_{(\Omega',\Omega)} : D_f(\Omega) &\longrightarrow D_f(\Omega') \\ D^\alpha_{\Omega-f} &\longmapsto j^{(f)}_{(\Omega',\Omega)}(D^\alpha_{\Omega-f}) \coloneqq \mathrm{d}_{\Omega'-f}\Big(j_{(\Omega',\Omega)}\Big(\mathrm{d}^{-1}_{\Omega-f}(D^\alpha_{\Omega-f})\Big)\Big) \end{aligned}$$

that is, considering the definitions of the isomorphisms $\mathrm{d}_{\Omega-f}$, $\mathrm{d}_{\Omega'-f}$ and $j_{(\Omega',\Omega)}$, one obtains:

$$\begin{aligned} j^{(f)}_{(\Omega',\Omega)} : D_f(\Omega) &\longrightarrow D_f(\Omega') \\ D^\alpha_{\Omega-f} &\longmapsto j^{(f)}_{(\Omega',\Omega)}(D^\alpha_{\Omega-f}) = D^\alpha_{\Omega'-f}. \end{aligned}$$

We have yet, as Definition 3.22 establishes, that the ordered pair

$$\Big(\mathbb{D}'_f\Big(\Gamma(\mathbb{R}^n)\Big), j^{(f)}\Big(\Gamma^2(\mathbb{R}^n)\Big)\Big)$$

is a bonded family and also an extension of the bonded family

$$\Big(\mathbb{C}\Big(\Gamma(\mathbb{R}^n)\Big), j\Big(\Gamma^2(\mathbb{R}^n)\Big)\Big).$$

Now in possession of the bonded family

$$\Big(\mathbb{D}'_f\Big(\Gamma(\mathbb{R}^n)\Big), j^{(f)}\Big(\Gamma^2(\mathbb{R}^n)\Big)\Big)$$

and knowing, as we saw in 5.33, that the family of $S$-groups $\mathbb{D}'(\Gamma(\mathbb{R}^n))$ is an extension of the family of $S$-groups $\mathbb{D}'_f(\Gamma(\mathbb{R}^n))$, we conclude also in this case that the extension, to $\mathbb{D}'(\Gamma(\mathbb{R}^n))$, of the bonding $j^{(f)}(\Gamma^2(\mathbb{R}^n))$ of the bonded family

$$\Big(\mathbb{D}'_f\Big(\Gamma(\mathbb{R}^n)\Big), j^{(f)}\Big(\Gamma^2(\mathbb{R}^n)\Big)\Big),$$

is well-defined (again according to Definition 3.22).

With the notation

$$j'\Big(\Gamma^2(\mathbb{R}^n)\Big) = \Big\{j'_{(\Omega',\Omega)}\Big\}_{(\Omega',\Omega)\in\Gamma^2(\mathbb{R}^n)}$$

for the referred extension, we have, again appealing to Definition 3.22, that its members, $j'_{(\Omega,\Omega)}$, are the functions defined by

$$\begin{aligned} j'_{(\Omega',\Omega)} : D(\Omega) &\longrightarrow D(\Omega') \\ D_\Omega^\alpha &\longmapsto j'_{(\Omega',\Omega)}(D_\Omega^\alpha) \coloneqq e_{\Omega'}\left( j^{(f)}_{(\Omega',\Omega)}\left(e_\Omega^{-1}(D_\Omega^\alpha)\right)\right) \end{aligned}$$

or yet, with the definitions of the isomorphisms $e_\Omega$, $e_{\Omega'}$ and $j^{(f)}_{(\Omega',\Omega)}$ in mind,

$$\begin{aligned} j'_{(\Omega',\Omega)} : D(\Omega) &\longrightarrow D(\Omega') \\ D_\Omega^\alpha &\longmapsto j'_{(\Omega',\Omega)}(D_\Omega^\alpha) = D_{\Omega'}^\alpha. \end{aligned}$$

From this and with a new resort to Definition 3.22, we can say that the ordered pair

$$\left(\mathbb{D}'\Big(\Gamma(\mathbb{R}^n)\Big), j'\Big(\Gamma^2(\mathbb{R}^n)\Big)\right)$$

is a bonded family and also an extension of the bonded family

$$\left(\mathbb{D}'_f\Big(\Gamma(\mathbb{R}^n)\Big), j^{(f)}\Big(\Gamma^2(\mathbb{R}^n)\Big)\right).$$

## 5.35 The Restrictions $\Theta'(\Delta(\mathbb{R}^n))$ and $\Theta^f(\Delta(\mathbb{R}^n))$

According to Definition 2.28(b), being $K \subseteq \mathbb{R}^n$ a non-empty compact set,

$$C_K^\infty(\mathbb{R}^n) \coloneqq \Big\{\phi \in C^\infty(\mathbb{R}^n) : \text{ support of } \phi \subseteq K\Big\}.$$

Yet, from the Definition 2.28(c), being $\Omega \subseteq \mathbb{R}^n$ an open set ($\Omega \in \Gamma(\mathbb{R}^n)$),

$$C_\Omega^\infty(\mathbb{R}^n) \coloneqq \bigcup_{K\subseteq\Omega} C_K^\infty(\mathbb{R}^n) = \Big\{\phi \in C^\infty(\mathbb{R}^n) : \text{ support of } \phi \text{ is a compact included in } \Omega\Big\}.$$

Also, by Definition 2.30, a distribution of domain $\Omega \in \Gamma(\mathbb{R}^n)$ (or a distribution in $\Omega$), is a linear functional

$$\begin{aligned} \Lambda : C_\Omega^\infty(\mathbb{R}^n) &\longrightarrow C \\ \phi &\longmapsto \Lambda(\phi) \end{aligned}$$

such that: for each compact $K \subseteq \Omega$, there exist a constant $C_K < \infty$ and an integer $N_K \geqslant 0$, generally depending on $K$, such that

$$|\Lambda(\phi)| \leqslant C_K \max\Big\{\left|\left(\partial^\alpha_{\mathbb{R}^n}(\phi)\right)(x)\right| : x \in \Omega \text{ and } |\alpha| \leqslant N_K\Big\}$$

for every $\phi \in C_K^\infty(\mathbb{R}^n)$; a distribution $\Lambda$ is said to be of finite order if the condition above is attended for a choice of $N_K = N$ independent of $K$.

Taking now into account the definitions above remembered, one easily concludes that:

**(a)** If $\Omega', \Omega \in \Gamma(\mathbb{R}^n)$ are such that $\Omega' \subseteq \Omega$, i.e., for $(\Omega', \Omega) \in \Delta(\mathbb{R}^n)$, one has

$$C^\infty_{\Omega'}(\mathbb{R}^n) \subseteq C^\infty_{\Omega}(\mathbb{R}^n);$$

**(b)** If $\Lambda : C^\infty_\Omega(\mathbb{R}^n) \longrightarrow C$ is a distribution of domain $\Omega$ and $\Omega' \in \Gamma(\mathbb{R}^n)$ such that $(\Omega', \Omega) \in \Delta(\mathbb{R}^n)$, then, the restriction of the function $\Lambda$ to the subset $C^\infty_{\Omega'}(\mathbb{R}^n)$ of its domain $C^\infty_\Omega(\mathbb{R}^n)$, that is, the function $\Lambda|_{C^\infty_{\Omega'}(\mathbb{R}^n)}$ defined by

$$\begin{aligned} \Lambda|_{C^\infty_{\Omega'}(\mathbb{R}^n)} : C^\infty_{\Omega'}(\mathbb{R}^n) &\longrightarrow C \\ \phi &\longmapsto \Big(\Lambda|_{C^\infty_{\Omega'}(\mathbb{R}^n)}\Big)(\phi) := \Lambda(\phi), \end{aligned}$$

also is a distribution of domain $\Omega' \subseteq \Omega$. In other terms, we have: if $\Lambda \in D'(\Omega)$, then, $\Lambda|_{C^\infty_{\Omega'}(\mathbb{R}^n)} \in D'(\Omega')$ for every $\Omega' \in \Gamma(\mathbb{R}^n)$ such that $\Omega' \subseteq \Omega$.

**(c)** If $\Lambda \in D'(\Omega)$ is of finite order, then, for every $\Omega' \in \Gamma(\mathbb{R}^n)$ such that $\Omega' \subseteq \Omega$, $\Lambda|_{C^\infty_{\Omega'}(\mathbb{R}^n)}$ also is of finite order, that is, if $\Lambda \in D'_f(\Omega)$, then $\Lambda|_{C^\infty_{\Omega'}(\mathbb{R}^n)} \in D'_f(\Omega')$.

As we will see, the family of functions defined ahead is particularly important regarding our goals described in 5.32(a) and (b).

**Definition.** $\Theta'(\Delta(\mathbb{R}^n))$ is the family of functions indexed by the set $\Delta(\mathbb{R}^n)$, defined by:

$$\Theta'\Big(\Delta(\mathbb{R}^n)\Big) := \Big\{\Theta'_{(\Omega',\Omega)}\Big\}_{(\Omega',\Omega)\in\Delta(\mathbb{R}^n)}$$

where $\Theta'_{(\Omega',\Omega)}$, for each $(\Omega', \Omega) \in \Delta(\mathbb{R}^n)$, is the following function

$$\begin{aligned} \Theta'_{(\Omega',\Omega)} : D'(\Omega) &\longrightarrow D'(\Omega') \\ \Lambda &\longmapsto \Theta'_{(\Omega',\Omega)}(\Lambda) := \Lambda|_{C^\infty_{\Omega'}(\mathbb{R}^n)} \end{aligned}$$

(the result (b) above ensures that $\Theta'_{(\Omega',\Omega)}$ assumes values in $D'(\Omega')$).

Among the properties of the family $\Theta'(\Delta(\mathbb{R}^n))$ we highlight the ones presented below, where $\Omega, \Omega', \Omega'' \in \Gamma(\mathbb{R}^n)$ such that $\Omega'' \subseteq \Omega' \subseteq \Omega$ are arbitrarily chosen.

**(d)** $\Theta'_{(\Omega'',\Omega')}\Big(\Theta'_{(\Omega',\Omega)}(\Lambda)\Big) = \Theta'_{(\Omega'',\Omega)}(\Lambda)$ for every $\Lambda \in D'(\Omega)$.

This results directly from the definitions of the functions of the family $\Theta'(\Delta(\mathbb{R}^n))$.

**(e)** $\Theta'_{(\Omega',\Omega)}(\Lambda) \in D'_f(\Omega')$ for every $\Lambda \in D'_f(\Omega)$.

Follows from (c) above and the definition of $\Theta'_{(\Omega',\Omega)}$.

**(f)** $\Theta'_{(\Omega',\Omega)}(\Lambda_f) = \Lambda_{\rho_{(\Omega',\Omega)}(f)}$ for every $f \in C(\Omega)$.

In fact, from the definition of distribution induced by a function $f \in C(\Omega)$, formulated in 2.31 (page 70), we have that,

$$\begin{aligned} \Lambda_f : C^\infty_\Omega(\mathbb{R}^n) &\longrightarrow C \\ \phi &\longmapsto \Lambda_f(\phi) := \textstyle\int_\Omega f(x)\phi(x)\,\mathrm{d}x \end{aligned}$$

and

$$\begin{aligned} \Lambda_{\rho_{(\Omega',\Omega)}(f)} : C^\infty_{\Omega'}(\mathbb{R}^n) &\longrightarrow C \\ \phi &\longmapsto \Lambda_{\rho_{(\Omega',\Omega)}(f)}(\phi) = \textstyle\int_{\Omega'} \big(\rho_{(\Omega',\Omega)}(f)\big)(x)\phi(x)\,\mathrm{d}x\,. \end{aligned}$$

On the other hand, by the definition of $\Theta'_{(\Omega',\Omega)}$,

$$\Theta'_{(\Omega',\Omega)}(\Lambda_f) = \Lambda_f|_{C^\infty_{\Omega'}(\mathbb{R}^n)}$$

and, hence, for each $\phi \in C^\infty_{\Omega'}(\mathbb{R}^n)$,

$$\Big(\Theta'_{(\Omega',\Omega)}(\Lambda_f)\Big)(\phi) = \Big(\Lambda_f|_{C^\infty_{\Omega'}(\mathbb{R}^n)}\Big)(\phi) = \Lambda_f(\phi) = \int_\Omega f(x)\phi(x)\,\mathrm{d}x\,.$$

Now, since $\phi \in C^\infty_{\Omega'}(\mathbb{R}^n)$, then, $\phi(x) = 0$ for $x \notin \Omega'$ and, hence, since $\Omega' \subseteq \Omega$, we obtain that

$$\Big(\Theta'_{(\Omega',\Omega)}(\Lambda_f)\Big)(\phi) = \int_{\Omega'} f(x)\phi(x)\,\mathrm{d}x = \int_{\Omega'} \Big(\rho_{(\Omega',\Omega)}(f)\Big)(x)\phi(x)\,\mathrm{d}x = \Lambda_{\rho_{(\Omega',\Omega)}(f)}(\phi)$$

for every $\phi \in C^\infty_{\Omega'}(\mathbb{R}^n)$, that is,

$$\Theta'_{(\Omega',\Omega)}(\Lambda_f) = \Lambda_{\rho_{(\Omega',\Omega)}(f)}.$$

**Remark.** As we saw in item 2.31 (page 70), for each $\Omega \in \Gamma(\mathbb{R}^n)$, the function $\eta_\Omega$ (there, in 2.31, denoted by $\eta$) defined by

$$\begin{aligned} \eta_\Omega : C(\Omega) &\longrightarrow D'_f(\Omega) \subseteq D'(\Omega) \\ f &\longmapsto \eta_\Omega(f) := \Lambda_f, \end{aligned}$$

is an isomorphism from the group $C(\Omega)$ onto the group $C'(\Omega) = \{\Lambda_f : f \in C(\Omega)\} \subseteq D'_f(\Omega)$ of the distributions in $\Omega$ induced by the functions $f \in C(\Omega)$. Hence, the property (f) above can be written as follows:

$$\Theta'_{(\Omega',\Omega)}\Big(\eta_\Omega(f)\Big) = \eta_{\Omega'}\Big(\rho_{(\Omega',\Omega)}(f)\Big)$$

or yet, with the identifications

$$f \equiv \eta_\Omega(f) \quad \text{and} \quad \rho_{(\Omega',\Omega)}(f) \equiv \eta_{\Omega'}\Big(\rho_{(\Omega',\Omega)}(f)\Big)$$

promoted by the isomorphisms $\eta_\Omega$ and $\eta_{\Omega'}$, and using the habitual abuse of language (referred, for instance, in 3.24), as

$$\Theta'_{(\Omega',\Omega)}(f) = \rho_{(\Omega',\Omega)}(f) \quad \text{for every} \quad f \in C(\Omega).$$

**(g)** $\Theta'_{(\Omega',\Omega)}$ is a homomorphism from the group $D'(\Omega)$ into the group $D'(\Omega')$.

This results from the fact that, for $\Lambda^{(1)}, \Lambda^{(2)} \in D'(\Omega)$ arbitrarily fixed,

$$(\Lambda^{(1)} + \Lambda^{(2)})|_{C^\infty_{\Omega'}(\mathbb{R}^n)} = \Lambda^{(1)}|_{C^\infty_{\Omega'}(\mathbb{R}^n)} + \Lambda^{(2)}|_{C^\infty_{\Omega'}(\mathbb{R}^n)}.$$

**(h)** $\Theta'_{(\Omega',\Omega)}\Big(D^\alpha_\Omega(\Lambda)\Big) = j'_{(\Omega',\Omega)}(D^\alpha_\Omega)\Big(\Theta'_{(\Omega',\Omega)}(\Lambda)\Big)$ for every $D^\alpha_\Omega \in D(\Omega)$ and every $\Lambda \in D'(\Omega)$.

In fact, from the definition of distribution derivative (Definition 2.33) we have, for $\alpha \in \mathbb{N}^n$, $\Omega \in \Gamma(\mathbb{R}^n)$ and $\Lambda \in D'(\Omega)$, that $D^\alpha_\Omega(\Lambda)$ is the distribution (in $\Omega$) given by:

$$\begin{aligned} D^\alpha_\Omega(\Lambda) : C^\infty_\Omega(\mathbb{R}^n) &\longrightarrow C \\ \phi &\longmapsto \Big(D^\alpha_\Omega(\Lambda)\Big)(\phi) = (-1)^{|\alpha|}\Lambda\Big(\partial^\alpha_{\mathbb{R}^n}(\phi)\Big). \end{aligned}$$

The definition of $\Theta'_{(\Omega',\Omega)}$, in turn, tells us that

$$\Theta'_{(\Omega',\Omega)}\Big(D^\alpha_\Omega(\Lambda)\Big) = D^\alpha_\Omega(\Lambda)|_{C^\infty_{\Omega'}(\mathbb{R}^n)},$$

and, therefore

$$\left(\Theta'_{(\Omega',\Omega)}\Big(D^\alpha_\Omega(\Lambda)\Big)\right)(\phi) = \Big(D^\alpha_\Omega(\Lambda)\Big)(\phi) = (-1)^{|\alpha|}\Lambda\Big(\partial^\alpha_{\mathbb{R}^n}(\phi)\Big)$$

for every $\phi \in C^\infty_{\Omega'}(\mathbb{R}^n)$.

Now, as we know, $\Theta'_{(\Omega',\Omega)}(\Lambda) \in D'(\Omega')$ and, hence, resorting to the definition of $D^\alpha_{\Omega'}$, we have

$$\begin{aligned} \left(D^\alpha_{\Omega'}\Big(\Theta'_{(\Omega',\Omega)}(\Lambda)\Big)\right)(\phi) &= (-1)^{|\alpha|}\Big(\Theta'_{(\Omega',\Omega)}(\Lambda)\Big)\Big(\partial^\alpha_{\mathbb{R}^n}(\phi)\Big) = \\ &= (-1)^{|\alpha|}\Big(\Lambda|_{C^\infty_{\Omega'}(\mathbb{R}^n)}\Big)\Big(\partial^\alpha_{\mathbb{R}^n}(\phi)\Big) = \\ &= (-1)^{|\alpha|}\Lambda\Big(\partial^\alpha_{\mathbb{R}^n}(\phi)\Big) \end{aligned}$$

for every $\phi \in C^\infty_{\Omega'}(\mathbb{R}^n)$.

From the last two expressions, we obtain that

$$\left(\Theta'_{(\Omega',\Omega)}\Big(D^\alpha_\Omega(\Lambda)\Big)\right)(\phi) = \left(D^\alpha_{\Omega'}\Big(\Theta'_{(\Omega',\Omega)}(\Lambda)\Big)\right)(\phi)$$

for every $\phi \in C^\infty_{\Omega'}(\mathbb{R}^n)$, i.e., that

$$\Theta'_{(\Omega',\Omega)}\Big(D^\alpha_\Omega(\Lambda)\Big) = D^\alpha_{\Omega'}\Big(\Theta'_{(\Omega',\Omega)}(\Lambda)\Big).$$

However, as we saw in 5.34 (page 308),

$$D^{\alpha}_{\Omega'} = j'_{(\Omega',\Omega)}(D^{\alpha}_{\Omega}).$$

Therefore,

$$\Theta'_{(\Omega',\Omega)}\Big(D^{\alpha}_{\Omega}(\Lambda)\Big) = j'_{(\Omega',\Omega)}(D^{\alpha}_{\Omega})\Big(\Theta'_{(\Omega',\Omega)}(\Lambda)\Big).$$

**(i)** The family $\Theta'(\Delta(\mathbb{R}^n))$ is a restriction for the bonded family

$$\left(\mathbb{D}'\Big(\Gamma(\mathbb{R}^n)\Big), j'\Big(\Gamma^2(\mathbb{R}^n)\Big)\right).$$

In fact, reporting ourselves to the definition of restriction for a bonded family, Definition 3.13, we see that conditions (a), (b) and (c) there presented correspond, exactly and respectively, to properties (g), (d) and (h) above of the family $\Theta'(\Delta(\mathbb{R}^n))$. Regarding condition (d) of the referred definition, this one is trivially satisfied by the family $\Theta'(\Delta(\mathbb{R}^n))$, once that the derivatives $D^{\alpha}_{\Omega} \in D(\Omega)$ (corresponding, in this case, to the homomorphisms $\Phi \in H(\gamma)$ referred to in condition (d)) are endomorphisms on the group $D'(\Omega)$.

Let us now consider, for each $\Omega \in \Gamma(\mathbb{R}^n)$, the subgroup of the abelian group $D'(\Omega)$, constituted by the finite order distributions (in $\Omega$):

$$D'_f(\Omega) \subseteq D'(\Omega).$$

The property (e) tells us that, for each $\Omega' \in \Gamma(\mathbb{R}^n)$ such that $(\Omega',\Omega) \in \Delta(\mathbb{R}^n)$, the restriction (in the usual sense of restriction of a function to a subset of its domain) of the restriction (in the sense of being member of the family $\Theta'(\Delta(\mathbb{R}^n))$) $\Theta'_{(\Omega',\Omega)} : D'(\Omega) \longrightarrow D'(\Omega)$ to $D'_f(\Omega)$, i.e., the function

$$\begin{aligned} \Theta'_{(\Omega',\Omega)}|_{D'_f(\Omega)} : D'_f(\Omega) &\longrightarrow D'(\Omega') \\ \Lambda &\longmapsto \Big(\Theta'_{(\Omega',\Omega)}|_{D'_f(\Omega)}\Big)(\Lambda) := \Theta'_{(\Omega',\Omega)}(\Lambda), \end{aligned}$$

assumes values in $D'_f(\Omega')$, that is,

$$\Big(\Theta'_{(\Omega',\Omega)}|_{D'_f(\Omega)}\Big)(\Lambda) \in D'_f(\Omega') \quad \text{for every} \quad \Lambda \in D'_f(\Omega).$$

Thus, for each function $\Theta'_{(\Omega',\Omega)} : D'(\Omega) \longrightarrow D'(\Omega')$ of the family $\Theta'(\Delta(\mathbb{R}^n))$, the corresponding function

$$\Theta'_{(\Omega',\Omega)}|_{D'_f(\Omega)} : D'_f(\Omega) \longrightarrow D'_f(\Omega')$$

only involves finite order distributions. We hence obtain the family $\Theta^f(\Delta(\mathbb{R}^n))$ defined ahead, whose members are functions that associate finite order distributions with finite order distributions.

**Definition.** $\Theta^f(\Delta(\mathbb{R}^n))$ is the family of functions indexed by the set $\Delta(\mathbb{R}^n)$ and defined by:

$$\Theta^f\Big(\Delta(\mathbb{R}^n)\Big) := \Big\{\Theta^f_{(\Omega',\Omega)}\Big\}_{(\Omega',\Omega)\in\Delta(\mathbb{R}^n)}$$

where $\Theta^f_{(\Omega',\Omega)}$, for each $(\Omega',\Omega) \in \Delta(\mathbb{R}^n)$, is the following function

$$\begin{aligned}\Theta^f_{(\Omega',\Omega)} : D'_f(\Omega) &\longrightarrow D'_f(\Omega')\\ \Lambda &\longmapsto \Theta^f_{(\Omega',\Omega)}(\Lambda) := \Theta'_{(\Omega',\Omega)}(\Lambda),\end{aligned}$$

that is,

$$\Theta^f_{(\Omega',\Omega)} := \Theta'_{(\Omega',\Omega)}|_{D'_f(\Omega)}.$$

Taking into account the definition above, as well as the one for the semigroup $D_f(\Omega) \subseteq D(\Omega)$ (in 2.35), namely,

$$D_f(\Omega) := \Big\{D^\alpha_{\Omega-f} : \alpha \in \mathbb{N}^n\Big\},$$

where

$$D^\alpha_{\Omega-f} := D^\alpha_\Omega|_{D'_f(\Omega)},$$

that is,

$$\begin{aligned}D^\alpha_{\Omega-f} : D'_f(\Omega) &\longrightarrow D'_f(\Omega)\\ \Lambda &\longmapsto D^\alpha_{\Omega-f}(\Lambda) := D^\alpha_\Omega(\Lambda),\end{aligned}$$

the reader will find no trouble to verify that, from the properties (g), (d) and (h) above regarding the family $\Theta'(\Delta(\mathbb{R}^n))$ result, respectively, the properties (j), (k) and (l) below, concerning the family $\Theta^f(\Delta(\mathbb{R}^n))$.

**(j)** $\Theta^f_{(\Omega',\Omega)}$, for each $(\Omega',\Omega) \in \Delta(\mathbb{R}^n)$, is a homomorphism from the group $D'_f(\Omega)$ into the group $D'_f(\Omega')$;

**(k)** $\Theta^f_{(\Omega'',\Omega')}(\Theta^f_{(\Omega',\Omega)}(\Lambda)) = \Theta^f_{(\Omega'',\Omega)}(\Lambda)$ for every $\Omega, \Omega', \Omega'' \in \Gamma(\mathbb{R}^n)$ such that $\Omega'' \subseteq \Omega' \subseteq \Omega$ and every $\Lambda \in D'_f(\Omega)$;

**(l)** $\Theta^f_{(\Omega',\Omega)}(D^\alpha_{\Omega-f}(\Lambda)) = j^{(f)}_{(\Omega',\Omega)}(D^\alpha_{\Omega-f})(\Theta^f_{(\Omega',\Omega)}(\Lambda))$ for every $(\Omega',\Omega) \in \Delta(\mathbb{R}^n)$ and every $\Lambda \in D'_f(\Omega)$.

These three properties, (j), (k) and (l), of the family $\Theta^f(\Delta(\mathbb{R}^n))$, are, respectively, the first three conditions required by Definition 3.13 for this family to be a restriction for the bonded family

$$\bigg(\mathbb{D}'_f\Big(\Gamma(\mathbb{R}^n)\Big), j^{(f)}\Big(\Gamma^2(\mathbb{R}^n)\Big)\bigg).$$

The fourth and last condition required in the quoted definition is, such as for the case of the family $\Theta'(\Delta(\mathbb{R}^n))$, also trivially satisfied here. In short, we have:

**(m)** the family $\Theta^f(\Delta(\mathbb{R}^n))$ is a restriction for the bonded family

$$\left(\mathbb{D}'_f\Big(\Gamma(\mathbb{R}^n)\Big), j^{(f)}\Big(\Gamma^2(\mathbb{R}^n)\Big)\right).$$

## 5.36 The $S$-Space $\mathscr{D}'_f(\mathbb{R}^n)$

In item 5.34 we proved that

$$\left(\mathbb{D}'_f\Big(\Gamma(\mathbb{R}^n)\Big), j^{(f)}\Big(\Gamma^2(\mathbb{R}^n)\Big)\right)$$

is a bonded family that is an extension of the bonded family

$$\left(\mathbb{C}\Big(\Gamma(\mathbb{R}^n)\Big), j\Big(\Gamma^2(\mathbb{R}^n)\Big)\right)$$

of the $S$-space

$$\mathscr{C}(\mathbb{R}^n) = \left(\mathbb{C}\Big(\Gamma(\mathbb{R}^n)\Big), j\Big(\Gamma^2(\mathbb{R}^n)\Big), \rho\Big(\Delta(\mathbb{R}^n)\Big)\right).$$

We also saw, in 5.35, that the family $\Theta^f(\Delta(\mathbb{R}^n))$ is a restriction for the bonded family

$$\left(\mathbb{D}'_f\Big(\Gamma(\mathbb{R}^n)\Big), j^{(f)}\Big(\Gamma^2(\mathbb{R}^n)\Big)\right).$$

Now, taking into account that $\Theta^f_{(\Omega',\Omega)} = \Theta'_{(\Omega',\Omega)}|_{D'_f(\Omega)}$ and that $\eta_\Omega(f) = \Lambda_f \in D'_f(\Omega)$ for every $f \in C(\Omega)$, we have

$$\Theta^f_{(\Omega',\Omega)}\Big(\eta_\Omega(f)\Big) = \Theta'_{(\Omega',\Omega)}\Big(\eta_\Omega(f)\Big) \quad \text{for every} \quad f \in C(\Omega).$$

But, by 5.35(f), $\Theta'_{(\Omega',\Omega)}(\Lambda_f) = \Lambda_{\rho_{(\Omega',\Omega)}(f)}$, that is, $\Theta'_{(\Omega',\Omega)}(\eta_\Omega(f)) = \eta_{\Omega'}(\rho_{(\Omega',\Omega)}(f))$ and, hence, we have that, for every $f \in C(\Omega)$,

$$\Theta^f_{(\Omega',\Omega)}\Big(\eta_\Omega(f)\Big) = \eta_{\Omega'}\Big(\rho_{(\Omega',\Omega)}(f)\Big)$$

or, in the usual abuse of language,

$$\Theta^f_{(\Omega',\Omega)}(f) = \rho_{(\Omega',\Omega)}(f) \quad \text{for every} \quad f \in C(\Omega).$$

Therefore, as Definition 3.23(a) establishes, the family $\Theta^f(\Delta(\mathbb{R}^n))$ is a prolongation of $\rho(\Delta(\mathbb{R}^n))$ to the bonded family

$$\left(\mathbb{D}'_f\Big(\Gamma(\mathbb{R}^n)\Big), j^{(f)}\Big(\Gamma^2(\mathbb{R}^n)\Big)\right).$$

In short, the conclusions above allow us to conclude, according to the definitions of $S$-space (Definition 3.15) and $S$-space extension (Definition 3.23(b)), that the triplet $\mathscr{D}'_f(\mathbb{R}^n)$ defined by

$$\mathscr{D}'_f(\mathbb{R}^n) := \left(\mathbb{D}'_f\Big(\Gamma(\mathbb{R}^n)\Big), j^{(f)}\Big(\Gamma^2(\mathbb{R}^n)\Big), \Theta^f\Big(\Delta(\mathbb{R}^n)\Big)\right)$$

is a $S$-space and also a strict and closed extension (Definition 3.23(c-1)) of the $S$-space

$$\mathscr{C}(\mathbb{R}^n) = \left(\mathbb{C}\Big(\Gamma(\mathbb{R}^n)\Big), j\Big(\Gamma^2(\mathbb{R}^n)\Big), \rho\Big(\Delta(\mathbb{R}^n)\Big)\right)$$

of continuous functions.

We will henceforth refer to $\mathscr{D}'_f(\mathbb{R}^n)$ as the $S$-space of the finite order distributions.

## 5.37 The $S$-Space $\mathscr{D}'(\mathbb{R}^n)$

As we know (from item 5.34), the ordered pair $(\mathbb{D}'(\Gamma(\mathbb{R}^n)), j'(\Gamma^2(\mathbb{R}^n)))$ is a bonded family that is an extension of the bonded family $(\mathbb{D}'_f(\Gamma(\mathbb{R}^n)), j^{(f)}(\Gamma^2(\mathbb{R}^n)))$. Furthermore, we also know (5.35(i)) that the family $\Theta'(\Delta(\mathbb{R}^n))$ is a restriction for the bonded family $(\mathbb{D}'(\Gamma(\mathbb{R}^n)), j'(\Gamma^2(\mathbb{R}^n)))$. We have yet, due to the definition of the members of the family $\Theta^f(\Delta(\mathbb{R}^n))$, that

$$\Theta'_{(\Omega',\Omega)}(\Lambda) = \Theta^f_{(\Omega',\Omega)}(\Lambda) \quad \text{for every} \quad \Lambda \in D'_f(\Omega),$$

which means, with Definition 3.23(a) in mind, that the family $\Theta'(\Delta(\mathbb{R}^n))$ is a prolongation of $\Theta^f(\Delta(\mathbb{R}^n))$ to the bonded family $(\mathbb{D}'(\Gamma(\mathbb{R}^n)), j'(\Gamma^2(\mathbb{R}^n)))$.

Thus, and keeping in mind the definitions of $S$-space and $S$-space extension, the triplet $\mathscr{D}'(\mathbb{R}^n)$ defined by

$$\mathscr{D}'(\mathbb{R}^n) := \left(\mathbb{D}'\Big(\Gamma(\mathbb{R}^n)\Big), j'\Big(\Gamma^2(\mathbb{R}^n)\Big), \Theta'\Big(\Delta(\mathbb{R}^n)\Big)\right)$$

is a $S$-space and also an extension of the $S$-space

$$\mathscr{D}'_f(\mathbb{R}^n) := \left(\mathbb{D}'_f\Big(\Gamma(\mathbb{R}^n)\Big), j^{(f)}\Big(\Gamma^2(\mathbb{R}^n)\Big), \Theta^f\Big(\Delta(\mathbb{R}^n)\Big)\right)$$

of the finite order distributions. We will say that $\mathscr{D}'(\mathbb{R}^n)$ is the $S$-space of the distributions.

We now inquire: Is the extension $\mathscr{D}'(\mathbb{R}^n)$ of the $S$-space $\mathscr{D}'_f(\mathbb{R}^n)$ locally closed and coherent?

In order to answer this question, we will resort to two results from Schwartz' Distributions Theory, stated below in A and B, whose proofs can be found, for instance, at the book by W. Rudin[30] referred to at the footnote 14 at page 69.

[30] At the referred W. Rudin's book, results A and B are, exactly, the theorems 6.21 and 6.26 translated to our notation and terminology.

**A** − *Let* $\overline{\Gamma}(\Omega) \subseteq \Gamma(\Omega)$ *be a cover of* $\Omega \in \Gamma(\mathbb{R}^n)$ *and* $\{\Lambda_{\Omega'}\}_{\Omega' \in \overline{\Gamma}(\Omega)}$ *be a family of distributions such that:* $\Lambda_{\Omega'} \in D'(\Omega')$ *for every* $\Omega' \in \overline{\Gamma}(\Omega)$ *and*

$$\Theta'_{(\Omega'' \cap \Omega', \Omega')}(\Lambda_{\Omega'}) = \Theta'_{(\Omega'' \cap \Omega', \Omega'')}(\Lambda_{\Omega''})$$

*for every* $\Omega', \Omega'' \in \overline{\Gamma}(\Omega)$ *such that* $\Omega'' \cap \Omega' \neq \varnothing$*. Then, there exists a single distribution* $\Lambda \in D'(\Omega)$ *such that*

$$\Theta'_{(\Omega', \Omega)}(\Lambda) = \Lambda_{\Omega'} \quad \textit{for every} \quad \Omega' \in \overline{\Gamma}(\Omega).$$

**B** − *Let* $\Lambda \in D'(\Omega)$ *and* $K \subseteq \Omega \in \Gamma(\mathbb{R}^n)$ *be a compact set. Then, there exist* $f \in C(\Omega)$ *and* $\alpha \in \mathbb{N}^n$ *such that*

$$\Lambda(\phi) = (-1)^{|\alpha|} \int_{\Omega} f(x)(\partial^{\alpha}_{\mathbb{R}^n}(\phi))(x)\, \mathrm{d}x$$

*for every* $\phi \in C^{\infty}_{K}(\mathbb{R}^n)$*.*

Result A tells us, exactly, in terms of Definition 3.19, that the $S$-space of the distributions, $\mathscr{D}'(\mathbb{R}^n)$, is coherent.

On the other hand, it results from B that $\mathscr{D}'(\mathbb{R}^n)$ is a locally closed extension of the $S$-space $\mathscr{D}'_f(\mathbb{R}^n)$. In order to prove this last claim, let $\Omega \in \Gamma(\mathbb{R}^n)$, $K \subseteq \Omega$ compact, and $\Lambda \in D'(\Omega)$ be all arbitrarily fixed. Let now $\Omega' \in \Gamma(\Omega)$ be such that $\Omega' \subseteq K$. Since

$$C^{\infty}_{\Omega'}(\mathbb{R}^n) \subseteq C^{\infty}_{K}(\mathbb{R}^n),$$

then, taking B into account, we have, for some $f \in C(\Omega)$ and $\alpha \in \mathbb{N}^n$, that

$$\Lambda(\phi) = (-1)^{|\alpha|} \int_{\Omega'} f(x)(\partial^{\alpha}_{\mathbb{R}^n}(\phi))(x)\, \mathrm{d}x$$

for every $\phi \in C^{\infty}_{\Omega'}(\mathbb{R}^n)$. But, since

$$\Lambda(\phi) = \Big(\Theta'_{(\Omega', \Omega)}(\Lambda)\Big)(\phi)$$

for $\phi \in C^{\infty}_{\Omega'}(\mathbb{R}^n)$ (once that $\Theta'_{(\Omega', \Omega)}(\Lambda) = \Lambda|_{C^{\infty}_{\Omega'}(\mathbb{R}^n)}$), comes that

$$\Big(\Theta'_{(\Omega', \Omega)}(\Lambda)\Big)(\phi) = (-1)^{|\alpha|} \int_{\Omega'} f(x)(\partial^{\alpha}_{\mathbb{R}^n}(\phi))(x)\, \mathrm{d}x \tag{5.37-1}$$

for every $\phi \in C^{\infty}_{\Omega'}(\mathbb{R}^n)$.

On the other hand, from the definition of derivative of a distribution (Definition 2.33), we have that

$$\left(D^{\alpha}_{\Omega'}\Big(\Theta'_{(\Omega', \Omega)}(\Lambda_f)\Big)\right)(\phi) = (-1)^{|\alpha|} \Big(\Theta'_{(\Omega', \Omega)}(\Lambda_f)\Big)\Big(\partial^{\alpha}_{\mathbb{R}^n}(\phi)\Big) \tag{5.37-2}$$

for every $\phi \in C^{\infty}_{\Omega'}(\mathbb{R}^n)$. But, according to the definitions of $\Theta'_{(\Omega',\Omega)}$ and $\Lambda_f$,

$$\Big(\Theta'_{(\Omega',\Omega)}(\Lambda_f)\Big)(\varphi) = \Big(\Lambda_f|_{C^{\infty}_{\Omega'}(\mathbb{R}^n)}\Big)(\varphi) = \Lambda_f(\varphi) = \int_{\Omega'} f(x)\varphi(x)\,\mathrm{d}x$$

for $\varphi \in C^{\infty}_{\Omega'}(\mathbb{R}^n)$ and, hence, since $\varphi = \partial^{\alpha}_{\mathbb{R}^n}(\phi) \in C^{\infty}_{\Omega'}(\mathbb{R}^n)$ if $\phi \in C^{\infty}_{\Omega'}(\mathbb{R}^n)$, we obtain

$$\Big(\Theta'_{(\Omega',\Omega)}(\Lambda_f)\Big)\Big(\partial^{\alpha}_{\mathbb{R}^n}(\phi)\Big) = \Lambda_f\Big(\partial^{\alpha}_{\mathbb{R}^n}(\phi)\Big) = \int_{\Omega'} f(x)(\partial^{\alpha}_{\mathbb{R}^n}(\phi))(x)\,\mathrm{d}x$$

for $\phi \in C^{\infty}_{\Omega'}(\mathbb{R}^n)$. Taking now this result into (5.37-2) one concludes that

$$\left(D^{\alpha}_{\Omega'}\Big(\Theta'_{(\Omega',\Omega)}(\Lambda_f)\Big)\right)(\phi) = (-1)^{|\alpha|}\int_{\Omega'} f(x)(\partial^{\alpha}_{\mathbb{R}^n}(\phi))(x)\,\mathrm{d}x$$

for every $\phi \in C^{\infty}_{\Omega'}(\mathbb{R}^n)$. This last equation, along with equation (5.37-1), show us that

$$\Big(\Theta'_{(\Omega',\Omega)}(\Lambda)\Big)(\phi) = \left(D^{\alpha}_{\Omega'}\Big(\Theta'_{(\Omega',\Omega)}(\Lambda_f)\Big)\right)(\phi)$$

for every $\phi \in C^{\infty}_{\Omega'}(\mathbb{R}^n)$, that is,

$$\Theta'_{(\Omega',\Omega)}(\Lambda) = D^{\alpha}_{\Omega'}\Big(\Theta'_{(\Omega',\Omega)}(\Lambda_f)\Big),$$

or yet, taking 5.35(f) into account, that

$$\Theta'_{(\Omega',\Omega)}(\Lambda) = D^{\alpha}_{\Omega'}\Big(\Lambda_{\rho_{\Omega',\Omega}(f)}\Big).$$

In short, we have then proved the following result:

- Let $\Omega \in \Gamma(\mathbb{R}^n)$, $K \subseteq \Omega$ compact, and $\Lambda \in D'(\Omega)$ be all arbitrarily fixed. For each $\Omega' \in \Gamma(\Omega)$ such that $\Omega' \subseteq K$, there exist $f \in C(\Omega)$ and $\alpha \in \mathbb{N}^n$ such that
$$\Theta'_{(\Omega',\Omega)}(\Lambda) = D^{\alpha}_{\Omega'}\Big(\Lambda_{\rho_{(\Omega',\Omega)}(f)}\Big).$$

Let now $\Omega \in \Gamma(\mathbb{R}^n)$ and $\Lambda \in D'(\Omega)$ be arbitrarily chosen. As we know, for each $x \in \Omega$, there exist $\Omega' \in \Gamma(\Omega)$ and a compact $K \subseteq \mathbb{R}^n$ such that

$$x \in \Omega' \quad \text{and} \quad \Omega' \subseteq K \subseteq \Omega.$$

Thus being, and remembering ourselves that $D^{\alpha}_{\Omega'}(\Lambda_{\rho_{(\Omega',\Omega)}(f)}) \in D'_f(\Omega')$ for every $f \in C(\Omega)$ and every $\alpha \in \mathbb{N}^n$, the result above allow us to conclude that

$$\Theta'_{(\Omega',\Omega)}(\Lambda) \in D'_f(\Omega')$$

and, therefore, in terms of Definition 3.23(c-2), that the $S$-space $\mathscr{D}'(\mathbb{R}^n)$ is a locally closed extension of the $S$-space $\mathscr{D}'_f(\mathbb{R}^n)$.

## 5.38 The $\widetilde{\mathscr{C}}(\mathbb{R}^n)$ and $\overline{\mathscr{C}}(\mathbb{R}^n)$-Distributions

The considerations of this item involve, in an essential way, the concepts of domain of distributions of 1st and 2nd species, as well as the notions associated with such domains, introduced, all of them, through Definition 5.10. Hence, we suggest the reader, before going forward, to review this definition.

Now, keeping Definition 5.10(b) in mind, we observe that the $S$-space of the continuous functions,

$$\mathscr{C}(\mathbb{R}^n) = \Big(\mathbb{C}\big(\Gamma(\mathbb{R}^n)\big), j\big(\Gamma^2(\mathbb{R}^n)\big), \rho\big(\Delta(\mathbb{R}^n)\big)\Big),$$

satisfies all the required conditions for the $S$-space $\mathscr{G}(I)$ of the referred definition, that is, $\mathscr{C}(\mathbb{R}^n)$ is an abelian, surjective, with identity, and coherent $S$-space, and its family of $S$-groups,

$$\mathbb{C}\big(\Gamma(\mathbb{R}^n)\big) = \Big\{\mathbb{C}(\Omega) = \big(C(\Omega), \partial(\Omega)\big)\Big\}_{\Omega\in\Gamma(\mathbb{R}^n)},$$

is such that its groups are two by two disjoints, that is,

$$C(\Omega) \cap C(\Omega') = \varnothing \quad \text{for} \quad \Omega \neq \Omega'.$$

This, then, allows us to specialize (particularize), that is, to translate the definition in question to the case where $\mathscr{G}(I)$ is taken as the $S$-space $\mathscr{C}(\mathbb{R}^n)$. Hence proceeding, we get that the ordered pair

$$\big(\mathscr{C}(\mathbb{R}^n), \widehat{\mathscr{C}}(\mathbb{R}^n)\big),$$

where

$$\widehat{\mathscr{C}}(\mathbb{R}^n) = \Big(\widehat{\mathbb{C}}\big(\Gamma(\mathbb{R}^n)\big) = \Big\{\widehat{C}(\Omega) = \big(\widehat{C}(\Omega), \widehat{\partial}(\Omega)\big)\Big\}_{\Omega\in\Gamma(\mathbb{R}^n)}, \widehat{j}\big(\Gamma^2(\mathbb{R}^n)\big), \widehat{\rho}\big(\Delta(\mathbb{R}^n)\big)\Big)$$

is a $S$-space, will be said:

- a domain of distribution of 1st species (determined by $\mathscr{C}(\mathbb{R}^n)$) if and only if $\widehat{\mathscr{C}}(\mathbb{R}^n)$ is isomorphic to the strict and closed extension, $\widetilde{\mathscr{C}}(\mathbb{R}^n)$, of $\mathscr{C}(\mathbb{R}^n)$, defined in the 1st TESS;
- a domain of distribution of 2nd species (determined by $\mathscr{C}(\mathbb{R}^n)$) if and only if $\widehat{\mathscr{C}}(\mathbb{R}^n)$ is isomorphic to the locally closed and coherent extension, $\overline{\mathscr{C}}(\mathbb{R}^n)$, described in the 2nd TESS, of $\widetilde{\mathscr{C}}(\mathbb{R}^n)$ (defined in the 1st TESS);

It then results that $(\mathscr{C}(\mathbb{R}^n), \widetilde{\mathscr{C}}(\mathbb{R}^n))$ and $(\mathscr{C}(\mathbb{R}^n), \overline{\mathscr{C}}(\mathbb{R}^n))$ are, respectively, domains of distribution of 1st and 2nd species, both determined by the $S$-space $\mathscr{C}(\mathbb{R}^n)$. On the other hand, from items 5.36 and 5.37, we know that:

- the $S$-space $\mathscr{D}'_f(\mathbb{R}^n)$ of the finite order distributions is a strict and closed extension of $\mathscr{C}(\mathbb{R}^n)$ and, therefore, by the 1st TESS, it is isomorphic to $\widetilde{\mathscr{C}}(\mathbb{R}^n)$;
- the $S$-space $\mathscr{D}'(\mathbb{R}^n)$ of the distributions is a locally closed and coherent extension of the strict and closed extension, $\mathscr{D}'_f(\mathbb{R}^n)$, of $\mathscr{C}(\mathbb{R}^n)$.

Let us consider now the Proposition 4.34 and observe that the results above mentioned show us that, with the choices

$$\begin{aligned} \mathscr{G}(I) &= \mathscr{C}(\mathbb{R}^n), \\ \widetilde{\mathscr{G}}(I) &= \widetilde{\mathscr{C}}(\mathbb{R}^n), \\ \overline{\mathscr{G}}(I) &= \overline{\mathscr{C}}(\mathbb{R}^n), \\ \overset{\circ}{\mathscr{G}}(I) &= \mathscr{D}'_f(\mathbb{R}^n) \quad \text{and} \\ \overset{*}{\mathscr{G}}(I) &= \mathscr{D}'(\mathbb{R}^n), \end{aligned}$$

all hypotheses of the referred proposition are satisfied, which allow us to affirm its thesis, namely,

$$\mathscr{D}'(\mathbb{R}^n) = \overset{*}{\mathscr{G}}(I) \simeq \overline{\mathscr{G}}(I) = \overline{\mathscr{C}}(\mathbb{R}^n),$$

that is, the $S$-space $\mathscr{D}'(\mathbb{R}^n)$ of the distributions is isomorphic to $\overline{\mathscr{C}}(\mathbb{R}^n)$.

Hence, since $\mathscr{D}'_f(\mathbb{R}^n) \simeq \widetilde{\mathscr{C}}(\mathbb{R}^n)$ and $\mathscr{D}'(\mathbb{R}^n) \simeq \overline{\mathscr{C}}(\mathbb{R}^n)$, we have that

$$\Big(\mathscr{C}(\mathbb{R}^n), \mathscr{D}'_f(\mathbb{R}^n)\Big) \quad \text{and} \quad \Big(\mathscr{C}(\mathbb{R}^n), \mathscr{D}'(\mathbb{R}^n)\Big)$$

are, also, domains of distribution of 1st and 2nd species, respectively, determined by the $S$-space $\mathscr{C}(\mathbb{R}^n)$.

In an obvious sense, we can say that the domains $(\mathscr{C}(\mathbb{R}^n), \widetilde{\mathscr{C}}(\mathbb{R}^n))$ and $(\mathscr{C}(\mathbb{R}^n), \mathscr{D}'_f(\mathbb{R}^n))$ are isomorphic, essentially equal, as well as the domains $(\mathscr{C}(\mathbb{R}^n), \overline{\mathscr{C}}(\mathbb{R}^n))$ and $(\mathscr{C}(\mathbb{R}^n), \mathscr{D}'(\mathbb{R}^n))$.

Considering now the Definitions 5.10(c) and (d), let us remark that the concepts associated to a domain of distribution, for instance, of 1st species (defined in part (c) of Definition 5.10), namely, $\widetilde{\mathscr{G}}(I)$-distribution, domain, derivative, and restriction of $\widetilde{\mathscr{G}}(I)$-distributions, are "defined unless isomorphism" in the following sense: a $\widetilde{\mathscr{C}}(\mathbb{R}^n)$-distribution, for instance, is, according to Definition 5.10(c-1), any member of any of the groups $\widehat{C}(\Omega)$ of the $S$-space $\widehat{\mathscr{C}}(\mathbb{R}^n)$ of the domain $(\mathscr{C}(\mathbb{R}^n), \widehat{\mathscr{C}}(\mathbb{R}^n))$, indistinctly if we take $\widehat{\mathscr{C}}(\mathbb{R}^n) = \widetilde{\mathscr{C}}(\mathbb{R}^n)$ (in which case $\widehat{C}(\Omega) = \widetilde{C}(\Omega)$) or $\widehat{\mathscr{C}}(\mathbb{R}^n) = \mathscr{D}'_f(\mathbb{R}^n)$ (where $\widehat{C}(\Omega) = D'_f(\Omega)$) or any other $S$-space $\widehat{\mathscr{C}}(\mathbb{R}^n)$ isomorphic to $\widetilde{\mathscr{C}}(\mathbb{R}^n)$; in the same way, derivative of $\widetilde{\mathscr{C}}(\mathbb{R}^n)$-distribution is defined in 5.10(c-3) as being any endomorphism of the semigroups $\widehat{\partial}(\Omega)$ of the $S$-space $\widehat{\mathscr{C}}(\mathbb{R}^n)$, without caring if we take $\widehat{\mathscr{C}}(\mathbb{R}^n) = \widetilde{\mathscr{C}}(\mathbb{R}^n)$ (where $\widehat{\partial}(\Omega) = \widetilde{\partial}(\Omega)$) or $\widehat{\mathscr{C}}(\mathbb{R}^n) = \mathscr{D}'_f(\mathbb{R}^n)$ (where $\widehat{\partial}(\Omega) = D_f(\Omega)$) or any other $S$-space isomorphic to $\widetilde{\mathscr{C}}(\mathbb{R}^n)$. In the same sense and by the same reasons, the corresponding concepts referring to a domain of 2nd species (these ones defined in part (d) of Definition 5.10) — $\overline{\mathscr{G}}(I)$-distribution, domain, derivative and restriction of $\overline{\mathscr{G}}(I)$-distributions — are defined unless isomorphism.

In other terms, the concepts associated with domains of distribution of 1st and 2nd species — $\widetilde{\mathscr{G}}(I)$ and $\overline{\mathscr{G}}(I)$-distributions, domain, derivative and restriction of $\widetilde{\mathscr{G}}(I)$ and $\overline{\mathscr{G}}(I)$-distributions — were not defined through the exhibition of specific mathematical objects, built at the heart of some mathematical theory (Analysis, for example), to which we attributed the denominations referred in Definition 5.10(c) and (d). Instead, the definitions in question characterize the referred concepts through a set of requirements, bundled up in the notions of $S$-space, abelian, surjective, and with identity $S$-space, $S$-space extension, strict, closed, locally closed, and coherent extensions, among others; in few words, they are axiomatic definitions.

In this perspective, as a consequence, the Definitions 5.10(b-1) and (b-2) of the notions of domain of distribution of 1st and 2nd species, respectively, also are axiomatic. Hence, strictly speaking, we must refer, for instance, to the ordered pairs $(\mathscr{C}(\mathbb{R}^n), \mathscr{D}'_f(\mathbb{R}^n))$ and $(\mathscr{C}(\mathbb{R}^n), \mathscr{D}'(\mathbb{R}^n))$, as being models (of the axiomatically defined notions) of the domains of distribution of 1st and 2nd species determined by $\mathscr{C}(\mathbb{R}^n)$, respectively, and not as we did earlier, referring to the pairs above as if they were themselves domains of distribution.

Due to the 1st and 2nd TESS, all these axiomatic notions introduced through Definition 5.10(b), (c) and (d) are consistent, i.e., they do have models, and they are categoric, that is, they admit, essentially, a single model. Hence, for example, the 2nd species domain of distribution determined by the $S$-space $\mathscr{C}(\mathbb{R}^n)$, that is, the axiomatic structure $(\mathscr{C}(\mathbb{R}^n), \widehat{\mathscr{C}}(\mathbb{R}^n))$ defined in 5.10(b-2) admits as a model the pair

$$\Big(\mathscr{C}(\mathbb{R}^n), \mathscr{D}'(\mathbb{R}^n)\Big),$$

while the 1st species domain of distribution determined by the same $S$-space $\mathscr{C}(\mathbb{R}^n)$, that is, the axiomatic structure $(\mathscr{C}(\mathbb{R}^n), \widehat{\mathscr{C}}(\mathbb{R}^n))$ defined in 5.10(b-1) has the pair

$$\Big(\mathscr{C}(\mathbb{R}^n), \mathscr{D}'_f(\mathbb{R}^n)\Big)$$

as one of its models. These models are, unless isomorphism, unique and, hence, for this reason and in this sense, we conclude that:

**(a)** the finite order distributions and its associated concepts — domain, derivative and restriction of finite order distributions — are, essentially, the $\widetilde{\mathscr{C}}(\mathbb{R}^n)$-distributions and its corresponding associated concepts, in the sense of, unless isomorphism, composing the only model of the axiomatic structure defined in 5.10(b-1) and (c) with $\mathscr{G}(I) = \mathscr{C}(\mathbb{R}^n)$;

**(b)** the distributions (of all orders, finite and infinite) and its associated concepts are, essentially, the $\overline{\mathscr{C}}(\mathbb{R}^n)$-distributions and its corresponding associate concepts, in the sense of, unless isomorphism, composing the only model of the axiomatic framework formulated by Definition 5.10(b-2) and (d) with $\mathscr{G}(I) = \mathscr{C}(\mathbb{R}^n)$.

Reporting ourselves to the Preliminaries 5.32, we see that conclusions (a) and (b) above are, precisely, what there we had established, informally, as the goals of this section.

# Schwartz' Distributions Axiomatics

## 5.39 Initial Remarks

As we saw in item 5.38, to each $S$-space $\mathscr{G}(I)$, abelian, surjective, with identity, coherent and whose family of $S$-groups is such that its groups are two by two disjoints, correspond two axiomatic structures, namely, the 1st and 2nd species domains of distribution determined by $\mathscr{G}(I)$ or, equivalently, the $\widetilde{\mathscr{G}}(I)$ and $\overline{\mathscr{G}}(I)$-distributions and its associated concepts, characterized by Definitions 5.10(b), (c) and (d).

Furthermore, regarding these axiomatic notions, we also saw that:

- the axiomatics formulated in items 5.12 and 5.17 are equivalents to the to axiomatic definition of $\widetilde{\mathscr{G}}(I)$-distributions and its associated concepts, given in 5.10(b-1) and (c), that, in turn, is categoric (due to the 1st TESS);
- the axiomatics described in 5.21 and 5.28 are equivalents to the axiomatic definition of $\overline{\mathscr{G}}(I)$-distributions and its associated concepts, given in 5.10(b-2) and (d), that, in turn, is categoric (due to the 2nd TESS);

Taking now into account the conclusions (a) and (b) obtained at the end of item 5.38, namely,

- the finite order distributions $S$-space, $\mathscr{D}'_f(\mathbb{R}^n)$, is, unless isomorphism, the only model of the 1st species domain of distribution determined by the $S$-space $\mathscr{C}(\mathbb{R}^n)$, or, equivalently, of the $\widetilde{\mathscr{C}}(\mathbb{R}^n)$-distributions axiomatic and its associated concepts;
- the distributions $S$-space, $\mathscr{D}'(\mathbb{R}^n)$, is, unless isomorphism, the only model of the 2nd species domain of distribution determined by the $S$-space $\mathscr{C}(\mathbb{R}^n)$, or, equivalently, of the $\overline{\mathscr{C}}(\mathbb{R}^n)$-distributions axiomatic and its associated concepts,

it then results that:

- the axiomatics described in 5.12 and 5.17, translated (particularized) to the case where $\mathscr{G}(I)$ is taken as the $S$-space $\mathscr{C}(\mathbb{R}^n)$, are equivalents, categoric and have the finite order distributions and its associated concepts as its unique model (unless isomorphism);

- the axiomatics formulated in items 5.21 and 5.28, specialized (translated) to the case where $\mathscr{G}(I) = \mathscr{C}(\mathbb{R}^n)$, are equivalents, categoric and have as model, unique unless isomorphism, the distributions (of all orders, finite and infinite) and its associated concepts.

In what follows, in the next four items, these particularizations, to the case where $\mathscr{G}(I) = \mathscr{C}(\mathbb{R}^n)$, of the $\widetilde{\mathscr{G}}(I)$ and $\overline{\mathscr{G}}(I)$-distributions axiomatics formulated in 5.12 and 5.21, respectively, as well as its simplified versions described in 5.17 and 5.28, are explicitly presented; due to obvious reasons, instead of the denominations "$\widetilde{\mathscr{C}}(\mathbb{R}^n)$-Distributions Axiomatic" and "$\overline{\mathscr{C}}(\mathbb{R}^n)$-Distributions Axiomatic" for the referred particularizations, we will use "Finite Order Distributions Axiomatic" and "Distributions Axiomatic", respectively.

## 5.40 Finite Order Distributions Axiomatic

- **Primitive Terms:** finite order distribution, domain, addition, derivative, and restriction of finite order distributions;

- **Precedent Theories:** Classical Logic, Set Theory, and $S$-Spaces Theory;

- **Axioms:** The ones formulated ahead, regarding the $S$-space

$$\mathscr{C}(\mathbb{R}^n) = \Bigg( \mathbb{C}\Big(\Gamma(\mathbb{R}^n)\Big) = \Big\{\mathbb{C}(\Omega) = \Big(C(\Omega), \partial(\Omega)\Big)\Big\}_{\Omega\in\Gamma(\mathbb{R}^n)}, \\ j\Big(\Gamma^2(\mathbb{R}^n)\Big), \\ \rho\Big(\Delta(\mathbb{R}^n)\Big) = \Big\{\rho_{(\Omega',\Omega)}\Big\}_{(\Omega',\Omega)\in\Delta(\mathbb{R}^n)}\Bigg).$$

**Axiom 1** Every function $f \in C(\Omega)$, for each $\Omega \in \Gamma(\mathbb{R}^n)$, is a finite order distribution.

**Axiom 2** To each finite order distribution, $\Lambda$, corresponds a single open set of $\mathbb{R}^n$ (i.e., an element of $\Gamma(\mathbb{R}^n)$), denominated the domain of $\Lambda$, in such a way that, if $\Lambda \in C(\Omega)$, then, $\Omega$ is the domain of $\Lambda$.

**Axiom 3** The addition is an operation that to each ordered pair $(\Lambda, \Lambda')$ of finite order distributions with the same domain $\Omega$, associates a single finite order distribution with the same domain $\Omega$, denoted by $\Lambda + \Lambda'$ and denominated the sum of $\Lambda$ with $\Lambda'$, in such a way that if $\Lambda, \Lambda' \in C(\Omega)$, the sum $\Lambda + \Lambda'$ is the one obtained through the addition of the group $C(\Omega)$.

**Notation.** We will denote by $D'_f(\Omega)$ the class of finite order distributions with domain $\Omega$. From Axiom 1 and Axiom 2 we have that $C(\Omega) \subseteq D'_f(\Omega)$.

**Axiom 4** The derivatives of finite order distributions with domain $\Omega$, whose class is denoted by $D_f(\Omega)$, are functions with domain and codomain equal to $D'_f(\Omega)$, one (unique) per multi-index $\alpha \in \mathbb{N}^n$ denoted by $D^\alpha_{\Omega-f} : D'_f(\Omega) \longrightarrow D'_f(\Omega)$, such that:

**(4-1)** $D^\alpha_{\Omega-f}(\Lambda + \Lambda') = D^\alpha_{\Omega-f}(\Lambda) + D^\alpha_{\Omega-f}(\Lambda')$ for every $\alpha \in \mathbb{N}^n$ and every $\Lambda, \Lambda' \in D'_f(\Omega)$;

**(4-2)** for each $\alpha \in \mathbb{N}^n$, $D^\alpha_{\Omega-f}(f) = \partial^\alpha_\Omega(f)$ for every $f \in C^{|\alpha|}(\Omega)$;

**(4-3)** the class $D_f(\Omega)$, with the usual operation of composition of functions, is a semigroup and the function

$$\begin{aligned} \mathrm{d}_{\Omega-f} : \partial(\Omega) &\longrightarrow D_f(\Omega) \\ \partial^\alpha_\Omega &\longmapsto \mathrm{d}_{\Omega-f}(\partial^\alpha_\Omega) := D^\alpha_{\Omega-f} \end{aligned}$$

is an isomorphism.

**Axiom 5** Let $\Omega \in \Gamma(\mathbb{R}^n)$ and $\partial^\alpha_\Omega \in \partial(\Omega)$. If $f \in C(\Omega)$ is such that $f \notin C^{|\alpha|}(\Omega)$, then, $D^\alpha_{\Omega-f}(f) \notin C(\Omega)$.

**Axiom 6** Let $\Omega \in \Gamma(\mathbb{R}^n)$ and $\Lambda \in D'_f(\Omega)$. There exist $\alpha \in \mathbb{N}^n$ and $f \in C(\Omega)$ such that

$$\Lambda = D^\alpha_{\Omega-f}(f)$$

(every finite order distribution is a derivative of some function $f \in C(\Omega)$).

**Axiom 7** The restrictions of finite order distributions are functions $\Theta^f_{(\Omega',\Omega)} : D'_f(\Omega) \longrightarrow D'_f(\Omega')$, one (unique) for each $(\Omega', \Omega) \in \Delta(\mathbb{R}^n)$, that associate to each finite order distribution with domain $\Omega$, a single finite order distribution with domain $\Omega' \subseteq \Omega$, in such a way that:

**(7-1)** $\Theta^f_{(\Omega',\Omega)}(\Lambda+\Lambda') = \Theta^f_{(\Omega',\Omega)}(\Lambda)+\Theta^f_{(\Omega',\Omega)}(\Lambda')$ for every $\Lambda, \Lambda' \in D'_f(\Omega)$;

**(7-2)** $\Theta^f_{(\Omega',\Omega)}(f) = \rho_{(\Omega',\Omega)}(f)$ for every $f \in C(\Omega) \subseteq D'_f(\Omega)$;

**(7-3)** $\Theta^f_{(\Omega'',\Omega')}(\Theta^f_{(\Omega',\Omega)}(\Lambda)) = \Theta^f_{(\Omega'',\Omega)}(\Lambda)$ for every $\Omega, \Omega', \Omega'' \in \Gamma(\mathbb{R}^n)$ such that $\Omega'' \subseteq \Omega' \subseteq \Omega$ and every $\Lambda \in D'_f(\Omega)$;

**(7-4)** $\Theta^f_{(\Omega',\Omega)}(D^\alpha_{\Omega-f}(\Lambda)) = D^\alpha_{\Omega'-f}(\Theta^f_{(\Omega',\Omega)}(\Lambda))$ for every $\alpha \in \mathbb{N}^n$ and every $\Lambda \in D'_f(\Omega)$.

## 5.41 Finite Order Distributions Axiomatic: Simplified Version

- **Primitive Terms:** finite order distribution, domain, derivative, and restriction of finite order distributions;
- **Precedent Theories:** Classical Logic, Set Theory, and $S$-Spaces Theory;
- **Axioms:** The ones formulated ahead, regarding the $S$-space of continuous functions

$$\mathscr{C}(\mathbb{R}^n) = \Big(\mathbb{C}\Big(\Gamma(\mathbb{R}^n)\Big), j\Big(\Gamma^2(\mathbb{R}^n)\Big), \rho\Big(\Delta(\mathbb{R}^n)\Big)\Big).$$

**Axiom 1** Every function $f \in C(\Omega)$, for each $\Omega \in \Gamma(\mathbb{R}^n)$, is a finite order distribution.

**Axiom 2** To each finite order distribution, $\Lambda$, corresponds a single open set of $\mathbb{R}^n$ (i.e., an element of $\Gamma(\mathbb{R}^n)$), denominated the domain of $\Lambda$, in such a way that, if $\Lambda \in C(\Omega)$, then, $\Omega$ is the domain of $\Lambda$.

**Notation.** We will denote by $D'_f(\Omega)$ the class of finite order distributions with domain $\Omega$. From Axiom 1 and Axiom 2 we have that $C(\Omega) \subseteq D'_f(\Omega)$.

**Axiom 3** The derivatives of finite order distributions with domain $\Omega$, whose class is denoted by $D_f(\Omega)$, are functions with domain and codomain equal to $D'_f(\Omega)$, one (unique) per multi-index $\alpha \in \mathbb{N}^n$ denoted by $D^\alpha_{\Omega-f} : D'_f(\Omega) \longrightarrow D'_f(\Omega)$, such that:

**(3-1)** for each $\alpha \in \mathbb{N}^n$, $D^\alpha_{\Omega-f}(f) = \partial^\alpha_\Omega(f)$ for every $f \in C^{|\alpha|}(\Omega)$;

**(3-2)** the class $D_f(\Omega)$, with the usual operation of composition of functions, is a semigroup and the function

$$\begin{aligned} \mathrm{d}_{\Omega-f} : \partial(\Omega) &\longrightarrow D_f(\Omega) \\ \partial^\alpha_\Omega &\longmapsto \mathrm{d}_{\Omega-f}(\partial^\alpha_\Omega) := D^\alpha_{\Omega-f} \end{aligned}$$

is an isomorphism.

**Axiom 4** Let $\Omega \in \Gamma(\mathbb{R}^n)$ and $\Lambda \in D'_f(\Omega)$. There exist $\alpha \in \mathbb{N}^n$ and $f \in C(\Omega)$ such that

$$\Lambda = D^\alpha_{\Omega-f}(f)$$

(every finite order distribution is a derivative of some function $f \in C(\Omega)$).

**Axiom 5** For $\Omega \in \Gamma(\mathbb{R}^n)$, $f, g \in C(\Omega)$ and $\alpha \in \mathbb{N}^n$, one has:

$$D^\alpha_{\Omega-f}(f) = D^\alpha_{\Omega-f}(g) \quad \text{if and only if} \quad f - g \in N(\partial^\alpha_\Omega).$$

**Axiom 6** The restrictions of finite order distributions are functions $\Theta^{f}_{(\Omega',\Omega)} : D'_f(\Omega) \longrightarrow D'_f(\Omega')$, one (unique) for each $(\Omega',\Omega) \in \Delta(\mathbb{R}^n)$, that associate to each finite order distribution with domain $\Omega$, a single finite order distribution with domain $\Omega' \subseteq \Omega$, in such a way that:

**(6-1)** $\Theta^{f}_{(\Omega',\Omega)}(f) = \rho_{(\Omega',\Omega)}(f)$ for every $f \in C(\Omega) \subseteq D'_f(\Omega)$;

**(6-2)** $\Theta^{f}_{(\Omega'',\Omega')}(\Theta^{f}_{(\Omega',\Omega)}(\Lambda)) = \Theta^{f}_{(\Omega'',\Omega)}(\Lambda)$ for every $\Omega, \Omega', \Omega'' \in \Gamma(\mathbb{R}^n)$ such that $\Omega'' \subseteq \Omega' \subseteq \Omega$ and every $\Lambda \in D'_f(\Omega)$;

**(6-3)** $\Theta^{f}_{(\Omega',\Omega)}(D^{\alpha}_{\Omega-f}(\Lambda)) = D^{\alpha}_{\Omega'-f}(\Theta^{f}_{(\Omega',\Omega)}(\Lambda))$ for every $\alpha \in \mathbb{N}^n$ and every $\Lambda \in D'_f(\Omega)$.

## 5.42 Distributions Axiomatic

- **Primitive Terms:** distribution, domain, addition, derivative, and restriction of distributions;
- **Precedent Theories:** Classical Logic, Set Theory, and $S$-Spaces Theory;
- **Axioms:** The ones formulated ahead, regarding the $S$-space of continuous functions

$$\mathscr{C}(\mathbb{R}^n) = \left(\mathbb{C}\left(\Gamma(\mathbb{R}^n)\right), j\left(\Gamma^2(\mathbb{R}^n)\right), \rho\left(\Delta(\mathbb{R}^n)\right)\right).$$

**Axiom 1** Every function $f \in C(\Omega)$, for each $\Omega \in \Gamma(\mathbb{R}^n)$, is a distribution.

**Axiom 2** To each distribution, $\Lambda$, corresponds a single open set of $\mathbb{R}^n$ (i.e., an element of $\Gamma(\mathbb{R}^n)$), denominated the domain of $\Lambda$, in such a way that, if $\Lambda \in C(\Omega)$, then, $\Omega$ is the domain of $\Lambda$.

**Notation.** $D'(\Omega)$ will denote the class of distributions with domain $\Omega$. From Axiom 1 and Axiom 2 we have that $C(\Omega) \subseteq D'(\Omega)$.

**Axiom 3** The addition is an operation that to each ordered pair $(\Lambda, \Lambda')$ of distributions with the same domain $\Omega$, associates a single distribution with the same domain $\Omega$, denoted by $\Lambda + \Lambda'$ and denominated the sum of $\Lambda$ with $\Lambda'$, in such a way that if $\Lambda, \Lambda' \in C(\Omega)$, the sum $\Lambda + \Lambda'$ is the one obtained through the addition of the group $C(\Omega)$.

**Axiom 4** The derivatives of distributions with domain $\Omega$, whose class is denoted by $D(\Omega)$, are functions with domain and codomain equal to $D'(\Omega)$, one (unique) per multi-index $\alpha \in \mathbb{N}^n$ denoted by $D^{\alpha}_{\Omega} : D'(\Omega) \longrightarrow D'(\Omega)$, such that:

**(4-1)** $D^\alpha_\Omega(\Lambda + \Lambda') = D^\alpha_\Omega(\Lambda) + D^\alpha_\Omega(\Lambda')$ for every $\alpha \in \mathbb{N}^n$ and every $\Lambda, \Lambda' \in D'(\Omega)$;

**(4-2)** for each $\alpha \in \mathbb{N}^n$, $D^\alpha_\Omega(f) = \partial^\alpha_\Omega(f)$ for every $f \in C^{|\alpha|}(\Omega)$;

**(4-3)** the class $D(\Omega)$, with the usual operation of composition of functions, is a semigroup and the function

$$\begin{aligned} \mathrm{d}_\Omega : \partial(\Omega) &\longrightarrow D(\Omega) \\ \partial^\alpha_\Omega &\longmapsto \mathrm{d}_\Omega(\partial^\alpha_\Omega) := D^\alpha_\Omega \end{aligned}$$

is an isomorphism.

**Axiom 5** Let $\Omega \in \Gamma(\mathbb{R}^n)$ and $\partial^\alpha_\Omega \in \partial(\Omega)$. If $f \in C(\Omega)$ is such that $f \notin C^{|\alpha|}(\Omega)$, then, $D^\alpha_\Omega(f) \notin C(\Omega)$.

**Axiom 6** The restrictions of distributions are functions $\Theta'_{(\Omega',\Omega)} : D'(\Omega) \longrightarrow D'(\Omega')$, one (unique) for each $(\Omega', \Omega) \in \Delta(\mathbb{R}^n)$, that associate to each distribution with domain $\Omega$, a single distribution with domain $\Omega' \subseteq \Omega$, in such a way that:

**(6-1)** $\Theta'_{(\Omega',\Omega)}(\Lambda + \Lambda') = \Theta'_{(\Omega',\Omega)}(\Lambda) + \Theta'_{(\Omega',\Omega)}(\Lambda')$ for every $\Lambda, \Lambda' \in D'(\Omega)$;

**(6-2)** $\Theta'_{(\Omega',\Omega)}(f) = \rho_{(\Omega',\Omega)}(f)$ for every $f \in C(\Omega) \subseteq D'(\Omega)$;

**(6-3)** $\Theta'_{(\Omega'',\Omega')}(\Theta'_{(\Omega',\Omega)}(\Lambda)) = \Theta'_{(\Omega'',\Omega)}(\Lambda)$ for every $\Omega, \Omega', \Omega'' \in \Gamma(\mathbb{R}^n)$ such that $\Omega'' \subseteq \Omega' \subseteq \Omega$ and every $\Lambda \in D'(\Omega)$;

**(6-4)** $\Theta'_{(\Omega',\Omega)}(D^\alpha_\Omega(\Lambda)) = D^\alpha_{\Omega'}(\Theta'_{(\Omega',\Omega)}(\Lambda))$ for every $\alpha \in \mathbb{N}^n$ and every $\Lambda \in D'(\Omega)$.

**Axiom 7** Let $\Omega \in \Gamma(\mathbb{R}^n)$. If $\Lambda \in D'(\Omega)$ and $x \in \Omega$, there exists $\Omega' \in \Gamma(\Omega)$ such that $x \in \Omega'$ and

$$\Theta'_{(\Omega',\Omega)}(\Lambda) = D^\alpha_{\Omega'}(f)$$

for some $\alpha \in \mathbb{N}^n$ and $f \in C(\Omega')$.

**Definition.** Let $\Omega \in \Gamma(\mathbb{R}^n)$ and $\overline{\Gamma}(\Omega) \subseteq \Gamma(\Omega)$ be a cover of $\Omega$. Let also

$$\left\{\Lambda_{\Omega'}\right\}_{\Omega' \in \overline{\Gamma}(\Omega)}$$

be a family of distributions indexed by $\overline{\Gamma}(\Omega)$ such that $\Lambda_{\Omega'} \in D'(\Omega')$ for every $\Omega' \in \overline{\Gamma}(\Omega)$. We will say that the family $\{\Lambda_{\Omega'}\}_{\Omega' \in \overline{\Gamma}(\Omega)}$ is **coherent** if and only if

$$\Theta'_{(\Omega'' \cap \Omega',\Omega')}(\Lambda_{\Omega'}) = \Theta'_{(\Omega'' \cap \Omega',\Omega'')}(\Lambda_{\Omega''})$$

for every $\Omega', \Omega'' \in \overline{\Gamma}(\Omega)$ such that $\Omega'' \cap \Omega' \neq \varnothing$.

**Axiom 8** Let $\Omega \in \Gamma(\mathbb{R}^n)$. If $\{\Lambda_{\Omega'}\}_{\Omega' \in \overline{\Gamma}(\Omega)}$ is a coherent family of distributions, then, there exists a single distribution $\Lambda \in D'(\Omega)$ such that

$$\Theta'_{(\Omega',\Omega)}(\Lambda) = \Lambda_{\Omega'} \quad \text{for every} \quad \Omega' \in \overline{\Gamma}(\Omega).$$

## 5.43 Distributions Axiomatic: Simplified Version

- **Primitive Terms:** distribution, domain, derivative, and restriction of distributions;
- **Precedent Theories:** Classical Logic, Set Theory, and $S$-Spaces Theory;
- **Axioms:** The ones formulated ahead, regarding the $S$-space of continuous functions

$$\mathscr{C}(\mathbb{R}^n) = \Big(\mathbb{C}\big(\Gamma(\mathbb{R}^n)\big), j\big(\Gamma^2(\mathbb{R}^n)\big), \rho\big(\Delta(\mathbb{R}^n)\big)\Big).$$

**Axiom 1** Every function $f \in C(\Omega)$, for each $\Omega \in \Gamma(\mathbb{R}^n)$, is a distribution.

**Axiom 2** To each distribution, $\Lambda$, corresponds a single open set of $\mathbb{R}^n$ (i.e., an element of $\Gamma(\mathbb{R}^n)$), denominated the domain of $\Lambda$, in such a way that, if $\Lambda \in C(\Omega)$, then, $\Omega$ is the domain of $\Lambda$.

**Notation.** $D'(\Omega)$ will denote the class of distributions with domain $\Omega$. From Axiom 1 and Axiom 2 we have that $C(\Omega) \subseteq D'(\Omega)$.

**Axiom 3** The derivatives of distributions with domain $\Omega$, whose class is denoted by $D(\Omega)$, are functions with domain and codomain equal to $D'(\Omega)$, one (unique) per multi-index $\alpha \in \mathbb{N}^n$ denoted by $D^\alpha_\Omega : D'(\Omega) \longrightarrow D'(\Omega)$, such that:

**(3-1)** for each $\alpha \in \mathbb{N}^n$, $D^\alpha_\Omega(f) = \partial^\alpha_\Omega(f)$ for every $f \in C^{|\alpha|}(\Omega)$;

**(3-2)** the class $D(\Omega)$, with the usual operation of composition of functions, is a semigroup and the function

$$\begin{aligned} \mathrm{d}_\Omega : \partial(\Omega) &\longrightarrow D(\Omega) \\ \partial^\alpha_\Omega &\longmapsto \mathrm{d}_\Omega(\partial^\alpha_\Omega) \coloneqq D^\alpha_\Omega \end{aligned}$$

is an isomorphism.

**Axiom 4** For $\Omega \in \Gamma(\mathbb{R}^n)$, $f, g \in C(\Omega)$ and $\alpha \in \mathbb{N}^n$ one has:

$$D^\alpha_\Omega(f) = D^\alpha_\Omega(g) \quad \text{if and only if} \quad f - g \in N(\partial^\alpha_\Omega).$$

**Axiom 5** The restrictions of distributions are functions $\Theta'_{(\Omega',\Omega)} : D'(\Omega) \longrightarrow D'(\Omega')$, one (unique) for each $(\Omega', \Omega) \in \Delta(\mathbb{R}^n)$, that associate to each distribution with domain $\Omega$, a single distribution with domain $\Omega' \subseteq \Omega$, in such a way that:

**(5-1)** $\Theta'_{(\Omega',\Omega)}(f) = \rho_{(\Omega',\Omega)}(f)$ for every $f \in C(\Omega) \subseteq D'(\Omega)$;

**(5-2)** $\Theta'_{(\Omega'',\Omega')}(\Theta'_{(\Omega',\Omega)}(\Lambda)) = \Theta'_{(\Omega'',\Omega)}(\Lambda)$ for every $\Omega, \Omega', \Omega'' \in \Gamma(\mathbb{R}^n)$ such that $\Omega'' \subseteq \Omega' \subseteq \Omega$ and every $\Lambda \in D'(\Omega)$;

**(5-3)** $\Theta'_{(\Omega',\Omega)}(D^\alpha_\Omega(\Lambda)) = D^\alpha_{\Omega'}(\Theta'_{(\Omega',\Omega)}(\Lambda))$ for every $\alpha \in \mathbb{N}^n$ and every $\Lambda \in D'(\Omega)$.

**Axiom 6** Let $\Omega \in \Gamma(\mathbb{R}^n)$. If $\Lambda \in D'(\Omega)$ and $x \in \Omega$, there exists $\Omega' \in \Gamma(\Omega)$ such that $x \in \Omega'$ and

$$\Theta'_{(\Omega',\Omega)}(\Lambda) = D^\alpha_{\Omega'}(f)$$

for some $\alpha \in \mathbb{N}^n$ and $f \in C(\Omega')$.

**Definition.** Let $\Omega \in \Gamma(\mathbb{R}^n)$ and $\overline{\Gamma}(\Omega) \subseteq \Gamma(\Omega)$ be a cover of $\Omega$. Let also

$$\left\{\Lambda_{\Omega'}\right\}_{\Omega' \in \overline{\Gamma}(\Omega)}$$

be a family of distributions indexed by $\overline{\Gamma}(\Omega)$ such that $\Lambda_{\Omega'} \in D'(\Omega')$ for every $\Omega' \in \overline{\Gamma}(\Omega)$. We will say that the family $\{\Lambda_{\Omega'}\}_{\Omega' \in \overline{\Gamma}(\Omega)}$ is **coherent** if and only if

$$\Theta'_{(\Omega'' \cap \Omega',\Omega')}(\Lambda_{\Omega'}) = \Theta'_{(\Omega'' \cap \Omega',\Omega'')}(\Lambda_{\Omega''})$$

for every $\Omega', \Omega'' \in \overline{\Gamma}(\Omega)$ such that $\Omega'' \cap \Omega' \neq \varnothing$.

**Axiom 7** Let $\Omega \in \Gamma(\mathbb{R}^n)$. If $\{\Lambda_{\Omega'}\}_{\Omega' \in \overline{\Gamma}(\Omega)}$ is a coherent family of distributions, then, there exists a single distribution $\Lambda \in D'(\Omega)$ such that

$$\Theta'_{(\Omega',\Omega)}(\Lambda) = \Lambda_{\Omega'} \quad \text{for every} \quad \Omega' \in \overline{\Gamma}(\Omega).$$

## 5.44 Final Remarks

In possession of the axiomatics above, that, as we saw, categorically define the Schwartz' distributions, the remarks made at the Prologue of this monograph, regarding the possibility of an approach to the theory of distributions, that is mathematically rigorous and appealing from the perspective of a mathematics user, become now more concrete. In fact, a fast look at the referred axiomatics reveals that one does not require any sophisticated mathematical background in order to understand them. On the contrary, even at the second year of a good (bachelor) undergraduate mathematics degree, a diligent student would be able to take any of the axiomatics above as a definition, and explore it through principles of logical inference.

Finally, regarding the purely mathematical aspect, the definition of the axiomatic structure of $S$-space, the formulation of the extension problems associated to it and, mainly, the obtainment of the theorems (1st and 2nd TESS) that solves them, are the points of this monograph that deserve to be highlighted and that, maybe, can be useful in other mathematical questions, such as they were when leading us to the, conceptually simple, distributions axiomatic when particularized to the case of the $S$-space of continuous functions, $\mathscr{C}(\mathbb{R}^n)$. After all, as well observed José Sebastião e Silva[31],

> *"... qu'on ne sait jamais si d'autres applications ne seront pas possibles."*

[31] Reference cited at the footnote 6 at page 16.